 \documentclass[reqno]{amsart}
\usepackage[top=1in, bottom=1in, left=1.35in, right=1.35in]{geometry}

 \usepackage{amsmath,amssymb}
 \usepackage{algpseudocode}
 \usepackage{algorithm}
 \usepackage{algorithmicx}
 \usepackage{caption}
 \usepackage{subcaption}
 \usepackage{amsthm}
 \numberwithin{equation}{section}
 \usepackage{amsfonts}
 \usepackage{graphicx}
\usepackage{forest}
 \usepackage{enumitem}
\usepackage{comment}

\newenvironment{widequote}
  {\begin{list}{}{\leftmargin=1.7em \rightmargin=1.7em}\item\relax}
  {\end{list}}

\usepackage{tikz}
\usetikzlibrary{
  arrows.meta,     
  positioning,     
  shapes.geometric,
  fit              
}

  \graphicspath{{converted_eps/}}

\usepackage[hypertexnames=false]{hyperref}

\usepackage[hang,flushmargin]{footmisc}
\usepackage{xcolor}
\usepackage{fancyhdr}
\usepackage{amsthm}
\usepackage{listings}
    \usepackage{cite}
    \usepackage{tabularx}  
    \usepackage{epsfig}
    \usepackage{float}
    \usepackage{listings}
    \usepackage{appendix}

 \theoremstyle{plain}
  \newtheorem{thm}{Theorem}
 \newtheorem{prop}{Proposition}[section]
 \newtheorem{lem}[prop]{Lemma}
 
 \theoremstyle{definition}
 \newtheorem{definition}[prop]{Definition}
 
  \newtheorem*{example*}{Example}
 \theoremstyle{remark}
 \newtheorem{remark}[prop]{\normalfont\bfseries Remark}
\newtheorem*{remark*}{\normalfont\bfseries Remark}      

\theoremstyle{definition}
\newtheorem{assumption}{Assumption}

\newcommand\appo[1]{ \td{\hat{#1}} }
\newcommand\appow[1]{ \wt{\wh{#1}} }

\newcommand\tbf[1]{\textbf{{#1}}}
\newcommand\mbf[1]{\mathbf{{#1}}}

\newcommand\mw[1]{\mathrm{{#1}}}
\newcommand \nlinf[1]{ {\| #1 \|}_{L^{\infty}} }
\newcommand \xxa[1]{ {\| #1 \|}_{X_1} }
\newcommand \xxb[1]{ {\| #1 \|}_{X_2} }

\let \hyr = \hyperref
         
\let \mf = \mathsf
\let \nolim = \nolimits         
\let \tts = \textstyle         
\let \tf = \tfrac
 \let\pa=\partial
 \let\al=\alpha
 \let\b=\beta
 \let\d=\delta
 \let\g=\gamma
 \let\e=\varepsilon
  
 \let \kp = \kappa
 \let\lam=\lambda

  \let \sig = \sigma
 \let\f=\frac
 \let\inf = \infty
 \let \les = \lesssim
  \let \gtr = \gtrsim
 \let\om=\omega
 
 \let \th = \theta
 \let \pr = \prime
 \let \vp = \varphi
 \let\G= \Gamma
\let\B = \Big
 \let\D=\Delta
 
 \let\S=\Sigma
 \let\Om=\Omega
 \let\td = \tilde
 \let\wt=\widetilde
 \let\wh=\widehat
 \let \olin = \overline

 \let\teq \triangleq
 
 \let\pa=\partial
 \let \bsh = \backslash
 \let \vs = \vspace
\let\sq = \square

\newcommand{\BT}{\mathbb{T}}
\newcommand{\BZ}{\mathbb{Z}}

\def \nmod { \mf{n}_{\hat W_2}}
\def \FM { \mf{FM}}

\def \phyy { \mf{phy} }
\newcommand\jan[1]{  { \la{#1}\ra } }
\def \cNF {\mathsf{NF}}
\def \fNS {\hat f_{\mathsf{NS}}}
\def \CCs { \mathsf{C} } 

\def \Ein {E_{\mf{in}}}
\def \Einm {E_{\mf{in}, m}}
\def \TTy { \mathsf{T}_1}
\def \numm {\mathsf{num}}

\def\uds{\underset}
\def \nlocc{ \mf{nloc}}

\let \dag = \dagger

\usepackage[normalem]{ulem}

\def \dap { \mf{p}}

  \def\cA{{\mathcal A}}
 \def\cB{{\mathcal B}}
 \def\cC{{\mathcal C}}
 
 \def\cE{{\mathcal E}}
 \def\cF{{\mathcal F}}

 \def\cI{{\mathcal I}}
 
 \def\cK{{\mathcal K}}
 \def\cL{{\mathcal L}}
 
 \def\cN{{\mathcal N}}

 \def\cR{{\mathcal R}}
 
 \def\cT{{\mathcal T}}

 \def\cX{{\mathcal X}}

 \def\cN{{\mathcal N}}

 \def\na{\nabla}
 \def\la{\langle}
 \def\ra{\rangle}
\def\lt{\left}
\def\rt{\right}
\def\one{\mathbf{1}}

 \def \weg{\wedge}

\def \tot {\mathsf{tot}}

\def\lemmainvelPII{2.3}
\def\seclinevo{3} 
 \def\secerridea{3.6}
 \def\secvelcomp{4}
 \def\secholfar{4.5}

\def\secudx{B.4}
 \def\secpiecepol{C.2}
\def\secestapprfar{C.3}

\def\seccutoffnear{{D.2}}
\def\appsolurep{C.1}
\def\secresid{C.4}
\def\secvellocest{4.7}

\def\secbddphol{8.4}
\def\secbddg{8.4.2}

\def\suppsecnondphol{8.5.2}
\def\suppsecnonest{8.5}
\def\suppseclinffar{8.6}
\def\suppsecholfar{8.7}
\def\suppsecbdother{8.9}

\def\secholconst{5}

\def\secCIIItonorm{8.5}

\def \appker{B}
\def \appexplcit{D}
\def \applinfest{E}
\def \secsolucomp{3.5}

\def \supptwosharp{5}
\def \supptwovel{6-7}
\def \supponenon{8}

\def \kin{ k_0 }
\def \ain{ \al_0 }

 \newcommand{\beq}{\begin{equation}}
 \newcommand{\eeq}{\end{equation}}
  \newcommand{\bal}{\begin{aligned} }
  \newcommand{\eal}{\end{aligned}}
    \newcommand{\bga}{\begin{gathered} }
  \newcommand{\ega}{\end{gathered}}
 \newcommand{\ben}{\begin{eqnarray}}
 \newcommand{\een}{\end{eqnarray}}
 \newcommand{\beno}{\begin{eqnarray*}}
 \newcommand{\eeno}{\end{eqnarray*}}

  \newcommand{\bseq}{\begin{subequations}}
 \newcommand{\eseq}{\end{subequations}}

 \newcommand{\Den}{\mathrm{Den}}

 \newcommand{\ee}{\mathbf{e}}

 \newcommand{\nn}{\mathbf{n}}
 \newcommand{\uu}{\mathbf{u}}

 \newcommand{\xx}{\mathbf{x}}
 \newcommand{\yy}{\mathbf{y}}
\newcommand{\zz}{\mathbf{z}}

 \newcommand{\R}{\mathbb{R}}
 
\newcommand{\UU}{\mathbf{U}}
\newcommand{\VV}{\mathbf{V}}

 \newcommand{\vom}{\boldsymbol\omega}

\newcommand{\sgn}{\mathrm{sgn}}
 \newcommand{\supp}{\mathrm{supp}}

 \author{Jiajie Chen and Thomas Y. Hou}
\address{Jiajie Chen, Courant Institute, New York University, New York, NY 10003.   
\newline 
Current address: Department of Mathematics, University of Chicago, Chicago, IL 60637.
\newline
Email: jiajiechen@uchicago.edu
}

\address{
Thomas Y. Hou, Applied and Computational Mathematics, Caltech, Pasadena, CA 91125.
\newline 
Email: hou@cms.caltech.edu
}

\title[Stable blowup of 2D Boussinesq equations]{Stable nearly self-similar blowup of the 2D Boussinesq and 3D Euler equations with smooth data I: Analysis}
\begin{document}

\begin{abstract}
 Inspired by numerical evidence of a potential 3D Euler singularity \cite{luo2014potentially,luo2013potentially-2}, we prove finite-time, nearly self-similar blowup of the 2D Boussinesq and 3D axisymmetric Euler equations with smooth initial data of finite energy and boundary. 
The proof encounters several essential difficulties. One of the essential difficulties is to control a number of nonlocal terms that do not seem to offer any damping effect. Another essential difficulty is that the strong advection normal to the boundary introduces a large growth factor for the perturbation when using weighted $L^2$ or $H^k$ estimates. We overcome these difficulties by combining weighted $L^\infty$ estimates and weighted $C^{1/2}$ estimates, and by developing sharp functional inequalities using the symmetry properties of the kernels and some techniques from optimal transport. Moreover, we decompose the linearized operator into a leading order operator and a finite-rank operator. We design the leading order operator to obtain sharp stability estimates. The contribution from the finite-rank operator to linear stability is estimated by constructing approximate space-time solutions. These ingredients enable us to establish the nonlinear stability of the approximate self-similar profile and to prove stable nearly self-similar blowup. 
\end{abstract}

 \setcounter{tocdepth}{1}
 \hypersetup{bookmarksdepth=2}

 \maketitle
 \tableofcontents

\section{Introduction}

The question whether the 3D incompressible Euler equations can develop a finite time singularity from smooth initial data of finite energy is one of the most outstanding open questions in the theory of nonlinear partial differential equations and fluid dynamics. The main difficulty is due to the presence of the vortex stretching term $\vom \cdot \nabla \uu$ in the vorticity equation:
\begin{equation}
   \vom_{t} + \uu \cdot \nabla \vom = \vom \cdot \nabla \uu,
   \quad  \uu = \na \times (-\D)^{-1} \vom,
  \label{eqn_eu_w}
\end{equation}
where $\vom = \nabla \times \uu$ is the \emph{vorticity vector} of the fluid, and the velocity vector $\uu$ is related to $\vom$ via the \emph{Biot-Savart law}. The velocity gradient $\nabla \uu$ formally has the same scaling as vorticity $\vom$. Thus the vortex stretching term has a nonlocal quadratic nonlinearity in terms of vorticity. 
However, the advection can lead to dynamic depletion of the nonlinear vortex stretching, which could prevent a finite time blowup, see e.g., \cite{constantin1996geometric,deng2005geometric,hou2006dynamic,constantin1986,lei2009stabilizing,hou2008dynamic}. Interested readers may consult the excellent surveys \cite{constantin2007euler,gibbon2008three,hou2009blow,kiselev2018,majda2002vorticity} and the references therein.

The blowup analysis presented in this paper is inspired by the computation of Luo-Hou \cite{luo2014potentially,luo2013potentially-2} in which they presented some convincing numerical evidence that the 3D axisymmetric Euler equations with smooth initial data and boundary develop a potential finite time singularity.
Inspired by the recent breakthrough of Elgindi \cite{elgindi2019finite} (see also \cite{elgindi2019stability})  on the blowup of the axisymmetric Euler equations without swirl for $C^{1,\alpha}$ initial velocity, we have proved asymptotically self-similar blowup of the 2D Boussinesq equations and the nearly self-similar blowup of the 3D axisymmetric Euler equations with $C^{1,\alpha}$ velocity and boundary in \cite{chen2019finite2}. The blowup analysis presented in \cite{chen2019finite2} takes advantage of the $C^{1,\alpha}$ velocity with a small $\al$ in an essential way and does not generalize to prove the Hou-Luo blowup scenario with smooth initial data.
The results presented in this paper provide the first rigorous proof of stable nearly self-similar blowup of the 2D Boussinesq and 3D Euler equations with smooth data and boundary.

The main results of this paper are stated by the two informal theorems below. 
A quantitative version of Theorem \ref{thm1a} is stated in Theorem \ref{thm:main} in Section \ref{sec:main_result}, and the full statement is given in Theorem \ref{thm:main_stab_full}.  The precise statement of Theorem \ref{thm1b} is given by Theorem \ref{thm:euler} in Section \ref{sec:euler}.

\begin{thm}\label{thm1a}
Let $\th$, $\uu$ and $\om$ be the density, velocity and vorticity in the 2D Boussinesq equations \eqref{eq:bous_sys}, respectively.
There is a family of smooth initial data $(\th_0, \om_0)$ with $\th_0(x,y)$ even,  $\om_0(x,y)$ odd in $x$, $\th_0(0, y) \equiv 0$, and finite-energy  velocity  $\uu_0$  such that the solution of the 2D Boussinesq equations develops a singularity in finite time $T<+\infty$. 
The blowup solution $(\th(t), \om(t))$ is nearly self-similar in the sense that $(\th(t), \om(t))$ with suitable dynamic rescaling is  close to an approximate blowup profile $ (\bar \th, \bar \om)$ up to time $T$. 
Moreover, within the same symmetry class, 
the blowup is stable for initial data $ (\na \th_0, \om_0)$  close to $(\na \bar \th, \bar \om)$ in some weighted $C^{1/2}$ norm.

 \end{thm}

\begin{thm}\label{thm1b}

Consider the 3D axisymmetric Euler equations in the cylinder $(r,z) \in [0, 1] \times \BT$. Let $u^\th$ and $\om^\th$ be the angular velocity and angular vorticity, respectively. The solution of the 3D Euler equations \eqref{eq:euler10}-\eqref{eq:euler20} develops a nearly self-similar blowup (in the sense described in Theorem \ref{thm1a}) in finite time for some smooth initial data 
$\om_0^{\th}$, $u_0^{\th}$ supported away from the symmetry axis $r=0$. The initial velocity has finite energy; $u_0^{\th}$ is even and $\om_0^{\th}$ is odd and periodic in $z$.  
Moreover, within the same symmetry class, the blowup is stable for $ (u^\th_0, \om^\th_0)$
close to the approximate blowup profile $({\bar u}^\th, {\bar \om}^\th)$ after proper rescaling subject to some constraint on the initial support size.

\end{thm}

\subsection{A novel framework of analysis with computer assistance}\label{sec:frame_novel}

One of our main contributions is to introduce a novel framework of analysis that enables us to obtain sharp stability estimates. In our analysis, we combine sharp functional inequalities, weighted energy estimates, a finite-rank perturbation method, and approximate space-time solutions constructed numerically with rigorous error control. We follow the framework in \cite{chen2019finite,chen2019finite2,chen2021HL} to establish finite time blowup of the 2D Boussinesq and 3D Euler equations by proving the nonlinear stability of an approximate steady state of the dynamic rescaling formulation. A very important first step is to construct an approximate steady state  with sufficiently small residual errors. 
We achieve this by decomposing the solution into a semi-analytic part capturing the far field behavior of the solution and a numerically computed part with compact support. The approximate steady state gives an approximate self-similar profile. See  Section \ref{sec:ASS} for further discussion. We remark that there have been some recent exciting developments in using a physics-informed neural network (PINN) to construct an approximate steady state of the 2D Boussinesq equations, see \cite{Tristan2022}.

Establishing linear stability of the approximate steady state is the most crucial step in our blowup analysis. One essential difficulty is that the advection normal to the boundary for smooth initial data introduces a large growth factor if we use weighted $L^2$ or $H^k$ energy estimates similar to \cite{chen2019finite,chen2019finite2,chen2021HL,elgindi2019finite}. To overcome this difficulty, we use weighted $L^\infty$ estimates and the outgoing property of the flow to extract damping from the local terms.
To close the energy estimates, we use a combination of the weighted $L^\infty$ estimates and the weighted $C^{1/2}$ estimates.

To estimate the nonlocal terms, we derive sharp $C^{1/2}$ estimates for $\na \uu$ using the symmetry properties of the kernels and some techniques from optimal transport \cite{villani2021topics}. We decompose the Biot-Savart law
into two parts. The main part captures the most singular part of the Biot-Savart law, and we apply the sharp functional inequalities for its $C^{1/2}$ estimate. 
The terms from the second part are more regular. We can approximate them by a finite-rank operator and obtain sharp estimates by constructing approximate space-time solutions with rigorous error control.

We use the 2D Boussinesq equations to give a high-level description of the linear stability analysis. 
Let $\bar{\omega}$, $\bar{\theta}$ be an approximate steady state. 
We denote $W = (\omega, \theta_x, \theta_y)$, decompose $W = \overline{W} +\widetilde{W}$ with $\overline{W} = (\bar{\omega},\bar{\theta}_x, \bar{\theta}_y)$, and denote by $\cL$ the linearized operator around $\overline{W}$ that governs the perturbation $\widetilde{W}$ in the dynamic rescaling formulation (see Section \ref{sec:linop}),
\begin{equation}
\label{eq:model_nloc_0}
\widetilde{W}_t = \cL(\widetilde{W}),
\end{equation}
where the coefficients of $\cL$ depend on the approximate steady state $\overline{W}$. 
We further decompose the linearized operator $\cL$ into a leading order operator $\cL_0$ plus a finite-rank perturbation $\cK$, i.e  $\cL = \cL_0 + \cK$. We design the leading order operator $\cL_0$ to obtain sharp stability estimates using weighted energy estimates and sharp functional inequalities.

In Part I of our paper, we perform the weighted energy estimates. In our analysis, we decompose 
$\widetilde{W} = \widetilde{W}_1 + \widetilde{W}_2$.  The first term $\widetilde{W}_1$ captures the main part of the perturbation, which is essentially governed by the leading order operator $\cL_0$ with a weak coupling to $\widetilde{W}_2$ through nonlinear interaction. The second term $\widetilde{W}_2$ captures the contribution from the finite-rank operator. We mainly perform stability analysis for $\widetilde{W}_1$ since $\wt W_2$ is driven by $\wt W_1$ (see \eqref{eq:model2_decoup_0}).

We use the following toy model to illustrate the main ideas by considering $\cK$ as a rank-one operator $\cK(\widetilde{W}) = a(x) P(\widetilde{W})$ for some operator $P$ satisfying 
(i) $P(\widetilde{W})$ is constant in space; (ii) $\|P(\widetilde{W})\| \leq c \|\widetilde{W}\| $. Given initial data $\widetilde{W}_0$, we decompose \eqref{eq:model_nloc_0} as follows 
\beq\label{eq:model2_decoup_0}
\bal
\partial_t \widetilde{W}_1(t) &= \cL_0 \widetilde{W}_1,  \quad \widetilde{W}_1(0) = \widetilde{W}_0 ,  \\
\partial_t \widetilde{W}_2(t) &= \cL \widetilde{W}_2 + a(x) P(\widetilde{W}_1(t)), \quad \widetilde{W}_2(0) = 0.
\eal
\eeq
It is easy to see that $\widetilde{W} = \widetilde{W}_1 + \widetilde{W}_2$ solves \eqref{eq:model_nloc_0} with initial data $\widetilde{W}_0$. The second part $\widetilde{W}_2$ is driven by the rank-one forcing term $a(x) P(\widetilde{W}_1(t))$. Using Duhamel's principle and the fact that $P(\widetilde{W}_1(t))$ is constant in space, we obtain 
\beq\label{eq:model_nloc4_0}
\widetilde{W}_2(t) = \int_0^t P(\widetilde{W}_1(s)) e^{\cL (t-s) } a(x) d s.
\eeq
Since the leading operator $\cL_0$ has the desired stability property by construction, $\widetilde{W}_1(t) = e^{\cL_0 (t)} \widetilde{W}_0$ decays in $L^{\inf}(\vp)$ ($\vp$ is a singular weight) and we can control $P(\widetilde{W}_1(s))$. By checking the decay of $e^{\cL (t)}a(x)$ in the energy space for large $t$, we can obtain the stability estimate of $\wt W_2$. A crucial idea in the estimate of $\wt W_2$ is to bridge the energy estimates and numerical PDEs via an approximate solution in space and time. 
Note that $e^{\cL (t) } a(x)$ is equivalent to solving the linear evolution equation $v_t = \cL(v)$ with initial data $v_0 = a(x)$. Due to the rapid decay of the linearized equation, we solve this initial value problem using a numerical scheme up to a modest time. The stability property of $\widetilde{W}_1$ allows us to control the numerical error in computing $e^{\cL (t) } a(x)$ and obtain  sharp stability estimates for $\widetilde{W}_2$. 
Without the flexibility to design and estimate \( \cL_0=\cL-\cK\), where \(\cK\) is a finite-rank operator, we would not be able to establish linear and nonlinear stability of the approximate self-similar profile.

We establish nonlinear stability using the stability lemma (Lemma \ref{lem:PDE_nonstab}), which depends on various constants in the estimates. The constants in the weighted energy estimates depend on the approximate steady state constructed numerically in Section \ref{sec:ASS} and the singular weights. The approximate steady state is represented using piecewise polynomials and explicit basis functions. We can obtain rigorous bounds for its high-order derivatives. These bounds in turn provide rigorous bounds for lower-order derivatives, pointwise values, and various integrals involving the approximate steady state by using standard numerical analysis.
In Part II of our paper \cite{ChenHou2023b}, we provide sharp and rigorous upper bounds for these constants and the residual error of the approximate steady state. In  Part II \cite[Section {\secvelcomp}]{ChenHou2023b}, we also estimate the velocity in the regular case by bounding various integrals with computer assistance. These sharp estimates of the constants enable us to prove that the inequalities in our stability lemma hold for our approximate self-similar profile. This completes the stability analysis of the approximate self-similar profile and proves the nearly self-similar blowup of the 2D Boussinesq and 3D Euler equations. See Section \ref{proof-sketch} for more discussions of our blowup analysis.

We remark that we have used fine properties of the approximate steady state,  especially its outgoing property \eqref{eq:outgoing}, in an essential way to establish its linear stability. Without obtaining relatively sharp upper bounds for the constants in the energy estimates that depend on the approximate self-similar profile, we would not have been able to apply Lemma \ref{lem:PDE_nonstab} to prove nonlinear stability. Thus, it seems quite difficult to prove stable blowup from smooth initial data without using any computer assistance.

We note that in the blow-up analysis for the critical nonlinear Schr\"odinger equation (see, e.g., \cite{merle2005blow}), one exploits an explicit ground state. Since no explicit ground state is available for the 3D Euler equations, the numerically constructed approximate steady state with a small residual error plays the same key role as the explicit ground state in blow-up analyses of other nonlinear PDEs, including the nonlinear Schr\"odinger equation \cite{merle2005blow} and the Keller--Segel system \cite{collot2022KSblow}.

To pass from the 2D Boussinesq to the 3D axisymmetric Euler equations, we follow the same ideas as in our previous work \cite{chen2019finite2} by controlling the support of the solution to be in a small region close to the boundary and away from the symmetry axis and performing localized 
elliptic estimates. The asymptotic scaling properties of the 3D Euler equations are the same as those of the 2D Boussinesq equations up to asymptotically small terms, after making appropriate changes of variables. We will provide some additional estimates to control these asymptotically small terms and prove the blowup of the 3D Euler equations.

\subsection{Comparison with finite co-dimension stability argument}

Our argument shares some similarities with the blowup analysis using a topological argument \cite{masmoudi2008blow} and finite co-dimension stability analysis; see e.g. \cite{merle2022blow, merle2022implosion2}. In the finite co-dimension stability analysis, one performs $H^k$ estimates on a linearized operator $\cL$ with large $k$ to obtain coercivity of $\cL - \cK$ for some compact operator $\cK$, and estimates $\cK$ using spectral theory and semigroup methods. 
A topological or fixed-point argument is then used to avoid potential unstable modes generated by $\cK$.

There are a few differences between this method and ours. First,
we do not require the energy space to be a Hilbert space. 
Second, we develop a constructive method to analyze the finite-rank operator $\cK$ (which approximates a compact operator)  by solving decoupled linear PDEs in space-time with rigorous error control. This yields explicit, computable bounds for $\cK$, 
whereas semigroup methods 
in \cite{merle2022blow, merle2022implosion2} produce \textit{implicit, uncomputable constants} that cannot be used to close the nonlinear stability estimates around a numerical profile. 
 Third, our method determines conditions on $\cK$ for \textit{full} stability. 
Moreover, if one attempts to establish nonlinear stability of 
 a numerical approximate self-similar profile in the $H^k$ norm with large $k$, 
 it would be extremely difficult to construct the approximate profile with a sufficiently small residual error in this norm.

\subsection{Review of literature}
There has been a lot of effort in studying 3D Euler singularities using various simplified models. 
Several 1D models, including the Constantin-Lax-Majda (CLM) model \cite{CLM85}, the De Gregorio (DG) model \cite{DG90}, the gCLM
model \cite{OSW08} and the Hou-Li model \cite{hou2008dynamic}, have been introduced to study the effect of advection and vortex stretching in 3D Euler. Singularity formation has been established for the CLM model in \cite{CLM85}, for the DG model with smooth data in \cite{chen2019finite} and with $C^{1-}$ data in \cite{chen2021regularity}, and for the gCLM model with various parameters in \cite{Elg17,chen2020singularity,chen2020slightly,chen2019finite,Elg19,Cor10}. In \cite{choi2014on}, the authors proved the blowup of the Hou-Luo model proposed in \cite{luo2014potentially}. In \cite{chen2021HL}, Chen-Hou-Huang proved the asymptotically self-similar blowup of the Hou-Luo model by extending the method of analysis for the finite time blowup of the De Gregorio model by the same authors in \cite{chen2019finite}. Inspired by their work on the vortex sheet singularity \cite{caflisch1989a}, Caflisch and Siegel have studied complex singularity for the 3D Euler equations, see \cite{caflisch1993,caflisch2009}.

In \cite{choi2015finite,kiselev2018finite}, 
the authors proposed several simplified models to study the Hou-Luo blowup scenario \cite{luo2014potentially,luo2013potentially-2} and established finite time blowup of these models. In these works, the velocity is determined by a simplified Biot-Savart law in a form similar to the key lemma in the seminal work of Kiselev-Sverak \cite{kiselev2013small}. 
 In \cite{elgindi2017finite,elgindi2018finite}, Elgindi and Jeong proved finite time blowup for the 2D Boussinesq and 3D axisymmetric Euler equations in a domain with a corner using $\mathring{C}^{0,\al}$ data.
There has been some recent progress in searching for potential Euler and Navier-Stokes singularities in the interior domain, see  \cite{HouHuang2021,HouHuang2022,Houeuler2021,Hounse2021}. 

There have been some interesting recent results on the potential instability of the Euler blowup solutions, see \cite{vasseur2020blow,lafleche2021instability}. In a recent paper \cite{Chen2022}, we showed that the blowup solutions of the 2D Boussinesq and 3D Euler equations with $C^{1, \alpha}$ velocity considered in \cite{elgindi2019finite,chen2019finite2} are also unstable using the notion of stability introduced in \cite{vasseur2020blow,lafleche2021instability}. The blowup analysis in \cite{elgindi2019finite,chen2019finite2} is based on the stability of a self-similar blowup profile using the dynamic rescaling formulation. In comparison, the linear stability in  \cite{vasseur2020blow,lafleche2021instability} is performed by directly linearizing the 3D Euler equations around a particular blowup solution with a fixed blowup time $T$ in the original physical variables.

\vs{0.05in}
\paragraph{\bf Organization} The rest of the paper is organized as follows. Sections \ref{sec:lin}--\ref{sec:EE} are devoted to the blowup analysis for 2D Boussinesq and Section \ref{sec:euler} is devoted to the blowup analysis for 3D Euler. In Section \ref{sec:lin}, we discuss the main ingredients and ideas in establishing linear stability of an approximate profile using various simplified models. In Section \ref{sec:sharp}, we develop sharp  H\"older estimates using optimal transport. In Section \ref{sec:finite_rank}, we introduce the $L^{\inf}$-based finite-rank perturbation method. Section \ref{sec:EE} is devoted to energy estimates and Section \ref{sec:ASS} is devoted to the construction of an approximate self-similar profile using the dynamic rescaling formulation. 
Some technical estimates, derivations, parameters, and the full stability inequalities are deferred to the appendices.

\section{Main ideas for stability analysis}\label{sec:lin}

In this section, we outline the main ingredients in our stability analysis. We mainly focus on the 2D Boussinesq equations and generalize the analysis to 3D Euler in Section \ref{sec:euler}.

We begin by introducing the dynamic rescaling formulation (Section \ref{sec:dyn_res})
and stating an abbreviated version of the main stability and blowup result (Section \ref{sec:main_result}). The overall structure of the proof is outlined in Section~\ref{proof-sketch}, including the outgoing flow property and finite-rank perturbation.
Properties of the approximate steady state are
described in Section~\ref{sec:property}. 
The remaining subsections explain the main ideas behind the linear stability estimates and the analytic low-rank corrections used to achieve cubic vanishing of the perturbation and error.

Denote by $\om^{\th}$, $u^{\th}$ and $\phi^\th$  the angular vorticity, angular velocity, and angular stream function, respectively.
The 3D axisymmetric Euler equations are given below:
\beq\label{eq:euler10}
\pa_t (ru^{\th}) + u^r (r u^{\th})_r + u^z (r u^{\th})_z = 0, \quad 
\pa_t (\f{\om^{\th}}{r}) + u^r ( \f{\om^{\th}}{r} )_r + u^z ( \f{\om^{\th}}{r})_z = \f{1}{r^4} \pa_z( (r u^{\th})^2 ),
\eeq
where the radial velocity $u^r$ and the axial velocity $u^z$ are given by the Biot-Savart law:
\beq\label{eq:euler20}
-(\pa_{rr} + \f{1}{r} \pa_{r} +\pa_{zz}) {\phi^\th} + \f{1}{r^2} {\phi^\th} = \om^{\th}, 
\quad  u^r = -\phi^\th_z \, , \quad u^z = \phi^\th_r + \f{1}{r} {\phi^\th} , 
\eeq
with the no-flow boundary condition ${\phi^\th}(1, z ) = 0$ on the solid boundary $r = 1$ 
and a periodic boundary condition in $z$. For 3D axisymmetric Euler blowup at the boundary $r=1$, 
we know that it has scaling properties asymptotically the same as those  of the 2D Boussinesq equations \cite{majda2002vorticity,chen2019finite2}. 
We therefore study the 2D Boussinesq equations on the upper half space:
\bseq\label{eq:bous_sys}
\begin{align}
\om_t +  \uu \cdot \na \om  &= \th_{x},  \label{eq:bous1}\\
\th_t + \uu \cdot  \na \th & =  0 , \label{eq:bous2} 
\end{align}
where the velocity field $\uu = (u , v)^T : \R_+^2 \times [0, T) \to \R^2$ is determined via the Biot-Savart law
\beq\label{eq:biot}
 - \D \phi = \om , \quad  u =  - \phi_y , \quad v  = \phi_x,
\eeq
\eseq 
where $\phi$ is the stream function with the no-flow boundary condition $\phi(x, 0 ) = 0$ at $y=0$.

\subsection{Dynamic rescaling formulation}\label{sec:dyn_res}

Following \cite{chen2019finite,chen2019finite2,chen2021HL}, we use the dynamic rescaling formulation. 
\footnote{
The dynamic rescaling formulation was introduced in \cite{mclaughlin1986focusing,landman1988rate} and is closely related to the modulation technique. 
It has been used in many singularity formation problems, including nonlinear
Schr\"odinger equations \cite{kenig2006global,merle2005blow}, compressible Euler and NSE equations \cite{buckmaster2019formation,buckmaster2019formation2,merle2022implosion2}, incompressible fluids \cite{chen2019finite2,elgindi2019finite}, and related
models \cite{chen2019finite,chen2020slightly,chen2021regularity,chen2020singularity}. 
}
Let $ \om(x, t), \th(x,t) , \uu(x, t)$ be the solutions of \eqref{eq:bous_sys}. 
It is easy to show that 
\beq\label{eq:rescal1}
\bal
  \td{\om}(x, \tau) &= C_{\om}(\tau) \om(   C_l(\tau) x,  t(\tau) ), \quad   \td{\th}(x , \tau) = C_{\th}(\tau)
  \th( C_l(\tau) x, t(\tau)),  \\
    \td{\uu}(x, \tau) &= C_{\om}(\tau)  C_l(\tau)^{-1} \uu(C_l(\tau) x, t(\tau)) , 
\eal
\eeq
are the solutions to the dynamic rescaling equations
 \beq\label{eq:bousdy0}
\bal
\td{\om}_{\tau}(x, \tau) + ( c_l(\tau) \xx + \td{\uu} ) \cdot \na \td{\om}  &=   c_{\om}(\tau) \td{\om} + \td{\th}_x , \qquad 
\td{\th}_{\tau}(x , \tau )+ ( c_l(\tau) \xx + \td{\uu} ) \cdot \na \td{\th}  = c_{\th} \td \th,
\eal
\eeq
where $\td \uu = (\td u, \td v)^T = \na^{\perp} (-\D)^{-1} \td{\om}$, $\xx = (x, y)^T$, 
\beq\label{eq:rescal2}
\bal
  C_{\om}(\tau) = \exp\big( \int_0^{\tau} c_{\om} (s)  d s  \big), \quad  C_l(\tau) = \exp\big( \int_0^{\tau} -c_l(s) ds    \big) , \quad  C_{\th}(\tau)  =  \exp\big( \int_0^{\tau} c_{\th} (s)  d s \big),
\eal
\eeq
$  t(\tau) = \int_0^{\tau} C_{\om}(s ) d s $ and the rescaling parameters $c_l(\tau), c_{\th}(\tau), c_{\om}(\tau)$ satisfy \cite{chen2019finite2}
\beq\label{eq:rescal3}
c_{\th}(\tau) = c_l(\tau ) + 2 c_{\om}(\tau).
\eeq

We have the freedom to choose the time-dependent scaling parameters $c_l(\tau)$ and $c_{\om}(\tau)$ 
 according to some normalization conditions. 
 \footnote{These two parameters relate to the two-parameter family of scaling symmetries in Boussinesq and 3D Euler.
 } 
 After we choose them, the resulting rescaling equations \eqref{eq:bousdy0} are equivalent to the original equation via \eqref{eq:rescal1}, \eqref{eq:rescal2}.

To simplify notation, we still use $t$ for the rescaled time in \eqref{eq:bousdy0} and write $\td \om, \td \th$ as $\om, \th$
\beq\label{eq:bousdy1}
\bal
&\om_t + (c_l \xx + \uu) \cdot \na \om = \th_x + c_{\om} \om ,\quad  \th_t + (c_l \xx + \uu)\cdot \na \th =  c_{\th} \th .
\eal
\eeq
where $\uu= (u,v)$. Following \cite{chen2021HL}, we impose the following normalization conditions on $c_{\om}, c_l$
\beq\label{eq:normal}
c_l = 2 \f{\th_{xx}( 0) }{\om_x( 0 ) }, \quad c_{\om} = \f{1}{2} c_l + u_x(0), \quad c_{\th} = c_l + 2 c_{\om} .
\eeq
We consider initial data with the following symmetries 
\beq\label{eq:Sym}
   \om(x, y) \mw{ \ is \ odd \ in \  } x,  
  \quad  \th(x, y) \mw{ \ is \ even \ in \ } x, 
   \quad \th(0, y) \equiv   0 , 
    \tag{Sym}
\eeq
which are preserved by \eqref{eq:bous_sys} and \eqref{eq:bousdy1}-\eqref{eq:normal}. For smooth data, 
conditions \eqref{eq:normal} enforce 
\beq\label{eq:normal1}
\theta_{xx}(t, 0)=\theta_{0, xx}(  0), \quad \omega_x(t, 0)=\omega_{0, x}( 0)
\eeq
for all time, where $\th_0, \om_0$ are initial data to \eqref{eq:bousdy1}.
\footnote{We can derive the ODEs of $\theta_{xx}(t,0)$ and $\om_x(t, 0)$
\[
\tf{d}{dt} 
\om_x(t, 0) = (c_{\om} - c_l - u_x(0)) \om_x(t,0) + \th_{xx}(t,0), 
\quad \tf{d}{dt} \th_{xx}(t, 0) =  (c_{\th} - 2 (c_l + u_x(0))) \th_{xx}(t, 0),
\]
where we use $v|_{y =0} =0, v_x(t, 0) = 0$. Under the conditions \eqref{eq:normal}, the right hand sides vanish.
}

\subsection{Main stability result}\label{sec:main_result}

In this section, we state an abbreviated version of the main stability and blowup result for the 2D Boussinesq equations. The full stability theorem is stated later in  Theorem~\ref{thm:main_stab_full} in Section~\ref{sec:EE_stab_thm}, near the end of the energy estimates.

\begin{thm}[\bf Nonlinear stability and blowup -- abbreviated version]\label{thm:main}
Let $(\bar{\th},\bar{\om}, \bar \uu,  \bar c_l, \bar c_{\om})$ be the approximate self-similar profile constructed in Section~\ref{sec:ASS} and $E_* = 5 \cdot 10^{-6}$.  Assume that the initial data $(\om_0, \th_0)$ satisfy the  symmetry \eqref{eq:Sym} with $ \om_0  \in  C^{k, \al}, \th_0  \in C^{k + 1, \al}$ for some $k \geq 5, \al \in (0, 1)$. Assume also that $(\om_0, \th_0)$ are sufficiently close to the profile $(\bar \om, \bar \th)$ in the energy $\Ein$, introduced in 
\eqref{energy} in Section~\ref{sec:EE_stab_thm}, which controls weighted $C^{1/2}$ norms of $\om_0$ and $\nabla \th_0$. 

Then the corresponding solution $(\om, \th)$ to the dynamic rescaling equations \eqref{eq:bousdy1} exists globally in the rescaled time and remains close to the approximate profile for any $t \geq 0$:
\[
 \| \om(t) - \bar \om \|_{L^{\inf}}, \  \| \th_x(t) - \bar \th_x \|_{L^{\inf}} ,
\  \| \th_y(t) - \bar \th_y \|_{L^{\inf}} < 200 E_* , \  | u_x(t, 0 ) - \bar u_x( 0)| , \   | \bar c_{\om} - c_{\om}(t)| < 100 E_*.
\]
 Consequently, we can choose smooth initial data $\om_0, \th_0 \in C_c^{\inf}$ with finite energy
$\uu_0 \in L^2$ such that the corresponding solution of the physical equations \eqref{eq:bous_sys} blows up in finite time $T$.

\end{thm}

\vs{0.05in}
\paragraph{\bf{Nearly self-similar blowup}}

We give an \emph{informal} description of the nearly self-similar blowup solution 
and its blowup rate. 
The vorticity  $\om_{\phyy}$  in the physical space \eqref{eq:bous_sys} has the  form 
\[
  \om_{\phyy}(x, t(\tau)) = C_{\om}^{-1}(\tau) \om_{\mf{ss}}( C_l(\tau)^{-1} x, \tau ),
  \quad || \om_{ss}(\tau) - \bar \om ||_{L^{\infty}} \ll 1 ,
\]
where $\om_{\mf{ss}}$ is the self-similar variable, corresponding to $\td \om$ in \eqref{eq:rescal1}. 
We can generalize \eqref{eq:rescal2} to $C_{\om}(\tau) = C_{\om}(0) \exp( \int_0^{\tau} c_{\om}(s) ds ) , C_{\th}(\tau) = C_{\om}(0)^2 \exp( \int_0^{\tau} c_{\th}(s) ds )$  with an extra factor $C_{\om}(0)$, while keeping the same formulas for $C_l,  t(\tau)$ in  \eqref{eq:rescal2}. 
 Using estimates in Theorem \ref{thm:main} and $c_l(\tau) \equiv \bar c_l$ 
 satisfied by the blowup solution \eqref{eq:normal},\eqref{eq:normal1},  \emph{informally}, we 
 approximate the relation \eqref{eq:rescal1}, \eqref{eq:rescal2}:
 \[
 \bal
  || \om_{\phyy}(t(\tau) )||_{L^{\infty}} & \approx C_{\om}^{-1}(\tau) \| \bar \om \|_{L^{\infty}} \, , 
 & \quad C_{\om}(\tau) & \approx C_{\om}(0) e^{ \bar c_{\om} \tau} \, ,
 & \quad C_l(\tau) & =  e^{ - \bar c_l \tau} \, , \\
 T = t(\infty) & \approx \f{ C_{\om}(0) }{ |\bar c_{\om}|}
\approx  \f{ \| \bar \om \|_{L^{\infty}} }{  \| \om_{\phyy}(0) \|_{L^{\infty}}   |\bar c_{\om}| } \, ,
& \quad T - t(\tau) & \approx \f{ C_{\om}(\tau) }{ |\bar c_{\om}| } \, , \quad 
 &  C_l(\tau)^{-1} & \approx C (T- t(\tau))^{ \bar c_l / \bar c_{\om}} \, ,   \\
\eal
 \]
 which implies 
 \[
   \om_{\phyy}(x, t(\tau) ) \approx  \f{1 }{ ( T-t(\tau) ) |\bar c_{\om}|   } 
 \bar \om\big( \f{ C \cdot x }{ (T- t(\tau))^{ -\bar c_l / \bar c_{\om}} } \big)  \, , 
 \]
 with $-\bar c_l / \bar c_{\om} \approx 2.92 $ (see \eqref{eq:bous_expo}),  for some $C>0$ depending on $\bar c_l, \bar c_{\om},  C_{\om}(0)$. 
 Here, $\approx$ denotes an  approximate relation.  The blowup time is approximately 
 proportional to $\| \om_{\phyy}(0) \|_{L^{\infty}}^{-1}$. 
 The exact relation can be inferred from Theorem \ref{thm:main}, \eqref{eq:rescal1}, \eqref{eq:rescal2}.  If instead $(\bar \om, \bar \th, \bar c_l, \bar c_{\om})$ solves the profile equations exactly, i.e. \eqref{eq:bousdy1} without $\pa_t \om, \pa_t \th$ terms, 
and $\om_{\mf{ss}}(\cdot, 0)= \bar \om, 
\th_{\mf{ss}}(\cdot, 0) = \bar \th $,  the above approximations become \emph{exact}. Since we only prove that $\om_{\mf{ss}}(\tau)$ remains sufficiently close to $\bar \om$ and do not prove convergence of $\om_{\mf{ss}}(\tau)$ as $\tau \to \infty$, Theorem \ref{thm:main} does not imply \textit{asymptotically} self-similar blowup.

\begin{remark}[Regularity and singular weights]

The energies $\Ein$ for the initial data and $E_4(t)$, defined in \eqref{energy4}, involve singular weights near the origin, such as $|x|^{-a}$ with $a>2$. In Section~\ref{sec:finite_pertb}, we introduce a finite-rank decomposition whose main component, denoted by $W_1$, agrees with the full perturbation $W$ initially, which vanishes $ O(|x|^3)$ near $x=0$. The local well-posedness argument in Section~\ref{sec:reg_van} shows that, the conditions $\nabla^i W_1(t,0)=0$ are preserved for $i\leq 2$. Thus $W_1$ vanishes to cubic order 
$O(|x|^3)$ near the origin, which makes these singularly weighted norms finite and justifies their use in the stability estimates.

\end{remark}

\subsection{The main steps in the proof of Theorem \ref{thm:main}}
\label{proof-sketch}

We follow the framework in \cite{chen2019finite,chen2019finite2,chen2021HL} to establish finite time blowup by proving the nonlinear stability of an approximate steady state to \eqref{eq:bousdy1}. 
The proof of  Theorem \ref{thm:main} combines five main ingredients: the construction of an
approximate steady state, weighted coercivity estimates, sharp weighted nonlocal estimates, a finite-rank perturbation method for the linear stability, and nonlinear stability estimates.

In this paper, we use the bar notation, $\bar{\cdot}$, for
the approximate steady state and its related quantities, e.g., $\bar \om$, $\bar \th$, $\bar c_{\om}, \bar c_l$. 
We use boldface to denote vectors, e.g., $\xx=(x,y)$, $\uu=(u,v)$, $\UU=(U,V)$. We use the ``hat'' notation for approximate quantities, such as the approximate space-time solutions $\hat F_i$, $\hat W_2$ and the finite-rank approximations of velocity $\hat \uu$ in Sections \ref{sec:finite_pertb}, \ref{sec:appr_vel}.

\vs{0.05in}

\paragraph{\bf Dependency tree of the proof} 
To clarify the structure of the proof, we present the dependency tree of the proof in Figure 
\ref{fig:tree}. In this figure, \textit{Thm, Lem, App, Sec, P1, P2, Supp1, Supp2} are short for Theorem, Lemma, Appendix, Section, Paper I (the current paper), Paper II \cite{ChenHou2023b}, the Supplementary material for Paper I \cite{ChenHou2023aSupp} and Paper II \cite{ChenHou2023bSupp}, respectively.

\begin{figure}[t]
\centering
\begin{forest}
 for tree ={l sep=3pt,
    s sep = 1pt,
    grow' = east,
    parent anchor=east,
    child anchor=west,
    }    
      [Thm \ref{thm:main}: proved \\ by  
      nonlinear \\ stability \\ Lem \ref{lem:PDE_nonstab}{,} P1  \\
        \& inequalities \\ in App \ref{sec:ineq}{,} P1, text width= 2.4cm  
               [Linear \\
               stability,   text width=1.7cm    
                [Sharp H\"older estimates: \\ Sec \ref{sec:sharp} \& App \ref{app:sharp}{,} P1, text width= 4cm
                      [ Compute the sharp \\ constant:
                          Sec {\supptwosharp} Supp2, text width= 4cm
                      ]
                  ]
                 [
                     Finite-rank perturbation{:} \\ Sec \ref{sec:finite_rank}{,} P1,  text width=4cm
                 ]
                 [  
                     Linear energy estimates{:}\\ Sec \ref{sec:EE_intro}-\ref{sec:rank1}{,} P1, text width=4cm
                 ]
                 [ 
                     Estimate nonlocal terms: 
                     Sec \ref{sec:u_comp_idea} (ideas){,} P1  \\ 
                        \& Sec {\secvelcomp} \& App {\appker}{,} P2, text width=4cm
                      [ A few more similar cases \\
                      and explicit formulas: Sec {\supptwovel} Supp2, text width=4cm]
                  ]
              ]   
              [Nonlinear \\ estimates,  text width=1.6cm 
                 [Construct finite-rank \\
                    part $\hat W_2${,} Sec \ref{sec:finite_pertb}-\ref{sec:appr_vel}{,}  \\
                    Estimate $\hat W_2${,} Sec \ref{sec:W2}{,} P1 
                    ,text width=4cm          
                        [Construct approximate  \\ space-time solution $\hat F_i$
                         in $\hat W_2${:}  \\ 
                        Sec {\seclinevo} \& App {\appexplcit}{,} P2, 
                        text width=4.7cm] 
                       [Estimate residual  
                       error \\ for $\hat W_2$: 
                        Sec {\secerridea} (ideas){,} \\
                        App {\secresid}{,} {\appexplcit} \& {\applinfest}{,} P2, text width=4.7cm]                            
                  ]
                  [ Estimate nonlinear terms:
                    Sec \ref{sec:rank1}-\ref{sec:non}{,} P1, text width=4cm
                        [Estimate similar nonlinear
                          \& error terms: Sec {\supponenon} Supp1, text width=4.8cm
                        ]
                  ]    
              ]              
              [Approximate \\ profile  \\ , text width= 2 cm                                                                          
                    [Construct profile{:}  Sec \ref{sec:ASS} \\ (methods){,} P1{;} 
                    App {\appsolurep} \\ \& App {\appexplcit}   (estimates){,} P2, text width=4.9cm ] 
                    [Estimate profile: \\
                     App {\secpiecepol}{,} {\secestapprfar}{,} {\appexplcit}{,} P2, text width=4.9cm]    
                    [Estimate residual error{:} 
                    Sec {\secerridea} (ideas){,} 
                    App {\secresid}{,} {\appexplcit} \& {\applinfest}{,} P2, 
                    text width=4.9cm
                    ]
              ]                                         
        ]      
\end{forest}
\vspace{-0.05in}
\captionsetup{width=0.9\textwidth} 
 \captionof{figure}{
Level 2: main steps. Level 3: main estimates and methods. Level 4: compute and estimate explicit functions or constants, and additional similar estimates. 
}
\label{fig:tree}

\end{figure}

Below, we briefly outline the steps of the proof and ideas. Each step is relatively self-contained.

\vs{0.05in}
\paragraph{\bf Step 1: Approximate steady state}\label{para:step1}

We construct the approximate steady state by numerically solving \eqref{eq:bousdy1},\eqref{eq:normal}, with details in Section \ref{sec:ASS}. Its properties are discussed in Section \ref{sec:property}, and the residual error estimate is given in Part II \cite[Appendix {\secresid}]{ChenHou2023b}.

\vs{0.05in}
\paragraph{\bf Step 2: Weighted coercivity estimates of local terms}\label{para:step2}

The key stability mechanism arises from the \textit{outgoing property} \eqref{eq:outgoing} of the flow, discussed in Section \ref{sec:analy_para}:
\beq\label{eq:sec2_outgo}
  \bar c_l x + \bar u(x, y) \geq c_1 x,  \quad c_1 \approx 0.47, \quad  \bar c_l y + \bar v(x, y) \geq c_2 y, \quad c_2 \approx  3, 
\eeq
for all $(x, y) \in \R^2_{++}$.
 This property means that the perturbation is transported from the origin to spatial infinity, generating a strong stability mechanism.  
We exploit this mechanism through singularly weighted $L^\infty$-type estimates:
the commutator between the weight and the transport operator generates damping terms for the local part of the linearized operator,
 which can be made \textit{arbitrarily strong}
in any bounded region away from the origin; see Section \ref{sec:outgo}.

To close the top-order estimate of the nonlocal terms, we use 
weighted $C^{1/2}$ estimates and derive similar damping terms. See Sections \ref{sec:model_local} and \ref{sec:idea_EE_calpha} for further discussion.

\vs{0.05in}

Having identified the local damping mechanism, we turn to 
weighted $L^{\infty}$ and $C^{1/2}$ nonlocal estimates. The most singular part is estimated in \hyr[para:step3]{Step~3} and the more regular parts in \hyr[para:step4]{Step~4}.

\vs{0.05in}

\paragraph{\bf Step 3: Sharp H\"older estimates}\label{para:step3}
\footnote{
Although we have standard $C^{\al}$ estimates for the Riesz transform $\na \na^{\perp} (-\D)^{-1}  \om$, the associated constants are usually not given explicitly and are not sharp enough for our purposes. 
  Here, we do not assume that $\al$ is small.
  }
A key difficulty is to control the zero-order singular integral operator $\na \uu=\na\na^\perp(-\D)^{-1}\om$ in singular weighted norms. We observe crucially that we can use optimal transport techniques 
\footnote{
Optimal transport has been applied to establish many sharp functional inequalities and to study functional inequalities in detail. We refer the reader to Villani's book \cite{villani2021topics} for background and further discussion.
}
to establish sharp $C^{1/2}$ estimates for $\na \uu$ with localized kernels, e.g. $u_x(x, a, b)$ defined in  \eqref{eq:ux_local}, which captures the most singular part in the nonlocal estimates. We refer to Section \ref{sec:sharp} and Appendix \ref{app:sharp} for details.
The constants in these bounds are some explicit $L^1$ integrals \textit{independent} of the weights $\vp, \psi_i$. In Section {\secholconst} of the supplementary material II in Part II \cite{ChenHou2023bSupp}, we estimate these explicit integrals using integral formulas, numerical quadrature, and obtain the constants for these bounds rigorously.

\vs{0.05in}
\paragraph{\bf Step 4: Estimate of more regular nonlocal terms.}\label{para:step4}

After estimating the most singular part in \hyr[para:step3]{Step~3}, the remaining nonlocal terms involve velocities with \emph{desingularized} kernels. 
In contrast to the estimates in \hyr[para:step3]{Step~3}, these estimates \emph{depend on the weight}.  A typical term has the form
\footnote{
The weights $ \vp, \psi_1$ below and in \eqref{eq:idea_norm_form} are used 
for weighted $L^{\infty}$ and $C^{1/2}$ estimates for $\om$. See \eqref{wg:linf_decay} 
and \eqref{wg:hol}.
}
\[
  I(\om) = \psi_1 (u_x(\om) - \hat u_x(\om)) - u_{x,a,b}(\om \psi_1) ,   
   \quad  u_{x,a,b}(f) \teq u_x(f, a, b),
\]
where $\hat f$ is the finite-rank approximation for the relevant velocity quantities $f=\uu$ or $\na\uu$  constructed in Section \ref{sec:appr_vel}. We subtract the most singular part $u_{x,a,b}(\om\psi_1)$, which is estimated by the method in \hyr[para:step3]{Step~3}.
The remainder term $I(\om)$ is more regular: 
it gains almost one derivative from $\om$ and is log-Lipschitz. We rewrite $I(\om)$ as an integral of $\om$ and bound $|I(\om)(x)|$ or its H\"older quotient $|I(\om)(x) - I(\om)(z) | \cdot |x-z|^{-1/2}$
in the form
\beq\label{eq:idea_norm_form}
        \CCs_1 ||\om\vp||_{L^\infty}
        + \CCs_2 [\om \psi_1]_{C_{x_1}^{1/2}} +        \CCs_3 [\om\psi_1]_{C_{x_2}^{1/2}} , 
\eeq
where $[ \cdot ]_{C_{x_i}^{1/2}}$ is the directional H\"older seminorm introduced in \eqref{hol:semi}. The coefficient functions $\CCs_i$, depending on $x$ or $(x,z)$, bound explicit weighted $L^1$ integrals involving the weights $\psi_1,\vp$.

We estimate these integrals using the scaling symmetries of the kernels,
symmetrization of the integrands, and rigorous numerical quadrature. We explain the main ideas of the integral estimates in Section \ref{sec:u_comp_idea} and refer to Part II \cite[Section {\secvelcomp}]{ChenHou2023b} for full details. The $C^{1/2}$ seminorm in \eqref{eq:idea_norm_form} 
 is mainly used to control the most singular part  in the $C^{1/2}$ estimate in \hyr[para:step3]{Step~3}.

\vs{0.05in}
\paragraph{\bf Step 5: Finite-rank perturbation}\label{para:step5}

To estimate finite-rank operators, such as those introduced in \hyr[para:step4]{Step~4}, 
we develop an $L^{\inf}$-based finite-rank perturbation method in Section \ref{sec:finite_rank}.

The method decomposes the perturbation so that the main component $W_1$ can be controlled by energy estimates for the modified linearized operator $\cL-\cK$, where $\cK$ is a finite-rank operator. By designing $\cK$ to approximate the nonlocal terms, we obtain sharper linear stability estimates for  $\cL - \cK$. The contribution of $\cK$ is captured by 
the variable $\widehat W_2$. Instead of estimating $\widehat W_2$ by energy estimates, we use Duhamel's formula \eqref{eq:bous_W2_appr} and construct and estimate finitely many approximate space-time solutions $\hat F_i$ associated with the range of $\cK$.

\vs{0.05in}
\paragraph{\bf Step 6: Energy estimates}

In Section \ref{sec:EE}, we perform weighted $L^{\infty}$ and $C^{1/2}$ stability estimates and construct the full energy $E_4$ \eqref{energy4}. 
To obtain  top-order $C^{1/2}$ coercivity estimates, we combine the weighted coercivity estimates of local terms in \hyr[para:step2]{Step~2} and sharp $C^{1/2}$ estimates in \hyr[para:step3]{Step~3}. While the outgoing property \eqref{eq:sec2_outgo}  generates sufficiently strong damping terms only  for the local parts, a key observation is that the nonlocal operator is \emph{essentially} local up to a \emph{compact operator}, since the Biot--Savart kernel is smooth away from the diagonal \(y=x\). 
We exploit this by showing that the commutator between the weight and nonlocal
operator is more regular. Combining the outgoing property, the \emph{essential local structure}, and the profile asymptotics, in a suitable weighted $C^{1/2}$ space, we  obtain a decomposition of the linearized operator 
\[
 \cL = \cA + \cC, 
\]
 where $\cA$ is coercive and $\cC$ is compact. The coercivity of $\cA$ captures the stability of the leading-order terms at the level of regularity
\footnote{
This is similar to the $H^k$ coercivity estimate for the linearized operator up to 
lower-order terms \cite{merle2022implosion2,chen2024Euler,chen2024vorticity,buckmaster2022smooth}.
}
and is the cornerstone of the nonlinear stability analysis.

The associated stability conditions \eqref{eq:PDE_nondiag} with $|x- z| \to 0$
provide the explicit 
conditions ensuring the coercivity of $\cA$ in this top-order $C^{1/2}$ estimate. \footnote{
If $|x-z| \gtr 1$, we get $ I= |f( x ) - f( z )|  /  |x - z |^{1/2} \les   
\| f \|_{L^{\infty}}$, and hence $I$ is lower order in regularity.
}
This part of the argument does not require the estimates for the more regular parts of $\uu,\na \uu$ 
from \hyr[para:step4]{Step~4}, nor the estimates for the approximate space-time solution $\widehat F_i$ from \hyr[para:step5]{Step~5}.
Indeed, these more regular terms are log-Lipschitz, contribute to the  compact operator $\cC$, 
and their $C^{1/2}$ estimates vanish as $|x-z|\to 0$. 

To estimate the more regular terms, we use the method in \hyr[para:step5]{Step~5} and approximate these terms using some finite-rank operator. We quantify the approximation error of the more regular nonlocal terms using methods in \hyr[para:step4]{Step~4}. We further bound the (semi)norms in the nonlocal estimates, e.g. those in \eqref{eq:idea_norm_form}, using the energy $E_4$. With these estimates, we derive the coefficients in the stability conditions \eqref{eq:PDE_nondiag} 
associated with the energy estimates. These coefficients depend \textit{explicitly} on the profile and weights and are computable. The full inequalities, which include the energy estimates of several norms 
outlined in Figure \ref{Fig:order}, are given in Appendix \ref{sec:ineq}.

\vs{0.05in}
\paragraph{\bf Step 7: Verify inequalities for nonlinear stability}

The preceding steps reduce the nonlinear stability analysis to a finite list of explicit inequalities
collected in Appendix \ref{sec:ineq}. The inputs to these inequalities are the 
coefficient functions for weighted $L^{\infty}$ and $C^{1/2}$ estimates 
of nonlocal terms from \hyr[para:step3]{Step~3} and \hyr[para:step4]{Step~4}, the estimates for the approximate space-time solutions $\widehat F_i$ 
from \hyr[para:step5]{Step~5}, the bounds for the local residual error and the residual operator from Part II \cite{ChenHou2023b} 
(see ideas in Appendix \ref{app:idea_num_est}) and the bootstrap threshold $E_*$.

All these quantities depend only on the approximate steady state and the chosen weights. 
In Part II \cite{ChenHou2023b},  we obtain rigorous computer-assisted bounds 
for these quantities and use them to verify the inequalities in Appendix \ref{sec:ineq}. The corresponding codes are available in~\cite{ChenHou2023code}.

After verifying the stability conditions \eqref{eq:PDE_nondiag}, we obtain nonlinear stability estimates $E_4(t)< E_*$ for all $t>0$ using  Lemma \ref{lem:PDE_nonstab}, which implies the bounds 
and blowup results in Theorem \ref{thm:main}.

\subsubsection{Small and large parameters in stability analysis}\label{sec:analy_para}

We highlight the main methods that produce the small and large parameters needed to close the stability estimates.

\vs{0.05in}

\paragraph{\bf Small residual error}

The residual error plays a crucial role as a small parameter for the stability analysis.
Our key stability estimates can be captured by the following inequality
\bseq\label{eq:idea_EE}
\beq\label{eq:idea_EE_0}
  \tfrac{d}{dt} E(t) \leq - \lam E(t) + C E(t)^2 + \e,
\eeq
for the energy $E(t)$, where $\lam, C >0$ are stability constants, and $\e$ estimates the residual error of the profile. To obtain nonlinear stability estimates, sufficient 
\emph{quantitative} conditions are 
\beq\label{eq:idea_EE_close}
   4  C \e  < \lam^2,
\eeq
\eseq
and $E(0) < \f{2 \e}{\lam}$. By constructing a profile with residual error $\e$ small enough so that \eqref{eq:idea_EE_close} holds, the nonlinear terms in the stability analysis can be treated \emph{perturbatively}.

\vspace{0.05in}

\paragraph{\bf Outgoing property} 

We use the outgoing property \eqref{eq:sec2_outgo} to generate \textit{arbitrarily strong} damping terms for the local part of the linearized operator. See Step 2 and Section \ref{sec:outgo}.

\vspace{0.05in}
\paragraph{\bf Finite-rank perturbation}

We construct finite-rank approximations $\hat \uu(\om)$ for $\uu$ 
and $\wh{(\na \uu)}(\om)$ for the nonsingular part of the integral defining
$\na\uu(\om)$, which are more regular than $\na\uu(\om)$. In the weighted $L^\infty$ estimates of the velocity remainders
$\uu_A=\uu-\hat\uu, \, (\na\uu)_A=\na\uu-\wh{(\na\uu)},$
and in the $C^{1/2}$ estimates of $\uu_A$, we refine the approximation $\hat \uu, 
\wh{ (\na \uu) }$ following Section \ref{sec:appr_vel} and
increase the rank $N$. This yields \emph{sufficiently small} coefficient bounds in these estimates,  such as $\CCs_1$ in
\eqref{eq:idea_norm_form}. Hence, after subtracting the approximations of these
more regular nonlocal terms, the stability of $\cL-\cK_N$ for sufficiently large
$N$ is reduced to estimates similar to the top-order $C^{1/2}$ estimates and the
 model problem with localized velocity in Section \ref{sec:model_trunc}. In particular,
 the velocity remainder $\uu_A, (\na \uu)_A$ are small enough to be treated as perturbations of
the local damping term. In practice, the anisotropic structure (see Sections
\ref{sec:aniso}, \ref{sec:aniso_est}) allows us to use only a small rank $N$
and to treat the nonlocal terms in $\cL-\cK_N$ perturbatively.

\begin{remark}[\bf Nonlocal terms are not small]
All the nonlocal terms in the linearized equations are not small. 
Without the finite-rank perturbation method, we cannot treat nonlocal terms perturbatively. Without sharp $C^{1/2}$ estimates in \hyr[para:step3]{Step 3}, the stability constants become much weaker, requiring a much smaller residual error of the profile to close the stability estimates.
\end{remark}

\subsubsection{Exact vs. approximate profiles.}

Following this work, a natural question is:

 \begin{widequote}
If one has an \emph{exact} profile rather than only an \emph{approximate}
profile, what properties of the profile are needed to construct blowup from
\emph{smooth}, \emph{finite-energy} initial data?
 \end{widequote}

For an \emph{approximate profile}, where $\e \neq 0$, the stability constants $C,\lambda$ in
\eqref{eq:idea_EE_0} must be computable: even with linear stability $\lambda>0$, a computable 
and effective constant $C$ is needed to verify
\eqref{eq:idea_EE_close}.  For an \emph{exact} smooth profile, 
however, the residual vanishes $\varepsilon=0$.
Nonlinear stability follows from \eqref{eq:idea_EE_0} for sufficiently small
initial perturbations, without requiring quantitative bounds on all constants.
One may then use \emph{non-quantitative} functional analytic arguments to
control the compact part of the nonlocal operators.
Thus the outgoing condition \eqref{eq:sec2_outgo}, together with suitable
regularity and decay of the exact profile, becomes the essential structural
input for proving finite-codimension stability.

The methods developed in this work 
indicate the main ingredients:
singularly weighted estimates capturing the outgoing condition, a finite-rank or compact decomposition of the nonlocal terms, and a splitting of the perturbation
into a main part estimated by energy estimates and a smoother part estimated
by Duhamel's formula. Building on these ideas, the first author recently proved
finite-codimension stability of exact $C^\alpha$ self-similar vorticity profiles
and sharp asymptotically self-similar blowup for 3D axisymmetric Euler without
swirl from $C_c^\alpha$ vorticity and $C^{1,\alpha}\cap L^2$ velocity
\cite{chen2026eulerI,chen2026eulerII} for any $\al \in (0, \f13 )$.

This suggests a potential route toward an analytic
proof of blowup with smooth data: 
\begin{itemize}[leftmargin=1em]
\item construct an exact smooth profile with the outgoing condition similar to \eqref{eq:sec2_outgo};
\item prove finite-codimension stability of the profile using the analytic
mechanisms above.
\end{itemize}
We remark that constructing such an exact smooth blowup profile remains a fundamental open
question. One possible approach is through a fixed-point argument; see
\cite{hou2025nonuniqueness,chen2026eulerI,chen2026eulerII}.

\subsubsection{Numerical analysis and computer-assisted proofs}
 In our energy estimates, we need to derive rigorous and tight piecewise bounds of various quantities involving the approximate steady state, singular weights, and several explicit functions. 
 The full details are developed in Part II \cite{ChenHou2023b}. 
  In Appendix \ref{app:idea_num_est}, 
 we provide a self-contained discussion of some ideas and methods 
 for rigorous bounds using  numerical analysis.

Computer-assisted proofs have played an important role in  several PDE problems, especially in computing explicit tight bounds of complicated integrals \cite{castro2014remarks,cordoba2017note,gomez2014turning} or bounding the norms of linear operators \cite{castelli2018rigorous,enciso2018convexity}. We refer to \cite{gomez2019computer} for an excellent survey on computer-assisted proofs in PDEs.  Note that our computer-assisted approach to establish stability differs from existing approaches, e.g. \cite{castro2020global}.  We \textit{do not} use direct computation to quantify the spectral gap of the linearized operator, since this operator \emph{cannot be} approximated by a finite-rank operator.

In the rest of this section, we discuss the main ideas for our stability analysis.

\subsection{Properties of the approximate steady state}\label{sec:property}

\begin{figure}[ht]
   \centering
      \includegraphics[width=\textwidth]{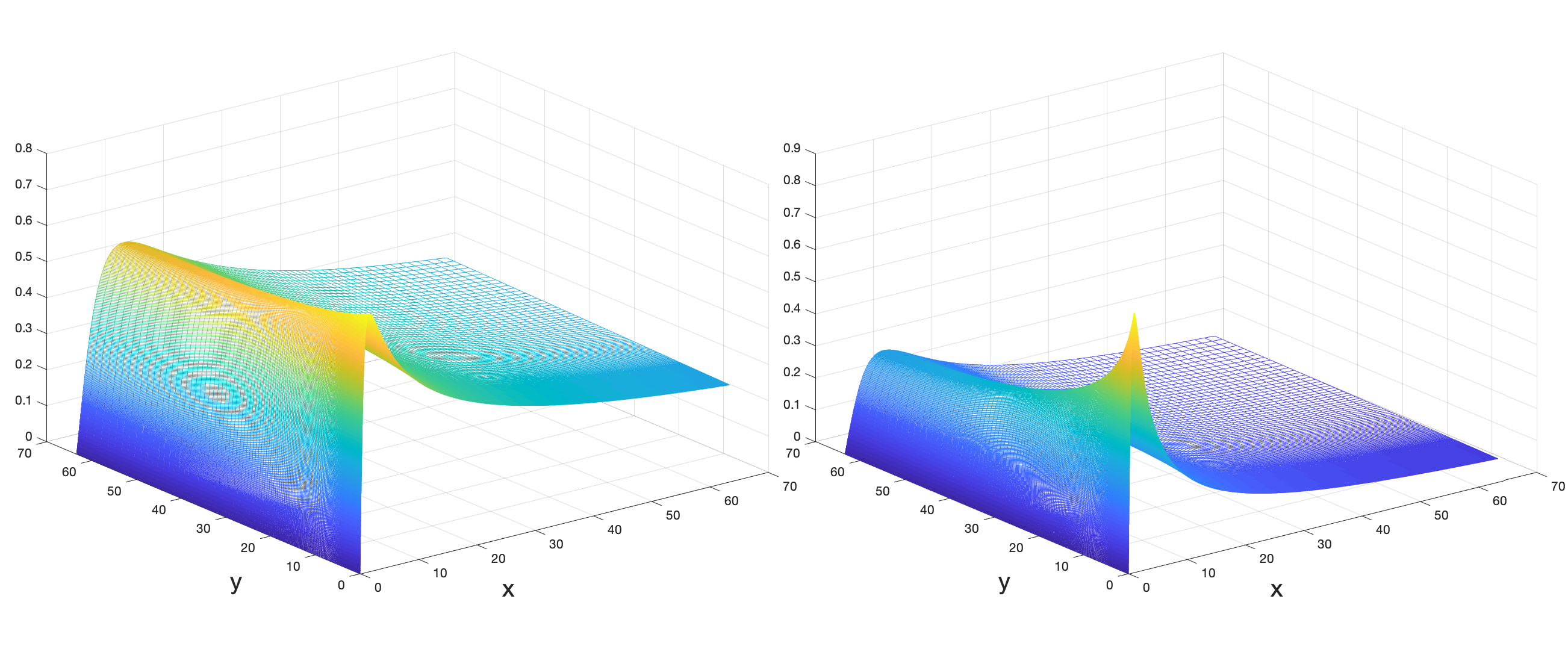}
      \vspace{-0.45in}
      \caption{Approximate blowup profile in the near field.
      Left: $\bar \om$; right: $\bar \th_x$. The plots illustrate the shape
      and anisotropy of the profile.}
      \label{fig:prof}
\end{figure}

\begin{figure}[ht]
   \centering
      \includegraphics[width=\textwidth]{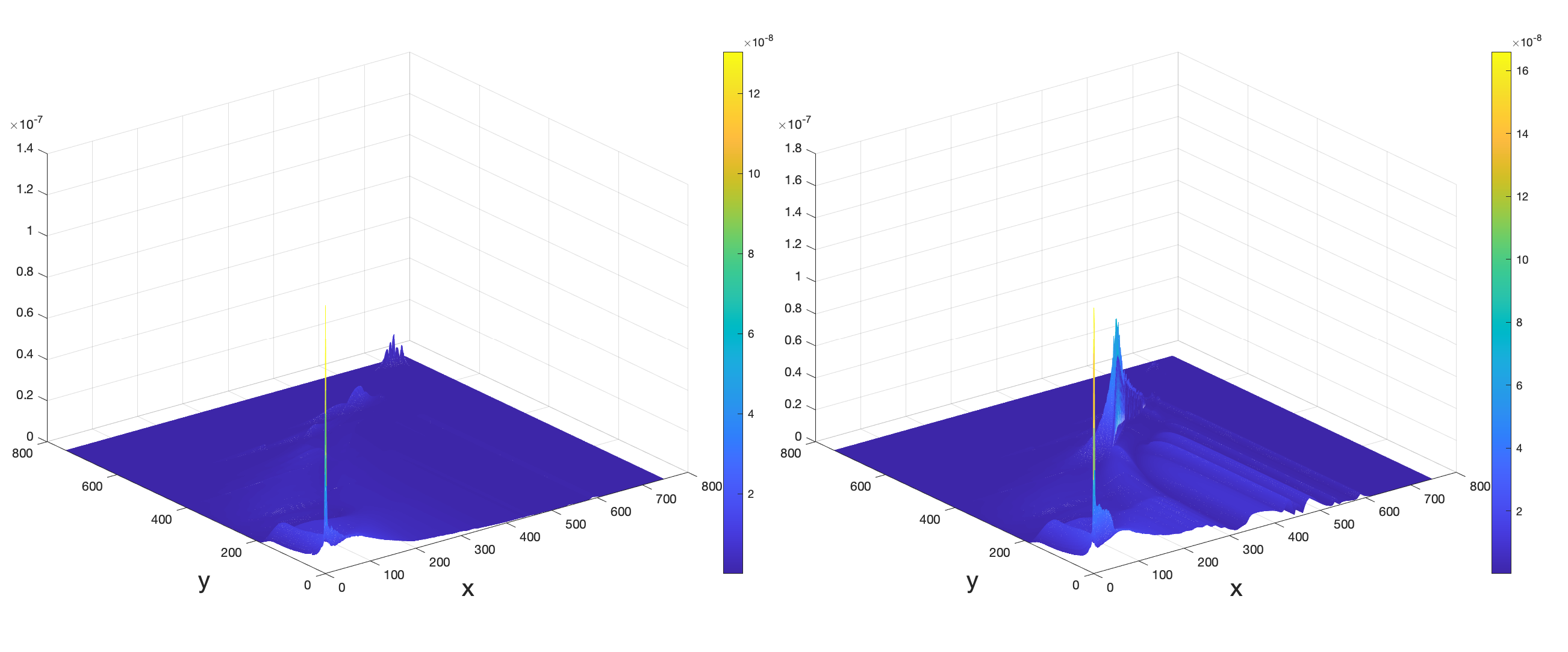}
      \captionsetup{width=0.93\textwidth}
      \vspace{-0.45in}
      \caption{Rigorous piecewise bounds for the weighted residual errors of
      the approximate steady state. Left: bound of $\bar \cF_{loc, 1}\vp_1$ in the
      $\om$-equation. Right: bound of $\bar \cF_{loc, 2}\vp_2$ in the
      $\th_x$-equation. These weighted bounds and $\vp_i$ are used in the energy estimates 
      in Section \ref{sec:EE}.
      }
      \label{fig:prof_err}
\end{figure}

We record the structural properties of the approximate steady state used in the
stability analysis. The construction of the approximate steady state $(\bar \om, \bar \th, \bar c_{\om}, \bar c_l)$ of \eqref{eq:bousdy1},\eqref{eq:normal} and its associated approximate stream function $\bar \phi^N$ is discussed in detail in Section \ref{sec:ASS}. 
    We then derive the residual $\bar \cF_i$ \eqref{eq:bous_err} 
    and its \emph{essential local} part $\bar \cF_{loc,i}$ \eqref{eq:bous_errM}.
\footnote{
The superscript $N$ in $\bar \phi^N$ stands for \emph{numerics}. Due to numerical errors, 
$  -\D \bar \phi^N  \neq \bar \om$. 
The local residual error $\bar \cF_{loc,i}$ \eqref{eq:bous_errM} differs from $\bar \cF_i$ \eqref{eq:bous_err}  by contributions from the error $\bar \e= \bar \om + \D \bar \phi^N$ 
and a near-origin correction. In Section \ref{sec:comb_vel_err}, we rigorously combine the estimate of the \emph{nonlocal error} $\bar \cF_i - \bar \cF_{loc, i}$ and the energy estimates. 
} 
In Figure \ref{fig:prof}, we plot $\bar \om, \bar \th_x$
\footnote{
We plot the variable $\bar \th_x$ rather than $\bar \th$ since $\bar \th$ grows in the far-field. 
}
on the grid points. In Figure \ref{fig:prof_err}, we plot rigorous piecewise  bounds of the weighted $L^{\inf}(\vp_i)$ norms of $\bar \cF_{loc,i}, i=1,2$. Since $\vp_1, \vp_2$ are very singular near $x=0$, with leading order $ |x|^{-2.4} |x_1|^{-\f{1}{2}}$, $c |x|^{-\f{5}{2}} |x_1|^{-\f{1}{2}}$, the weighted $L^{\inf}(\vp_i)$ norms 
  of $\bar \cF_{loc,i}$ are relatively large near the origin. By contrast, the unweighted errors $ \bar \cF_{loc, i}$  are very small near the origin and less than $2 \cdot 10^{-12}$, since we use a uniform fine grid near the origin.

\vspace{0.05in}
\paragraph{\bf Key outgoing property}

The advection in \eqref{eq:bousdy1} satisfies the \emph{crucial outgoing property}
\footnote{
\label{foot:outgo_bd}
In the energy estimates, we do not use inequality \eqref{eq:outgoing} 
or \eqref{eq:bous_aniso} with uniform constants $c_i$ for all $x$. Instead, we use the piecewise bounds for  $\bar c_l x + \bar u, \bar c_l y + \bar v$ to derive damping terms, 
and piecewise bounds for $\na \bar \th, \na \bar \om$, which are much sharper. These piecewise bounds are incorporated into the nonlinear stability inequalities in Appendix \ref{sec:ineq},
which are further verified with computer assistance. 
}
\begin{equation}\label{eq:outgoing}
\bar c_l x + \bar u(x, y) \geq c_1 x,  \quad c_1 \approx 0.47, \quad 
\bar c_l y + \bar v(x, y) \geq c_2 y, \quad c_2 \approx  3, 
\tag{Outgo}
\end{equation}
for all $(x, y) \in \R^2_{++}$. For $(x, y)\in \R^2_{++}$ near the origin, 
the advection is \textit{anisotropic}, and we have
\[
\bar c_l x + \bar u(x, y) \approx 0.47 x, \quad \bar c_l y + \bar v(x, y) \approx 5.54 y . 
\]

Property \eqref{eq:outgoing} is the key stability mechanism for the linearized operator; see Section~\ref{sec:outgo}.

\vspace{0.05in}
\paragraph{\bf Scaling exponents}

The scaling constants and the velocity gradient near the origin satisfy
\footnote{
The displayed decimal values are rounded. In the computer-assisted verification, 
$\bar c_l , \bar c_{\om}, \bar u_x(0)$ in \eqref{eq:bous_expo} 
and $\bar \al_1$ in \eqref{solu:decay} are computed from the represented approximate profile with rigorous interval-arithmetic enclosures.
}
\footnote{
By scaling symmetry, if
$(\om_*,\th_*,c_{l,*},c_{\om,*})$
\emph{exactly} solves the profile equations, i.e. \eqref{eq:bousdy1} without 
$\pa_t\om,\pa_t\th$ terms, then so does $ (a\om_{*,b},a^2 \th_{*,b},a \, c_{l,*},a \, c_{\om,*}) $
for any $a,b>0$, where $\om_{*,b}(x)\teq \om_*(bx), \,  \th_{*,b}(x)\teq b^{-1}\th_*(bx)$. 
Similarly, multiplying the approximate profile by an amplitude factor $a$,
one may renormalize $\bar c_l$, $\bar c_{\om}$, or $\bar u_x(0)$ in \eqref{eq:bous_expo} while leaving the ratio $\bar c_l/\bar c_{\om}$ \emph{invariant}. Thus, 
the ratio $\bar c_l/\bar c_{\om}$ characterizes the blowup rate. 
}
\beq\label{eq:bous_expo}
\bar c_l \approx 3.00649798, \quad \bar c_{\om} \approx -1.02942517,  \quad  \bar u_x(0) \approx -2.532674 , \quad  \bar v_x(0) = 0.
\eeq
Note that the ratio $\bar c_l / \bar c_{\om} \approx -2.9205600$ is very close to the one reported by Hou-Luo \cite{luo2014potentially,luo2013potentially-2}.

\vspace{0.05in}
\paragraph{\bf Regularity, representation, and decay }
The variables $\bar \om, \bar \phi^N$ are odd in $x$ and $\bar \th$ is even in $x$. The approximate steady state $(\bar \om, \bar \th, \bar \phi^N)$ is represented by piecewise 6th order polynomials $ \bar \om_2, \bar \th_2, \bar \phi_2$ supported in $[0, D_1]^2$ with $D_1 \approx 10^{15}$, semi-analytic parts $\bar \om_1, \bar \th_1, \bar \phi_1$ that capture the far-field behavior of the solutions, an analytic part $\bar \phi_3$ that captures $\bar \phi^N$ near $x=0$ to reduce the round-off error, and a small rank-one correction $\bar \phi_{\mf{cor}}$ such that $\bar \om - (-\D) \bar \phi^N = O(|x|^2)$ near $0$:
\[
\bar \om = \bar \om_1 + \bar \om_2, \quad \bar \th = \bar \th_1 + \bar \th_2, \quad \bar \phi^N = \bar \phi_1 + \bar \phi_2 + \bar \phi_3 + \bar \phi_{\mf{cor}}.
\]
See \eqref{eq:ASS_decomp1} and Section \ref{sec:ASS}. Since each basis function in the representation is $C^{4,1}$,
\footnote{
The least regular basis function is the piecewise 6th-order B-spline $B_{1,i}(x) B_{2,j}(y)$,  
which is a linear combination of terms of the form $(x - a)^5_+ (y - b)_+^5$ for some $a, b$ and has $ C^{4, 1}$ regularity.
}
we have $\bar \om, \bar \th, \bar \phi^N \in C^{4,1}$. From classical embedding inequalities, for the \textit{exact velocity} $\bar \uu \teq \uu(\bar \om) = \na^{\perp}(-\D)^{-1} \bar \om $,
 we have $\na \bar \uu \in C^{4, \al}$ 
 \footnote{ 
 Due to numerical errors, we do not have $\bar \phi^N = (-\D)^{-1} \bar \om$. Thus, the regularities of the \emph{exact velocity} $\uu(\bar \om)$ and the approximate velocity $\na^{\perp} \bar \phi^N$ are different. We estimate this error rigorously in Section \ref{sec:comb_vel_err}.
 }
 for any $\al \in [0, 1)$ and $ \f{ \bar \uu }{|\xx|}  \in L^{\infty}$. For $i\leq 4$, 
 from the representation in Section \ref{sec:ASS}, the solution satisfies the regularity and decay rate 
\footnote{
For sufficiently large $|\xx|$, we have $\bar \om = \bar \om_1, \bar \th = \bar \th_1$ with formulas in \eqref{eq:ASS_semi1},\eqref{eq:ASS_decomp1}, since $\bar \om_2, \bar \th_2$ have compact support.
}
\beq\label{solu:decay} 
|\na^i \bar \om(\xx)| \les \la \xx \ra^{\bar \al_1-i}, \quad  |\na^i \bar \th(\xx)| \les \la \xx \ra^{1 + 2 \bar \al_1 - i}, 
\quad |\na^i \bar \phi^N| \les \la \xx \ra^{2 + \bar \al_1- i}, 
\quad \bar \al_1 \in [-0.343, -0.342] ,
\tag{Decay}
\eeq
where $ \bar \al_1$ is close to $ \bar c_{\om} / \bar c_l$, and the implicit constants are independent of $\xx$. The decay properties are used in the far-field weighted estimates. See Section \ref{sec:ASS} for more details.

\vspace{0.05in}
\paragraph{\bf Anisotropy of the profile}

The profile is anisotropic in the sense that the $y$-derivatives are much
smaller than the $x$-derivatives, especially in the near field
$(x,y)\in [0,1]^2$:
 \beq\label{eq:bous_aniso}
|\bar \th_y| \leq c_3|\bar \th_x| , \quad c_3 \approx 0.16, \quad 
|\bar \om_y| \leq c_4 |\bar \om_x|, \quad c_4 \approx 0.23.
\eeq

Although $c_3,c_4$ are not small enough to make all such terms perturbative, this anisotropy helps identify the leading interactions and provides smallness for several coupling terms.
See Footnote \ref{foot:outgo_bd} for more discussions of the constants $c_3, c_4$ and \eqref{eq:bous_aniso}.

\subsection{Linearized equations}\label{sec:linop}

We decompose a solution $(\om_{\tot}, \th_{\tot}, c_{\om, \tot}, c_{l, \tot})$ 
to \eqref{eq:bousdy1} as the approximate profile $(\bar \om, \bar \th, \bar c_{\om}, \bar c_l)$
and its perturbation $( \om, \th, c_{\om}, c_l )$:
\[
  (\om_{\tot} , \, \th_{\tot}, \, c_{\om, \tot} , \, c_{l, \tot} ) \teq ( \om + \bar \om, \, \th + \bar \th, \, c_{\om} + \bar c_{\om}, \, c_l + \bar c_l)  .
\]

Since $\om, \na \th$ have similar regularity, we study the system of $(\om, \th_x, \th_y)$ and denote 
\beq\label{eq:eta_xi}
\eta = \th_x, \quad \xi = \th_y.
\eeq

Linearizing \eqref{eq:bousdy1} around $(\bar \om, \bar \th, \bar \uu, \bar c_l, \bar c_{\om})$, we obtain the equations for the perturbation
\beq\label{eq:bous_lin11}
\bal
\om_t &= - (\bar c_l x +\bar \uu) \cdot \na \om + \th_x + \bar c_{\om} \om 
 -  \uu \cdot \na \bar \om + c_{\om} \bar  \om + \bar F_{\om} + N(\om)  , \\
\th_t & = - (\bar c_l x + \bar \uu) \cdot \na \th + \bar c_{\th} \th + c_{\th} \bar \th - \uu \cdot \na \bar \th 
+ \bar F_{\th} + N(\th) ,\quad 
\uu 
=\na^{\perp}(-\D)^{-1} \om, \quad  \bar c_{\th} = \bar c_l + 2 \bar c_{\om},
\eal
\eeq
where $\bar \uu= (\bar u, \bar v) = \nabla^{\perp}(-\D)^{-1} \bar \om$ denotes the \emph{exact velocity}. Here, we define the residual errors \eqref{eq:bous_err} $\bar F_{\om}, \bar F_{\th}$ associated with the approximate steady state $\bar \om, \bar \th, \bar c_l, \bar c_{\om}$, and related functions $\bar \cF_i$
\beq\label{eq:bous_err}
\bga
 \bar F_{\om}   = - (\bar c_l x + \bar \uu ) \cdot \na \bar \om + \bar \th_x + \bar c_{\om} \bar\om \,  , \quad  \bar F_{\th} = - (\bar c_l x + \bar \uu ) \cdot \na \bar \th + \bar c_{\th} \bar \th , 
 \\
\olin \cF_1  \teq \bar F_{\om}, \qquad \olin \cF_2 \teq \pa_x \bar F_{\th}, \qquad \olin \cF_3 \teq \pa_y \bar F_{\th} \, , 
 \ega
 \eeq
We have the following conditions by linearizing \eqref{eq:normal}
\footnote{
In \eqref{eq:normal}, we can write the full parameters as 
$c_l = \bar c_l + \td c_l, c_{\om} = \bar c_{\om} + \td c_{\om}$ with 
perturbation $\td c_l ,\td c_{\om}$. In \eqref{eq:normal_pertb}, we abuse notation and $c_l, c_{\om}$ denote the perturbations $\td c_l, \td c_{\om}$.
}
\beq\label{eq:normal_pertb}
c_{\om} =  u_x(0), \quad  c_l \equiv 0 , \quad c_{\th} = c_l + 2 c_{\om} 
= 2 c_{\om} = 2 u_x(0).
\eeq
Moreover, we define   the nonlinear terms $N(\om), N(\th)$ and related functions $\cN_i$: 
\beq\label{eq:bous_non}
\bal
    \cN_1 & =  N(\om)  \teq - \uu \cdot \na \om + u_x(0) \om ,  \quad 
  \quad N(\th) = - \uu \cdot \na \th + 2 u_x(0) \th , \\
 \cN_2 & = \pa_x N(\th) = - \uu \cdot \na \eta - u_x \eta - v_x \xi + 2 u_x(0) \eta,  \\
 \cN_3  & = \pa_y N(\th) =- \uu \cdot \na \xi - u_y \eta - v_y \xi + 2 u_x(0) \xi. 
 \eal
\eeq

Taking derivatives of the $\th$-equation in \eqref{eq:bous_lin11} and using 
\eqref{eq:bous_non} and \eqref{eq:bous_err} for the nonlinear $\cN_i$ and error terms $\olin \cF_i$, we obtain 
\beq\label{eq:bous_lin12}
\bal
 \pa_t \eta & =  - (\bar c_l x +\bar \uu) \cdot \na \eta +
(2 \bar c_{\om} - \bar u_x) \eta - \bar v_x \xi - \uu_x \cdot \na \bar \th 
- \uu \cdot \na \bar \th_x  + 2 c_{\om} \bar \th_x 
+  \cN_2 + \olin \cF_2  , \\
  \pa_t \xi & = -  (\bar c_l x +\bar \uu) \cdot \na \xi +
(2 \bar c_{\om} + \bar u_x) \xi  - \bar u_y \eta - \uu_y \cdot \na \bar \th 
- \uu \cdot \na \bar \th_y  + 2 c_{\om} \bar \th_y 
+ \cN_3 + \olin \cF_3 , \\
\eal
\eeq
where we have used $ c_{\th} = c_l + 2 c_{\om}$ by \eqref{eq:normal_pertb}.  Due to the normalization conditions \eqref{eq:normal1} and the odd symmetries of $\th_x, \om$ we have the following vanishing conditions near the origin
\beq\label{eq:normal_vanish}
 \om = O(|x|^2), \quad 
\th_x = O( |x|^2), \quad \th_y = O( |x|^2).
\eeq

Note that numerical evidence for the linear stability of the above system has been reported by Liu \cite{liu2017spatial} (see \cite[Section 3.4]{liu2017spatial}), who showed that the eigenvalues of the discretized linearized operator have negative real parts bounded away from $0$.

\subsection{Main terms of the system}\label{sec:idea_main_term}
We identify the main terms in the system \eqref{eq:bous_lin11},\eqref{eq:bous_lin12}.

\subsubsection{Anisotropy in the $x, y$ directions}\label{sec:aniso}

 In the near field, the $y$-derivatives of the solution, e.g., $\bar \om_y, \bar \th_y, \bar \th_{xy}$, are relatively small. In \eqref{eq:bous_lin12}, $\xi = \th_y$ enjoys much better stability estimates than those of $\eta =  \th_x$ due to the flow structure: compression along the $x$-direction and outward flow along the $y$-direction. Indeed, since $\bar u_x(0) \approx -2.5$  near the origin and $\bar c_{\om} \approx -1$  \eqref{eq:bous_expo}, we have 
\[
( 2 \bar c_{\om} - \bar u_x) \eta \approx 0.5 \eta, \quad ( 2 \bar c_{\om} + \bar u_x ) \xi \approx  - 5.5 \xi .
\]
These sign conditions also hold in the near field. The first term contributes a growing term in the $\eta$-equation, while the second contributes a strong damping term in the $\xi$-equation.

\subsubsection{Weak coupling}\label{sec:weak_coup}
Note that $\bar v_x \approx 0$ near $0$ \eqref{eq:bous_expo} and $\bar v(x, 0) = 0$ due to the boundary condition, $\bar v_x$ is quite small in the near field and near the boundary. 
Therefore, $\xi$ is weakly coupled to the equation of $\eta$ in such a region, which is the most difficult region of the analysis. See the discussion in Section \ref{sec:aniso_est}. Using the above analysis and dropping the smaller terms and the $\xi$ equation, we identify the main terms in the linear part of the system \eqref{eq:bous_lin11} 
\beq\label{eq:bous_main}
\bal
\pa_t \om &= - (\bar c_l x +\bar \uu) \cdot \na \om + \eta + \bar c_{\om} \om 
- u \bar \om_x +  c_{\om} \bar  \om + \cR_{\om}, \\
 \pa_t \eta & =  - (\bar c_l x +\bar \uu) \cdot \na \eta +
(2 \bar c_{\om} - \bar u_x) \eta  - u_x \bar \th_x 
- u \bar \th_{xx}  + 2 c_{\om} \bar \th_x  + \cR_{\eta}, 
\eal
\eeq
where $\cR_{\om}, \cR_{\eta}$ denote the remaining terms in the equations. While $\cR_{\omega}$ and $\cR_{\eta}$ are not negligible, we first analyze \eqref{eq:bous_main} by setting $\cR_{\omega} = 0$ and $\cR_{\eta} = 0$, as this simplified case is more tractable than \eqref{eq:bous_lin11} and \eqref{eq:bous_lin12}. The resulting analysis guides our choice of energies and parameters.

\subsection{Coercivity estimates for the local parts }\label{sec:model_local}

Following \cite{chen2019finite,chen2019finite2,chen2021HL}, we perform weighted energy estimates and obtain the damping terms for the energy estimates from the local terms, especially $(\bar c_l x + \bar \uu) \cdot \na f$ in \eqref{eq:bous_main}. We design the energy space $X$ so that the local part of the linearized equations is \emph{stable} in $X$ and the nonlocal terms can be effectively estimated in $X$.

\subsubsection{Outgoing flow generates sufficiently strong damping}\label{sec:outgo}

We use property \eqref{eq:outgoing} to generate \textit{sufficiently strong} damping terms in $|x| \in [R_1, R_2]$ for any $0< R_1< R_2$. Denote $W = (\om, \eta, \xi)$. The linearized equations \eqref{eq:bous_lin11},\eqref{eq:bous_lin12} without $\cN_i, \bar \cF_i$ can be written 
schematically as 
\[
  \pa_t W + b(\xx) \cdot \na W = \bar M \cdot W + \cL_{nloc} \om,\quad 
  b(\xx) = \bar c_l \xx + \bar \uu ,
\]
where $\bar M(\xx) \in \R^{3 \times 3}$ denotes the coefficient matrix for the local terms, e.g. 
$\bar c_{\om} \om, (2 \bar c_{\om} - \bar u_x) \eta$, and $\cL_{nloc} \om$ denotes the nonlocal terms, including the $c_{\om}$ terms, which depend only on $\om$.

Multiplying both sides by a weight $\vp$ and commuting $\vp$ with the transport term, we obtain 
\[
  \pa_t (W \vp) + b(x) \cdot \na ( W \vp ) = (\bar M + d_{\vp}(x) ) \cdot (W \vp)  
  + \vp \cL_{nloc} \om , \quad 
    d_{\vp}(x) \teq (b \cdot \na \vp) \vp^{-1},
\] 
where $d_{\vp}(x) \cdot W \vp$ comes from the commutator. If $\bar M + \f{b \cdot \na \vp}{\vp}$ is negative (in an appropriate sense), we obtain damping terms. For this purpose, we extract damping terms from $d_{\vp}(x)$ and design 
\[
  \vp = \vp_s \vp_0^m , \quad \vp_s = |\xx|^{\al} + |\xx|^{\b}, 
  \quad \vp_0 \asymp 1 ,
\]
with $\al < 0,  m > 0$. Since $\al <0$, the weight $\vp$ is singular near $\xx =0$. We choose the exponents $\al,\b$ in $\vp_s$ to capture the stability mechanisms near $\xx=0$ and for large $|\xx|$. They are also chosen according to the asymptotics of $W$ in these regions, so that $W\vp$ is well-defined. We consider both $\b <0$ and $\b >0$. We have
\[
  d_{\vp} = ( b \cdot \na \vp_s ) \vp_s^{-1}
  +  m ( b \cdot \na \vp_0) \vp_0^{-1} \teq d_1( \vp_s ) + m d_2(\vp_0). 
\]

We can choose $\vp_0 \in C^1$ freely to obtain $ d_2 \leq  0$. 
For example, for any $0< R_1 < R_2$, we design  
\[
\bga
  \vp_0( \xx )=1, \quad |\xx| \leq \f{R_1}{2},
  \quad \vp_0( \xx ) = \f{1}{2},  
  \quad |\xx| \geq 2 R_2,  \\
\quad  \pa_i \vp_0( \xx )  < 0, \quad |\xx| \in [R_1, R_2].
\ega
\]
From property \eqref{eq:outgoing}, for some constant $C(R_1, R_2) > 0$ depending on $\vp_0$ and profile, we get
\[
  m d_2(\vp_0) = m ( b \cdot \na \vp_0) \vp_0^{-1} 
  \leq  - m C(R_1, R_2) \one_{[R_1, R_2]}( | \xx | ).
\]

Taking $m > 0$ sufficiently large, we generate \textit{sufficiently strong} damping terms in $|x| \in [R_1, R_2]$. Since $d_{\vp} W \vp$ appears only on the diagonal of $(\bar M + d_{\vp}) W \vp$, 
if we neglect the nonlocal terms $\vp \cL_{nloc} \om$, by applying Lemma \ref{lem:PDE_stab}
with energy $\max_i ( || W_i \vp||_{L^{\infty}} )$, we achieve the diagonally dominant stability condition \eqref{eq:PDE_diag} in $|x| \in [R_1, R_2]$, with $a_{ii} = -(\bar M_{ii} + d_{\vp}) >0, a_{ij} = \bar M_{ij}, i \neq j$.

In the actual estimates, we use different weights $\vp_i$ for different
components $W_i$. The same mechanism applies componentwise. Near the origin,
the damping comes from the singular factor $\vp_s$ and the local
matrix $\bar M(\xx)$; this is analyzed in Section~\ref{sec:model1}.

Stability for large $|\xx|$ is much easier to establish since the coefficients decay \eqref{solu:decay}. Assume 
\[
\vp_i \sim c_i |x|^{\b_i}, \quad  \na \vp_i \sim c_i \b_i |x|^{\b_i - 1 } 
\]
for large $|\xx|$.  Using the asymptotics \eqref{solu:decay}, we obtain 
\beq\label{eq:idea_damp_far}
\lim_{|\xx| \to \infty} \bar M + \rm{diag}( d_{\vp_1}, d_{\vp_2}, d_{\vp_3} ) 
= \rm{diag} (\bar c_{\om} + \bar c_l \b_1, 2 \bar c_{\om} + \bar c_l \b_2,  2 \bar c_{\om} + \bar c_l \b_3 )
+ \ee_{12} ,
\eeq
where $\ee_{12}$ is the matrix with $\ee_{12, 12} = 1, \ee_{12, ij} = 0$ for $(i, j) \neq (1, 2)$.
This off-diagonal term comes from  $\eta = \th_x $ in \eqref{eq:bous_lin11}. 
For stability in the far-field, we choose 
\[
        \b_1<|\bar \al|,
        \qquad
        \b_2, \, \b_3<2|\bar \al|,
        \qquad
        \bar \al= \tf{\bar c_{\om}}{ \bar c_l}\approx - \tf13.
\]
All exponents of $L^{\infty}$ weights $\vp_i, \vp_{g,i}$ (defined in \eqref{wg:linf_decay}, \eqref{wg:linf_grow} ) satisfy this constraint, e.g. 
\[
\vp_i: \   - \tf16   < |\bar \al| ,  \quad \tf17  < 2 |\bar \al| , 
\quad \, \vp_{g, i} : \tf{1}{16} < |\bar \al| ,  \quad  \al_{g, n}   < 2 |\bar \al| ,
\]
The corresponding stability constraint for the $C^{1/2}$ weights $\psi_i$  can be analyzed similarly.

It remains to explain why the nonlocal terms do not destroy the damping mechanism above. We decompose
 $\cL_{nloc} = \cL_{\na \uu} + \cL_{\uu} + \cL_{c_{\om}}$ in \eqref{eq:bous_lin11},\eqref{eq:bous_lin12}, where $\cL_{f} \om$ is the linear operator with the nonlocal term $f(\om)$, e.g. $\cL_{\na \uu}$ contains $\uu_x \cdot \na \bar \th$ in \eqref{eq:bous_lin12}. For weighted estimates, we have
\beq\label{eq:idea_Lnloc}
  \vp \cL_{nloc} \om
  = \cL_{\na \uu} (\om \vp) + ( [\vp, \cL_{\na \uu}] \om   + \vp \cL_{\uu} \om  + \vp \cL_{c_{\om}} \om ),
\eeq
where the commutator is given by
\[
   [\vp, \cL_{\na \uu}] \om  = \vp \cL_{\na \uu} \om - \cL_{\na \uu}(\om \vp) .
\]
As shown later in Section \ref{sec:nonloc_wg_idea}, the commutator gains one derivative.
The last three nonlocal terms in 
\eqref{eq:idea_Lnloc} are more regular than  $\om$, and can be approximated by finite-rank operators. Formally, $\cL_{\na \uu} (\om \vp)$ has size $O(1)$ relative to $\om \vp$, with size \textit{independent of}  $\vp$.  By designing the weight $\vp$ to generate a damping term with a large coefficient, we treat $\cL_{\na \uu} (\om \vp)$ as a perturbation to the local damping terms. This further yields the top order stability in Section \ref{sec:nonloc}.

\vs{0.05in}

\paragraph{\bf Nonlocal terms are small near $0$ }
By approximating $\uu$ with low-rank operators $\hat\uu$ from their Taylor expansion near $0$
(see Section \ref{sec:appr_vel}), the resulting nonlocal terms in \eqref{eq:bous_lin11}–\eqref{eq:bous_lin12}, 
e.g.$(\uu - \hat\uu) \cdot \na \bar\th_x$, vanish to higher order than the local terms
and thus can be treated perturbatively near $0$ \textit{without requiring large} local damping terms. See \eqref{eq:velA_reg},\eqref{eq:reg_W1} for the vanishing orders.

\vs{0.05in}

\paragraph{\bf Design of weights}

 Since we can compute  the coefficients $b, \bar M$, and only need to generate $O(1)$ damping to control the off-diagonal term $\bar M_{ij} W_j$ via \eqref{eq:PDE_diag} and Lemma \ref{lem:PDE_stab}, we explicitly  design  weights using combinations of powers $|\xx|^{a} x_1^{b}$  (for certain $a,b$) in Appendix \ref{app:para_wg}, and do not need to include large parameters in the weights. 
 These weights $\rho$ capture the outgoing stability mechanism. Moreover, 
  for $|\xx|$ not too large and $x_i \geq 0$, they satisfy 
 \[
   \pa_i \rho(\xx) <0, \quad  \rho(\xx)^{-1} \cdot (b_i(\xx) \pa_i \rho(\xx) )<0 
 \]

\subsubsection{Model for the local term} \label{sec:model1}

To analyze linear stability near $0$, from the above discussions, we can drop the nonlocal terms and focus on the local terms in the main system \eqref{eq:bous_main}. We also drop the remaining terms in \eqref{eq:bous_main} and approximate \eqref{eq:bous_main} near $0$ by the following model in $\R_2^{++}$
\beq\label{eq:model_loc}
\bal
\om_t + ( a_1 x \pa_x + a_2 y \pa_y) \om  & = -\om + \eta,  \\
\quad \eta_t + ( a_1 x \pa_x +a_2 y \pa_y ) \eta & = a_3 \eta, \quad a_1 = 0.5, \ a_2 = 5.5, \ a_3 = 0.5,
\eal
\eeq
with $\om, \eta$ being odd in $x$, where we have used \eqref{eq:bous_expo} to obtain approximations
\[
\bar c_l x + \bar u \approx (\bar c_l + \bar u_x(0)) x \approx 0.5 x, \quad 
\bar c_l y + \bar v \approx 5.5 y, \quad 2 \bar c_{\om} - \bar u_x(0) \approx 0.5.
\]

\paragraph{\bf{Weighted $L^{\infty}$ space}}

Since $a_3 >0$, $a_3 \eta$ in \eqref{eq:model_loc} contributes to a growing term. We consider  weighted $L^{\infty}$ estimates to take advantage of the transport structure. Suppose that $\vp = r^{-\g}, r =  (x^2 + y^2)^{1/2}$. Multiplying the $\eta$-equation by $\vp$ and performing a direct calculation yields 
\beq\label{eq:model_linf}
  \pa_t (  \eta \vp) +  ( a_1 x  \pa_x  + a_2 y \pa_y ) ( \eta \vp) 
  = ( a_3  \vp + a_1 x \pa_x \vp + a_2 y \pa_y \vp ) \eta \teq a(\g) \eta \vp, 
\eeq
Since $ x\pa_x \vp = -\g x^2 r^{-\g-2} , y\pa_y \vp = -\g y^2 r^{-\g-2}$ and $a_2 \geq a_1$, we get 
\beq\label{eq:damp_toy_linf}
\bal
& a(\g)= a_3 +  \f{ a_1 x\pa_x \vp + a_2 y \pa_y \vp }{\vp} 
= a_3 - \g \f{ a_1 x^2 + a_2 y^2}{x^2 + y^2} \leq a_3 - a_1 \g.
\eal
\eeq

Since $a_1 = 0.5, a_2 = 5.5, a_3 = 0.5$, to obtain a damping factor $a(\g) \leq 0$, we can choose $\g \geq 1$. Notice that for the system \eqref{eq:bous_lin11}, \eqref{eq:bous_lin12}, $\om, \eta$ vanish at least quadratically near $0$ by \eqref{eq:normal_vanish}. Therefore, we can choose $1 < \g \leq  2$ to derive the damping terms in the $\eta$ equation.

Performing $L^{\inf}$ estimate with weight $\vp = r^{-\g}$ and $\g > 1$ on both equations
in  \eqref{eq:model_loc}, we get 
\beq\label{eq:model_linf2}
\f{d}{dt} || \om \vp ||_{\infty} \leq (-1 - a_1 \g) || \om \vp ||_{\infty}
+ || \eta \vp ||_{\infty} , \quad  \f{d}{dt} || \eta \vp ||_{\infty} 
\leq (a_3 - a_1 \g )  || \eta \vp ||_{\infty}.
\eeq
It is easy to further obtain that $\max( || \om \vp ||_{\infty}, || \eta \vp ||_{\infty} ) $ decays exponentially fast.

From \eqref{eq:damp_toy_linf}, since $a_2$ is much larger than $a_1$, as the ratio $\lam = y / x$ increases, we get a much larger damping factor 
\[
a(\g, \lam) = a_3 - \g \f{a_1 + a_2\lam^2}{1 + \lam^2} , \quad a(\g, 0) = a_3 - a_1 \g
= \f{1}{2} - \f{1}{2} \g,  \ a(\g, \inf) = a_3 - \g a_2 = \f{1}{2} - \f{11}{2} \g \;.
\]

If one performs a weighted $L^2$ stability estimate of \eqref{eq:model_loc} similar to \cite{chen2019finite,chen2019finite2,chen2021HL,elgindi2019finite} with a singular weight $\vp = x^{-\al} y^{-\b}$, integration by parts yields that 
the $y$-advection contributes 
\[ 
I = - \int \eta a_2 y \pa_y \eta \vp  = \f{a_2(1-\b)}{2} \int \eta^2 \vp ,
\]
to the energy estimate of $\int \eta^2 \vp$. 
Since $ \eta(x, 0) $ need not vanish,
\footnote{
The Boussinesq system \eqref{eq:bous_lin11},\eqref{eq:bous_lin12} does not preserve 
 $\eta(x, 0) = 0$ due to the nonlocal terms.
} 
one requires $\b < 1$ for a well-defined energy; however, this results in a large anti-damping term $I$.
\footnote{
  If $\b$ is close to $1$, in estimating the nonlocal terms, e.g. $u_x$, since $u_x(x, 0) \neq 0$, we expect an estimate with a large constant: 
$ || u_x \vp^{1/2}||_2 \leq C (1-\b)^{-1/2} || \om \vp^{1/2}||_2$. If $1-\b$ is not small, $I$ is a 
large anti-damping term in the energy estimate. One requires a very singular weight to extract a damping term for $\eta$, e.g. $\al \geq 14$ if $\b =0$. 
}

A potential approach is to perform an $H^k$ estimate \cite{merle2022implosion2} with large $k$, but the coercivity of the linearized operator in $H^k$ estimate is unclear, even for top-order terms.\footnote{
 Taking a partial derivative $\pa_x^i \pa_y^j$ plays a role similar to that of a singular weight $ x^{-i}y^{-j}$. To derive a damping term for $ \pa_x^k \eta$ in the model \eqref{eq:model_loc}, one needs $k \geq 14$. Due to the mixed derivative terms $\pa_x^i \pa_y^j \eta$, it is not clear if $H^k$ estimates can be closed for \eqref{eq:model_loc} with $k$ not very large.
This approach can lead to many more terms in the system \eqref{eq:bous_lin11}, \eqref{eq:bous_lin12}, e.g. mixed derivative terms and $\pa_x^i  (u_x \bar \th_x) $, which can be difficult to control. Moreover, due to the boundary, $\pa_y$ does not commute with the nonlocal operator $ \na^{\perp}(-\D)^{-1}$ for the velocity. Thus, in the $H^k$ estimates, it is not clear if one can obtain stability for the leading order operator. 
}
Moreover, constructing a profile with a small residual in $H^k$ norm for large $k$, e.g. $k \geq 14$, is extremely challenging.

This is a new difficulty absent in \cite{chen2019finite,chen2021HL,chen2019finite2,elgindi2019finite}.
The Boussinesq system \eqref{eq:bous_lin11},\eqref{eq:bous_lin12} is a generic 2D system. We overcome it by using $L^{\inf}$-type estimates, the outgoing property \eqref{eq:outgoing},
and developing sharp estimates of the nonlocal terms in appropriate functional spaces.

\subsubsection{Vanishing order of the perturbation}\label{sec:vanish}

From \eqref{eq:normal_vanish}, the natural symmetry class only gives quadratic
vanishing of the perturbations $\om,\eta,\xi$ near $\xx=0$, i.e.,
$O(|\xx|^2)$. 
On the other hand, to obtain larger damping factors near $0$, we choose a more singular weight (larger $\g$) for the stability analysis, which is effectively of order $|\xx|^{-3}$ near the origin.
See Sections~\ref{sec:hol_singu}, \ref{sec:model_trunc} and the model problem
\eqref{eq:model_loc}, \eqref{eq:damp_toy_linf} for the motivation. Thus, to make
the weighted $L^\infty$ energy norms finite and to extract the required damping
near $\xx=0$, the main perturbation and the residual error must vanish
cubically.

We achieve this by a finite-rank decomposition of the perturbation and analytic
low-rank corrections to the residual error. More precisely, we perform the
energy estimates on the main part
$W_1=(\om_1,\eta_1,\xi_1)$, which satisfies $W_1=O(|\xx|^3)$ near $\xx=0$;
the corrected residual error has the same cubic vanishing. The cubic vanishing
order is preserved by the $W_1$ equations. We discuss the method and mechanism 
via a model in Section~\ref{sec:model_low_rank}, define $W_1$ in Section~\ref{sec:finite_pertb},
and establish its vanishing order in Section~\ref{sec:reg_van}. See also
Sections~\ref{sec:decoup}, \ref{sec:decoup_modi} for the corresponding
decomposition of the full Boussinesq system.
In the rest of this section, the singular weighted estimates are applied to the main
component $W_1$, which vanishes as $O(|\xx|^3)$ near the origin.

\begin{remark}[\bf $O(|\xx|^2)$ vs $O(|\xx|^3)$ vanishing]

The cubic vanishing condition is \emph{not the decisive source} of stability. For perturbations vanishing only $O(|\xx|^2)$, one can still follow the analysis of \eqref{eq:model_loc} and Section~\ref{sec:nonloc} to establish linear stability of \eqref{eq:bous_lin11}--\eqref{eq:bous_lin12} near $0$ and the top-order $C^{1/2}$ estimates. The advantage of imposing cubic vanishing on $W_1$ is that it allows us to use more singular weights near the origin and obtain stronger coercivity estimates 
for the local parts, leaving more room to absorb the nonlocal and nonlinear terms for nonlinear stability.

\end{remark}

\subsection{Analytic low-rank correction for cubic vanishing}
\label{sec:model_low_rank}

To overcome the apparent mismatch of vanishing order discussed above, we perform \emph{analytic}  low-rank correction of the perturbation and residual error, using an idea
similar to Taylor expansion. 

We first explain, in a simplified model, how this correction restores the cubic vanishing order needed in the singular weighted estimates. Consider the scalar equation
\beq\label{eq:lin_model}
   \pa_t W=\cL W,
   \qquad
   \cL W=-\bar V(\xx)\cdot\na W+\bar c(\xx)W+\bar F(\xx)P(W),
\eeq
where $\xx=(x,y)\in\R\times\R_+$, $\bar V = (\bar V_1, \bar V_2) \in \R^2$, and $P(W)$ is a rank-one functional. We assume
that $W$, $\bar V_1$, and $\bar F$ are odd in $x$, while $\bar c, \bar V_2$ are even in
$x$, and $\bar V_2(x, 0) \equiv 0$. This is the symmetry structure of the terms in the linearized
equation \eqref{eq:bous_lin11} for $\om$ or $\eta$. Suppose that near the origin
\beq\label{eq:van_F}
   \bar F(\xx)=O(|\xx|^2),
   \qquad
   \pa_{xy}\bar F( 0)\neq 0.
\eeq
Even if the initial data satisfy $W_0=O(|\xx|^3)$, the forcing term
$\bar F P(W)$ may generate an $O(|\xx|^2)$ contribution. Thus
\eqref{eq:lin_model} does not preserve cubic vanishing in general. This
is the obstruction caused by the nonlocal terms in the full linearized equation.

To remove the bad quadratic mode, we choose a smooth  function
$\chi_\dag \in C_c$, odd in $x$, such that
\beq\label{eq:chi_xy}
   \chi_\dag(\xx)=xy+O(|\xx|^3)
   \qquad \text{near } \xx = 0.
\eeq
We decompose $W = W_1 + W_2$ with $W_i$ solving 
 \bseq\label{eq:lin_W}
 \begin{align}
    \pa_t W_1 + \overline V(\xx) \cdot \na W_1 & = \bar c(\xx) W_1 + ( \overline F - \overline F_{xy}(0) \chi_\dag ) P(W_1), 
    \label{eq:lin_W1} \\
  \pa_t W_2 + \overline V(\xx) \cdot \na W_2 & = \bar c(\xx) W_2 +  \overline F P(W_2) + \overline F_{xy}(0) \chi_\dag \cdot P(W_1),   \label{eq:lin_W2}   
\end{align}
 \eseq
 with initial data 
 \[
         W_1 |_{t = 0} = W_0, \quad W_2 |_{t = 0} = 0.
 \]

The correction is chosen so that, by Taylor expansion and \eqref{eq:chi_xy},
\beq\label{eq:F_cor1}
   \bar F(\xx)-\pa_{xy}\bar F( 0 )\chi_\dag(\xx) =O(|\xx|^3), 
   \qquad \text{near } \xx= 0.
\eeq
Thus, if $W_1(0, \xx)=O(|\xx|^3)$, 
the cubic vanishing $W_1(t, \xx) = O(|\xx|^3)$ near $\xx= 0$ is preserved.

The identity \eqref{eq:F_cor1} is an analytic identity, not a
numerical cancellation. Indeed, if
\beq\label{eq:ass_f_bar}
   \bar F(\xx)= \tts{ \sum}_i c_i f_i(\xx),
\eeq
where the basis functions $f_i$ are explicitly known, for example piecewise
polynomials or explicitly defined global functions, then
\beq\label{eq:basis_derivative_F}
   \pa_{xy}\bar F( 0)= \tts{\sum}_i c_i\pa_{xy}f_i( 0).
\eeq
Thus the coefficient of the correction is obtained by differentiating the
represented function. The cancellation in \eqref{eq:F_cor1} is
imposed on the function $\bar F$ itself, not on sampled grid values. This is the
main reason that the cubic vanishing condition used later is exact.

\subsubsection{Duhamel formula} 
\label{sec:duhamel_computer}

We next explain how $W_2$ is estimated. 
Recall $\cL$ from \eqref{eq:lin_model}. 
Following the
method outlined in Section~\ref{sec:frame_novel},  equation
\eqref{eq:lin_W2} admits the Duhamel representation
\beq\label{eq:duhamel_W2_model}
   W_2(t)
   =
   \int_0^t \pa_{xy}\bar F( 0)P(W_1(s))e^{\cL(t-s)}\chi_\dag\,ds .
\eeq
We do not need to construct and estimate the \emph{full} semigroup $e^{\cL s}$ generated by $\cL$, but only $e^{\cL s}\chi_\dag$.

A purely abstract semigroup estimate is not sufficient for our purpose, since
the constants are \emph{not computable}. Instead, we construct an approximate
space-time solution $\widehat g(t,\xx)$ to
\beq\label{eq:lin_g}
   \pa_t g=\cL g,
   \qquad
   g(0,\xx)=\chi_\dag .
\eeq
The function $\widehat g$ is represented by explicit basis functions in space
and time as in \eqref{eq:ass_f_bar}, and its residual $(\pa_t-\cL)\widehat g$ is estimated rigorously in the weighted norms. See Section \ref{sec:how_numerics_enter} 
for further discussion of rigorous bounds. Using $\widehat g$, we define the approximate 
solution for $W_2$
\beq\label{eq:What2_model}
   \widehat W_2(t,\xx)
   =
   \int_0^t \pa_{xy}\bar F( 0 )P(W_1(s))\widehat g(t-s,\xx)\,ds .
\eeq
Then $\widehat W_2$ satisfies
\beq\label{eq:What2_eq_model}
   \pa_t \widehat W_2
   =
   \cL\widehat W_2+\pa_{xy}\bar F( 0 )\chi_\dag P(W_1)+ \cR(W_1),
\eeq
where the \emph{explicit residual} accounts for the error in the initial data and in solving the PDE:
\beq\label{eq:R_W1_model}
\bal
   \cR(W_1)
   &=
   \pa_{xy}\bar F( 0 )P(W_1(t))\big(\widehat g(0)-\chi_\dag\big)  
   +\int_0^t \pa_{xy}\bar F( 0 )P(W_1(s))
       (\pa_t-\cL)\widehat g(t-s)\,ds .
\eal
\eeq
Thus the error in replacing $e^{\cL t}\chi_\dag$ by $\widehat g$ is completely
explicit. We modify the $W_1$-equation \eqref{eq:lin_W1}
accordingly by including the term $-\cR(W_1)$: 
\[
\pa_t   W_1 =  \cL     W_1  - \overline F_{xy}(0) \chi_\dag  P(W_1) 
- \cR(   W_1, t) . 
\]
Hence the decomposition $W = W_1 + \hat W_2 $ remains exact after $\widehat g$ is fixed.

For the singularly weighted estimates on $W_1$, $\cR(W_1)$
\eqref{eq:R_W1_model} must have the same cubic vanishing. Thus, we require the approximate
space-time solution to satisfy the vanishing conditions
\beq\label{eq:g_cond}
   \widehat g(t,\xx)=O(|\xx|^2),\qquad
   \widehat g(0,\xx)-\chi_\dag(\xx)=O(|\xx|^3),\qquad
   (\pa_t-\cL)\widehat g(t,\xx)=O(|\xx|^3) , 
\eeq
near $\xx= 0$. 
The next subsection explains how to establish these conditions rigorously.

\subsubsection{Approximate space-time solution and corrections}\label{sec:numer}

To establish \eqref{eq:g_cond}, the key idea is to perform \emph{analytic} low-rank  correction near $0$ using suitable Taylor expansion. The construction of $\wh g$ 
and corrections for 2D Boussinesq are discussed in detail in  Section \ref{sec:low_rank_cor} and \cite[Section 3]{ChenHou2023b}.

\vs{0.1in}
\paragraph{\bf Step 1. Numerical construction}

We construct the approximate solution to \eqref{eq:lin_g} numerically at discrete times $t_k = k h$ using 
representation \eqref{eq:ass_f_bar} 
\beq\label{eq:g_0form}
   \wh g^{(0)}(t_k, x, y) = 
   \sum\nolimits_{ij} 
   a_{k, ij} B_{1,i}( x) B_{2, j}(y) ,
\eeq
with piecewise polynomial basis functions $B_{1,i}(x), B_{2,j}(y)$ 
\footnote{
 The subscripts $1$ and $2$ in $B_{1,i}$ and $B_{2,j}$ indicate the $x$- and $y$-directions, respectively.
 }
and $t_k \leq T$, where $T$ is the final time up to which we construct the numerical solution.
By interpolating $a_{k,ij}$ in $t$ using piecewise cubic polynomials
$a_{ij}(t)$ (see \cite[Section 3.4]{ChenHou2023b}), whose coefficients are
\emph{linear combinations} of $a_{k,ij}$, 
we extend $\wh g^{(0)}(t_k, x, y)$ to a globally defined function $ \wh g^{(0)}(t, x, y)$ with 
compact support in time $t\in [0, T]$. 
In each time interval $[t_k, t_{k+1}]$, 
$ \wh g^{(0)}(t, x, y)$ satisfies 
\beq\label{eq:g_1form}
\wh g^{(0)}(t, x, y) = \sum\nolimits_{ij} 
a_{ ij}(t) B_{1,i}(x) B_{2,j}(y) .
\eeq
We require that $B_{1,i}(x)$ is odd in $x$ so that $ \wh g^{(0)}$ is odd in $x$ 
and vanishes $O(|\xx|)$ near $\xx=0$. 
\footnote{
To construct the approximation for $\xi$ to the linearized Euler equations, we choose basis $B_{1,i}(x)$ for $\xi$ even in $x$.  
}
We choose the basis functions $B_{1,i}(x), B_{2,j}(y) \in C^{4, 1}$ so that $  \wh g^{(0)} \in C^{4, 1}$.

\vs{0.1in}
\paragraph{\bf Step 2. First analytic correction}

While the \emph{exact} solution to \eqref{eq:lin_g} vanishes $O(|\xx|^2)$ dynamically, due to numerical error, the numerical solution $ \wh g^{(0)}$ in \eqref{eq:g_1form} \emph{need not}
preserve  $O(|\xx|^2)$ vanishing. The first analytic correction is used to restore this
quadratic vanishing.

 Let $\chi_{\dag,1}$ be a $C^{4,1}$ function that is odd in $x$, has a compact support, and satisfies 
\[
 \chi_{\dag, 1}(\xx) = x + O(|\xx|^3) ,
\]
near $x=0$. We perform the first \emph{analytic} correction on $ \wh g^{(0)}$:
\beq\label{eq:g_cor1}
 \wh g^{(1)}(t, \xx ) \teq  \wh g^{(0)}(t, \xx ) - \pa_x  \wh g^{(0)}(t, 0) \cdot \chi_{\dag, 1}( \xx)  , 
\eeq
and obtain $ \wh g^{(1)} \in C^{4, 1}$. We obtain the vanishing order $ \wh g^{(1)}(t, \xx ) = O(|\xx|^2)$
for the same reason as \eqref{eq:F_cor1} and the analytic structure \eqref{eq:g_0form}.  The value $\pa_x  \wh g^{(0)}(t, 0)$ is of the order of the numerical error and is expected to be very small. 
We achieve the first condition in \eqref{eq:g_cond}.

Here, the distinction between \emph{analytic expression} and \emph{numerical
evaluation} is essential.
The correction uses the exact analytic derivative of the explicit representation of $\wh g^{(0)}$. 
From \eqref{eq:g_1form}, 
we derive the exact analytic derivative, determined by the basis functions and coefficients:
\[
\pa_x \wh g^{(0)}(t,0,0)= \sum\nolimits_{ij} 
a_{ij}(t) (\pa_x B_{1,i})(0) B_{2,j}(0),
\]
What is not needed is a separate numerical evaluation of that derivative.

\vs{0.1in}
\paragraph{\bf Step 3. Second analytic correction}

While the \emph{exact} solution to \eqref{eq:lin_g} satisfies 
$g(0, \xx) - \chi_\dag \equiv 0, (\pa_s - \cL) g(s, \xx) \equiv 0$, hence the cubic vanishing in \eqref{eq:g_cond},
due to numerical error, $ \wh g^{(1)}$ constructed in the first two steps \emph{need not} satisfy \eqref{eq:g_cond}. The second analytic correction is used to restore the
cubic vanishing of the error satisfied by the exact solution. 

We introduce the error parts:
\beq\label{def:err}
  \cE_0( \xx) \teq \wh g^{(1)}(0, \xx) - \chi_\dag ,  \quad  \cE(s, \xx) \teq  (\pa_s -\cL)  \wh g^{(1)}(s, \xx) .
\eeq
Since $\wh g^{(1)}(t, \xx)$ has spatial $C^{4, 1}$ regularity by \eqref{eq:g_1form} and 
\eqref{eq:g_cor1}, we obtain the spatial regularity $\cE(s , \cdot) \in C^{3, 1}$ for any $s \geq 0$ 
and $ |\pa^{\b} \cE(s)|$ is bounded in $s$ 
for any $|\beta| \leq 2$.

By definition of $\cL$ in \eqref{eq:lin_model},
the conditions in \eqref{eq:van_F}, and $ \wh g^{(1)}(t, \xx) = O(|\xx|^2)$ near $\xx=0$, 
using a Taylor expansion, we obtain the following leading order terms near $\xx=0$
\beq\label{eq:err_2nd}
  \cE_0( \xx) = \pa_{xy} \cE_0(0) x y + O(|\xx|^3), \quad  \cE(s, \xx)
  = \pa_{xy}\cE(s, 0) x y + O(|\xx|^3) . 
\eeq
The scalar $ \pa_{xy} \cE_0(0)$ and the function $\pa_{xy}\cE(s, 0)$
 depend on $\wh g^{(1)}, \overline F, \overline V, \bar c(x)$ and are of the order of the numerical error. 

  We introduce a function $\chi_{\dag,2} \in C^{4,1}$, that is odd in $x$, has a compact support, and satisfies 
\footnote{ 
For the Boussinesq system, in Sections \ref{sec:finite_pertb}, \ref{sec:appr_vel} and more details in \cite[Section 3]{ChenHou2023b}, the nonlocal correction term $P(f)$ near $\xx=0$ is given by 
$P(f) = c_{\om}(f)=- \pa_{xy}(-\D)^{-1} f(0)$. We choose 
$\chi_{\dag, 2} = f_{\chi,1}$ or $\chi_{1,2}$ defined in \eqref{eq:cutoff_near0_all}. 
We have  $\chi_{\dag, 2}(\xx) = \D( \f{1}{6} x y^3 \kp(x) \kp(y))$ with $\kp(x)$ being smooth cutoff function supported near $x=0$ and $\kp(x)=1$ near $x=0$, and obtain 
$\pa_{xy}(-\D)^{-1} \chi_{\dag, 2}(0) = 0$ and $\chi_{\dag, 2}(\xx) = x y + O(|\xx|^3)$ near $\xx=0$.
}
\beq\label{def:chi2}
\chi_{\dag, 2}(\xx) = x y + O(|\xx|^3), \quad  P(\chi_{\dag, 2}) = 0.
\eeq

To further eliminate the quadratic error, we perform the second \emph{analytic} correction
\beq\label{def:g2}
  \wh g^{(2)}(t, \xx) = \wh g^{(1)}(t, \xx) + a(t) \chi_{\dag, 2}.
\eeq
From \eqref{eq:g_0form}, \eqref{eq:g_cor1}, and \eqref{def:g2},
$\wh g^{(2)}(t, \xx)$ has the representation \eqref{eq:ass_f_bar} with $C^{4, 1}$ basis functions. 

\vs{0.1in}

\paragraph{\bf Constructing and solving ODE for $a(t)$ analytically}

We choose $a(t)$ to eliminate the quadratic parts near $\xx=0$ in \eqref{eq:err_2nd} by imposing the following \emph{analytic} condition on $a(t)$:
\beq\label{eq:van_g2}
  \pa_{xy} (   \wh g^{(2)} - \chi_\dag )|_{s=0, \xx=0 } = 0,
  \quad \big( \pa_{xy}  (\pa_s - \cL) \wh g^{(2)}\big)(s, 0) = 0.
\eeq

Due to symmetry, using \eqref{eq:err_2nd}, \eqref{def:g2}, we rewrite the above analytic conditions equivalently:
\beq\label{eq:ODE1}
 a(0) + \pa_{xy} \cE_0(0) = 0,
 \quad  \pa_{xy} ( \pa_s - \cL) a(s) \chi_{\dag, 2} |_{ \xx=0} + \pa_{xy} \cE(s, 0) = 0.
\eeq

For any function $f \in C^3$, odd in $x$, and $f = O(|\xx|^2)$, from \eqref{eq:lin_model}, we derive 
\[
   \pa_{xy} ( \cL  f )(0)=  (\bar c(0) - \pa_x \overline V_1(0) - \pa_y \overline V_2(0) ) \pa_{1 2} f(0)
 +  \pa_{1 2}\overline F(0) \cdot P( f ).
\]

Using the above formula,  $\pa_{xy} \chi_{\dag, 2}(0) = 1, P(\chi_{\dag, 2}) = 0$ by \eqref{def:chi2}, 
we rewrite equations \eqref{eq:ODE1} as
\beq\label{eq:ODE2}
\bal
 0  & = \f{d}{d s} a(s) - a(s) \pa_{xy} \cL \chi_{\dag, 2}(0) 
  + \pa_{xy} \cE(s, 0) \\
 & =  \f{d}{d s} a(s) - a(s) \big( \bar c(0 ) - \pa_x \overline V_1(0) - \pa_y \overline V_2(0)   \big)
 + \pa_{xy} \cE(s, 0), \\
 a(0) & = - \pa_{xy} \cE_0(0).
 \eal
\eeq

The above relation gives a \emph{linear} ODE in $a(s)$, which we solve \emph{exactly and analytically}
\beq\label{def:a}
  a(t) = - e^{\bar \lam t}  \pa_{xy} \cE_0(0)
- \int_0^t e^{ \bar \lam (t- s) }  \pa_{xy} \cE(s, 0 )  d s,
\quad \bar \lam = \bar c(0 ) - \pa_x \overline V_1(0) - \pa_y \overline V_2(0).
\eeq

Since $\cE(s,{\bf x}) $ has $C^{3,1}$ spatial regularity in ${\bf x}$, and $\na^{\beta} \cE(s, {\bf x})$ is bounded in $s$ for $|\beta| \leq 2$, 
we can define the above integral \emph{analytically}. 

The formula \eqref{def:a} is the key point in the second correction:  once the defect coefficients are determined, $a(t)$ is obtained by solving the ODE \eqref{eq:ODE2} \emph{analytically}, 
rather than by another numerical approximation. With the above correction, we construct $\hat g = \hat g^{(2)}$  and establish \eqref{eq:g_cond}.

\subsubsection{Vanishing order of error for stream function}

To approximate the velocity profile $\bar \uu=\na^\perp(-\D)^{-1}\bar\om$, we construct an approximate stream function $\bar\phi^N$ by solving the Poisson equation $-\D\phi=\bar\om$.  Due to numerical error, there is a local residual error 
\[
\bar\e=\bar\om+\D\bar\phi^N .
\]
We represent $\bar\phi^N$ using basis functions similar to  \eqref{eq:ass_f_bar}. Their odd symmetry in $x_1$ implies only $\bar\e=O(|x_1|)$, whereas controlling the nonlocal error $\uu(\bar\e)=\na^\perp(-\D)^{-1}\bar\e$ requires quadratic vanishing. As in 
the previous parts in Section \ref{sec:model_low_rank}, e.g., \eqref{eq:g_cor1}, we use an analytic low-rank  correction to ensure $\bar\e=O(|x|^2)$; see the end of Section \ref{sec:ASS_multi}. We further decompose $\bar\e=\bar\e_1+\bar\e_2$, with a rank-one component $\bar\e_2$, to improve the quadratic vanishing of $\bar\e$ to cubic vanishing of $\bar\e_1$. We then apply nonlocal functional inequalities to $\bar\e_1$ to estimate $\na^\perp(-\D)^{-1}\bar\e_1$; see Section \ref{sec:comb_vel_err}.

\subsubsection{How numerical values enter the weighted estimates}
\label{sec:how_numerics_enter}

We briefly explain
how numerical information enters the estimates once the exact local cancellation has been
enforced analytically.  

We focus on the region near $\xx= 0$, where the weight is singular and the vanishing order is essential. Suppose that a corrected
function $f$ satisfies $f(\xx)=O(|\xx|^3)$ near $x=0$. Then
\[
   |\xx|^{-3}|f(\xx)|
   \les
   \|\na^3 f\|_{L^\infty(|\xx|\le 1)} , \quad 
   \mbox{for} \  |\xx|\le 1 \, .
\]
Therefore the singular weighted estimate reduces to an ordinary $C^3$ bound on
the corrected function. Such $C^3$ bounds are obtained directly from the basis representation. If
\beq\label{eq:f_basis_C3}
   f(\xx)= \tts{\sum}_{i,j}a_{ij} B_{1,i}(x)B_{2, j}(y),
\eeq
then
\beq\label{eq:f_basis_C3_derivative}
   \na^3 f(\xx)
   =
   \tts{\sum}_{i,j}a_{ij}\na^3\big(B_{1,i}(x)B_{2,j}(y)\big).
\eeq

Rigorous bounds on the coefficients $a_{ij}$ and on the derivatives of the basis functions give rigorous piecewise $C^3$ estimates for $f$. Near $x=0$, numerical errors may affect the coefficients $a_{ij}$, but they do not affect the \emph{exact local cancellation}, because the cancellation is imposed after the coefficients have been fixed and is verified at the level of the represented functions.  Together with the high-order regularity bounds obtained from the basis representation, the exact vanishing allows us to control the singularly weighted estimates rigorously. In particular, the cubic vanishing order required by the singular weights is stable under the computer-assisted construction, which enables robust weighted energy estimates around the numerically constructed profiles, such as $\bar\om,\bar\th$ in \eqref{eq:bous_lin11}, and the approximate solution $\hat g$.

In Appendix \ref{app:idea_num_est}, we provide a self-contained outline of some key methods  for  rigorous weighted estimates of functions and residual error.
 As examples, we  discuss the weighted estimates of $\hat g$ in \eqref{eq:g_cond}.  
 More detailed methods for rigorous numerical computations are developed in \cite{ChenHou2023b}.

\subsection{Ideas for Weighted H\"older estimates}\label{sec:idea_EE_calpha}

Since   $ \na \uu = \na \na^{\perp}(-\D)^{-1}\om $ in \eqref{eq:bous_lin12} is not bounded from $L^{\infty} \to L^{\infty}$, to close the estimates, we perform weighted $C^{\al}$ estimates.

The basic idea is to regard a weighted $C^\al$ estimate as a weighted $L^\infty$ estimate for
differences on $\R^2_+ \times \R^2_+$. More precisely, 
for a variable $f$, we estimate the weighted difference $ \big( f\vp(x)-f\vp(z)\big) g(x-z)$, 
where $\vp$ is the spatial weight and $g(h)$ is a H\"older weight, e.g. $g(h)=|h|^{-\al}$.  The following  identity explains where the damping terms in the weighted $C^{\al}$ estimate come from.

\begin{lem}\label{lem:hol_comp}
Suppose that $f$ satisfies
\beq\label{eq:hol_comp0}
\pa_t f  + b(x) \cdot \na f = c(x) f(x) + \cR, \quad x  \in \R^2_+.
\eeq
Given some weights $g(x_1, x_2)>0$ even in $x_1, x_2$, 
$g \in C^1(\R^2 \bsh \{0 \})$, and $\vp$, we denote the operator $\d$ and function $F$
\[
\d(p)(x, z) \teq p(x) - p(z), \quad
F(x, z, t)  = \d(f \vp)(x, z) g(x-z),
\quad d(x) = c(x) + \f{ b \cdot \na \vp}{\vp}, \quad x, z \in \R^2_+.
\]
Then we have
\beq\label{eq:hol_comp}
\bal
\pa_t F + ( b(x) \cdot \na_x + b(z)  \cdot \na_z ) F
& = \left( d(x) + \f{ (b(x) - b(z)) \cdot (\na g) (x-z) }{g(x-z)}  \right) F \\
+ & (d(x) - d(z))g(x-z)  (f \vp)(z)  + \d(\cR \vp)(x, z) g(x-z) \;.
\eal
\quad
\eeq
\end{lem}

The proof is a direct calculation and is deferred to Appendix~\ref{app:lem:hol_comp}. We will
use the first term on the right-hand side of \eqref{eq:hol_comp} as the main damping term. This damping has two sources: 
\[
        d(x) = c(x) + \f{ b \cdot \na \vp}{\vp} ,  \quad         \f{ (b(x)-b(z))\cdot (\na g)(x-z)}{g(x-z)} .
\]
 The first one comes from the singular spatial weight $\vp$. The second one comes from the H\"older weight $g(x-z)$. 
The term involving $d(x)-d(z)$ is treated as a lower-order
source term and will be controlled by a stronger weighted $L^\infty$ estimate.

\vs{0.05in}
\paragraph{\bf Weighted $C^{\al}$  stability of model problem }
To illustrate the damping mechanism in the weighted $C^\alpha$ estimate, we apply Lemma~\ref{lem:hol_comp} to the $\eta$-equation in the model problem \eqref{eq:model_loc}.
Denote 
\[
\bal
& \vp_2   = |x|^{-\g_2},  \  g(h) = |h|^{-\al}, \ h \in \R^2, \quad 
b(x) =(a_1 x_1, a_2 x_2), \quad F = ( (\eta \vp_2)(x) -  ( \eta \vp_2)(z) )g(x-z)
  \eal
\]
for $x = (x_1, x_2), z = (z_1, z_2) \in \R^2_{+}$.  Using the identity \eqref{eq:damp_toy_linf}, we obtain 
\[
d(x) = a_3 + \f{ b(x) \cdot \na \vp_2}{\vp_2} = a(\g_2).
\]

Since $h_i \pa_i g(h) =  - \al \f{h_i^2}{ |h|^2} g $, the H\"older weight $g$ gives the additional damping
  \beq\label{eq:damp_toy_hol}
\bal
 & \f{ (b(x) - b(z)) \cdot (\na g)(x-z) }{g(x-z)}  = - \al \f{  a_1(x_1-z_1)^2 + a_2(x_2-z_2)^2}{ |x-z|^{2}}  \teq e(\al, x, z).  
 \eal
\eeq
Since $e(\al, x, z)<0$, $ e(\al, x, z) F$ in \eqref{eq:hol_comp} is also a damping term. 
Using Lemma \ref{lem:hol_comp} with $\cR = 0$, we obtain 
\[
\bal
\pa_t F + ( b(x) \cdot \na_x + b(z) \cdot \na_z) F
 & = ( a(\g_2)(x) + e(\al, x, z) ) F  + ( a(\g_2)(x) - a(\g_2)(z) ) g(x-z) (\eta \vp_2)(z).
\eal
\]

The last term is the main point that requires some care. Although $a(\g_2)(x)$ is bounded, it is
not sufficiently regular in the same weighted H\"older norm. Thus we do not estimate this term
only by the $C^\al$ norm of $\eta\vp_2$. Instead, we use the stronger weighted $L^\infty$
estimate $\| \eta |x|^{-\g} \|_{\infty}$ with $\g=\g_2+\al$.
This choice is natural because $\vp_2 |x|^{-\al}=|x|^{-\g_2-\al}=|x|^{-\g}.$
Moreover, since $a(\g_2)(x)|x|^\al$ and $|x|^\al$ are $C^\al$, and since
$a(\g_2)(x)\in L^\infty$, we have
\[
        \left|
        \big(a(\g_2)(x)-a(\g_2)(z)\big)g(x-z)|z|^\al
        \right|
        \le C_\al .
\]

Combining this bound with 
\[
a(\g_2) \leq a_3 - a_1 \g_2, \quad  e(\al,x, z) \leq -a_1 \al
\]
 by \eqref{eq:damp_toy_linf}, \eqref{eq:damp_toy_hol}, we obtain
\[
\tf{d}{dt} || F ||_{L^{\inf}(x, z)}
\leq (a_3 - a_1 (\g_2 + \al) ) || F ||_{L^{\inf}(x, z)}
+ C_{\al} || \eta |z|^{-\g} ||_{\inf} ,
\]
and yield the weighted $C^\alpha$ estimate for $\eta$. Similarly, we obtain a weighted $C^{\al}$ estimate for $\om$. Since $a_1 = a_3 = \f{1}{2}$, choosing $\g_2 + \al  = \g> 1$ and combining the above estimate and \eqref{eq:model_linf2}, we establish a stability estimate for model 
\eqref{eq:model_loc} in a combination of weighted $L^{\inf}$ and $C^{\al}$ spaces.

\subsubsection{Anisotropic estimates and the most difficult scenario}\label{sec:aniso_est}

Denote $ f=(\om, \eta, \xi)$. Motivated by the above analysis, we will perform  $L^{\inf}$ estimates for $f_{i} \vp_i$ and $C^{1/2}$ estimates for $f_{i} \psi_i$.
We choose the $L^{\inf}$ weight $\vp_i$ to be at least $|x|^{-\f{1}{2}}$ more singular near $x =0$ than the H\"older weight $\psi_i$. Since damping coefficients similar to $a(\g_2)$ in \eqref{eq:damp_toy_hol} in the model problem are not $C^{1/2}$ near $x=0$, 
\footnote{
If $| \na^i f | \les |x|^{p_1 -i}, i\leq 1  $, 
then $f |x|^{-p_1} \les 1$. We require $p_2 \leq p_1 - \f{1}{2}$ to ensure that 
$f |x|^{-p_2 } \in C^{1/2}$ locally. 
}
the ratio $\f{\psi_i}{ \vp_i}$ provides a vanishing power $|x|^{\f{1}{2}}$, which compensates the low regularity of the damping coefficients.  See \eqref{wg:linf_decay},\eqref{wg:hol}.

From the above analysis of the model problem, \eqref{eq:damp_toy_linf}, and \eqref{eq:damp_toy_hol}, we see that the estimate is anisotropic in $x$ and $y$. In the near field in $\R^2_{++}$, from \eqref{eq:damp_toy_linf}, if $y/x$ is not small, we get a much larger damping term. Moreover, the solution $(\bar \om, \na \bar \th)$ decays in the far-field.
Thus the most difficult region in $\R^2_{++}$ for the analysis is a sector near the boundary, 
 for instance
\[
        \S_S=\{(x,y): |(x,y)|\le 2,\; y/x\le 0.1\}.
\]

From \eqref{eq:damp_toy_hol}, we also get a much larger damping factor if $|x_2 - z_2|  / |x_1 - z_1|$ is large. Thus, the most difficult case in the H\"older estimate is the estimate of pairs $(x, z)$ with 
\[
|x_2 - z_2| \ll |x_1 - z_1|.
\]

\subsection{Estimate the nonlocal terms $\na \uu$}\label{sec:nonloc}

In this section, we discuss the idea of weighted $C^{1/2}$ estimates of the nonlocal term $\na \uu$ and consider the stability of a model problem with nonlocal terms. The key point is that the most singular part of $\na \uu$ must be estimated
directly by sharp $C^{1/2}$ estimates, while the more regular part can be
treated perturbatively by the finite-rank method in Section~\ref{sec:finite_rank}.
Thus,  the sharp H\"older estimates for $\na \uu$ in Section~\ref{sec:sharp} are essential.

\subsubsection{Weighted estimates}\label{sec:nonloc_wg_idea}

We focus on the term $u_x\bar\th_x$ in the main system \eqref{eq:bous_main}. This is the
most difficult nonlocal term to estimate. Other nonlocal terms, such as
$u\bar\om_x$ and $u\bar\th_{xx}$ in \eqref{eq:bous_main}, involve 
the more regular nonlocal term  $u$  or contain smaller coefficients.
\footnote{
In the leading order system for the $C^{1,\al}$ singular solution \cite{chen2019finite2},
$u_x\bar\th_x$ is also the main nonlocal term. Other nonlocal terms involving
$\uu,\na\uu$ contain a small factor $\al$.
}

We first isolate the most singular part of $u_x$. Denote by $u_x(x,a,b)$ the localized
version of $u_x$:
\beq\label{eq:ux_loc_idea}
\bal
 u_{x} (x, a, b)
 &\teq
 - \f{1}{\pi} P.V. \int_{|x_1-y_1| \leq a,\; |x_2-y_2| \leq b}
 K_1(x-y) W(y)\,dy,
 \qquad
 K_1(s)=\f{s_1s_2}{|s|^4},
\eal
\eeq
where $W$ is the odd extension of $\om$ in $y$ from $\R^2_+$ to $\R^2$; see
\eqref{eq:ext_w_odd}. We decompose
\[
        u_x = u_{x,S}+u_{x,R}, \qquad u_{x,S}(x)=u_x(x,a,b),
\]
where $u_{x,S}$ captures the singularity near $y=x$, and $u_{x,R}$ is more regular.
The H\"older estimate for the regular part $u_{x,R}\psi$ is obtained using the method
discussed in Section~\ref{sec:u_comp_idea}.

The top-order nonlocal term is the singular part $u_{x,S}\psi$. We observe that the commutator 
\beq\label{eq:ux_commu0}
\bal
{[}u_x(\cdot, a, b), \psi {]}(\om)
&\teq
u_x(\om)(x,a,b)\psi(x)-u_x(\om\psi)(x,a,b)  \\
&=
-\f{1}{\pi}
\int_{|x_1-y_1| \leq a,\; |x_2-y_2| \leq b}
K_1(x-y)W(y)\big(\psi(x)-\psi(y)\big)\,dy
\eal
\eeq
is more regular than the original singular integral, because the difference
$\psi(x)-\psi(y)$ cancels part of the kernel singularity. Therefore we use the
decomposition
\beq\label{eq:ux_commu1}
u_x(\om)(x,a,b)\psi
=
u_x(\om\psi)(x,a,b)+[u_x(\cdot,a,b),\psi](\om)
\teq I_1+I_2 .
\eeq

The first term $I_1$ is the most singular part. It is estimated using the sharp
H\"older inequalities in Section~\ref{sec:sharp}.  
Indeed, if $\om\vp\in L^\infty$ for a stronger weight $\vp$, then
$K_1(x-y)(\psi(x)-\psi(y))$ has only a $|x-y|^{-1}$ singularity, and $I_2$ is
log-Lipschitz. The $C^{1/2}$ estimate of $I_2$ can be bounded in terms of
$\|\om\vp\|_{L^\infty}$ by the methods in Section~\ref{sec:u_comp_idea}. If the
localization parameters $a$ and $b$ are small, this estimate also gives a small factor of order $(\max(a,b))^{1/2}$.

\begin{remark}[\bf Constant independent of $\psi$]
It is crucial to commute $\psi$ with the nonlocal operator as in \eqref{eq:ux_commu1}.
After this commutation, the sharp $C^{1/2}$ estimate for
$u_x(\om\psi)(x,a,b)$ has a constant independent of $\psi$. This allows the large
damping generated by the singular weight, discussed in Section~\ref{sec:outgo},
to dominate the top-order nonlocal terms. This is one of the key mechanisms that leads to stability at the top order.
\end{remark}

\subsubsection{Top-order coercivity estimates}\label{sec:hol_singu}

We  isolate the most singular regime in the top-order weighted $C^{1/2}$ estimate
and explain the coercivity estimates for the leading part of the linearized operator.
This discussion is meant to connect the weighted H\"older identity in Lemma \ref{lem:hol_comp},  the anisotropic scenario from Section~\ref{sec:aniso_est},
and the sharp H\"older estimates from Section~\ref{sec:sharp}.

We consider the model obtained by dropping the nonlinear and error terms in
\eqref{eq:bous_lin11}, \eqref{eq:bous_lin12}:
\beq\label{eq:model_hol_sing0}
\bal
\om_t
&= -(\bar c_l \xx+\bar \uu)\cdot\na \om+\bar c_{\om}\om+\eta,\\
\eta_t
&= -(\bar c_l \xx+\bar \uu)\cdot\na \eta
 +(2\bar c_{\om}-\bar u_x)\eta
 -(u_x-\hat u_x)\bar\th_x-(v_x-\hat v_x)\bar\th_y .
\eal
\eeq
We omit the weak coupling term $\bar v_x\xi$ and the more regular nonlocal terms,
such as $\uu$, since these terms are lower order in the limiting singular
$C^{1/2}$ estimate discussed below.

Here $\hat u_x,\hat v_x$ are the leading finite-rank approximations of $u_x,v_x$
near $\xx=0$:
\beq\label{eq:model_hol_appr}
\hat u_x =
\big(u_x(0)+C_{u_x}(\xx)\cK_{00}\big)\chi(\xx),
\qquad
\hat v_x = C_{v_x}(\xx)\cK_{00}\chi(\xx).
\eeq
They satisfy
\beq\label{eq:model_hol_vanish}
u_x-\hat u_x = O(|\xx|^{5/2}),
\qquad
v_x-\hat v_x = O(|\xx|^{5/2}),
\eeq
for $\om=O(|\xx|^{5/2})$ near $\xx=0$; see Section~\ref{sec:reg_van}.
The approximation \eqref{eq:model_hol_appr} is a simplified version of the
finite-rank approximation of the velocity developed in Section~\ref{sec:appr_vel}.

\vspace{0.1in}
\paragraph{\bf{Goal of the estimates and heuristic}}

We focus on the most difficult H\"older configuration identified in
Section~\ref{sec:aniso_est}: the pair $(\xx,\zz)$ satisfies $x_2=z_2$.
We also restrict to the most singular scenario where $|\xx-\zz| \to 0$. Lemmas~\ref{lem:holx_ux} and \ref{lem:holx_uy} provide the sharp H\"older estimates of the nonlocal terms. The diagonal local terms identified in Lemma \ref{lem:hol_comp} provide damping.
The goal of the following estimates is to prove that this damping dominates the nonlocal terms and other 
off-diagonal local terms  and to verify the linear stability condition
\eqref{eq:PDE_diag} for the model \eqref{eq:model_hol_sing0}, in this singular scenario $x_2=z_2, |\xx-\zz| \to 0$.


 In this regime, heuristically, the H\"older difference is essentially a half derivative $D$ in $x_1$-direction. 
If $D$ acts on a Lipschitz term, its estimate vanishes as $|x-z| \to 0$. If $D$ acts on the nonlocal terms, e.g. $u_x - \hat u_x$, since $\hat u_x$ is more regular and $[u_x]_{C_x^{1/2}}$ can be bounded using Lemma \ref{lem:holx_ux}, we treat it as  $2.55 [ \om \psi_1 ]_{C_x^{1/2} }$. If $D$ applies to the local term, we use the energy to bound it.

To keep the exposition focused, we first describe the case where $\xx$ is away from
the origin. The case where $\xx$ is close to $0$ is more technical because the weights
are singular, but it is actually easier from the viewpoint of stability: the vanishing
orders in \eqref{eq:model_hol_vanish} provide additional smallness. See
Remark~\ref{rem:dp_near0} for the estimates near $\xx=0$.

Following Section~\ref{sec:aniso_est}, we choose the H\"older weights
\eqref{wg:hol} so that
\beq\label{wg_prop:hol}
\psi_1 \sim |\xx|^{-2},
\qquad
\psi_2 \sim |\xx|^{-5/2}
\eeq
near $\xx=0$. 
For $f\in C^{2,1}$,  $|f(\xx)|\les |\xx|^3$ ensures
$f\psi_i\in C^{1/2}$. See Sections~\ref{sec:vanish} and \ref{sec:reg_van}. Denote
\[
   \bar b(\xx)=\bar c_l\xx+\bar\uu(\xx)
   =(\bar c_l x+\bar u,\,\bar c_l y+\bar v).
\]

\vspace{0.1in}
\paragraph{\bf{Weighted $C^{1/2}$ equations}}

For a pair $(\xx,\zz)$, applying Lemma~\ref{lem:hol_comp} to $\om\psi_1$
and $\eta\psi_2$ gives
\bseq\label{eq:model_hol_sing}
\beq
\bal
 & \pa_t \d(\om\psi_1)g_1
 +\big(\bar b(\xx)\cdot\na_{\xx}+\bar b(\zz)\cdot\na_{\zz}\big)
   \big(\d(\om\psi_1)g_1\big)
 = \bar e_1(\xx,\zz)\d(\om\psi_1)g_1+B_1(\xx,\zz),\\
 & \pa_t \d(\eta\psi_2)g_2
 +\big(\bar b(\xx)\cdot\na_{\xx}+\bar b(\zz)\cdot\na_{\zz}\big)
   \big(\d(\eta\psi_2)g_2\big)
 = \bar e_2(\xx,\zz)\d(\eta\psi_2)g_2+B_2(\xx,\zz),
\eal
\eeq
where $g_i=g_i(\xx-\zz)$ and $\d f=f(\xx)-f(\zz)$. The terms $B_i$ are the
non-damping terms:
\beq
\bal
B_1(\xx,\zz)
&\teq \d(\eta\psi_1)g_1(\xx-\zz)
+\d(\bar d_1)g_1(\xx-\zz)(\om\psi_1)(\zz),\\
B_2(\xx,\zz)
&\teq
-\d\!\left((u_x-\hat u_x)\bar\th_x\psi_2
+(v_x-\hat v_x)\bar\th_y\psi_2\right)g_2(\xx-\zz)
+\d(\bar d_2)g_2(\xx-\zz)(\eta\psi_2)(\zz).
\eal
\eeq
\eseq
Here
\beq\label{eq:model_hol_sing1}
\bal
\bar d_1 &= \bar c_{\om}+\f{\bar b\cdot\na\psi_1}{\psi_1},
\qquad
\bar d_2 = 2\bar c_{\om}-\bar u_x+\f{\bar b\cdot\na\psi_2}{\psi_2},\\
\bar e_i(\xx,\zz)
&= \bar d_i(\xx)
+\f{(\bar b(\xx)-\bar b(\zz)) \cdot (\na g_i)(\xx-\zz)}
       {g_i(\xx-\zz)} .
\eal
\eeq
The coefficients $\bar e_i$ 
contain both the
damping from the spatial weights $\psi_i$ and from the H\"older weights
$g_i$, as explained in Section~\ref{sec:idea_EE_calpha}. We aim to control the
terms $B_i$ by this damping.

\vspace{0.1in}
\paragraph{\bf{Estimates of $B_i$}}

We briefly discuss the estimates of $B_i$.  For the $u_x$ term in $B_2$, write
\[
\d((u_x-\hat u_x)\bar\th_x\psi_2)g_2
=
\d((u_x-\hat u_x)\psi_1 A)g_2,
\qquad
A(\xx)=\bar\th_x(\xx)\cdot \tf{\psi_2}{\psi_1}(\xx).
\]
When $\xx$ is away from $0$ and $|\xx-\zz| \to 0$ , the part where the difference
falls on $A$ is regular and vanishes.
The leading contribution comes from $\d( (u_x-\hat u_x)\psi_1) g_2 A(\xx)$.   Using the commutator decomposition \eqref{eq:ux_commu1}, the regularity of
$\hat u_x$, and the fact that the commutator $[u_x,\psi_1](\om)$ is more regular,
the leading singular part is estimated by applying Lemma~\ref{lem:holx_ux} to
$u_x(\om\psi_1)$: 
\bseq\label{eq:model_hol_sing_S2}
\beq
|\d(u_x(\om\psi_1))g_2 A(\xx)|
\leq
2.55[\om\psi_1]_{C_x^{1/2}}
\big|\bar\th_x \cdot \tf{\psi_2}{\psi_1}(\xx)\big|,
\eeq
where we used $g_1(\xx-\zz) = g_2(\xx-\zz) = |x_1 -z_1|^{-1/2}, \forall \, x_2 = z_2$ \eqref{wg:hol}. Similarly, Lemma~\ref{lem:holx_uy} gives
\beq
|\d(v_x(\om\psi_1))g_2|
\big|\bar\th_y \cdot \tf{\psi_2}{\psi_1}(\xx)\big|
\leq
2.53\|\d(\om\psi_1)g_1\|_{L^\infty}
\big|\bar\th_y \cdot \tf{\psi_2}{\psi_1}(\xx)\big| .
\eeq
\eseq

The remaining terms in $B_i$ \eqref{eq:model_hol_sing} are more regular, and their weighted estimates vanish in this limiting top-order estimate. For example,  
for  $\xx$ away from the origin, we have
\[
  \d(\bar d_i)g_i(\xx-\zz)(f_i\psi_i)(\zz) \to 0, \quad  |\d((u-\hat u)\bar\om_x\psi_1)g_1| \to 0, 
\qquad \text{as } |\xx-\zz|\to0.
\]
In the energy estimates in Section~\ref{sec:EE}, these terms are retained and controlled by the energy.

Finally, the coupling term in the $\om$ equation satisfies, in the limit
$|\xx-\zz|\to0$ with $x_2=z_2$,
\beq\label{eq:model_hol_sing_S1}
|\d(\eta\psi_1)g_1|
\leq
\big|\tf{\psi_1}{\psi_2}(\xx)\big| \cdot [\eta\psi_2]_{C_x^{1/2}}.
\eeq

\paragraph{\bf{Summary of the estimates}}

\begin{figure}[t]
   \centering
      \includegraphics[width = 0.65\textwidth]{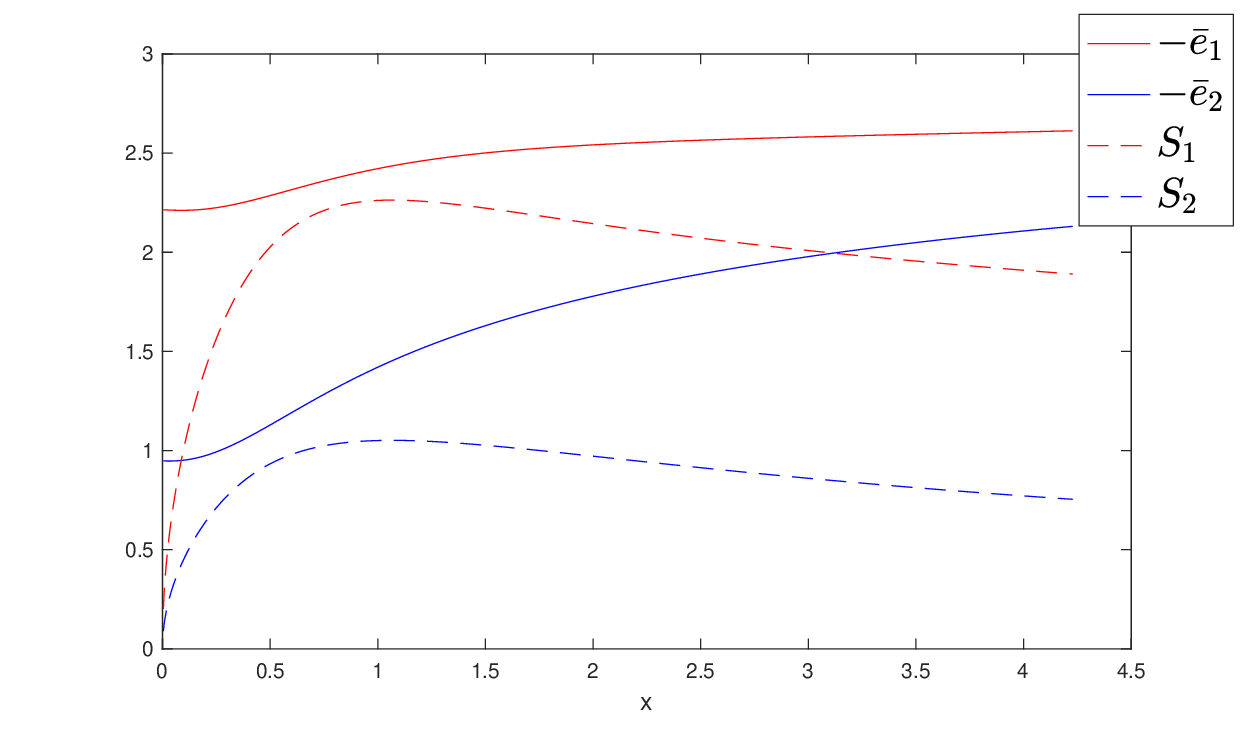}
      \vspace{-0.2in}
        \captionsetup{width=0.9\textwidth}
      \caption{Coefficients $\bar e_1,\bar e_2$ of damping terms and estimates of the bad terms $B_1,B_2$.}
      \label{fig:damp_hol_sing}
\end{figure}

For $x_2=z_2$ with $|x_1-z_1| \to 0$, the damping coefficients
$\bar e_i$ in \eqref{eq:model_hol_sing1} simplify to
\[
\bar e_i(\xx)
\teq
\bar d_i(\xx)- \lim_{z_1 \to x_1} \f12\f{\bar b_1(\xx)-\bar b_1(\zz)}{x_1-z_1}
=
\bar d_i(\xx)-\f12\pa_1\bar b_1(\xx)
=
\bar d_i(\xx)-\f12(\bar c_l+\bar u_x(\xx)).
\]
We choose
\beq\label{wg:hol_mu1}
\mu_1=0.668
\eeq
to balance the two H\"older components $\d(\om\psi_1)g_1$ and
$\d(\eta\psi_2)g_2$ in the energy. In the limit $x_2=z_2$,
$|x_1-z_1|\to0$, estimate 
\eqref{eq:model_hol_sing_S1} for $B_1$ 
and \eqref{eq:model_hol_sing_S2} for $B_2$ imply
\beq\label{eq:model_hol_sing2}
\bal
  \limsup_{z_1 \to x_1 }
   |B_1(\xx, \zz)|
  & \leq S_1(\xx)\cdot \mu_1[\eta\psi_2]_{C_x^{1/2}},
  && \limsup_{z_1 \to x_1 } \mu_1|B_2(\xx,\zz)|
  \leq S_2(\xx)[\om\psi_1]_{C_g^{1/2}},\\
 S_1(\xx) & =\B|\f{\psi_1}{\mu_1\psi_2}(\xx)\B|,
&&
S_2(\xx)  =\mu_1(2.55|\bar\th_x|+2.53|\bar\th_y|)
\f{\psi_2}{\psi_1}(\xx).
\eal
\eeq
Therefore, by \eqref{eq:PDE_diag} and Lemma~\ref{lem:PDE_stab}, stability follows
if the damping dominates the bad terms:
\beq\label{eq:model_hol_sing3}
\bar e_i(\xx)+S_i(\xx)\leq -c_i
\eeq
for some $c_i>0$. In Figure~\ref{fig:damp_hol_sing}, we illustrate the estimates by plotting
$S_1,S_2$ in \eqref{eq:model_hol_sing2}, together with the damping coefficients
$-\bar e_1,-\bar e_2$, on the boundary $y=0$ for $x<5$. The figure visualizes
the dominance of the damping terms in the most singular scenario.
\footnote{
The plotted bounds can be verified rigorously using piecewise sharp bounds of the
profile and interval arithmetic. They are included only to illustrate the estimates
and are not used in the proof.
}

The estimates above are limiting estimates for the full top-order $C^{1/2}$ stability
estimate. As $|\xx-\zz|\to0$, the more regular terms omitted in the model
\eqref{eq:model_hol_sing} vanish at the level of the $C^{1/2}$ quotient. Thus,
in this subsection we track only the leading contributions and verify the stability
of the top-order terms, without using the estimates for the more regular terms.

For $|\xx-\zz|$ sufficiently small but with a different ratio
$|x_2-z_2|/|x_1-z_1|$, Section~\ref{sec:aniso_est} shows that the H\"older weight
gives additional damping through $\f{ (\bar b(x)- \bar b(z)) \cdot \na g_i}{g_i}$,
as in \eqref{eq:model_hol_sing1}. Using the sharp H\"older estimates in
Section~\ref{sec:sharp}, one verifies that the damping terms also dominate in
these cases. The case where $|\xx-\zz|$ is not sufficiently small, where the
omitted more regular terms must be retained, is handled in the full energy estimates
in Section~\ref{sec:EE}. For larger $|\xx-\zz|$, the H\"older estimates for
$\na\uu$ have better constants in Lemmas~\ref{lem:holx_ux}--\ref{lem:holy_ux};
see the discussion below \eqref{eq:EE_tau1}.

\begin{remark}
The size of the stability factors $-c_i-(\bar e_i+S_i)$ is closely related to the
choice of weights. See the discussion of the stability factors in \eqref{eq:EE_target1}
and the related remark.
\end{remark}

\subsubsection{A model problem for localized velocity} \label{sec:model_trunc}

We next introduce a simplified model to explain the weak coupling between the weighted
$L^\infty$ and $C^{1/2}$ estimates in the full energy argument. This model also
motivates the localization of velocity in \eqref{eq:ux_loc_idea} and the finite-rank
approximation of the more regular nonlocal terms. Consider
\beq\label{eq:model_trunc}
\om_t = -d(\xx)\om + a(\xx)u_x(\om,\e)(\xx).
\eeq
Here $u_x(\om,\e)$ denotes the localized velocity $u_x(\om,a,b)$ in
\eqref{eq:ux_loc_idea} with $a=b=\e$. In this model, we omit the transport terms
for simplicity. In the weighted energy estimate, their contribution is absorbed into
the damping terms, similarly to $-d(\xx)\om$; see Sections~\ref{sec:model1} and
\ref{sec:hol_singu}. The nonlocal term $u_x(\om,\e)$ models the singular
nonlocal terms involving $\na\uu$ in \eqref{eq:bous_lin11}, \eqref{eq:bous_lin12}.
We also omit the more regular nonlocal terms, such as $\uu,c_\om$, and the
nonsingular part $u_{x,R}=u_x-u_x(\om,\e)$, which are estimated by the methods
in Sections~\ref{sec:model2} and \ref{sec:finite_pertb}. Assume
\beq\label{eq:model3_1}
-d(\xx)+ M_0 |a(\xx)|\leq -c_1<0,\quad
-d(\xx)\leq -c_2<0,\quad
d(\xx), \, |a(\xx)|, \, [d]_{C^{1/2}}, \, [a]_{C^{1/2}}\leq c_3,
\eeq
for some constants $c_i>0$, where
$ [ \, \cdot\, ]_{C^{1/2}}$ is the H\"older semi-norm  
with $g(h) = |h|^{-1/2}$ \eqref{hol:g}  and 
\[
M_0 = C_1(\infty) +  (\tf12 C_1(\infty)+C_2(\infty) ) . 
\]
The constants $C_i(\cdot)$ are defined in
Lemmas~\ref{lem:holx_ux}, \ref{lem:holy_ux}. These lemmas imply the uniform
estimate
\beq\label{eq:model3_2}
[ u_x(\om,a,b) ]_{C^{1/2}} \leq M_0 [\om]_{C^{1/2}}
\eeq
for any $a,b$. The first condition in \eqref{eq:model3_1} is the analogue of
the stability condition \eqref{eq:model_hol_sing3} in the $C^{1/2}$ estimate
when $|\xx-\zz| \to 0$. It states that the local damping $d(\xx)$ dominates the estimate of nonlocal terms. \footnote{
In Sections \ref{sec:wg_hol_g} and \ref{sec:EE_hol}, by exploiting the anisotropic flow 
and the damping term from the H\"older weight $g_i$, we significantly weaken the 
first stability constraint  in \eqref{eq:model3_1}.
}

We show that the model is linearly stable if $\e$ is sufficiently small. We use
$C$ to denote absolute constants independent of $c_i$ and $\e$. From the formula
for the localized velocity in \eqref{eq:ux_loc_idea}, we have
\beq\label{eq:model3_linf}
|u_x(\om,\e)(\xx)| \leq C\e^{1/2}[\om]_{C^{1/2}},
\qquad
|a(\xx)u_x(\om,\e)(\xx)| \leq Cc_3\e^{1/2}[\om]_{C^{1/2}} .
\eeq
Thus the $L^\infty$ estimate is almost closed: the only coupling from the
$C^{1/2}$ norm to the $L^\infty$ norm carries the small factor $\e^{1/2}$.

We next derive the $C^{1/2}$ estimate. Denote
$\d f(\xx,\zz)=f(\xx)-f(\zz)$. Using \eqref{eq:model3_1}, we have
\[
\bal
-\f{\d(d\om)(\xx,\zz)}{|\xx-\zz|^{1/2}}
&=
-d(\xx)\f{\d\om(\xx,\zz)}{|\xx-\zz|^{1/2}}
-\f{d(\xx)-d(\zz)}{|\xx-\zz|^{1/2}}\om(\zz), \quad
\B|\f{d(\xx)-d(\zz)}{|\xx-\zz|^{1/2}}\om(\zz)\B|
 \leq c_3\|\om\|_{L^\infty}.
\eal
\]
For the nonlocal term, using \eqref{eq:model3_1}, \eqref{eq:model3_2}, and
\eqref{eq:model3_linf}, we get
\[
\bal
\B|\f{\d(a u_x(\om,\e))(\xx,\zz)}{|\xx-\zz|^{1/2}}\B|
&\leq
|a(\xx)|\B|\f{\d(u_x(\om,\e))(\xx,\zz)}{|\xx-\zz|^{1/2}}\B|
+\B|\f{a(\xx)-a(\zz)}{|\xx-\zz|^{1/2}}\B|\,|u_x(\om,\e)(\zz)|\\
&\leq
|a(\xx)|M_0 [\om]_{C^{1/2}}
+Cc_3\e^{1/2}[\om]_{C^{1/2}}.
\eal
\]
It follows that
\beq\label{eq:model3_4}
\bal
\pa_t\f{\d\om(\xx,\zz)}{|\xx-\zz|^{1/2}}
&=
-d(\xx)\f{\d\om(\xx,\zz)}{|\xx-\zz|^{1/2}}
+B(\xx,\zz),\\
|B(\xx,\zz)|
&\leq
\big(|a(\xx)| M_0+Cc_3\e^{1/2}\big)[\om]_{C^{1/2}}
+c_3\|\om\|_{L^\infty}.
\eal
\eeq

To apply Lemma~\ref{lem:PDE_stab}, we use the mixed energy
\beq\label{eq:model_trunc_EE}
E=\max\left(\|\om\|_{L^{\inf}},\tau^{-1}\|\om\|_{C^{1/2}}\right),
\eeq
and will choose a large $\tau$.
Estimates \eqref{eq:model3_linf} and \eqref{eq:model3_4} reduce the stability conditions \eqref{eq:PDE_diag} to
\beq\label{eq:model3_stab}
L^{\infty}:\quad
c_2-Cc_3\e^{1/2}\tau\geq \lam,
\qquad
C^{1/2}:\quad
d(\xx)-|a(\xx)|M_0-Cc_3\e^{1/2}-\tau^{-1}c_3\geq \lam,
\eeq
for some $\lam>0$. Since $d(\xx)-|a(\xx)|M_0\geq c_1$ by \eqref{eq:model3_1},
if $\e$ is sufficiently small so that
\[
Cc_3\e^{1/2}<\min( \f{c_1}{2}, \f{c_2}{2},\f{c_2c_1}{2c_3} ),
\]
then
\[
d(\xx)-|a(\xx)|M_0-Cc_3\e^{1/2}\geq \tf{1}{2} c_1.
\]
We can then choose
\beq\label{eq:model3_5}
\tau=\sqrt{\f{2c_2}{C\e^{1/2}c_1}}\sim \e^{-1/4},
\eeq
which is large, so that $\tau^{-1}$ is small and the stability conditions
\eqref{eq:model3_stab} hold for some $\lam>0$.

\vs{0.05in}
\paragraph{\bf{Insights}}

The $L^\infty$ estimate is almost closed because the $L^{\infty}$ estimate of the localized velocity has a small factor $\e^{1/2}$; see \eqref{eq:model3_linf}. Thus, we may regard $\|\om\|_{L^\infty}$ as the more stable component. The $C^{1/2}$ estimate is more delicate.
The damping assumption in \eqref{eq:model3_1}
controls the top-order contribution. The remaining coupling from the $L^\infty$ norm
to $C^{1/2}$ seminorm is made perturbative by choosing the mixed energy with a small coupling weight
$\tau^{-1}$ for $\|\om\|_{C^{1/2}}$. In the notation of
\eqref{eq:PDE_diag},  $\tau^{-1}$ and $c_3\tau^{-1}$ here correspond to
the weights $\mu_i$ and the off-diagonal estimate
$\sum_{j\neq i}a_{ij}\mu_i\mu_j^{-1}$. Our energy estimates for
\eqref{eq:bous_lin11} follow the same principle; see Step 5 in
Section~\ref{sec:EE_intro}.

\vs{0.05in}
\paragraph{\bf Obtain small $\e$ by approximation}

In the full problem, smallness is not obtained by literally localizing $u_x$ at an
arbitrarily small scale. Instead, we construct a finite-rank approximation $\hat u_x$
to the nonsingular part of $u_x$. The remaining error $u_x-\hat u_x$ plays the role
of the localized velocity $u_x(\om,\e)(\xx)$ in the model problem. Increasing the rank
$N$ improves the approximation and therefore produces the analogue of a smaller
localization parameter $\e$.

\subsection{Other notations}

We introduce the directional H\"older seminorms  $C_x^{\al}, C_y^{\al}$
\beq\label{hol:semi}
[ \om ]_{ C_x^{\al}(D)} \teq \sup_{ x, z \in D, \, x_2 = z_2, x_1 \neq z_1} \f{ | \om(x) - \om(z)| }{ |x_1 - z_1|^{\al}} , \  
[ \om ]_{ C_y^{\al}(D)} \teq \sup_{ x, z \in D,  \, x_1 = z_1, x_2 \neq z_2} \f{ | \om(x) - \om(z)| }{ |x_2 - z_2|^{\al}} .
\eeq

Given a weight $g(h): \R^2\bsh \{ 0 \} \to \R_+ $ that is $-\al$-homogeneous , i.e. $g(\lam h_1, \lam h_2) = \lam^{-\al} g(h)$, e.g., $g(h) = |h|^{-\al}$, we define the weighted H\"older seminorm 
\beq\label{hol:g}
[ \om ]_{C^{\al}_g(D)} 
= \sup\nolimits_{x, z \in D, x \neq z } | (\om(x) - \om(z) ) g(x-z) |.
\eeq
In this paper, we mostly use $D = \R^2_{++}$. In this case, we drop $D$ to simplify notation. 
If $g(h) \equiv |h|^{-\al}$, we simplify $[ \om ]_{C^{\al}_g(D)}$ as $[ \om ]_{C^{\al}(D)}$.

We introduce  $\la a \ra  \teq (|a|^2 + 1)^{1/2}$ and  define the inner product in the \emph{first quadrant} $\R^2_{++}$
\beq\label{inner}
\la f, g \ra =  \int_{\R^2_{++}}  f( \xx) g(\xx ) \, d \xx .
\eeq

We use $\xx$, $(x,y)$, and $(x_1,x_2)$ interchangeably to denote a vector in $\R^2$.

\section{Sharp H\"older estimates via optimal transport}\label{sec:sharp}

In this section, we establish the sharp directional $C^{1/2}$ estimate for $\na \uu$ using the symmetry properties of the kernels and some techniques from optimal transport. Note that novel functional inequalities for related Biot-Savart laws play an important role in
\cite{kiselev2013small,elgindi2019finite}. The sharp H\"older estimates below
are likewise an important ingredient in our stability analysis.

Obtaining a sharp estimate of $\na\uu(\xx)-\na\uu(\zz)$ for arbitrary
$\xx,\zz\in\R^2_{++}$ by a direct argument is difficult, since the
corresponding kernel has a complicated sign structure and loses some useful symmetries of $\na\na^\perp(-\D)^{-1}$. Instead, we estimate the
directional seminorms $C_x^{1/2}$ and $C_y^{1/2}$ in \eqref{hol:semi}. This
choice is adapted to the anisotropy of the linearized equations
\eqref{eq:bous_lin11}, \eqref{eq:bous_lin12}: the damping in the $y$ direction
is much stronger, so the sharp estimate in the $x$ direction is the most
important; see Section~\ref{sec:aniso_est}. The directional seminorms also
reduce the dimension of the pair $(\xx,\zz)$ from four to three, since one
assumes either $x_1=z_1$ or $x_2=z_2$. In these cases, the kernel in
$\na\uu(\xx)-\na\uu(\zz)$ has better symmetry and simpler sign properties,
which allow us to reduce the two-dimensional estimates to estimating many 1D integrals.  Once the directional estimates
$[\na\uu]_{C_x^{1/2}}$ and $[\na\uu]_{C_y^{1/2}}$ are obtained, the full
$C^{1/2}$ estimate follows by the triangle inequality.

Below, we first introduce the kernels and
localized nonlocal operators, and state the sharp directional estimates in
Section~\ref{sec:holder_vel}. We introduce the connection between 
optimal transport and sharp H\"older estimates in Section \ref{sec:idea_opt}. The rest of this section and Appendix~\ref{app:sharp} are devoted to the proof 
of four directional estimates, which we summarize in Table \ref{tab:hol}. The estimate of
\([u_x]_{C^{1/2}_x}\) is the most important one and contains the main idea of
the optimal-transport argument. Other estimates follow the same strategy,
with additional technical modifications caused by the boundary and by the local
terms in \(u_y\) and \(v_x\).
\begin{table}[t]
\centering
\begin{tabular}{c|c|c|c|c}
\hline
Quantity & Direction &  lemma & seminorm(s) in the upper bound & Proof \\
\hline
$u_x$ & $C^{1/2}_x$ &  Lemma \ref{lem:holx_ux} &  $[\omega]_{C^{1/2}_x}$ & Section~\ref{sec:holx_ux} \\
$u_x$ & $C^{1/2}_y$ & Lemma \ref{lem:holy_ux}  & $[\omega]_{C^{1/2}_x}, \  [\omega]_{C^{1/2}_y}$ & Section~\ref{sec:holy_ux} \\
$u_y,\ v_x$ & $C^{1/2}_x $ & Lemma \ref{lem:holx_uy}  &  $[\omega]_{C^{1/2}_x},\ [\omega]_{C^{1/2}_y}$ & Appendix~\ref{app:sharp} \\
$u_y,\ v_x$ & $C^{1/2}_y$ & Lemma \ref{lem:holy_uy}  & $[\omega]_{C^{1/2}_x},\ [\omega]_{C^{1/2}_y}$ & Appendix~\ref{app:sharp} \\
\hline
\end{tabular}
\caption{Outline of sharp H\"older estimates.}
\label{tab:hol}
\end{table}

The proofs have two levels. Conceptually, they use the cancellation of sign-changing kernels and match the positive and negative parts via a 1D transportation map. Technically, the main work is to verify the kernel sign structure,
decompose the domain according to the sign structure, and compare the transportation costs for localized and non-localized kernels.

\vs{0.05in}

\paragraph{\bf Analytic estimates and certified constants} The estimates in this section 
and Appendix \ref{app:sharp} are \emph{analytic}. The constants are given by \emph{explicit integrals} obtained from the optimal transport arguments.  In the supplementary material II in Part II \cite[Section {\secholconst}]{ChenHou2023bSupp}, we estimate these explicit integrals using some integral formulas 
and numerical quadrature with computer assistance, and obtain rigorous upper bounds.

\vs{0.05in}

\paragraph{\bf Kernels}
The kernels associated with $\na \uu$ are given by
\beq\label{eq:kernel_du}
K_1(y) \teq \f{y_1 y_2}{|y|^4}, \quad K_2(y) \teq  \f{1}{2}\f{ y_1^2 - y_2^2}{ |y|^4} , 
\quad G(y) = - \f{1}{2} \log |y|,
\eeq
where $\f{1}{\pi} G$ is the Green function of $-\D$ in $\R^2$. Note that $\pa_1 \pa_2 G = K_1,  \pa_1^2 G  = K_2 $.

Denote by $ K_{i,s}$ the symmetrized kernel 
\beq\label{eq:kernel_du_sym}
K_{i,s}(x, y) \teq K_i( x- y) - K_i (x_1 + y_1, x_2 - y_2)  - K_i (x_1 - y_1, x_2 + y_2) + K_i(x+y).
\eeq

Consider the odd extension of $\om$ in $y$ from $\R_2^+$ to $\R_2$
\beq\label{eq:ext_w_odd}
W(y) = \om(y)  \ \mathrm{ for } \ y_2 \geq 0, \quad W(y) = - \om(y_1, -y_2) \quad  \mathrm{ for } \ y_2 < 0. 
\eeq
$W$ is odd in both $y_1$ and $y_2$ variables. Clearly, $u_x$ can be written as 
\[
u_x(x) =  - \f{1}{\pi} \int_{\R^2}  K_1(x- y) W(y) dy = - \f{1}{\pi} \int_{\R_2^{++}} \om(y) K_{1,s}(x,y) dy.
\]

We localize the singular integral because in the energy estimates only the
near-diagonal part of the kernel contributes to the top-order H\"older seminorm.
The far part has a smoother kernel and will be treated later as a more regular
nonlocal term.  The localization also gives sharper constants when the
localization scale is small compared with \(|x-z|\).

For any $a ,b_1, b_2>0, x_2 > 0 $, we consider the localized velocity 
\beq\label{eq:ux_local}
\bal
Q_{a, b_1, b_2}(x) & \teq  [x_1 -a, x_1 +a] \times [x_2 - b_1, x_2 + b_2], \\
u_{x}(x, a, b_1, b_2)  & \teq  - \f{1}{2\pi} P.V. \int_{ y \in Q_{a,b_1, b_2}(x) } \f{ 2 (x_1 - y_1) (x_2 - y_2)}{ | x-y|^4 } W(y) dy , \\ 
u_y(x, a )  &\teq \f{1}{2 \pi} P.V. \int_{  y \in Q_{a, a}(x)   }
\f{  (x_1 - y_1)^2 - (x_2 - y_2)^2 }{ |x-y|^4} W(y) dy + \f{\om(x)}{2}, \\  
v_x(x, a )  &\teq \f{1}{2 \pi} P.V. \int_{  y \in Q_{a, a}(x)  }
\f{  (x_1 - y_1)^2 - (x_2 - y_2)^2 }{ |x-y|^4} W(y) dy - \f{\om(x)}{2} .\\ 
\eal
\eeq
If $b_1 = b_2 = b$, we simplify  $u_x(x, a, b_1, b_2) $ as $u_x(x, a, b)$; if $b_1 = b_2 = a$, we further simplify $u_x(x, a, b_1, b_2) $ as $u_x(x, a)$.

\subsection{H\"older estimates of the velocity}\label{sec:holder_vel}

We have the following estimates for $\na \uu$. We discuss the ideas in Section \ref{sec:idea_opt} and the proof in Sections \ref{sec:holx_ux}, \ref{sec:holy_ux}, and Appendix \ref{app:sharp}. We localize the velocity in \eqref{eq:ux_local} to obtain improvement of the constant $C_1( \f{b}{|x-z|} )$ when $|x-z|$ is large.

\vspace{0.1in}
\begin{lem}[Estimate of ${[} u_x {]}_{C_x^{1/2}}$]\label{lem:holx_ux}
For any $ b_1, b_2 >0, a \geq \f{1}{2} |x_1 - z_1| , x=  (x_1, x_2),  z = (z_1, x_2) \in \R_2^+ $, and $D $ covering $( Q_{a, b_1, b_2}(x) \cup Q_{a, b_1,b_2}(z) ) \cap \R_2^+$ \eqref{eq:ux_local}, we have 
\[
 \f{|  u_x( x,  a, b_1  , b_2  ) -  u_x( z, a, b_1 , b_2 ) | }{ |x-z|^{1/2}} 
 \leq C_1 \lt( \f{b}{|x-z|} \rt) [ \om ]_{C_x^{1/2}(D)} ,
  \]
where $b = \max(b_1, b_2)$ and $C_1(b)$ is an increasing function given by 
\[
\bal
C_1(b) & \teq \f{4}{ \pi} \B| \int_0^b d s_2 \int_{ f(s_2)}^{\infty} | T(s_1, s_2) - s_1|^{1/2} 
\D(s_1, s_2) d s_1 \B| ,\\
\D(s_1, s_2) &= \f{ (s_1 + 1/2) s_2 }{   ( (s_1+1/2)^2 + s_2^2)^2 } 
- \f{ (s_1 - 1/2) s_2 }{   ( (s_1 - 1/2)^2 + s_2^2)^2 } .\\
\eal
\]
Here, $f(s_2)$ is the unique solution in $[0, \infty)$ satisfying $\D( f(s_2), s_2) = 0$ and $T(s_1, s_2)$ is the unique solution in $[0, f(s_2))$ that solves 
\[
\int_{ T(s_1, s_2)}^{s_1} \D( t , s_2) d t = 0,
\]
for $s_1 > f(s_2)$. In particular, $T(s_1, s_2)$ can be obtained explicitly by solving a cubic equation and $C_1(b) \leq 2.55$ for any $b> 0$.

\end{lem}

If $a < \f{1}{2} |x_1 - z_1|$, the singular region is small. We can simply apply the triangle inequality to estimate each term. We localize the seminorm to region $D$ since we only use the seminorm to control $w(x) - w(z)$ for $x, z \in D$. The same reasoning applies to the following lemmas.

\begin{remark}[\bf Generalization]
The above Lemma can be further generalized to the localized velocity $ \om \ast K_1(s) \one_{ -a_1 \leq s \leq a_2, -b_1 \leq s_2 \leq b_2 } $, i.e., we do not need $a_1 = a_2$ in \eqref{eq:ux_local}. The proof follows from the same argument. Yet, we  only use the special case $a_1 = a_2$ in our later estimates.
\end{remark}

The upper bounds in the following Lemmas involve $[\om]_{C_x^{1/2}}$ and $[\om]_{C_y^{1/2}}$. We will further bound it using the energy norm.

\begin{lem}[Estimate of ${[} u_x{]}_{C_y^{1/2}}$]\label{lem:holy_ux}
For any $a , b \geq |x-z| , x=  (x_1, x_2),  z = (x_1, z_2) \in \R_2^+ $, and $D$ covering $( Q_{a, b}(x) \cup Q_{a, b}(z) ) \cap \R_2^+$ \eqref{eq:ux_local}, we have 
\[
 \f{|  u_x( x,  a ,b  ) -  u_x( z, a, b ) | }{ |x-z|^{1/2}} 
 \leq \f{1}{2} C_1( \f{ a}{|x-z| } ) [ \om]_{C_y^{1/2}(D) }
 +  C_2( \f{a}{ |x-z|}) [\om ]_{C_x^{1/2}(D)}   ,
\]
where $C_1(a)$ is defined in the previous Lemma and $C_2(a)$ is given by 
\[
C_2(a) =  \f{ \sqrt{2} }{\pi} \int_0^a \int_0^{\infty} y_1^{1/2} \B| \f{y_1 (1/2-y_2)}{ (y_1^2 + (1/2-y_2)^2)^2 }   +\f{y_1 (1 /2 + y_2) }{ (y_1^2 + (1 / 2+ y_2)^2)^2 } \B|  dy.
\] 
In particular, $C_2(a) \leq \f{4.26}{\pi}$. 
\end{lem}

Next we estimate the other kernel. We remark that for $u_y(x, a)$ and $v_x(y, a)$, the estimates are different due to the local term related to $\om$ \eqref{eq:ux_local}.

\begin{lem}[$C_x^{1/2}$ estimate of $ u_y, v_x$]\label{lem:holx_uy}
For any $a \geq 2 |x-z| , x=  (x_1, x_2),  z = (z_1, x_2) \in \R_2^+ $, $D$ covering $( Q_{a}(x) \cup Q_{a}(z) ) \cap \R_2^+$ \eqref{eq:ux_local}, and $\tau > 0$, we have 
\[
\bal
 \f{|  v_x( x,  a   ) -  v_x( z, a  ) | }{ |x-z|^{1/2}} 
 &\leq \f{1}{\pi} C_1(\tau) \max( [ \om ]_{C_x^{1/2}(D) }  , \tau^{-1} [ \om ]_{C_y^{1/2}(D)}  ), \\
  \f{|  u_y( x,  a   ) -  u_y( z, a  ) | }{ |x-z|^{1/2}} 
 &\leq \f{1}{\pi} C_2(\tau) \max( [ \om ]_{C_x^{1/2}(D) }  , \tau^{-1} [ \om ]_{C_y^{1/2}(D)}  ),
 \eal
\]
for some constant $C_1(\tau), C_2(\tau) > 0$ with $ \f{1}{\pi} C_1(0.582) \leq 2.53$ and $ \f{1}{\pi} C_2( 0.582) \leq  1.55$.
\end{lem}

In general, the localized $u_y$ is not in $C_y^{1/2}$ due to the presence of the boundary and the discontinuity of $W$ cross $y = 0$. Thus, we consider the estimate without localizing the kernel.

\begin{lem}[$C_y^{1/2}$ estimate of $ u_y, v_x$]\label{lem:holy_uy}
For $x = (x_1, x_2), z = (x_1, z_2) \in \R_2^+$, and any $\tau > 0$, we have 
\[
\bal
 \f{|  v_x( x,  \infty  ) -  v_x( z, \infty  ) | }{ |x-z|^{1/2}} 
 &\leq  \f{C_3(\tau)}{\pi} \max( [ \om ]_{C_x^{1/2} } , \tau^{-1} [ \om ]_{C_y^{1/2} }  ) , \\
 \f{|  u_y( x,  \infty  ) -  u_y( z, \infty  ) | }{ |x-z|^{1/2}} 
 &\leq  \f{ C_4(\tau)}{\pi}  \max( [ \om ]_{C_x^{1/2} }
 ,  \tau^{-1} [ \om ]_{C_y^{1/2} } ),
 \eal
\]
for some constant $C_3(\tau), C_4(\tau) >0$. We have $ \f{1}{\pi} C_3( 0.582) \leq 2.60, \f{1}{\pi} C_4(0.582) \leq 2.61 $.
\end{lem}

In the proof of the above two lemmas, we provide computable upper bounds for $C_i(\tau),i \leq 4$.
The parameter \(\tau\) is chosen to match \emph{precisely} the anisotropy built into the  H\"older seminorm $ [ \om \psi_1 ]_{C_{g_1}^{1/2}} $ with  weight $g_1$ used in the energy estimates. See \eqref{eq:wg_hol_g}. 
In practice, we choose $\tau = g_1( 0, 1) / g_1(1, 0) \approx 0.582$.

\subsection{Connection to optimal transport and ideas of the proof}\label{sec:idea_opt}

A key observation is that the H\"older estimate is related to an optimal transport problem. We illustrate the ideas by proving a sharp H\"older estimate of the Hilbert transform. 
The Hilbert transform can be seen as an approximation of $u_x(\om)$, which is exact if $\om(x, y)$ is constant in $y$ \cite{luo2013potentially-2,choi2014on}.

We estimate $ \f{1} {|x-z|^{1/2}}  | H f(x) - H f(z) |$ by $[f ]_{C^{1/2}}$. Due to translation and scaling symmetry, we can assume $x = 1, z= -1$ without loss of generality. Then we need to estimate 
\beq\label{eq:hol_hil}
 S = H f(1) - H f(-1) = \f{1}{\pi} \int_{\R} ( \f{1}{ 1- y} + \f{1}{1 + y} ) f(y) dy 
= \f{2}{\pi} \int_{\R} \f{1}{1 - y^2} f(y) dy.
\eeq
The kernel $k(y) = \f{2}{1- y^2}$ is positive on $(-1, 1)$ and negative for $|y| > 1$, and satisfies $ \int k(y) dy = 0$. 

Denote $k^{\pm}(y) = \max( \pm k(y), 0)$. An estimate of $S$ using $[f]_{C^{1/2}}$ is equivalent to estimating the transportation cost of moving the positive region of $k(y)$ with measure $k^+(y) dy $ to its negative region with measure $k^-(y) dy$ with distance function $c(x, y) = |x - y|^{1/2} $.

For example, if $k(y) = \d_{1}(y) + \d_2(y) - \d_3(y) -\d_4(y)$, where $\d_a(x) $ is the Dirac function centered at $a$, then we get 
\[
|\int k(y) f(y) dy |= | f(1) + f(2) - f(3) - f(4)|
\leq | f(2) - f(3)|  + |f(1) - f(4)| \leq (\sqrt{1} + \sqrt{3} )  || f||_{C^{1/2}}. 
\]

The above estimate can be interpreted as moving the mass from $2$ to $3$ and $1$ to $4$ with cost function $|x- y|^{1/2}  || f||_{C^{1/2}}$. Using the language of optimal transport, 
to obtain sharp estimate of $S$ \eqref{eq:hol_hil}, we are seeking a measurable map $T$ such that $T \# k^+ dy  = k^- dy$, where $(T \# \mu)(A) = \mu( T^{-1}(A))$ for a measurable set $A$, and the following cost 
\[
 C(T) = \int_{ k(y) \geq 0} | T( y ) - y|^{1/2}  k(y) dy  || f||_{C^{1/2}}
\]
is as small as possible. 
Based on the above discussion, we have the following transportation lemma, which is used repeatedly in the H\"older estimate.
\begin{lem}[Transportation Lemma]\label{lem:trans}
Suppose that there exists $c \in (a, b)$ such that $k < 0$ on $(a,c)$, $k > 0 $ on $(c, b)$, 
$k \in C( (a,c) \cup (c, b))$, $k |x-c|^{\al} \in L^1_{loc}$ with $P.V. \int_a^b k(x) dx = 0$.
For $\al \in (0, 1), g \in C^{\al}(a, b)$, we have
\[
 \B| P.V. \int_a^b k(x) g (x) dx \B| \leq \int_c^b |k(x)| |x - T(x)|^{\al} dx [g]_{C_x^{\al}}
 =  \int_a^c |k(x)| |x - T(x)|^{\al} dx [g]_{C^{\al}},
\]
where integrals across $c$ are understood in the principal-value sense, and $T(x)$ solves  $P.V. \int_{x}^{T(x)} k(s) ds = 0$ with $ (T(x) - c) ( x-c) \leq 0$.

\end{lem}

We use the fact that $|h|^{\al}$ is concave for $\al \in (0, 1)$ to design the map $T$. In our later estimates of $\na \uu$, we use the above Lemma with $\al = \f{1}{2}$

\begin{proof}

Since $P.V.\int_a^b k=0$, we may replace $g$ by $g-g(c)$ and thus assume that $g(c) = 0$. The resulting $g$-integrals in Lemma \ref{lem:trans} and in the following proof are absolutely convergent.  

First we study how to construct the map $T$.  
 For $x_1 < x_2 < x_3< x_4$ and $\al \in (0, 1)$, we have
\[
\bal
 &|\int (\d_{x_1} + \d_{x_2} - \d_{x_3}  - \d_{x_4} ) g(x) dx |
 = |g(x_1) + g(x_2) - g(x_3) - g(x_4)| \\
 \leq & \min( |x_1 - x_3|^{\al} + |x_2 - x_4|^{\al}, |x_1 -x_4|^{\al} + |x_2 - x_3|^{\al}  ) [g]_{C_x^{\al}} \\
 =& ( |x_1 -x_4|^{\al} + |x_2 - x_3|^{\al} )  [g]_{C^{\al}}.
 \eal
\]
The above estimate indicates that to find an optimal map $T$ moving $ (\d_{x_1} + \d_{x_2}) dx $ to 
$ (\d_{x_3} + \d_{x_4}) dx $ with cost $|x-y|^{\al}$, we should choose $T(x_1) = x_4, T(x_2) = x_3$, which implies that $T(x)$ is decreasing in $x$. Due to conservation of mass and the sign properties of $k$, a natural construction of $T : (a, c) \to (c, b)$ is given by 
\beq\label{eq:trans_lem1}
\int_{x}^{T(x)} k(s) d s = 0, 
\eeq
for $x < c$, which implies $ T^{\prime}(x) k( T(x)) = k(x)$ for smooth $k$.
The idea of the above map is to move the mass in the positive region to its closest negative region that has not been occupied due to the monotonicity of $T$. Using a change of variable $y = T(x) : (a, c) \to (c, b]$, we get
\[
 \int_{c}^b k(x) g(x) d x = - \int_a^c k( T(x)) g(T(x)) T^{\prime}(x) dx 
 = - \int_a^c k(x)  g(T(x)) dx.
\]
It follows 
\[
|\int_a^b  k g d x  | = | \int_a^c k(x) (g(T(x)) - g(x)) dx  |
\leq \int_a^c |k(x)| |T(x) - x|^{\al} dx [g]_{C_x^{\al}}.
\]
Similarly, we can define $T : (c, b) \to (a, c)$ by solving $\int_{T(x)}^x k(s) ds = 0$, which is also equivalent to \eqref{eq:trans_lem1}. The first inequality in Lemma \ref{lem:trans} follows from the same argument.
\end{proof}

In the estimates below, the main work is therefore to identify the sign structure of the kernel
and to construct an efficient one-dimensional transport map.

\subsubsection{$C^{1/2}$ estimate of the Hilbert transform}

We use the Hilbert transform as an example to illustrate Lemma \ref{lem:trans}.
We apply Lemma \ref{lem:trans} with $k(y) = \f{1}{1-y^2}$, up to the harmless
prefactor \(2/\pi\), and $g(y) = f(y)$ to estimate \eqref{eq:hol_hil}. For any $x > 0$, we construct the transportation map $T(y)$ by solving 
\[
 0 =\int_x^{T(x)} k(y) dy = \int_x^{T(x)} \f{1}{1-y^2} dy = 0,
\]
which implies 
\[
\f{ x+1}{1-x} = \f{ T+1}{ T-1}, \quad T(x) = \f{1}{x},
\]
where we have used $T(x) > 1$ if $x<1$ and $T(x) < 1$ if $x>1$ due to the sign of $\f{1}{1-y^2}$.  This map also applies to $y < 0$. Applying this map to \eqref{eq:hol_hil}, we yield 
\beq\label{eq:hol_hil2}
\bal
|S| & = |\f{2}{\pi} \int_{y > 1} k(y) ( f( y) - f( T(y)) ) dy 
+ \f{2}{\pi} \int_{y < -1} k(y) ( f(y) - f( T(y)) ) dy | \\
&\leq \f{4}{\pi} \int_{y>1} |k(y)| | y-  T(y)|^{1/2} dy  [f]_{C^{1/2}} 
 = \f{4}{\pi} \int_{y > 1} \f{1}{ |y^2 - 1|}  \f{ |y^2 - 1|^{1/2}}{ |y|^{1/2}} dy [ f ]_{C^{1/2}} \\
 & = [ f]_{C^{1/2}} \f{4}{\pi} \int_{y>1} \f{1}{ |y^2 - 1|^{1/2} y^{1/2}} dy = C [ f ]_{C^{1/2}}, \quad \f{ C}{\sqrt{2}} \approx 2.37 .
\eal
\eeq

Since $x, z$ are arbitrary, we yield $[H f]_{C^{1/2}} \leq \f{C}{\sqrt{2}} [ f ]_{C^{1/2}}$. The equality achieves if $|f(y) - f( \f{1}{y})| = |y-\f{1}{y}|^{1/2} [f]_{C^{1/2}}$ for all $y>0$ and $y<0$. Since the Hilbert transform satisfies $H (H f) = - f$, the sharp constant in $[ H f]_{C^{1/2}} \leq C_* [ f ]_{C^{1/2}}$ satisfies $C_* \geq 1$.

Below, we apply Lemma \ref{lem:trans} to 
prove  Lemmas \ref{lem:holx_ux}, \ref{lem:holy_ux} for the H\"older estimate of $u_x$, and refer to Table \ref{tab:hol} for the organization of estimates for other cases.  
We need two steps. Firstly, we identify the sign of the kernel $K$ in the integral of $ \na \uu( \xx ) - \na \uu( \zz )$. Next, we fix a variable in the 2D integral in one direction, e.g. fix $x=a$, and then apply Lemma \ref{lem:trans} to estimate the 1D integral in the other direction, e.g., on the line $ \{ (a, y): y \in \R\}$. One may generalize Lemma \ref{lem:trans} to 2D and construct the 2D optimal transport map directly. Yet, the domain where the kernel $K$ is positive or negative is complicated. To avoid this difficulty, we build the 2D transport map using the 1D Lemma \ref{lem:trans} repeatedly. The odd symmetry of the kernel $K_1(s)$ in $s_1$ enables us to apply this approach to obtain sharp estimate of $u_x$ effectively. See Remark \ref{rem:holx_ux1}.

\subsection{Estimate of $[u_x]_{C_x^{1/2}}$}\label{sec:holx_ux}

In this subsection, we prove Lemma~\ref{lem:holx_ux} for $u_x$. The proof has three steps. First, we identify the sign structure of the kernel
$K(s)$ \eqref{eq:holx_ux_ker} in $u_x(\xx)-u_x(\zz)$ when $x_2=z_2$. Second, for each fixed vertical
variable $s_2$, we apply the 1D transport lemma in the horizontal
variable $s_1$. Third, for the localized operator, we compare the resulting
transportation cost with the cost in the non-localized case.

We begin with the non-localized case, corresponding to $a=b_1 = b_2 = \infty$ in Lemma 
\ref{lem:holx_ux}.  In the $C_x^{1/2}$ estimate of $u_x$, we have $x_2 = z_2$. 
Using the scaling symmetry and translation invariance, we only need to estimate the following 
\beq\label{eq:holx_ux_int1}
u_x(\f{1}{2}, x_2) - u_x(-\f{1}{2}, x_2) 
= - \f{1}{\pi} P.V. \int_{ \R^2} K(s) W( s_1,  x_2 - s_2) ds 
\eeq
for any $x_2$, where $W$ is an odd extension of $\om$ from $\R_2^+$ to $\R^2$, and $K(s)$ is given by 
\beq\label{eq:holx_ux_ker}
 K(s) = K_1( s_1 + \f{1}{2}, s_2) - K_1(s_1 - \f{1}{2}, s_2) = \f{ (s_1 + \f{1}{2}) s_2}{  ( (s_1 + \f{1}{2})^2 + s_2^2 )^2}
-  \f{ (s_1 - \f{1}{2}) s_2}{  ( (s_1 - \f{1}{2})^2 + s_2^2 )^2} .
\eeq

Since $K(s)$ is odd in $s_2$, we consider $s_2 \geq 0$ without loss of generality. We will only use the H\"older seminorm of $W$, $[ W]_{C_x^{1/2}}$, to estimate the above quantity. Note that $ [W]_{C_x^{1/2}(\R^2)} = [ \om]_{C_x^{1/2}(\R^2_+)}$. Without loss of generality, we can assume that $x_2 = 0$.

We first identify the sign structure of the kernel. A direct calculation yields 
\[
K(s) = \f{ s_2 \D_1(s_1, s_2)}{  ( (s_1 + \f{1}{2})^2 + s_2^2 )^2 
( (s_1 - \f{1}{2})^2 + s_2^2 )^2 }, 
\quad \D_1 = s_2^4 - 2 s_1^2 s_2^2 - 3 s_1^4 + \f{1}{2} s_1^2 + \f{1}{2} s_2^2 + \f{1}{16}.
\]

For a fixed $s_2$,  $\D_1(s_1, s_2) = 0$ implies 
\beq\label{eq:ux_thres}
s_1 = f(s_2)  = \B(  \f{ \f{1}{2} - 2 s_2^2 + \sqrt{ 16 s_2^4 + 4 s_2^2 + 1 }  }{ 6 }  \B)^{1/2}.
\eeq

\begin{figure}[h]
   \centering
      \includegraphics[width = 0.4\textwidth  ]{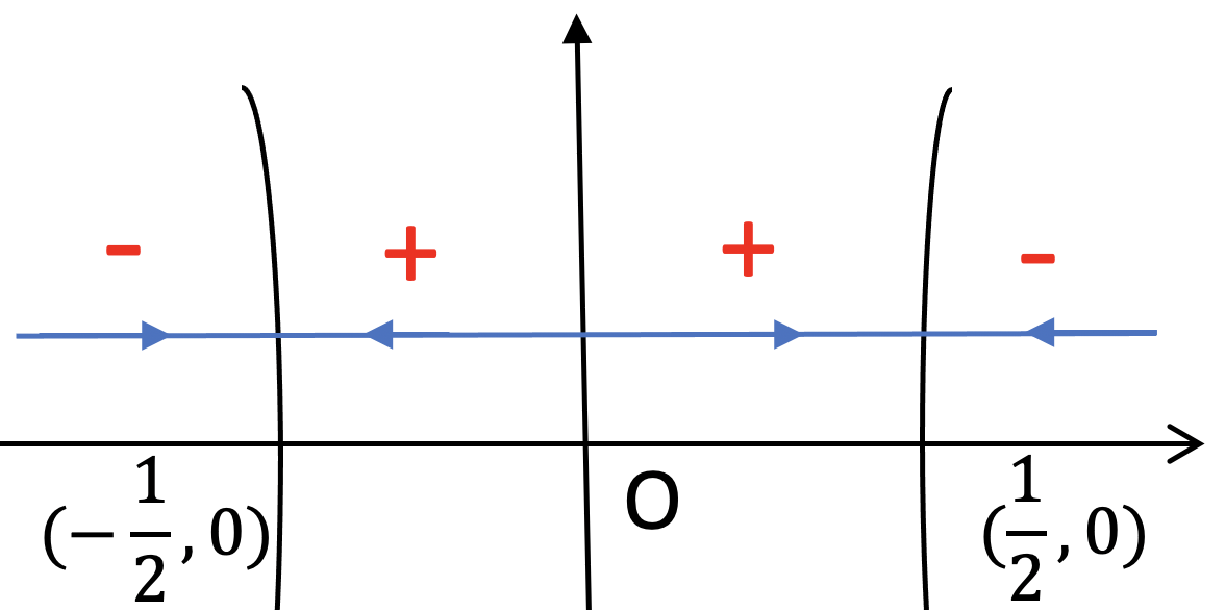}
      \caption{
      The black curve gives a schematic illustration of the sign-changing curve $s_1=f(s_2)$ in \eqref{eq:ux_thres},  and $\pm$ indicates $\sgn(K(s))$ in different regions. The blue arrows indicate the direction of 1D transportation plan. 
      }
                  \label{fig:ux_OT}
 \end{figure}

Moreover, for $s_1 , s_2 \geq 0$, it is easy to see that $\D_1(s_1, s_2) \geq 0$ if and only if 
\beq\label{eq:holx_ux_sign1}
K(s_1, s_2) \geq 0 \mathrm{ \ for \ } s_1  \in [0, f(s_2)],
\quad 
K(s_1, s_2) \leq 0 \mathrm{ \ for \ } s_1 \geq f(s_2).
\eeq
See Figure \ref{fig:ux_OT} for an illustration of $\sgn(K(s))$.

Note that the sign changes if $s_2 \leq 0$ since $K$ is odd in $s_2$. Since $K_1$ is odd in $s_1$, we get
\[
\bal
\int_0^{\inf} K(s_1, s_2) d s_1 &= \int_0^{\inf} K_1( s_1 + \f{1}{2}, s_2) -
K_1( s_1 - \f{1}{2}, s_2 )  d s_1 \\
&= \int_{1/2}^{\inf } K_1(s_1, s_2)  d s_1- \int_{-1/2}^{\inf} K_1(s_1, s_2) d s_1
=  - \int_{-1/2}^{1/2} K_1(s_1, s_2) d s_1 = 0.
\eal
\]

To estimate the integral in \eqref{eq:holx_ux_int1},  we first fix $s_2$ and then apply Lemma \ref{lem:trans} to estimate 
\beq\label{eq:holx_ux_int2}
 I(s_2) = \int_{ \R} K(s_1, s_2) W(s_1, x_2 - s_2) ds_1
 = (\int_{\R_- } + \int_{R_+}) K(s_1, s_2) W(s_1, x_2 - s_2) ds_1 \teq I_-(s_2) + I_+(s_2).
\eeq

We do so for the following reason. Near the singularity, from Taylor expansion of \eqref{eq:ux_thres}: $s_1 = \f{1}{2} + O(s_2^4)$, the curve $\G = \{ s: s_1 = f(s_2)\}$ is close to a straight line in the vertical direction. See Figure \ref{fig:ux_OT} for an illustration. Similar to the idea below \eqref{eq:trans_lem1}, an effective plan in 2D is to move the mass in the positive region to its closest possible negative region that has not been occupied. 
Thus, we expect that an effective 2D transport plan $(x, y) \to T(x, y)$ is orthogonal to the curve $\G$ and thus almost parallel to the $x$ direction.

\begin{remark}[\bf Shape of the curve]\label{rem:holx_ux1}
The fact that near the singularity $s = (\pm \f{1}{2}, 0)$, the curve $\G$ is almost vertical is due to the odd symmetry of $K_1(s_1, s_2)$ in $s_1$. In fact, from \eqref{eq:holx_ux_ker}, 
for $s$ close to $(\f{1}{2}, 0)$, we have $K(s) \approx - K_1(s_1 - \f{1}{2}, s_2)$, whose sign is determined by $\sgn(s_1 -\f{1}{2})$.
\end{remark}

Since $K(s_1, s_2)$ is even in $s_1$, we can estimate $I_+(s_2), I_-(s_2)$ in the same way.  To apply Lemma \ref{lem:trans}, we first construct $T(\cdot, s_2)$ on $[0, \infty)$ by solving 
\beq\label{eq:cubic_T0}
\int_{T(s_1, s_2)}^{s_1} K( t, s_2) d t = 0.
\eeq
We will show later that this equation has a unique solution of $T$ on $[0 ,\infty]$ for $s_1 > 0$.  Then applying Lemma \ref{lem:trans} to $I_+(s_2)$ and using $[ W]_{C_x^{1/2}} = [\om]_{C_x^{1/2}}$, we get 
\[
|I_+(s_2)| \leq  [\om]_{C_x^{1/2}}  \B|
 \int_{ f(s_2)}^{\inf} K(s_1, s_2) | T(s_1, s_2) - s_1|^{1/2} d s_1 \B| \teq M(s_2).
\]
See the blue arrows in Figure \ref{fig:ux_OT} for an illustration of this transportation plan.
Since $K(s_1, s_2)$ is even in $s_1$, the estimate of $I_-(s_2)$ in \eqref{eq:holx_ux_int2} is the same: $|I_-(s_2)| \leq M(s_2)$. Since $K(s)$ is odd in $s_2$, from \eqref{eq:cubic_T0}, we get $T(s_1, -s_2) = T(s_1, s_2)$. Therefore, the estimate of $I( s_2)$ in \eqref{eq:holx_ux_int2} is the same as $I( -s_2)$: $|I( s_2) |\leq 2 M(|s_2|)$. Integrating the estimate of $I(s_2)$ over $s_2$, we yield 
\[
\B| - \f{1}{\pi} P.V. \int_{ \R^2} K(s) W( s_1,  x_2 - s_2) ds \B|
\leq \f{4}{\pi}  [\om]_{C_x^{1/2}}  \B| \int_0^{\infty} ds_2\int_{ f(s_2)}^{\inf} K(s_1, s_2) | T(s_1, s_2) - s_1|^{1/2} d s_1 \B|,
\]
which along with \eqref{eq:holx_ux_int1} prove Lemma \ref{lem:holx_ux} in the case of $a = b_1 = b_2 = \infty$.

\subsubsection{Formula of $T$}

From \eqref{eq:holx_ux_ker} and $ \f{xy}{ (x^2 + y^2)^2 } = -\f{1}{2} \pa_x \f{y}{ x^2 + y^2}$,  equation \eqref{eq:cubic_T0} is equivalent to 
\[
 s_2 \B( \f{ 1 }{ (T+ \f{1}{2})^2 + s_2^2}  -  \f{ 1 }{ (T - \f{1}{2})^2 + s_2^2} \B)
=  s_2  \B( \f{ 1 }{ (s_1 + \f{1}{2})^2 + s_2^2   } -  \f{ 1 }{ (s_1 - \f{1}{2})^2 + s_2^2   } \B)
\]
where we have simplified $T(s_1, s_2)$ as $T$. For $s_2 \neq 0$, expanding the identity yields 
\beq\label{eq:cubic_T10}
0 = - 1 - 8 T s_1 + 16 T s_1 ( T^2 + T s_1 + s_1^2 ) - 8 s_2^2 ( 1 - 4  s_1 T + 2 s_2^2) .
\eeq

The above equation is cubic in $T$, and thus can be solved explicitly. In Appendix \ref{app:OT_map}, we show that it has a unique real root and derive its solution formula.

\begin{remark}[\bf Effectiveness of the estimates]\label{rem:holx_ux2}
When $\omega(x, y)$ is constant in $y$, we have $u_x(\omega)(x, 0) = H \omega(x)$, as observed in \cite{choi2014on,luo2013potentially-2}. Thus, the optimal constant in Lemma \ref{lem:holx_ux} must exceed that of the Hilbert transform \eqref{eq:hol_hil2}. Here, we obtain the upper bound $C_1(b) \leq 2.55$, very close to $C/\sqrt{2} \approx 2.37$ for the Hilbert transform \eqref{eq:hol_hil2}, reflecting the effectiveness of using 1D transport maps to construct the 2D transport map in this setting.
\end{remark}

\subsubsection{Localized estimate of $u_x$}

Next, we estimate
$u_x(\xx,a,b_1,b_2)-u_x(\zz,a,b_1,b_2)$ with $x_2=z_2$ using
$[W]_{C_x^{1/2}} = [\om]_{C_x^{1/2}}$. The argument parallels the non-localized case. 
We first show that the localized kernel has the same one-sign-change structure as
the full kernel, and construct the corresponding transport map in the $x$ direction. 
In Section~\ref{sec:comp_cost}, we then compare the localized and non-localized transportation costs and show that the localized cost is smaller.

Without loss of generality, we assume $ x_1 = \f{1}{2}, z_1 = -\f{1}{2}, x_2 = 0$. Denote 
\[
I_a \teq [-a, a], \quad I_b = [-b_1, b_2], \quad Q \teq I_a \times I_b, \quad b = \max(b_1, b_2), 
\]

Since we assume $a \geq \f{1}{2} |x_1 - z_1|$ in Lemma \ref{lem:holx_ux}, we have 
\beq\label{eq:holx_ux_a}
a  \geq 1/2.
\eeq

The  kernel associated with $u_x( \f{1}{2}, x_2) - u_x( -\f{1}{2}, x_2)$ \eqref{eq:ux_local} becomes 
\beq\label{eq:ux_ker_loc}
\bal
K_{a, b}(s_1, s_2) = & \one_{ s_2 \in I_b} \B( \f{ (s_1 + \f{1}{2}) s_2}{  ( (s_1 + \f{1}{2})^2 + s_2^2 )^2} \one_{ s_1 + \f{1}{2} \in I_a }
-  \f{ (s_1 - \f{1}{2}) s_2}{  ( (s_1 - \f{1}{2})^2 + s_2^2 )^2 } \one_{ s_1 - \f{1}{2} \in I_a }  \B)  \\
 = &   ( ( K_1(s_1 + 1/2, s_2) - K_1(s_1 - 1/2, s_2)) \one_{ s + 1/2 \in Q } \\
 & - K_1( s_1 - 1/2, s_2) ( \one_{ s - 1/2 \in Q} - \one_{s  + 1/2 \in Q } ) ) .
\eal
\eeq

Since $a \geq \f{1}{2}$ and $K_1(s)$ is odd in $s_1$, for fixed $s_2$,  we have
\[
\int_{0 }^{\inf} K_{a, b}(s_1, s_2) d s_1
= (\int_{ 1/2}^{  a} - \int_{-1/2}^{a} ) K_1(s ) d s_1
= - \int_{-1/2}^{1/2}  K_1(s) d s_1 = 0, \quad \int_{-\inf}^0 K_{a, b}(s_1, s_2) d s_1 = 0. 
\]

Similar to the case without localization, for each $s_2$, we consider the transportation from the positive part of $K_{a, b}$ to its negative part. Firstly, we identify the sign of $K_{a, b}$. 
We restrict to $s_2 \in [-b_1, b_2]$ and $s_2 \neq 0$ since otherwise $ K_{a, b} = 0$. We focus on $s_1 , s_2 \geq 0$ and the estimate for $s_1 <0$ or $s_2<0$ is the same. Since $a \geq \f{1}{2}$, we have 
\beq\label{eq:holx_ux_s1_low}
s_1 \pm 1/2 \geq -a , \quad \mathrm{for  } \quad s_1 \geq 0.
\eeq
Thus, for $s_1 \geq 0$, we can neglect the constraint $ s_1 \pm \f{1}{2} \geq -a$ in the localization in \eqref{eq:ux_ker_loc}.

\paragraph{\bf{Case 1}: $ a \in [1/2, 1]$ } 

Clearly, $K_{a, b}(s_1, s_2) \geq 0$ for $s_1 \leq \f{1}{2}$ since both kernels in \eqref{eq:ux_ker_loc} are non-negative. For $s_1 > \f{1}{2}$, since $a \leq 1$ and $s_1 + \f12 > 1$, we get 
\[
K_{a, b}(s_1, s_2) = - K_1(s_1 - \f{1}{2}, s_2) \one_{ s-1/2 \in Q} \leq 0.
\]

In this case, we denote $s_c(s_2) = \f{1}{2}$. 

\paragraph{\bf{Case 2}: $a \in (1,  f(s_2) + \f{1}{2} )$}

Recall $f(s_2)$ from \eqref{eq:ux_thres}. For $s_1 > a - \f{1}{2} > \f{1}{2}$, we get 
\[
K_{a, b} =  -K_1(s_1 - \f{1}{2}, s_2) \one_{ s-1/2 \in Q} \leq 0.
\]

For $s_1 \leq a - \f{1}{2} < f(s_2) $, using \eqref{eq:holx_ux_sign1} and $s_1 + \f{1}{2} \leq a$, we obtain
\[
K_{a, b} = K_1( s_1 + \f{1}{2}, s_2) - K_1( s_1 - \f{1}{2}, s_2) \geq 0.
\]

We denote $s_c(s_2) = a - \f{1}{2}$.

\paragraph{\bf{Case 3}: $a \geq f(s_2) + \f{1}{2} $}

For $ s_1 < f(s_2) $, using \eqref{eq:holx_ux_sign1} and $ s_1 \pm \f{1}{2} < a$, we get 
\[
K_{a, b} =  K_1( s_1 + \f{1}{2}, s_2) - K_1( s_1 - \f{1}{2}, s_2) \geq 0.
\]

For $s_1 \geq f(s_2) > \f{1}{2}$, since $K_1( s_1 + \f{1}{2}, s_2) - K_1( s_1 - \f{1}{2}, s_2) \leq 0$ \eqref{eq:holx_ux_sign1} and 
\[ \one_{ s - 1/2 \in Q} - \one_{s  + 1/2 \in Q } 
= \one_{ s_1 - 1/2 \leq a} - \one_{ s+1/2 \leq a} \geq 0, 
\]
we get 
\[
K_{a, b} \leq  - K_1( s_1 - 1/2, s_2) ( \one_{ s - 1/2 \in Q} - \one_{s  + 1/2 \in Q } ) ) \leq 0.
\]

We denote $s_c(s_2) = f(s_2)$. In summary, for fixed $s_2$, we define
\beq\label{eq:ux_thres_loc}
\bal
s_c(s_2) &= \f{1}{2},  &&  \mathrm{\ if \ } a  \in [\f{1}{2}, 1], \qquad \qquad s_c(s_2) = a - \f{1}{2}, \
\mathrm{\ if \ }  a \in (1, f(s_2) + \f{1}{2} ),  \\
s_c(s_2) &= f(s_2), && \mathrm{\ if \ }  a \geq f(s_2) +1/2 ,\eal
\eeq
which satisfies
\beq\label{eq:holx_ux_sign2}
 K_{a, b}( s_1, s_2) \geq 0 , \ s_1 \in [0, s_c], \quad K_{a, b}(s_1, s_2) \leq 0, \ s_1 \in [s_c, \infty], \quad  s_c(s_2) \leq f(s_2),
\eeq
where the last inequality follows from the definition of $s_c$ and $f(s_2) \geq \f{1}{2}$ \eqref{eq:ux_thres}.

In each case $i=1,2,3$, we construct the transport map by solving 
\beq\label{eq:cubic_T4}
\int_{T_i(s_1, s_2)}^{s_1} K_{a, b}( x, s_2) d x = 0, \quad T_i \leq a + \f{1}{2}.
\eeq

We add the restriction $T_i \leq a + \f{1}{2}$ since $ K_{a, b}(s) = 0$ for $s_1 > a + \f{1}{2}$ by definition \eqref{eq:ux_ker_loc}. Applying Lemma \ref{lem:trans} in the $s_1$ direction 
and using $K_{a,b}(s) = 0$ for $|s_2| \geq b$, we yield
 \[
\bal
I_i \teq  \B|\int_{s_1 \geq 0, s_2 \geq 0} K_{a, b}( s_1, s_2) \om( s_1, -s_2) ds \B|
\leq  \int_{ 0}^{b} \int_{ 0}^{ s_c(s_2)} | K_{a, b}(s)| | T_i(s) - s_1|^{1/2} ds \cdot  [\om]_{C_x^{1/2}}.
\eal
\]

\subsubsection{Comparison of the cost}\label{sec:comp_cost}

The guiding principle of this subsection is that localization can \emph{only decrease} the transportation cost. Indeed, after localization, the positive and negative parts of the singular kernel are supported in a smaller region, so the transportation distance needed to match them is expected to decrease. We now make this principle rigorous 
by proving
\beq\label{eq:ux_comp}
I_i \leq  \int_0^b \int_{ f(s_2)}^{\inf} |K(s)| \cdot | T(s) - s_1|^{1/2} ds \cdot [\om]_{C_x^{1/2}} ,
\eeq
where $T$ is defined in \eqref{eq:cubic_T0}. It suffices to prove, for any $s_2$, 
\beq\label{eq:ux_comp1}
J_i \teq \int_0^{ s_c(s_2)} | K_{a,b}(s)| | T_i(s) - s_1|^{1/2} ds 
\leq \int_0^{ f(s_2)} | K(s)| | T(s) - s_1|^{1/2} ds .
\eeq

 We focus on $s_2 \in [-b_1, b_2] $ and $s_2 \neq 0$. The proof is based on two comparisons. First, the localized kernel is no larger than the full kernel on $ [0, s_c(s_2)]$:
\beq\label{eq:ux_comp_ker}
|K_{a,b}(s)| \leq |K(s)|, \quad  s_1 \in [0, s_c(s_2)].
\eeq
Second, we use \eqref{eq:cubic_T0} and \eqref{eq:cubic_T4} to show that the localized transportation 
$T_i$ map does not transport mass farther than the non-localized map $T$:
\beq\label{eq:ux_comp2}
   s_1 \leq s_c(s_2)  \leq T_i(s) \leq T(s) ,  \quad s_1 \in [0, s_c(s_2)]
\eeq
and thus $T_i(s) - s_1 \leq T(s) - s_1$. Clearly, inequality \eqref{eq:ux_comp1} follows from \eqref{eq:ux_comp_ker} and \eqref{eq:ux_comp2}.

\vspace{0.1in}
\paragraph{\bf{Compare the kernels}}
From \eqref{eq:holx_ux_sign2} and \eqref{eq:holx_ux_sign1}, since $s_c(s_2) \leq f(s_2)$, we get $K_{a, b}(s) , K(s) \geq 0$ for $s_1 \in [0, s_c(s_2)]$. Hence, for fixed $s_2 \in[-b_1, b_2]$ ($b_1, b_2$ may differ), \eqref{eq:ux_comp_ker} is equivalent to 
\[
0 \leq K(s) - K_{a, b}(s) = 
K_1(s_1 + \f{1}{2}, s_2) (1 - \one_{s_1 + 1/2 \in I_a})
- K_1( s_1 - \f{1}{2}, s_2) (1 - \one_{s_1 - 1/2 \in I_a}) \teq I.
\]

From the definition of $s_c(s_2)$ \eqref{eq:ux_thres_loc} and \eqref{eq:holx_ux_s1_low}, for $s_1 \in [0, s_c(s_2)]$, we have 
\[
s_1 \pm 1/2 \geq -a, \quad s_1 - 1/2 \leq a , \quad 1 - \one_{s_1 - 1/2 \in I_a} = 0,
\]
which along with $K_1(s_1 + \f{1}{2}, s_2) \geq0$ \eqref{eq:ux_ker_loc} implies \eqref{eq:ux_comp_ker}
\[
I = K_1(s_1 + \f{1}{2}, s_2) (1 - \one_{s_1 + 1/2 \in I_a}) \geq 0.
\]

\begin{remark}
In the above derivations, we consider $s_2 \geq 0$. If $s_2 \leq 0$, one needs to track the sign to prove inequality \eqref{eq:ux_comp_ker}. 
\end{remark}

\paragraph{\bf{Compare the maps}}
To prove \eqref{eq:ux_comp2}, our idea is to use the equations \eqref{eq:cubic_T0}, \eqref{eq:cubic_T4} and the sign of the kernels $K_{a, b}, K$ to compare $T_i$ and $T$.

We fix $s_2 > 0$ in the following derivations. To simplify the notation, we simplify $T(s_1, s_2)$ as $T(s_1)$ in some places.
Since $T_i, T$ \eqref{eq:cubic_T0}, \eqref{eq:cubic_T4} are decreasing and $s_c(s_2)$ is a fixed point for $T_i(\cdot, s_2)$, for $s_1 \leq s_c(s_2)$, we get 
\beq\label{eq:ux_comp3}
T_i( s_1, s_2) \geq T_i( s_c(s_2), s_2) = s_c(s_2) , 
\quad T(s_1, s_2) \geq T(  f(s_2), s_2) = f(s_2)\geq s_c(s_2).
\eeq
Moreover, from \eqref{eq:cubic_T0}, \eqref{eq:cubic_T4}, we have 
\beq\label{eq:holx_ux_conv}
T_i(T_i(s_1)) = s_1, \quad T( T(s_1) ) = s_1.
\eeq

Denote 
\beq\label{eq:holx_ux_ker_pm}
K^{\pm} = K_1( s_1 \pm \f{1}{2}, s_2 ), \quad 
K_{a, b}^+ = K_1(s_1 + \f{1}{2}, s_2) \one_{ s_1 + \f{1}{2} \leq a}, \quad 
K_{a, b}^{-} =  K_1(s_1 - \f{1}{2}, s_2) \one_{ s_1 - \f{1}{2} \leq a}.
\eeq
We remark that $K^-$ is not non-negative but $K^+$ is positive. By definition, we have 
\beq\label{eq:holx_ux_ker_sign}
\bal
& K_{a, b} = K_{a, b}^+ - K_{a, b}^-, \quad K = K^+ - K^- .  \\
& K^+(s) \geq 0,  \ s_1, s_2 \geq 0, \quad K^-(s) \geq 0, \ s_1 \geq 1/2, s_2 \geq 0.
\eal
\eeq

Next, we study each case in the order of $3,2,1$ to prove \eqref{eq:ux_comp2}.

\paragraph{\bf{Case 3}:  $a \geq f(s_2) + \f{1}{2}$ }

In this case, recall  $s_c(s_2) = f(s_2)$ from \eqref{eq:ux_thres_loc}.

For $s_1 \leq a - \f{1}{2}$, we get $K_{a, b} = K$ \eqref{eq:ux_ker_loc}. Hence, equations \eqref{eq:cubic_T0} and \eqref{eq:cubic_T4} are the same for $s_1 \leq a-1/2$, and we get 
\beq\label{eq:holx_ux_case3_1}
T_3(s_1, s_2) = T(s_1, s_2) , \quad s_1 \in [ 
T( a - \f{1}{2}),  a - \f{1}{2}] . 
\eeq

It follows \eqref{eq:ux_comp2} for $s_1 \in [T(a-1/2), f(s_2)]$. We recall that from \eqref{eq:ux_comp3}, $a-1/2 \geq f(s_2)$ and $T(a-1/2) = T_3(a-1/2)$, we have 
\beq\label{eq:holx_ux_case3_2}
 T(a - 1/2) \leq T(f(s_2)) = f(s_2) \leq a-1/2 , \quad T(s_1) , T_3(s_1)  \geq a-1/2, \quad s_1 \leq T(a-1/2).
\eeq

Next, we compare $T(s_1), T_3(s_1)$ for $s_1 < T( a- 1/2) \leq f(s_2)$. From \eqref{eq:cubic_T0},\eqref{eq:cubic_T4}, and $T(a-1/2) = T_3(a-1/2) \leq a-1/2$, we have
\[
\bal
&\int_{ T(a - 1/2)}^{ a - 1/2} K(s) d s_1 = \int_{ T_3(a - 1/2)}^{ a - 1/2} K_{a, b}(s) d s_1 
= \int_{ T(a - 1/2)}^{ a - 1/2} K_{a, b}(s) d s_1 = 0 . \\
\eal
\]

Moreover, from \eqref{eq:ux_ker_loc} and \eqref{eq:cubic_T4}, we have 
\[
\bal
K_{a, b}(t, s_2) &= - K_{a,b}^-(t, s_2) = -K^-(t, s_2) , \quad t \in [a-1/2, T_3(s_1)] \subset [1-1/2, a+1/2],  \\
K_{a,b}(t, s_2) &= K(t, s_2) , \quad t \leq T(a-1/2) \leq a-1/2.
\eal
\]
Plugging the above identities in \eqref{eq:cubic_T0}, \eqref{eq:cubic_T4} for $s_1 \leq T(a-1/2)$, we yield 
\[
\bal
0 &= \int_{s_1}^{T(s_1) } K(t, s_2) d t   =  \int_{s_1}^{ T( a - 1/2)} K(t, s_2) + \int_{ a -1/2}^{ T}  K(t,s_2) d t  ,  \\
0 & = \int_{s_1}^{T_3(s_1)} K_{a, b}(t, s_2) d t
= \int_{s_1}^{ T( a - 1/2)} K(t, s_2) - \int_{a -1/2}^{ T_3(s_1)} K^-(t, s_2) d t  .
\eal
\]

Note that $K = K^+ - K^-$. Calculating the difference between the two identities yields
\[
0 =  \int_{a -1/2}^{T(s_1)} (K^+ - K^-) d s + \int_{a -1/2}^{T_3(s_1)} K^-  
= \int_{a-1/2}^{T(s_1)} K^+ d s_1 + \int_{T(s_1)}^{T_3(s_1)} K^- d s_1.
\]

Recall $s_2 \neq 0$ and from \eqref{eq:holx_ux_case3_2}, we obtain $T_3(s_1) ,T(s_1 ) \geq a-1/2 \geq 1/2$.
From \eqref{eq:holx_ux_ker_sign},  we yield $K^+ > 0$ and $K_- > 0$ for $s_1 \geq \min(a-1/2, T_3, T)$.  Since $T$ is decreasing and 
\[
T(s_1) \geq T( T(a-1/2)) =a-1/2, \quad s_1 \leq T(a-1/2),
\]
the first integral is non-negative. We prove  $T(s_1) \geq T_3(s_1)$ for $s_1 \leq T(a-1/2)$, which along with \eqref{eq:holx_ux_case3_1} implies  \eqref{eq:ux_comp2}.

The proof in the case 1,2 is completely similar. 
\vspace{0.1in}
\paragraph{\bf{Case 2}:  $a \in (1, f(s_2) + \f{1}{2}) $}

Recall  $s_c(s_2) = a - \f{1}{2} \leq f(s_2)$ from \eqref{eq:ux_thres_loc} and \eqref{eq:ux_comp3}. For any $s_1 \leq a - \f{1}{2} \leq f(s_2)$, 
using \eqref{eq:cubic_T0}, \eqref{eq:cubic_T4}, and an argument similar to that in case 3, we yield 
\[
\bal 
0 &= \int_{s_1}^{T_2(s_1)} K_{a, b}(t, s_2)  d t
= \int_{s_1}^{a - 1/2} K(t, s_2) d t - \int_{a-1/2}^{T_2(s_1)} K^-(t, s_2) d t,  \\
 0 &= \int_{s_1}^{T(s_1)} K(t, s_2) d t
 = \int_{s_1}^{a - 1/2} K(t, s_2) d t + \int_{a-1/2}^{T(s_1) }  (K^+ - K^-)(t, s_2)  d t ,
\eal
\]
where we have used  $K_{a,b}(t ,s_2) = -K^-(t, s_2)$ for $t \geq a-1/2$ \eqref{eq:ux_ker_loc}, \eqref{eq:holx_ux_ker_pm} in the second equality.

Comparing the difference between two identities yields 
\[
0 = \int_{a-1/2}^{T(s_1)} K^+ - K^-(t, s_2)  d t  + \int_{a-1/2}^{T_2(s_1) } K^-(t, s_2) d t =\int_{a-1/2}^{T(s_1) } K^+(t, s_2) d t + \int_{T(s_1)}^{T_2(s_1)}  K^- (t,s_2)d t .
\]

Recall from \eqref{eq:ux_comp3} that $T(s_1), T_2(s_1) \geq s_c(s_2) = a -1/2$ for $s_1 \leq a-1/2$. For $s_2 \neq 0$ and  $s_1 \geq \min( T, T_2, a-1/2)  = a -1/2 > 1/2$, we have $K^- > 0, K^+ > 0$ \eqref{eq:holx_ux_ker_sign}. We obtain $T(s_1) \geq T_2(s_1)$, which implies \eqref{eq:ux_comp2}.

\vspace{0.1in}
\paragraph{\bf{Case 1}: $a \in [1/2, 1]$}

In this case, $s_c(s_2) = \f{1}{2} < f(s_2)$. From \eqref{eq:ux_ker_loc}, 
we yield 
\[
\bal
0 &\leq K_{a, b} = \one_{s_1 + 1/2 \leq a} K_1(s_1 + 1/2, s_2) - K_1(s_1 -1/2, s_2) \\
&\leq K_1(s_1+1/2, s_2) - K_1(s_1-1/2, s_2) = K, \quad s_1 \in [0, 1/2],  \\
K_{a, b} &= -K_1(s_1 -1/2, s_2) = - K^-(s_1, s_2) , \quad K^-(s_1, s_2) \geq 0,   \quad s_1 \in [1/2, a + 1/2].
\eal
\]

For any $s_1 < \f{1}{2}$ and $s_2 \neq 0$, from \eqref{eq:ux_comp3}, we get $T_1(s_1) \geq 1/2 , T(s_1) \geq f(s_2) > 1/2$. Using \eqref{eq:cubic_T0}, \eqref{eq:cubic_T4} and the above estimates for $K_{a,b}$, we yield 
\[
0 = \int_{s_1}^{1/2} K_{a, b}(t, s_2) d t - \int_{1/2}^{T_1(s_1)} K^-(t, s_2) d t
= \int_{s_1}^{1/2} K(t, s_2) d t + \int_{1/2}^T (K^+ - K^-)(t, s_2) d t.
\]

It follows 
\[
0 = \int_{s_1}^{1/2} ( K - K_{a, b})(t, s_2) d t + \int_{1/2}^{T} K^+ + \int_T^{T_1} K^-(t,s_2) d t  \teq II_1 + II_2 + II_3 \;.
\]

From \eqref{eq:ux_comp_ker} and $K_{a, b}, K >0$ on $t \in [0, s_c(s_2)] = [0, 1/2]$, we get $II_1 \geq 0$. Note that $K^- , K^+ > 0$ for $s_1 > 1/2$ \eqref{eq:holx_ux_sign1}, \eqref{eq:holx_ux_sign2}. Since $T_1 , T  > 1/2$, we must obtain $T(s_1) \geq T_1(s_1)$, which implies \eqref{eq:ux_comp2}.

We have proved \eqref{eq:ux_comp2} in all three cases, which implies $|T(s) - s_1| \geq |T_i(s) - s_1|$. 
This shows that localization does not increase the transportation cost. Combining this estimate and \eqref{eq:ux_comp_ker}, we prove \eqref{eq:ux_comp} that the non-localized bound is larger and conclude the proof of Lemma \ref{lem:holx_ux}.

\subsection{Estimate of $[u_x]_{C_y^{1/2}}$}\label{sec:holy_ux}

Recall from Lemma~\ref{lem:holy_ux} that \(b_1=b_2=b\) and
\(x_1=z_1\) in this case.  The proof has an additional difficulty due to the
discontinuity of the odd extension \(W\) across the boundary \(y_2=0\).  We
therefore split the integral into two regions.  In the region \(y_2\ge m\),
$W$ belongs to \(C^{1/2}_y\), and the transport argument from
Section~\ref{sec:holx_ux} applies after a rotation.  In the region
\(y_2\le m\), $W$ is not H\"older continuous in $y$, so we
instead control the contribution by its \(C^{1/2}_x\) seminorm and a
monotonicity estimate for the truncated kernel.  By translation and scaling,
we assume without loss of generality that
\[
z_2 = m + 1/2, \quad  x_2 = m - 1/2 ,  \quad x_1 = z_1 = 0, \quad m \geq 1/2 .
\]
We introduce the localized kernel   
\[
K_{1, a, b}(s) \teq  K_1(s) \one_{ |s_1| \leq a, |s_2| \leq b} (s) .
\]
Then we have 
\[
u_x(z) - u_x(x) = \f{1}{ \pi} \int_{\R^2} W(y) \B( K_{1, a, b}( y_1, y_2 - (m-1/2))
- K_{1, a, b}( y_1 ,  y_2 - (m+1/2) ) \B) dy  , 
\]
where $W$ is the odd extension of $\om$ in $\R^2$ \eqref{eq:ext_w_odd}. Note that $W$ is not H\"older in the $y$-direction near $y_2 = 0$, we cannot use the same method as that in the estimate of $[u_x]_{C_x}^{1/2}$. On the other hand, since $W \in C_y^{1/2}(  \R \times [m, \infty) )$, we can apply the previous method to obtain 
\beq\label{eq:ux_holy_pf1}
\B| \f{1}{\pi} \int_{ y_2 \geq m} 
W(y) ( K_{1, a, b}( y_1 , y_2- (m-1/2) )
- K_{1, a, b}( y_1 , y_2- (m+1/2)  ) dy )  \B|
\leq \f{1}{ 2} C_1(a) [ \om ]_{ C_y^{1/2}} .
\eeq
Rotating the coordinate by $90$ degree, we obtain the case studied in Section \ref{sec:holx_ux}.

It remains to estimate 
\beq\label{eq:ux_holy_pf2}
\bal
I(b) & = \f{1}{  \pi} \int_{ y_2 \leq m } W(y) ( K_{1, a, b}( y_1, y_2 - (m-1/2) )
- K_{1, a, b}( y_1, y_2  - (m+1/2)  ) dy )  \\
& = \f{1}{ \pi} \int_{y_2 \leq 0} W(y_1, y_2 + m) ( K_{1, a, b}( y_1 , y_2 + 1/2)
- K_{1, a, b}( y_1 , y_2 - 1/2 ) dy ) . 
\eal
\eeq
Since $W$ is not H\"older continuous across $y=0$, we use $[W]_{C_x^{1/2}}$ to control $I$. Our idea is to compare the integral $I(b)$ with the case $ b = \infty$, $I(\infty)$. 
The following elementary monotonicity lemma allows us to compare the localized boundary contribution with the non-localized one.

\begin{lem}\label{lem:comp_trunc}
Suppose $f , f g \in L^1$ and $g \geq 0$ is monotone increasing on $[0, \infty]$.  
For any $0 \leq k \leq b \leq c$, we have 
\[
\int_{b - k}^{b+k}  |f(x - k) | g(x) dx 
\leq \int_{b-k}^{c-k} | f(x-k) - f(x+k)| g(x) dx 
+ \int_{c-k}^{c + k}  |f( x -k) | g(x) dx. 
\]
\end{lem}

\begin{proof}
Denote by $R, L$ the right and the left side of the above inequality, respectively. We have 
\[
\bal
R - L &\geq \int_{b-k}^{c-k} \B( |f(x-k) | - |f(x+ k)| \B) g(x) dx 
+ \int_{c-k}^{c + k}  |f( x -k) | g(x) dx - \int_{b-k}^{b+k} |f(x-k)| g(x) dx \\
& =\int_{b +k}^{c+k} |f(x-k)| g(x) dx 
- \int_{b-k}^{c-k} |f(x+k)| g(x) dx 
= \int_{b}^c |f(x) | ( g(x+k) - g(x-k) ) dx.
\eal
\]
Since $g$ is increasing on $ [0, \infty)$, we prove $R \geq L$.
\end{proof}

Now, we are in a position to estimate $I$. Since $K_{1, a,b}(y_1, y_2)$ is odd in $y_1$, we yield 
\[
\bal
|I | & \leq  \f{1}{\pi} \int_{y_2 \leq 0, y_1 \geq 0} \sqrt{2 y_1} \B| K_{1, a, b}(y_1, y_2 + 1/2) 
- K_{1, a,b}(y_1, y_2 - 1/2)\B| dy \cdot [ \omega ]_{C_x^{1/2}}  .\\
\eal
\]

For a fixed $y_1$ with $|y_1| \leq a$ and $b \geq 1/2$, using the definition of $K_{1, a, b}$ \eqref{eq:ux_local}, the even symmetry of $K_{1, a, b}(y_1, y_2 + 1/2)  - K_{1, a,b}(y_1, y_2 - 1/2)$ in $y_2$,  Lemma \ref{lem:comp_trunc} with $k =1/2$ and $c = \infty$, we get 
\[
\bal
&\int_{ y_2 \leq 0} | 
K_{1, a, b}(y_1, y_2 + 1/2) 
- K_{1, a,b}(y_1, y_2 - 1/2)| d y_2  \\
= & \int_{ y_2 \geq 0} | 
K_{1, a, b}(y_1, y_2 + 1/2) 
- K_{1, a,b}(y_1, y_2 - 1/2)| d y_2  \\
= & \int_0^{b-1/2 } |K_1(y_1, y_2 + 1/2) - 
K_1(y_1, y_2 - 1/2) | d y_2
+ \int_{b-1/2}^{b+1/2}  |  K_1(y_1, y_2 - 1/2) | dy_2  \\
\leq  &  \int_0^{ \infty } |K_1(y_1, y_2 + 1/2) - 
K_1(y_1, y_2 - 1/2) | d y_2 .
\eal
\]
For  $y_2 > 2  |y_1| + 1$,  since $|K_1(y_1, y_2 \pm \f12)| \les |y_2|^{-2}$, $ \int_{c-1/2}^{c+1/2}  |  K_1(y_1, y_2 - 1/2) | dy_2$ vanishes as $c \to \infty$.  Since $K_{1, a, b}(y) =0$ for $|y_1| \geq a$, integrating the above inequality in $y_1$ from $0$ to $a$,  we prove 
\[
|I | \leq  \f{1}{\pi} \int_0^a \int_0^{\inf} \sqrt{2 y_1} 
|K_1(y_1, y_2 + 1/2) - 
K_1(y_1, y_2 - 1/2) | d y \cdot [\om]_{C_x^{1/2}} .
\]

Combining the above estimate, \eqref{eq:ux_holy_pf1}, 
and \eqref{eq:ux_holy_pf2}, we prove Lemma \ref{lem:holy_ux}.

\section{$L^{\inf}$-based finite-rank perturbation method}\label{sec:finite_rank}

In this section, we introduce the $L^{\inf}$-based finite-rank perturbation
method to control the more regular nonlocal terms in the linearized equations.
The idea is to decompose the linearized operator into a stable main part
$\cL-\cK$ and a finite-rank operator $\cK$.

 In Section~\ref{sec:model2}, we explain the method in a rank-one
model. In Section~\ref{sec:finite_pertb}, we apply the method to the linearized Boussinesq system by first
decomposing \(W=W_1+W_2\) and then replacing \(W_2\) by the approximate
space-time solution \(\widehat W_2\). In Section~\ref{sec:appr_vel}, we construct the
finite-rank approximations of the regular velocity terms. In Section~\ref{sec:reg_van}, we  record
the qualitative regularity and vanishing properties needed to justify the
singular weighted estimates in Section~\ref{sec:EE}.

The decomposition, the representation formulas, the design of the finite-rank
operators, the analytic corrections enforcing the vanishing order, and the
qualitative regularity estimates for the perturbation in
Section~\ref{sec:reg_van} are all \emph{analytic}. The only computer-assisted
ingredients in this section are the construction of the approximate space-time
solutions $\hat F_i, \hat F_{\chi,i}$ representing $\hat W_2$ in
Section~\ref{sec:decoup_modi}, and the quantitative bounds on 
$\hat F_i, \hat F_{\chi,i}$ and  their residual errors.

\subsection{A toy model with a nonlocal term}\label{sec:model2}

We use a model problem to illustrate the ideas of stability analysis of a 
linearized equation perturbed from a simpler linearized equation.  Consider
\beq\label{eq:model_nloc}
 f_t = \cL_0 f + a(x) P(f) \teq \cL f, \quad P(f) = \int_{\R_2^{++} } f g dx,
\eeq
in $\R_2^{++}$, where $a, g$ are some given time-independent functions
and $a(x)$ may not be small. Operator $\cL_0$ models the local terms in \eqref{eq:bous_lin11}, \eqref{eq:bous_lin12}, and the rank-one operator $a(x) P(f)$ models the nonlocal terms. 
We assume that $\cL_0$ is linearly stable in $L^{\inf}(\vp)$ with some singular weight $\vp$, which can be studied following Section \ref{sec:model_local}, and $a(x) \in L^{\inf}(\vp)$.
We want to understand the long time behavior and the stability of the above model using the information of $\cL_0$.

A natural attempt is to use Duhamel's principle and the semi-group $e^{\cL_0 t}$ to represent the solution to \eqref{eq:model_nloc}. However, since $a(x)$ is not small, $a(x) P(f)$ cannot be absorbed perturbatively. 
\footnote{Another attempt is to project $f$ onto some space $Y$ orthogonal to $g(x)$ or $a(x)$ to eliminate the nonlocal term. However, the projection is not compatible with our $L^{\inf}$-based estimates.} 
Thus, the model isolates the central difficulty: the leading operator is stable, but
the finite-rank perturbation is not small and cannot be absorbed perturbatively in the stability estimate.

\subsubsection{Rank-one perturbation}\label{sec:rank_one}

Following the ideas in Section \ref{sec:frame_novel}, we decompose \eqref{eq:model_nloc} as follows 
\beq\label{eq:model2_decoup}
\bal
\partial_t f_1(t) &= \cL_0 f_1,  \quad f_1(0) = f_0 ,  \\
\partial_t f_2(t) &= \cL f_2 + a(x) P(f_1(t)), \quad f_2(0) = 0, 
\eal
\eeq
for initial data $f_0$, and then represent $f_2$ using Duhamel's principle 
\beq\label{eq:model_nloc4}
f_2(t) = \int_0^t P(f_1(s)) e^{\cL (t-s) } a(x) ds.
\eeq
Thus \(f_2\) is driven only by the scalar coefficient \(P(f_1(s))\) and the 
\emph{solution-independent} function \(e^{\cL s }a\). If $e^{\cL t} a(x)$ decays in $L^{\inf}(\vp)$ for large $t$, we can establish  $L^{\inf}(\vp)$ 
stability estimate of $\cL$.

Note that by choosing zero initial data for $f_2$ and using the fact that $P(f_1(s))$ is independent of space, we can solve $f_2$ for an arbitrary forcing coefficient $P(f_1(s))$.

A related idea appears in the $T(1)$ and $T(b)$ theorems in harmonic analysis \cite{david1984boundedness,david1985operateurs}, where one tests the operator on special functions to obtain boundedness. 
 In our setting, we estimate the finite-rank contribution by testing the decay of $e^{\cL t} a$
for the finitely many initial data $a(x)$ in the range of the finite-rank operator.
\footnote{
The analogy is only conceptual. In the $T(1)$ and $T(b)$ theorems, one tests
$T$ on $1$ or $b$ to obtain $L^2$ boundedness of a singular integral operator.
Our idea is also related in spirit to the Sherman-Morrison formula \cite{sherman1950adjustment} which connects the invertibility of $A \in R^{n\times n}$ and its rank-one perturbation.
}

 \subsubsection{Decay of $e^{\cL t} a$ and constructing approximate solution to $f_2$}\label{sec:model2_solu}

Though the operator $\cL$ and $a(x)$ \eqref{eq:model_nloc} are given, it is difficult to prove decay of $e^{\cL t} a$ in the weighted norm analytically since $\cL$ is nonlocal. The operator $\cL$ for the Boussinesq system \eqref{eq:bous_lin11} is even more complicated.

An alternative approach is to solve \eqref{eq:model_nloc} numerically from initial data $ a(x)$ to obtain an approximate solution $\hat g(t, x)$. Then by showing the error $ e^{\cL t} a - \hat g(t, x)$ is small and verifying the decay of $ \hat g(t)$, we obtain the decay estimates of $e^{\cL t} a$. 
However, it is difficult to estimate the error in the weighted norm rigorously. The bound obtained from standard \textit{a-priori} error estimate is neither explicit nor sharp enough for our purposes. 
\footnote{
The \textit{a-priori} error estimate gives $| e^{\cL t} a(x)- \hat g(t, x) | \leq C_1 (h^m + k^n) e^{C_2t},$
for some constants $C_1, C_2$ depending on $a(x)$ and $\cL$, where $h, k$ are the mesh size and time step in the computation, and $m, n$ are the  numerical scheme order. However, 
$C_1, C_2$ are difficult to estimate and may be large. 
Moreover, since we seek decay estimates for $ e^{\cL t} a$ with moderately large $t$, say $t\geq 10$, the factor $e^{C_2t}$ can render the estimate impractical.
}

Instead, we seek \textit{a-posteriori} error estimate. Firstly, we solve \eqref{eq:model_nloc} numerically and obtain a numerical solution $\hat g(t_k, x)$ at time $t_k$, which is represented by piecewise polynomials and thus defined globally in $x$. Then we interpolate the solution $\hat g(t_k, x)$ in time $t$ using piecewise cubic polynomials to obtain solution $\hat g(t, x)$ defined on $[0, T] \times \R_2^+$. Once \(\widehat g\) is fixed, we estimate the error  \textit{a-posteriori} and treat it as a small perturbation in the stability estimates.

For this purpose,  we introduce the residual error and a residual operator $\cR$ related to the nonlocal term $P(f_1(t))$ in \eqref{eq:model2_decoup}
\beq\label{eq:model_nloc5}
\bal
e(t, x)  & \teq  (\pa_t - \cL) \hat g, \quad  e_0(x) = \hat g(0) - g_0 , \quad g_0 = a(x) \\
\cR(f_1, t)  & = P(f_1(t)) e_0(x) + \int_0^t P(f_1(s)) e(t-s,x) ds.
\eal
 \eeq
 Since $\hat g$ is defined everywhere in space and time, we can estimate $e(t, x)$ and $e_0(x)$. Using the approximate solution $\hat g(t, x)$ for $e^{\cL t} g_0$, we construct the approximate solution to $f_2$  \eqref{eq:model_nloc4}
\beq\label{eq:model_nloc6}
\hat f_2(t) = \int_0^t P(f_1(s)) \hat g(t-s) ds.
\eeq
By definition and \eqref{eq:model_nloc5}, we have 
\[
\bal
(\pa_t -\cL) \hat f_2 
&= P(f_1(t)) \hat g(0)  + \int_0^t P(f_1(s)) (\pa_t - \cL) \hat g(t-s) ds  \\
&= P(f_1(t)) ( a(x) + e_0) + \int_0^t  P(f_1(s)) e(t - s, x) d s
= P(f_1(t)) a(x) + \cR(f_1, t).
\eal
\]

If the error $e(t, x)$ and $e_0$ are small, we can show that the norm of the residual operator is small in some suitable functional space (e.g., by imposing bootstrap assumptions on $f_1$)
\beq\label{eq:model_nloc7}
|| \cR(f_1, t) ||_X \leq \e || f_1||_X, \quad \e \ll 1. 
\eeq

Now, we modify the decomposition \eqref{eq:model2_decoup} as follows ($f = f_1 + \hat f_2 $)
\bseq\label{eq:model2_decoup2}
\begin{align}
\pa_t f_1  &= \cL_0 f_1 - \cR(f_1, t) , \quad  f_1(0) = f_0,  \label{eq:model2_decoup2:a} \\
\pa_t \hat f_2 &= \cL \hat f_2 + a(x) P(f_1(t)) + \cR(f_1 , t) .  \label{eq:model2_decoup2:b} 
\end{align}
\eseq
Thus the numerical approximation enters the \emph{exact} $f_1$-equations \eqref{eq:model2_decoup2:a} only through the residual operator \(\cR\). We remark that the solution \eqref{eq:model_nloc6} constructed by the numerical solution $\hat g$ solves equation  
\eqref{eq:model2_decoup2:b} \textit{exactly}. Now, due to the smallness \eqref{eq:model_nloc7}, we can apply the stability estimate of $\cL_0$ and treat $\cR(f_1, t)$ as perturbation to obtain stability estimate of $f_1$.

\begin{remark}
Using the decomposition \eqref{eq:model2_decoup},\eqref{eq:model2_decoup2}, constructing an approximating solution $\hat g$ to $e^{\cL t} a(x)$, and testing its decay, we replace a difficult nonlocal term in the original problem \eqref{eq:model_nloc} by a small error term $\cR(f_1, t)$ in \eqref{eq:model2_decoup2} that can be treated as a small perturbation. Moreover, we do not need to assume any specific form about the rank-one operator $a(x) P(f_1(t))$.
\end{remark}

\subsection{Finite-rank perturbations to the linearized operators}\label{sec:finite_pertb}

We generalize the idea from Section \ref{sec:model2} to the linearized Boussinesq
system. The construction below has four parts. First, we reorganize the nonlocal terms in the linearized equations \eqref{eq:bous_lin11},\eqref{eq:bous_lin12} so that they have better vanishing order near the origin.  Second, we introduce the near-origin corrections $\cK_{1i}$ and $\cNF_i$,
based on Taylor expansion, and the finite-rank approximations $\cK_{2i}$ for
the more regular nonlocal terms. Third, we decompose the perturbation into a
main part $W_1$ and a finite-rank part $W_2$, arranging the correction terms so
that  $W_1$-equation preserves \emph{cubic vanishing}, compatible with the singular weights used later. 
Finally, following Section~\ref{sec:model2_solu}, we replace $W_2$ by an approximate
component $\hat W_2$ with a controlled residual.

We perform the linear stability analysis on $\cL-\cK_1 - \cK_2$ in Section~\ref{sec:EE}; this plays the role
of the stability estimate for $\cL_0$ in the model problem of
Section~\ref{sec:model2}. The finite-rank contribution is then captured by a
second component via Duhamel's formula, as in \eqref{eq:model2_decoup}.

We begin with the first step. We introduce
 \beq\label{eq:u_tilde}
\td u = u - u_x(0)x, \quad \td v \teq  v -v_y(0) y = v + u_x(0) y. 
\eeq
Since $\om_x(0) = 0$ \eqref{eq:normal_vanish}, $ \om_x(0) = -\D \phi_x(0)$, $\phi(x, 0) = 0$, and $\phi$ is odd in $x$, for $\xx= (x,y)$, we yield 
\[
\phi = \phi_{xy}(0) x y +  O(|\xx|^4  ), \quad \td u = O(|\xx|^3), \quad \td v = O(|\xx|^3), \quad  \na \td \uu = O(|\xx|^2) ,
\]
near $\xx=0$,  for perturbations regular enough. Recall $c_{\om} = u_x(0)$ \eqref{eq:normal_pertb}. Using $\td u_x = u_x - u_x(0), \td u_y = u_y, \td v_x = v_x, \td v_y = v_y + u_x(0)$,  
\[
\bal
 - \uu \cdot  \na \bar \om + c_{\om} \bar \om &=  
 - \td \uu \cdot \na \bar \om + c_{\om}(\bar \om - x \bar \om_x + y \bar \om_y), \\
- \uu \cdot \na \bar \th_x - \uu_x \cdot \na \bar \th + 2 c_{\om} \bar \th_x
&= - \td \uu \cdot \na \bar \th_x - \td \uu_x \cdot  \na \bar \th  + c_{\om}( \bar \th_x - x \bar \th_{xx} + y \bar \th_{xy} ) , \\
- \uu \cdot \na \bar \th_y  - \uu_y \cdot \na \bar \th
+ 2 c_{\om} \bar \th_y  
&= - \td \uu \cdot \na \bar \th_y - \td \uu_y \cdot \na \bar \th + c_{\om}( 3 \bar \th_y - x \bar \th_{xy} + y \bar \th_{yy} ) , \\
 \eal
\]
and denoting
\beq\label{eq:f_cw}
 \bar f_{c_{\om}, 1} = \bar \om - x \bar \om_x + y \bar \om_y, 
\ \bar  f_{c_{\om}, 2} =  \bar \th_x - x \bar \th_{xx} + y \bar \th_{xy}  , \
 \ \bar f_{c_{\om}, 3} =  3 \bar \th_y - x \bar \th_{xy} + y \bar \th_{yy} ,
\eeq
we can rewrite \eqref{eq:bous_lin11}, \eqref{eq:bous_lin12} as follows 
\begin{align}\label{eq:lin}
\pa_t \om &= - (\bar c_l x  + \bar \uu) \cdot \na \om  +  \bar c_{\om} \om + \eta - \td \uu \cdot \na \bar \om + c_{\om} \bar f_{c_{\om}, 1} + \cN_1 + \olin \cF_1
\teq \cL_1  + \cN_1 + \olin \cF_1 
, \notag \\
\pa_t \eta &= - (\bar c_l x + \bar \uu) \cdot \na \eta  + (2 \bar c_{\om}- \bar u_x) \eta 
- \bar v_x \xi  - \td \uu_x \cdot  \na \bar \th - \td \uu \cdot \na \bar \th_x
+ c_{\om} \bar f_{c_{\om}, 2}  + \cN_2 +  \olin \cF_2 \notag  \\
& \teq \cL_2 + \cN_2 + \olin \cF_2, \notag \\
\pa_t \xi & = - (\bar c_l x + \bar \uu) \cdot \na \xi  + 
(2 \bar c_{\om} + \bar u_x) \xi - \bar u_y \eta 
- \td \uu_y \cdot  \na \bar \th - \td \uu \cdot \na \bar \th_y
+  c_{\om} \bar f_{c_{\om}, 3}  + \cN_3 + \olin \cF_3  \notag \\
&  \teq \cL_3  + \cN_3 + \olin \cF_3  . 
\end{align}

The nonlocal terms $c_{\om} \bar f_{c_{\om}, i}$, $ - \td \uu \cdot \na f$ for $f = \bar \om, \bar \th_x, \bar \th_y$, and $\na \td \uu_R$, the nonsingular part of the integral, are more regular than $\om$. We will choose finite-rank operators to approximate them.

\subsubsection{Correction near the origin}

As discussed in Section \ref{sec:vanish},
the use of singular weights for stability analysis 
requires the main component \(W_1\) to vanish cubically at the origin.  However, the linearized nonlocal
terms and the residual/nonlinear terms may generate quadratic Taylor modes.
The role of the corrections below is to remove precisely these low-order modes.
The operators \(\cK_{1i}\) correct the linear terms involving \(c_\omega\), while
\(\cNF_i\) correct the low-order Taylor coefficients from the nonlinear
terms and the residual.
This correction follows the same idea in Section \ref{sec:model_low_rank}.

We consider the corrections 
\beq\label{eq:appr_near0}
\bal
\cK_{1i}(\om)  &\teq c_{\om}(\om) \bar f_{c_{\om}, i},  \\
\cNF_{1}(\om, \eta, \xi) &=  ( c_{\om} \om_{xy}(0)    + \pa_{x y} \olin \cF_1(0) ) f_{\chi, 1}, \quad    f_{\chi, 1}  \teq  \D ( \tf16 x y^3 \cdot \chi_{\cNF}  )  ,  \\
\cNF_{2}(\om, \eta, \xi) & =   ( c_{\om} \eta_{xy}(0) + \pa_{xy} \olin \cF_2 (0) ) f_{\chi, 2}, \quad
f_{\chi, 2} \teq x y \cdot \chi_{\cNF} ,  \\
\cNF_{3}(\om, \eta, \xi) & = (  c_{\om} \xi_{xx}(0) + \pa_{xx} \olin \cF_3(0) ) f_{\chi, 3} , \quad
f_{\chi, 3} \teq \tf{1}{2} x^2  \cdot \chi_{\cNF},
\eal
\eeq
where $\chi_{\cNF}$ is some cutoff function with $\chi_{\cNF} = 1 + O(|x|^4)$ near $x=0$ constructed in \eqref{eq:cutoff_near0_all}. The form of $f_{\chi, 1}$ allows us to get $ \uu( f_{\chi, 1}) = \na^{\perp}  (-\D)^{-1}( f_{\chi, 1})$ analytically, and we have $f_{\chi, 1} = xy + h.o.t$. The operator $\cK_{1i}$ is a correction to the linear part, and $\cNF_i$ is a correction to the nonlinear term and the residual in \eqref{eq:lin}, respectively. 
These corrections are chosen so that the main component $W_1$ introduced
below preserves cubic vanishing near the origin.

To estimate  $ \om_{xy}(t, 0), \th_{xxy}(t, 0)$ in the operators \eqref{eq:appr_near0}, we derive their ODE using \eqref{eq:bousdy1}
\beq\label{eq:bous_ODE_xy}
\bal
\tf{d}{dt} \om_{xy}(t, 0) & = ( -2  c_l +  c_{\om}) \om_{xy}(t, 0) - \om_x^2(t, 0) + \th_{xxy}(t, 0),  \\
 \tf{d}{dt} \th_{xxy}(t, 0) & = ( -2  c_l + 2 c_{\om} - u_x(t, 0)) \th_{xxy}(t, 0) - 2 \om_x(t, 0) \th_{xx}(t, 0) .  \\ 
\eal
\eeq
The ODE is for the \emph{full perturbation} rather than $W_1$. Since $\om_x(t, 0), \th_{xx}(t, 0)$ are \emph{time-independent} \eqref{eq:normal1}, to estimate $\om_{xy}(t, 0), \th_{xxy}(t, 0)$, using \eqref{eq:bous_ODE_xy}, we only need to control $c_l, c_{\om}, u_x(t, 0)$ rather than some higher order norm of $\om, \th$, e.g. $|| \om ||_{C^2}, || \th||_{C^3}$.
This is important since it keeps the low-rank correction near the origin compatible with the
\emph{low-regularity} weighted energy framework.

\subsubsection{Approximation of the velocity}\label{sec:appr_vel_short}

In Section \ref{sec:appr_vel}, for  $f = \uu, \na \uu$,  we will construct finite-rank approximations $\appo{f}$ for $\td f$ (see \eqref{eq:u_appr}) so that the remainder $\td f-\appo{f}$ has smaller constants in the weighted estimates by the norms $ \| \om \vp \|_{L^{\inf}}, [ \om \psi_1 ]_{C_x^{1/2}}, [ \om \psi_1 ]_{C_y^{1/2}}$.  We remark that the approximations of $u, v, u_x, u_y, v_x$ are constructed \emph{separately}. 
Thus they need not satisfy the differentiation identities 
satisfied by the true velocity; for example,
we do not necessarily have
\[
\pa_x^i \pa_y^j \hat u = \wh{ \pa_x^i \pa_y^j u},  \quad 
\pa_x^i \pa_y^j \hat v = \wh{ \pa_x^i \pa_y^j v},
\]
for $i+j = 1$. Below, we only record how these approximations enter the finite-rank operators:
\beq\label{eq:appr_vel}
\bal
\cK_{21} = - \appo{\uu} \cdot \na \bar \om , \quad 
\cK_{22} = - \appow{\uu_x} \cdot \na \bar \th - \appo{ \uu} \cdot \bar  \na \th_x, 
\quad  \cK_{23} = - \appow{ \uu_y} \cdot \na \bar \th - \appo{  \uu} \cdot \na \bar  \th_y \, .
\eal
\eeq
These operators approximate the contributions from the regular nonlocal terms in
\eqref{eq:lin}.

\subsubsection{Decomposition of the system}\label{sec:decoup}

Denote $W_1 =(\om_1, \eta_1, \xi_1), W_2 = (\om_2, \eta_2 , \xi_2)$. 
We decompose the perturbation $(\om, \eta,\xi)$ in \eqref{eq:lin} into  $W_1$ and $W_2$, where \(W_1\) is the component estimated by the
weighted energy estimates, while \(W_2\) captures the finite-rank contribution and is represented by Duhamel's formula. Recall the notations \eqref{eq:bous_non} and \eqref{eq:bous_err}. Following Section \ref{sec:model2} and \eqref{eq:model2_decoup}, we decompose \eqref{eq:lin} as follows 
\beq\label{eq:bous_decoup}
 \bal
\pa_t W_{1, i}  & = (\cL_i -\cK_{1i} - \cK_{2 i}) W_1
+ \cN_i(W_1 + W_2) + \olin \cF_i - \cNF_i( W_1 + W_2)  ,\\
\pa_t W_{2, i} & =  \cL_i   W_2 + \cK_{1i }(W_1) + \cK_{2i}(W_1) + \cNF_i(W_1 + W_2) , \\
W_1 |_{t =0 }  &=  (\om_0, \eta_0, \xi_0), \quad W_2 |_{t = 0} = (0, 0, 0) , \\
 \eal
\eeq
with $\om_0, \eta_0, \xi_0$ being the initial perturbation with vanishing order $O(|x|^3)$. 
We have
\[
\pa_t (W_1 + W_2) = \cL_i (W_1 + W_2) + \cN_i(W_1 + W_2) + \olin \cF_i,
\]
which are the same equations as \eqref{eq:lin}. Since $W_1 + W_2$ has initial data $(\om_0, \eta_0, \xi_0)$, $W_1 + W_2$ solves \eqref{eq:lin} with the given initial data. 
Using the definitions \eqref{eq:appr_near0} and a Taylor expansion near $x=0$, we obtain that the vanishing conditions $\om_1, \eta_1, \xi_1 = O(|x|^3)$ are preserved.

\begin{remark}[$\mathbf{\pa_y W_{i, 2} \neq \pa_x W_{i, 3}}$]\label{rem:lin_th_deri}
Although $W_{1, 2} + W_{2, 2} = \th_x, W_{1, 3} + W_{2, 3} = \th_y$, since the finite-rank operators $\cK_{i j}$ we choose do not satisfy similar partial derivative relations, the solution 
to \eqref{eq:bous_decoup} does not satisfy $\pa_y W_{i, 2} = \pa_x W_{i, 3}$ for $i=1$ or $i=2$.
\end{remark}

Let us motivate the decomposition \eqref{eq:bous_decoup} and clarify the roles of the finite-rank terms.
At the linear level, we choose finite-rank operators $\cK_{1i}, \cK_{2i}$ to approximate  more regular non-local terms in $\cL_i$. Then $\cL_i - \cK_{1i}- \cK_{2i}, \cL_i$ serve as the $\cL_0, \cL$ operators in the model problem  \eqref{eq:model_nloc}, respectively. The decomposition of the solutions $W_1, W_2$ is similar to \eqref{eq:model2_decoup}. Since we want to perform energy estimate on $W_1$ using more singular weights, 
we design $ \cK_{1i}, \cNF_i$ to correct the linear, nonlinear and residual terms near the origin.
The ranges of these operators are finite-dimensional.  
Although $\cNF_i(W_1 + W_2)$  involves nonlinear factors, e.g. $c_{\om}(\om_1 + \om_2) \pa_{xy} ( \om_1 + \om_2) (0)$, since these factors are constant in space, we can still apply Duhamel's formula in \eqref{eq:model_nloc4} to $\cNF_i(W_1 + W_2)$, i.e.,
\[
\int_0^t e^{ \cL (t-s) }  \B( c_{\om} \pa_{xy}(\om_1(s) + \om_2(s) )(0)  \bar f \B) ds
= \int_0^t c_{\om} \pa_{xy}(\om_1(s) + \om_2(s) )(0)  e^{\cL( t-s)} \bar f ds,
\]
and obtain the formula of $W_2$ in \eqref{eq:bous_decoup}. 

\vspace{0.1in}
\paragraph{\bf{Avoid the loss of derivatives}}

Note that in the equation of $W_1$ in \eqref{eq:bous_decoup}, it contains the nonlinear terms $\uu(W_1 + W_2) \cdot \na (W_1 + W_2) $. In general, the term $\uu \cdot \na W_2$ can lead to loss of derivatives. Note that $W_2$ in \eqref{eq:bous_decoup} is driven by the forcing terms of the following forms 
\beq\label{eq:bous_W2_1}
\sum\nolimits_{1\leq i \leq N} a_i(W_1, W_2)(t) f_i
\eeq
for some $N$, time-dependent scalars (independent of $x$)  $a_i(W_1, W_2)(t)$, and time-independent functions $f_i$, e.g. $ c_{\om}( W_1) \bar f_{c_{\om}, i}$ in \eqref{eq:appr_near0}. By choosing smoother $f_i$ in the approximation, we can obtain solution $W_2$ smooth enough for our energy estimates and overcome the above difficulty.

\subsubsection{Constructing the approximate solution of $W_2$ and modifying the decomposition}\label{sec:decoup_modi}

Following Section~\ref{sec:model2_solu}, we replace the exact evolution of the
finite-rank forcing terms in the $W_2$ equation by approximate space-time
solutions. 

There are two types of finite-rank terms. The first has coefficients 
$a_i(W_1)(t)$, given by linear functionals of $W_1$, and range functions $\bar F_i$; the
second is $\cNF_i$ in \eqref{eq:appr_near0} and has coefficients
$a_{\mf{nl},i}(W)(t)$, given by nonlinear functionals of $W$, and range functions
$\bar F_{\chi,i}$.

For the first type, we represent 
\beq\label{eq:bous_finite_rank1}
\cK_{1 j}(W_1) + \cK_{2  j}(W_1) 
= \sum\nolim_{1\leq i \leq n_1} a_i(W_1)(t) \cdot \olin F_{ij},
\quad j  \in \{1,2,3\} ,
\eeq
 where $\vec{\cK}_{i \cdot} = (\cK_{i 1} , \cK_{i2}, \cK_{i3} )$, $\bar F_i(x) = (\bar F_{i, 1}(x), \bar F_{i, 2}(x), \bar F_{i, 3}(x) ) : \R^2_+ \to \R \times \R \times \R$, and $a_i(W_1)(t)$ is some linear functional on $W_1$. For example, the formula \eqref{eq:appr_near0} can be written as 
\[
a(W_1)(t)  \cdot (\bar f_{ c_{\om}, 1}, \bar f_{c_{\om}, 2}, \bar f_{ c_{\om}, 3 } ), \quad a(W_1)(t) =
c_{\om}(W_1)=  u_x(W_1(t))(0) ,
\]
where we have used \eqref{eq:normal_pertb} for $c_{\om}$. Recall the operators $\cNF_i$ and functions $f_{\chi, 1}$ from \eqref{eq:appr_near0}. Writing \eqref{eq:appr_near0} as vectors, we have 
\beq\label{eq:appr_near02}
\bal
\cNF( W ) &=  (c_{\om} \pa_{xy}\om(0) + \pa_{xy} \olin \cF_1(0) ) \cdot ( f_{\chi, 1}, 0, 0)
+  (c_{\om} \pa_{xy}\eta(0) + \pa_{xy} \olin \cF_2(0) ) \cdot ( 0, f_{\chi, 2},  0) \\
& \quad +  (c_{\om} \pa_{xx} \xi(0) + \pa_{xx} \olin \cF_3(0) ) \cdot ( 0, 0, f_{\chi, 3}) \\
& \teq \tts{\sum}_{1\leq i \leq 3} \  a_{ \mf{nl}, i}( W(t)) \bar  F_{\chi, i},
\eal
\eeq
where $\mf{nl}$ is short for \emph{nonlinear},
\beq\label{eq:appr_near02_op}
\bal
  a_{\mf{nl}, \cdot} &= 
(c_{\om} \pa_{xy}\om(0) + \pa_{xy} \olin \cF_1(0) , \, c_{\om} \pa_{xy}\eta(0) + \pa_{xy} \olin \cF_2(0) , \, c_{\om} \pa_{xx} \xi(0) + \pa_{xx} \olin \cF_3(0) ),  \\
\bar F_{\chi, i} & \teq  f_{\chi, i} \ee_i, 
\eal
\eeq
and $\ee_1, \ee_2, \ee_3$ are the standard basis for $\R^3$.   Denote by $\hat F_i(t, x)$ and $ \hat F_{\chi, i}(t, x)$ the approximation of $e^{\cL t} \bar F_i$ and $e^{\cL t} \bar F_{\chi, i} $. Following \eqref{eq:model_nloc5} and \eqref{eq:model_nloc6}, and using the idea in \eqref{eq:bous_W2_1}, we construct the approximate solution to $W_2$ in \eqref{eq:bous_decoup} as follows
\beq\label{eq:bous_W2_appr}
\hat W_2(t) = \sum_{i \leq n_1} \int_0^t a_i(W_1(s)) \hat F_i( t-s) ds
+ \sum_{ i \leq 3} \int_0^t a_{\mf{nl}, i}(W_1(s) + \hat {W_2}(s))  \hat F_{\chi, i}(t - s) ds.
\eeq

We introduce the residual operator
\begin{align}\label{eq:bous_err_op}
& \cR_l(W_1)  \teq \sum\nolim_{1\leq i \leq n_1} \B( a_i( W_1(t)) ( \hat F_i(0, \cdot) - \bar F_i )
   + \int_0^t a_i(W_1(s) )  ( (\pa_t - \cL ) \hat F_i )(t-s) ds \B), \notag  \\
 & \cR_{ \mf{nl} }( W)  \teq \sum\nolim_{1\leq i \leq 3} 
\B( a_{ \mf{nl}, i}( W(t)) ( \hat F_{\chi, i}(0, \cdot) - \bar F_{\chi, i} )
   + \int_0^t a_{ \mf{nl},i}(W(s) )  ( (\pa_t - \cL ) \hat F_{\chi,i} )(t-s) ds \B), \notag 
     \\
& \cR(W_1, \hat W_2)  \teq \cR_{l}(W_1) + \cR_{\mf{nl} }(W_1 + \hat W_2) ,  
\end{align}
where $\cR_l$ measures the initial error and evolution error of the
approximate solutions $\widehat F_i$, while \(\cR_{nl}\) measures the
corresponding error for \(\widehat F_{\chi,i}\). Note that $\cR(W_1, \hat W_2)$ is a vector in $\R^3$.

\begin{remark}
Given $W_1, \hat F_i, \hat F_{\chi, i}$, 
the formula  \eqref{eq:bous_W2_appr} does not determine $\hat W_2$  completely since  $a_{\mf{nl},i}(W_1+\hat W_2)$ in \eqref{eq:appr_near02_op} depends on $\hat W_2$. At the linear level, $\hat W_2$  \eqref{eq:bous_W2_appr} is determined. 
We show that $a_{\mf{nl},i}(W_1+\hat W_2)$ \eqref{eq:appr_near02_op} is much smaller than the linear part (the first part in \eqref{eq:bous_W2_appr}) and control $a_{ \mf{nl}, i}(W(s))$ using a bootstrap condition \eqref{eq:W2_non_boot}. Then we can still use \eqref{eq:bous_W2_appr} to estimate $\hat W_2$ effectively. See more discussion on local-in-time regularity  in Section \ref{sec:reg_van}.
\end{remark}

Similar to \eqref{eq:model2_decoup2}, using the above operators, we modify the decomposition \eqref{eq:bous_decoup} as follows 
\bseq\label{eq:bous_decoup2}
\begin{align}
\pa_t W_{1, i}  & = (\cL_i -\cK_{1i} - \cK_{2 i}) (W_1)
+ \cN_i(W_1 + \hat  W_2) + \olin \cF_i - \cNF_i( W_1 + \hat W_2) - \cR_i(W_1, \hat W_2) , \label{eq:bous_decoup2:W1} \\
\pa_t \hat W_{2, i} & =  \cL_i  \hat W_2 + \cK_{1i }(W_1) + \cK_{2i}(W_1) + \cNF_i(W_1 + \hat W_2) 
+ \cR_i( W_1, \hat W_2) ,   \label{eq:bous_decoup2:W2} \\
 W_1 |_{t =0 } & =   (\om_0, \eta_0, \xi_0), \quad \hat W_2 |_{t = 0} = (0, 0, 0) ,  \notag 
 \end{align} 
where $\cR = (\cR_1, \cR_2, \cR_3)$. The above decomposition is a nonlinear generalization of \eqref{eq:model2_decoup2}.  The residual operator $\cR$ is the  Boussinesq analogue of the residual operator 
$\cR(f_1,t)$ in \eqref{eq:model_nloc5}--\eqref{eq:model2_decoup2}, with additional nonlinear
terms $a_{\mf{nl},i}(W_1+\hat W_2)$ coming from $\cNF_i$. With the definition \eqref{eq:bous_err_op}, we solve the $\hat W_2$ equation using the formula \eqref{eq:bous_W2_appr} \textit{exactly}. It is easy to see that 
\beq
  W_1 + \hat W_2 = (\om, \eta, \xi)
\eeq
\eseq
solves \eqref{eq:lin} from initial data $(\om_0, \eta_0, \xi_0)$. 
If the errors in \eqref{eq:bous_err_op}, e.g. $\hat F_i(0, \cdot) - \bar F_i, (\pa_t - \cL) \hat F_i$,  and the residual error $\bar \cF_i$ \eqref{eq:bous_err} are small, we expect the following estimates for $\cR$
\[
\| \cR(W_1, \hat W_2) \|_X \leq \e ( \| W_1 \|_X +  \| W_1 + \hat W_2 \|^2_X   + \bar \e)
\]
in some suitable weighted space $X$ with very small $\e, \bar \e$, where 
the last two terms come from the estimate of $a_{ \mf{nl},i}(W_1 + \hat W_2)$ in \eqref{eq:appr_near02}. Then, the residual operator $\cR$ can be treated as a small perturbation in \eqref{eq:bous_decoup2}.  
At the linear level, $\hat W_2$ is almost decoupled from the $W_1$ equation.

\begin{remark}[\bf Smallness of error]
Since $\hat F_i$ solves $\pa_t F_i = \cL F_i$ numerically, and the initial data and coefficients of $\cL$ are smooth enough, in principle, by choosing a high order numerical scheme with sufficiently small mesh size and time step, one can make the errors $\hat F_i(0, \cdot) - \bar F_i, \, (\pa_t - \cL) \hat F_i$ very small. 
\end{remark}

\paragraph{\bf Cubic vanishing}

For initial perturbations $\om_0,\eta_0,\xi_0=O(|\xx|^3)$,
the vanishing orders of the error in \eqref{eq:bous_err_order}, 
\eqref{eq:appr_near0},  and equation \eqref{eq:bous_decoup2} imply that $\om_1,\eta_1,\xi_1=O(|\xx|^3)$ is preserved. In
Section~\ref{sec:reg_van}, we discuss the regularity of $W_1,\hat W_2$ in
\eqref{eq:bous_decoup2} and justify these cubic vanishing orders.

\subsubsection{Analytic low-rank corrections and construction of $\hat F_i$}
\label{sec:low_rank_cor}

We now explain how the approximate space-time solutions in
\eqref{eq:bous_W2_appr}, \eqref{eq:bous_err_op} are constructed. 
We have discussed the construction  for a model problem in Section~\ref{sec:model_low_rank}; here we record the corresponding construction for the full Boussinesq system. The formulas
for the initial data $\bar F_i$ and $\bar F_{\chi,i}$ are given in
Appendix~\ref{app:init_data}.

For each initial datum $\bar G\in\{\bar F_i,\bar F_{\chi,i}\}$, let $\hat G$ be the corrected approximate solution associated with $\pa_tG=\cL G$ and initial data $\bar G$. Writing $\hat G=(\hat\om,\hat\eta,\hat\xi)$, for $\xx= (x, y)$, we require
\footnote{
 The regularity of $(\bar \th_x, \bar \th_y)$ differs from that of the $(\eta, \xi)$ components of $\hat G$ since we use representations based on the 6th order B-spline to construct
    $(\bar \om, \bar \th)$ and $\hat G$. The basis for $\hat G$ is one order more regular than that of $\na  \bar \th$.
} 
\begin{subequations}\label{eq:bous_err_order}
\begin{align}
  \hat G(t,\cdot) & =O(|\xx|^2),\qquad
 & &  \hat \om , \ \hat \eta \ \mathrm{odd\ in}\ x,\qquad
  \hat \xi \ \mathrm{even\ in}\ x,  \quad  \hat \xi(t, 0, y) = 0,
    \label{eq:bous_err_order:a} \\
 \hat G(0,\xx)-\bar G(\xx) & = O(|\xx|^3), \qquad
 & &  (\pa_t-\cL)\hat G(t,\xx)=O(|\xx|^3), 
\label{eq:bous_err_order:b}  \\
 \hat G(t) & \in C^{4, 1}, \quad  & & \supp(\hat G(t)) \subset [-L_2, L_2] \times [0, L_2], \quad L_2 < 10^{16} 
\label{eq:bous_err_order:c}  .
\end{align}
\end{subequations}

\paragraph{\bf Construction of approximate space-time solutions.}

For each initial datum  $\bar G\in\{\bar F_i,\bar F_{\chi,i}\}$, we first construct a raw approximate solution $\hat G^{(0)}$ by numerically solving $\pa_t G=\cL G$ at
discrete times $t_k$ up to a finite time $T_i$.  We then set the approximate
solution $\hat G^{(0)}$ to zero for $t>T_i$ and interpolate  $\hat G^{(0)}(t_k)$ in time using piecewise cubic polynomials to obtain $\hat G^{(0)}(t)$ defined for any $ t \geq 0$. 
\footnote{
The jump discontinuity of $\hat G^{(0)}(t,\cdot)$ at $t=T_i$ contributes boundary terms in the residual operator \eqref{eq:bous_err_op}, via terms such as $(\pa_t-\cL)\hat F_i$. These boundary terms are estimated rigorously; see the formulas in \cite[Section \secsolucomp]{ChenHou2023b}.
}
We represent $\hat G^{(0)}(t,\xx)$ by explicit, compactly supported basis functions in space, rather than merely by its values on grid points. The spatial basis consists of piecewise 6th order B-splines with $C^{4,1}$ regularity and imposes the symmetries \eqref{eq:bous_err_order:a}. This representation gives the compact support and smoothness properties in \eqref{eq:bous_err_order:c},
and symmetries in \eqref{eq:bous_err_order:a}.
\footnote{
These symmetries are consistent with \eqref{eq:Sym}. For $\th(x,y)$ even in
$x$ with $\th(0,y)=0$, we have $\na \th(0,y)=0$.
}

\vs{0.05in}

\paragraph{\bf Analytic corrections and vanishing order.}
The raw numerical solution $\hat G^{(0)}$ need not satisfy the exact vanishing conditions in
\eqref{eq:bous_err_order}. We therefore apply analytic low-rank corrections
near $\xx=0$, following the methods explained in Section~\ref{sec:numer}, and denote the corrected approximate solution by $\hat G$. The first correction restores the quadratic vanishing of the approximate solution in \eqref{eq:bous_err_order:a}, and the second correction enforces the cubic
vanishing of the initial error and residual error in
\eqref{eq:bous_err_order:b}. These corrections are performed at the level of
the explicitly represented functions. Hence the improved vanishing in
\eqref{eq:bous_err_order} is an \emph{exact analytic property} of the
corrected functions, even though the coefficients of these functions are
obtained numerically. 

The full construction of $\hat G$ with \emph{analytic corrections} is developed in detail  in Part II \cite[Section {\seclinevo}]{ChenHou2023b}, especially \cite[Section 3.5]{ChenHou2023b}.

\vs{0.05in}

\paragraph{\bf Ideas of Residual estimates.}

We decompose the residual errors in \eqref{eq:bous_err_order:b} into local and
nonlocal parts:
\[
 (\pa_t-\cL)\hat G
        =
        \cE_{loc}(\hat G)+\cE_{nloc}(\hat G), \quad \hat G = (\hat \om, \hat \eta, \hat \xi) .
\]
Here  $\cE_{loc}(\hat G)$ depends locally on represented basis functions for $\hat G$ and their derivatives, e.g.  error associated with $( \bar c_l \xx + \bar \uu ) \cdot \na \hat G$ in \eqref{eq:lin}.  The nonlocal error part depends on $\hat G$ through the nonlocal velocity, 
such as the error associated with $\td \uu_x(\hat \om)\cdot \na \bar \th$ in
\eqref{eq:lin}. 
\footnote{
For the approximate functions, this error comes from the
error in the numerical stream-function solver: if $\psi^N$ denotes the
numerical stream function, then the error $-\D\psi^N-\om$ contributes to
$\cE_{nloc}(\hat G)$. 
}
We estimate $\cE_{loc}(\hat G)$ directly from the basis representation; some
ideas are outlined in Appendix~\ref{app:idea_num_est}. The detailed local
error estimates for $\hat F_i$ and $\hat F_{\chi,i}$ are developed in
Part~II~\cite[Section~{\seclinevo}]{ChenHou2023b}, using numerical analysis. The nonlocal error estimates are combined with the estimate of nonlocal perturbation, e.g. 
$\td \uu_x \cdot \na \bar \th$ \eqref{eq:lin}, in the energy estimate. These estimates are used to control the residual operator $\cR$ in \eqref{eq:bous_err_op}, and enter later in the estimates of $\hat W_2$; see Sections~\ref{sec:W2} and \ref{sec:comb_vel_err}.

Similarly, for the residual error $\bar \cF_i$ in \eqref{eq:bous_err}, we
decompose $\bar \cF_i=\bar \cF_{loc,i}+\bar \cF_{nloc,i}$. The local part 
$\bar \cF_{loc,i}$ is estimated using basis representation and numerical analysis, with details in Part~II~\cite[Appendix~{\secresid}]{ChenHou2023b}. The nonlocal part is combined with the energy estimates in Section~\ref{sec:comb_vel_err}.

\vs{0.05in}

\paragraph{\bf Illustration of decay of approximate solution}

\begin{figure}[t]
   \centering
      \includegraphics[width = \textwidth  ]{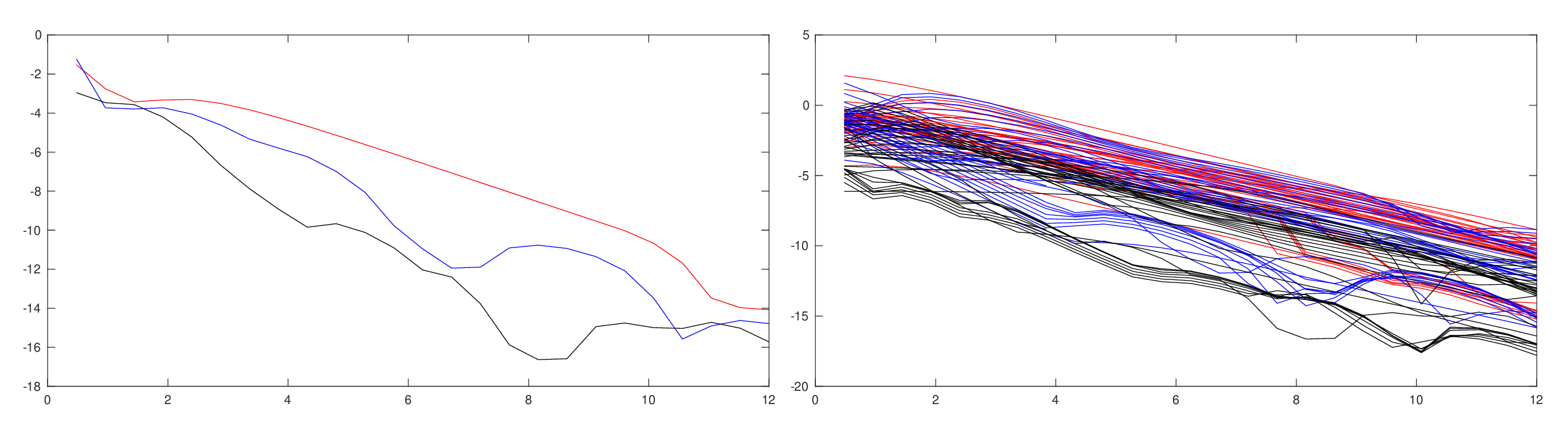}
      \caption{Plot of weighted discrete norm $\log || \hat f(t_i) \rho ||_{\ell^{\inf}}$, where $\hat f= \hat \om, \hat \eta, \hat \xi$, $\rho(\xx) = |\xx|^{-2} + 1$, and the $\ell^{\infty}$ norm  is taken over the grid points. Here $t_i=0.48i$ and $0\le t_i\le 12$. 
      Red, blue, and black curves represent $\hat \om, \hat \eta, \hat \xi$, respectively. Left: the first mode related to $\bar f_{c_{\om}, i}$ \eqref{eq:appr_near0}. Right: all modes. 
}
            \label{fig:lin_decay}
 \end{figure}

In Figure \ref{fig:lin_decay}, we plot $\log \| \hat f(t) \rho \|_{\ell^{\infty}} ,\rho = |x|^{-2} + 1$ with discrete $\ell^{\inf}$ norm computed over the grid points to illustrate the decay of $|| e^{\cL t} f \cdot \rho||_{\ell^{\inf}}$.
Here $\hat f$ is an approximate space-time solution from initial data in Appendix \ref{app:init_data}.
Over the time interval $[0, 12]$, the discrete norm  decreases by a factor about $e^{-10} \approx 4.5 \cdot 10^{-5}$. In Part II \cite{ChenHou2023b}, we use the method in Section \ref{sec:model2_solu} to estimate the decay rigorously. The observed exponential decay  is consistent with the numerical evidence of linear stability reported by Liu \cite[Section 3.4]{liu2017spatial}. By constructing approximate solution to $e^{\cL t} f$, we exploit this spectral property rigorously.

\subsection{Approximating the regular part of the velocity}\label{sec:appr_vel}

Recall that the nonlocal velocity $\uu, \na \uu$
appearing in \eqref{eq:lin} can be represented as singular integral
operators acting on the vorticity:
\[
f(\om)(x)
=
\int_{\R^2_+} K_f(x-y)\om(y)\,dy,
\qquad
f\in \{u,v,u_x,v_x,u_y\},
\]
where $K_f$ denotes the corresponding kernel. To estimate the nonlocal
terms efficiently, we approximate these operators by some finite-rank
operators $ \hat f(\om)$. 
Since $\na \cdot \uu = u_x + v_y = 0$, 
we have $v_y = -u_x$.  For each $f\in \{u,v,u_x,v_x,u_y\}$, we construct
$\hat f(\om)$ so that
\beq\label{eq:vel_appr_err}
\B|
\int_{\R^2_+} K_f(x-y)\om(y)\,dy - \hat f(\om)(x) \B|
\leq
C_1(x,\g)\,
\max \B(
\|\om\vp_1\|_{L^\infty}, \, s_f \max_{i=1,2}\g_i [\om\psi_1]_{C_{x_i}^{1/2}}
 \B).
\eeq
Here $s_f=0$ for $f=u,v$, so that no $C^{1/2}$ seminorm is used for the
velocity itself, and $s_f=1$ for $f=u_x,v_x,u_y$. 
The parameter $\gamma$ is related to the H\"older seminorm 
for $\om$ in the energy $E_1$ \eqref{energy1}. The kernel $K_f$ is the kernel associated with
$\pa_x^i\pa_y^j(-\D)^{-1}\om$, $i+j\leq 2$, and the H\"older seminorms
$C_x^{1/2}, C_y^{1/2}$ are defined in \eqref{hol:semi}.
Given the weights $\vp_1, \psi_1, \gamma_i$ and kernels $K_f$, by refining the approximation $\hat f$, we aim to obtain a relatively small constant $C_1(x,\g)$.

\vs{0.05in}

\paragraph{\bf Design of the approximation.}

The approximation is designed differently in three regions,

(1) In Section~\ref{sec:appr_vel_near0}, we approximate the velocity near
$\xx=0$ by its leading-order Taylor expansion. This expansion depends
only on two scalar functionals of $\om$: $u_x(\om)(0)$ and
$\cK_{00}(\om)$.

(2) In Section~\ref{sec:u_appr_1st}, we interpolate the nonsingular part
of the velocity along the boundary $y=0$ at finitely many points
$(x_i,0)$.

(3) In Section~\ref{sec:u_appr_far}, we approximate the velocity 
in the  far-field by the leading-order term of the symmetrized kernel
expansion.

The first two ingredients are combined into a near-field
approximation $\hat f_1$, while the third produces a far-field
approximation $\hat f_2$. The final approximation is $\hat f=\hat f_1+\hat f_2$ in 
\eqref{eq:u_appr}.

\vs{0.05in}

We first explain the basic interpolation idea for the near-origin and boundary approximation.  Since $K_f(z)$ is smooth away from $z=0$, a natural approach is to approximate the nonsingular part of $K(x - y)$ by interpolating $K(x - y)$ on finitely many points $x_i$: 
\[
K(x- y) \one_{|x-y| \geq \e} \approx \sum\nolimits_{1\leq i \leq n} 
\chi_i(x) K(x_i - y) \one_{|x_i-y| \geq \e} ,
\]
for some $\e > 0$, where $\{\chi_i\}_{i=1}^n$ is a partition of unity on the domain $D$ of approximation, i.e. $\sum_i \chi_i(x) =1$ for all $x\in D$, with $\chi_i$ localized near $x_i$. This leads to the finite-rank approximation
\beq\label{eq:u_appr_idea}
\hat f_{\rm idea}(\om)(x) = 
\sum\nolimits_{1\leq i \leq n} 
\chi_i(x) \int K(x_i - y) \one_{|x_i-y| \geq \e} \om(y) dy.
\eeq

Due to the decay of the coefficients 
of the nonlocal terms in
\eqref{eq:lin}, such as $\na\bar\om$ in $\uu\cdot\na\bar\om$ and
$\na\bar\th$ in $\na\uu\cdot\na\bar\th$, 
these terms are small for large $|x|$. Therefore, it is enough to approximate
$\uu$ and $\na\uu$ and obtain improved estimates with smaller $C_1(\xx, \g)$ \eqref{eq:vel_appr_err} in the near field. This is especially important
near the boundary, where the damping factor is relatively small; see
Section~\ref{sec:aniso_est}.

\vs{0.05in}

\paragraph{\bf Qualitative vs quantitative regularity} 

Although we consider sufficiently regular solutions, e.g.,
$\om+\bar\om\in C^{k,\al}$ and $\th+\bar\th\in C^{k+1,\al}$ with
$k\geq 5$ and $\al\in(0,1)$, these classical regularity bounds are used only
\emph{qualitatively}. They justify the Taylor expansions, singular-integral
representations, 
and vanishing orders of perturbation used in the weighted energy estimates; see Section~\ref{sec:reg_van}.

In contrast, the \emph{quantitative} stability estimates use only weighted $L^\infty$ and
$C^{1/2}$ norms of the perturbation  $\om$. 
Assume $\om \vp \in L^{\infty}$ for some weight $\vp$.
Since $\uu =  \na^{\perp}(-\D)^{-1}\om$ is log-Lipschitz, we can approximate $\uu(x)$ in $C^{\b} \cap L^{\inf}$, for any $\b < 1$, by interpolating its values at sufficiently many points. The $C^{\b} \cap L^{\inf}$ norm of the approximation error is bounded by $c || \om \vp ||_{\inf}$ with a small constant $c$. 
Similarly, for $\na \uu = \na \na^{\perp}(-\D)^{-1}\om$, 
the nonsingular part $K(z) \one_{|z| \geq \e} \ast \om $ of $\na \uu$ is Lipschitz and can be approximated in $C^{1/2} \cap L^{\inf}$. 
To control the singular part of $K(z) \one_{|z|\leq \e} \ast \om$ of $\na \uu$, we need H\"older regularity of $\om$. This motivates the form of estimates in \eqref{eq:vel_appr_err}.

\subsubsection{Approximation near $0$}\label{sec:appr_vel_near0}

Throughout Section \ref{sec:appr_vel_near0}, $\om$ denotes a \emph{generic input} to the Biot--Savart law, not 
the full perturbation $\om$ in \eqref{eq:lin} or \eqref{eq:bous_lin11}. Assume that $\om$ has vanishing order $O(|\xx|^{2+\alpha})$ for some $\alpha>0$ near $\xx=0$ and $\om$ is odd in $x$. We will apply the following derivation to $\om = \om_1$, the first component $W_{1,1}$ of $W_1$ in \eqref{eq:bous_decoup2:W1}, and to the error  $\om= \hat{\e}_1,\bar \e_1$  in \eqref{eq:uerr_dec1}.

To obtain sharper estimates near $\xx=0$, we extract the leading Taylor expansion of the velocity 
generated by $\om$, 
which is determined by two functionals $u_x(\om)(0)$ and $\cK_{00}(\om)$. After subtracting the corresponding terms, the velocity remainder has improved vanishing near the origin. We then localize this expansion and use it as the near-origin part of the finite-rank approximation.

Recall that $ -\D\phi =\om $, $\phi(x, 0) \equiv 0$, and $\phi, \om$ are odd in $x$. 
For $\phi \in C^{4,\al}$ locally and $\na \om(0) = 0$, we obtain
\bseq\label{eq:stream_vanish}
\beq\label{eq:stream_vanish:a}
\pa_x^i \pa_y^{2j} \phi(0) = 0, \, i + 2j \leq 4, \quad 0 = \om_x(0)  =  - \phi_{xxx}(0) - \phi_{xyy}(0) , 
 \quad  \om_{xy}(0) = -\phi_{xxxy}(0) - \phi_{xyyy}(0),
\eeq
For $\om$ further satisfying $\na^2 \om(0)=0$, we obtain 
$\phi_{xxxy}(0) + \phi_{xyyy}(0) =0$. Thus, using Taylor expansion, we expand the stream function $\phi$ as  
\beq
\phi(x, y) = \phi_{xy}(0) xy + \f{1}{6} ( \phi_{xxxy}(0) x^3 y + \phi_{x yyy}(0) x y^3 ) + h.o.t. 
= \phi_{xy}(0) xy +  \f{1}{6}  \phi_{xxxy}(0) ( x^3 y - x y^3)  + h.o.t. 
\eeq
\eseq

We can represent $\phi_{xxxy}(0)$ as an integral of $\om$
\beq\label{eq:u_appr_near0_K}
\phi_{xxxy}(0 ) = \f{2}{\pi} \int_{\R^2_{++}}  K_{00}(y) \om(y)  dy, \quad 
K_{00}(y) \teq  \f{ 24 y_1 y_2 (y_1^2 - y_2^2) }{ |y|^8} ,
\quad 
\cK_{00}(\om) \teq \f{1}{\pi} \int_{\R^2_{++}} K_{00}(y) \om(y) dy.
   \eeq
For $\om = O(|x|^{2+\al})$ with a suitable decay, the above integral is well-defined. By definition, we have $\phi_{xxxy}(0) = 2 \cK_{00}(\om)$.  When no confusion arises, we suppress the dependence of $u_x(0)$,
$\cK_{00}$, and the approximation terms on the input $\om$.

Recall $\uu = \na^{\perp} \phi = (-\pa_y, \pa_x) \phi$.  Note that $u_x(0) = - \phi_{xy}(0)$. Using the above formulas, near $0$, 
we obtain the leading order expansion  for $\na \uu$ and $\uu$ 
\beq\label{eq:u_appr_TL}
\bal
 u & = -\phi_y =  u_x(0) x - ( \tf{1}{3} x^3 -  x y^2 ) \cK_{00} + h.o.t., & &  \\
 v &= \phi_x = -u_x(0) y + ( x^2 y - \tf{1}{3} y^3 ) \cK_{00} + h.o.t., & & \\
u_x & = u_x(0) - (x^2 - y^2) \cK_{00} + h.o.t., 
 & \quad v_x  &=   2 x y  \cK_{00} + h.o.t. ,  \\ 
u_y &= - \phi_{yy} = 2 xy \cK_{00} + h.o.t. , & \quad v_y & = - u_x ,
 \eal
 \eeq
where $h.o.t.$ denotes the terms with higher vanishing power near $\xx= 0 $ 
and are estimated below.

We now introduce the spatial coefficients multiplying $u_x(0)$ and $\cK_{00}$ in \eqref{eq:u_appr_TL} by
\beq\label{eq:u_appr_near0_coe}
\bal
C_{u0} &= x , \quad C_{v0} = -y, \quad C_{u_x 0} = 1, \quad C_{u_y 0} = C_{v_x 0} = 0,  \\
C_u & = - ( \f{1}{3} x^3 - x  y^2 ), \quad C_v = x^2 y - \f{1}{3} y^3,  \\
C_{u_x} & = -(x^2 - y^2), \quad C_{v_x} = 2x y,  \quad C_{ u_y} =  2 xy, 
\eal
\eeq
and rewrite the above leading order formulas as 
\bseq\label{eq:u_appr_near0}
\beq
  f(x, y)= u_x(0) C_{f0}(x, y)  + \cK_{00} C_f(x, y) + h.o.t., \quad  f = u, v, u_x, v_x, u_y. 
\eeq

The above expansion identifies the leading-order terms of $\uu$ and $\na \uu$
near the origin in terms of $u_x(\om)(0)$ and $\cK_{00}(\om)$.  For $v_y$, we use $v_y = -u_x$ to estimate it.
We localize the above leading order terms to construct the approximation term near $0$ in the next subsection.

We record the vanishing  of the higher-order terms (h.o.t.) in the expansion above. 
For $g $ odd in $x$, using the Biot-Savart law  $\phi(g) = (-\D)^{-1} g$, standard $C^{\al}$ estimate, 
and \eqref{eq:u_est_direct}, we obtain 
\footnote{
When $k=0$, the  $\dot C^{\al}$ estimates of $\na^2 \phi(g)$  follow standard estimates, see e.g., 
\cite[Lemma 4.6]{majda2002vorticity}. Then use the equation $- \D \pa_x^k \phi = \pa_x^k \om$ and $-\D \phi = \om$ to obtain homogeneous $C^{\al}$ estimate for $\na^{k+2} \phi, k \geq 0$. 
 }
  \beq\label{eq:stream_local_reg1}
   \|  \na^{k+2} \phi(g) \|_{\dot C^{\al}(\R^2_+)} \les_{\al} \| \na^k g \|_{\dot C^{\al}} 
   , \quad 
    | \na^2 \phi(g)(\xx) | = | \na \uu(g)(\xx) | 
    \les_{\al, \ell} \| g \|_{C^{\al}} +  \| \, |\xx|^{\ell} g \|_{L^{\infty}}.
\eeq

For $\om \in C^{2,\al}$ 
odd in $x$ with $\na^i \om(0) =0$ for $i\leq 2$, we have $\phi = \phi(\om) \in C^{4,\al}$ locally. Using a Taylor expansion \eqref{eq:stream_vanish}, and the above estimates with $g = \om$, 
we obtain 
  \beq\label{eq:stream_local_reg}
    | \phi -   \phi_{xy}(0)  xy - \tf{1}{6} \phi_{xxxy}(0)(x^3 y - x y^3) | \les 
    |\xx|^{4 + \al} || \na^4 \phi||_{\dot C^{\al}} 
    \les  |\xx|^{4 + \al} || \om   ||_{C^{2,\al}}
    .
  \eeq
  Taking derivatives, we obtain the corresponding vanishing orders for the h.o.t. terms: 
\beq\label{eq:stream_local_reg2}
  | f(\om) - ( u_x(\om)(0) \cdot C_{f0}(\xx)  + \cK_{00}(\om) C_f(\xx) ) | 
  \les |\xx|^{k+\al}  \| \om \|_{C^{2,\al}} ,
\eeq
\eseq
where $k = 3$ for $ f = u, v$ and $k =2$ for  $f = u_x, v_x, u_y$.  This is the quantitative reason for subtracting the $u_x(\om)(0)$ and $\cK_{00}(\om)$
modes from $\uu, \na \uu$: the remainder has higher vanishing order and is therefore compatible
with the singular weight used  in Section \ref{sec:EE}.

\subsubsection{Approximation along the boundary}\label{sec:u_appr_1st}

We next construct the near-field approximation along the boundary, combining
the Taylor expansion near $\xx=0$ from Section~\ref{sec:appr_vel_near0} with
interpolation at finitely many boundary points $(x_i,0)$. 
We construct these boundary approximation terms \emph{only} for $f=u,u_x$; they are not needed for $v,u_y,v_x$.

To localize the approximation and interpolate smoothly along
the boundary, we introduce some cutoff functions.  Let $\chi$ be the cutoff function  in \eqref{eq:cutoff_exp} 
and $\td \chi = 1 -\chi$. They satisfy
\bseq\label{eq:vel_chi_bd}
\beq\label{eq:vel_chi}
\td \chi(x) \teq 1 -\chi(x) , \quad  \chi(x) = 0,  \  \td \chi(x) = 1, \quad x \leq 0, \quad \chi(x) = 1, \ \td \chi(x) = 0, \quad x \geq 1 .
\eeq

Given $0< x_0 < x_1 < ...< x_n < x_{n+1}$ and $y_0>0$, we construct the cutoff functions 
\beq\label{eq:vel_cutoff_bd}
\bal
\chi_0(x, y)  &=  \td \chi( \f{x - x_0}{ x_1- x_0} ) \td \chi( \f{y-y_0}{y_0} ),   \\
\chi_i(x, y)  & =  \B( \chi( \f{x - x_{i-1}}{x_i - x_{i-1}} ) \one_{x \leq x_i}  + \td \chi( \f{ x - x_i}{ x_{i+1} - x_i} ) \one_{x \geq x_i} \B) \td \chi( \f{y-y_0}{y_0} ), \quad 1\leq i \leq n,
\eal 
\eeq
and extend them as even functions in $x$.

By construction,  $\{\chi_i\}_{0\leq i\leq n}$ form a partition
of unity in the $x$ direction near the boundary:
\beq\label{eq:vel_cutoff_bd:c}
\sum\nolimits_{0\leq i\leq n}\chi_i(x,y)=1,
 \   \mbox{on} \, [-x_n, x_n] \times [0,  y_0],
\quad \supp(\chi_i) \subset [-x_{i+1}, x_{i+1}] \times [0, 2 y_0].
\eeq
\eseq
The factor $\td\chi((y-y_0)/y_0)$ localizes these cutoff functions to the boundary
layer $y\leq 2y_0$.

We approximate $\uu,\na\uu$  along the boundary using  the following modification of \eqref{eq:u_appr_idea}:
\beq\label{eq:u_appr_1st}
\bal
\hat f_1(\om) (x, y) &\teq   C_{f0}(x, y) u_x(\om)(0) + \hat f_{10}(\om)(x, y), \\
 \hat f_{10}(\om)(x, y) &\teq  C_f(x, y)  \B( \cK_{00}(\om) \chi_0(x, y) + \sum\nolimits_{ 1 \leq i \leq n} \f{ \fNS(\om)(x_i, 0)  }{C_f( x_i, 0)} \chi_i(x, y)  \B) , 
\eal 
\eeq

Here, $\fNS(\om)$ is the non-singular part of the integral related to the velocity
\beq\label{eq:u_appr_1st2}
\bal
\fNS(\om)(x_i, 0) &= \int_{y \in \R^2, \max(|y_1 - x_i|, |y_2|) \geq t_i } K_f( (x_i , 0) - y)  
 \om_{\mf{ext}}(y)  dy 
-  C_{f0}(x_i, 0) u_x(\om)(0), 
\eal
\eeq
where $  \om_{\mf{ext}}$ denotes the odd extension of $\om$ in $y$ from $\R^2_+$ to $\R^2$.
Here $K_f$ is the kernel for $f=u,v,u_x,v_x,u_y$, $\mathsf{NS}$ is short for
\emph{nonsingular}, and $C_{f0},C_f$ are the spatial coefficients in
\eqref{eq:u_appr_near0}; see also \eqref{eq:u_appr_near0_coe}. For each fixed
$x_i$, $\fNS(\om)(x_i,0)$ is a scalar linear functional of $\om$; hence each
summand in \eqref{eq:u_appr_1st} is rank-one in $\om$. Thus $\hat f_1(\om)$ is
a finite-rank operator of $\om$. 
The subscript $1$ indicates that $\hat f_1$ is the first part of the
approximation.

We use the full boundary approximation above \emph{only} for $f=u,u_x$. For
$f=v,v_x,u_y$, the associated coefficients in the nonlocal terms are
relatively small, e.g. $\bar\omega_y$ in $v\bar\omega_y$ and $\bar\theta_y$
in $v_x\bar\theta_y$ in \eqref{eq:lin}. Moreover, the nonlocal term $u_y$
\emph{only} appears in the $\xi$-equation \eqref{eq:lin}, where the stability estimates are stronger.
Thus, for $f=v,v_x,u_y$, we keep only the near-origin Taylor approximation
from Section~\ref{sec:appr_vel_near0}, without adding the boundary
approximation terms $\fNS(\om)(x_i,0)$ \eqref{eq:u_appr_1st2}:
\beq\label{eq:u_appr_1st_other}
\hat f_1(\om)(x,y)
=
C_{f0}(x,y)u_x(\om)(0)+C_f(x,y)\cK_{00}(\om)\chi_0(x,y),
\qquad f=v,v_x,u_y .
\eeq

For $f=u,u_x$, $\hat f_1$ is given by \eqref{eq:u_appr_1st}. By
definition, for $f = u, u_x$, we have the following identities for the approximation error 
at $(x, y) = (x_i, 0)$
\bseq\label{eq:u_appr_1st_xi}
\beq\label{eq:u_appr_1st_xi:a} 
f(\om)(x_i,0)-\hat f_1(\om)(x_i,0)
= f(\om)(x_i,0)- ( C_{f0}(x_i,0)u_x(\om)(0) +  \fNS(\om)(x_i,0) ) , \  1 \leq i\leq n,
\eeq
whereas for \emph{all five} components $f=u,v,u_x,v_x,u_y$ and $\xx = (x, y)$, near $\xx=0$, we have
\beq\label{eq:u_appr_1st_xi:b}
f(\om)(\xx)-\hat f_1(\om)(\xx)
=
f(\om)(\xx)-C_{f0}(\xx)u_x(\om)(0)-C_f(\xx)\cK_{00}(\om),
\quad \xx \in [0, x_0] \times [0, y_0].
\eeq
\eseq

Thus \(\widehat f_1\) has two roles.  Near the origin, 
by  \eqref{eq:u_appr_1st_xi:b}, it matches the leading Taylor expansion of $f$ 
as in \eqref{eq:u_appr_near0} for \emph{all five cases}.  Along the boundary, for \(f=u,u_x\), 
by \eqref{eq:u_appr_1st_xi:a} and \eqref{eq:u_appr_1st2}, it interpolates
the nonsingular part of the Biot--Savart integral at the selected points
\((x_i,0)\), which is similar to \eqref{eq:u_appr_idea}. Here, we consider a weighted version of \eqref{eq:u_appr_idea} with vanishing coefficient $C_f(x)$  \eqref{eq:u_appr_near0_coe} so that the approximation has the correct vanishing order near $x=0$.

\subsubsection{Approximation in the far-field}\label{sec:u_appr_far}

We next construct the far-field approximation. The idea is to
replace the global integral  in $u_x(\om)(0)$ by a truncated integral
determined by the main part of the symmetrized kernel in the far field. This
improves the estimate when $|\xx|$ is large.

Recall $u_x(\om)(0) = - \f{4}{\pi} \int_{\R^2_{++}} \f{y_1 y_2}{|y|^4} \om(y) dy$.  For a large radius $R_i$, we define the truncated integral 
\beq\label{eq:u_appr_In}
I_i(\om) \teq  - \f{4}{\pi}\int_{ y \in \R^2_{++}, \, \max(y_1, y_2) \geq R_i} \f{y_1 y_2}{|y|^4} \om(y) dy.
\eeq
Similar truncated integrals appear in the velocity approximation in \cite{kiselev2013small}. It arises as the leading order part of the symmetrized kernel: for $u_x = K \ast \om$, $\om$ odd in $x$, and $K(s) =-\f{1}{\pi} \f{s_1 s_2}{|s|^4}$ , we have
\beq\label{eq:u_appr_2nd_main}
\bal
K^{sym}(x, y) & \teq
K(x - y)  + K(x+y) - K( x_1 - y_1, x_2 + y_2)
- K(x_1 + y_1, x_2 - y_2) \\
&
= -\f{4}{\pi}\f{y_1 y_2}{|y|^4} + l.o.t.
\eal
\eeq
for $ \max |y_i| \geq 2 \max |x_i|$, where the lower order term (l.o.t.) decays faster in $y$ with $|l.o.t.| \les |\xx| \cdot |\yy|^{-3}$. For $f = u, v, v_y$, using a similar argument, we obtain the leading order part of the associated kernel $K_f$ 
\beq\label{eq:u_appr_2nd_main2}
K_f^{sym} = - C_{f0} \f{4}{\pi} \f{y_1 y_2}{|y|^4}  + l.o.t.
\eeq
for $ \max |y_i| \geq 2 \max |x_i|$, where $C_{f0}$ is defined in \eqref{eq:u_appr_near0_coe}. When $|y| / |x|$ is small, the function $-\f{4}{\pi}\f{y_1 y_2}{|y|^4}$ does not approximate $K^{sym}(x, y)$ well, and thus we truncate it
\beq\label{eq:u_appr_2nd_main3}
\int_{ \max |y_i| \geq 2 \max |x_i|  } - \f{y_1 y_2}{|y|^4}  \om(y) dy.
\eeq

The truncation in \eqref{eq:u_appr_2nd_main3} depends on the evaluation point
$x$, and hence the above operator is not finite-rank. We approximate \eqref{eq:u_appr_2nd_main3} by finite-rank operators by replacing the $x$-dependent cutoff with finitely many fixed truncations, and then interpolating between them using a partition of unity in the \(x\)-variable.
Concretely, we approximate $\one_{\max |y_i| \geq 2 \max |x_i|}$ by
\[
g(x, y) =\tts{\sum}_i \one_{  |y_1| \vee |y_2| \geq R_i} \chi^R_i(x)
\]
where $  \{ \chi^R_i \}_{ i\geq 1}$ is a partition of unity on the approximation domain  and $\chi^R_i(x)$ is localized to the region where $|x|$ is comparable to $R_i$. Thus, for $x$ close to $R_i$, we obtain 
\[
g(x, y) \approx \one_{  |y_1| \vee |y_2| \geq R_i}  \approx \one_{ |y_1| \vee |y_2| \geq |x|}.
\]

We now make this approximation precise.  Recall the cutoff functions  $\chi, \td \chi$ from \eqref{eq:vel_chi}. Given $R_0 < R_1 < ..< R_m $, we define  
\bseq \label{eq:vel_cutoff_R}
\begin{align}
\chi^R_i(x, y) & = \td \chi( \f{x - R_{i}}{ R_{i+1} - R_i } )
 \td \chi( \f{y - R_{i}}{ R_{i+1} - R_i } )
 - \td \chi( \f{x - R_{i-1}}{ R_{i} - R_{i-1} } )
 \td \chi( \f{y - R_{i-1}}{ R_{i} - R_{i-1} } ), \ & & 1 \leq i \leq  m - 1 , \notag \\
 \chi^R_m(x, y) & = 1 - \td \chi( \f{x - R_{m-1}}{ R_{m} - R_{m-1} } )
 \td \chi( \f{y - R_{m-1}}{ R_{m} - R_{m-1} } ) , & & \\
   S_i^R & =  1  - \td \chi( \f{x - R_{i-1}}{ R_{i} - R_{i-1} } )
 \td \chi( \f{y - R_{i-1}}{ R_{i} - R_{i-1} } ),  
 & & 1 \leq i\leq m. \notag
\end{align}
By construction, the cutoff functions $\{\chi_i^R\}_{1\leq i\leq m}$ form a
partition of unity in the far-field region:
\beq
\tts{\sum}_{1\leq i\leq m}\chi_i^R(x,y)=1,
 \qquad \forall \, \max(x,y)\geq R_1.
\eeq

Moreover, these functions satisfy
\beq\label{eq:vel_cutoff_R:c}
\bal
& \supp (\chi_i^R )  \subset [0,R_{i+1}]^2 \setminus [0,R_{i-1}]^2,
\quad 1\leq i\leq m-1, & \quad
& \chi_m^R=1 \quad \text{for } \max(x,y)\geq R_m, \\
& \supp (S_i^R) \cap \R^2_{++}  \subset \R^2_{++}\setminus [0,R_{i-1}]^2,
& \quad & S_i^R=1 \quad \text{for } \max(x,y)\geq R_i.
\eal
\eeq
\eseq

Recall the functional $I_i(\om)$ from \eqref{eq:u_appr_In}. We now define the second approximation by
\beq\label{eq:u_appr_2nd}
\bal
\hat f_2(\om)(x, y) & =C_{f0}(x, y) (1 -  \chi_{\tot}(x, y)) \cdot \big( \tts{\sum}_{ 1\leq i \leq m } \chi^R_i(x, y)  ( I_i(\om) - u_x(\om)(0)) \big), \\
\chi_{\tot}(x, y) &= \tts{\sum}_{ 0 \leq i \leq n}  \chi_i(x, y),
\eal
\eeq
where $\chi_{\tot}$ denotes the sum of the cutoff functions  \eqref{eq:vel_cutoff_bd:c} for the \emph{first} approximation in Section \ref{sec:u_appr_1st}. 
By definition \eqref{eq:vel_chi_bd}, we have 
\[
 \chi_{\tot}(x, y) = 1 
  \quad \mbox{on} \   [0, x_n] \times   [ 0,  y_0],
 \qquad  \chi_{\tot}(x, y)= \chi_i(x, y) =  0   \quad \mbox{for} \  x \geq x_{n+1}, \, \mbox{or}  \  y \geq 2 y_0. 
\]
Thus $1-\chi_{\tot}$ turns on the second approximation only away from the near-origin and boundary approximation region $[0, x_n] \times   [ 0,  y_0]$. 

For $x = R_i$ and $y \in [0, R_i]$  satisfying $\chi_{\tot}(x, y) = 0$, from \eqref{eq:u_appr_1st2}, we get 
\[
\hat f_1(\om)(x, y) = C_{f0}(x, y) \cdot u_x(\om)(0), \quad \hat f_2(\om)(x, y) = C_{f0}(x, y) \cdot  (I_i(\om) - u_x(\om)(0)), 
\]
and therefore
\[
\bal
f - \hat f_1(\om) - \hat f_2(\om)
&=  f - C_{f0} I_i(\om)
= \int  K_f(x- y)  \om(y ) dy - C_{f0} I_i(\om) \\
&= \int_{\R^2_{++}}( K_f^{sym}(x, y) + \f{4}{\pi} C_{f0}(x) \f{y_1 y_2}{|y|^4} \one_{\max(|y_1|, |y_2|) \geq R_i} ) \om(y) dy,
\eal
 \]
 where $K_f^{sym}$ is the symmetrized kernel \eqref{eq:u_appr_2nd_main}. Thus, for large $|x|$, the approximation $\hat f_1 + \hat f_2$ can be seen as a smooth finite-rank interpolation of the far-field main term in \eqref{eq:u_appr_2nd_main2}, \eqref{eq:u_appr_2nd_main3}.

\begin{remark}[\bf Approximate part of the integral]
The lower-order terms in \eqref{eq:u_appr_2nd_main} and
\eqref{eq:u_appr_2nd_main2} are small only in the far-field regime
$|y|\gg |x|$. Thus, the finite-rank approximation \eqref{eq:u_appr_2nd} approximates only the far-field
main part of the integral defining $f$. 
Nevertheless, it gives a better estimate for
$f-\hat f_1-\hat f_2$ than for $f$ itself, which is sufficient for our
purpose.
\end{remark}

Combining \eqref{eq:u_appr_1st} and \eqref{eq:u_appr_2nd}, we define the full approximation of
$f=\uu,\nabla \uu$, and the corresponding approximation of
$\td f = f - C_{f0} u_x(0)$ by 
\beq\label{eq:u_appr}
 \hat f(\om) \teq \hat f_1(\om) + \hat f_2(\om),
 \quad \appo{ f}(\om) \teq \hat f(\om) - C_{f0} u_x(\om)(0) = \hat f_{10}(\om) + \hat f_2(\om),
\eeq
where the notations $\td u, \td v$ are introduced in \eqref{eq:u_tilde}. Clearly, $\hat f , \appo{f}$ are finite-rank operators of $\om$.
Using the operators $\appo{f}$ for $f =u, v, \uu_x, \uu_y$, 
we construct the finite-rank operator $\cK_{2i}$ in \eqref{eq:appr_vel}.

Since $\chi_i$ in \eqref{eq:u_appr_1st} has a compact support by \eqref{eq:vel_cutoff_bd:c}, 
we obtain $|\chi_i C_{f}| \les 1$ 
for $C_{\bullet}$ in \eqref{eq:u_appr_near0_coe}. Using $|\chi_i^R|\leq 1, |\chi_{\tot}| \les 1$, and the definitions of $\hat f_{10}(\om) $ 
in \eqref{eq:u_appr_1st} and $\hat f_2(\om)$ in \eqref{eq:u_appr_2nd}, we estimate 
\beq\label{eq:u_appr_growth}
\bal
|\appo{f}(\om)(\xx) | & = | (\hat f_{10}(\om)+ \hat f_2(\om)) (\xx)| \\
& \les |\cK_{00}(\om)| + |C_{f0}(\xx)|  (   \tts{\sum}_{i\leq n} |\fNS(\om)(x_i, 0)| 
+ \tts{\sum}_{i \leq m} |I_i(\om) - u_x(\om)(0)| ).
\eal
\eeq

\begin{remark}[\bf Different approximation centers and parameters]

We use \emph{different number} of boundary approximation terms, i.e.  $n$ in \eqref{eq:u_appr_1st},
and different approximation centers $x_i$ and size $t_i$ of the singular regions  in 
\eqref{eq:u_appr_1st}, \eqref{eq:u_appr_1st2} for $\hat u(\om)$ and $\wh{u_x} (\om)$. 
Hence, the family of boundary cutoff functions 
in \eqref{eq:vel_chi_bd} are also different for $\hat u$ and $\wh {u_x}$. 
Moreover, in the approximation \eqref{eq:u_appr_1st_other} for $u_y,v_x, v$, parameter $x_0$ in \eqref{eq:vel_cutoff_bd} for cutoff function $\chi_0$ in \eqref{eq:u_appr_1st_other} 
\emph{differs from} that for $u, u_x$. 

On the other hand, for the approximation in the far-field, we choose 
the \emph{same number} of approximation terms and \emph{the same radii} $R_i$ in \eqref{eq:u_appr_In}
for \emph{all five cases} $u, v, u_x, u_y, v_x$. 

We collect these parameters $x_i, t_i, R_i$ in  Appendix~\ref{app:appr_vel_para}.

To obtain sharp estimate of the constants in \eqref{eq:vel_appr_err} and in the energy estimates 
for $\uu - \hat \uu, \na \uu - \widehat{ \na \uu }$, which are important for reducing the number of approximation terms and obtaining sharp energy estimates, we estimate the integrals with computer assistance. See Section 4 of Part II \cite{ChenHou2023b} for details. 
In total, fewer than $50$ approximation terms are used.

\end{remark}

\begin{example*}
We give two concrete examples of the finite-rank approximations. From \eqref{eq:u_appr} and \eqref{eq:u_appr_near0_coe}, the approximations for $\wh{u_x}$ and $\hat v$ are
\[
\bal
  \wh{u_x}(\om) & = u_x(\om)(0)    - (x^2 - y^2) (\cK_{00}(\om) \chi_{u_x, 0} - \tts{\sum}_{1\leq i\leq n_{u_x}} \f{\wh{u_x}_{, \mf{NS} }(\om)(x_i,0) }{x_i^2}  \chi_{u_x, i}(x, y)) \\
  & \quad + (1 - \chi_{tot}) \cdot ( \tts{\sum}_{1\leq  i \leq m_{u_x} } \chi_i^R ( I_i(\om) - u_x(\om)(0) ) \, ), 
 \qquad \quad \chi_{tot} = \sum_{ 0 \leq i \leq n_{u_x}}  \chi_{u_x, i} , \\
 \hat v(\om)  & = -u_x(\om)(0) y + (x^2 y - \tf13 y^3  ) \cK_{00}(\om) 
 \chi_{v, 0}
 - y ( 1-  \chi_{v, 0})  \cdot ( \tts{\sum}_{1\leq  i \leq m_{u_x} } \chi_i^R ( I_i(\om) - u_x(\om)(0) ) \, ) ,
\eal
\]
where $\cK_{00}(\om),  \wh{u_x}_{, \mf{NS} }(\om), I_i(\om)$ are the functionals defined in \eqref{eq:u_appr_near0_K}, \eqref{eq:u_appr_1st2},  and \eqref{eq:u_appr_In}.

For $\max( x, y) \geq R_m$, since $R_m, m$ chosen in \eqref{para:appr_2nd} 
and $x_{\bullet}$ in Appendix \ref{app:appr_vel_para} satisfies $R_m > x_{n+1}$, the support properties of $\chi_i$ in \eqref{eq:vel_cutoff_bd:c} 
and $\chi_i^R$ in \eqref{eq:vel_cutoff_R:c} imply $\chi_i = 0, \chi_{\mf{tot}} = 0, 
\chi_i^R =0, i \leq m-1, \chi_m^R = 1$. Thus, 
the approximations $\hat f_1$ \eqref{eq:u_appr_1st}, $\hat f_2$ \eqref{eq:u_appr_2nd}, 
$\hat f, \appo{f}$ \eqref{eq:u_appr}  reduce to 
\begin{align}\label{eq:eg_appr_far}
 \hat f_1(\om) & = C_{f0}  u_x(\om)(0),
  & \   \hat f_2(\om) &= C_{f0} ( I_m(\om)- u_x(\om)(0) ), \\
  \hat f(\om ) & = \hat f_1(\om) + \hat f_2(\om) =  C_{f0}  I_m(\om) ,
 &  \appo{f}(\om) & = C_{f0} ( I_m(\om)- u_x(\om)(0) ),
 \  f \in \{ u, v, u_x, u_y, v_x \}, \notag 
\end{align}
where $R_m = 2^4 \cdot 13 $,
and $I_m(\om) 
= - \f{4}{\pi}\int_{ \max(y_1, y_2) \geq R_m} \f{y_1 y_2}{|y|^4} \om(y) dy $ by \eqref{eq:u_appr_In}.

\end{example*}

\subsubsection{Reformulation of the approximations}

The formulas  \eqref{eq:u_appr_1st} and\eqref{eq:u_appr_2nd} are useful for estimating the approximation error, because
only a few local cutoff functions are active at each point.  For constructing
the approximate space-time solutions $\hat F_i, \hat F_{\chi,i}$, however, it is better to rewrite the same
finite-rank operators \eqref{eq:u_appr_1st}, \eqref{eq:u_appr_2nd} in a telescoping form.  The telescoping form exposes cancellations between neighboring scalar functionals and uses smoother spatial
coefficients.

Recall $\chi_i$ defined in \eqref{eq:vel_cutoff_bd}. We introduce 
\[
S_j = \tts{\sum}_{ 0\leq i \leq j} \, \chi_i , 
\quad \mbox{for} \  0 \leq j \leq n ,
\quad \chi_j = S_j - S_{j-1} , \, \mbox{for} \  1 \leq j \leq n .
\]
 By definition, we have $S_0 = \chi_0$ and 
\[
S_j(x, y) = 1 ,\quad  \forall \, (x, y) \in  [0,  x_{j}] \times [0, y_0] , \qquad  S_j(x, y) = 0, \ \forall \, x \geq x_{j+1}.
\]

By rewriting the summation, we can rewrite the operator $\hat f_{10}(\om)$ in \eqref{eq:u_appr_1st} as follows 
\begin{align}
&  \hat f_{10}(\om)  =  C_f   \B( \cK_{00}(\om) S_0  + \sum_{ 1 \leq i \leq n} \f{  \fNS(\om)(x_i, 0) }{C_f( x_i, 0)}  ( S_i - S_{i-1} ) \B)    \label{eq:u_appr_1st_sum}   \\
& = C_f(x, y) \B( S_0 ( \cK_{00}(\om) - \f{ \fNS(\om)(x_1, 0) }{ C_f( x_1, 0)  }  )   
+ S_n \f{ \fNS(\om)(x_n, 0) }{ C_f( x_n, 0)  }      +
\sum_{1 \leq i \leq n -1 } S_i \B( \f{ \fNS(\om)(x_i, 0) }{ C_f( x_i, 0)  } - \f{ \fNS(\om)(x_{i+1}, 0) }{ C_f( x_{i+1}, 0)  }  \B)  \B) . \notag
\end{align}

For notational simplicity, we suppress the dependence of $C_f$ and $S_i$ on
$(x,y)$. The second line of \eqref{eq:u_appr_1st_sum} rewrites the coefficients
in terms of differences of neighboring scalar functionals. This makes explicit
the cancellation among $\cK_{00}(\om)$ and
$\fNS(\om)(x_i,0)/C_f(x_i,0)$ with $1\leq i\leq n$.

Similarly, using \eqref{eq:vel_cutoff_R} and
\[
\chi_i^R=S_i^R-S_{i+1}^R,\qquad 1\leq i\leq m-1,
\qquad S_m^R=\chi_m^R ,
\]
we rewrite the second  approximation \eqref{eq:u_appr_2nd} as
\begin{align}\label{eq:u_appr_2nd_sum}
\hat f_2(\om) & = C_{f0}  (1 -  \chi_{tot} )  \B( \sum\nolimits_{ 1\leq i \leq m -1} (S_i^R - S_{i+1}^R)  ( I_i(\om) - u_x(\om)(0))  + S_m^R (I_m(\om) - u_x(\om)(0)) \B) \notag \\
& = C_{f0} (1 - \chi_{tot}) \B( S_1^R ( I_1(\om) - u_x(\om)(0) ) 
+ \sum\nolimits_{ 2\leq i \leq m} S_i^R (I_i(\om) - I_{i-1}(\om) )  \B), 
\end{align}
For notational simplicity, we suppress the dependence of $\hat f_2$,
$\chi_{\tot}$, and $S_i^R$ on $(x,y)$. The second line of
\eqref{eq:u_appr_2nd_sum} rewrites $\hat f_2$ in terms of the increments
$I_i(\om)-I_{i-1}(\om)$, which satisfy sharper estimates than the individual
quantities $I_i(\om)-u_x(\om)(0)$.

We use the two representations of $\hat f_{10}$ and $\hat f_2$ for different
purposes. In the estimates of the velocity error $f-\hat f_1-\hat f_2$ in
Section~4 of Part~II \cite{ChenHou2023b}, we use the original formulas
\eqref{eq:u_appr_1st} for $\hat f_{10}$ and \eqref{eq:u_appr_2nd} for
$\hat f_2$. Indeed, by the local support of the cutoffs $\chi_i$ in
\eqref{eq:vel_chi_bd} and $\chi_i^R$ in \eqref{eq:vel_cutoff_R}, for each
fixed $\xx$ at most two neighboring summands are nonzero.

On the other hand, in the numerical construction of $\wh W_2$ in
\eqref{eq:bous_W2_appr}, we use the reformulated formulas
\eqref{eq:u_appr_1st_sum} and \eqref{eq:u_appr_2nd_sum}. In these formulas,
the spatial coefficients are $C_f(x,y)S_i$ and
$C_{f0}(1-\chi_{\tot})S_i^R$, rather than $C_f(x,y)\chi_i$ and
$C_{f0}(1-\chi_{\tot})\chi_i^R$. The cutoff factors $S_i,S_i^R$ are
numerically more regular than $\chi_i,\chi_i^R$: for example, as $x$
increases, $S_i$ has one transition from $1$ to $0$, whereas $\chi_i$ has two
transitions, from $0$ to $1$ and then from $1$ to $0$. Thus the corresponding
spatial coefficients can be approximated by the numerical basis with smaller
approximation error. The reformulated formulas also exploit the cancellations
above, leading to smaller residual errors, i.e. $\cR(\cdot)$ in \eqref{eq:bous_err_op}, in solving the linearized equations.

In the nonlinear stability estimates, we estimate $\hat f_1$ using both
representations and take the sharper bound. For $\hat f_2$, we use
\eqref{eq:u_appr_2nd_sum}.

\subsection{Qualitative regularity and vanishing order of perturbation}\label{sec:reg_van}

In Section \ref{sec:finite_pertb}, we decompose the full linearized equation \eqref{eq:lin} into equations for $W_1$ and $\hat W_2$ in \eqref{eq:bous_decoup2}. Here we establish the \emph{qualitative} regularity and vanishing properties  for $W_1, \hat W_2$, which are needed for the later weighted energy estimates 
and justify the finiteness of the singularly weighted norms and functionals in Section~\ref{sec:EE}. 
The constants below are irrelevant to the later quantitative stability estimates. Readers uninterested in qualitative regularity may proceed directly to Section~\ref{sec:EE}.

Although $W_1$ and $\hat W_2$ are coupled in \eqref{eq:bous_decoup2}, 
we view  $W_1$-equation as a linear equation for $W_1$ with forcing depending on
  the full perturbation $W=W_1+\hat W_2$ and residual  terms; see \eqref{eq:reg_W1_pf1}.  Below, we first discuss the qualitative regularity of the profile and $W$ in Section \ref{sec:reg_van_pf}, then the regularity of the finite-rank operators and residual terms in Section \ref{sec:reg_finite}, and finally the regularity and cubic vanishing of $W_1$ and related terms in Sections \ref{sec:reg_W12}, \ref{sec:reg_other}.

\subsubsection{Regularity of the profile and full perturbation}\label{sec:reg_van_pf}

We decompose the full solution $(\om_{\tot}, \th_{\tot})$ ($\tot$ stands for \textit{total}) to the self-similar Boussinesq equations \eqref{eq:bousdy1}-\eqref{eq:normal} as 
\beq\label{eq:whole_solu}
 \om_{\tot} = \bar \om +  \om, \quad \th_{\tot} = \bar \th +  \th ,
 \quad \uu_{\tot} = \bar \uu + \uu,  \quad \uu = (u, v),
\eeq
where $ (\om,  \th, \uu)$ is the \emph{full} perturbations
and $\bar \om, \bar \th, \bar \uu$ is the profile.
Consider initial data satisfying
\beq\label{eq:reg_init}
\om_{\tot, 0} \in C^{\kin, \ain}, \, \th_{\tot,0} \in C^{\kin+1, \ain},  
\quad  \pa_x\om_{\tot, 0}(0) \neq 0 , \quad  \om_{\tot,0}, \, \na \th_{\tot,0}\in C_c, 
\eeq
for some $\kin \geq 5, \ain \in (0, 1)$, and satisfying the symmetries \eqref{eq:Sym}.

 The local well-posedness theory
\footnote{
Using \eqref{eq:normal}, we can rewrite the self-similar Boussinesq equations \eqref{eq:bousdy1}-\eqref{eq:normal} as a closed system for $\om_{\tot}, \th_{\tot}$ involving $\om_{\tot,x}(0), \th_{\tot,xx}(0)$. For $( \om_{\tot,0}, \th_{\tot,0} )  \in \cX^k, \cX^k \teq C^{k,\al} \times C^{k+1, \al}, k \geq 1, \al > 0$ or $\cX^k \teq H^k \times H^{k + 1} , k \geq 3$ 
with $\pa_x \om_{\tot, 0}(0) \neq 0$, the local well-posedness theory (LWP) follows from the standard argument in \cite{chae1999local,chae1997local} with $ \int_0^t || \na (\uu_{\tot},\th_{\tot})||_{L^{\infty}}$ controlling the breakdown of $\cX^k$ regularity. In fact, compared to the original Boussinesq 
\eqref{eq:bous_sys}, the extra terms $c_l \xx \cdot \na ( \om_{\tot}, \th_{\tot})$ form a Lipschitz vector field, and $c_{\om} \om_{\tot}, c_{\th} \th_{\tot}$ are bounded in the same functional space $\cX^k$. These terms contribute to bounded terms in the energy estimates \cite{chae1999local,chae1997local}. Moreover, using \eqref{eq:normal1} and \eqref{eq:normal}, 
we can bound $c_l, c_{\om}, c_{\th}$ by $|| ( \om_{\tot} , \th_{\tot}) ||_{\cX^k}$ via embedding inequalities. Thus, we obtain the same type of energy estimates as in \cite{chae1999local,chae1997local} and further prove LWP. 
} 
implies that 
$\om_{\tot}(t) \in C^{\kin,\ain} $, $\th_{\tot}(t) \in C^{\kin+1,\ain}$, as long as
\[
M(t) \teq  \tts{\int_0^t }
        \|\na(\uu_{\tot},\th_{\tot})\|_{L^\infty}\,ds
\]
remains finite.  We next record the compact-support property. Denote $V\teq c_l\xx+\uu_{\tot}$. From \eqref{eq:bousdy1}, we have
\[
  \pa_t \om_{\tot} +  V \cdot \na \om_{\tot} = \th_{\tot, x} + c_{\om}\om_{\tot},
  \quad \pa_t \na \th_{\tot} + V \cdot \na (\na \th_{\tot}) = - (\na V)^T \na \th_{\tot} 
  + c_{\th} \na \th_{\tot}
\]

Since 
$\frac{|V(\xx)|}{1 + |\xx|} \les 1 + \| \na \uu_{\tot} \|_{L^{\infty}}$ 
and $\om_{\tot,0}, \na \th_{\tot,0}$ have compact support, estimating the growth of  characteristics obtains that $M(t)$ controls the growth of the supports of $\om_{\tot}, \na \th_{\tot}$ on $[0, t]$. 
Thus, $\om_{\tot}(t)$ and $\na \th_{\tot}(t)$ remain compactly supported as long as $M(t) < +\infty$. 
From \eqref{eq:whole_solu} and $\kin \geq 5$, 
 $\om, \th$ have the same regularity as $\bar \om, \bar \th \in C^{4,1}$, respectively. Moreover, we have 
\[
\om(t) = -\bar \om, \quad  \na \th(t) = -\na \bar \th , \quad  |\xx| > R, \quad
\mbox{for \ some \ } R =  R( M(t), \om_{\tot,0}, \th_{\tot,0} ) > 0 .
\]

Denote $\uu(\bar \om) = \na^{\perp}(-\Delta)^{-1} \bar \om$.  
We summarize the regularity of $\om, \th$ and of the profile $\bar \om, \bar \th$ as discussed in Section~\ref{sec:property} and \eqref{solu:decay}.
\bseq\label{eq:solu_reg1}
\begin{align}
  \om, \ \th,  \  \bar \om, \ \bar \th \in C^{4, 1} , 
  \quad \bar \phi^N \in C^{4, 1} & , & 
  \quad  \tf{1}{|\xx|} \uu(\bar \om)  \in L^{\infty}, 
  \  \na \uu(\bar \om) \in C^{4, \al} & , \quad \forall \, \al \in (0, 1),  \label{eq:solu_reg1:a} \\
    |\na^i (\om(t), \bar \om)(\xx)| \les C_t \la \xx \ra^{ \bar \al_1 - i } &, &
  \qquad   |\na^i (\th(t), \bar \th)(\xx)| \les C_t \la \xx \ra^{ 1 + 2 \bar \al_1 - i } & , 
\quad i \leq 3  ,
  \label{eq:solu_reg1:decay} 
\end{align}
where $C_t =  C( M(t), \om_{\tot,0}, \th_{\tot,0} )$, and we recall $ \bar \al_1 < 0$ 
from \eqref{solu:decay}. The above regularity and \eqref{solu:decay} implies that the residual error $\bar \cF_i$ in \eqref{eq:bous_err} and $ \bar f_{c_{\om}, i}$ in \eqref{eq:f_cw} satisfy 
\beq\label{eq:solu_reg1:b}
  \bar f_{c_{\om}, i} , \quad  \bar \cF_i \in C^{2, 1}  ,
  \quad | \bar f_{c_{\om}, i}| + | \bar \cF_i| \les \la \xx \ra ^{ \bar \al_1} , \quad \bar \al_1 < 0 .
\eeq
\eseq

Finally, we record the local expansion near the origin.  Recall the symmetries of $(\om_{\tot}, \th_{\tot})$ from \eqref{eq:Sym}. The profile $(\bar \om, \bar \th)$ has the same properties as $(\om_{\tot}, \th_{\tot})$, and hence so does $(\om, \th)$.  In addition to \eqref{eq:reg_init}, we consider initial perturbations with cubic vanishing:
\beq\label{eq:van_init_cubic}
 \om_0 = \om_{\tot, 0} - \bar \om = O(|x|^3), \quad  \na \th_0 = \na (\th_{\tot,0} - \bar \th) = O(|x|^3) ,
\eeq
which implies $\om_{ 0, x}(0) = \th_{0, xx}(0) = 0$. The symmetries and $\om_x(0) = 0, \th_{xx}(0) = 0$ for the perturbation obtained by linearizing \eqref{eq:normal1} yield the expansion near $\xx = (x, y) = 0$ 
\beq\label{eq:van_pertb_sq}
   \om =   \om_{xy}(0)  xy + O(|\xx|^3), \
  \    \th =   \th_{xxy}(0)\cdot \tf12 x^2 y + O(|\xx|^4),
 \quad 
  \na \om(0) = 0, \  \na^2 \th(0) = 0 , 
\eeq

Since the cubic vanishing \eqref{eq:van_init_cubic} of the \emph{full} perturbation 
$W= (\om, \na \th)$ is \textit{not} preserved by its evolution \eqref{eq:lin}, we \emph{do not} perform singularly weighted estimates to $W$.  Instead, the decomposition
\(W=W_1+\widehat W_2\) in Section \ref{sec:finite_pertb} is designed so that \(W_1\) preserves
cubic vanishing.

\subsubsection{Regularity of the finite-rank terms}\label{sec:reg_finite}

Denote $\eta = \th_x, \xi = \th_y$ as in \eqref{eq:eta_xi} and $W = (\om, \eta, \xi)$. From \eqref{eq:solu_reg1:a}, we have $W \in C^{3, 1}$. In Section \ref{sec:decoup_modi}, we decompose 
\beq\label{eq:reg_decomp}
W  =(\om, \th_x, \th_y) =  W_1 + \hat W_2, \quad  W_1 = (\om_1, \eta_1, \xi_1),
\quad \hat W_2 = (\hat \om_2, \hat \eta_2, \hat \xi_2),
\eeq
where $\hat W_2$ is represented by \eqref{eq:bous_W2_appr} and $(W_1, \hat W_2)$ solve \eqref{eq:bous_decoup2}. 
We now discuss the regularity of
the finite-rank and residual terms in \eqref{eq:bous_decoup2}. The finite-rank terms
\footnote{
$\cNF(W)$ in $\cK_F$ depends nonlinearly on $W$; we slightly abuse terminology
by calling $\cK_F$ a finite-rank operator.

}
\beq\label{eq:KF_reg1}
\cK_F = \vec{\cK}_1(W_1) + \vec{\cK}_2(W_1) + \cNF(W) = 
\tts{\sum}_{ 1\leq i\leq n_1 } \, 
a_i(W_1) \bar F_{i } + 
\tts{ \sum}_{1\leq i\leq 3}   \,  a_{ \mf{nl}, i}(W) F_{\chi, i}
\eeq
in \eqref{eq:bous_decoup2} consists of $\cK_{1i}$ from \eqref{eq:appr_near0}, $\cNF_i$ from \eqref{eq:appr_near0},\eqref{eq:appr_near02}, and $\cK_{2 i}$ from \eqref{eq:appr_vel},
where $\bar  F_{i } ,\bar F_{\chi, i} \in \R^3$. The approximations $\appo{\uu}, \appow{\na \uu}$ entering $\cK_{2 i}$ are constructed in \eqref{eq:u_appr_1st},\eqref{eq:u_appr_2nd},\eqref{eq:u_appr}. 

We record three qualitative facts about these terms: the scalar coefficients are
well-defined for the regularity class under consideration, the range functions are sufficiently smooth and decaying, and the velocity approximations preserve the vanishing orders needed near the origin.

We first estimate  $a_{\bullet}(\cdot)$  in
\eqref{eq:KF_reg1}. The terms $a_i(W_1(t))$ are linear functionals of $W_1$ and include
\begin{enumerate}[ leftmargin=3em, itemsep=0pt, topsep=1pt]
    \item $c_{\om}(\om_1)$ from $\cK_{1i}$ in \eqref{eq:u_appr_near0};
    \item $\cK_{00}(\om_1)$ and $\fNS(\om_1)$, up to constants independent of
    $\om_1$, from \eqref{eq:u_appr_1st}, \eqref{eq:u_appr_near0_K}, and
    \eqref{eq:u_appr_1st2};
    \item $I_i(\om_1)-u_x(\om_1)(0)$ from \eqref{eq:u_appr_In} and
    \eqref{eq:u_appr_2nd}.
\end{enumerate}
The nonlinear terms $a_{\mf{nl},i}(W(t))$ involve $\na^2 W(0)$ and
$\na^2\bar\cF(0)$; see \eqref{eq:appr_near0} and \eqref{eq:appr_near02}.

For $\om_1 \in C^{2, \al}, \al \in (0,1)$ with $\na^2 \om_1(0) = 0$, we have 
\[
  |\om_1(\xx)| \les \min( 1, |\xx|^{2 + \al}) || \om_1 ||_{C^{2,\al}}.
\]
For any $\ell >0$ and parameters $x_i$ for approximation \eqref{eq:u_appr_1st} chosen in Appendix \ref{app:appr_vel_para}, it follows 
\footnote{
\label{foot:LWP_bound_functional}
  We can rewrite $\cK_{00}(\om_1), I_i-u_x(0), \fNS(x_i,0)$ as $\int g(y) \om_1(y) dy$ for suitable kernels $g$. 
The kernels $g$ for $\cK_{00}$ in  \eqref{eq:u_appr_near0_K} and $I_i -u_x(0)$ in \eqref{eq:u_appr_In} satisfy $|g(y)|\les |y|^{-4} $, while the kernel $g$ for $\fNS(x_i, 0)$
satisfies $|g(y)| \les |y|^{-2}$. 
For $|y|\geq 1$ sufficiently large, one uses the odd extension of $\om_1(y)$ in $y_2$, its odd symmetry in $y_1$, and symmetrizes the integral. For example, after symmetrization, we rewrite  $I_i -u_x(0)$ in \eqref{eq:u_appr_In}  as 
\[
I_i -u_x(0) = \f{4}{\pi} \int_{ \max(y_1, y_2) \leq R_i, y_i \geq 0} \f{y_1 y_2}{|y|^4} \om(y) d y
\teq \int_{\R^2_{++}} g^{\mf{sym}}(y) \om(y) d y .
\]
The resulting kernel $g^{\mf{sym}}$ has bounded support and thus satisfies $|g^{\mf{sym}}(y)| \les C(x_i) |y|^{-4}$. Estimate for $u_x(\om_1)(0)$ in \eqref{eq:KF_reg2:a} follows from 
\eqref{eq:stream_local_reg1}. See Appendix \ref{app:euler_stream_sym} for the symmetrization argument and estimates. These symmetrized kernels are integrable against $W_1$ for $W_1 (1 + |x|^{-2 -\al})\in L^{\infty}$.
}
\bseq\label{eq:KF_reg2}
\beq\label{eq:KF_reg2:a}
\bal
    |c_{\om}(\om_1) |  & = |u_x(\om_1)(0)|  \les_{\al, \ell}  \| \om_1 \|_{C^{\al}} + \| \, |\xx|^{\ell}  \om_1 \|_{L^{\infty}}  ,  \\
    |\cK_{00}( \om_1 )  | 
     &+   |I_i(\om_1) - u_x(\om_1)(0)|
     \les || \om_1 ||_{C^{2, \al}}, \qquad  |\fNS(\om_1)(x_i, 0)|   \les  || \om_1  ||_{C^{\al}}  , \\
    |a_{ \mf{nl}, i}(W)|  & \les
     1 + |\na^2 W(0)| \cdot |c_{\om}(\om)    |
\les 
     1 + |\na^2 W(0)| \cdot 
(  \| \om \|_{C^{\al}} + \| \, |\xx|^{\ell}  \om \|_{L^{\infty}}  ) . 
\eal
\eeq

We next turn to the range functions in \eqref{eq:KF_reg1}. The functions in
$\cK_{1i}$ are $\bar f_{c_{\om},i}\in C^{2,1}$ by
\eqref{eq:appr_near0}, \eqref{eq:f_cw}, and \eqref{eq:solu_reg1:b}. The
range functions in the approximations of $\uu$ and $\na\uu$ are built from
the coefficients $C_f,C_{f0}$, the cutoff functions
$\chi_i,\chi_i^R,\chi_{\tot}$ in \eqref{eq:u_appr_1st} and
\eqref{eq:u_appr_2nd}, and their products. Since
$\chi_i,\chi_i^R,\chi_{\tot}$, defined in \eqref{eq:vel_chi_bd}, \eqref{eq:vel_cutoff_R}, and
\eqref{eq:u_appr_2nd}, are $C^\infty$ cutoff functions, and since
$C_f,C_{f0}$, defined in \eqref{eq:u_appr_near0_coe}, are polynomials in
$x,y$, these range functions are smooth.

The functions $\bar F_{ij}$ in $\cK_{2i}$, defined in \eqref{eq:appr_vel} and \eqref{eq:u_appr}, 
are products of $(\na\bar\om,\na^i\bar\th)\in C^{2,1}$, $i\leq 2$, with the smooth
coefficients $C_f,C_{f0}$ and cutoff functions $\chi_i,\chi_i^R,\chi_{\tot}$  described above. 
 Hence $\bar F_{ij}\in C^{2,1}$. Moreover, by definition of $f_{\chi,i}$ in
\eqref{eq:appr_near0} and \eqref{eq:cutoff_near0_all}, we have
$f_{\chi,i}\in C_c^\infty$. Combining these observations with
\eqref{eq:KF_reg2:a}, for $W_1\in C^{2,\al}$ with $\na^2W_1(0)=0$, we obtain
\beq\label{eq:KF_reg2:b}
\bal
  || \cK_{1 i}(W_1) ||_{C^{2, 1}} +   || \cK_{2 i}(W_1) ||_{C^{2, 1}} 
  & \les \| W_1 \|_{C^{2,\al}} + \| \, |\xx|^{\ell} \om_1 \|_{L^{\infty}} ,  \\
  \la \xx \ra^{10} \cdot |\cNF(W)(\xx)| +  \| \cNF(W) \|_{C^{2, 1}}  & \les | \na^2 W(0) | 
\cdot  ( \| \om \|_{C^{\al}} + \| \, |\xx|^{\ell} \om \|_{L^{\infty}} ) + 1.
\eal 
\eeq
\eseq

Using estimates \eqref{eq:u_appr_growth} for the approximation 
 $\appo{f}$ of $\uu, \na \uu$, \eqref{solu:decay} for $\bar \om, \bar \th$, 
 and estimates \eqref{eq:solu_reg1:b} for the range functions $  \bar f_{c_{\om}, i}$ of 
 $\cK_{1 i}$, and \eqref{eq:KF_reg2:a} for the functionals, we bound 
 \bseq\label{eq:KF_decay}
 \beq
| \cK_{1 i}(W_1)(\xx) | + | \cK_{2 i}(W_1)(\xx)| \les  \la \xx \ra^{\bar \al_1} 
\cdot ( \| W_1 \|_{C^{2,\al}} + \| \, |\xx|^{\ell} \om_1 \|_{L^{\infty}} ).
 \eeq

Using estimates similar to \eqref{eq:u_appr_growth} for the range functions in  $\appo{f}$, 
and \eqref{solu:decay} for $\bar \om, \bar \th$, we obtain
\beq\label{eq:KF_decay:b}
 |\bar F_{ij}| \les \la \xx \ra ^{\bar \al_1}, \quad f_{\chi, i} \in C_c^{\infty} , 
 \quad j \in \{ 1, 2, 3\}.
\eeq
\eseq

\paragraph{\bf Vanishing order of velocity and its approximation} 

We now discuss the vanishing orders associated with the velocity approximations. These estimates will be used later for the terms involving $\td \uu, \na \td \uu$, and their approximation errors. 
Recall 
the velocity approximations $\appo{f}$ from \eqref{eq:u_appr_1st}, \eqref{eq:u_appr_2nd}, \eqref{eq:u_appr}. By the definitions of $\hat f_2$ \eqref{eq:u_appr_2nd}, $\chi_i^R$ \eqref{eq:vel_cutoff_R}, and $\td \chi$ \eqref{eq:vel_chi}, 
we have $\hat f_2(\xx) = 0$ for $\max(|x_1|, |x_2|) \leq R_0$. By definition of $\hat f_{10}$ \eqref{eq:u_appr_1st} and $C_f(\xx)$ \eqref{eq:u_appr_near0_coe}, we have $|C_{\na \uu}| \les |\xx|^2, | C_{\uu} | \les |\xx|^3$. Thus, using \eqref{eq:u_appr} and \eqref{eq:KF_reg2:a}, the velocity approximations satisfy
\footnote{
The functionals in $\appo{f}$ \eqref{eq:u_appr} only 
involve $\cK_{00}, I_i-u_x(0), \fNS$. We do not need the bound for $u_x(0)$ in \eqref{eq:KF_reg2:a}.
}
\bseq\label{eq:vel_reg}
\beq\label{eq:vel_reg:a}
 |\appo{\uu}(\om_1) | \les |\xx|^3 || \om_1 ||_{C^{2,\al}} ,
  \quad   |\appow{ \na \uu}(\om_1)| \les |\xx|^2  || \om_1 ||_{C^{2,\al}} .
\eeq

Using $|\na \bar \th| \les |\xx|, |\bar \om| \les |\xx| $, \eqref{eq:vel_reg:a}, and the definition of $\cK_{2i}$ in \eqref{eq:appr_vel}, we obtain 
\beq\label{eq:reg_K2}
 |\vec \cK_{2}(W_1) | \les |\xx|^3 || \om_1 ||_{C^{2,\al}}.
\eeq

Next, we consider vanishing of the exact velocity. 
Assume that  $f\in C^{2,\al}$, $f |\xx|^\ell\in L^\infty$ for some
$\ell>0$, $f$ is odd in $x$, and  $\na f(0)=0$. 
Denote  $\uu(f)=\na^{\perp}\phi(f)$. 
Using identities in \eqref{eq:stream_vanish:a}, we obtain $\na^3 \phi(f)(0)=0$ and $\na^2 \uu(f)(0) = 0$. 
 Combining this with $v_y(f)(0)=-u_x(f)(0), u_y(f)(0) = v_x(f)(0) = 0$, 
 using estimates \eqref{eq:stream_local_reg1}, and Taylor expansion near $0$,  we estimate 
\beq\label{eq:vel_reg:b}
\bal
\td \uu(f) \teq \uu(f) - u_x(f)(0) (x, - y),  \quad 
|\td\uu(f)|
 & \les |\xx|^3\|f\|_{C^{2,\al}} ,
 \quad 
   |\na\td\uu(f)| & \les  |\xx|^2 \|f\|_{C^{2,\al}} .
\eal 
\eeq
\eseq
We apply this estimate to $f=\om$ and $f=\om_1$, using
$\na\om(0)=0$ from \eqref{eq:van_pertb_sq} and $\na\om_1(0)=0$.

In \eqref{eq:vel_rem}, we define velocities $\uu_A, \uu_{x, A}, \uu_{y, A} $ by subtracting the  
finite-rank approximations:
\bseq\label{eq:velA_reg}
\beq 
  \uu_A(f) = \td \uu(f) -  \appo{\uu }(f), 
\quad \uu_{x, A}(f)  =  \td \uu_x(f) - \appow{\uu_x}(f), \quad 
\uu_{y, A}(f)  = \td \uu_y(f) - \appow{\uu_y}(f).
\eeq

For $\om_1$ with $\na^i \om_1(0)=0, i\leq 2$, applying \eqref{eq:stream_local_reg2} near $\xx=0$, and applying \eqref{eq:vel_reg:a}, \eqref{eq:vel_reg:b}, and triangle inequality away from $\xx=0$, we obtain \textit{improved} vanishing orders for $\uu_A, (\na \uu)_A$ compared to $\td \uu, \na \td \uu$ 
\footnote{
\label{fn:uA_reg}
From \eqref{eq:u_appr_1st},\eqref{eq:u_appr_2nd},\eqref{eq:u_appr}, near $\xx=0$, we have $f_A = f - \hat f = f - (u_x(0) C_{f0}(\xx) + \cK_{00}(\om_1) C_f(\xx)) $
for $f = \uu, \na \uu$, which satisfies \eqref{eq:u_appr_near0}. 
We can improve the powers in \eqref{eq:velA_reg:b} to $|\xx|^{4-}$ for $\uu_A$ and $|\xx|^{3-}$ for $(\na \uu)_A$, since $\om_1 \in C^{2, 1}$ by Lemma \ref{lem:W1_reg}.
Yet, later we will replace 
$|| \om_1 ||_{C^{2, \al}}$ in \eqref{eq:velA_reg:b} by the energy $E_4$ (see Remark \ref{rem:E4_reg}). 
Since the norm $ [ \om_1 \psi_1 ]_{C_{g_1}^{1/2}}$  in $E_4$ (with $\psi_1, g_1$ defined in \eqref{wg:hol}) is similar to $|| \om_1 |\xx|^{-2} ||_{C^{1/2}}$ and $|| \om_1 ||_{C^{2, 1/2}}$ in terms of vanishing order near $0$, we can control \emph{only} the weaker powers $|\xx|^{7/2}, |\xx|^{5/2}$ in \eqref{eq:velA_reg:b} by $E_4$.
}
\beq\label{eq:velA_reg:b}
|\uu_A(\om_1)| \les |\xx|^{7/2}  \| \om_1 \|_{C^{2,1/2}} ,
\quad 
 | (\na \uu(\om_1))_A|  \les  |\xx|^{5/2}  || \om_1 ||_{C^{2,1/2}}  .
\eeq

\eseq

\subsubsection{Regularity of $\hat W_2, W_1$ and residual error}\label{sec:reg_W12}

We estimate the qualitative regularity and vanishing properties of $\hat W_2$ and $W_1$. Recall that $\hat W_2$ is built from the approximate space-time solutions $\hat F_i$ and $\hat F_{\chi,i}$ in
\eqref{eq:bous_W2_appr} and \eqref{eq:bous_err_op}. These functions satisfy properties \eqref{eq:bous_err_order}
\begin{align}\label{eq:Fi_supp}
   \hat F_i, \ \hat F_{\chi,i} & \in C_c^{4,1}, 
   \quad   \supp(\hat F_{\bullet}) \subset [-L_2, L_2] \times [0, L_2] , \ 
   \Longrightarrow  \
  & & \tf{1}{|\xx| } \uu(\hat F_{\om}) \in L^{\infty}, \  \na \uu(\hat F_{\om} ) \in C^{4, \g} , \,  \notag \\
\hat F_i,\hat F_{\chi,i} & = O(|\xx|^2),\qquad
   \hat F_i(0,\xx)-\bar F_i(\xx)=O(|\xx|^3),\qquad
  & & (\partial_t-\cL)\hat F_i=O(|\xx|^3) ,
\end{align}
for any $\g  < 1$, where $\hat F_{\om}$ is the $\om$-component of 
the solution  $\hat F_{i}$ or $\hat  F_{\chi,i, 1} $ for some $i$.

Since $\hat F_i, \hat F_{\chi,i}, f_{\chi,i}$ have compact support by \eqref{eq:Fi_supp}, \eqref{eq:KF_reg2:b}, using the decay of  $\bar F_{ij}$ in \eqref{eq:KF_decay:b}, the definitions of $\hat W_2$ \eqref{eq:bous_W2_appr} and $\cR_l, \cR_{ \mf{nl} }, \cR$ \eqref{eq:bous_err_op},  estimates \eqref{eq:KF_reg2} for the functionals, 
and the vanishing in \eqref{eq:Fi_supp}, for any $\al \in (0, 1)$, we obtain  the following qualitative estimates
\begin{align}\label{eq:W2_reg}
 || \hat W_2(t) ||_{C^{4, 1}} & + |\xx|^{-2} |\hat W_2(t, \xx)|  \les_{\al}  1 +  \sup\nolim_{s \leq t} \big( || W_1(s) ||_{C^{2,\al}} 
+ |\na^2 W(s, 0)|  \cdot ( \| \om(s) \|_{C^{\al}} + \| \om(s) |\xx|^{\ell} \|_{L^{\infty}} ) \big) \notag , \\
    || \cR_l(t) ||_{C^{2, 1}}  & + ( |\xx|^{-3} + |\xx|^{ -\bar \al_1}  )   | \cR_l(t, \xx) |    \les_{\al} \sup\nolim_{s \leq t}  || W_1(s) ||_{C^{2,\al}}, \notag  \\
 || \cR_{ \mf{nl}}(t) ||_{C^{2, 1}}  & + ( |\xx|^{-3} +  |\xx|^{ - \bar \al_1}  )   | \cR_{ \mf{nl}}(t, \xx) |  \les_{\al} 
1 + \sup\nolim_{s \leq t}  |\na^2 W(s, 0)|  \cdot ( \| \om(s) \|_{C^{\al}} + \| \om(s) |\xx|^{\ell} \|_{L^{\infty}} ) ,  \notag \\
& \hspace{1.1in} \supp(\hat W_2)    \subset  [-L_2, L_2] \times [0, L_2].  
\end{align}

Recall $\bar \al_1 < 0$ from \eqref{solu:decay}. Estimate $ (|\xx|^{-3} + |\xx|^{ - \bar \al_1}) g \in L^{\infty}$ 
implies the cubic vanishing $|g| \les |\xx|^{3}$  and far-field decay $|g| \les |\xx|^{ \bar \al_1}$. The implicit constants may depend on the pre-computed approximate solution 
$\hat F_i, \hat F_{\chi, i}$, the profile $\bar \om, \bar \th$, and the parameters used in the finite-rank construction. 

Recall the decomposition \eqref{eq:reg_decomp} of $(\om, \th)$
Since $(\om_{\tot}, \th_{\tot})$ and  perturbation $(\om, \th)$ satisfy the symmetries \eqref{eq:Sym}, and $\hat F_i$ satisfies the symmetries 
\eqref{eq:bous_err_order:a}, we have 
\beq\label{eq:Sym_W12}
  \omega_1,\eta_1,\hat\omega_2,\hat\eta_2
   \quad\text{are odd in }x,\qquad
   \xi_1,\hat\xi_2
   \quad\text{are even in } x, 
   \qquad
   W_1(0,y)=\hat W_2(0,y)=0 .
\eeq

With the above preparations, we have the following local-in-time regularity for $W_1$.

The next lemma justifies the use of singularly weighted estimates for $W_1$:
as long as the standard continuation criterion holds, $W_1$ remains $C^{2,1}$
and vanishes to $O(|x|^3)$ near the origin.

\begin{lem}[\bf Regularity of $W_1$]\label{lem:W1_reg}
Suppose that the initial data satisfy 
$\om_{\tot, 0} \in C^{\kin, \ain}$, $ \th_{\tot,0} \in C^{\kin+1, \ain}  $ for some $\kin \geq 5, \ain \in (0, 1)$, and conditions \eqref{eq:reg_init},\eqref{eq:van_init_cubic},\eqref{eq:Sym}.
Assume also that 
\beq\label{eq:W1_BKM}
M(\bar t) = \tts{\int}_0^{\bar t} || \na(\uu_{\tot}, \th_{\tot})(t) ||_{L^{\infty}} d t < +\infty ,
\eeq
for some $\bar t > 0$. Then, for $ 0\leq t \leq \bar t$,
\bseq\label{eq:reg_W1}
\beq\label{eq:reg_W1:a}
W_1(t) \in C^{2, 1}, \quad  \na^i W_1(0, t) = 0, i \leq 2,  \quad |W_1(t, \xx ) | \les 
C_{\bar t}
\min( x_1 |\xx|^2, 1).
\eeq
 Moreover, if $\bar\alpha_1$ denotes the far-field exponent of the profile
introduced in Section~\ref{sec:fit_profile}, then $\hat W_2$ has compact
support and
\beq\label{eq:reg_W1:b}
 |\na^i \om_1(t, \xx)| \les C_{\bar t} |\xx|^{\bar \al_1 - i},
 \quad  |\na^i(\eta_1, \xi_1)(t, \xx)| \les C_{\bar t}  |\xx|^{2 \bar \al_1 - i}, \quad i \leq 2.
\eeq
\eseq
Here $C_{\bar t}=C(M(\bar t),\om_{\tot,0},\th_{\tot,0})$ may change from line
to line.
\end{lem}

\begin{proof}

By the regularity and vanishing results established in Section  \ref{sec:reg_van_pf}, 
the full perturbation $(\om, \th)$
satisfies the regularity \eqref{eq:solu_reg1} and vanishing order 
\eqref{eq:van_pertb_sq} on $ [0, \bar t]$ with $\bar t = t(\om_{\tot,0}, \th_{\tot,0} ) > 0$. 
Since $W = W_1 + \hat W_2$ and $\hat W_{2,0}(x) \equiv 0$ \eqref{eq:bous_W2_appr},  we deduce that $W_{1,0}(x) = O(|\xx|^3)$ and $W_{1, 0} \in C^{2, 1}$.

We derive a closed equation for $W_1$ whose forcing is regular and vanishes cubically near the
origin on $[0, \bar t]$. We define a forcing $F_{W, i}$ and fix $\ell$
\[
  F_{W, i}  \teq \cN_i( W ) + \bar \cF_i - \cNF_i(W) - \cR_{ \mf{nl}}(W), \quad  i \leq 3  ,
  \quad \ell \in (0, |\bar \al_1|).
\]

From the definitions of $\cN_i, \bar \cF$ \eqref{eq:bous_non}, and of $\cNF$ \eqref{eq:appr_near0}, using Taylor expansion, the regularity of $W$ \eqref{eq:solu_reg1} and of $\cR_{ \mf{nl}}$ \eqref{eq:W2_reg}, 
and the decay estimates for $\bar \cF_i, \om, \th$ in \eqref{eq:solu_reg1:decay}, \eqref{eq:solu_reg1:b}, $\cR_{ \mf{nl}}$ in \eqref{eq:W2_reg}, and $ \cNF_i(W)$ in \eqref{eq:KF_reg2:b}, 
for $t \in [0, \bar t ]$,   we obtain that $  F_{W, i}(t, \cdot ) \in C^{2, 1}$ and satisfies  
\beq\label{eq:LWP_force}
\bal
  \| F_{W, i}(t, \cdot ) \|_{C^{2, 1}} & \leq C_{\bar t}, \\
   |F_{W, i}(t, \xx)| & \les \min( |\xx|^3, |\xx|^{\bar \al_1}) ( 1+ \| W \|_{C^{3, 1}} + \| \, |\xx|^{\ell} \om \|_{L^{\infty}} )^2 
  \leq C_{\bar t} \min( |\xx|^3, |\xx|^{\bar \al_1})  .
\eal 
\eeq
 Recall $\cL$ from \eqref{eq:lin}. From \eqref{eq:bous_decoup2},\eqref{eq:bous_err_op}, $W_1$ solves a linear equation with a forcing $F_{W_1, i}$ 
\beq\label{eq:reg_W1_pf1} 
  \pa_t W_1 = (\cL - \cK_1 - \cK_2 - \cR_l ) W_1 + F_{W}.
\eeq

The regularity \eqref{eq:solu_reg1:a} implies that the coefficients in $\cL$ \eqref{eq:lin} 
  are in $C^{2, 1}$ and that $\na (\bar c_l \xx + \bar \uu) \in C^{2, 1}$. 
  We fix $\al \in (0, 1)$ and $\ell \in (0, 10^{-1})$ and consider the following class of functions 
  \[
      Z \teq \{ G \in C^{2,\al}  : G |\xx|^{\ell} \in L^{\infty},\,
        G = (G_{\om}, G_{\eta}, G_{\xi})\ \mbox{satisfies symmetries} \ \eqref{eq:Sym_W12},  \, \na^i G(0) = 0, i \leq 2 \} .
  \]

For $G \in Z$,  using \eqref{eq:vel_reg:b}, we obtain 
\footnote{
  We add $S = (\bar c_l x + \bar \uu)\cdot \na G$ in \eqref{eq:reg_W1_pf2} so that the 
 transport term $S$ is cancelled in the expression \eqref{eq:reg_W1_pf2}.
}
\bseq\label{eq:reg_W1_pf2}
\beq
\na^i ( (\cL - \cK_1)(G) + (\bar c_l x + \bar \uu)\cdot \na G ) (0) =0, \quad i \leq 2.
\eeq

Moreover, for $G \in Z$, from \eqref{eq:reg_K2}, \eqref{eq:W2_reg}, and \eqref{eq:reg_W1_pf1}, we obtain 
\beq
  \na^i \vec \cK_2(G)(0) =0, \quad  \na^i \cR_l(G)(0) = 0, \quad \na^i F_{W}(0) = 0, \quad i \leq 2 .
\eeq
\eseq

Thus all terms on RHS of \eqref{eq:reg_W1_pf1}, except for the 
transport term $(\bar c_l x + \bar \uu)\cdot \na G$, preserves the vanishing: $\na^i G(0)=0, i\leq 2$.
We use \eqref{eq:reg_W1_pf2}, 
estimates for $F_{W, i}$ in \eqref{eq:LWP_force},  estimates for $\cK_{1i}, \cK_{2i}, \cR_l$ in \eqref{eq:KF_reg2} and in \eqref{eq:KF_decay}, the nonlocal estimates \eqref{eq:stream_local_reg1}, and perform $C^{2,\al}$ estimates on $W_1$ and weighted $L^{\infty}$ estimates on 
$ W_1 |\xx|^{\ell}$. Applying the standard local well-posedness arguments to \eqref{eq:reg_W1_pf1} in the class $Z$, for example by the particle-trajectory method and iteration  \cite{chae1999local,elgindi2017finite}, we obtain local well-posedness of $W_1$-equation  \eqref{eq:reg_W1_pf1} in $Z$.
\footnote{
Let $\bar X(t , \cdot )$ be the flow map generated by $\bar c_l \xx + \bar \uu$. One may solve \eqref{eq:reg_W1_pf1} using the scheme 
\[
  \pa_t (G_{n+1} \circ \bar X ) = [ (\cL +  (\bar c_l \xx + \bar \uu) \cdot \na ) G_n  - (\cK_1 + \cK_2 + \cR_l) G_n + F_W ] \circ \bar X
\]
with $G_0(t, x) \equiv W_1(0, x) \in Z$. There is no transport term $(\bar c_l \xx + \bar \uu) \cdot \na G_n$ on the right side. 
Since $\xx = 0$ is a stationary point, i.e. $\bar X(t, 0) \equiv 0$, conditions \eqref{eq:reg_W1_pf2} imply that the iteration preserves $\na_x^i  G_n |_{x =0} = 0$ for $i \leq 2 $.  
}
Since  \eqref{eq:reg_W1_pf1} is a linear equation on $W_1$ with a forcing $F_{W_1}(t) \in Z$ for $t \in [0, \bar t]$ by \eqref{eq:LWP_force}, we obtain $W_1(t) \in Z$ for $t \in [0, \bar t]$.

Using  \eqref{eq:W2_reg}, \eqref{eq:solu_reg1:a} and \eqref{eq:reg_decomp},  we have $W= W_1 + \hat W_2 \in C^{3,1}$ and $\hat W_2 \in C^{2, 1}$ on $[0, \bar t]$.
Thus, $W_1(t) \in C^{2, 1}$ as well.
 Since $\na^i W_1(t, 0) = 0$ for $i \leq 2$, 
 from the symmetries \eqref{eq:Sym}, we obtain $|W_1(t, x) | \les C_{\bar t} |\xx|^2 |x_1|$. 
From \eqref{eq:W2_reg}, $\hat W_2$ has compact support. Since $W = W_1 + \hat W_2$ satisfies the far-field estimates in \eqref{eq:solu_reg1:decay}, we derive \eqref{eq:reg_W1:b} from 
\eqref{eq:solu_reg1:decay}. This completes the proof.
\end{proof}

\begin{remark}[\bf Bounded energies $E_i$]
Under the assumptions of Lemma \ref{lem:W1_reg}, 
\eqref{eq:reg_W1} implies 
\[
W_1 \cdot  (|\xx|^{-3} + 1)\in L^{\infty}, \quad W_1 \cdot ( |\xx|^{-5/2} + 1 ) \in C^{1/2}.
\]
From \eqref{eq:reg_W1}, one can further verify the boundedness of the energies $E_i$, defined later in \eqref{energy1}, \eqref{energy2}, \eqref{energy3}, \eqref{energy4}.
 Using \eqref{eq:velA_reg}, we also obtain $\uu_A \psi_u , (\na \uu)_A  \psi_1 
  \in \dot{C}^{1/2}$, where $\psi_u, \psi_1$ are the weights defined in \eqref{wg:hol}, 
  \eqref{wg:elli} for weighted nonlocal estimates.
\end{remark}

\begin{remark}[\bf Energy $E_4$ controls regularity]\label{rem:E4_reg}

We \emph{only} use the higher regularity of $(\om, \na \th)$ in the above argument \emph{qualitatively}: it justifies the Taylor expansion, vanishing order, and regularity of the auxiliary quantities. A key point is that this \emph{qualitative} higher regularity can be propagated by 
\emph{quantitatively} controlling much lower-regularity norms, which in our setting corresponds to the energy $E_4$. This is in the spirit of the BKM continuation criterion.

For $\na^i W_1(0) =0, i\leq 2$, any $\al \in (0, 1]$ and $\b >0$, we bound the functionals  \beq\label{eq:functional_w_est}
| \cK_{00}(\om_1)| + |c_{\om}(\om_1)| + |I_i(\om_1) - u_x(\om_1)(0)| +| \fNS(\om_1)(x_i, 0)|
\les || \om_1( |\xx|^{-2-\al} + |\xx|^{\b} ) ||_{L^{\infty}} . 
\eeq
This requires less regularity of $\om_1$ than the estimates in \eqref{eq:KF_reg2:a}. 
 See Footnote \ref{foot:LWP_bound_functional} and Section \ref{sec:rank1}. Since energy $E_4(W_1, W)$ defined in \eqref{energy4} 
controls this weighted norm, $\na^2 W(0)$, and $|c_{\om}(\om)|$, we obtain 
estimates   \eqref{eq:KF_reg2}, \eqref{eq:vel_reg:a}, \eqref{eq:reg_K2}, \eqref{eq:W2_reg} with 
norm $|| W_1 ||_{C^{2, \al}}$
and functionals $|\na^2 W(0)|, c_{\om}(\om)$ replaced by energy $E_4$. Similarly, by symmetrizing integrals in $\td \uu(\om_1), \na \td \uu(\om_1)$, $\uu_A(\om_1), (\na \uu(\om_1))_A$, we can bound their weighted 
$L^{\infty}$ and $C^{1/2}$ norms by $E_4$, and obtain estimates \eqref{eq:vel_reg:b} and \eqref{eq:velA_reg:b} with norms $|| W_1 ||_{C^{2, \al}}, | W_1 ||_{C^{2, 1/2}}$ replaced by $E_4$.

Since $\vp_{g, i} \gtr 1$ by \eqref{wg:linf_grow}, 
 the definitions of $E_3$ \eqref{energy3} and $E_4$ \eqref{energy4} imply 
 $
    || \na \th_1 ||_{L^{\infty}} = || (\eta_1, \xi_1) ||_{L^{\infty}} \les E_3 \les E_4.
 $
Using the decomposition \ref{eq:whole_solu} and \eqref{eq:reg_decomp}, we obtain
\[
|( \na \uu(\om_{\tot}), \na \th_{\tot})| 
\les 1 + |(\na \uu( \hat \om_2 ), \hat \eta_2, \hat \xi_2 ) | 
+ | ( \na \uu( \om_1 ), \eta_1, \xi_1) | \les 1 + E_4^2(W_1, W).
\] 
Thus, $\int_0^t E_4^2(s) ds$ controls the continuation criterion \eqref{eq:W1_BKM}.
We will establish in Theorem \ref{thm:main} that $E_4(t)$ remains small for any $t \geq 0$, which justifies the derivations in this Section and Section \ref{sec:EE}.
\end{remark}

\subsubsection{ Qualitative regularity of quantities in the energy estimates}\label{sec:reg_other}

For later use in Section~\ref{sec:EE}, we collect the qualitative regularity and vanishing
orders of several auxiliary terms defined there: $\UU$, $\UU_A$, and
$\UU_{\FM}$ in \eqref{eq:U}; $\bar \e$ in \eqref{eq:wb_vel_err};
$\bar \e_1$, $\bar \e_2$, $\hat \e_1$, and $\hat \e_2$ in
\eqref{eq:uerr_dec1}; $\hat \e$ in \eqref{eq:lin_evo_main1};
$\td \cN_i$ in \eqref{eq:bous_nonM}; $\bar \cF_{loc,i}$ in
\eqref{eq:bous_errM}; $\cR_{loc,i}$ in \eqref{eq:W2_errM}; and
$\hat W_{2,\cdot,M}$ in \eqref{eq:W2_M2}. Under the assumptions of Lemma \ref{lem:W1_reg}, and using the above 
estimates, we obtain the following regularity and vanishing orders of these terms near $0$, with regularity controlled by energy $E_4$  and universal constants:
\beq\label{eq:reg_other}
\bga
\bar \e , \ \bar \e_1 , \ \hat \e, \ \hat \e_1  \in C^{2, 1}, \quad  \bar \e, \ \hat \e   = O(|x|^2),  \quad \bar \e_1, \ \hat \e_1 = O(|x|^3),  \\
\UU_{\FM} \in C^3, \quad  \UU_{\FM} = O(|x|^3), \\ 
\UU_A, \ \uu_A(\bar e_1), \ \uu_A(\hat \e_1) = O(|x|^{7/2}),  \quad  (\na \UU)_A,\ (\na \uu)_A(\bar \e_1), \
(\na \uu)_A(\hat \e_1) 
 = O(|x|^{5/2}),  \\
\bar \cF_{loc, i}, \ \cR_{loc, i} \in C^{2, 1}, \quad  \td \cN_i, \ \bar \cF_{loc,i}, \ \cR_{loc, i} = O(|\xx|^3), 
\quad 
 \hat W_{2, \cdot ,M} \in C_c^{3, 1}, 
   \  \hat W_{2, \cdot ,M} = O(|\xx|^3) .
 \ega
\eeq
The proof follows from previous arguments. In Section \ref{sec:EE}, we provide quantitative estimates.

\section{Energy estimates}\label{sec:EE}

In this section, we perform the energy estimates for the perturbation $W_1$ in
\eqref{eq:bous_decoup2}. These estimates are the core of the stability argument. 
Our goal is to control several weighted $L^\infty$ and H\"older norms
of $W_1$, and to establish the nonlinear stability estimates
\eqref{eq:PDE_diag}, \eqref{eq:PDE_nondiag} in Lemmas~\ref{lem:PDE_stab} and
\ref{lem:PDE_nonstab}. We then use these estimates to prove
Theorem~\ref{thm:main}, \ref{thm:main_stab_full}.

The main difficulties are the coupling between
weighted $L^\infty$ and $C^{1/2}$ estimates and  the nonlocal estimates with
finite-rank approximations. 
The lowest-order weighted $L^\infty$ estimates do not close by themselves.  A
key point is to develop methods that make their dependence on the
$C^{1/2}$ norms sufficiently weak, so that the coupled stability estimates can 
be established.

The estimates involve several parameters, which play an
important role in controlling the coupling and relative sizes of different norms in the energies.
They are determined sequentially by the stability conditions. In particular,
after the key linear estimates are established, later quantities are added to the
energy through damped estimates of the form \eqref{eq:damp_sys}, with 
coupling parameters chosen small enough to preserve the linear stability estimates.

 To clarify the structure of this section, we first isolate the main equations, then outline the estimates in Section~\ref{sec:EE_intro} and the dependencies
among estimates and parameters in Figure~\ref{Fig:order}.

\vs{0.02in}

\paragraph{\bf Analytic estimates and computer assistance}
Before entering the estimates, we clarify the role of computer assistance in
this section. The weighted $L^\infty$ and $C^{1/2}$ energy estimates, 
ODE estimates of several functionals, 
modified decomposition for nonlocal error estimates, the nonlinear estimates, and the bootstrap argument are derived analytically. 
These analytic estimates reduce nonlinear stability to finite inequalities involving computable coefficient functions.

Computer assistance is used to bound explicit integrals 
for nonlocal estimates, derive piecewise bounds of explicit coefficients, track coefficients in the energy estimates, and rigorously verify inequalities for nonlinear stability 
in Appendix \ref{sec:ineq}.

\subsection{The main equation}\label{sec:EE_intro}

We first isolate the key linear part of the $W_1= (\om_1, \eta_1, \xi_1)$ equation. After choosing approximations for the velocity $ \uu$ and using the decomposition in Section \ref{sec:finite_pertb}, the main equations \eqref{eq:bous_decoup2}, corresponding to the leading part of \eqref{eq:lin}, for $W_1$ read
\beq\label{eq:lin_main}
\bal
\pa_t \om_1 + (\bar c_l x + \bar \uu  ) \cdot \na \om_1 &= \bar c_{\om} \om_1 + \eta_1  -  \uu_{A}(\om_1) \cdot \na \bar \om ,  \\ 
\pa_t \eta_1 + (\bar c_l x  + \bar \uu  ) \cdot \na \eta_1  & = (2 \bar c_{\om} - \bar u_x) \eta_1 
-  \bar v_x \xi_1 -  \uu_A(\om_1) \cdot \na \bar \th_x -   \uu_{x, A}(\om_1) \cdot \na \bar \th , \\
\pa_t \xi_1 + (\bar c_l x + \bar \uu  ) \cdot \na \xi_1 & = 
(2 \bar c_{\om} + \bar u_x ) \xi_1 -  \bar u_{y} \eta_1 - \uu_A(\om_1) \cdot \na \bar \th_y -  \uu_{y, A}(\om_1) \cdot \na \bar  \th,
\eal
\eeq
Here, $\bar \uu=(\bar u,\bar v)$. Unless otherwise specified, $\uu_A$, $\uu_{x,A}, \uu_{y,A}$ are the velocity remainders depending on $\om_1$, rather than $\om$, after subtracting the finite-rank approximations in 
\eqref{eq:u_appr_1st}, \eqref{eq:u_appr_2nd}, \eqref{eq:u_appr}:
\bseq\label{eq:vel_rem}
\beq
\uu_A(f) \teq \td \uu(f) -  \appo{\uu}(f), 
\quad \uu_{x, A}(f)  \teq  \td \uu_x(f) - \appow{\uu_x}(f), \quad 
\uu_{y, A}(f)  \teq \td \uu_y(f) - \appow{\uu_y}(f).
\eeq
  From the definitions \eqref{eq:u_appr} and \eqref{eq:u_tilde}, we also have
\beq
\uu_A(f) = \uu(f) -  \hat {\uu}(f), \quad 
 \uu_{x, A}(f)  =   \uu_x(f) - \wh {\uu_x}(f),  \quad 
\uu_{y, A}(f)  \teq \uu_y(f) - \wh{\uu_y}(f).
\eeq
\eseq
These terms have the vanishing orders \eqref{eq:velA_reg}. Note that we do not have $ \pa_x ( \uu_A ) = \uu_{x, A}$ since we choose the approximations for $\uu, \uu_x, \uu_y$ separately. Similarly, we do not have $\pa_y \eta_1 = \pa_x \xi_1$.

Equation~\eqref{eq:lin_main} corresponds to the linear operator $\cL_i - \cK_{1 i} - \cK_{2 i}$ in \eqref{eq:bous_decoup2}. It is obtained by subtracting
$\cK_{2i}$ in \eqref{eq:appr_vel} from $\cL$ in  \eqref{eq:lin}, which
replaces the ``tilde" velocity $(\td \uu, \na \td \uu)$  by the remainders  
$( \uu_A, (\na \uu)_A)$.  We also subtract $\cK_{1i}$ \eqref{eq:appr_near0} from $\cL$ to remove 
$c_{\om}\bar f_{\chi,i}$ terms in \eqref{eq:lin}.
To simplify the presentation, we temporarily omit the residual terms $\cR_i$, the correction
$\cNF(W_1+\hat W_2)$, nonlinear $\cN_i$, and error term $\bar \cF_i$ in \eqref{eq:bous_decoup2}.
In Sections~\ref{sec:linf_decay}--\ref{sec:rank1}, we first establish linear stability estimates for \eqref{eq:lin_main}. 
In Section~\ref{sec:comb_vel_err}, we extend these estimates 
to the full system,  which reintroduces the residual, nonlinear, and nonlocal error terms and modifies the decomposition.

We adopt the notation from \eqref{eq:bous_decoup2} and introduce $\bar W$
\beq\label{eq:EE_W1}
W_{1, 1} = \om_1, \quad W_{1, 2} = \eta_1, \quad W_{1, 3} = \xi_1 , \quad \bar W = (\bar \om, \bar \th_x, \bar \th_y).
\eeq

We consider initial data satisfying \eqref{eq:Sym}, \eqref{eq:reg_init}, \eqref{eq:van_init_cubic}. Then the perturbation $W_1$ preserves  $W_1 = O(|x|^3)$ and satisfies \eqref{eq:Sym_W12} and Lemma \ref{lem:W1_reg}; 
see Section \ref{sec:reg_van}.

We now introduce the notation for weighted estimates: $\cT_d(\rho), d_{i, L}(\rho)$  denote the coefficients of the damping terms in weighted estimates  and $b(x)$ denote the coefficient of the advection 
\beq\label{eq:dp}
\bal
 b(x) & = \bar c_l x + \bar \uu  , \quad  &   \cT_d(\rho) & =  \rho^{-1} \B(  (\bar c_l x + \bar \uu  ) \cdot \na \rho \B) 
 = \rho^{-1}  ( b \cdot \na \rho ) , \\
 d_{1, L}(\rho)  &  = \cT_d( \rho ) + \bar c_{\om} , \quad 
&  d_{2, L}(\rho) & = \cT_d( \rho) + 2 \bar c_{\om} - \bar u_x, \qquad
d_{3, L }(\rho)  = \cT_d( \rho) + 2 \bar c_{\om} + \bar u_x .
\eal
\eeq
The subscript \textit{L} is short for linear.
In the equation of $W_{1, i}$, we treat the terms other than the local terms of $W_{1, i}$ in \eqref{eq:lin_main} as forcing, or ``bad" terms 
\footnote{
Here ``bad'' means that these terms are not part of the local damping terms.  In the weighted estimates, they are compared against the damping generated by the local terms and bounded perturbatively.
}
:
\beq\label{eq:bad}
\bal
 B_1(W_1) & \teq \eta_1  -  \uu_{A}(\om_1) \cdot \na \bar \om ,  \\
B_2(W_1) & \teq -  \bar v_x \xi_1 -  \uu_A(\om_1) \cdot \na \bar \th_x -   \uu_{x, A}(\om_1) \cdot \na \bar \th , \\
 B_3(W_1) &\teq -  \bar u_y \eta_1  - \uu_A(\om_1) \cdot \na \bar \th_y -  \uu_{y, A}(\om_1) \cdot \na \bar \th.
\eal
\eeq

With the above notations, we simplify \eqref{eq:lin_main} as 
\[
\pa_t W_{1 ,i} + b \cdot \na W_{1, i} = d_{i,L}(1) W_{1, i} + B_i,
\]
where $d_{i, L}(1)$ denotes the coefficient in \eqref{eq:dp} with weight $\rho \equiv 1$ and $\cT_d(1) = 0$. Using the above equations, 
for a \emph{time-independent} weight $\rho$, the weighted quantity $W_{1, i} \rho$ satisfies
\beq\label{eq:lin_sim}
\pa_t (W_{1, i}\rho) + b \cdot  \na (W_{1,i}  \rho ) 
= d_{i,L}( \rho) W_{1, i} \rho + B_i  \rho.
\eeq

In each weighted estimate, we compare the damping coefficient $d_{i,L}( \rho)$ with 
the forcing $B_i  \rho$.

\vs{0.1in}

\paragraph{\bf Dependency of estimates and parameters}\label{sec:EE_depend}

The estimates and parameter choices are interdependent. We organize them according to the order in Figure~\ref{Fig:order}. This figure summarizes the dependencies among weights, parameters, and estimates, and outlines estimates in Sections~\ref{sec:linf_decay}--\ref{sec:non}. We discuss the ideas of nonlocal estimates in Section~\ref{sec:u_comp_idea} and state the nonlinear stability and blowup theorem in Section~\ref{sec:EE_stab_thm}. Below, we explain the stability estimates step by step.

In Step~1, we choose different powers $|x|^p$ in $\rho_{ij}$ \eqref{wg:u}, \eqref{wg:elli} to capture the vanishing order near $0$ and decays of $\uu_A, (\na \uu)_A$; 
 the precise coefficients of these powers are less important and are chosen \emph{independent} of later estimates.

In Step~2, we perform top-order weighted $C^{1/2}$ estimates, neglecting
the more regular terms and following the estimates in Section~\ref{sec:hol_singu}.
We choose the H\"older weights $\psi_i, g_i$ to obtain coercive estimates
for the $C^{1/2}$ estimates in the near-diagonal regime
$|\xx-\zz|\to0$; see Section~\ref{sec:wg_hol_g}.

In Step~3, we perform the weighted $L^\infty(\vp_i)$ estimates and choose the
decaying weights $\vp_i$ in \eqref{wg:linf_decay}. Following
Section~\ref{sec:model_local}, we design these weights to derive 
sufficient strong local damping terms.  At this stage, we
also estimate the \emph{unapproximated} nonlocal terms $\td\uu$ and $\na\td\uu$ in
$\cL$ \eqref{eq:lin}, using the ideas of Section~\ref{sec:u_comp_idea}, to
identify the region $D$ where these terms are not small enough to be treated
perturbatively to the local damping terms.

In Step~4, we construct the finite-rank approximations
$(\hat \uu,\wh{\na \uu})$ for $(\uu,\na \uu)$ in
Section~\ref{sec:appr_vel}. The purpose is to improve the bounds for
$\uu_A=\uu-\hat\uu$ and $(\na\uu)_A=\na\uu-\wh{\na\uu}$ in the region $D$
identified in Step~3, until they are small compared with the local damping
terms.

Step~5 is the key linear stability step. We perform the weighted
$L^\infty(\vp_i)$ estimates with the energy $E_1$ \eqref{energy1} in
Section~\ref{sec:linf_decay}, and the weighted $C^{1/2}$ estimates with the
energy $E_2$ \eqref{energy2} in Section~\ref{sec:EE_hol}. These two estimates
are coupled, since neither estimate closes on its own. To establish the stability condition \eqref{eq:PDE_diag}, we use two large-parameter mechanisms from Section~\ref{sec:analy_para} and the ideas from Section~\ref{sec:outgo}. By adjusting the weights, we generate larger damping
away from $0$ and increase the diagonal coefficients $a_{ii}$ in
\eqref{eq:PDE_diag}. By increasing the approximation rank while keeping the
weights fixed, we improve the estimates of the more regular nonlocal terms and
decrease the off-diagonal coefficients $a_{ij}$, $j\neq i$. This allows us to
close the $L^\infty$ estimates with only weak coupling to the $C^{1/2}$ norm,
encoded by the small weight $\tau_1^{-1}$ in $E_1$ \eqref{energy1}; the same weak coupling then
allows us to close the $C^{1/2}$ stability estimates. This mechanism is
analogous to \eqref{eq:model3_linf} for the model problem. We also adjust the weights
$\vp_2,\vp_3$ and the finite-rank approximation to achieve stability with a
lower-rank approximation.

\vspace{0.05in}

\paragraph{\bf Damped system}

After the key linear stability estimate in Step~5, the remaining estimates are
\emph{perturbative} in nature. The energy $E_2$ already controls $W_1$ and the
nonlocal terms in \emph{any bounded region}. Since the coefficients in \eqref{eq:lin_main} decay in the far field, the
estimates in Steps~6--7 are perturbative and follow from the same damped
transport estimate with forcing:
\beq\label{eq:damp_sys}
    \pa_t f   + b \cdot \na f   = d \cdot f   + F  ,
  \tag{Damp}
\eeq
Here $f$ denotes a quantity to be added to the energy, and $F$ consists of terms
controlled by previously established energies. After performing the weighted
estimate as in Section~\ref{sec:model_local}, we have a damping coefficient
$d(x)\leq -a<0$. The forcing term satisfies
\beq\label{eq:damp_force}
  |F(t, x) | \leq C_1(x)  E_{\mw{stab}} + C_2(x) E_4^2,
  \quad  || C_i ||_{L^{\infty}} \leq c_i < \infty, 
\eeq
where $E_{\mw{stab}}$ denotes an energy for which stability has already been established, and $E_4$ is the final energy \eqref{energy4}. 
The second term in $F$ is quadratic and is treated perturbatively.

Choosing $\mu>0$ with $\mu c_1<a$, up to the quadratic $E_4^2$ terms, 
Lemma~\ref{lem:PDE_stab} establishes linear stability for the following enlarged energy with stability factor $\lam$
\footnote{
In \eqref{lem:PDE_stab} for the norm $|| f||_{L^{\infty}}$, we have $a_{ii} = - d(x) \geq a, \mu_i = \mu$, $C_1(x) = \sum_{j\neq i} a_{ij} \mu_j^{-1}$ for the bound of $B_i$.
} :
\[
E_{\mw{new}} = \max( E_{\mw{stab}}, \mu  || f||_{L^{\infty}}),
\quad 
\lam = \min( \lam_{\mf{stab}}, a - \mu c_1  ) > 0,
\]
where $\lam_{\mf{stab}}$ denotes the stability factor for energy $E_{\mw{stab}}$. 
We use this perturbative estimate to add further quantities to the energy after
Step~5, while preserving linear stability.

In Step~7, we also estimate the ODEs of some functionals 
to improve nonlinear constants, e.g. $a_{ij,2}$ in \eqref{eq:PDE_nondiag}, 
using the same damped estimate with $f(t,x)$ independent of $x$. Such technical
improvements would not be necessary if the residual error could be reduced
significantly.
\footnote{
This is much easier to achieve for 1D PDEs, especially when there are no
nonlocal terms.
}

In Step 8, the nonlinear terms are much smaller and are treated perturbatively.

\begin{figure}[t]
\centering

\begin{tikzpicture}[
  >=Stealth,            
  node distance=3.5mm and 10mm,  
  box/.style={          
    draw, 
    rectangle, 
    rounded corners, 
    align=center,
    inner sep= 4pt,
    text width=9cm
  },
  smallbox/.style={     
    draw, 
    rectangle, 
    rounded corners, 
    align=center,
    inner sep=4pt,
    text width=4.5cm
  }
]

\node[box] (profi) {%
1. Approximate profile 
$\bar \om, \bar \th, \bar u, \bar c_{\om},\bar c_l $,
weights $\rho_{ij}$, \\
 sharp H\"older estimates constants from Sec \ref{sec:sharp}
};

\node[box, below=of profi] (wghol) {%
2. Top order weighted $C^{1/2}$ estimates by \textit{neglecting} more regular terms: 
 choose weights $\psi_i$, $g_i(h), i\leq 3$, 
 parameters $\mu_1, \mu_2$ in $E_2$ \eqref{energy2}.
};

\node[box, below=of wghol] (wglinf1) {%
3. Weighted $L^{\infty}$ estimates with decaying weights 
on \textit{local} terms Sec \ref{sec:linf_decay}: choose weights $\vp_i, i \leq 3$.

};

\node[box, below=of wglinf1] (finite_rank) {
4. Construct finite-rank approximation $\hat \uu, \wh{\na \uu}$ of $\uu, \na \uu$
in region $D$ (Sec \ref{sec:appr_vel}) to improve bounds for $  \uu - \hat \uu,\na \uu - \wh{\na \uu}$
};

\node[smallbox, right=of wglinf1] (wg_u1) {
Diagnostic weighted bounds of $\td \uu, \na \td \uu$ in $\cL$ \eqref{eq:lin} using 
$|| \om_1 \vp_1||_{L^{\infty}}$ and weighted $C^{1/2}$ norm of $\om_1$. Ideas in Sec \ref{sec:u_comp_idea}.
Identify region $D$ where $\td \uu, \na \td \uu$ are
not perturbative relative to local damping terms,
i.e. where \eqref{eq:linf_diag_decay}, \eqref{eq:PDE_diag} fails. 
};

\node[box, below=of finite_rank] (lin_stab){
5. {\bf Key linear stability} (Sec \ref{sec:linf_decay}, \ref{sec:EE_hol}): determine rank and coupling weight $\tau_1$ between $L^{\infty}$ and $C^{1/2}$ norms in $E_1$ \eqref{energy1} to close coupled weighted $L^{\infty}(\vp_i)$ estimate 
with $E_1$ \eqref{energy1} and weighted $C^{1/2}$ estimates with $E_2$ \eqref{energy2}
for linearized equation perturbed by operators $\hat \uu, \wh {\na \uu}$ \eqref{eq:lin_main}.
Further adjust $\vp_2, \vp_3$.
 };

 \node[smallbox, right=of lin_stab] (wg_u2) {
Compute weighted bound of $\uu - \hat \uu, \na \uu - \wh{\na \uu}$ using $L^{\infty}(\vp_1)$ and weighted $C^{1/2}$ norm of $\om_1$.
Ideas in Sec \ref{sec:u_comp_idea}.
};

\node[box, below=of lin_stab] (wg_linf2){
6. Weighted $L^{\infty}$ estimate with growing weight 
$\vp_{g, i}$ (Sec \ref{sec:linf_grow}) to bound $|| (\om_1, \eta_1, \xi_1) ||_{L^{\infty}}$
for nonlinear estimates. Determine $\vp_{g, i}, \tau_2, \mu_4$ in energy $E_3$ \eqref{energy3}.

};

\node[smallbox, right=of wg_linf2] (wg_u3){
Idea: Estimate damped system with a force \eqref{eq:damp_sys}.
};

\node[box, below=of wg_linf2] (non) {
7. Nonlinear estimates:
Estimate terms in energy $E_4$ \eqref{energy4} using ODEs (Sec \ref{sec:rank1}) and determine $\mu_{ij}, i \geq 5$. Estimate $\wh W_2$ using $E_4$ (Sec \ref{sec:W2}), and then estimate nonlinear and error terms 
(Sec \ref{sec:comb_vel_err},\ref{sec:non}).

};

\node[smallbox, right=of non] (W2_solu) {
Construct approximate space-time solution $\hat F_i$ and 
finite-rank parts $\wh W_2, \cR_i$ following Sec \ref{sec:decoup_modi}.  \\
Idea: Estimate a damped system with forcing \eqref{eq:damp_sys}.
};

\node[box, below=of non] (stab){
8. Nonlinear stability: determine required smallness of residual error $\bar \cF$ of the profile (refine profile construction if $\bar \cF$ is not small enough) and choose $E_*$ to 
establish the stability estimates by verifying \eqref{eq:PDE_nondiag}.

};

 \draw[->, thick] (profi) -- (wghol);
 \draw[->, thick] (wghol) -- (wglinf1);
 \draw[->, thick] (wglinf1) -- (finite_rank);
 \draw[->, thick] (finite_rank) -- (lin_stab);
 \draw[->, thick] (lin_stab) -- (wg_linf2);
 \draw[->, thick] (wg_linf2) -- (non);
 \draw[-, ultra thick] (wg_linf2)-- (wg_u3);
  \draw[->, thick] (non) -- (stab);
  \draw[->, ultra thick] (W2_solu) -- (non);

\draw[->, ultra thick] (wghol.east) -- (wg_u1.west);
\draw[->, ultra thick] (wghol.east) -- (wg_u2.west);
\draw[->, ultra thick] (wglinf1.east) -- (wg_u1.west);
\draw[->, ultra thick] (wglinf1.east) -- (wg_u2.west);
\draw[->, ultra thick] (wg_u1.west) -- (finite_rank.east);
\draw[->, ultra thick] (finite_rank.east) -- (wg_u2.west);
 \draw[->, ultra thick] (wg_u2.west) -- (lin_stab.east);
 \draw[->, ultra thick] (finite_rank.east) -- (W2_solu.west);

\end{tikzpicture}
      \captionsetup{width=0.97\textwidth} 
  \captionof{figure}{
Dependency of estimates, parameters, and weights.
An arrow $a\to b$ means that the choices of parameters, weights, or estimates
in box $b$ depend on those in box $a$.
   }
  \label{Fig:order}
\end{figure}

\vspace{0.05in}
\paragraph{\bf Linear stability factor}
The  top-order estimates in the most singular regime in Section~\ref{sec:hol_singu}, together with the constants $c_1,c_2$ in \eqref{eq:model_hol_sing3}, indicate the size of the damping available in the
weighted $C^{1/2}$ estimates with weights $\psi_i$. 
Guided by this benchmark,
we choose the target stability factor $\lambda_*$, the lower bound for
$\lambda$ in \eqref{eq:PDE_diag}, to be on the order of $\min_i c_i$:
\beq\label{eq:EE_target1}
 \lambda_* \in [0.035,0.08].
\eeq

\begin{remark}
Using methods in  Sections \ref{sec:outgo} and \ref{sec:analy_para}, we could obtain a larger linear stability factor by generating stronger damping and increasing the approximation rank. 
Yet, we require weights that are smaller for large $|x|$ or decay faster at infinity, for instance by choosing a larger $m$ or a more negative $\b$ in the weight
$\vp_s$ in Section~\ref{sec:outgo}. This increases the constants in nonlinear
estimates through embeddings such as $ |f| \leq \vp^{-1}\|f\vp\|_{L^\infty}$, and thus requires a smaller residual error to close the nonlinear stability estimates. Thus, rather than maximizing the linear stability factor, 
we  balance the linear damping, nonlinear estimates, and
computational cost.
\footnote{
One can also obtain a larger stability factor by performing the estimates in a
space with a more singular weight $|\xx|^{-k}$ near $\xx=0$. This is analogous
to performing $H^k$ estimates with larger $k$ in
\cite{merle2022implosion2,chen2024Euler,chen2024vorticity}. However, this would
require the residual error to be small in a more singular weighted norm, which
is very challenging.
}
\end{remark}

\subsubsection{Choosing H\"older weights}\label{sec:wg_hol_g}

We use $\vp_{\bullet}$ to denote weights for $L^{\infty}$ estimates 
and $\psi_{\bullet}$ for $C^{1/2}$ estimates. 
In the energy estimates, we choose the H\"older weights $\psi_{\bullet}$ first since they determine the top-order
weighted $C^{1/2}$ estimates and set the benchmark for the range of stability factor used later. The $L^\infty$ weights will be motivated and chosen later, see  \eqref{eq:wg_linf_decay},
\eqref{eq:wg_linf_grow}. We choose
\beq\label{eq:wg_hol}
\bal
 \psi_{1}  & = p_{11} |x|^{-2} + p_{12} |x|^{-1/2} + p_{13} |x|^{-1/6},  \\
 \quad \psi_{2} & = p_{2 1} |x|^{-5/2} + p_{2 2} |x|^{-1} + p_{2 3} |x|^{-1/2} + p_{2 4} |x|^{1/6} , \quad  \psi_3(x) \equiv \psi_2(x),
 \eal
\eeq
where $p_{ij}$ are given in \eqref{wg:hol}. We first choose different powers to control the perturbation in both the near-field and the far-field, while obtaining damping factors in these two regions from the local parts of the equation. This choice, as well as the subsequent choice of the coefficients $p_{ij}$, is guided by the outgoing properties of the flow and the weighted H\"older estimates developed in Sections \ref{sec:outgo} and
\ref{sec:model1}.
Next, we apply the estimates in Section \ref{sec:hol_singu} to treat the nonlocal terms in \eqref{eq:bad} perturbatively. Finally, we adjust the coefficients so that \eqref{eq:model_hol_sing3} holds with $c_1, c_2$ as large as possible. 
 In this way, the top-order nonlocal terms are controlled perturbatively by the damping generated from the local terms and the outgoing flow.

\begin{remark}[\bf Notation $\phi$]
The reader should not confuse the weights $\psi_i, \vp_i$ with the notation for the stream function $\phi = (-\D)^{-1}\om$. In the energy estimates, we rarely use the stream function.
\end{remark}

Next, we determine the weights $g_i$ in the H\"older seminorm
$C_{g_i}^{1/2}$ 
\eqref{hol:g}. They play a different role from the spatial weights $\psi_i$ and define anisotropic variants of the
homogeneous $\dot C^{1/2}$ seminorm. We use Lemmas \ref{lem:holx_ux}--\ref{lem:holy_uy} and the triangle
inequality  for H\"older estimate of $\na \uu$:
\[
\bal
   g(x-z)|f(x) - f(z) |  \leq g(x-z) (|f(x_1, x_2) - f(z_1, x_2)| + | f(z_1, x_2) - f( z_1, z_2)| ) .
   \eal
\]

In general, a direct use of the triangle inequality may introduce an extra
factor $\sqrt{2}$.  For example, if $[f]_{C_x^{1/2}}\leq A$ and
$[f]_{C_y^{1/2}}\leq A$, then one obtains
$[f]_{C^{1/2}}\leq \sqrt{2}A$.  The choice
\[
    g(h_1,h_2)=\big(|h_1|^{1/2}+c|h_2|^{1/2}\big)^{-1}
\]
would remove this loss
\footnote{
In the same example, choosing $c=1$ gives
$g(x-z)|f(x)-f(z)|\leq A$.
}
but it leads to an unbounded contribution from the damping factor
\beq\label{eq:wg_hol_gtry}
 d_g = g(x-z)^{-1}(b(x) - b(z)) \cdot (\nabla g)(x-z)
 \eeq
since $|\partial_1 g(0,h_2)|=\infty$ and
$|\partial_2 g(h_1,0)|=\infty$.  We therefore use the regularized
anisotropic weight
\beq\label{eq:wg_hol_g}
 g_{i}(h) = g_{i0}(h) g_{i0}(1, 0)^{-1} ,  \quad 
 g_{i0}(h) =  ( \sqrt{ |h_1| + q_{i1} |h_2| } + q_{i3} \sqrt{ |h_2| + q_{i2} |h_1| }   )^{-1}, 
\eeq
for some small $q_{i1},q_{i2} >0 $. 
We divide $g_{i0}(1, 0)$ to normalize $g_i(1, 0)= 1$. To exploit the anisotropy of the flow, see Section \ref{sec:aniso_est}, we choose
$q_{i3}<1$.  The parameters $q_{ij}$ are given in \eqref{wg:hol}.

For general pairs $(x,z)$, estimating
$g(x-z)(\na\uu(x)-\na\uu(z))$ still gives a larger constant than in the
one-dimensional cases $x_1=z_1$ or $x_2=z_2$.  However, when $\frac{|x_2-z_2|}{|x_1-z_1|}$
is not too small, the damping contribution $d_g$ in
\eqref{eq:wg_hol_gtry} becomes stronger.  This extra damping compensates for
the larger constant in the full H\"older estimate and allows the nonlocal
terms to remain perturbative.

For $\eta$ and $\xi$, we set $g_3(h)=g_2(h)$.  To choose the parameters
$q_{ij}$ in \eqref{eq:wg_hol_g}, we first find the point
$x_*\in \R^2_{++}$ where the damping is weakest for sufficiently small
$|x-z|$ with $x_2=z_2$; equivalently, $x_*$ corresponds to the maximizer of the
left-hand side of \eqref{eq:model_hol_sing3}.  Near $x_*$, we then perform the
H\"older estimates for other values of the ratio
$|x_2-z_2|/|x_1-z_1|$, while keeping $|x-z|$ small.  As in Section
\ref{sec:hol_singu}, the more regular nonlocal terms and the contributions
from the approximate space-time solution $\hat F_i$  vanish in this local leading-order
analysis.  Thus the choice of $q_{ij}$ is determined by 
the approximate profile $(\bar \om, \bar \th)$ and the local damping
mechanism alone.  We choose $q_{ij}$ so that this damping remains at least as
strong as in the case $x_2=z_2$.

\subsection{Ideas of estimating the nonlocal terms}\label{sec:u_comp_idea}

In the energy estimates, we need weighted $L^\infty$ and $C^{1/2}$ estimates for the velocity remainders $\uu_A,\uu_{x,A},\uu_{y,A}$ in \eqref{eq:vel_rem},\eqref{eq:lin_main} assuming
$\om_1\vp_1,\om_1\vp_{g,1}\in L^\infty$ and $\om_1\psi_1\in C^{1/2}$ for suitable weights
$\vp_1, \vp_{g, 1}, \psi_1$.  This subsection discusses the
main ideas and methods behind these estimates; the full implementation,
including the rigorous computation, is given in Part~II \cite{ChenHou2023b}.
The main point is that the nonlocal terms can be reduced to explicit weighted
kernel integrals, while the genuinely singular parts are controlled by the
weighted $C^{1/2}$ norm. The constants in these estimates 
depend only on the weights, not on the perturbation.

We consider $\om_1$ odd in $x$. For $f = \uu_A, (\na \uu)_A$, the nonlocal velocity can be written as
\[
 I(f)(x) = \int_{ \R^2} K_{f}(x, y) \Om_1(y) dy ,
\]
for some kernel $K_f$, where $\Om_1$ is the odd extension of $\om_1$ in $y$ from $\R^2_{+}$ to $\R^2$ \eqref{eq:ext_w_odd}. 
The formulas of $\na \uu$ are given in \eqref{eq:ux_local}. 
When $f= (\na \uu)_A$, integrals are understood in the principal-value sense.

We first consider $f = \uu_A$. The kernel involves $\na^{\perp} \log |x-y|$ and has 
a locally integrable singularity of order $|x-y|^{-1}$ near $y=x$. To obtain a sharp 
weighted estimate of $  \uu_A$ with some weight $\rho$ singular  near $\xx=0$, 
we exploit the oddness of $\Om_1$ in $y_1, y_2$, symmetrize the kernel 
\bseq\label{eq:u_comp_idea1}
\beq
 K_f^{sym}(x, y)  = K_f(x,  y ) - K_f(x, -y_1, y_2) - K_f(x, y_1, -y_2) + K_f(x, - y), \\
 \eeq
 and then apply the $L^{\inf}$ estimate 
 \begin{align}
| \rho(x) I(\uu_A)(x)| & = \rho(x)\B|  \int_{\R^2_{++}} K_{\uu_A}^{sym}(x,y) \om_1(y) dy \B| 
\leq  \rho(x)\int_{\R^2_{++}} | K_{\uu_A}^{sym}(x, y)| \vp_1^{-1}(y) dy 
\cdot  || \om_1 \vp_1 ||_{L^{\inf}} \notag \\
&  \teq \rho(x) C(\uu, x) || \om_1 \vp_1||_{L^{\inf}} ,
\end{align}
\eseq
where $C(\uu, x)$ denotes the last integral on the second line. This estimate is sharp in the following sense: for fixed $x$, equality can be achieved by $\om_1 \vp_1(y) = C \cdot \sgn(K_{\uu_A}^{sym}(x, y))$ for some constant $C$. 
Thus, once the weight $\vp_1$ is fixed, the bound $C(\uu, x)$ 
is an integral of explicit kernels and weights, and is independent of the perturbation $\om_1$. 
We can estimate it effectively for {\textit{all}} $x$ using the scaling symmetry of the kernel and 
numerical integration with rigorous error estimates.

For $f = (\na \uu)_A$, 
the kernel has a nonintegrable singularity of order $|x-y|^{-2}$ near $y=x$.
We decompose the integral $I(f)$ into the nonsingular part ($\mf{NS}$) and the singular part (S) with singular region $R$ centered around $x$ with radius $r(x)$ 
\beq\label{eq:u_comp_idea2}
\bal
I(f) &= I_{ \mf{NS}} (f) + I_S(f),  \quad  R(x) = \{ y: \max_{i=1,2} |x_i - y_i| \leq r(x )\}. \\
 I_{\mf{NS} }(f) & = \int_{ R^c} K_f(x, y) \Om_1(y) dy , \quad 
I_S(f) = \int_{ R} K_f(x, y) \Om_1(y) dy .
\eal
\eeq

The nonsingular part is estimated by the weighted $L^{\infty}$ method above.
For the singular part, 
we use the commutator decomposition,
e.g., \eqref{eq:ux_commu1}. We apply the above $L^{\inf}$ estimate \eqref{eq:u_comp_idea1} to the regular term. The singular term related to $\na \uu(\om_1 \psi_1)$ is estimated using $ || \om_1 \psi_1||_{C^{1/2}}$. For example, we have the following estimate 
\[
\bal
\B| \int_{ |s_1|, |s_2| \leq \tau  }  \f{ s_1 s_2}{|s|^4} (\om_1 \psi_1)(x-s) ds \B| 
= & \B| \int_{  0\leq s_1 , |s_2| \leq \tau }  \f{ s_1 s_2}{|s|^4} ( (\om_1 \psi_1)(x-s) - (\om_1 \psi_1)(x_1+s_1, x_2 - s_2 ) d s \B| \\
\leq & C \tau^{1/2}   [\om_1 \psi_1]_{C_x^{1/2}} , 
\eal
 \] 
 where $C$ is independent of $\tau$.  Thus, for a singular weight $\rho$,  we estimate $\rho I(f)$  as follows 
\footnote{
The odd extension $\Om_1$ of $\om_1$ is not continuous across $y=0$. 
 We carefully use the $[\om_1 \psi_1]_{C_y^{1/2}}$ norm  to control the integral in the singular region. See details in Supplementary material II of Part II \cite[Section 6.1]{ChenHou2023bSupp}. 
}
\beq\label{eq:u_comp_idea3}
| \rho(x) I(f)(x)| \leq  \CCs_1(x) || \om_1 \vp_1 ||_{L^{\inf}}
+  \CCs_2(x) [ \om_1 \psi_1]_{C_x^{1/2}} + \CCs_3(x) [ \om_1 \psi_1 ]_{C_y^{1/2}},
\quad 
\eeq
for explicit computable functions $\CCs_i(x)$. We further bound $\om_1$-norms using energy $E_1$ \eqref{energy1}.

\subsubsection{Ideas of $C^{1/2}$ estimates}\label{sec:idea_int_hol}

The weighted $C^{1/2}$ estimate is more involved. 
We outline the ideas below and refer full details to Part II \cite[Section 4]{ChenHou2023b}. 
The basic strategy is to isolate the most singular part, estimate it by the
sharp H\"older estimates from Section~\ref{sec:sharp}, and treat the remaining
more regular part by optimizing two estimates.

We focus on the estimate of the more singular term $(\na \uu)_A$ \eqref{eq:vel_rem}.  
To estimate $(\na \uu)_A(\om) \psi$ with a weight $\psi$, we decompose it into the regular and singular parts following \eqref{eq:ux_commu0}, \eqref{eq:ux_commu1}, e.g.
\[
\bga
    u_{x, A}(\om) \psi = u_x(\om) \psi - \hat u_x(\om) \psi
   = I_S  + I , \\
   I_S  =   u_x( \om \psi)(x, a, b),\quad 
  I =  [u_x(\cdot, a, b), \psi](\om)  + (u_x(\om) - u_x(\om)(x, a, b)) \psi - \hat u_x \psi
  \ega
\]
where $[u_x(\cdot, a, b), \psi](\om)$ is the commutator defined in \eqref{eq:ux_commu0}. For the singular part $I_S$, we use the sharp H\"older estimates in Lemmas \ref{lem:holx_ux}-\ref{lem:holy_uy}.

It remains to estimate the more regular part $I$. 
Due to the commutator and removal of the singularity, 
this term is Log-Lipschitz and is similar to $J(\om)(s, 0)$ below with $a=0$
\[
J(\om)(s, a) = \int_{ a \leq \max_i |s_i - y_i|  } K( s, y) \om(y) dy,   \quad  |K(s, y) | \leq C_1 |s-y|^{-1},
\quad |\pa K(s, y)| \leq C_2 |s-y|^{-2} .
\]
 For $x, z \in \R_2^{++} $,  we further decompose 
 \[
       I(s) = I_R( s, a) + I_S(s, a).
 \]
Here, $I_R( s, a)$ is the regular part of $I$ with a distance $|a|$ away from the singularity and is similar to $J(f)(s, a)$, while $I_S(s, a)$ is the singular part similar to $J(f)(s,0) - J(f)(s, a)$. 
The regular part $I_R( s, a)$ is estimated through its piecewise Lipschitz norm, which is 
bounded by the norm $|| \om \vp||_{\inf}$ and is reduced to 
estimating explicit integrals depending on the weight $\vp$. 
See the previous part in Section \ref{sec:u_comp_idea} for more details.  For $I_S(s)$, we estimate its piecewise $L^{\inf}$ norm by estimating integrals depending on the weights. This gives bounds similar to the following
\[
\bal
 |\pa I_R(s, a)| & \leq   C_1(s)\log (C_3 a^{-1}) || \om \vp||_{\inf} , \quad   |I_S(s, a)| \leq C_2(s) |a| \cdot || \om \vp||_{L^{\inf}} ,  \\
 \f{ |I(x) - I(z)|}{|x-z|^{1/2}} & \leq  \f{  |I_R(x, a) - I_R(z, a)| + |I_S(x, a) - I_S(z, a)|}{ |x-z|^{ 1/2} }   
 \\
 &  \leq \B( C_4(x, z) \log ( \f{C_3}{a}  ) |x-z|^{ \f{1}{2} } + \f{ C_5(x, z) |a| }{  |x-z|^{1/2}}  \B) || \om \vp ||_{L^{\inf}}.
 \eal
\]

The sharp estimate is obtained by $a \sim |x-z|$, which balances two terms in the upper bound.
In practice, we perform several decompositions with different $a$ and optimize the resulting estimates.
Since we further bound 
the $L^{\infty}$ norm $|| \om_1 \vp_1||_{L^{\inf}}$ and 
$C^{1/2}$ seminorm 
 $[ \om_1 \psi_1 ]_{C_{x_i}^{1/2}}$ by  energy $E_1$ \eqref{energy1}, and these (semi)norms contribute to $E_1$ with different weights, we also use a small portion of the H\"older norm to improve the estimate of $I$, and then optimize among these bounds.

\subsubsection{Scaling symmetry and rescaled integral}\label{sec:scal_sym}

It remains to compute the explicit weighted kernel integrals in the above estimates. 
This is complicated by two singularities.
First, the weight $\rho(x)$ is singular near $0$, which significantly amplifies the 
 error in computing the integrals, e.g. \eqref{eq:u_comp_idea1}. 
Second, the kernel $K_{\uu}(x,y), K_{\na \uu}$ are singular near $y = x$. If there are only a few $x$, one can design a mesh adapted to the singularity $y = x$ and then apply standard quadrature. However, 
we need estimates of the integrals for \emph{all} $x$.  The key observation is that the kernels enjoy a scaling symmetry, which allows us to confine the relevant singularity to an $O(1)$ region.

Denote $f_{\lam}(x) \teq f(\lam x)$. We consider the kernels about $\na u$, which are singular of order $2$:
 \[
 K(\lam x, \lam y) = \lam^{-2} K(x, y)
 \]
For a rescaling parameter $\lam$ to be chosen, 
taking $y = \lam \hat{y}, x = \lam \hat{x}$, we get 
\[
 \rho(x) \int_{\R_2^{++}} K^{sym}(x, y) \om(y) dy 
= \rho(\lam \hat{x}) \int_{\R_2^{++} } K^{sym}( \lam \hat{x}, \lam \hat{y}) \om( \lam \hat{y})  \lam^2 d \hat{y} 
= \rho_{\lam}(\hat{x})\int_{\R_2^{++} } K^{sym}(  \hat{x},  \hat{y}) \om_{\lam}(  \hat{y})   d \hat{y} .
\]

Applying the $L^{\infty}$ estimate gives
\[
| \rho(x) \int_{\R_2^{++} } K^{sym}(x, y) \om(y) dy | 
\leq || \om_{\lam} \vp_{  \lam}||_{L^{\infty}} \rho_{\lam}(\hat{x})
\int_{\R_2^{++} }  | K^{sym}(\hat{x}, \hat{y}) | \vp_{ \lam}(\hat{y})^{-1} d\hat{y}.
 \]
Since $|| \om_{\lam} \vp_{\lam}||_{L^{\infty}} = || \om \vp ||_{L^{\infty}}$, it suffices to compute the rescaled integral. 
By choosing $\lam$ so that $|\hat x|\asymp 1$, both the rescaled point $\hat x$ and the kernel
singularity $\hat y=\hat x$ lie in an $O(1)$ region. This allows us to use an adaptive mesh concentrated in a fixed bounded region, instead of re-meshing across many scales of $x$. 
Moreover,  $\hat{x}$ can be chosen away from $0$, e.g. $|\hat x| \asymp 1$, so that $\rho_{\lam}(\hat{x})$ is nonsingular. For example, $|x|^{-2} = \lam^{-2} |\hat{x}|^{-2}$ by taking $\lam = |x| / |\hat{x}|$ with $|\hat{x}|\asymp 1$. 
This rescaling method is the key device for overcoming the singularities 
of weights and kernels  in the estimate of integrals, and it reduces a continuum of singular integral estimates at many scales to a finite family of estimates on normalized configurations. See Part II \cite[Section 4]{ChenHou2023b} for details.

\vs{0.05in}
\paragraph{\bf Weighted $L^{\infty}$ and $C^{1/2}$ estimates}

Let $\hat \uu, \wh {\na \uu}$ be the finite-rank approximation constructed in 
\eqref{eq:u_appr} and $\uu_A, (\na \uu)_A$ be the associated velocity remainder \eqref{eq:vel_rem}. 
Let $\psi_u , \psi_1$ be the weights defined in \eqref{wg:hol}, \eqref{wg:elli}. Using the above methods, 
for $\om$ odd in $x_1$, we establish weighted 
directional $C^{1/2}$ estimates 
for $\psi_u \uu_A , \psi_1 (\na \uu)_A$ in the following form
\footnote{
Since $\uu, \na \uu$ have different scaling and asymptotics, we apply $\psi_u$ to estimate $\uu_A$ 
and $\psi_1$ to $(\na \uu)_A$.
}
\footnote{
Since $\uu$ is about one-order more regular than $\om$, we do not need 
the $C^{1/2}$-norm of $\om$ in the upper bound \eqref{eq:holest:u}.
}
\bseq\label{eq:holest}
\begin{align}
   \f{ | \d(\psi_u f_A , x, z ) | }{ |x-z|^{1/2}} 
  & \leq \CCs_{ 1/2, f, i}(x, z, \vp,  \psi_u)   || \om \vp||_{\inf} ,
  \hspace{1.85in} f = u, v , 
  \label{eq:holest:u} \\
   \f{ | \d(\psi_1 f_A, x, z ) | }{ |x-z|^{1/2}} 
  & \leq \CCs_{ 1/2, f, i}(x, z, \vp, \psi_1,  \g) \max(  \| \om \vp \|_{\inf} , 
   \max_{j=1,2} \g_j  [ \om \psi_1]_{C_{x_j}^{ 1/2 }(\R^2_{+}) } ), 
   \  f = \pa u, \pa v, 
   \label{eq:holest:du}
\end{align}
\eseq
with $\pa = \pa_x, \pa_y$, $x_1 = z_1$ for $i=2$ and $x_2 = z_2$ for $i=1$. Here, the 
coefficient functions $\CCs_{1/2,f, i}(x, z, \bullet)$ 
are \emph{computable} and bounded in any compact domain of $\R_2^{++}$, depend on $\g$, weights, and the approximations. Using $v_y = -u_x$, we do not need extra  estimates for $v_y$.

We choose the parameter $\g_j$ in \eqref{eq:holest:du} according to the coupling weight in $E_1$ \eqref{energy1}
\[
  \g_1 =  \tau_1^{-1}, \quad \g_2 = \tau_1^{-1} g_1(0,1) . 
\] 
If in addition  $\vp = \vp_1$, using the definition \eqref{hol:g} and \eqref{eq:hol_to_E1},
we bound the norm in \eqref{eq:holest:du}: 
\beq\label{eq:holest2}
\max(  || \om \vp_1||_{\inf} , \max_{j=1,2} \g_j  [ \om \psi_1]_{C_{x_j}^{1/2}(\R_2^{+}) } )  
\leq \max(  || \om \vp_1||_{\inf} , \tau_1^{-1} [ \om \psi_1]_{C_{g_1}^{1/2}} )
\leq E_1.
\eeq
Since we fix $\gamma$ and $\psi_1$, when no confusion arises, we simplify 
$\CCs_{ 1/2, f, i}(x, z, \vp,  \psi_1,\g)$ as  
$\CCs_{ 1/2, f, i}(x, z, \vp)$.

As in \eqref{eq:u_comp_idea3}, we establish weighted $L^{\infty}$ estimates in the following form 
\beq\label{eq:idea_linf}
 |\rho_{ij} f_{ij}(\om_1)(x)| \leq \CCs_{ij, 1}(x ; \vp) || \om_1 \vp||_{\inf} 
 + \CCs_{ij, 2}(x ; \vp) [ \om_1 \psi_1]_{C_x^{1/2}}
 + \CCs_{ij, 3}(x ; \vp ) [ \om_1 \psi_1]_{C_y^{1/2}},
\eeq
with explicit computable coefficient functions. 
For the perturbation, the $\om_1$-norms are further bounded by energy $E_1$ \eqref{energy1} 
or energy $E_3$ \eqref{energy3}. Here $\vp = \vp_1, \vp_{g, 1}, \vp_{elli}$ 
 are defined in \eqref{wg:linf_decay}, \eqref{wg:elli},  $\rho_{ij}$ is given in  \eqref{wg:elli}, and $f_{ij}$ denotes $f_{01} = u_A, f_{10} = v_A,  f_{11} = u_{x, A},  f_{20} = v_{x, A},   f_{0 2} = u_{y, A} $. 
See more discussion of estimates \eqref{eq:idea_linf} with $\vp = \vp_1, \vp_{g, 1}$ to Sections \ref{sec:linf_decay}, \ref{sec:linf_grow}, respectively.

In summary, the output of these integral estimates 
 are  the weighted $L^{\infty}$ \eqref{eq:idea_linf} and $C^{1/2}$ estimates 
  \eqref{eq:holest} of velocity remainder $\uu_A, (\nabla \uu)_A$  by weighted $L^{\infty}, C^{1/2}$ norms of $\om_1$ 
with explicit computable coefficient functions. 
These estimates are then inserted into the weighted energy estimates in Sections~\ref{sec:linf_decay}--\ref{sec:linf_grow}. See full computation and estimates details  in \cite[Section 4]{ChenHou2023b}.

\subsection{Weighted $L^{\inf}$ estimate with decaying weights}\label{sec:linf_decay}

We first establish weighted $L^\infty$ estimates using decaying weights
$\vp_i$. They provide stronger damping in the
far field (corresponds to more negative $\b_i$ in \eqref{eq:idea_damp_far}) 
and make the linear stability  easier to establish. We choose
\footnote{
From Lemma \ref{lem:W1_reg} and \eqref{eq:Sym_W12}, for $W_1 = (\om_1, \eta_1, \xi_1)$, we can divide $W_1$ by $|x_1|^{1/2}$ and $ || W_{1, i} \vp_i||_{L^{\infty}}$ is bounded.
}
\beq\label{eq:wg_linf_decay}
\bal
\vp_1 & = ( p_{41} |x|^{-2.4} + p_{42} |x|^{-1/2} ) |x_1|^{-1/2} + p_{43} |x|^{-1/6},  \ \qquad \qquad \qquad \vp_4 =  \psi_1 |x_1|^{-1/2} , \\
\vp_{i-3}& = ( p_{i1}|x|^{-5/2} + p_{i2} |x|^{-3/2} + p_{i3} |x|^{-1/6} ) |x_1|^{-1/2}
+ p_{i4} |x|^{-1/4} + p_{i5} |x|^{1/7} , \  i=5, 6 ,\\
\eal
\eeq
with parameters $p_{ij}$ given in \eqref{wg:linf_decay}. 
 See the weighted $L^{\inf}$ estimate in Sections \ref{sec:outgo}, \ref{sec:model1} for more motivations.  Recall $\psi_i$ from \eqref{eq:wg_hol} for weighted $C^{1/2}$ estimates. We apply $\psi_1 , \vp_1$ for $\om$, $\psi_2, \vp_2$ for $\eta$, and $\psi_3, \vp_3$ for $\xi$.
The additional weight $\vp_4$ is used in Section \ref{sec:hol_Q1} 
for the estimate of $\| \vp_4 \om_1 \|_{L^{\infty}}$ needed for the H\"older estimate.
We explain the choice of  $\vp_i$ in Section \ref{sec:wg_linf_decay}.

Before deriving the weighted $L^\infty$ inequalities, we first record the
nonlocal estimates in a form that can be inserted directly into the energy estimates. 
The following estimates bridge the nonlocal estimates in Section~\ref{sec:u_comp_idea} and the weighted \(L^\infty\) energy estimates. Using the weights $\vp_1, \psi_1$ and following the methods in Section~\ref{sec:u_comp_idea}, 
we obtain computable coefficient functions $\CCs_{ij,k}(x)$ for the kernel integral estimates of $\uu_A, (\na \uu)_A$ defined in \eqref{eq:vel_rem}. 
 In particular, we have
\bseq\label{eq:u_linf_est1}
\beq
 |\rho_{ij} f_{ij}(\om_1)(x)| \leq \CCs_{ij, 1}(x) || \om_1 \vp_1||_{\inf} 
 + \CCs_{ij, 2}(x) [ \om_1 \psi_1]_{C_x^{1/2}}
 + \CCs_{ij, 3}(x) [ \om_1 \psi_1]_{C_y^{1/2}},
\eeq
where $f_{ij}(\om_1)$ denotes different velocity terms obtained from $\om_1$: \footnote{
This indexing is consistent with the stream function $\phi$ representation 
\[
 u = - \pa_y \phi, \quad  v = \pa_x \phi , \quad  u_x = - \pa_{xy }\phi,  
 \quad v_x = \pa_{x x} \phi ,  \quad u_y = - \pa_{yy} \phi ;
 \]
the indices $i,j$ record the number of $x$ and $y$ derivatives on $\phi$.
}
\beq
 f_{01} = u_A, \quad f_{10} = v_A, \quad f_{11} = u_{x, A}, \quad 
 f_{20} = v_{x, A},  \quad f_{0 2} = u_{y, A} .
\eeq
\eseq

The functions $\CCs_{ij,k}(x)$ \emph{encodes} the piecewise bounds in the kernel integral estimates 
and are bookkeeping devices. They only depend on $\vp_1, \psi_1$, the kernels,  and the finite-rank approximation $\hat{\uu}, \wh {\na \uu}$, but \emph{not on} $\om_1$. We use the \emph{sans-serif font} $\CCs_{\bullet}$ to distinguish these functions from generic constants.  We refer more discussion to Section~\ref{sec:piece_bd}.

We use the weight $\rho_{ij}$ to capture the vanishing order near $0$ and decays of $\uu_A, \na \uu_A$
\beq\label{wg:u}
 \rho_{10} = \rho_{01} = |x|^{-3} + |x|^{-7/6}, \quad \rho_{20} = \rho_{11} = \rho_{0 2} = \psi_1 .
\eeq
See also definitions of these weights in Appendix \ref{app:para_wg}. For $ (\na \uu)_A$, we choose $\rho_{ij} = \psi_1, i+j=2$, since we need to estimate $(\na \uu)_A \psi_1$ using the H\"older norm of $\om_1 \psi_1$ and $\na \uu$ and $\om_1$ are of the same order. To control $\uu_A$, we do not need to use the H\"older seminorm. Thus 
\beq\label{eq:u_linf_est11}
\CCs_{ij, 2}(x) =  \CCs_{ij, 3}(x) = 0, \quad i+j =1.
\eeq

\paragraph{\bf Singular factor $|x_1|^{-1/2}$} 
There is one additional point caused by the factor $|x_1|^{-1/2}$ in \eqref{eq:wg_linf_decay}.  Note that $\vp_i$ contains the singular term $|x_1|^{-1/2}$ \eqref{eq:wg_linf_decay} and $\bar \om_x, \bar \th_{xx}$ \emph{do not} vanish on $x_1 = 0$. 
To bound $u_A \bar \om_x \vp_1, u_A \bar \th_{xx} \vp_2$ in the energy estimate of \eqref{eq:lin_main}, since $\om_1, u$ are odd in $x_1$, we use the odd symmetry of $u_A$ in $x_1$ and $u_A(0, x_2) = 0$ to absorb the singularity. Since $u_A$ is almost $1$ order more regular than $\om_1$, we can develop an estimate 
\[
|u_A(x)| \leq C(\vp_1)(x) \cdot \| \om_1 \vp_1 \|_{\infty} ,\quad
C(\vp_1)(x) \les |x_1|  \cdot ( 1 + | \log |x_1| | ) C_2(\vp_1)(x),
\]
with locally bounded coefficient $C_2(\vp_1)(x)$. Thus, $\rho_4 u_A$ can be bounded by $|| \om_1 \vp_1 ||_{\infty}$ locally, where $\rho_4$ given in \eqref{wg:elli} captures the singular weight $|x_1|^{-1/2}$.  
This estimate follows similar methods in Section \ref{sec:u_comp_idea}; 
see Part II \cite[Appendix {\secudx}]{ChenHou2023b}. 
We optimize this estimate and \eqref{eq:u_linf_est1} for $u_A$. 
 Then the coefficient $\CCs_{01, 1}(x)$ bounding $u_A \rho_{10} = u_A \rho_4 \cdot \f{\rho_{10}}{\rho_4}$ vanishes along $x_1= 0$.

\subsubsection{Piecewise upper bounds}\label{sec:piece_bd}

We first obtain the nonlocal estimates with computable pointwise upper bounds 
for $\xx$ in each mesh cell. We discretize a very large domain $[0, D]^2$ in $\R_2^{++}$ using the same mesh $y_i$ in Section \ref{sec:ASS_rep} for computing the profiles. Using the method in Section {\secvelcomp} in Part II \cite{ChenHou2023b}, we can obtain piecewise bounds of $ \uu_A,  (\na \uu)_A$ and $\rho_{10} \uu_A, \rho_{20} (\na \uu)_A$ in each grid $[y_i, y_{i+1}] \times [y_j, y_{j+1}]$. In particular, the coefficients $ \CCs_{ij, 1}(\xx), \CCs_{ij, 2}(\xx), \CCs_{ij, 3}(\xx)$ in the upper bound \eqref{eq:u_linf_est1} are piecewise constants. We store these bounds as $n \times n$ matrices. 

The estimates in the far-field $ x \notin [0, D]^2$ are much easier and very small,  
since the coefficients of nonlocal terms in \eqref{eq:lin_main} have fast decay 
and we choose $D$ very large. 
We estimate the asymptotics bounds for $\uu_A, (\na \uu)_A$ 
using a scaling argument similar to that in Section \ref{sec:scal_sym}. 
See Supplementary Material II in Part II \cite[Section 7]{ChenHou2023bSupp}. The same ideas apply to all other estimates.

\vs{0.1in}
\paragraph{\bf{Bookkeeping notations for nonlocal estimates}}

The weighted estimates below contain several repeated combinations
of the same piecewise coefficient bounds $\CCs_{jk, i}$.  To avoid obscuring the structure of
the estimates with long expressions, we introduce the following bookkeeping
notation.  These definitions \emph{do not introduce new estimates}; they only package
the nonlocal estimates and coefficients $\CCs_{kl , i }(x), i \in \{ 1,2,3\}$ in \eqref{eq:u_linf_est1} following Section \ref{sec:u_comp_idea}.  We define functions 
\bseq\label{eq:EE_oper}
\begin{align}
\cT_{u}(f)   & \teq  \CCs_{01, 1} |f_x| + \CCs_{10, 1} |f_y|  ,  \label{eq:EE_oper:a}  \\ 
\CCs(\uu_x, \na \bar \th,  i ) & \teq \CCs_{11, i } | \bar \th_x| + \CCs_{20, i} |\bar \th_y|,  
\qquad  
\CCs(\uu_y, \na \bar \th, i)   \teq \CCs_{02, i} |\bar \th_x| + \CCs_{11, i} |\bar \th_y| ,   \
 1\leq i \leq 3, \label{eq:EE_oper:b} \\
\CCs( f, \na \bar \th, a )  & \teq \CCs(f,  \na \bar \th, 1) a_1 + \CCs(f, \na \bar \th, 2) a_2 + \CCs(f, \na \bar \th, 3) a_3 ,  \quad \ a  = (a_1, a_2, a_3) \in \R^3_+,
\label{eq:EE_oper:c}
\end{align}
where the vector input $f = \uu_x$ or $f = \uu_y$. We introduce 
\beq\label{eq:EE_oper:d}
\CCs(f_{ij}, a ) = \CCs_{ij, 1} a_1 + \CCs_{ij, 2} a_2 + \CCs_{ij, 3} a_3, 
\quad  a = (a_1, a_2, a_3), \  i + j = 1, 2 ,
\eeq
\eseq
where $f_{ij}$ denotes different velocity term  \eqref{eq:u_linf_est1}.
We use \emph{sans-serif font} $\CCs$ to distinguish $\CCs(\uu_x, \na \bar \th, i )$ from 
generic functions and constants. We use $\cT_{u}$ to estimate transport-type terms: $\uu_A \cdot\na f$. The calligraphic $T$ is short for \emph{transport}. We use $\CCs(\uu_x, \na \bar \th, i)$ 
and $\CCs(\uu_y, \na \bar \th, i)$ to estimate  $\uu_{x, A} \cdot \na \bar \th$ and $\uu_{y, A} \cdot\na  \bar \th$, respectively. 
The functions $\CCs(\bullet, \na \bar \th, i)$  with $i =1,2,3$ in \eqref{eq:EE_oper:b} track the bounds associated 
with  $ || \om_1 \vp_1||_{\inf} $, $[ \om_1 \psi_1]_{C_x^{1/2}}$, $ [ \om_1 \psi_1]_{C_y^{1/2}}$ in nonlocal estimates \eqref{eq:u_linf_est1}, respectively.
We use weight $a$ in \eqref{eq:EE_oper} to encode relative sizes of $3$ different norms 
and $\CCs(f_{ij}, a), \CCs(f, \na \bar \th, a) $ to combine their corresponding coefficient bounds.
Note that $\CCs(f_{ij}, a), \CCs(f, \na \bar \th, a) $ are linear in $a$.

\vs{0.05in}

\paragraph{\bf $L^{\infty}(\vp_i)$ estimates}

 As in \eqref{eq:lin_sim}, we perform weighted $L^{\inf}$ estimates for $W_{1} =( \om_1, \eta_1, \xi_1)$ in \eqref{eq:lin_main}
\beq\label{eq:linf_decay1}
  \pa_t ( W_{1, i} \vp_i) + (\bar c_l x + \bar \uu ) \cdot \na ( W_{1, i} \vp_i) 
 =  d_{i,L}(\vp_i) W_{1, i} \vp_i + B_i(x) \vp_i  \, , 
\eeq
where we used $d_{i, L}(\cdot)$ in \eqref{eq:dp}  to denote the  damping  coefficient, and $B_i(x)$ are the bad terms defined in \eqref{eq:bad}. 
We estimate $B_i(x)$  using the above pointwise estimates for the nonlocal terms
\beq\label{eq:linf_decay2}
\bal
 | (B_1 \vp_1) (x)|  & \leq \f{\vp_1}{\vp_2} || \eta_1 \vp_2||_{\inf}
+  \f{\vp_1}{\rho_{10}} \cT_u( \bar \om) || \om_1 \vp_1 ||_{\inf},   \\
 | (B_2\vp_2)(x)|  & \leq \f{\vp_2}{\vp_3} |\bar v_x| \cdot || \xi_1 \vp_3 ||_{\inf}
+ \f{\vp_2}{\rho_{10}} \cT_{u}(\bar \th_x) || \om_1 \vp_1||_{\inf} \\
  + \f{\vp_2}{\psi_1} & ( \CCs(\uu_x, \na \bar \th, 1) || \om_1 \vp_1||_{\inf} 
 + \CCs(\uu_x, \na \bar \th,  2) [\om_1 \psi_1]_{C_x^{1/2}} 
+ \CCs(\uu_x,  \na \bar \th, 3)  [\om_1 \psi_1]_{C_y^{1/2}} ),  \\
 | (B_3 \vp_3)(x)|  & \leq \f{\vp_3}{\vp_2} |\bar u_y| || \eta_1 \vp_2 ||_{\inf}
+ \f{\vp_3}{\rho_{10}} \cT_{u}(\bar \th_y) || \om_1 \vp_1||_{\inf} \\
  + \f{\vp_3}{\psi_1} & ( \CCs(\uu_y, \na \bar \th, 1) || \om_1 \vp_1||_{\inf} + \CCs(\uu_y, \na \bar \th , 2) [\om_1 \psi_1]_{C_x^{1/2}}
+ \CCs(\uu_y, \na \bar \th,  3)  [\om_1 \psi_1]_{C_y^{1/2}} ).
\eal
\eeq

  We use $B_2$ to illustrate the above weighted estimates. We estimate $B_2$ defined in \eqref{eq:bad} as
  \[
    |B_2 \vp_2| 
    \leq |\bar v_x \f{\vp_2}{\vp_3}| \cdot || \xi_1 \vp_3 ||_{L^{\infty}}
    +  \f{\vp_2}{\rho_{10}} ( |u_A \rho_{10} | \cdot |\bar \th_{x x}|
    + |v_A \rho_{10}| \cdot |\bar \th_{ x y }| )
    + \f{\vp_2}{\psi_1} ( | u_{x, A} \psi_1| \cdot | \bar \th_x| 
    + | v_{x, A} \psi_1 | \cdot  |\bar \th_y| )
    \teq I_1 + I_2 + I_3.
  \]
Since $u_{x, A}, v_{x, A}$ correspond to $f_{11}, f_{2 0}$ in \eqref{eq:u_linf_est1}, 
using \eqref{eq:u_linf_est1} and \eqref{eq:EE_oper:b}, we bound  terms in $I_3$:
\[
\bal
   | u_{x, A} \psi_1| \cdot |\bar \th_x| 
    + | v_{x, A} \psi_1| \cdot |\bar \th_y|
  & \leq ( \CCs_{11, 1} || \om_1 \vp_1 ||_{L^{\infty}}
  + \CCs_{11, 2} [\om_1 \psi_1]_{C_x^{1/2}}
  + \CCs_{11, 3} [\om_1 \psi_1]_{C_y^{1/2}}  ) |\bar \th_x|  \\
  & \qquad + ( \CCs_{20, 1} || \om_1 \vp_1 ||_{L^{\infty}}
  + \CCs_{20, 2} [\om_1 \psi_1]_{C_x^{1/2}}
  + \CCs_{20, 3} [\om_1 \psi_1]_{C_y^{1/2}}  ) |\bar \th_y|  \\
& = \CCs(\uu_x, \na \bar \th,  1)  || \om_1 \vp_1||_{L^{\infty}}
+ \CCs( \uu_x, \na \bar \th,  2 ) [\om_1 \psi_1]_{C_x^{1/2}} 
+ \CCs(\uu_x, \na \bar \th, 3) [\om_1 \psi_1]_{C_y^{1/2}} ,
  \eal
\]
with computable coefficients $\CCs_{ij,k}$ and $\CCs(\uu_x, \na \bar \th, i)$ depending on the profile and weights. 

We bound  transport-type terms in $I_2$ using  
estimate \eqref{eq:u_linf_est1}, $\rho_{10} = \rho_{01}$ by \eqref{wg:u}, and \eqref{eq:EE_oper:a}
\[
|u_A \rho_{10} | \cdot |\bar \th_{x x}|    + |v_A \rho_{10}| \cdot |\bar \th_{ x y }|
\leq  \CCs_{01,1} \nlinf{ \om_1 \vp_1 }  \cdot |\bar \th_{x x}|  + \CCs_{10, 1} \nlinf{ \om_1 \vp_1} \cdot |\bar \th_{x y}| 
= \cT_u( \bar \th_x) \nlinf{\om_1 \vp_1}. 
\]

Combining estimates of $I_j$ yields the upper bound of $B_2 \vp_2$ in \eqref{eq:linf_decay2}. 
 Other terms  in \eqref{eq:linf_decay2} are estimated similarly. Thus, the estimates do not depend on the 
 perturbation; they depend only on the precomputed coefficient functions $\CCs_{ij,k}$ in \eqref{eq:EE_oper} and the profile.

\subsubsection{Coupling weights between the $L^{\inf}$ and the H\"older norm}\label{sec:wg_couple}

The weighted $L^{\inf}$ estimates are not closed, since estimates \eqref{eq:linf_decay2} involves $[\om_1 \psi_1]_{C_{x_i}^{1/2}}$, which controls $(\na \uu)_A$ via \eqref{eq:u_linf_est1}.
To apply the linear stability Lemma \ref{lem:PDE_stab}, we need to combine weighted $L^{\inf}$ and $C^{1/2}$ (semi)norms into one coupled energy with coupling weight $\tau_1$ so that the diagonal dominated condition \eqref{eq:PDE_diag} holds.

Recall the weighted H\"older seminorm from \eqref{hol:g}. We introduce the first energy for $W_1$
\beq\label{energy1}
E_1(W_1) =\max(  \max_i  || W_{1,i} \vp_i ||_{\inf}, \tau_1^{-1}  \max( [ \om_1 \psi_1 ]_{C^{1/2}_{g_1}},  \sqrt{2} || \om_1 \vp_4 ||_{\inf} ) 
) , \  \vp_4 = \psi_1 |x_1|^{-1/2} .
\eeq
When no confusion arises, we write $E_1$ for $E_1(W_1)$ and $E_1(t)$ for $E_1(W_1(t))$.  The coupling weight $\tau_1$ is chosen below. Note that the H\"older seminorm $[\ \cdot \ ]_{ C_{g_1}^{1/2} }$  \eqref{hol:g} is only defined in $\R_2^{++}$. We add the extra $L^{\inf}$ norm $ || \om_1 \vp_4||_{\inf}$ to control 
the seminorm in \emph{upper half space} $[\om_1 \psi_1]_{C_x^{1/2}(\R_2^+)}$. See more discussions in Section \ref{sec:wg_linf_decay}. Since 
\[
g_1(h_1, 0) =| h_1|^{- \f{1}{2}}, \quad g_1(0, h_2) = g_1(0, 1) |h_2|^{-\f{1}{2}}
\]
by \eqref{eq:wg_hol_g},  using \eqref{eq:hol_w_cross} and the estimate in Section \ref{sec:wg_linf_decay}, we obtain 
\begin{gather}
E_1(t) \geq \tau_1^{-1 } \max( [ \om_1 \psi_1]_{C_{g_1}^{1/2}}, \sqrt{2} || \om_1 \psi_1 |x_1|^{- \f{1}{2}}||_{L^{\inf}}  )  \geq  \tau_1^{-1} \max(  [ \om_1 \psi_1]_{ C_x^{1/2}(\R_2^+) }, \  g_1(0, 1) [ \om_1 \psi_1 ]_{C_y^{1/2}}), \notag \\ 
 [ \om_1 \psi_1]_{ C_x^{1/2}(\R_2^+) } \leq \tau_1 E_1(t),
 \qquad
 [ \om_1 \psi_1 ]_{C_y^{1/2}} \leq \tau_1 g_1(0,1)^{-1} E_1(t). 
  \label{eq:hol_to_E1}
\end{gather}
This is the basic conversion from the H\"older seminorms in
\eqref{eq:u_linf_est1} to the coupled energy $E_1$.

Using $E_1(t)$, the above estimates, and the notation \eqref{eq:EE_oper}, we can simplify the estimate \eqref{eq:u_linf_est1} for $\na \uu_A$ by combining $3$ (semi)norms :
$\nlinf{\om_1 \vp_1},  [ \om_1 \psi_1]_{ C_{x}^{1/2}(\R_2^+) },  [ \om_1 \psi_1]_{ C_y^{1/2}(\R_2^+) }$:
\beq\label{eq:u_linf_est12}
\bal
  |\rho_{ij} f_{ij}| 
  & \leq  (  \CCs_{ij , 1}  + \CCs_{ij, 2} \tau_1 + \CCs_{ij, 3} \tau_1 g_1(0,1)^{-1} ) E_1  = E_1 \cdot \CCs( f_{ij}, (1, \tau_1, \tau_1 g_{1}(0,1)^{-1} )).
 \eal 
\eeq

We substitute \eqref{eq:hol_to_E1} in the $C^{1/2}$ seminorm terms in the bound for $B_i$ in \eqref{eq:linf_decay2} and simplify 
\[
  \CCs( f, \na \bar \th, 2) [\om_1 \psi_1]_{C_x^{1/2}} 
+ \CCs( f, \na \bar \th,  3)  [\om_1 \psi_1]_{C_y^{1/2}} 
\leq  \tau_1 \cdot \CCs(f, \na \bar \th,  (0, 1,  g_1(0, 1)^{-1})) E_1(t), \quad f = \uu_x, \uu_y\;,
\]
where  $\CCs(f, \na \bar\th, a )$ is defined in \eqref{eq:EE_oper}. 
By definition of $E_1$ in \eqref{energy1}, the norms $\nlinf{\om_1 \vp_1}, \nlinf{ \eta_1 \vp_2}$, 
$\nlinf{\xi_1 \vp_3 }$  are all bounded by $E_1$ with \emph{unit constant}.

Applying Lemma \ref{lem:PDE_stab} to the system \eqref{eq:linf_decay1}, \eqref{eq:lin_main}  with energy $E_1$
and using estimates \eqref{eq:linf_decay2}, we obtain the linear stability constraints \eqref{eq:PDE_diag} for $\om_1 \vp_1, \eta_1 \vp_2, \xi_1 \vp_3$:
\bseq\label{eq:linf_diag_decay}
{\small
\begin{align}
  -d_{1, L}(\vp_1) - \f{\vp_1}{\vp_2} - \f{\vp_1}{\rho_{10}} \cT_u(\bar \om) \hspace{1.9in} & \geq \lam, \label{eq:linf_diag_decay:w} \\
\big( - d_{2, L}(\vp_2) - \f{\vp_2}{\vp_3} |\bar v_x| - \f{\vp_2}{ \rho_{10}} \cT_u(\bar \th_x)
 - \f{\vp_2}{ \psi_1} \CCs(\uu_x, \na \bar \th,  1)  \big)
 - \tau_1  \f{\vp_2}{\psi_1} \CCs( \uu_x, \na \bar \th,  (0, 1, g_1(0, 1)^{-1}))  & \geq \lam, 
 \label{eq:linf_diag_decay:eta} \\
 \big( - d_{3, L}(\vp_3)  - \f{\vp_3}{\vp_2} |\bar u_y| - \f{\vp_3}{ \rho_{10}} \cT_u(\bar \th_y)
 - \f{\vp_3}{ \psi_1} \CCs(\uu_y, \na \bar \th,  1)  \big)
 - \tau_1  \f{\vp_3}{\psi_1} \CCs( \uu_y, \na \bar \th , (0, 1, g_1(0, 1)^{-1}))  & \geq \lam 
  \label{eq:linf_diag_decay:xi} .
\end{align}
}
\eseq
In \eqref{eq:linf_diag_decay:w}, \eqref{eq:linf_diag_decay:eta}, \eqref{eq:linf_diag_decay:xi},  the first parenthesis contains the damping term and all bad terms that can be bounded using \emph{only} the weighted $L^{\inf}$ part of $E_1$ \eqref{energy1}, which 
is \emph{independent of} the $\tau_1$. 
The second term is
the loss coming from the $C^{1/2}$ seminorms in the estimate of $(\na \uu)_A$ in \eqref{eq:linf_decay2} and \emph{depends on} $\tau_1$ and the coupled  $C^{1/2}$-seminorm.
Thus \eqref{eq:linf_diag_decay} is precisely where the weighted $L^{\inf}(\vp_i)$ stability estimate
of $\om_1, \eta_1, \xi_1$ places an upper constraint on $\tau_1$.

Inequalities \eqref{eq:linf_diag_decay} impose the first constraints on 
the coupling weight $\tau_1$. The \(L^\infty\) estimates prefer \(\tau_1\) not
too large due to the last term in \eqref{eq:linf_diag_decay:eta},\eqref{eq:linf_diag_decay:xi}, because the H\"older contribution in the nonlocal estimates is
multiplied by \(\tau_1\). In contrast, later, the H\"older estimates prefer \(\tau_1\) large,
because the $L^{\infty}$ component enters the H\"older estimates with a relative weight \(\tau_1^{-1}\).

\vs{0.05in}

\paragraph{\bf Making the $L^{\infty}$-$C^{1/2}$ coupled terms perturbative}

We discuss the choice of $\tau_1$ below. The $L^{\infty}$-$C^{1/2}$ estimates are coupled through  the estimate of 
nonlocal remainder $(\na \uu)_A$ defined in \eqref{eq:vel_rem}.  Since the integral defining $\na \uu(\om_1)(x)$ is more regular away from the singularity $y=x$, by increasing the rank of the approximation $\wh{\na \uu}$, we can systematically reduce the coefficient functions  $\CCs_{ij, k}(x), k\leq 3$, in the estimate \eqref{eq:u_linf_est1}  of the velocity remainder $(\na \uu)_A$  for $x$ in any fixed compact domain.\footnote{
Recall the estimates in the model problem in the localized velocity in Section \ref{sec:model_trunc}. The non-singular part $(\na \uu)(x) - (\na \uu)(x, \e )$
can be approximate by a finite-rank operator with error bounded by $\nlinf{\om_1 \vp_1}$ in compact domain. The $L^{\infty}$ pointwise estimate of singular part $(\na \uu)(x, \e, \e )$ 
by $C^{1/2}$ norm contains the \emph{small} factor $\e^{1/2}$. 
}
This reduction is crucial since smaller $\CCs_{ij,k}$ yields smaller
coefficients $\CCs(\pa \uu, \na \bar \th, \cdot)$ in \eqref{eq:linf_diag_decay}, thereby allowing
larger choices of $\tau_1$ while still satisfying \eqref{eq:linf_diag_decay}.
For any prescribed large value of $\tau_1$,  the functions $\CCs_{ij, k}$ can be reduced by the above mechanism, at the expense of using a sufficiently high-rank approximation,  so that the coupled terms are absorbed by the remaining diagonal damping. This yields the diagonal-dominance inequality
\eqref{eq:linf_diag_decay}.

 We want to obtain an overall stability factor $\lam_*$ in the range \eqref{eq:EE_target1} and thus choose $\lam \approx  \lam_*$. We choose the largest $\tau_1$ for which \eqref{eq:linf_diag_decay} holds, with the available finite-rank approximation.  We choose $\tau_1$ large for the same reason as $\tau$ in \eqref{eq:model3_5} and \eqref{eq:model_trunc_EE} for the model problem. 
 This choice makes it easier to establish the weighted $C^{1/2}$ stability estimates.
\footnote{
Consider Lemma \ref{lem:PDE_stab} with energy $E_1$ \eqref{energy1}.
In the stability condition \eqref{eq:PDE_diag} for the norm $\tau_1^{-1} [\om_1 \psi_1]_{C_{g_1}^{1/2}}$ in $E_1$, we have $\mu_i = \tau_1^{-1}$. Choosing a larger $\tau_1$ yields smaller 
$\mu_i$ and $\sum_{ j\neq i } a_{ij} \mu_i \mu_j^{-1}$ in this condition.
\label{fn:hol_tau}
}
In our estimate, we choose 
\beq\label{eq:EE_tau1}
\tau_1 = 5.
\eeq

Although $\tau_1$ is only moderately large, it is enough for us to control the bad terms 
perturbatively in the H\"older estimate when $|x-z|$ is not very small. 
There are three reasons. First, we get the small factor $\tau_1^{-1}$ when estimating the more regular terms using the weighted $L^{\inf}$ norm. See Footnote
\ref{fn:hol_tau}. Second, as $|x-z|$ increases, the localized estimates in the $C^{1/2}$ estimates 
in Lemmas \ref{lem:holx_ux}-\ref{lem:holx_uy} give smaller constants. Third, for $|x-z|$ not too small, the H\"older estimate of nonlocal terms $(\na \uu)_A$ using \eqref{eq:u_linf_est12} and triangle inequality can provide estimates better than Lemmas \ref{lem:holx_ux}-\ref{lem:holx_uy}. Note that from \eqref{eq:linf_diag_decay}, choosing a larger $\tau_1$ requires better estimates of the nonlocal terms, e.g., smaller coefficients $\CCs_{jk, i}$ in \eqref{eq:u_linf_est1}, \eqref{eq:EE_oper}.
This would require higher-rank approximations of $\uu$ and $\na \uu$, and hence higher computational cost. 
Thus we use the moderate value \eqref{eq:EE_tau1}.

\begin{figure}[t]
   \centering
      \includegraphics[width = \textwidth  ]{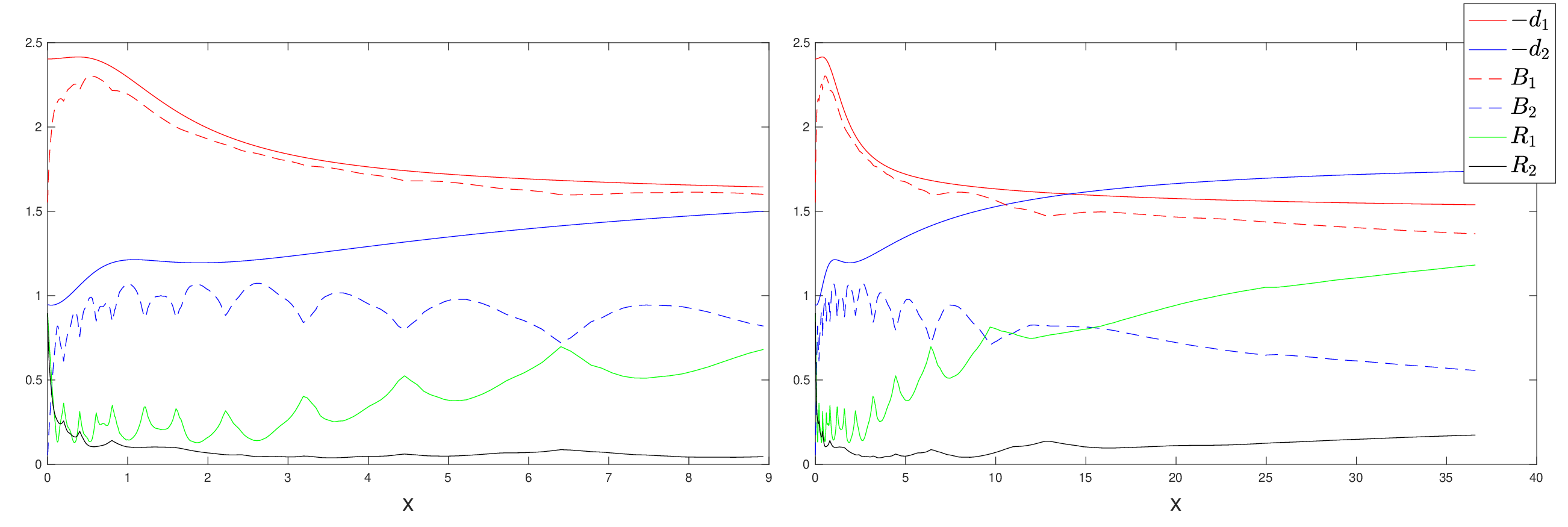}
      \captionsetup{width=0.96\textwidth} 
            \caption{Weighted $L^{\inf}(\vp_i)$ stability estimates.
Left: estimates in $x \in [0, 1.8]$; Right: estimates in $ x \in [0, 35]$. The red curves show the coefficient of the damping term $-d_1(\vp_1)$ and the estimate of $B_1$; the blue curves show $ -d_2(\vp_2), B_2$ in the estimates of $\eta_1$. The green and black curves are the stability factors in the estimates of $\om_1$ and $\eta_1$. A dashed curve lying below the solid one indicates $B_i < -d_i$ and that condition \eqref{eq:PDE_diag} for linear stability holds.
      }
                  \label{fig:Linf1}
 \end{figure}

Below, we visualize the $L^{\infty}(\vp_i)$  estimates.  Figure \ref{fig:Linf1} plots the rigorous piecewise lower bounds for the damping: $-d_{2, L}(\vp_2) > 0$ \eqref{eq:linf_diag_decay}, estimates of the bad terms, i.e. 
the terms subtracted from the damping $-d_{i, L}$ in \eqref{eq:linf_diag_decay}, and the lower bound of the stability factors (LHS of \eqref{eq:linf_diag_decay}) along the boundary. 
The $\xi_1$ variable enjoys much better estimates, so we do not plot its estimates.

We choose approximations $\hat \uu, \wh {\na \uu }$ for $\uu, \na \uu$ along the boundary in Section \ref{sec:appr_vel} so that the weighted estimates of $\uu - \hat \uu, \na \uu - \wh {\na \uu }$ are small. Near the approximation center $x_i$ in \eqref{eq:u_appr_1st}, we have better estimates of the bad terms. In Figure \ref{fig:Linf1}, the points $x_i$ are near the local minimum of the blue dashed curve. Since the coefficients of  nonlocal terms \eqref{eq:lin_main}, e.g., $\bar \om_x, \bar \th_x$, decay at large $x$, 
in addition to the approximations near boundary \eqref{eq:u_appr_1st}, we only add $7$ approximation terms \eqref{eq:u_appr_2nd}, \eqref{eq:vel_cutoff_bd} in a large domain $[0, 200]^2$.  No further approximations are needed for large $x$.

We choose $\vp_1$ slightly weaker than $\vp_2$ near the origin \eqref{eq:wg_linf_decay} so that $\vp_1 / \vp_2 || \eta_1 \vp_2||_{\inf}$ is small, and we can obtain large stability factors 
$( \geq 0.7 )$ for both $\om$ and $\eta$ near $0$. This allows us to control a larger weighted residual error near the origin. 

\begin{figure}[t]
   \centering
      \includegraphics[width = \textwidth  ]{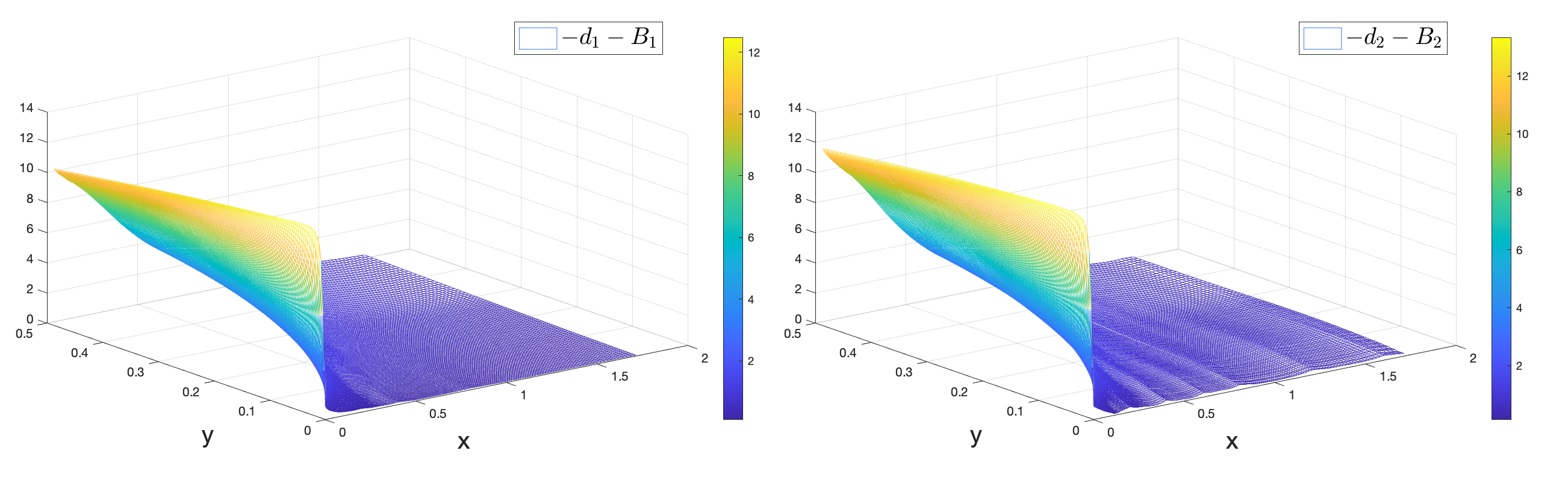}
      \caption{Weighted $L^{\inf}$ stability estimates in the near-field. Left: stability factors in the estimate of $\om_1$. Right: stability factor in the estimate of $\eta_1$.
      }
                  \label{fig:Linf12}
 \end{figure}

In Figure \ref{fig:Linf12}, we plot rigorous piecewise lower bounds of the stability factors, e.g., the left hand side of \eqref{eq:linf_diag_decay}, in $L^{\inf}(\vp_i)$ estimates of $\om_1, \eta_1$ in the near-field. Due to the anisotropy of the flow, the damping terms and the stability factors are larger if the angle $y/x$ is large. See Section \ref{sec:aniso_est}.

To verify inequality \eqref{eq:linf_diag_decay},  we derive its piecewise bounds  using the estimates of the approximate steady state and the weights in Appendix C and Appendix A.1 of Part II \cite{ChenHou2023b}.

\subsubsection{Choosing the weights $\vp_i$}\label{sec:wg_linf_decay}

We discuss the order of choosing parameters and weights for $\vp_i$ in Figure \ref{Fig:order}. 
The formulas are given in \eqref{wg:linf_decay}. We first choose $\vp_1$ with different powers to capture the vanishing order of $\om_1$ near $0$ and its decay at infinity. 
We include the power $|x_1|^{-1/2}$ in $\vp_i$ \eqref{eq:wg_linf_decay} to control $|| \om_1 |x_1|^{-1/2} \psi_1||_{\inf}$ for the $C^{1/2}$ estimate. 
See Section \ref{sec:hol_Q1}. 
The energy estimates of $|| \om_1 |x_1|^{-1/2} \psi_1||_{\inf}$ require bounding  
$ W_{1, i} |x_1|^{-1/2} \psi_1$ and other singular weighted quantities. We then adjust the parameters in $\vp_1$ to obtain a good damping factor $d_{1, L}(\vp_1)(x)$ defined in \eqref{eq:dp} 
for $\om_1$. 
With $\vp_1$ and $\psi_1$ chosen in \eqref{eq:wg_hol},
we estimate nonlocal terms $\uu_A, (\na \uu)_A$ in \eqref{eq:u_linf_est1}. 

Next, using estimates for $(\na \uu)_A, \uu_A$, which are independent of $\vp_2, \vp_3$, and formally setting $\xi_1 = 0$, we choose exponents of different powers in $\vp_2$ and adjust coefficients in $\vp_2$.
The goals are to (a) obtain better stability factors in the weighted $L^{\inf}$ estimates of $\om_1, \eta_1$, and 
(b)  allow a larger $\tau_1$ in \eqref{eq:linf_diag_decay}.  Since equations for $\xi_1$ and $\eta_1$ are similar, we choose the same combination of powers in $\vp_2$ and $\vp_3$ \eqref{eq:wg_linf_decay}.  As $\xi_1$ is weakly coupled with $(\om_1,\eta_1)$ and has much stronger stability estimate (see Section \ref{sec:idea_main_term}), we determine the parameters in $\vp_3$ \emph{only} after obtaining $\vp_1, \vp_2$.

\subsubsection{Weighted $L^{\inf}$ estimate related to the H\"older norm}\label{sec:hol_Q1}

Recall that $W_1 = (\om_1, \eta_1, \xi_1 )$ from \eqref{eq:EE_W1}, $\om_1(x), \eta_1(x)$ are odd in $x_1$ 
and $\xi_1$ is even in $x_1$. To simplify our energy estimates, using the symmetry of $W_{1,i}(x)$ in $x_1$, we only perform H\"older estimates in $\R^2_{++}$, which controls
\[
 I_x  \teq \f{ ( W_{1,i} \psi_i)(x_1, x_2) - ( W_{1,i} \psi_i) (z_1, x_2)}{ |x_1 - z_1|^{1/2}} , 
 \quad I_y \teq  \f{ ( W_{1,i} \psi_i)(x_1, x_2) - ( W_{1,i} \psi_i) (x_1, z_2)}{ |x_2 - z_2|^{1/2}} , 
\]
for $  (x_1,x_2),  (z_1, x_2), (x_1 ,z_2)  \in \R^2_{++}$. 
Due to the symmetry in the $x$ direction, we have 
\[
|I_y(x_1, x_2, z_2)| = | I_y(-x_1, x_2, z_2)|  .
\]

Hence the $C^{1/2}$ estimate in $\R^2_{++}$  controls $C_y^{1/2}$ seminorm in $\R^2_+$. The $x_1$ direction is slightly different. 
Since $\na \uu$ is nonlocal, to control the weighted H\"older norm of $\na \uu(\om_1)$ in $\R^2_{++}$, we need to control $[\om_1 \psi_1]_{C_{x_i}^{1/2}(\R^2_+)}$.  Yet, the above estimate does not directly control $I_x(x_1, x_2, z_1)$ with $x_1 <0< z_1 $ or $[\om_1 \psi_1]_{C_{x}^{1/2}(\R^2_+)}$. 
In fact, for an odd function $f$, it is easy to obtain that 
\[
[f]_{C_x^{1/2}(\R^2_{+})} \leq \sqrt{2} [ f]_{C_x^{1/2}(R^2_{++})}, 
\]
which leads to an extra factor $\sqrt{2}$ into the H\"older estimate. 
Instead, we use a weighted $L^{\infty}$ estimate to control $I_x(\om_1 \psi_1)$.
Let $F \teq \om_1 \psi_1$. For $x_1 < 0< z_1$, since $F$ is odd in $x_1$,
\[
|F(x_1, x_2) - F(z_1, x_2)| = | F(-x_1, x_2) + F(z_1, x_2)|, \quad  |x_1 - z_1| = |x_1| + |z_1|
\]
 we have 
\beq\label{eq:hol_w_cross}
\bal
\f{ | F( -x_1, x_2) + F(z_1, x_2) | }{ |x_1 - z_1|^{1/2}} 
& \leq  \max \B( \f{ 2 | F(-x_1, x_2) | }{  |2 x_1|^{1/2}}  ,  
\f{ 2 | F (z_1, x_2) | }{ |2 z_1|^{1/2}}  
\B) \f{  |x_1|^{1/2}  +  |z_1|^{1/2} }{  ( 2|x_1| + 2|z_1|)^{1/2}  }  \\
& \leq \max \B( \f{ 2 | F(-x_1, x_2) | }{  |2 x_1|^{1/2}}  , 
\f{ 2 | F(z_1, x_2) | }{ |2 z_1|^{1/2}}  
\B),
\eal
\eeq
where we  used the Cauchy-Schwarz inequality in the last inequality. Thus it suffices to control $|| \om_1 \psi_1 |x_1|^{-1/2} ||_{L^{\inf}}$. From \eqref{eq:hol_w_cross} and  $ [\om_1 \psi_1]_{C_{g_1}^{1/2}}  \geq  [\om_1 \psi_1]_{C_x^{1/2}(\R^2_{++}) }$ \eqref{energy1}, we include the norm 
\[
\tau_1^{-1} \sqrt{2} || \om_1 \vp_4||_{\inf},  \quad \vp_4 \teq \psi_1 |x_1|^{-1/2}
\]
with a specific weight $\tau_1^{-1} \sqrt{2}$ in the energy $E_1$ \eqref{energy1} so that $E_{1}(t) \geq \tau_1^{-1} [\om_1 \psi_1]_{C_x^{1/2}(\R_2^+)}$. See \eqref{eq:hol_to_E1}. Thus, by controlling the H\"older seminorm on $\R^2_{++}$ and the $L^\infty(\vp_4)$ norm, energy $E_1$ effectively  controls the H\"older seminorm on $\R^2_{+}$ 
needed in the nonlocal estimates.

   We then perform $L^{\inf}(\vp_4)$ estimate of $\om_1$ using the estimates of nonlocal terms and derivations in \eqref{eq:linf_decay2} and Section \ref{sec:linf_decay} 
\[
\bal
  & \pa_t (  \om_1 \vp_4) + (\bar c_l x + \bar \uu ) \cdot \na ( \om_1 \vp_4)
 = d_{4,L}(x) (  \om_1 \vp_4) +  B_1(x) \vp_4,  \\
& d_{4,L}(x)  \teq d_{1,L}(\vp_4) =   \cT_d( \vp_4) + \bar c_{\om} , \quad 
|B_1(x) \vp_4|  \leq \f{\vp_4}{\vp_2} || \eta_1 \vp_2||_{\inf}
+  \f{\vp_4}{\rho_{10}} \cT_u( \bar \om) || \om_1 \vp_1 ||_{\inf} ,
 \eal
\]
where  $\cT_d, \,  d_{1,L}( \cdot )$ are defined in \eqref{eq:dp}, and $\cT_{\uu}(\bar \om)$ is defined in \eqref{eq:EE_oper}. 
Applying Lemma \ref{lem:PDE_stab} with energy $E_1$ in \eqref{energy1}, 
we obtain condition \eqref{eq:PDE_diag} for $|| \om_1 \vp_4||_{\inf}$ with weight $ \tau_1^{-1} \sqrt{2}$
\footnote{
 Consider Lemma \ref{lem:PDE_stab} 
 with the energy $\max( c_1 \max || W_{1,i} \vp_i ||_{L^{\infty}}, c_2 || \om_1 \vp_4||_{L^{\infty}}  )
 $, where $ c_1 = 1$ and $c_2 = \tau_1^{-1} \sqrt{2}$ are related to $E_1$ \eqref{energy1}. 
 Since $c_2 / c_1 = c_2 $, the factor $\mu_i$ in the inequality \eqref{eq:PDE_diag} for
 $|| \om_1 \vp_4||_{L^{\infty}}$ is $c_2$. Thus, we have the factor $c_2 =  \tau_1^{-1} \sqrt{2} $ in  \eqref{eq:linf_diag_add}.
}
\beq\label{eq:linf_diag_add}
- d_{4,L}  -  \f{\sqrt{2} }{\tau_1 }  ( \f{\vp_4}{\vp_2}  +  \f{\vp_4}{\rho_{10}} \cT_u( \bar \om) ) \geq \lam  ,
\eeq
for some $\lam > 0$. From \eqref{wg:hol}, \eqref{wg:linf_decay}, the ratio is bounded: $\vp_4 / \vp_2 = |x_1|^{-1/2} \psi_1 / \vp_2 \les 1$.

Compared with the weighted $L^{\infty}(\vp_1)$ estimate \eqref{eq:linf_diag_decay:w}, the damping margin here is larger because the bad terms are multiplied by the smaller coupling parameter 
$\f{\sqrt 2}{\tau_1}$ with $\tau_1 = 5$ \eqref{eq:EE_tau1}.

\subsection{Weighted H\"older estimates}\label{sec:EE_hol}

Recall the weights $\psi_i$ for $\om_1, \eta_1, \xi_1$ \eqref{wg:hol} and $g_i(h)$ \eqref{eq:wg_hol_g} in the weighted $C^{1/2}$ seminorms, the notation $W_{1,i}$ \eqref{eq:EE_W1}, and the simplified equation \eqref{eq:lin_sim}. The goal of the weighted H\"older estimate is to control $ [ W_{1,i} \psi_i ]_{C_{g_i}^{1/2} }$,
where $[ \ \cdot  \ ]_{C_{g_i}^{1/2}}$ is defined in \eqref{hol:g}. Together with the weighted $L^{\inf}$ estimate, we control the second energy 
\beq\label{energy2}
\bal
E_{2,0}( W_1)  & \teq    \tau_1^{-1} 
\max( [ W_{1,1}\psi_1 ]_{C_{g_1}^{1/2}} , \mu_1
[ W_{1, 2}\psi_2 ]_{C_{g_2}^{1/2}} ,   \mu_2 [ W_{1, 3}\psi_3 ]_{C_{g_3}^{1/2}}
),  \\
E_2(W_1)  &\teq  \max( E_1(W_1), E_{2, 0}( W_1) ),  \qquad \mu_{h} \teq \tau_1^{-1}(1,  \mu_1, \mu_2 ) .
 \eal
 \eeq
The energy $E_{1}$ is defined in \eqref{energy1}. When no confusion arises, we write $E_2$ for $E_2(W_1)$ and $E_2(t)$ for $E_2(W_1(t))$.  The coupling weights $\mu_1, \mu_2$ 
  \footnote{
  One can absorb $\mu_1,\mu_2$  into the definitions of  $\psi_2, \psi_3$. 
 Yet, in view of Lemma \ref{lem:PDE_stab}, it is easier to 
 adjust $\psi_i$ (form of the weight) and $\mu_i$ (size of weighted norm) separately for better stability estimates using the forms in \eqref{energy2}.
 }
are determined by analyzing the most singular scenario in Section \ref{sec:hol_singu} \eqref{wg:hol_mu1}. They are given in \eqref{wg:EE}.
 The energy $E_{1}$ is defined in \eqref{energy1}. 
 We have normalized the coefficient of the most singular power $|x|^{-2}$ in $\psi_1$ to be $1$. See \eqref{wg:hol}.

 Following the derivations in the weighted H\"older estimates in Section \ref{sec:model1} and using \eqref{eq:lin_sim} and Lemma \ref{lem:hol_comp}, we derive the following for $ W_{1,i} \psi_i $ and any $x, z \in \R^2_{++}$
\beq\label{eq:hol_est1}
\bal
\pa_t H_i + ( b(x) \cdot \na_x  + b(z)  \cdot \na_z ) H_i & = \big( ( d_{i, L}(\psi_i ) W_{1, i} \psi_i) (x)
- ( d_{i,L}(\psi_i) W_{1, i} \psi_i) (z) \big) g_i(x-z) \\ 
& \quad + d_{g_i} H_i + \big( (B_i \psi_i)(x) - (B_i \psi_i)(z) \big) g_i(x-z) \\
& \teq   I_1 + I_2 + I_3 \teq R_i,
 \eal
\eeq
Here, $b(x)$ is the coefficient of the advection \eqref{eq:dp}, $d_{g_i}$ is the damping factor from $g_i$ in the H\"older estimate, and $J_i, H_i$ are given below 
\beq\label{eq:hol_nota}
\bal
  J_i \teq W_{1, i} \psi_i  , \quad 
 H_i(x, z) =  ( J_i(x) - J_i(z) ) g_i(x-z) , \quad
  d_{g_i} \teq   \f{ (b(x) - b(z)) \cdot (\na g_i) (x-z) }{g_i(x-z)}  .
\eal
\eeq
The factor $B_i$ is the bad term defined in \eqref{eq:bad}, and $d_{i,L}$ is defined in \eqref{eq:dp} 
\beq\label{eq:hol_dp}
\bal
d_{1,L}(\psi_1) &= \cT_d( \psi_1) + \bar c_{\om} , 
\qquad  d_{2,L}(\psi_2) = \cT_d(\psi_2) + 2 \bar c_{\om } - \bar u_x, \qquad 
d_{3,L}(\psi_3) = \cT_d( \psi_3) + 2 \bar c_{\om} + \bar u_x .
\eal
\eeq

We note that the second term $I_2$ in \eqref{eq:hol_est1} is already a damping term. See Section \ref{sec:model1} and discussion below Lemma \ref{lem:hol_comp}. To further simplify the notation, we introduce 
\beq\label{eq:hol_dp2}
a_i(x) = d_{i,L}(\psi_i)(x) .
\eeq

The H\"older estimate has three parts.  First, Section~\ref{sec:hol_basic} introduces
bookkeeping notation that packages elementary H\"older estimate. Second, Sections~\ref{sec:EEhol_idea}--\ref{sec:hol_est_I3} estimate  \(I_1,I_2,I_3\) in \eqref{eq:hol_est1}. Third, Section~\ref{sec:EE_hol_sum}
assembles these estimates and discusses the verification methods.

Below, we use $C_{x_1}^{1/2}$ interchangeably with
$C_x^{1/2}$, and $C_{x_2}^{1/2}$ interchangeably with $C_y^{1/2}$.

\subsubsection{Basic H\"older estimates and bookkeeping notations}\label{sec:hol_basic}

Compared with the weighted $L^\infty$ estimates, the H\"older estimates involve
two points $x,z\in\R^2_{++}$. The main difficulty comes from the $C^{1/2}$ estimates of  $\na \uu$. Since the sharp $C^{1/2}$ estimates for $\na \uu$ in Section~\ref{sec:sharp} are directional, in the $x_1$ and $x_2$ directions, 
we connect $x$ and $z$ by rectangular paths and use the triangle
inequality to reduce the H\"older quotient
$|f(x)-f(z)|g(x-z)$ to directional $C_{x_1}^{1/2}$ and $C_{x_2}^{1/2}$ estimates.

To track the coefficients of estimates effectively, we introduce the
notations $\d,\d_i,\d_{\sq}$ below to package several elementary 
estimates: product decompositions, H\"older  estimates for products, and estimates along the two rectangular
paths connecting $x$ and $z$. This reduces many weighted H\"older bounds to
directional $C_{x_i}^{1/2}$ estimates, while isolating the more delicate nonlocal estimates.

We first introduce the directional  H\"older quotients and difference 
\begin{align}\label{eq:hol_dif}
& \d_1(f, x, z) \teq \f{ |f( x ) - f(z)|}{|x-z|^{1/2}} , 
 \  z_1 > x_1, \,  z_2 = x_2,
 \quad 
 \d_2(f, x, z) \teq \f{ |f( x ) - f(z)|}{|x-z|^{1/2}} , 
 \  z_2 > x_2, \,  z_1 = x_1,  \notag \\
  & \d(f)(x, z) \teq f(x)- f(z) .
\end{align}
Here, $\d_1(f), \d_2(f)$ relate to the $C_x^{1/2}, C_y^{1/2}$ estimates, respectively.

Given $f, g$ and $x, z$, we have two simple identities for the product :
\bseq\label{eq:hol_prod}
\beq
\bal
  (f g)(x) - (fg)(z)
  &= (f(x) - f(z)) g(x) + (g(x) - g(z)) f(z)
   = \d(f) g(x) + \d(g) f(z),  \\
    (f g)(x) - (fg)(z) 
  & = (f(x) - f(z)) g(z) + (g(x) - g(z)) f(x)
  = \d(f) g(z) + \d(g) f(x) .
  \eal
\eeq

From these two identities, using the triangle inequality, we get
\beq
\bal
\d(f g)(x, z) & = \d(f)(x, z) g(a) + \d(g)(x, z) f(b) , \quad (a,b) = (x, z), \rm{ \ or \ } (z, x),   \\
\d_i(f g , x, z) & \leq \d_i(f, g, x, z), \quad 
\mbox{for \ } i = 1,  \ x_2 = z_2, \quad \mbox{  or }  \quad  i =2,  \ x_1 = z_1 ,
\eal
\eeq
\eseq
where we slightly abuse the notations of $\d_i$ and define $\d_i(f, g, x, z)$ as a basic $C^{1/2}$ estimate 
\beq\label{eq:hol_prod0}
\d_i(f, g, x, z) = \min_{(a, b) = (x,z), (z,x)} \d_i(f, x, z) |g(a)|  + \d_i(g, x, z) |f(b)  |, \quad  x_{3-i} = z_{3-i} .
\eeq
If $g = 1$, we get $\d_i(f, g, x, z) = \d_i(f, x, z)$. 

\begin{figure}[t]
\centering
\begin{subfigure}{.49\textwidth}
  \centering
  \includegraphics[width=0.45\linewidth]{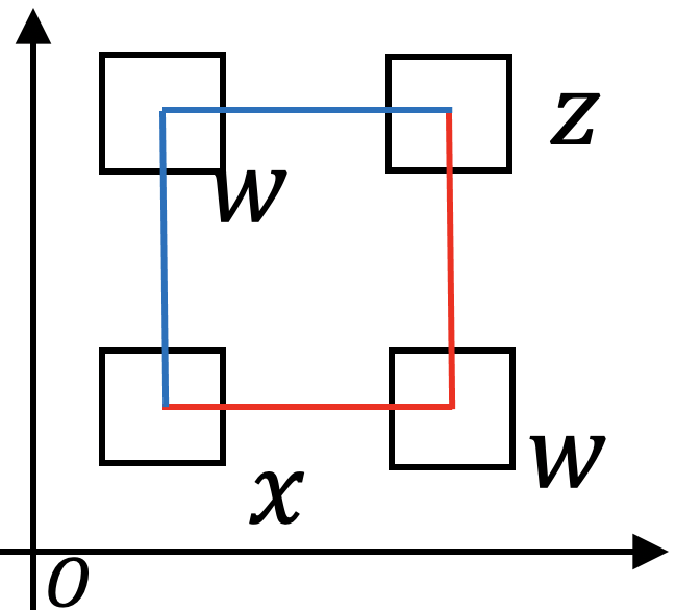}
    \end{subfigure}
    \begin{subfigure}{.49\textwidth}
  \centering
  \includegraphics[width=0.45\linewidth]{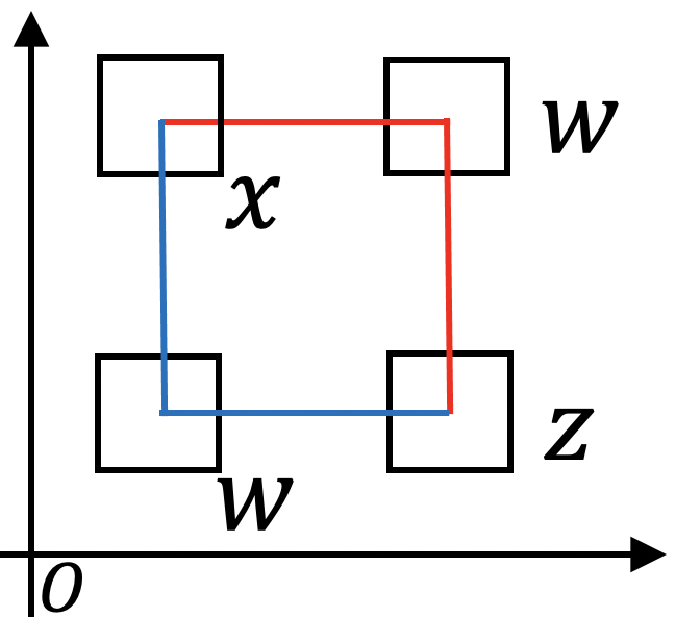}
    \end{subfigure}
\caption{Left, right figures correspond to the locations of $(x, z)$ in cases (1) $s \geq 0$, (2) $s< 0$, $s =(z_1-x_1)(z_2 - x_2)$. The red line and blue line represent two choices of the $C_{x_1}^{k}, C_{x_2}^{k}, k =\f{1}{2}$ or $1$ estimates used to estimate $f(x)-f(z)$. 
 }
\label{fig:hol_pa}
\end{figure}

We next record the path estimate used to pass from directional estimates to a two-point H\"older
estimate. Let $ h = z-x$ and set $w= (z_1, x_2)$. Then
\beq\label{eq:hol_path}
\bal
  | (f(z) - f(x) ) \rho(h) |  & =  | ( f(z) - f(w) + f(w) - f(x) ) \rho(h)| \\ 
 & \leq  ( \d_1( f, x, w ) |h_1|^{1/2} + \d_2(f, w, z ) |h_2|^{1/2} ) \rho(h) , 
  \quad  w =     (z_1, x_2) \  ,
 \eal
\eeq
where $\rho(h) = g(h)$ or $\rho \equiv 1$. The function $|h_i|^{1/2} g(h)$ is $0$-homogeneous and we apply the method in Section \ref{sec:hol_wg_est} to estimate it. 
Similarly,  we derive another estimate by choosing $w=  (x_1, z_2)$
\[
 | (f(z) - f(x) ) \rho(h) | 
  \leq  ( \d_2( f, x,  (x_1, z_2) ) |h_2|^{1/2} + \d_1(f,  (x_1, z_2), z ) |h_1|^{1/2} ) \rho(h) .
\]
We optimize two estimates 
with $\rho = g$ for $|\d(f)(x, z)| g(x-z)$. In Figure \ref{fig:hol_pa}, we illustrate the locations of $x, z$ and two estimates of $\d(f,x, z) g(h)$ based on 
the $C_{x_i}^{1/2}$ estimates ($\d_i(f, p, q)$). 

We introduce the notation \(\delta_\square\) to avoid repeatedly writing
the two possible rectangular-path estimates and to package the optimized rectangular-path estimate for product 
\beq\label{eq:hol_sq}
\bal
\d_{\sq}(f, g, x, z, s) 
& \teq \min( \d_1(f, g, x, (z_1, x_2) )  |s_1|^{1/2} + \d_2(f, g, (z_1, x_2), z ) |s_2|^{1/2},  \\
& \qquad  \quad \ \d_2(f, g, x, (x_1, z_2) )  |s_2|^{1/2} + \d_1(f, g, (x_1, z_2), z ) |s_1|^{1/2} ) , \\
\d_{\sq}(f, x, z, s) & \teq \d_{\sq}(f, 1, x, z, s), \\
\eal
\eeq
where $\d_i(\cdot)$ is defined in \eqref{eq:hol_prod0}, and we simplify the $\d_{\sq}$-notation if $g=1$. We use the ``box" notation $\sq$ since it mimics Figure \ref{fig:hol_pa} and indicates that we use $C_x^{1/2}, C_y^{1/2}$ path estimates to obtain the $C^{1/2}$ estimate. By definition, we get 
\beq\label{eq:hol_sq2}
|\d(f)(x,z ) | \leq \d_{\sq}(f, x, z, h), \quad | \d(f g)(x, z) | \leq \d_{\sq}(f, g, x, z, h).
\eeq

Note that 
\[
\d_i(f,x,z), \quad  \d_i(f,g, x,z), \quad \d_{\sq}(f, x,z, s), \quad \d_{\sq}(f, g, x, z, s)
\]
are symmetric in $x, z$. We introduce an extra variable $s$ to reduce bounding $\d_{\sq}(f,x,z,x-z) g(x-z)$ to estimating $\d_i(f, x, z)$ and $|x_j-z_j|^{1/2} g_i(x-z)$ separately. Function $\d_{\sq}(f,x,z, s) g_i(s)$ is $0$-homogeneous in $s$. See Section \ref{sec:EE_hol_sum}. We drop dependence of $\d, \d_{\sq}, \d_i$ on $x, z$ when  no confusion arises. These basic estimates reduce estimates for most terms to directional $C_{x_i}^{1/2}$ estimates.

\subsubsection{Estimate explicit coefficients}\label{sec:EEhol_idea}

In the weighted $C^{1/2}$ estimates, we estimate several products of explicit coefficients and perturbation quantities. We first explain how the explicit coefficients are estimated.
These weighted coefficients, denoted by $\bar p$ below, 
  depend only on the weights $\psi_i, \vp_i$, and the approximate steady state. Their derivatives can be estimated effectively using the method in  Appendix A.1 and Appendix C of Part II \cite{ChenHou2023b}. Away from $x_1 = 0$, we have $\bar p \in C^1$. 
\footnote{
    Although the weights $\psi_i, \vp_i$ are singular near $\xx=(x_1,x_2)=0$, 
    since the coefficient of the top-order nonlocal terms, i.e. $\na \bar \th$ in 
    $\uu_{x, A} \cdot \na \bar \th, \, \uu_{y, A} \cdot \na \bar \th$ \eqref{eq:lin_main}, 
    vanishes $O(x)$ near $x=0$, the estimates of nonlocal terms are much smaller near $x=0$. 
    See also Remark \ref{rem:dp_near0}. The more difficult part of the $C^{1/2}$ estimates comes from the region away from $0$, where 
  $\psi_i, \vp_i$ can be treated as smooth functions. 
  The singularity of the weights along $x_1 = 0$ and near $x=0$ mainly creates 
  \emph{technical difficulties} in deriving $C^{1/2}$ estimates for the coefficients and weights.
  }

We estimate the piecewise $C_x^{1/2}$ and $C_y^{1/2}$ seminorms of $\bar p(x)$ using the methods in  supplementary material I \cite[Section {\suppsecbdother}]{ChenHou2023aSupp} and Part II \cite[Appendix E]{ChenHou2023b}, and then combine these directional estimates with \eqref{eq:hol_sq}, \eqref{eq:hol_sq2} to estimate $g(x-z) \d( \bar p)$. For example, given $x, z$,  we have 
\[
\bal
 |\bar p(x) -\bar p(z)| \leq \d_i( \bar p, x, z) |x_i - z_i|^{1/2}
 = A_i(  x, z) |x_i - z_i|^{1/2}, \quad x_{3-i} = z_{3-i} \; ,
\eal
\]
for some constants $A_i$ depending on the weights and the approximate steady state. We discretize a very large  domain $[0, D]^2$ in $\R^2_{++}$ using the same mesh 
$y_i$ in Section \ref{sec:ASS_rep} and estimate these constants for $x \in Q_1, z \in Q_2$ for different grids uniformly. 

For instance, for the $C_x^{1/2}$ estimate, 
we track the piecewise bounds 
\[
A_1(x_1, x_2, z_1, x_2), \   z_2 = x_2,  \quad  \mbox{for} \ (x_1, x_2, z_1) \in  I_i \times I_j \times I_{i+ k}, \ I_{\ell} \teq [y_{\ell}, y_{\ell+1}], \ 0\leq k \leq m-1. 
\]
This gives an $n \times n \times m$ array of piecewise coefficient bounds. We have another estimate in \eqref{eq:hol_sq}, \eqref{eq:hol_sq2} by choosing another path from $x$ to $z$ and we optimize two estimates. We can restrict $z_1$ within $m$ grids from $x_1$ since for $z_1$ far-apart from $x_1$, $|z_1 - x_1|$ is not small, we can apply the triangle inequality to obtain the piecewise H\"older estimate.

In general, such an estimate has some overestimates, but it is well adapted to
the anisotropic regime of the problem.  In the worst cases for
the H\"older estimate, one often has $|x_2-z_2|\ll |x_1-z_1|$, and the directional estimate described above is effective.  See also Section \ref{sec:wg_hol_g}.

We now discuss the product estimate for a coefficient $\bar p$ and a perturbation $q$.
Using \eqref{eq:hol_dif}, we obtain 
\beq\label{eq:diff_dec1}
P  \teq  \d( \bar p  q) g( x- z )
 = (   \bar p(x) \d( q ) + \d( \bar p ) q(z) ) g(x- z) 
 \teq P_1 + P_2. 
 \eeq

The first term $P_1$ is controlled by the weighted H\"older seminorm of $q$.
The second term $P_2$ is more regular: the difference falls on the more regular, explicit
coefficient $\bar p$, while $q$ is controlled by the weighted $L^\inf$ norm. 
We also use the other product decomposition in
\eqref{eq:hol_prod} and optimize the two estimates using
\eqref{eq:hol_sq}, \eqref{eq:hol_sq2}. Below, we discuss different cases. In all cases, the estimate of $P_2$ is much smaller than that of $P_1$ when $|x-z|$ is small.

\subsubsection{Estimate of $I_1$}\label{sec:EE_hol_bad}

Recall $I_1, J_i, H_i , a_i$ from \eqref{eq:hol_est1}-\eqref{eq:hol_dp2} 
\[
H_i = \d(J_i) g_i(x-z), \quad J_i = W_{1, i} \psi_i,  \quad a_i(x) = d_{i, L}(\psi_i).
\]
We want to control the H\"older quantity  $H_i$.  We have 
\beq\label{eq:hol_I1_1}
\bal
I_1 & =  \d( a_i J_i) g_i( x- z )
 = ( a_i(x) \d( J_i) + \d( a_i) J_i(z) ) g_i(x- z) \\
 &= a_i(x) H_i +  \d(a_i) g_i(x- z) J_i(z) \teq I_{11} + I_{12}.
\eal
\eeq
The first term is a damping term. 
The second term is lower order, since the difference falls on the explicit coefficient
$a_i$.  We control $J_i(x)$ using the weighted $L^{\inf}$ norm in the energy $E_1$  \eqref{energy1}
\[
 |J_i(z)|= | (W_{1, i} \psi_i)(z) | \leq  || W_{1, i} \vp_i ||_{\inf}  \f{ \psi_i(z)}{\vp_i(z)}
 \leq E_1 \f{ \psi_i(z)}{\vp_i(z)}.
\]
Since $a_i$ is a given explicit function, we follow Section \ref{sec:EEhol_idea} and estimate $ \d(a_i) g_i(x- z)$ using the methods 
in supplementary material I \cite[Section {\suppsecbdother}]{ChenHou2023aSupp}, Part II \cite[ Appendix E]{ChenHou2023b}, and \eqref{eq:hol_sq},\eqref{eq:hol_sq2}. 
Away from $x=0$, when $|x-z|$ is small, since $a_i$ is locally $C^1$, $ \d(a_i) g_i(x-z)$ is very small. Hence
\beq\label{eq:hol_I12}
|I_{12}| \leq |\d(a_i)g_i(x-z) |  \f{ \psi_i(z)}{\vp_i(z)} E_1 .
\eeq

Similarly, we can also put the damping coefficient at the other point $z$ and estimate 
\beq\label{eq:hol_I1_2}
\bal
I_1  & = \td I_{11} + \td I_{12},\quad \td I_{11} = a_i(z) H_i, 
\quad \td I_{12} = \d(a_i) g_i(x- z) J_i(x) \\
 | \td I_{12}| 
 & \leq |\d(a_i)g_i(x-z) |  \f{ \psi_i(x)}{\vp_i(x)} E_1 .
\eal
\eeq

We choose one of the above estimates according to the relative size of the following terms 
\beq\label{eq:hol_I1_3}
 m_1 =  a_i(x) + \mu_{h,i} |\d(a_i)g_i(x-z) |  \f{ \psi_i(z)}{\vp_i(z)}, 
\quad  m_2 = a_i(z) + \mu_{h, i} |\d(a_i)g_i(x-z) |  \f{ \psi_i(x)}{\vp_i(x)}, 
\eeq
where $\mu_{h} \in \R^3_+$ is the weight of the H\"older seminorm in \eqref{energy2}. We use the decomposition \eqref{eq:hol_I1_1} and its estimates \eqref{eq:hol_I12} if $m_1 $ is smaller; we use \eqref{eq:hol_I1_2} if $m_2 $ is smaller. 
Since the left hand side of the stability condition \eqref{eq:PDE_diag}  
contains the negative of these quantities $-m_i$, this optimization 
provides a larger stability factor. 
\footnote{
Consider Lemma \ref{lem:PDE_stab} with $E_2$ \eqref{energy2}. 
In \eqref{eq:PDE_diag} for $[ W_{1, i} \psi_i ]_{C_{g_i}^{1/2}}$, 
$-a_i(p)>0$ in \eqref{eq:hol_I1_3} gives  the diagonal coefficient $a_{ii}$,  while $\cI_1(q) \teq \mu_{h,i} |\d(a_i)g_i(x-z) |  \f{ \psi_i(q)}{\vp_i(q)}$ in \eqref{eq:hol_I1_3}
contributes to the off-diagonal terms $a_{ij} \mu_i \mu_j^{-1}, i\neq j$. 
They contribute $\cI(p, q) \teq -a_i(p) - \cI_1(q)$ for $(p,q) = (x, z)$ or $(z, x)$ to LHS of \eqref{eq:PDE_diag}. 
Thus, choosing the smaller of $m_1,m_2$ in \eqref{eq:hol_I1_3} maximizes the
corresponding $\cI(x, z), \cI(z, x)$ and gives a larger stability factor in \eqref{eq:PDE_diag}. 
}
We  use similar optimizations several times to get larger stability factors, see, e.g. \eqref{eq:hol_prod0},\eqref{eq:hol_sq}. Following the methods in Section \ref{sec:EEhol_idea}, we track the piecewise bounds of the above functions and estimates, e.g. $a_i, \, \mu_{h, i} |\d(a_i)g_i(x-z) |  \f{ \psi_i(z)}{\vp_i(z)}$.

\begin{remark}[\bf Singular weights and regularity of $a_i$]\label{rem:dp_near0}
Since $ \psi_i$ \eqref{wg:hol} is singular, 
$a_i = d_{i, L}(\psi_i)$  in \eqref{eq:hol_dp}, \eqref{eq:dp} is not $C^{1/2}$ near $x=0$. 
However, the ratio $\psi_i/ \vp_i\lesssim |x|^{1/2}$ compensates this
loss in the weighted coefficient-difference term:
\[
 \d(a_i, x, z )|x-z|^{-1/2} \psi_i(x) / \vp_i(x)
\les  \d(a_i, x, z )|x-z|^{-1/2}  |x|^{1/2},
 \]
 which  remains bounded near $x=0$.
In Section {\secbddphol} of the supplementary material I \cite{ChenHou2023aSupp} (attached to this paper), we perform an improved estimates near $0$ and bound the explicit functions $\min( m_1, m_2)$ \eqref{eq:hol_I1_3} from above. Since the 
bad terms \eqref{eq:bad} are very small near $x=0$ due to the vanishing coefficients, e.g. $\bar \th_x$, this technical difficulty only has a tiny effect on the stability estimate. See Figure \ref{fig:hol1}. We optimize the improved near-origin estimate with the previous one.
\end{remark}

\subsubsection{Estimate of $I_3$}\label{sec:hol_est_I3}

The term \(I_3\) in \eqref{eq:hol_est1} contains all forcing terms \(B_i\) defined in \eqref{eq:bad}. 
We divide its estimate into local terms involving perturbation $(\om_1, \eta_1, \xi_1)$
in \eqref{eq:lin_main}, and nonlocal terms involving 
$\uu_A, (\na \uu)_A$ in \eqref{eq:lin_main}.   The local terms are handled by product estimates and 
the energy $E_2$ \eqref{energy2}; the nonlocal terms use weighted $L^{\infty}$ and $C^{1/2}$ nonlocal estimates, e.g. \eqref{eq:holest}.

\vs{0.1in}
\paragraph{\bf{Estimate of the local part}}
We focus on $\eta_1$ in $B_1$ \eqref{eq:bad} in the $C^{1/2}$ estimates of $\om_1 \psi_1$. Other terms $ - \bar v_x \xi_1$ in $B_2$, $-\bar u_y \eta_1$ in $B_3$ \eqref{eq:bad} can be estimated similarly. Note that the weights $\psi_1, \psi_2$ are different for $\om_1, \eta_1$. Using \eqref{eq:hol_prod}, we rewrite the difference as follows 
\[
 \d( \eta_1 \psi_1) g_1(x- z)
 = \d ( \eta_1 \psi_2 \cdot  \f{\psi_1}{\psi_2} ) g_1(x-z)
 = \B( \d( \eta_1 \psi_2 )  \f{ \psi_1 }{ \psi_2 }(x)
+  ( \eta_1 \psi_2 )(z) \d( \f{\psi_1}{ \psi_2} ) \B) g_1(x -z) 
\teq P_1 + P_2.
\]

The term $P_2$ is more regular since the difference falls on the more regular ratio $\psi_1/\psi_2$. We follow Section \ref{sec:EEhol_idea} and use \eqref{eq:hol_prod}-\eqref{eq:hol_sq2} to estimate $\d( \f{\psi_1}{ \psi_2} )  g_1(x -z)$. Using the weighted $L^{\inf}$ norm of $\eta_1$ and the energy $E_2$ \eqref{energy2}, we obtain 
\beq\label{eq:hol_eta_reg}
|P_2|
\leq || \eta_1 \vp_2 ||_{\inf}  \f{\psi_2(z) }{ \vp_2(z) } \B| \d_{\sq}( \f{\psi_1}{ \psi_2}, h )  g_1(x -z) \B|
\leq E_2  \f{\psi_2(z) }{ \vp_2(z) } \B| \d_{\sq}( \f{\psi_1}{ \psi_2}, h )  g_1( h) \B|, \quad h = x-z.
\eeq
Following Section \ref{sec:EEhol_idea}, we track the piecewise bound of the coefficient in this upper bound. 

For $P_1$, we use the weighted $C^{1/2}$ seminorm for $\eta_1 \psi_2$:
\beq\label{eq:hol_eta_grat}
\bal
|P_1|  & \leq 
\f{ \psi_1 }{ \psi_2 }(x)  \f{ g_1(x-z)}{g_2(x-z)} | \d( \eta_1 \psi_2) g_2( x-z)|
\leq \f{ \psi_1 }{ \psi_2 }(x) \B| \f{ g_1(x-z)}{g_2(x-z)} \B| [ \eta_1 \psi_2 ]_{C_{g_2}^{1/2}} \\
&\leq 
 \f{ \psi_1 }{ \psi_2 }(x) \B| \f{ g_1(x-z)}{g_2(x-z)} \B|   E_2  \tau_1 \mu_1^{-1} ,
 \eal
\eeq
where we used the energy $E_2$ \eqref{energy2} and $\mu_{h, 2}^{-1} = \tau_1 \mu_1^{-1}$ in the last inequality. Using another decomposition in \eqref{eq:hol_prod}, we get another estimate and we optimize them. 

 Although the bound \eqref{eq:hol_eta_grat} contains the factor $\tau_1$, 
 this does not amplify the final $C^{1/2}$ estimate,
even though $\tau_1$ is chosen to be moderately large and, in principle, can be chosen sufficiently
large; see Section \ref{sec:wg_couple}. 
Indeed, the seminorm $[\om_1\psi_1]_{C_{g_1}^{1/2}}$ enters $E_2$ \eqref{energy2} with weight $\tau_1^{-1}$.
Thus, in the estimate for $\tau_1^{-1} [\om_1\psi_1]_{C_{g_1}^{1/2}}$ of the energy, the relevant term is $\tau_1^{-1}P_1$. 
This prefactor $\tau_1^{-1} $ cancels  $\tau_1$ in the upper bound \eqref{eq:hol_eta_grat}.
The coupling weights $\mu_1, \mu_2$ in \eqref{energy2} measure the relative size 
among three $C^{1/2}$ seminorms and are chosen independently of $\tau_1$.

Since $g_1$ and $g_2$ are equivalent to $|h|^{-1/2}$ and homogeneous of order $-1/2$,  the quantity $\B| \f{ g_1(x-z)}{g_2(x-z)} \B|$ only depends on the ratio between $x_1 - z_1, x_2 - z_2$. We track this ratio for later verification.

For large $|x-z|$, we have a simple estimate based on triangle inequality 
\beq\label{eq:hol_loc_lg}
|\d( \eta_1 \psi_1) g_1(x- z)|
\leq || \eta_1 \vp_2||_{\inf}  ( \f{\psi_1}{\vp_2}(x) + \f{\psi_1}{\vp_2}(z)  ) g_1(x-z)
\leq E_2  ( \f{\psi_1}{\vp_2}(x) + \f{\psi_1}{\vp_2}(z)  ) g_1(x-z).
\eeq

\vs{0.1in}
\paragraph{\bf{Estimate of other local terms}}

Recall $W_1 = ( \om_1, \eta_1, \xi_1)$ from \eqref{eq:EE_W1} and the weights $\mu_{h}, \psi_i$ \eqref{wg:hol} in the energy $E_2$ \eqref{energy2}.
For $x_1 \leq z_1, i, j \in\{1,2,3\}$ and  $f \in C^{1/2}$, using the energy 
\[
||W_{1, i} \vp_i||_{L^{\inf}} \leq E_2, \quad  \mu_{h, i} [  W_{1, i} \psi_i]_{C_{g_i}^{1/2}} \leq E_2 
\]
from \eqref{energy2}, together with \eqref{eq:hol_prod}, \eqref{eq:hol_sq2}, we perform the $C^{1/2}$ estimate  for $f W_{1, i} \psi_i$ :
\beq\label{eq:hol_loc}
\bal
&  |\mu_{h, j} g_j(h) \d( f  W_{1, i}  \psi_i, x, z)|
 \leq  \mu_{h, j} g_j(h) ( | \d( f) W_{1, i} \psi_i(z)|
+ |  f(x) \d( W_{1, i} \psi_i)| ) \\
& \leq \f{\mu_{h, j} g_j(h)}{ \mu_{h, i} g_i(h)} | f(x) | E_2
+ \mu_{h, j} g_j(h) \d_{\sq}( f , x, z, h) \f{\psi_i}{\vp_i}(z) E_2 , \quad h = x-z.
\eal
\eeq
If $i = j$, the first term reduces to $|f(x)| E_2$. We only pick one decomposition in \eqref{eq:hol_prod} in which the coefficient multiplying $\d( W_{1, i} \psi_i)$ is evaluated at $x$, namely $f(x)$ in \eqref{eq:hol_loc},  to simplify the estimates. Note that we assume $x_1 \leq z_1$ 
and 
$\psi_2 = \psi_3, g_2 = g_3$ for the H\"older weights.
See \eqref{wg:hol}.

We apply \eqref{eq:hol_loc} to estimate $\bar v_x \xi_1 \psi_2$ in $B_2$ \eqref{eq:bad} with $f =\bar v_x, W_{1, 3} \psi_3 = \xi_1 \psi_2$,
and to estimate $\bar u_y \eta_1 \psi_3$ in $B_3$ \eqref{eq:bad} with $f = \bar u_y, W_{1, 2} \psi_2 = \eta_1 \psi_3$. 
These terms appear in the $I_3$-term in \eqref{eq:hol_est1} with $i=2,3$. 
We also apply \eqref{eq:hol_loc} in the nonlinear estimates in Section \ref{sec:non}.

\vs{0.1in}
\paragraph{\bf{Estimate of the nonlocal part}}
\footnote{
 From vanishing orders \eqref{eq:velA_reg} and $\psi_1 \sim |\xx|^{-2}$ \eqref{wg:hol} near $\xx=0$, 
 we have $(\na \uu)_A \psi_1 \sim |\xx|^{1/2}$ and it lies in $C^{1/2}$.
}
To control the nonlocal terms in $B_i$, we use the sharp $C_x^{1/2}$ and $C_y^{1/2}$ estimates in Section \ref{sec:sharp} for the most singular part and the estimates in Section 4 of Part II \cite{ChenHou2023b} for the more regular part. See Section \ref{sec:idea_int_hol} for the ideas. 

We focus on the estimate of $ -u_{x, A} \bar \th_x$ in $B_2$ \eqref{eq:bad}, which contributes to the largest part in the estimate. Using \eqref{eq:hol_prod0}, \eqref{eq:hol_sq}, \eqref{eq:hol_sq2}, we get
\beq\label{eq:hol_P40}
|\d( u_{x, A} \bar \th_x \psi_2 )| 
= |\d( u_{x, A} \psi_1 \cdot \f{\psi_2}{\psi_1} \bar \th_x  )| 
\leq \d_{\sq}( u_{x, A} \psi_1, \f{\psi_2}{\psi_1} \bar \th_x, h ) , \quad P_4 = \f{\psi_2}{\psi_1} \bar \th_x, 
\eeq
and it suffices to estimate $\d_i(  u_{x, A} \psi_1), \d_i( \f{\psi_2}{\psi_1} \bar \th_x)$ and the $L^{\inf}$ bounds of $ u_{x, A} \psi_1,P_4 $. 
For $u_{x, A}(p), p = x,z$, we use the $L^{\infty}$ estimate \eqref{eq:u_linf_est1} in Section \ref{sec:linf_decay}. The term $P_4$ is more regular: it has a vanishing order $|x|^{1/2}$ near $x=0$ and is in $C^{1/2}$. 
\footnote{
   Near $0$, since $\psi_1 \sim |x|^{-2},
  \psi_2 \sim |x|^{-5/2}$, $\bar \th_x \sim x_1$ from \eqref{wg:hol}, we have $\psi_2 \bar \th_x / \psi_1 \sim |x|^{1/2}$  and it lies in $C^{1/2}$.
}
We follow Section \ref{sec:EEhol_idea} to estimate it. 
We have 
\beq\label{eq:hol_P4}
| (u_{x, A} \psi_1) (p) \cdot \d_{i}( P_4) | \leq  C_i(x, z, p) E_1 , \qquad p = x, z, \quad i = 1,2, \quad x_{3-i} = z_{3-i} ,
\eeq
for some functions $C_{1}(x, z, p), C_{2}(x,z, p), p = x, z$ depending on the weights and the approximate profile.
We obtain piecewise upper bound of these functions as in Section \ref{sec:piece_bd}.

For $\d_i(u_{x,A}\psi_1)$, we apply the weighted estimates \eqref{eq:holest:du},
\eqref{eq:holest2} following ideas in \ref{sec:idea_int_hol}. Details are given in
Part II \cite[Section 4]{ChenHou2023b}. 
For $ i =1,2$ 
and $x_{3-i} = z_{3-i}$, we obtain 
\beq\label{eq:hol_nloc}
\bal
\d_i( u_{x, A}\psi_1 ) &\leq \tau_1^{-1} \CCs_{1/2, u_x, i}(x, z , \vp_1)  
\cdot
\max( \tau_1  || \om_1 \vp_1||_{\inf},   [\om_1 \psi_1]_{C_{g_1}^{1/2} }  ) 
\leq \tau_1 \cdot  ( \tau_1^{-1} \CCs_{1/2, u_x, i}(x, z , \vp_1) )   E_2, 
\eal
\eeq
with \emph{computable} coefficients bound $\CCs_{1/2,u_x, i}$ 
\emph{independent of} the perturbation $\om_1$.
\footnote{
We normalize the upper bound in \eqref{eq:hol_nloc} by the norm $\max( \tau_1 || \om_1 \vp_1||_{\inf},   [\om_1 \psi_1]_{C_{g_1}^{1/2} }  )$ so that $[\om_1 \psi_1]_{C_{g_1}^{1/2} }$  appears with \emph{unit prefactor}. 
Then, the estimate \eqref{eq:hol_nloc} is consistent with the sharp  $C^{1/2}$ estimates in Section \ref{sec:sharp}.
}
\footnote{
In the stability inequality \eqref{eq:PDE_diag} for the seminorm $A=   [ \eta_1 \psi_2 ]_{C_{g_2}^{1/2}}$ in $E_2$ with weight $\mu_{h,2}=\tau_1^{-1} \mu_1$, 
estimate for $\mu_{h,2} \d_i( u_{x, A}\psi_1 )$ via the bound \eqref{eq:hol_nloc}  contributes $ \tau_1^{-1} \mu_1  \times  \bar C \CCs_{1/2, u_x, i} $ to the term $a_{ij} \mu_i \mu_j^{-1}$. Here, 
 $  \bar C \CCs_{1/2, \cdot, i}  $ is the bound of $B_i$ by the energy $E_2$:  $\sum_{j \neq i} a_{ij}\mu_j^{-1}$ in \eqref{eq:PDE_stab_2} and \eqref{eq:PDE_diag}, while $\tau_1^{-1}\mu_1$ corresponds to $\mu_i$ in \eqref{eq:PDE_diag}.  
Note that $\bar C$ bounds the coefficient of nonlocal terms, e.g. $\bar \th_x$ in
\eqref{eq:hol_P40} and is $\tau_1$-independent; $\tau_1^{-1} \CCs_{1/2, u_x , i}$ decreases but \emph{bounded away from $0$} as $\tau_1 \to \infty$. 
Thus, the contributions of the bounds in \eqref{eq:hol_nloc}, \eqref{eq:hol_P40} to
condition \eqref{eq:PDE_diag} are not large when a larger $\tau_1$ is chosen.
}

We remark that when $|x-z|$ is small, the coefficients $\tau_1^{-1} \CCs_{1/2, u_x, i},i=1,2$  are very close to the constants provided by the sharp H\"older estimates in Section \ref{sec:sharp}, since the contribution of the more regular nonlocal part is small relative to the energy $E_1$.
Again, we track these piecewise upper bounds as in Section \ref{sec:piece_bd}. Combining the above estimates and the piecewise $L^{\inf}$ estimate of $P_4$ in \eqref{eq:hol_P40}, and using \eqref{eq:hol_sq2} and \eqref{eq:hol_prod}, we yield the estimate for $u_{x, A} \bar \th_x \psi_2$.

When $|x-z|$ is not small, we can apply the triangle inequality and the $L^{\inf}$ estimate of $ u_{x, A}$ in \eqref{eq:u_linf_est1} in Section \ref{sec:linf_decay} to obtain another bound. 
\footnote{
In practice, we only need to apply the above H\"older estimate when $|x-z|$ is small, e.g., when $x, z$ are within $40$ mesh grids designed in Section \ref{sec:ASS}. Beyond this range, the $L^{\inf}$ estimate already provides a better bound.
}

\vspace{0.05in}

\paragraph{\bf Remaining linear terms.}
The same decomposition, coefficient estimates, and directional H\"older
estimates apply to the remaining linear terms in 
$B_i$ in \eqref{eq:bad}, \eqref{eq:hol_est1}.  We record
the resulting estimates in \eqref{hol:lin}, with the modification in Section
\ref{sec:comb_vel_err} to track the nonlocal error. 
In Section \ref{sec:EE_hol_sum}, the reader may treat $\bar \uu^N = \bar \uu$ in \eqref{hol:lin}.

\subsubsection{Summarize the estimates}\label{sec:EE_hol_sum}

We now collect the estimates obtained in the previous subsections. 
The estimates above reduce the right hand side of \eqref{eq:hol_est1} to
diagonal damping terms and remainder terms bounded by the energy $E_2$ times
explicit coefficient functions.  We then assemble rigorous piecewise bounds
for these coefficients and verify the resulting diagonal dominance condition
\eqref{eq:PDE_diag} for linear stability.

For the right hand sides in \eqref{eq:hol_est1}, when $x,\; z$ are close, we obtain the following estimates 
\[
R_i = ( d_{g_i}(x, z) + a_i( p_{x, z}) ) H_i  + \td B_i, \quad  \td B_i = \hat I_{12} + I_3, \quad I_3 = \d ( B_i \psi_i ) g_i(x-z) , 
\]
for $i =1,2,3$, where 
\[
(p_{x, z}, \hat I_{12}) = (x, I_{12}) \quad  \mbox{or}  \quad (z, \td I_{12}),
\]
depending on the size of $m_1, m_2$ in \eqref{eq:hol_I1_3}.  Here, $\td B_i$ combines all terms that are not part of the diagonal damping: the term $\hat I_{12}$ \eqref{eq:hol_I1_1} or \eqref{eq:hol_I1_2} and $I_3$. We can estimate it as follows:
\beq\label{eq:hol_sum1}
\bal
& |\td B_i| \leq \CCs_i(x, z, h) E_2,  \quad h = x- z .  \\
\eal
\eeq

The coefficient $\CCs_i(x,z,h)$ is assembled from the bounds derived in
Sections \ref{sec:hol_basic}--\ref{sec:hol_est_I3}.  The
coefficient-difference part $I_{12}, \td I_{12}$ from $I_1$ is controlled by
\eqref{eq:hol_I12}, \eqref{eq:hol_I1_2} with optimization
\eqref{eq:hol_I1_3}.  The local terms in $I_3$ are controlled by the
product estimates \eqref{eq:hol_eta_reg}--\eqref{eq:hol_loc_lg}, their direct generalizations, and the
general local bound \eqref{eq:hol_loc}.   The nonlocal terms are controlled
by \eqref{eq:hol_P40}--\eqref{eq:hol_nloc}, using the sharp $C^{1/2}$ estimates from Section \ref{sec:sharp} for the most singular part and the
estimates from Part II \cite[Section 4]{ChenHou2023b} for the more regular
part. These coefficients depend only on weights $\vp_i, \psi_i$ chosen in \eqref{eq:wg_hol}, \eqref{eq:wg_linf_decay},
coupling weights $\tau_1, \mu_i$ in energies $E_1, E_2$ \eqref{energy1}, 
\eqref{energy2}, approximate profile $(\bar \om , \bar \th, \bar c_l , \bar c_{\om})$, and
finite-rank approximation constructed in Section \ref{sec:appr_vel}.  In particular, $\CCs_i(x,z,h)$ 
is a \emph{computable} coefficient function of $x,z, h$, and \emph{does not} depend on the perturbation.

The coefficient $\CCs_i(x,z,h)$ depend on $h$ through several $0$-homogeneous functions.
 For example, they involves $g_i(h) | h_j|^{1/2}$ via $ g_j(h)\d_{\sq}(x, z, h)$ \eqref{eq:hol_sq2} in \eqref{eq:hol_P40}, and $ \f{g_1(h)}{g_2(h)}$ in \eqref{eq:hol_eta_grat}. 
We can obtain piecewise upper bounds of the coefficients $\CCs_i(x, z, h)$ in the above estimates and \eqref{hol:lin}  following the discussions in Sections \ref{sec:piece_bd}, \ref{sec:EEhol_idea} and \ref{sec:EE_hol_bad}. 
Thus, the linear stability condition \eqref{eq:PDE_diag} 
associated with the $C^{1/2}$ seminorm in energy $E_2$ \eqref{energy2} reduces to
\beq\label{eq:hol_diag}
- d_{g_i}(x, z) - a_i(p_{x, z})  - \mu_{h,i} \CCs_i(x, z, h) > \lam, \quad \mu_{h} = \tau_1^{-1}(1,  \mu_1, \mu_2 ),
\eeq
uniformly in $x, z$ for some $\lam > 0$, where 
$h = x - z$ and $\mu_{h}$ is the coupling weight of the H\"older seminorm in \eqref{energy2}.

\vs{0.1in}
\paragraph{\bf{Verifying the stability conditions}}
It remains to verify \eqref{eq:hol_diag}. 
The fundamental verification method is to partition the range of $x, z$ and $h_1 / h_2$ into smaller intervals $I_x \times I_z \times I_h$, and then obtain a tight bound of \eqref{eq:hol_diag} in  $I_x \times I_z \times I_h$ using monotonicity. We discretize a large domain $[0, D]^2$ into small grid cells using the same mesh $\{ y_i\}_{i=0}^n$ as in Section \ref{sec:ASS}
\[
Q_{ij} = I_i \times I_j , \qquad I_i  =[y_i, y_{i+1}] ,   \qquad y_0 < y_1< .. < y_n, \quad n<800. 
\]

First, we fix the locations of $x, z$: $x \in Q_{ij}, z \in Q_{kl}$, and derive the piecewise bounds $\CCs_i(x, z, h)$ in $x, z$. The bound \eqref{eq:hol_loc_lg} 
using the triangle inequality (and \eqref{eq:check_hol_lg} similarly) involves $g_1(x-z)$ and is not homogeneous in $x-z$. We bound it using monotonicity of $g$ in Section \ref{sec:hol_wg_est}. 

We then estimate functions in $\CCs_i(x, z, h)$  \eqref{eq:hol_sum1}, \eqref{hol:lin} involving $h = x-z$, Since these functions are $0-$homogeneous in $h$, 
their values only depend on the ratio between $h_1$ and $h_2$. Thus, it suffices to consider the ratio between
\[
r_1 = x_1 - z_1,  \qquad r_2 = x_2 - z_2 . 
\]
Similar considerations apply to the damping factors $d_{g_i}$ defined in \eqref{eq:hol_est1}, which we estimate
\[
d_{g_i} 
 \leq C_1(x, z) f_1(x- z) + C_2(x, z)  f_2(x-z) 
 + C_3(x, z) f_3(x- z) + C_4(x, z)  f_4(x-z) 
\]
for some $0-$homogeneous functions $f_i$ 
\footnote{
These functions $f$ have the form $h_i g(h)^{-1} \pa_j g(h)$ for some $i, j\in \{1,2 \}$, where $g$ is some H\"older weight \eqref{eq:wg_hol_g}.
}
and piecewise constant bounds $C_i(x, z)$, see 
 the supplementary material I \cite[Section {\secbddg}]{ChenHou2023aSupp} (attached to this paper). 

We consider four different cases of $h$ depending on the sign of $r_1 / r_2$ and the size between $|r_1|, |r_2|$. We focus on $r_1 / r_2 >0$ and $r_2 \leq r_1 $ to illustrate the ideas. In this case,  we can normalize $r_1  = 1$ and $0 \leq r_2 \leq 1$, and reduce verifying \eqref{eq:hol_diag} to verifying an inequality in the single variable $r_2$. Since these functions have monotonicity properties, e.g. $g_i( 1, r_2 )$ is decreasing in $r_2$, these inequalities can be verified by partitioning  $r_2 \in [0, 1]$ into  smaller intervals 
\[
r_2 \in [b_i, b_{i+1}], \qquad 0 = b_0 < b_1 < ... < b_N = 1 .
\]

For fixed grids $x \in Q_1, z \in Q_2$, we can further bound the ratio $|z_2 -x_2| / |x_1 -z_1| \in [r_l, r_u]$ using the piecewise upper and lower bounds for $x_i, z_j$. Since the ratio $r_2 / r_1$ with $(r_1, r_2) = x-z$ lies within $[r_l, r_u]$, 
we do not need to verify \eqref{eq:hol_diag} for $h$ with $h_2 / h_1 = r_2 / r_1 \notin [r_l, r_u]$, or the case $h_2 / h_1 \notin  [r_l, r_u]$ can be marked as satisfying \eqref{eq:hol_diag} automatically.

When $|x-z|$ is far away, we have a much better estimate due to the improvement from the sharp H\"older estimates in Lemmas \ref{lem:holx_ux}-\ref{lem:holx_uy}. 
In practice, for $x \in Q_{i_1, i_2}, z \in Q_{j_1, j_2}$ with $\max( |i_1 - j_1|, |i_2 - j_2|) \geq 15$, we already have much better stability factors, the left side of \eqref{eq:hol_diag}.

\begin{figure}[t]
   \centering
      \includegraphics[width = \textwidth  ]{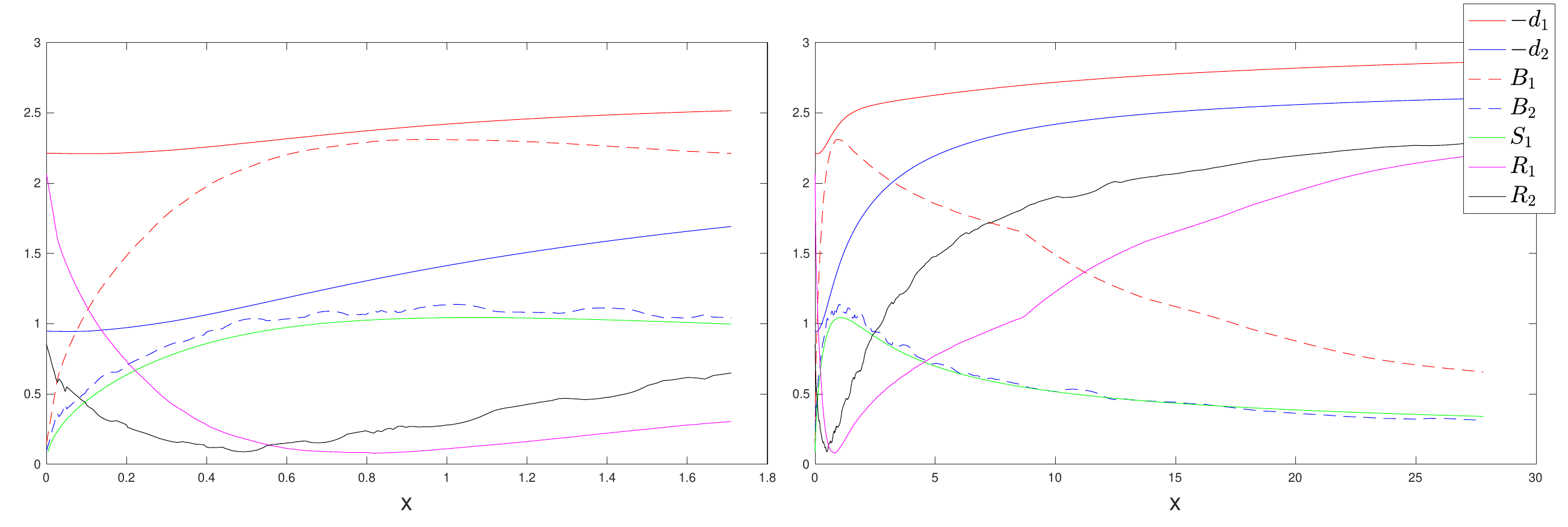}
      \captionsetup{width=0.95\textwidth} 
      \caption{Weighted H\"older estimates. Left figure: estimates near $0$, $x \in [0, 1.8]$; Right figure: estimates in a larger domain, $ x \in [0, 30]$.
      The red curves show the coefficient of the damping term for $ \d ( W_{1,1} \psi_1) g_1(x-z)$ and the estimate of the bad terms; the blue curves are for the H\"older estimate of $ W_{1,2} \psi_2$. The green curve is the same as the bound estimated in the most singular scenario based on the sharp inequalities in Section \ref{sec:hol_singu}. 
      The magenta curve is the stability factor in the H\"older estimate for $\om_1$,
      and the black curve is that for $\eta_1$.
      The stability factors are larger than $0.08$. The (red and blue) dashed curves below the solid one indicate that the weighted diagonal dominance condition \eqref{eq:PDE_diag} for stability holds.    
      }    
            \label{fig:hol1}
 \end{figure}

In Figure \ref{fig:hol1}, we visualize the rigorous piecewise  estimates along the boundary  with $r_2 / r_1 = 0$. Here, we consider  $x \in I_i \times I_0, z \in I_j  \times I_0, j -i = 2$, 
 where  $I_i =[y_i, y_{i+1}]$ is a small interval. This corresponds to one of the cases with the smallest damping. Other cases with small $j-i$ are similar. 
 The estimate of the bad terms in the $\eta_1$ equation is very close to the one in the most singular scenario in Section \ref{sec:hol_singu} based on the sharp inequalities. In some regions,
we have better estimates since $|x-z|$ is far away 
 \footnote{
Since the mesh spacing $y_{k+1}-y_k$ increases with $k$, the
separation $|x-z|$ can be large even when the index gap is small.  For example,
if $x\in I_i\times I_0$ and $z\in I_j\times I_0$ with $j-i=2$, then
$|x-z|$ can be large for large $i$.
 }
 and we have improved constants for the localized velocity from Lemmas \ref{lem:holx_ux}-\ref{lem:holy_ux}. 
For larger ratio $ |r_2 / r_1|$, we have larger stability factors than the case of $|r_2 / r_1|$ being very small due to the anisotropy of the flow. See Section \ref{sec:aniso_est}.

\begin{figure}[t]
   \centering
      \includegraphics[width = \textwidth  ]{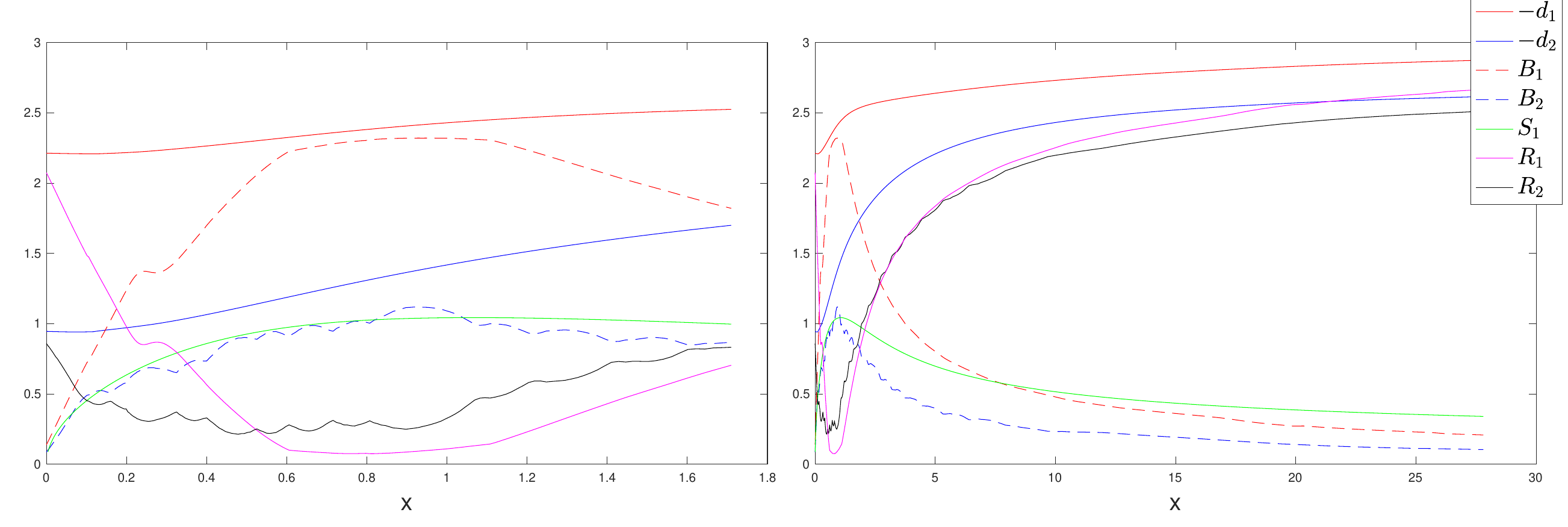}
      \caption{Weighted H\"older estimates with larger $|x-z|$              }    
            \label{fig:hol2}
 \end{figure}

In Figure \ref{fig:hol2}, we consider  $x \in I_i \times I_0, z \in I_j  \times I_0$ with $j - i = 10$. The stability factor for $\eta_1$ shown by the black curve becomes much larger and is larger than 0.3. 

\subsection{Weighted $L^{\inf}$ estimates with growing weights}\label{sec:linf_grow}

To close the nonlinear estimate in \eqref{eq:lin}, \eqref{eq:bous_decoup2}, we need to control $|| \om||_{\inf}, || \na \th||_{\inf}, || W_{1,i}||_{\inf}, ||\na \uu ||_{L^{\inf}}$,
which is not bounded by the energy $E_2$ \eqref{energy2} as $\vp_1$ \eqref{energy1} decays for large $|x|$ (see \eqref{eq:wg_linf_decay}). We further perform weighted $L^{\inf}(\vp_{g,i})$ estimates with the following weights $1\les \vp_{g, i}$ stronger than $\vp_i$ in the far-field
\beq\label{eq:wg_linf_grow}
\bal
\vp_{g,1} &= \vp_1 + p_{71} |x|^{1/16},  &   \vp_{g,2} & = \vp_2 + p_{81} |x|^{1/4} + p_{82} |x|^{ \al_{g,n}},  \\
\vp_{g,3} & = \vp_3 + p_{91} |x|^{1/4} + p_{92} |x|^{\al_{g, n}}, \qquad & \al_{g, n} &= 1/3 + 10^{-8}.
\eal
\eeq
The subscript \textit{``g"} is short for \emph{``grow"}. The $L^{\infty}(\vp_{g,i})$ estimates follow the analysis of the damped system \eqref{eq:damp_sys}.
Thus \emph{no new} stability mechanism is needed here: once the diagonal damping 
term in the estimate of $W_{1,i}\vp_{g,i}$ is negative, the remaining terms are treated \emph{perturbatively}.

The main growing terms $|x|^{1/16}$ in $\vp_{g, 1}$ and $|x|^{1/4}$ in $\vp_{g, 2}, \vp_{g, 3}$ are used to close the nonlinear estimate. The last term in $ \vp_{g, 2}, \vp_{g, 3}$ has a larger growth rate $\al_{g,n}$ but with a much smaller coefficient, $p_{92} \ll p_{91}, p_{82} \ll p_{81}$; they are used to close the nonlinear $C^{1/2}$ estimate. See \eqref{eq:wg_linf_grow_hol}.

We choose $p_{ij}$ so that the damping coefficients in the weighted $L^{\inf}(\vp_{g, i})$ estimate
\begin{gather}\label{eq:linf_grow_dp}
d_{1,L}(\vp_{g,1})  = \cT_d( \vp_{g,1}) + \bar c_{\om} , 
\quad  d_{2,L}(\vp_{g,2})  = \cT_d(\vp_{g,2}) + 2 \bar c_{\om } - \bar u_x, \quad 
d_{3,L}(\vp_{g,3}) = \cT_d( \vp_{g,3}) + 2 \bar c_{\om} + \bar u_x, 
\end{gather}
are negative and bounded by some $d_i$ with $d_i < 0$, where $\cT_d$ is defined in \eqref{eq:dp}.
For $|x|$ large enough, since $\bar \uu$ has sublinear growth implied by \eqref{solu:decay}, 
the main terms of $d_{i,L}$ are given by
\[
d_{1,L}(\vp_{g,1})= \bar c_l \al_1 + \bar c_{\om} + l.o.t., \qquad d_{2,L}(\vp_{g,2}) = \bar c_l \al_2 + 2 \bar c_{\om} + l.o.t., 
\qquad d_{3,L}(\vp_{g,3}) = \bar c_l \al_3 + 2 \bar c_{\om} + l.o.t,
\]
with $\bar c_l \approx 3, \bar c_{\om} \approx - 1$ \eqref{eq:bous_expo}, where $\al_i$ is the exponent of the last power in $\vp_{g,i}$ \eqref{eq:wg_linf_grow}. 
Following the analysis of far-field damping in \eqref{eq:idea_damp_far} and discussions therein,
we choose $\al_i$ so that these main terms are negative. 
We then choose $p_{71} = 1$ and determine the power $p_{81} |x|^{1/4}, p_{91} |x|^{1/4}$ in $\vp_{g,2}, \vp_{g,3}$ by setting $p_{82} = p_{92} =0$ to obtain a damping factor $d_{i,L}(\vp_{g,i}) < 0$ not too close to $0$. We further choose the last power $p_{82} |x|^{\al_{g, n}}, p_{92} |x|^{\al_{g, n}} $ with much smaller parameters $p_{82}, p_{92}$, while preserving $ d_{i,L}(\vp_{g,i}) <0$ not too close to $0$. The parameters are given in \eqref{wg:linf_grow}.

For some coupling weights  $\tau_2, \mu_4$ to be determined, we consider a new energy 
\beq\label{energy3}
\bal
E_3(W_1) & = \max\big( E_2(W_1),  \, \tau_2 \max( \mu_4 || \om_1 \vp_{g, 1} ||_{\inf} , || \eta_1 \vp_{g, 2}||_{\inf}, || \xi_1 \vp_{g, 3}||_{\inf} )  \big), \\
  \mu_g & \teq (\tau_2 \mu_4, \tau_2, \tau_2), \\
  \eal
\eeq
where $E_2(W_1)$ is defined in \eqref{energy2} and $\mu_g \in \R^3_+$ is the coupling weight 
for the norm $\nlinf{ W_{1,i} \vp_{g, i}}$.
 When no confusion arises, we write $E_3$ for $E_3(W_1)$ and $E_3(t)$ for $E_3(W_1(t))$.

To control $||\na \uu_A||_{L^{\inf}}$, we also use $|| \om_1 \vp_{g, 1}||_{ \inf }$ and the H\"older norm of $\om_1 \psi_1$.
Following the methods in Section \ref{sec:u_comp_idea} and using nonlocal estimates similar  to  \eqref{eq:u_linf_est1} in Section \ref{sec:linf_decay} but with the growing weight $\vp_{g,1}$ instead of $\vp_1$, we obtain
\beq\label{eq:u_linf_est2}
 |\rho_{ij} f_{ij}(x)| \leq \CCs_{g,ij, 1}(x) || \om_1 \vp_{g,1}||_{\inf} 
 + \CCs_{g,ij, 2}(x) [ \om_1 \psi_1]_{C_x^{1/2}}
 + \CCs_{g,ij, 3}(x) [ \om_1 \psi_1]_{C_y^{1/2}},
\eeq
where $g$ is short for ``grow", $f_{01} = u_A, f_{10} = v_A$, etc similar to those in \eqref{eq:u_linf_est1}, and $\rho_{ij}$ is the weight for $\uu$ in \eqref{wg:u}, \eqref{wg:elli}.
As in \eqref{eq:u_linf_est11}, we do not need the H\"older norm to control $\uu_A$:
\[
\CCs_{g,ij, 2}(x) =  \CCs_{g,ij, 3}(x) = 0, \qquad  i+j =1
\]
Since $\vp_{g, i}$ is growing, for large $|x|$, the above estimates provide stronger
far-field estimates than  \eqref{eq:u_linf_est12}, and we can obtain $|| f_{ij}||_{L^{\inf}} \leq C(\vp_{g,1}, \psi_1) \max( || \om_1 \vp_{g, 1} ||_{\inf}, [ \om_1 \psi_1]_{C^{1/2}} )$ from \eqref{eq:u_linf_est2}.

We optimize the estimates \eqref{eq:u_linf_est2}, \eqref{eq:u_linf_est12} for $\uu_A, (\na \uu)_A$ and use the energy $E_3 \geq E_1$ to obtain
\beq\label{eq:u_linf_est22}
\bal
|\rho_{ij} f_{ij}(x)| & \leq  \CCs_{g, ij}( \tau_2 \mu_4)(x) E_3(t) \\
 \CCs_{g, ij}( \kp )(x) & \teq \min( 
\CCs_{g, ij, 1}(x) \kp^{-1} + \CCs_{g, ij, 2}(x) \tau_1 
 + \CCs_{g, ij, 3}(x) \tau_1 g_1(0, 1)^{-1}  , \\
 & \qquad  \quad \,
 \CCs(f_{ij}, ( 1,\tau_1, \tau_1 g_1(0, 1)^{-1}  ) \, )   \, ).
  \eal
\eeq
Since $\tau_1$ and $g_1$ have been chosen, $ \CCs_{g, ij}$ depends only on 
the dummy variable $\kp$. 
To control the \emph{transport-type} term $\uu_A \cdot \na f $, we introduce  $\cT_{u, g}$ similar to $\cT_u$ in \eqref{eq:EE_oper} 
\beq\label{def:cT_ug}
\cT_{u, g}(f, \tau) = \CCs_{g, 01}( \tau ) |f_x | + \CCs_{g , 10}(\tau) |f_y |
\eeq

Performing weighted $L^{\inf}(\vp_{g,i})$ estimate 
on \eqref{eq:lin_main}, \eqref{eq:lin_sim} yields
\[
  \pa_t ( \mu_{g,i} W_{1, i} \vp_{g, i} ) + (\bar c_l x + \bar \uu ) \cdot \na ( \mu_{g,i} W_{1, i} \vp_{g,i}) 
 =  d_{i,L}(\vp_{g,i}) \cdot ( \mu_{g,i} W_{1, i} \vp_{g,i}) + \mu_{g,i} \vp_{g,i} B_{i}(x),
 \]
 where  $ d_{i,L}(\vp_{g,i})$  is the linear damping coefficient in  \eqref{eq:linf_grow_dp}.

Using the energy $E_3$ we can obtain pointwise bounds for $\om_1, \eta_1, \xi_1$, e.g. 
\[
\tau_2 \mu_4 |\vp_{g,1} \eta_1 |  
\leq  \tau_2 \mu_4 \vp_{g,1} \cdot (  \vp_2^{-1}  || \eta_1 \vp_2||_{\inf} \wedge \tau_2^{-1} \vp_{g,2}^{-1} ||\tau_2 \eta_1 \vp_{g,2} ||_{\inf} )
\leq \mu_4 \vp_{g,1} ( \tau_2 \vp_2^{-1} \wedge \vp_{g,2}^{-1}  ) E_3(t),
\]
where $a \wedge b = \min(a, b)$. Applying the above pointwise bounds and \eqref{eq:u_linf_est22},
\eqref{def:cT_ug}  for  $\uu_A,\na \uu_A$, we obtain the following
estimate for the bad terms:
\begin{align}\label{eq:linf_grow2}
&  \mu_{g,i} \vp_{g,i} |B_{i}(x)|
 \leq  \mu_{g, i}  A_{g,i}(x) E_3(t) , \hspace{2in} i =1,2,3,  \\
& \tau_2 \mu_4 A_{g,1}  = \mu_4 \vp_{g,1} ( \tau_2 \vp_2^{-1} \wedge  \vp_{g,2}^{-1}  ) 
+  \tau_2 \mu_4 \f{\vp_{g,1}}{ \rho_{01}}  \cT_{u, g}( \bar \om, \tau_2 \mu_4)  ,  \notag \\
& \tau_2 A_{g,2 } =  \vp_{g,2} | \bar v_x | ( \tau_2 \vp_3^{-1} \weg   \vp_{g,3}^{-1} ) 
+  \tau_2 \f{\vp_{g,2}}{\rho_{10}} \cT_{u, g} (\bar \th_x, \tau_2 \mu_4) 
+ \tau_2 \f{\vp_{g,2}}{ \psi_1}
\big( \CCs_{g,11}(\tau_2 \mu_4) |\bar \th_x| + \CCs_{g,20}(\tau_2 \mu_4) |\bar \th_y| \big) , \notag  \\
& \tau_2 A_{g,3}  =   \vp_{g,3} |\bar u_y|( \tau_2 \vp_2^{-1} \weg  \vp_{g,2}^{-1} ) 
+ \tau_2 \f{\vp_{g,3}}{\rho_{10}} \cT_{u, g} ( \bar \th_y, \tau_2 \mu_4) 
+ \tau_2 \f{\vp_{g,3}}{ \psi_1}
\big( \CCs_{g,02}(\tau_2 \mu_4) |\bar \th_x| + \CCs_{g,11}(\tau_2 \mu_4) |\bar \th_y| \big) . \notag
\end{align}

The linear stability condition \eqref{eq:PDE_diag} for $ || W_{1,i} \vp_{g,i}||_{\inf}$ with weights $
\mu_{g} = ( \tau_2 \mu_4,\tau_2, \tau_2 )$ reads 
\beq\label{eq:linf_diag_grow}
- d_{1,L}(\vp_{g,1}) - \mu_4 \tau_2  A_{g,1} \geq \lam, 
\quad  - d_{2,L}(\vp_{g, 2}) - \tau_2 A_{g,2} \geq \lam, \quad  - d_{3,L}(\vp_{g, 3}) - \tau_2 A_{g,3} \geq \lam, 
\eeq
 for some $\lam > 0$, where  $ d_{i,L}(\vp_{g,i})$ is given in  \eqref{eq:linf_grow_dp}.

We have chosen $p_{ij}$ and the weights 
$\vp_{g,i}$. We emphasize that the second bound in \eqref{eq:u_linf_est22}, which comes from 
\eqref{eq:u_linf_est1},  is independent of $\tau_2,\mu_4$; hence in any bounded region, it gives an upper bound for $\cT_{u,g}$ in \eqref{def:cT_ug} which is also independent of $\tau_2,\mu_4$. 
These upper bounds are provided by the energy $E_2$, for which stability has been established.
Thus, we obtain an upper bound of the forcing $\vp_{g,i} B_i$ by $E_2$ 
independent of $\tau_2, \mu_4$, as  seen from \eqref{eq:linf_grow2} after cancelling the common factor $\mu_{g,i}$.
\footnote{
 While $A_{g i}$ in \eqref{eq:linf_grow2} depends on $\tau_2, \mu_4$, 
 the above argument gives an upper bound independent of $\tau_2, \mu_4$. 
}
This is precisely an analog of \eqref{eq:damp_force} in the damped system. 
Thus, it remains to choose $\tau_2,\mu_4$ following the damped system \eqref{eq:damp_sys} so that, 
after  multiplication by the coupling weight $\mu_{g, i}$, the  bounds for $\mu_{g, i} \vp_{g, i}B_i$ in \eqref{eq:linf_grow2} are \emph{perturbative}. 
Since  $\na \bar \om, \na  \bar \th, \na^2 \bar \th,  \na \bar \uu$ decay 
by \eqref{solu:decay}, the asymptotics of the weights 
$\vp_{g,i}, \vp_i$ in \eqref{wg:linf_decay}, \eqref{wg:linf_grow} 
and the nonlocal estimates 
\footnote{
For $ |\om_1 | \les \la x \ra^{- 1/16} \nlinf{ \om_1 \vp_{g, 1} }$ 
 and $[ \om_1 \psi_1]_{C_{x_i}^{1/2}} \les E_2$ by \eqref{eq:wg_linf_grow}, \eqref{eq:hol_to_E1}, 
 standard estimates imply $ | \na \uu(\om_1)(x)|, |x|^{-1} |\uu(\om_1)(x)|
 \les_{\e} \la x \ra^{-1/16 + \e}$ for any $\e>0$. 
}
imply that the estimates $\mu_{g,i} A_{g,i}$ in \eqref{eq:linf_grow2} converge to $0$  as $\tau_2 \to 0$ uniformly for $\mu_4  \leq 1$,  e.g.
\[
\vp_{g,2} | \bar v_x | ( \tau_2 \vp_3^{-1} \weg   \vp_{g,3}^{-1} ) \to 0, \quad \tau_2 \f{\vp_{g,2}}{\rho_{10}} \cT_{u, g} (\bar \th_x, \tau_2 \mu_4)  \to 0.
\]

Thus, we can first choose  $\tau_2$ small enough to achieve the second and third stability condition in \eqref{eq:linf_diag_grow} with $\lam$ similar to that in \eqref{eq:linf_diag_decay}. 
With this fixed $\tau_2$, since  $\tau_2 \mu_4 A_{g,1} \to 0$ as $\mu_4 \to 0$, 
we can choose a small $\mu_4$ to achieve the first condition in \eqref{eq:linf_diag_grow}.

\begin{remark}[\bf Choices of $\tau_2, \mu_4$]

We do not simply set $\mu_4 = 1$ since it requires smaller $\tau_2$ to satisfy \eqref{eq:linf_diag_grow}, which leads to a weaker energy $E_3$ \eqref{energy3} and larger constants in later nonlinear estimates.
Parameters $\tau_2, \mu_4$ are given in \eqref{wg:EE}. The choices of $\vp_{g,i}$ and $\tau_2, \mu_4$ mainly affect the constants in the nonlinear estimates in the far-field since $\vp_i \geq \mu_{g,i}\vp_{g,i}$ for $|\xx|$ not very large.

\end{remark}

\subsection{Estimate of linear functionals}\label{sec:rank1}

In the previous sections, we established weighted $L^{\inf}$ and $C^{1/2}$
stability estimates for the linear equations \eqref{eq:lin_main} of $W_1$,
encoded in the energy $E_3$ \eqref{energy3}, provided that
\eqref{eq:linf_diag_decay}, \eqref{eq:linf_diag_grow},
\eqref{eq:linf_diag_add}, \eqref{eq:hol_diag} hold. 
To close the energy estimates for \eqref{eq:bous_decoup2}, it remains to
estimate the scalar functionals $a_i(W_1)$ and $a_{\mf{nl},i}(W)$ that enter
$\hat W_2$ \eqref{eq:bous_W2_appr} and the residual operators
$\cR$ \eqref{eq:bous_err_op}.
\footnote{
The error in $\cR$ \eqref{eq:bous_err_op} from the approximate solution,
e.g. $ \hat F_i(0) - \bar F_i$, is estimated in Part II
\cite[Section 3]{ChenHou2023b}.
}
Then, we further control $ \hat W_2$ and $\cR$. See Remark~\ref{rem:E4_reg}.

We first use energy $E_3$ \eqref{energy3} to directly bound the functionals for which the weighted
pointwise estimates are sufficient. See the \hyr[para:E3_pointwise]{Paragraph} below.
For the remaining functionals, especially
those involving the slowly decaying kernel in $c_{\om}$, we use sharper ODE
estimates. These ODE estimates follow the damped-system argument in
\eqref{eq:damp_sys}: the diagonal terms provide damping, while the remaining
terms are controlled by the energy with established stability and are treated
\emph{perturbatively} by choosing the new coupling weights in 
an enlarged energy.

We estimate the ODEs for $c_\omega(\omega_1)$, $c_\omega(\omega)$, 
and the Taylor coefficients \(\omega_{xy}(0)\), \(\th_{xxy}(0)\)
in Section \ref{sec:cw_w1}, \ref{sec:cw_w}, \ref{sec:wxx0}, respectively.
These estimates correspond to the ODE estimates in Step 7 of Figure~\ref{Fig:order}.
After these estimates, we construct the final energy $E_4$ in \eqref{energy4}.

\vs{0.05in}

\paragraph{\bf Controlling functionals by energy $E_3$}\label{para:E3_pointwise}

In \eqref{eq:KF_reg2}, we bound $a_i(W), a_{ \mf{nl},i}(W)$ using $C^{2,\al}$ norm of $W_1, W$; here, 
we obtain sharper bounds using the weighted energy $E_3$ in \eqref{energy3}.

From the definition of the energies \eqref{energy1}, \eqref{energy3}, we have the pointwise control:
\beq\label{eq:w_est1}
\bal
|W_{1, i}(x)| & \leq \vp_{Mi}^{-1}(x) E_3(t), & \quad  & W_{1, 1} = \om_1, \ W_{1,2}= \eta_1, \ W_{1, 3} = \xi_1, \\
\vp_{M1}(x) & \teq \max(\vp_1(x ), \tau_2 \mu_4 \vp_{g,1}(x) ), &  \quad 
 & \vp_{Mi}(x) \teq \max(\vp_i(x ), \tau_2  \vp_{g,i}(x) ),  \quad  i =2, 3 .
\eal
\eeq

We have two types of linear functionals $a_i(W_1)$. The first type is $c_{\om}(\om_1)$ from $\cK_{1i}(\om_1)$ \eqref{eq:appr_near0}. The second type is from $\cK_{i2}(\om_1)$ \eqref{eq:appr_vel} for the approximation of $\td \uu,  \na \td \uu$ \eqref{eq:u_appr_2nd}, \eqref{eq:u_appr_1st}, \eqref{eq:u_appr}. For $a_{ \mf{nl},i}$ defined in \eqref{eq:appr_near02_op}, we need to control $c_{\om}(W_1 + \hat W_2)$ and $\na^2( W_1 + \hat W_2)(0)$. 

For the second type of term, it is given by the integral 
$
 \int_{\R^2_{++}} \om_1(y) p(y) dy
$
for some function $p(y)$ that has a fast decay, e.g. it has a decay rate $|y|^{-4}$. 
We have two equivalent formulas \eqref{eq:u_appr_1st}, \eqref{eq:u_appr_1st_sum}
\eqref{eq:u_appr_2nd}, \eqref{eq:u_appr_2nd_sum} for approximating $\td \uu,  \na \td \uu$. Then we estimate it directly using \eqref{eq:w_est1}
\beq\label{eq:finite_mode_2nd}
| \int_{\R^2_{++}} \om_1(y) p(y) dy| 
\leq E_3(t) \int \vp_{M 1}^{-1}(y) |p(y)| dy .
 \eeq

With the estimates of the above functionals and ODEs in Sections \ref{sec:cw_w1}-\ref{sec:wxx0}, since the coefficients in \eqref{eq:u_appr_1st},\eqref{eq:u_appr_1st_sum},\eqref{eq:u_appr_2nd},\eqref{eq:u_appr_2nd_sum}, e.g. $C_f \chi_i, C_f S_i$, are given smooth functions with appropriate vanishing near $x=0$, we  estimate their derivatives and weighted norms following  Part II \cite[Appendix C, D]{ChenHou2023b} and
 supplementary material II of Part II \cite[Section {\secCIIItonorm}]{ChenHou2023bSupp} and estimate  $\hat \uu, \wh {\na \uu} $ \eqref{eq:u_appr}.

Recall the inner product \eqref{inner}. Next, we consider the $c_{\om}$-term  \eqref{eq:normal_pertb}. Controlling 
\beq\label{eq:green}
c_{\om}(\om) = u_x(0) = -\f{4}{\pi} \int_{ \R_2^{++}}  \f{y_1 y_2}{ |y|^4}  \om(y) dy 
= -\f{4}{\pi} \la \om, f_* \ra , \quad f_*(y) = \f{y_1 y_2}{|y|^4}. 
\eeq
effectively is nontrivial since the integrand $ \f{y_1y_2}{|y|^4}$ decays slowly (it is not in $L^1$) and our weight for $\om_1$ is very weak in the far-field. 
See $\tau_2 \mu_4 \vp_{g, 1}$ in \eqref{wg:linf_grow} 
and $\vp_1$ in \eqref{wg:linf_decay} \eqref{wg:EE}.

\begin{remark}[\bf Improved estimates of $c_{\om}$]
Using  the pointwise estimate  \eqref{eq:w_est1} directly gives
$|c_{\om}(\om_1)|\leq C_1E_3$, with $C_1$ about $170$--$300$.  Since this large constant enters the main nonlinear terms; see
\eqref{eq:cw_poor1}, it would force a smaller residual error to close the nonlinear stability estimates. We instead exploit the ODE for $c_{\om}(\om_1)$, whose linear damping
and faster-decaying forcing yield a much sharper bound. Similarly, we derive a sharper estimate of $c_{\om}(\om)$ in Section~\ref{sec:cw_w} via ODEs. 
The estimates below are technical but straightforward to implement with computer assistance.
\end{remark}

\subsubsection{Control $c_{\om}(\om_1)$}\label{sec:cw_w1}

Following \cite{chen2019finite2,chen2021HL}, we perform ODE estimates on $c_{\om}$. 
The estimates follow those of the damped system \eqref{eq:damp_sys} and are perturbative.
Using equations \eqref{eq:lin_main},\eqref{eq:bous_decoup2}, we derive equations for $ \la \om_1, q f_* \ra, \la \eta_1, q f_* \ra$ with $q =1, \chi_{ \mf{ode}}$ (to be defined in \eqref{eq:cw_chi0})
\beq\label{eq:lin_main_cw}
\bal
&\tf{d}{dt} \la \om_1, q f_* \ra 
 = \bar c_{\om} \la \om_1, q f_* \ra + \la \eta_1 , q f_* \ra + \la \G_1, q f_* \ra, 
  \\
& \tf{d}{dt} \la \eta_1, q f_* \ra 
  = 2 \bar c_{\om} \la \eta_1 , q f_* \ra + \la \G_2 , q f_* \ra,  \\
&   \G_1  \teq   \G_{1, M} + \cN_1 + \bar \cF_1 - \cNF_1 -\cR_1 ,
\quad \G_{1, M} = - (\bar c_l x + \bar \uu) \cdot \na  \om_1 - \uu_A \cdot \na \bar \om, \\
 &  \G_2   \teq  \G_{2, M}+ \cN_2 + \bar \cF_2 - \cNF_2 - \cR_2  , 
 \quad \G_{2, M} = -  (\bar c_l x + \bar \uu) \cdot \na \eta_1 + B_2(W_1),
\eal
\eeq
where $\la \cdot , \cdot \ra$  \eqref{inner} is the inner product on $\R^2_{++}$, $ B_2(W_1)$ denotes the bad term \eqref{eq:bad}, and $\cN_i, \bar \cF_i, \cNF_i, \cR_i$ are the nonlinear terms \eqref{eq:bous_non}, residual error \eqref{eq:bous_err}, rank-one correction \eqref{eq:appr_near0}, and residual operator \eqref{eq:bous_err_op}. The transport term $ \uu \cdot \na W_{1, i}$ in \eqref{eq:lin_main} is contained in $\cN_i$ \eqref{eq:bous_non}. We derive the ODE of $ \la \eta_1, q f_* \ra$ to control it in the ODE of $\la \om_1, q f_* \ra$. 

The main terms for $\G_i$ are given by $\G_{i, M}$ from the main linearized equations \eqref{eq:lin_main}. Using integration by parts, we get
\bseq\label{eq:ODE_IBP}
\beq
\bal
 - \int   ( \bar c_l y + \bar \uu) \cdot \na g(y) f_*(y) dy
=    \int g(y) \na \cdot (  ( \bar c_l y + \bar \uu ) f_*(y) )  dy 
=  \int g(y) \bar \uu \cdot  \na f_*(y) dy,
\eal
\eeq
where we have used an algebraic property of $f_*$ \eqref{eq:green}: 
\beq
\na \cdot (y f_*(y)) = 0, \quad \na \cdot ( \bar \uu f_*) = \bar \uu \cdot \na f_* ,
\quad 
f_*(y)=\tf{y_1y_2}{|y|^4}.
\eeq
\eseq

We will use \eqref{eq:ODE_IBP} with $g = \om_1, \eta_1, \om ,\eta$. The boundary term near the singularity at $0$ vanishes since $\om_1, \eta_1 , \om ,\eta = O(|x|)$. See \eqref{eq:Sym} and \eqref{eq:reg_W1}.

\begin{remark}[\bf Advantage of the ODEs]
In \eqref{eq:lin_main_cw} with $q=1$, to control $\la \om_1, f_*\ra, \la \eta_1, f_*\ra$, we bound integrals $\la \G_i, f_* \ra = \la (\om_1, \eta_1), g \ra$, where $g$ is either very small (of error size), e.g., $\bar \cF_i, \cN_i$, or decays faster than $f_*$. For example, 
$ \bar \uu \cdot \na f_*$ from \eqref{eq:ODE_IBP} decays faster than $f_*$ since $\bar \uu$ grows sublinearly, $O(|x|^{\g})$ with $\g \approx \f{2}{3}$. Since $\bar c_{\om} \approx -1$, $\bar c_{\om} \la \om_1, f_*\ra$ and  $2 \bar c_{\om} \la \eta_1, f_* \ra$ are damping terms.

\end{remark}

To bound the nonlocal terms in \eqref{eq:lin_main_cw} and in $B_2$ \eqref{eq:bad} involving $\uu_A,\na \uu_A$, we apply the estimates \eqref{eq:u_linf_est22}. For the local terms of $\om_1, \eta_1$
in \eqref{eq:lin_main_cw} and $B_2$ \eqref{eq:bad}, other than $\la \om_1, f_*\ra, \la \eta_1, f_*\ra$, we use the pointwise estimate \eqref{eq:w_est1}. The nonlinear and error terms in \eqref{eq:lin_main_cw} are treated perturbatively, using pointwise estimates \eqref{eq:u_linf_est22}, \eqref{eq:w_est1}, and integration by parts (IBP):
 \[
   \int \VV \cdot \na g(y) \cdot f(y) d y
   = - \int \na \cdot ( \VV f(y) ) g(y) dy 
   = -\int \VV \cdot \na f(y) \cdot g( y ) d y,
 \] 
 where $g$ is the perturbation, $\VV$ satisfies $\na \cdot \VV = 0$, and $f$ is some explicit function.

\vspace{0.1in}
\paragraph{\bf{Improvement} }

We further improve the above estimate by decomposing
\beq\label{eq:cw_chi0}
\la \om_1, f_* \ra = \la \om_1, \chi_{\mf{ode}} f_* \ra 
+  \la \om_1, (1-\chi_{\mf{ode}}) f_* \ra , 
\quad
\la \eta_1, f_* \ra = \la \eta_1, \chi_{\mf{ode}} f_* \ra 
+  \la \eta_1, (1-\chi_{\mf{ode}} ) f_* \ra , 
\eeq
where $\chi_{\mf{ode}}$ is a smooth cutoff function supported \textit{away} from the origin. We derive the ODEs for 
\[
\la \om_1, \chi_{\mf{ode}} f_* \ra, \quad  \la \eta_1, \chi_{\mf{ode}} f_* \ra
\]
using \eqref{eq:lin_main_cw} with 
$q = \chi_{\mf{ode}}$, and perform energy estimates on these terms. The only difference from 
the case $q=1$ is the advection term. Instead of having \eqref{eq:ODE_IBP} for $q=1$, we yield 
\[
- \int (\bar c_l y + \bar \uu) \cdot \na g(y) (f_* \chi_{\mf{ode}})(y) dy 
= \int g(y) \na \cdot ( (\bar c_l y + \bar \uu ) f_* \chi_{\mf{ode}} ) dy.
\]
Outside $\supp(1-\chi_{\mf{ode}})$, we get $\chi_{\mf{ode}} \equiv 1$ and yield the same integrand as \eqref{eq:ODE_IBP}. For $ \la \om_1, (1-\chi_{\mf{ode}}) f_* \ra $, where the integrand is supported near $0$, we 
apply the pointwise estimate \eqref{eq:w_est1} directly. 
We use this decomposition because the estimate via the ODE method is \textit{only} more effective than the pointwise estimate \eqref{eq:w_est1} for controlling the far-field part of the integral $\int \omega_1 f_*$ (where the ODE system's integrand 
decays faster). 
See the above remark. We choose 
\beq\label{eq:cw_chi}
\chi_{\mf{ode}}(x, y) =  1 - \chi_e(  (x- \nu_{31} ) / \nu_{31} ) \chi_e( (y - \nu_{32}) / \nu_{32} ), \quad \nu_{31} = 80, \quad \nu_{32} =  1200,
\eeq
where $\chi_e$ is the cutoff function defined in \eqref{eq:cutoff_exp}. 
Following the estimates discussed above, 
using the pointwise estimates \eqref{eq:u_linf_est22}, \eqref{eq:w_est1} and IBP, we control the main part in the ODEs \eqref{eq:lin_main_cw} 
\beq\label{eq:cw_W1_gamM}
\bal
 \tf{4}{\pi}  |\la \G_{i, M} ,  f_* \chi_{\mf{ode}} \ra| \leq  \mu_{5i, 1} E_3(t), \quad  i = 1,2,
 \eal
\eeq
for some constant $\mu_{5 i, 1}$. 

To explain the choice of weights, first ignore the remaining perturbative
terms $\la \G_i-\G_{i,M},\chi_{ode}f_*\ra$ in \eqref{eq:lin_main_cw}.
With $\mu_{ 5i, 1}$ determined in \eqref{eq:cw_W1_gamM}, we next choose coupling weights so that the forcing terms in the ODE system 
\eqref{eq:lin_main_cw} are perturbative relative to the damping term with coefficients  \(\bar c_\omega\) and \(2\bar c_\omega\). By the analysis of the damped system \eqref{eq:damp_sys} and the stability condition \eqref{eq:PDE_diag},
\footnote{
If $\G_i = \G_{i, M}$, the ODE \eqref{eq:lin_main_cw} and estimates \eqref{eq:cw_W1_gamM} imply 
\[
 \tf{d}{d t} A(t) = \bar c_{\om} A(t) + B(t) + f_1(t), \quad \tf{d}{dt} B(t) = 2 \bar c_{\om} B(t) + f_2(t), 
 \quad A(t) = \tf{4}{\pi} \la \om_1, \chi_{\mf{ode}} f_* \ra, \quad B(t) = \tf{4}{\pi} \la \eta_1, \chi_{\mf{ode}} f_* \ra
 \]
 where $ |f_i(t)| \leq \mu_{5 i,1} E_3(t)$. 
Consider Lemma \ref{lem:PDE_stab} with energy $E_{new}$ in \eqref{eq:cw_W1_mu1}. 
In the stability inequality \eqref{eq:PDE_diag} for $|B(t)|$, $a_{ii} = 2 \bar c_{\om}$, 
$\mu_i = \mu_{52}^{-1}$, the bound $\mu_{52,1}$ for $f_2(t)$ corresponds to 
the off-diagonal part $\sum_{j\neq i} a_{ij} \mu_j^{-1}$. Since $\bar c_{\om}<0$, \eqref{eq:PDE_diag} with $\lam = 0^+$ 
reduces to $2 |\bar c_{\om}| - \mu_{52}^{-1} \mu_{52,1} > 0$, which is the first bound 
in \eqref{eq:cw_W1_mu1:a}. 
The second bound in \eqref{eq:cw_W1_mu1:a} follows from $|B(t) + f_1(t)| \leq (\mu_{52} + \mu_{51,1} ) E_{new}$ and a similar argument.
} 
we can first choose a large $\mu_{52}$ and then $\mu_{51}$ so that 
\bseq\label{eq:cw_W1_mu1}
\beq\label{eq:cw_W1_mu1:a}
\mu_{52}^{-1} \mu_{52, 1} <  |2  \bar c_{\om}|  , \quad  \mu_{51}^{-1} ( \mu_{51, 1} + \mu_{52}  ) 
 < |\bar c_{\om}|.
\eeq

After multiplication by the small coupling weights $\mu_{52}^{-1}$ and
$\mu_{51}^{-1}$, the forcing bounds in \eqref{eq:cw_W1_gamM} become
perturbative. This yields linear stability for the enlarged energy $E_{new}$:
\beq
\bga
   E_{new} = \max(E_3, \mu_{51}^{-1} m_1, \mu_{52}^{-1} m_2  ), \\
    m_1 \teq | \tf{4}{\pi} \la \om_1, \chi_{\mf{ode}} f_* \ra| = |c_{\om}( \om_1 \chi_{\mf{ode}} ) |,
   \quad m_2 \teq | \tf{4}{\pi} \la \eta_1, \chi_{\mf{ode}} f_* \ra| = |c_{\om}( \eta_1 \chi_{\mf{ode}} ) | ,
   \ega
\eeq
\eseq
where we have used the notation \eqref{eq:green} for $c_{\om}(\cdot)$.
To close nonlinear estimates \eqref{eq:PDE_nondiag}, 
due to the remaining terms, we choose a slightly smaller weights. The weights $\mu_{51}, \mu_{52}$ chosen in \eqref{wg:EE} are mostly determined by the above linear stability constraints.

\subsubsection{Control  $c_{\om}(\om)$}\label{sec:cw_w}

Recall that $W_1 + \widehat W_2 = (\om, \eta, \xi)$ is the solution to \eqref{eq:lin} and 
\[
 W_1 = (\om_1, \eta_1, \xi_1),  \ \hat W_2 =(\hat \om_2, \hat \eta_2, \hat \xi_2),
 \quad  W_{1, 1} = \om_1, \ \hat W_{2, 1} = \hat \om_2, 
 \quad 
 c_{\om}( \om) = c_{\om}( \om_1 ) + c_{\om}( \hat \om_2 ). 
\]

We use similar ideas to estimate $c_{\om}(\om)$ by deriving the ODEs of $c_{\om}(\om_1)$ and $c_{\om}( \hat \om_2 )$ separately. We have derived the ODE of $c_{\om}(\om_1)$ in \eqref{eq:lin_main_cw} with $q=1$. 
Below, we derive the ODE for $c_{\om}( \hat \om_2 )$.

Recall the functionals $a_i(W_1), a_{ \mf{nl}, i}(W_1 + \hat W_2)$ 
and $\hat F_i, \hat F_{\chi ,i} \in \R^3$ from \eqref{eq:bous_W2_appr}. Denote
\beq\label{def:ode_bi}
b_i(s) \teq a_i(W_1(s)), \  1 \leq i \leq n_1,  \quad  b_{n_1 + j}(s) \teq a_{\mf{nl}, j}( W_1(s) + \hat W_2(s) ) ,
\ \hat F_{n_1+j} \teq \hat F_{\chi, j} , 
 \  j \leq 3  , 
 \eeq
where we abuse notation $\hat F_{n_1 +j}$ for $j=1,2,3$ to simplify notation. 
Recall that $\hat F_i, \hat F_{\chi ,i}$ are \emph{vector-valued} rather than 
scale-valued functions. Below, we use $\hat F_{i, 1} \in \R$ to indicate the first component 
of $\hat F_{i} \in \R^3$. Since $b_i$ is spatially independent, we use the formula \eqref{eq:bous_W2_appr} and linearity to get 
\beq\label{eq:lin_main_cw_w}
c_{\om}(\hat \om_2) =  c_{\om}(\widehat W_{2, 1}) = \sum\nolim_{i \leq \nmod} \int_0^{ t \weg T_i} b_i(t-s) c_{\om}(\widehat F_{i, 1} )( s) ds
 \teq \sum\nolim_{i \leq \nmod } I_i,  
 \quad \nmod =n_1 + 3.
\eeq
We restrict the integral to $s\leq T_i$ 
since we construct $\hat F_i$ with  $c_{\om}(\widehat F_{i, 1})(s)$ supported in $[0,T_i]$, although
$\widehat F_i$ itself has global-in-time support.
\footnote{
For $t > T_i$, we represent $\hat F_i(t) = \al_i(t) f(x)$ 
for some \emph{time-independent} function $f(x)$ with $c_{\om}(f)=0$. 
This term is the same as the analytic correction $a(t) \chi_2$ in \eqref{def:g2}
in Section \ref{sec:numer}. See details in Part II \cite[Section 3.5]{ChenHou2023b}. 
}
As chosen in \eqref{eq:appr_near0}, we label the first term as 
\beq\label{eq:ODE_cw_b1}
b_1(s) = a_1(W_1(s)) = c_{\om}( \om_1(s)), \quad 
b_1(s) \hat F_1  = c_{\om}( \om_1(s))\hat F_1, 
\quad  \om_1 = W_{1, 1}, 
\eeq
where $\hat F_1$ is the approximate solution from initial data $ (\bar f_{c_{\om}, 1}, 
\bar f_{c_{\om}, 2}, \bar f_{c_{\om}, 3})$.
 For each $I_i$ defined in \eqref{eq:lin_main_cw_w},  taking derivatives and using $\pa_t b_i( t-s ) = -\pa_s b_i( t-s )$ and integration by parts, we yield  
\begin{align}\label{eq:ode_cw_Ibi}
\f{d}{dt} I_i 
& = \one_{t < T_i} b_i(t-  t \weg T_i ) c_{\om}(\widehat F_{i, 1}(t \weg T_i ))
+ \int_0^{ t \weg T_i} - \pa_s b_i( t-s)  c_{\om} ( \widehat F_{i, 1}(s)) d s  \notag \\
& =  \one_{t < T_i} b_i( t- t \weg T_i ) c_{\om}(\widehat F_{i, 1}(t \weg T_i)) 
+ \int_0^{t\wedge T_i} b_i( t-s ) \pa_s c_{\om} (\hat F_{i, 1}(s)) ds  \notag \\
& \quad  
- b_i( t - t\weg T_i) c_{\om}( \widehat F_{i, 1}(t \weg T_i ))
 + b_i(  t ) c_{\om} ( \widehat F_{i, 1}(0) )  \notag \\
 & =   b_i(  t)  c_{\om} ( \hat F_{i, 1}(0) ) + \B( - \one_{t \geq T_i} b_i( t - T_i )  c_{\om}( \widehat F_{i, 1}(T_i)) + \int_0^{t\wedge T_i} b_i( t-s) \pa_s c_{\om} (\widehat F_{i, 1}(s)) ds \B)  \notag \\
& \teq  b_i(  t ) c_{\om} ( \widehat F_{i, 1}(0) ) + I_{b, i}.
\end{align}
The term $I_{b, i}$ is treated as a bad term. For $i=1$, the 
first mode $b_1(t) c_{\om} ( \hat F_{1,1}(0) )$ with $c_{\om} ( \hat F_{1,1}(0) ) \approx -2.5$ provides an additional damping term for $b_1(t) = c_{\om}( \om_1(t) )$.
To exploit this damping term, below, we combine the ODE estimate of $c_{\om}( \om_1)$ and $c_{\om}(\hat \om_2)$.
 Denote 
\[
\bar \lam_{c_{\om}} = c_{\om}( \widehat F_{1, 1}(0)) + \bar c_{\om}, \quad 
B_{c_{\om}, 1} 
= \tf{4}{\pi}  \la \G_1, f_* \ra ,
\]
where $\G_1$ is defined in \eqref{eq:lin_main_cw}.  Since $c_{\om}(g) = - \f{4}{\pi} \la g , f_* \ra$ \eqref{eq:green}, multiplying the ODE of $\la \om_1, f_* \ra$ in \eqref{eq:lin_main_cw} with $q=1$ by $-\f{4}{\pi}$ and combining it with derivations 
\eqref{eq:lin_main_cw_w}, \eqref{eq:ode_cw_Ibi} for $c_{\om}(\hat \om_2)$, we get
\[
\bal
\f{d}{dt} ( c_{\om}( \om_1) + c_{\om}(\hat \om_2))
 & = \bar c_{\om} c_{\om}( \om_1) - B_{ c_{\om}, 1} - \f{4}{\pi}\la\eta_1, f_*\ra
+ c_{\om}( \hat F_{1,1}(0)) b_1(t)
+ \sum_{i \geq 2} c_{\om}(\hat F_{i, 1}(0)) b_i(t)
+ \sum_{i \geq 1} I_{b, i} .\\
\eal
\]

Since $b_1(t) = c_{\om}(\om_1(t))$ by \eqref{eq:ODE_cw_b1}, using $\om= \om_1 + \hat \om_2$, formula of $ c_{\om}( \hat \om_2)$ in \eqref{eq:lin_main_cw_w}, we rewrite
\[
\bar c_{\om} c_{\om}(\om_1) + c_{\om}( \hat F_{1,1}(0)) b_1(t)  = \bar \lam_{ c_{\om} } c_{\om}( \om_1)
= \bar \lam_{ c_{\om} } c_{\om}(\om ) - \bar \lam_{ c_{\om} }
\tts{ \sum_{i\leq n_1 + 3}} \,
 \tts{ \int_0^{t\wedge T_i} } b_i( t-s)   c_{\om}(\hat F_{i,1}(s))  ds .
\]

Thus, using the above derivations and the formula of $I_{b, i}$ in \eqref{eq:ode_cw_Ibi}, we derive
\begin{align}\label{eq:ode_cw_all}
&\f{d}{dt} c_{\om}(\om)
 = \bar \lam_{c_{\om}}  c_{\om}(\om)
- B_{c_{\om}, 1} - \f{4}{\pi}\la\eta_1, f_*\ra
+ \sum\nolim_{i\geq 2} c_{\om}(\hat F_{i, 1}(0)) b_i( t) \\
&  + \sum\nolim_{i \geq 1} \B(- \one_{t \geq T_i} b_i(  t - T_i  )  c_{\om}( \hat F_{i, 1}(T_i)) + \int_0^{t\wedge T_i} b_i( t-s) ( \pa_s c_{\om} (\hat F_{i, 1}(s)) - \bar \lam_{c_{\om}} c_{\om}(\hat F_{i, 1}(s)) ) ds \B). \notag
\end{align}

We estimate \eqref{eq:ode_cw_all} following the \eqref{eq:damp_sys} system. 
We treat $\bar \lam_{c_{\om}}  c_{\om}(\om)$ as the damping term, and bound other terms perturbatively. We estimate $B_{ c_{\om}, 1}$ using the method in Section \ref{sec:cw_w1}.  
We bound $b_i(s), i\geq 2$ using \eqref{eq:finite_mode_2nd} and energy $E_3$ with established linear stability 
and bound $b_1(s) = c_{\om}(\om_1(s))$ using the enlarged energy of $E_3$ including 
$c_{\om}(\om_1)$ determined by Section \ref{sec:cw_w1}. The nonlinear modes $b_i$ are much smaller and bounded using bootstrap assumption.

 For $\f{4}{\pi}\la\eta_1, f_*\ra $, we use ODE \eqref{eq:lin_main_cw} and the method in Section \ref{sec:cw_w1}. The approximate terms 
\[
 c_{\om}( \hat F_{i, 1}(T_i)), \quad  \pa_s c_{\om} (\hat F_{i, 1}(s)) - \bar \lam_{c_{\om}} c_{\om}(\hat F_{i,1}(s))
 \]
are piecewise cubic polynomials in $s$ constructed numerically, and are estimated by the method in Section \ref{sec:W2}. We have estimated $b_1(t) = c_{\om}( \om_1(t))$ in Section \ref{sec:cw_w1} at the linear level. Using the above estimates and following the discussion around \eqref{eq:cw_W1_mu1}, at the linear level, we determine the coupling weights for $c_{\om}(\om)$ and $\la \eta_1, f_* \ra$ in the energy.

If we estimate $c_{\om}(\om_1), c_{\om}( \hat  \om_2)$ separately, we need to add a much 
 smaller coupling weight for $c_{\om}(\om)$ in the energy, leading to a constant about 
\textit{three times} larger for the nonlinear estimates.

\subsubsection{Control $\om_{xy}(0), \eta_{xy}(0), \xi_{xx}(0)$}\label{sec:wxx0}

First, note that $ \eta_{xy}(0) =  \xi_{xx}(0) = \th_{xxy}(0)$ as $\eta =\th_x, \xi = \th_y$ \eqref{eq:eta_xi}.
Recall the ODEs for the full solution in \eqref{eq:bous_ODE_xy}. Linearizing it around the profile 
$(\bar \om, \bar \th)$ and using the normalization conditions \eqref{eq:normal_pertb}, \eqref{eq:normal_vanish}, yields the ODEs for the perturbations 
\footnote{
For initial data $\om - \bar \om, \na \th - \na \bar \th  = O(|x|^2)$, we have  $ \om_x(0) = \bar \om_x(0),  \th_{xx}(0) = \bar \th_{xx}(0)$ 
dynamically by normalization \eqref{eq:normal1}.  The terms $ \om_x(0), \th_{xx}(0)$ in ODEs \eqref{eq:bous_ODE_xy}
are combined into the error $\pa_{xy} \bar \cF_1,   \pa_{xy} \bar \cF_2(0)$ in \eqref{eq:ODE_lin_xy}.
}
\beq\label{eq:ODE_lin_xy}
\bal
\tf{d}{dt} \om_{xy}(0) & = ( -2  \bar c_l + \bar c_{\om}) \om_{xy}(0) + \th_{xxy}(0) 
+ c_{\om} \bar \om_{xy}(0) + c_{\om} \om_{xy}(0) +  \pa_{xy} \bar \cF_1(0) , \\
 \tf{d}{dt} \th_{xxy}(0) & = ( -2  \bar c_l + 2 \bar c_{\om} -\bar u_x(0)) \th_{xxy}(0) +  c_{\om} \bar \th_{xxy}(0) + c_{\om} \th_{xxy}(0) +  \pa_{xy} \bar \cF_2(0), \\
 \eal
\eeq
where $\bar \cF_i$ is defined in \eqref{eq:bous_err}. The matrix involving $\om_{xy}(0), \th_{xxy}(0)$ has negative eigenvalues. We first estimate $\th_{xxy}(0)$ and then $\om_{xy}(0)$
following \eqref{eq:damp_sys}. Using the above ODEs, at the linear level, we determine the coupling weights for $\th_{xxy}(0), \om_{xy}(0)$ in the energy for stability,

To handle the nonlinear and error terms in \eqref{eq:lin_main_cw}, the ODE of $c_{\om}$, and \eqref{eq:ODE_lin_xy}, we choose the weights of the functionals in Sections \ref{sec:cw_w1}, \ref{sec:cw_w} slightly smaller than those determined by the linear estimates. The final energy \(E_4\) is obtained by adding scalar quantities that precisely control the functionals appearing in the finite-rank correction and the residual operator. Define
\beq\label{energy4}
\bal 
E_{4 ,0}(W_1, W) & \teq   \max(  \mu_{51}^{-1}  |c_{\om}( \om_1 \chi_{\mf{ode}})|, 
\mu_{52}^{-1} | c_{\om}(\eta_1 \chi_{\mf{ode}} ) |,  \mu_{6}^{-1} | c_{\om}(\om)|,  \mu_{62}^{-1} |c_{\om}(\eta_1)|,    \\ 
&  \qquad \qquad  \mu_7^{-1}  |\th_{xxy}(0)|, \mu_8^{-1} |\om_{xy}(0)| ), 
\qquad \quad \mu_6 =  61  , \  \mu_7 =  9.5, \ \mu_8 = 4.5 ,  \\
 E_4(W_1, W) & \teq  \max(  E_3(W_1),E_{4,0}(W_1, W) ) .
\eal 
\eeq 
When no confusion arises, we write 
$E_4,  E_4(t)$ for $E_4(W_1, W), E_{4}(W_1(t), W(t))$, respectively. 
Here, we define energy $E_3$ in \eqref{energy3}, $\chi_{\mf{ode}}$  in \eqref{eq:cw_chi}, $\mu_{ij}$  in \eqref{wg:EE}, and functional $c_{\om}(\cdot)$  in  \eqref{eq:green}. 

The variables $c_{\om}(\om_1 \chi_{\mf{ode}} ), c_{\om}(\eta_1 \chi_{\mf{ode}} ), c_{\om}(\eta_1)$ and parameters $\mu_{51}, \mu_{52}, \mu_{62}$ are \emph{intermediate variables} and are used \emph{only} in the ODEs in Sections \ref{sec:cw_w1}, \ref{sec:cw_w} along with \eqref{eq:cw_chi0} to control
\beq\label{eq:cw_est}
|c_{\om}(\om_1)| < \mu_5 E_4 , \quad  |c_{\om}(\om)|  < \mu_6 E_4  , \quad  \mu_5 = 76.
\eeq

To estimate the nonlinear mode $a_{ \mf{nl},i}(W_1 + \hat W_2)$ for \eqref{eq:appr_near0}, we impose the bootstrap bound 
\beq\label{eq:W2_non_boot}
 |c_{\om } \om_{xy}(0) + \pa_{xy} \bar \cF_1(0)| <  5 \mu_6 \cdot E_*, 
 \quad  |c_{\om } \th_{xxy}(0) + \pa_{xy} \bar \cF_2(0)| <  10 \mu_6 \cdot E_*, 
\eeq
where $E_*$ will be chosen in \eqref{eq:Ec}. Under the bootstrap assumptions $E_4(t) < E_*$, we 
will verify the stronger estimate  (see Appendix \ref{sec:ineq_cw} )
\beq\label{eq:W2_non_boot2}
\bal
  |c_{\om } \om_{xy}(0) + \pa_{xy} \bar \cF_1(0)| 
 & < \mu_8 \mu_6 E_* + |\pa_{xy} \bar \cF_1(0)|
 < 5 \mu_6 E_* \, ,  \\
   |c_{\om } \th_{xxy}(0) + \pa_{xy} \bar \cF_2(0)| 
 & < \mu_7 \mu_6 E_* + |\pa_{xy} \bar \cF_2(0)|
 < 10 \mu_6 E_* \, .
 \eal
\eeq

\subsection{Estimate $\wh W_2$ and the residual operator}\label{sec:W2}

Using the method in part II \cite[Section {\seclinevo}]{ChenHou2023b}, for initial data $\bar G_i(0) 
\in \{  \bar F_i(0), \bar F_{\chi, i} \}$ given in Appendix \ref{app:init_data} and the spatially independent factors $b_i(s) = a_i( W_1(s)), a_{\mf{nl}, i}( W(s)) $ in \eqref{eq:bous_err_op}, \eqref{eq:bous_W2_appr} (see also \eqref{def:ode_bi}), we construct an approximate space-time solution $\hat F_i, \hat W_2$ and its associated approximate stream function $\hat \phi^N$ and error $\hat \e$
\beq\label{eq:lin_evo_main1}
 \hat W_2(t) = \sum_{ i\leq \nmod} \int_0^t b_i(t- s) \hat G_i(s) ds , \  \hat \phi^N(t) = \sum_{ i\leq \nmod } \int_0^t b_i(t- s) \hat \phi^N_i(s) ds ,  \  \hat \e(t) =  \sum_{ i\leq \nmod } \int_0^t b_i(t- s) ( \hat G_i + \D \hat \phi^N_i )(s) ds .
\eeq
where $\nmod$  \eqref{eq:lin_main_cw_w} is the number of modes.
The residual error in $j$-th equation reads
\footnote{
$ D^2$ relates to leading order behaviors of $W$ near $x=0$, e.g. $D_i^2W_i(0)$,
   and \textit{does not} denote mix derivatives. 
}
\beq\label{eq:W2_err_est1}
\bal
    \cR_j(t)  &=  \cR_{loc, 0,  j}(t) -  ( \cR_j^{\bar \e} - D_j^2 \cR_j^{\bar \e}(0) \chi_{j 2} ) - ( \cR_j^{\hat \e} - D_j^2 \cR_j^{\hat \e}(0) \chi_{ j2}) \,  , \\
  \cR_{loc, 0, j} & = \sum\nolim_{i\leq \nmod} \int_0^t  b_i(t-s)  \cR_{\mf{num} ,i,  j}(s) ds , \\
\cR_j^{\bar \e} & = \cB_{op, j}( \uu( \bar \e), \hat W_2  ) , \quad  
   \cR_j^{\hat \e} = \cB_{op, j}( \uu(\hat \e), \bar W) , 
   \quad  D^2 = (\pa_{xy}, \pa_{xy}, \pa_{xx})   ,
   \quad j=1,2,3 .
   \\
\eal
\eeq
Here, $\cB_{op,j}$ denotes the bilinear operator defined in \eqref{eq:Blin}; 
we defer its formula since it  is not used here.
The functions  $\bar W, \chi_{j2}$ are given in \eqref{eq:EE_W1}, \eqref{eq:cutoff_near0_all}, and $\bar \e  = \bar \om - (-\D) \bar \phi^N$ will be defined in \eqref{eq:wb_vel_err}. 
The residual $\cR_{\mf{num}, i, j}= O(|x|^3)$ depends on the numerical solution $\hat \phi_i^N, \hat G_i$ \emph{locally}, and we absorbed the initial error $\hat G_i(0) - \bar G_i$ (see \eqref{eq:bous_err_op}) in $\cR_{\mf{num}}$. We remark that $\cR_j^{\bar \e}, \cR_j^{\hat \e}$ are \emph{nonlocal errors} since $\uu(\e)$ depends on $\e$ nonlocally. Moreover, in Part II \cite[Section 3]{ChenHou2023b}, we estimated 
\beq\label{eq:W2_res0}
   \int_0^{\infty} |\pa_x^k \pa_y^\ell f (s)| ds , \quad f  = \hat F_{i, j}, \ \hat \phi_i^N, \  
   \hat G_i, \  \hat G_i + \D \hat \phi_i^N , \ \cR_{\mf{num}, i, j}. 
\eeq

For later estimates, we add and subtract $f_{\chi, j}$ in the $\cR_j^{\bar \e}$ terms in \eqref{eq:W2_err_est1}
\footnote{
  The correction functions $f_{\chi, j}$ and $\chi_{j2}$ are different; see definitions \eqref{eq:cutoff_near0_all}. The subscript $2$ in $\chi_{j2}$ indicates the second type of cutoff function. We also introduced the first type, $\chi_{j1}$, in \cite{ChenHou2023b}.
}
\beq\label{eq:W2_err_est2}
\bal
  \cR_j(t) =  \cR_{loc, 0,  j}(t) - D_j^2 \cR_j^{\bar \e}(0) (f_{\chi, j} - \chi_{j 2} )
  - ( \cR_j^{\bar \e} - D_j^2 \cR_j^{\bar \e}(0) f_{\chi , j} ) -
( \cR_j^{\hat \e} - D_j^2 \cR_j^{\hat \e}(0) \chi_{j2} ) .
\eal 
\eeq

We control the spatial-independent factor $b_i(s)$ using the energy estimate  in Section \ref{sec:rank1} :
\[  
 |b_i(s)| < c_{i, 1}  E_4(s) + c_{i, 2} E_4(s)^2.
\]
In the bootstrap argument we will prove
$E_4(s)< E_*, \forall   s \geq 0$. 
Under such an assumption we have 
\beq\label{eq:W2_bi}
|b_i(s) | \leq c_i E_*, \quad  c_i = c_{i, 1}  + c_{i, 2} E_*  .
\eeq
for some threshold $E_*$ to be determined. Then we can control the local terms, e.g.
\begin{align}\label{eq:W2_res1}
| \pa_x^k \pa_y^\ell \cR_{loc, 0, j}(t)|
 &\leq  \sum_{i \leq \nmod } \sup_{s < t}|b_i(s)|  \int_0^{t} |\pa_x^k \pa_y^\ell \cR_{\mf{num},i, j}(s)| ds \leq  E_*   \sum_{i \leq \nmod} |c_i|  \int_0^{\infty} |\pa_x^k \pa_y^\ell \cR_{\mf{num},i, j}(s)| ds, 
 \, k + \ell \leq 2 , 
 \notag \\
|| \pa_x^k \pa_y^\ell \wh W_2(t) | 
&\leq  E_*  \sum_{i\leq \nmod } |c_i|  \int_0^{\infty} |\pa_x^k \pa_y^\ell \hat G_i(s)| ds, \quad \forall \, k + \ell \leq 5.
\end{align}
uniformly in time $t$ using monotonicity. The error and similar quantities are integrable in time since the approximate solution $\hat G_i(s)$ and residual error $ \cR_{\mf{num},i}$ can be decomposed into $A_i(t, x)+ B_i(t) f(x)$ with $A_i(t, x)$ compactly supported in time, and $B_i(t)$ decays exponentially fast in $t$.  \footnote{
The term $B_i(t) f(x)$ is the same as the analytic correction $a(t) \chi_2$ in \eqref{def:g2}
in Section \ref{sec:numer}. 
 $B_i(t)$ solves an ODE and admits an \emph{exact} formula as \eqref{def:a} with 
$\bar \lam < 0$ and a source  $\pa_{xy} \cE(s, 0)$ compactly supported in $s$. 
Thus, $B_i(t)$ decays exponentially. 
See  Part II \cite[Section 3.5]{ChenHou2023b} for more discussions. 
}
 Thus the estimate of $\widehat W_2$ and $\cR$ decouples into two independent pieces: the scalar coefficients \(b_i(t)\) \eqref{def:ode_bi}, controlled by \(E_4\) and estimate \eqref{eq:W2_bi}, and the precomputed time-integrated norms of the approximate space-time solutions and error, controlled in
Part~II \cite[Section 3]{ChenHou2023b}. 
This allows us to rigorously integrate the finite-rank component in the nonlinear stability argument.

\begin{remark}[\bf Linearity and decoupled estimates]\label{rem:W2_para}
  By the above decoupled structure, in practice, 
 we can pre-compute and estimate each mode, 
 e.g. $\hat G_i, \cR_{\mf{num}, i,j}$ in \eqref{eq:W2_res1}, completely \emph{in parallel}. 
Using linearity, triangle inequality, and \eqref{eq:W2_res1}, we further estimate $\cR$ and $\hat W_2$.
\end{remark}

\begin{remark}[\bf Effectiveness of estimates \eqref{eq:W2_res1}]

Although we estimate $ \pa_1^{\al} \pa_2^{\b} \wh W_2, \pa_1^{\al} \pa_2^{\b} \cR $ \eqref{eq:lin_evo_main1}, \eqref{eq:W2_err_est1} by applying triangle inequality and combining the estimates of different modes \eqref{eq:W2_res1}, such estimates do not necessarily yield a constant of $O(\nmod)$, since different solutions are large in different regions. When constructing approximations for the velocity in Section \ref{sec:appr_vel}, we apply some partition of unity. The coefficients of different approximations 
(serving as initial conditions for the approximate solution \eqref{eq:bous_W2_appr}, cf. Appendix \ref{app:init_data}) are large in different regions. This structure prevents large constants in the above estimates for $\cR, \wh W_2$ \eqref{eq:W2_res1}
\end{remark}

\subsection{Estimate of nonlocal error and modified decomposition}\label{sec:comb_vel_err}

The residual operator $\cR_j$ in \eqref{eq:W2_err_est1} and profile residual  $\bar \cF_i$ 
\eqref{eq:bous_err} in the numerical construction consists of local errors, nonlocal errors, and round-off errors. 
Below, we  decompose these errors and discuss their estimates.

\vspace{0.05in}

\paragraph{\bf Nonlocal error}

To construct the approximate stream function $\bar \phi^N$ and velocity $\na^{\perp} \bar \phi^N$ of $\bar \om$,
we numerically solve $-\D  \phi = \bar \om$  (see 
Section \ref{sec:ASS}). This induces a local error $ \bar \e= \bar \om + \D \bar \phi^N$
and nonlocal error $\na^{\perp}(-\D)^{-1} \bar \om - \na^{\perp} \bar \phi^N$. In the residual operator, we have a similar error $\hat \e$ in \eqref{eq:lin_evo_main1}. 

\vspace{0.05in}

\paragraph{\bf Local error}

The local error depends on the numerical bases (e.g., piecewise polynomials and semi-analytic functions) \emph{locally} in space and includes 
\[
\rm{all \  terms \ in \  \eqref{eq:W2_res0} }, \quad 
  \cR_{loc, 0, j} \ \eqref{eq:W2_err_est1}, \quad \bar \e = \bar \om + \D \bar \phi^N ,\quad  \hat \e
  \ \eqref{eq:lin_evo_main1}
   .
\]
We estimate them using numerical analysis and the methods in Part II \cite[Section {\secerridea}]{ChenHou2023b}.

We estimate the nonlocal error using \emph{functional inequalities} 
and discuss it below. Denote 
\beq\label{eq:wb_vel_err}
\bar \e = \bar \om - (-\D) \bar \phi^N,  \quad \bar \uu^N = \uu( -\D \bar \phi^N) = \na^{\perp} \bar \phi^N,  \quad \bar c_{\om}^N = \bar u^N_x(0) + \tf{1}{2} \bar c_l.
\eeq
We introduce $\bar c_{\om}^N$ since $\bar c_{\om}$ in \eqref{eq:normal} depends on $u_x(\bar \om)(0)$. Recall the  approximations $\hat \uu, \wh { \na \uu}$ \eqref{eq:u_appr}. To apply functional inequalities similar to \eqref{eq:u_linf_est1} for $\uu(\bar \e)$, the input must vanish to higher order than $O(|x|^2)$ 
near $x=0$, so that the singularly weighted norm is finite. 
Since  $ \bar \phi^N$ satisfies only $\bar \e  = O(|x|^2)$ 
(see \eqref{eq:poisson_err}), we correct $\bar \e$ near $0$. Similar consideration applies to $\hat \e$ \eqref{eq:lin_evo_main1}.

In this subsection, we rewrite \eqref{eq:bous_decoup2}, \eqref{eq:lin_main} so that every term is compatible with the singularly weighted estimates. First, we decompose the elliptic residuals so that 
the part for weighted estimates vanish cubically. Second, we combine the nonlocal velocity errors with the perturbation velocity into $\UU$ \eqref{eq:U}. Third, we rewrite the full  equations in \eqref{eq:lin_U}, where all residual errors vanish cubically, and the error and coefficients of the linearized operators are essentially local.

\vs{0.05in}
\paragraph{\bf Improved vanishing and decomposition}
For an error $\e = \bar \e$ or $\hat \e$, 
we choose a locally $C^4$ cutoff function $\chi_{ \e} = 1 +O(|x|^4)$ near $0$ and decompose $\e$ and $\uu(\e)$:
\footnote{
  We choose $\chi_{ \e}, \e_2$ in this specific form so that $ \e_2 =  \e_{xy}(0) xy + O(|x|^4)$ and we can obtain $\uu(\e_2 )$ explicitly. 
}
\beq\label{eq:uerr_dec1}
\bal
    \e   &=  \e_1 +  \e_2, \quad  \e_2 =  \e_{xy}(0) \D( \tf{1}{6} x y^3  \chi_{  \e}) , \\
 \uu(  \e_2) & = \na^{\perp}(-\D)^{-1}  \e_2 =  - \tf{1}{6}   \e_{xy}(0)\na^{\perp}  ( x y^3 \chi_{  \e}) ,  \\
   \uu(  \e)  & = \uu(  \e_1) + \uu(  \e_2)  = \uu_A(  \e_1) +  \hat \uu(  \e_1 ) + \uu(  \e_2) , \quad u_x(\e)(0) = u_x(\e_1)(0) . \\
 \eal
\eeq
Then $\e_1 = O(|\xx|^3)$. Here we use the decomposition $\uu = \hat \uu + \uu_A$ from 
\eqref{eq:vel_rem} and 
\[ 
u_x(\e_2)(0) = -\pa_{xy}(-\D)^{-1} \e_2(0)= \tfrac16 \e_{xy}(0) \cdot \pa_{xy} (x y^3 \chi_{\e}) |_{x, y=0} =0.
\]
We perform a similar decomposition for $\na \uu(  \e)$. 
We choose cutoffs $\chi_{\bar \e}$ for $\bar \e$ and $\chi_{\hat \e}$ for $\hat \e$ in \eqref{eq:cutoff_near0_all}; they have \emph{different parameters}.

Let $\om$ be the total perturbation without decomposition. Recall  from \eqref{eq:u_tilde}, \eqref{eq:u_appr}, \eqref{eq:vel_rem}:
\[
\uu(f) = \uu_A(f) + \hat \uu(f)= \uu_A(f) + u_x(f)(0)(x, -y) + \appo{ \uu}(f) .
\]
We combine the nonlocal terms of error and perturbations and define 
\footnote{
  Note that we do not include $\uu( -\D \hat \phi^N)$ in $\UU_A$.
}
\bseq\label{eq:U}
\beq\label{eq:U:a}
\bal
\UU & = \uu( \om + \bar \e), \quad  \td \UU = \UU - U_x(0)( x, -y), \quad  \UU_A = \uu_A(\om_1 + \bar \e_1 + \hat \e_1) ,  \\
 \UU_{\FM} & = \appo{\uu}(\om_1 + \bar \e_1 + \hat \e_1) + \uu( \bar \e_2 + \hat \e_2)  
+ \td \uu( - \D \hat \phi^N) , \\
\eal
\eeq
where the boldface font denotes vectors, e.g., $\uu = (u, v), \UU = (U, V)$. 
By definition, we have 
\footnote{
Since $u_x( \bar \e_2) (0) = u_x(\hat \e_2)(0) = 0$, we have $\uu(\bar \e_2 + \hat \e_2) = \td \uu(\bar \e_2 + \hat \e_2)$. Moreover, we write $- \D \hat \phi^N = \hat \om_2 - \hat \e$. Then 
 \[
 \bal
 \UU_A + \UU_{\FM} 
 &= (\uu_A + \appo{\uu} )(\om_1 + \bar \e_1 + \hat \e_1)
 + \uu( \bar \e_2 + \hat \e_2 ) + \td \uu( - \D \hat \phi^N )   = \td \uu( \om_1 + \bar \e_1 + \hat \e_1 ) 
+ \td \uu ( \bar \e_2 + \hat \e_2 ) 
+ \td \uu( \hat \om_2 - \hat \e ) \\
& = \td \uu( \om_1 + \hat \om_2 + \bar \e  + \hat \e- \hat \e )
= \td \uu( \om + \bar \e )
= \td \UU = \UU - U_x(0)(x, -y).
\eal
\]
}
\beq
\bal
\UU = \UU_A + U_x(0)(x, - y) + \UU_{\FM} .
\eal
\eeq
We use $\UU_A$ to denote the velocity remainder after subtracting the
finite-rank approximation, similar to $\uu_A$ in \eqref{eq:vel_rem}.
We use $\UU_{\FM}$ to denote \emph{finite-mode} velocity part after separating the affine term \(U_x(0)(x,-y)\). Each 
component of $\UU_{\FM}$ is  either given by a finite rank operator of $\om_1$,
or  has a finite-dimensional range. 
\footnote{
The range of $\td \uu(-\D \hat \phi^N)= \na^{\perp} \hat \phi^N 
+ \hat \phi^N_{xy}(0)(x, -y) $ is spanned by the numerical basis functions.
}
Similarly, we decompose $ \na \UU$ and define $ (\na \UU)_A, (\na \UU)_{\FM}$. 
\beq
\bal  
\na \UU & =   (\na \UU)_A  + U_x(0) \mw{diag}(1, -1) + (\na \UU)_{\FM},
\quad  (\na \UU)_A \teq (\na \uu)_A(\om_1 + \bar \e_1 + \hat \e_1), \\ 
\quad (\na \UU)_{\FM} & \teq \appow{(\na \uu)}( \om_1 + \bar \e_1 + \hat \e_1 )
+ \na \uu (\bar \e_2 + \hat \e_2) + \na \td \uu( - \D \hat \phi^N) .
\eal 
\eeq 
\eseq 
From \eqref{eq:normal}, \eqref{eq:normal1} and \eqref{eq:U},
for $\uu(f) \teq \na^{\perp}(-\D)^{-1} f$, 
we decompose 
the exact velocity for the profile $\bar \uu =\uu(\bar \om)$ 
into the numerical approximation $\bar \uu^N$ and error \eqref{eq:wb_vel_err}
\beq\label{eq:uerr_dec12}
\bar \uu \teq \uu(\bar \om) =\bar \uu^N + \uu(\bar \e), 
\quad \uu(\om + \bar \om) = \bar \uu^N + \UU, \quad  
\bar c_{\om} + c_{\om}(\om) = \bar c_{\om}^N + U_x(0) .
\eeq

The term $\UU_{\FM}$ is more regular. The approximation terms $\hat \uu(  \e_1 ), \td {\hat \uu}(\e_1),  \uu(  \e_2)$ for $\e = \bar \e, \hat \e$ only depend on $ \e$ via $ \e_{xy}(0)$ and finite many integrals ($<50$) \eqref{eq:u_appr} with smooth coefficients. We estimate the piecewise bounds of $\e$ following  Part II \cite[Section {\secerridea}]{ChenHou2023b}, then estimate these integrals and piecewise $C^3$ bounds of these terms.
We estimate $\td {\hat \uu }(\om_1)$ following Section \ref{sec:rank1} and $\uu(\hat \phi^N) = \na^{\perp} \hat \phi^N $ using the bounds of \eqref{eq:W2_res0} and method in \eqref{eq:W2_res1}. We 
isolate  $U_x(0) (x, -y)$ in \eqref{eq:U} since our estimate for such a term is larger than others. See Section \ref{sec:non_mainT}.

\vs{0.05in}

\paragraph{\bf Functional inequalities for $\uu_A(\e_1), (\na \uu)_A(\e_1),\UU_A$ }

Since we can bound $\hat \e_1$ piecewisely by $C_{\hat \e}(x)E_4$ 
using  \eqref{eq:W2_res0}, \eqref{eq:W2_res1}, \eqref{eq:W2_bi} and bound $\bar \e_1$ by $ C_{\bar \e}(x)$ with a small coefficient $ C_{\bar \e}(x)$ of error size, we obtain the weighted norm $ \nlinf{  \e_1 \vp_{elli} }, [ \e_1 \psi_1]_{C_{x_i}^{1/2}}$ for $ \e_1 = \bar \e_1, \hat \e_1$,
where $\vp_{elli}$ is a specific weight for nonlocal error estimates defined in \eqref{wg:elli}. Following the methods in Section \ref{sec:u_comp_idea} and using nonlocal estimates similar to  \eqref{eq:u_linf_est1} in Section \ref{sec:linf_decay} but with  $\vp_{elli}$ 
defined in \eqref{wg:elli}  instead of $\vp_i$, we obtain functional inequalities of the form
\footnote{
In practice, for $\e_1 = \bar \e_1,  \hat \e_1$, since we can estimate piecewise bounds for $\e_1$ and 
since $\e_1 $ is much smaller in certain regions, we develop estimates similar to, but sharper than those in \eqref{eq:u_linf_err} 
by replacing the global norm of $\e_1$ by localized norms. 
We also establish improved weighted $C^{1/2}$ estimates using localized norms. See   Part II \cite[
Lemma {\lemmainvelPII} and Section {\secvellocest}]{ChenHou2023b}. 
}
\beq\label{eq:u_linf_err}
 |\rho_{ij} f_{ij}(\e_1)(x)| \leq \CCs_{\e,ij, 1}(x) || \e_1 \vp_{elli}||_{\inf} 
 + \CCs_{\e,ij, 2}(x) [ \e_1 \psi_1]_{C_x^{1/2}}
 + \CCs_{\e,ij, 3}(x) [ \e_1 \psi_1]_{C_y^{1/2}},
\eeq
where $f_{ij}(\e_1)$ denotes different velocity terms of $\uu_A(\e_1), (\na \uu)_A(\e_1)$. See \eqref{eq:u_linf_est1}. 
For $\e_1 =\bar \e_1, \hat \e_1$, we apply the above weighted estimates to bound $\uu_A(\e_1)$ and $(\na \uu)_A(\e_1)$. 

We now combine the three contributions to estimate $\UU_A$ in \eqref{eq:U:a}.
Under the bootstrap assumption $\sup_{s \leq t} E_4(s) < E_*$, since  $\uu_A(f)$ is linear in $f$, combining the estimates \eqref{eq:u_linf_est22} for $\uu_A(\om_1), (\na \uu)_A(\om_1)$ and 
\eqref{eq:u_linf_err} for $\uu_A(\e_1), (\na \uu)_A(\e_1)$ with $\e_1 = \bar \e_1, \hat \e_1$,
at time $t$, we bound $\UU_A = \uu_A(  \om_1 + \hat \e_1 + \bar \e_1)$ \eqref{eq:U:a} with the form 
\beq\label{eq:comb_vel1}
|U_A(x)| \leq \CCs(x) E_3(t) + \CCs_{\hat \e}(x) \sup\nolim_{s \leq t } E_4(s) + \CCs_{\bar \e}(x) 
\leq (\CCs(x) +  \CCs_{\hat \e}(x) ) E_* + \CCs_{\bar \e}(x) ,
\eeq
for computable coefficient functions $\CCs, \CCs_{\hat \e}, \CCs_{\bar \e}$. 
Note that the norm of $\bar \e$ defined in \eqref{eq:wb_vel_err} is a \emph{computable}  constant.  The terms $\CCs(x) E_3, \CCs_{\hat \e}(x) E_4, \CCs_{\bar \e}(x)$ bound the nonlocal terms of $\om_1, \hat \e_1, \bar \e_1$, respectively. 
Under the bootstrap assumption $E_4(s) < E_*$ \eqref{eq:Ec}, the above upper bound is a concrete computable bound, and we use \eqref{eq:comb_vel1} to combine the estimate of $\uu(\om), \uu(\bar \e), \uu(\hat \e)$.

\subsubsection{Modified decomposition}

To simplify presentation of later modified decomposition and nonlinear terms, we first introduce several bilinear operators.

\paragraph{\bf{Notations and operators}}

Given a vector-valued function $\uu \in \R^2$ and a matrix-valued function $M \in \R^{2 \times 2}$, we introduce the bilinear operator $B_{op,i}( (\uu, M), G )$ for $(\uu, M)$ and $G= (G_1, G_2, G_3)$ 
\beq\label{eq:Blin_gen}
\bal
\cB_{op, 1}( (\uu, M), G )
  & = - \uu \cdot \na G_1 + M_{11}(0) G_1,  \\
\quad \cB_{op, 2}( (\uu, M), G ) & = - \uu \cdot \na G_2 +2 M_{11}(0) G_2 - M_{11} G_2 - M_{21} G_3, \\
\cB_{op, 3}( (\uu, M), G ) & = - \uu \cdot \na G_3 + 2 M_{11}(0) G_3 - M_{12} G_2 - M_{22} G_3.
\eal
\eeq
The operator is bilinear in the first input $(\uu, M)$ and the second input $G$. 
The input $G$ relates to $(\om, \eta, \xi)$ and $M$ relates to velocity gradient $\na \uu$ and its variants. If $M=  \na \uu$, namely, $M_{11} = u_x, M_{12} = u_y, M_{21} = v_x, M_{22} = v_y$, 
\footnote{
This is the usual indexing for $\na \uu$, which differs from the stream-function-based 
indexing in \eqref{eq:u_linf_est1}.
}
we drop $M$ to simplify the notation 
\beq\label{eq:Blin}
\bal
\cB_{op, 1}(\uu, G )  & = - \uu \cdot \na G_1 + u_x(0) G_1,  \\
  \cB_{op, 2}(\uu, G) & = - \uu \cdot \na G_2 +2 u_x(0) G_2 - u_x G_2 - v_x G_3,  \\
\cB_{op, 3}(\uu, G) & = - \uu \cdot \na G_3 + 2 u_x(0) G_3 - u_y G_2 - v_y G_3.
\eal
\eeq

We recall the nonlinear terms $\cN_i$ from \eqref{eq:bous_non}
\[
\cN_1 = - \uu \cdot \na \om + u_x(0) \om , \ 
\cN_2 = - \uu \cdot \na \eta - u_x \eta - v_x \xi + 2 u_x(0) \eta, \ 
 \cN_3 = - \uu \cdot \na \xi - u_y \eta - v_y \xi + 2 u_x(0) \xi.
\]

 Using the above bilinear operators, we simplify the nonlinear terms as
\beq\label{eq:bous_non_recall}
\bal
 \cN_i = \cB_{op, i}(\uu, (\om, \eta, \xi)) .
\eal
\eeq

Next, we use the above decomposition and \eqref{eq:W2_err_est1},\eqref{eq:W2_err_est2} to rewrite the decomposition \eqref{eq:bous_decoup2}. Using the operators \eqref{eq:Blin_gen}, \eqref{eq:Blin}, we first rewrite  \eqref{eq:bous_decoup2}  as
\beq\label{eq:lin_old}
\bal
 \pa_t W_{1, i}
+ \bar c_l \xx \cdot \na & W_{1, i}
=    \one_{i =1} \eta_1
+ \mf{a}_i  \bar c_l W_{1,i}
+ \cB_{op, i}( (\uu_A(\om_1), (\na \uu)_A(\om_1)) , \bar W) + \cB_{op, i}( \bar \uu, W_1 ) \\
& +  (\cN_i - D_i^2 \cN_i(0) \cdot f_{\chi, i})
+ (\bar \cF_i - D_i^2 \bar \cF_i(0)  \cdot  f_{\chi, i} )  - \cR_{loc, 0,  i} \\
& + D_i^2 \cR_i^{\bar \e}(0) (f_{\chi, i} -\chi_{ i 2} ) 
   + ( \cR_i^{\bar \e} - D_i^2 \cR_i^{\bar \e}(0) f_{\chi , i} ) +
( \cR_i^{\hat \e} - D_i^2 \cR_i^{\hat \e}(0) \chi_{i2} ) .
\eal 
\tag{Orig}
\eeq
where $\mf{a} = (1/2, 1, 1)$ and we used 
$\bar c_{\om} = \bar u_x(0) + \tf12 \bar c_l$ by \eqref{eq:normal}.  The first row rewrites the linearized equations from \eqref{eq:lin_main}
and corresponds to $(\cL_i - \cK_{1i} - \cK_{2i}) W_{1,i}$ in \eqref{eq:W2_err_est2}.
 The second row uses
\[
\cNF_i(W) = ( D_i^2 \cN_i(0) + D_i^2 \bar \cF_i(0)) \cdot f_{\chi,i}, \quad D^2 = (\pa_{xy}, \pa_{xy}, \pa_{xx}),
\]
from the definitions in  \eqref{eq:appr_near0},\eqref{eq:van_pertb_sq}. In the third row, we substitute the formula for $\cR_j$ from  \eqref{eq:W2_err_est2}  into \eqref{eq:bous_decoup2}, taking into account the \emph{minus} sign $-\cR_j$.

\paragraph{\bf{Modified nonlinear terms}} First, we combine the terms $I - D_i^2 I(0) f_{\chi, i}$ for three different $I$ :

(a) $I_{a,i} =\cB_{op, i}( \uu(\bar \e), \wh W_2) $ from the residual operator 
$\cR_i^{\bar \e}$ \eqref{eq:W2_err_est1}, \eqref{eq:lin_old};

 (b) $I_{b,i}$ being the terms involving $\uu(\bar \e)$ in the linearized equations of \eqref{eq:lin_main}, \eqref{eq:bad} 
\[
\bal
& - \uu(\bar \e) \cdot \na \om_1 + c_{\om}(\bar \e) \om_1 
= \cB_{op, 1}( \uu(\bar \e), W_1) , \\
& - \uu(\bar \e) \cdot \na \eta_1 + 2 c_{\om}(\bar \e) \eta_1 - u_x(\bar \e) \eta_1 - v_x(\bar \e) \xi_1  = \cB_{op, 2}( \uu(\bar \e), W_1) ,  \\
& - \uu(\bar \e) \cdot \na \xi_1 + 2 c_{\om}(\bar \e) \xi_1 - u_y(\bar \e) \eta_1 - v_y(\bar \e) \xi_1  = \cB_{op, 3}( \uu(\bar \e), W_1) ,
\eal
\]
which has vanishing order $O(|x|^3)$ due to \eqref{eq:reg_W1} and thus $ D_i^2 B_{op, i}( \uu(\bar \e), W_1)(0) = 0$; 

(c)  nonlinearity  $I_{c,i}=\cN_i = \cB_{op, i}(\uu, (\om, \eta, \xi))$ in \eqref{eq:bous_non_recall} and near-origin correction in \eqref{eq:lin_old}.

Using bilinearity  of $\cB_{op}$, $W_1 + \hat W_2 = (\om, \eta, \xi)$, 
and $\UU = \uu(\om + \bar \e)$ by \eqref{eq:U:a}, we combine 
\bseq\label{eq:bous_nonM}
\beq
\bal
 I_i  \teq I_{a, i} + I_{b, i} + I_{c, i}   & = \cB_{op, i}( \uu(\bar \e), \hat W_2)
+ \cB_{op, i}( \uu(\bar \e), W_1) + \cB_{op, i}( \uu(\om), (\om, \eta, \xi)  )  \\ 
  & = \cB_{op, i}( \uu(\bar \e) + \uu(\om) ,(\om, \eta, \xi) )
= \cB_{op, i}(\UU, (\om, \eta, \xi) ) ,  \\
\eal
\eeq
and define the modified nonlinear term 
\beq
  \td \cN_i  \teq   I_i - D_i^2 I_i(0) f_{\chi, i}
  = \cB_{op, i}( \UU, W_1) 
  +  \cB_{op, i}( \UU, \hat W_2)  - U_x(0) D_i^2 (W_{1,i} + \wh W_{2,i}) (0) f_{\chi, i}.
\eeq
\eseq

Here, we have used $W_1 + \hat W_2 = (\om, \eta, \xi), D^2 =  (\pa_{xy}, \pa_{xy} ,\pa_x^2)$, 
the fact that $B_{op, i}$ \eqref{eq:Blin} is bilinear, and $\UU = U_x(0)(x, -y) + O(|x|^2), \om, \eta, \xi = O(|x|^2)$ to obtain
 \[
 D_i^2 \cB_{op, i}(\UU, (\om, \eta, \xi ))(0) = U_x(0)  V_i, \quad V = ( \pa_{xy}\om(0)
 , \pa_{xy} \eta(0),  \pa_{xx} \xi (0) ).
 \]
 We remark that the \emph{full solution} $\eta = \th_x, \xi = \th_y$ satisfies $\eta_y =\xi_x$.

\vs{0.1in}
\paragraph{\bf{Modified residual error}}

We decompose residual error $\bar \cF_i$ \eqref{eq:bous_err} into the essentially local part and nonlocal part. Recall the general bilinear operator \eqref{eq:Blin_gen} $\cB_{op}(\uu_A, (\na \uu)_A, G)$. Since
\[
 \bar \uu = \uu(\bar \om) = ( \bar \uu - \uu_A(\bar \e_1)) + \uu_A(\bar \e_1)
 = ( \bar \uu^N + \uu(\bar \e_2) + \appo{\uu}( \bar \e_1) ) + \uu_A(\bar \e_1),
\]
by  \eqref{eq:uerr_dec1}, \eqref{eq:uerr_dec12}, we modify $\bar \cF_i$ \eqref{eq:bous_err}, \eqref{eq:lin_old}
by replacing $\bar \uu$ by the \emph{essential local part} $ \bar \uu - \uu_A(\bar \e_1)$
\beq\label{eq:bous_errM}
\bal
 \bar \cF_{loc, i} & \teq {II}_i - D_i^2 II_i(0) f_{\chi,i} = 
 {II}_i - D_i^2 \bar \cF_i(0) f_{\chi, i}, \\
   II_i  & = \bar \cF_i - \cB_{op, i}( (\uu_A(\bar \e_1), (\na \uu)_A(\bar \e_1) ), \bar W ), \\
II_1 & = - (\bar c_l x + \bar \uu - \uu_A(\bar \e_1) ) \cdot \na \bar \om + \bar \th_x 
+ \bar c_{\om} \bar \om, \\
II_2 & = - (\bar c_l x + \bar \uu - \uu_A(\bar \e_1) ) \cdot \na \bar \th_x + 2 \bar c_{\om} \bar \th_x - (  \bar u_x - u_{x, A}(\bar \e_1) )  \bar \th_x - ( \bar v_x - v_{x, A}(\bar \e_1) ) \bar \th_y, \\
II_3 & = - (\bar c_l x + \bar \uu - \uu_A(\bar \e_1) ) \cdot \na \bar \th_y + 2 \bar c_{\om} \bar \th_y -   (  \bar u_y - u_{y, A}(\bar \e_1) )  \bar \th_x - ( \bar v_y - v_{y, A}(\bar \e_1) ) \bar \th_y . \\
 \eal
\eeq
From our construction in Section \ref{sec:appr_vel} 
and \eqref{eq:velA_reg} and $\bar \e_1$ in \eqref{eq:uerr_dec1}, we have $\uu_A(\bar \e_1) = O(|x|^3), (\na \uu)_A(\bar \e_1) = O(|x|^2)$. Thus by definition of $\cB_{op, i}$ in \eqref{eq:Blin_gen}, 
 $D_i^2 II_i(0) = D_i^2 \bar \cF_i(0)$. 
From \eqref{eq:uerr_dec1}, \eqref{eq:uerr_dec12} and the discussion below, the above error $\bar \cF_{loc, i}$ essentially depends on the numerical construction \emph{locally}. 

Similarly, 
using $\cR_i^{\hat \e} = \cB_{op, i}( \uu(\hat \e), \bar W )$ \eqref{eq:W2_err_est1}, 
we modify the remaining part of $\cR_{\cdot}$ terms 
in \eqref{eq:lin_old} as (we have combined $ \cR_i^{\bar \e}$ into $\td \cN_i$ in \eqref{eq:bous_nonM})
\begin{align}
\cR_{loc, i} & \teq \cR_{loc, 0, i} - D_i^2 \cB_{op, i}( \uu(\bar \e), \hat W_2)(0)
\cdot (  f_{\chi, i}  -  \chi_{i2}  )  \label{eq:W2_errM} \\
& \quad - \B( \cB_{op, i}( \uu(\hat \e), \bar W ) -  D_i^2 \cB_{op, i}( \uu(\hat \e), \bar W )(0) \cdot \chi_{i 2}
- \cB_{op,i}( \uu_A(\hat \e_1), (\na \uu)_A(\hat \e_1), \bar W) \B). \notag 
\end{align}
The term $\cB_{op, i}( \uu(\hat \e), \bar W )  -  \cB_{op,i}( \uu_A(\hat \e_1), (\na \uu)_A(\hat \e_1), \bar W)$ is essentially local.

\vs{0.1in}
\paragraph{\bf{Modified bad terms}}

Recall $\UU_A$ from \eqref{eq:U}. We combine the terms $\uu_A(\hat \e_1 )$ in \eqref{eq:W2_errM}, $\uu_A(\bar \e_1)$ in \eqref{eq:bous_errM} and the bad terms \eqref{eq:bad} with $\bar \uu$ replaced by $\bar \uu^N$ to define the modified bad terms:
\beq\label{eq:bad_lin}
\bal
 B_{\mf{modi}, 1}(x) & 
 = \eta_1  -  \UU_{A} \cdot \na \bar \om ,  \\
B_{\mf{modi}, 2}(x) &  =  -  \bar v^N_x \xi_1 -  \UU_A \cdot \na \bar \th_x -   \UU_{x, A} \cdot \na \bar \th , \\
 B_{\mf{modi}, 3}(x)  & =  -  \bar u^N_y \eta_1  - \UU_A \cdot \na \bar \th_y -  \UU_{y, A} \cdot \na \bar \th.
\eal
\eeq

Note that we added $\cB_{op,i}( \uu(\bar \e), W_1)$ into $\td \cN_i$ \eqref{eq:bous_nonM},
subtracted $\cB_{op, i}( \uu_A(\bar \e_1), (\na \uu)_A(\bar \e_1), \bar W)$ in $\bar \cF_i$ \eqref{eq:bous_errM},
and added $\cB_{op, i}( \uu_A(\hat \e_1), (\na \uu)_A(\hat \e_1), \bar W)$ 
into $\cR_{loc, i}$ \eqref{eq:W2_errM} (hence subtracting it from $-\cR_{loc, i}$). Adding the 
corresponding terms with opposite signs 
\[
- \cB_{op,i}( \uu(\bar \e), W_1) + \cB_{op, i}( \uu_A(\bar \e_1), (\na \uu)_A(\bar \e_1), \bar W)    +  \cB_{op, i}( \uu_A(\hat \e_1), (\na \uu)_A(\hat \e_1), \bar W),
\]
to the RHS of the linear part in \eqref{eq:lin_old} (first row), which generates the linear part in \eqref{eq:lin_U} (corresponds to \eqref{eq:lin_U} without $\td \cN, \bar \cF_{loc}, \cR_{loc}$ terms). By combining  $\td \cN_i$ \eqref{eq:bous_nonM}, $\bar \cF_{loc, i}$ \eqref{eq:bous_errM}, $\cR_{loc,i}$ \eqref{eq:W2_errM}, we rewrite \eqref{eq:lin_old} equivalently as
\beq\label{eq:lin_U}
\bal
\pa_t \om_1 + (\bar c_l x + \bar \uu^N  ) \cdot \na \om_1 &= \bar c_{\om}^N \om_1 + B_{ \mf{modi}, 1} + \td \cN_1 + \bar \cF_{loc,1} - \cR_{loc, 1} ,  \\ 
\pa_t \eta_1 + (\bar c_l x  + \bar \uu^N  ) \cdot \na \eta_1  & = (2 \bar c_{\om}^N - \bar u_x^N) \eta_1 + B_{\mf{modi}, 2} + \td \cN_2 +  \bar \cF_{loc,2} - \cR_{loc, 2} 
 , \\
\pa_t \xi_1 + (\bar c_l x + \bar \uu^N  )  \cdot \na \xi_1 & = 
(2 \bar c_{\om}^N + \bar u^N_x ) \xi_1 + B_{ \mf{modi}, 3} + \td \cN_3+  \bar \cF_{loc, 3} - \cR_{loc, 3}.
\eal
\eeq

 The advantage of this modified equation is the following. The nonlocal elliptic errors
have been absorbed into $\UU$ \eqref{eq:U}, so that the error $\bar \cF_{loc, i}, \cR_{loc, i}$ in \eqref{eq:lin_U} 
are \emph{essentially local} and depend on the numerical basis \emph{locally}, and $\td \cN_i, \bar \cF_{loc, i}, \cR_{loc, i} $ have cubic vanishing  $O(|x|^3)$ near $0$. 
Moreover, the coefficients of the linearized operators, e.g. 
$\bar \uu^N, \bar c_{\om}^N$, are also \emph{locally}.
This allows us to estimate these functions 
rigorously using numerical analysis with tracking nonlocal errors.

The linear energy estimates 
for the main linear equation \eqref{eq:lin_main} in Sections \ref{sec:linf_decay}-\ref{sec:rank1} can be re-derived directly for the linear part in \eqref{eq:lin_U}  with coefficients $(\bar \uu, \na \bar \uu, \bar c_{\om})$
replaced by $(\bar \uu^N, \na \bar \uu^N, \bar c_{\om}^N)$, 
and with the weighted $L^{\infty}(\vp_i), L^{\infty}(\vp_{g, i})$ estimates 
and $C^{1/2}$ estimates of $\uu_A, (\na \uu)_A$, replaced by 
those of $\UU_A, (\na \UU)_A$. 
Then we obtain the  linear stability conditions 
\[
 L^{\infty}(\vp_i), i\leq 3 : \eqref{linf:WG1_lin}, 
 \quad  L^{\infty}(\vp_4): \eqref{linf:WG4_lin}, 
\quad  C^{1/2} : \eqref{hol:lin} , \quad 
L^{\infty}(\vp_{g, i}):  \eqref{linf:WG1g_lin}, 
 \]
corresponding to those established in Section \ref{sec:linf_decay}, \, \ref{sec:hol_Q1}, \, \ref{sec:EE_hol}, 
and \ref{sec:linf_grow}, respectively. We also modify the damping coefficients $d_{i,L}$ \eqref{eq:dp} to $d_{i,L}^{\mf{num}}$ \eqref{eq:dp_num}.

\begin{remark}[\bf Estimate nonlocal error]
The errors $\hat \e_1,  \bar \e_1$ \eqref{eq:uerr_dec1} are much smaller than the bootstrap threshold $E_*$ chosen in  \eqref{eq:Ec} for $|x|$ not very large, where the damping factor is relatively small. 
In this region, the mesh is much finer, making the numerical Poisson error $\bar \e, \hat \e$ in
\eqref{eq:wb_vel_err},\eqref{eq:lin_evo_main1}, and hence error $\hat \e_1, \bar \e_1$ \eqref{eq:uerr_dec1}, very small. For readability, the reader may temporarily regard  $\hat \e_1, \bar \e_1$ as $0$ and $\bar \uu^N = \bar \uu= \uu(\bar \om), \UU = \uu(\om), \UU_A = \uu_A(\om_1) $. The actual proof retains all error terms.
We combine the estimates of error \eqref{eq:uerr_dec1} and the perturbation to simplify the nonlocal error estimate significantly.  We do not apply standard estimates for the operator $\na^{\perp}(-\D)^{-1} \e$ to bound $\uu(\e)$ from bounds of $\e$. Such  error estimates are 
not sharp enough to close the estimate, and we need weighted estimates for the error. 
\end{remark}

\begin{figure}[h]
   \centering
      \includegraphics[width = \textwidth  ]{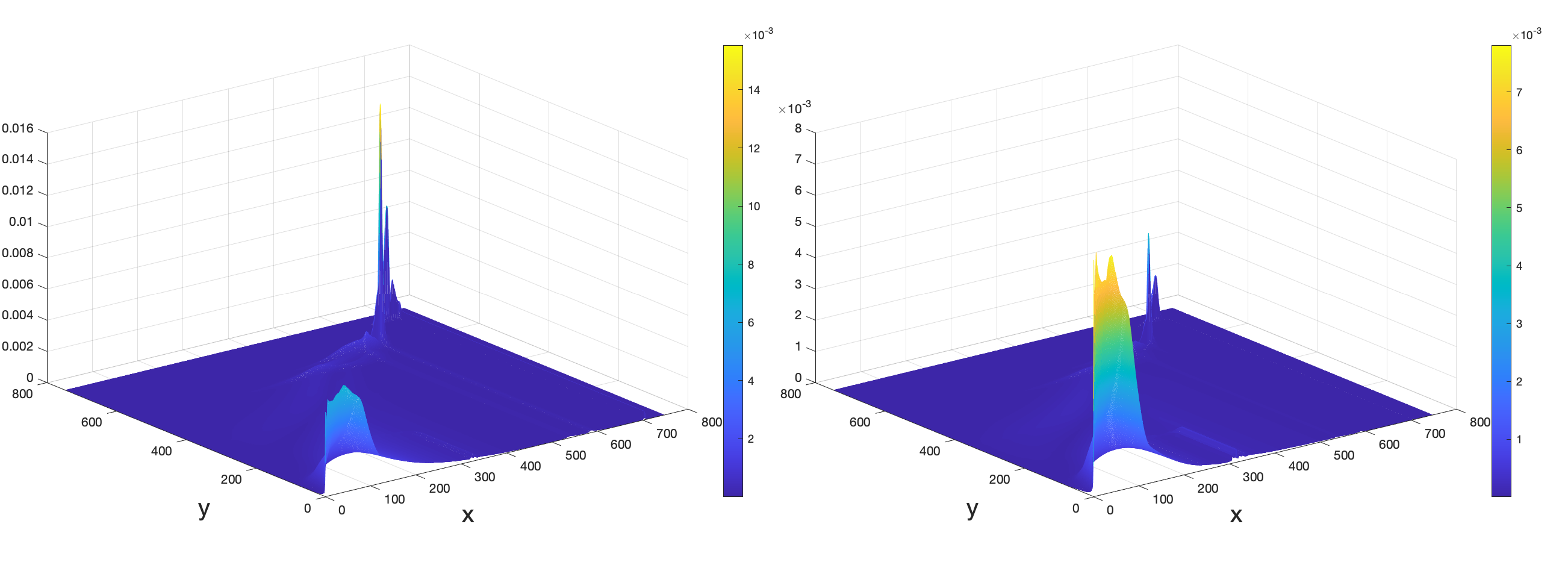}
      \captionsetup{width=0.95\textwidth} 
      \caption{Rigorous piecewise weighted bound of the local residual operator under adaptive mesh.  Left: bound of $\cR_{loc, 1} \vp_1$ in $\om$-equation. Right: bound of $ \cR_{loc, 2} \vp_2$ in $\eta$-equation.}
            \label{fig:err_evo}          
 \end{figure}

Using \eqref{eq:W2_res1} and following the estimates in Part II \cite[Section {\secerridea} \& Appendix {\secresid}]{ChenHou2023b}, 
we control the local residual operator $\cR_{loc, i}$ \eqref{eq:W2_errM}. In Figure \ref{fig:err_evo}, we plot the rigorous piecewise bound $C_{kl}$ for 
$| \vp_i \cR_{loc,i}| \leq C_{kl} E_4$ in $[y_k, y_{k+1}] \times [y_l, y_{l+1}]]$ with mesh $y_k$ in Section \ref{sec:ASS_rep}. 

In the near-field $|x| \leq  10^4 < y_{550} $, we have $ |\cR_{loc, 1} \vp_1| \leq 0.008 E_4, |\cR_{loc, 2} \vp_2| \leq 0.009 E_4$. We have $|\cR_{loc, 3}| \vp_3 \leq 0.002 E_4$ for $x$ in the 
range of mesh $[0, L]^2, L \approx 10^{13}$.  The near-field region with a large weighted error is about $[0, y_{100}]^2$ with $y_{100} \approx 0.4863$. In such a region, the error is much smaller than the remaining damping part in the weighted $L^{\inf}(\vp_i)$ estimate. See Figure \ref{fig:Linf1}. 
In the far-field region $|x| \geq 10^4$, where the weighted error is relatively large ($| \cR_{loc,1} \vp_1| \leq 0.016 E_4$, $|\cR_{loc, 2} \vp_2| \leq 0.009 E_4$), we have a large damping coefficient.
\footnote{
  To assess the relative size of $\cR$, under the nonlinear stability conditions \eqref{eq:PDE_nondiag}, the bounds $0.008$, $0.009$,$0.016$ (appearing as $a_{ij}$ in $a_{ij} E_*$) should be compared to 
  the stability factor $\lambda_*$ in \eqref{eq:EE_target1} (appearing as $a_{ii}$ in $a_{ii} E_*$). 
  The far-field error can be further reduced using a
finer mesh and a larger computational domain.  
  }
The estimate from the nonlinear modes $R_{ \mf{nl}}$ \eqref{eq:bous_err_op} is very small compared to the above bounds, and we have bounded it under the bootstrap assumption $E_4 < E_* = 5 \cdot 10^{-6}$, which will be discussed in Section \ref{sec:non}.

Using the piecewise weighted $L^{\inf}$ and $C^1$ bounds of $\cR_{loc, i}$ and the method in Part II \cite[Appendix E]{ChenHou2023b}, we derive the piecewise $C_{x_i}^{1/2}$ estimate of $\cR_{loc, i} \psi_i$. We combine this H\"older estimate of $\cR_{loc, i} \psi_i$ with the energy estimate in Section \ref{sec:EE_hol}. The resulting bound is very small compared to the least damping coefficients (near $x = (0.5, 0)$; see Figure \ref{fig:hol1}) since the estimates of $\cR_{loc,i}$ are much smaller near $x = (0.5, 0)$, and we have a small factor $ \tau_1^{-1}$
\footnote{
Consider Lemma \ref{lem:PDE_stab} with energy $E_2$ \eqref{energy2}. 
In the stability condition \eqref{eq:PDE_diag} for the H\"older norm, $\tau_1^{-1}$ corresponds to $\mu_i$ and leads to smaller ``bad terms'': $\sum_{j \neq i} |a_{ij}| \mu_i \mu_j^{-1}$. 
This property extends to Lemma \ref{lem:PDE_nonstab}.
}
for $ \cR_{loc, i} \psi_i$ in the H\"older energy estimate. This factor comes from the weight $\tau_1^{-1}$ of the H\"older norms in $E_2$ \eqref{energy2}

\begin{remark}[\bf Plots of estimates]\label{rem:plot}
To establish the nonlinear stability conditions \eqref{eq:PDE_nondiag}
associated with the energy $E_4$ \eqref{energy4}, we \emph{do not} rely on the plots. Instead, we
verify the required inequalities rigorously in Appendix~\ref{sec:ineq}. The
plots are included only to illustrate the size and spatial distribution of the
different terms in the nonlinear stability estimates.
\end{remark}

\subsection{Nonlinear estimates}\label{sec:non}

In Sections \ref{sec:linf_decay}--\ref{sec:rank1}, 
we established weighted  $L^{\inf}$ and $C^{1/2}$ estimates for $W_1$, ODE estimates for the
scalar functionals in the energy based on the main linear equation \eqref{eq:lin_main}. 
In Section \ref{sec:W2}, we estimate the local terms associated with 
$\hat W_2$, including $\cR_{loc, i}, \hat W_2$ using the estimates of scalar functionals. 
In Section \ref{sec:comb_vel_err}, we combine the nonlocal error and the nonlocal perturbation, perform the modified decomposition \eqref{eq:lin_U}, and directly generalize the linear estimates for \eqref{eq:lin_main} in Sections \ref{sec:linf_decay}--\ref{sec:rank1} to the linear part in \eqref{eq:lin_U}.

 Using the energy $E_4$ \eqref{energy4}, we control the $L^{\inf}$ norms of $W_{1,i}, \na \UU$, $\hat W_2$ (Sections \ref{sec:rank1}, \ref{sec:W2}). 
 In this section, we further control the remaining nonlinear terms $\td \cN_i$ in \eqref{eq:lin_U} 
and establish the \emph{full} nonlinear estimates for the energy $E_4$ 
of $W_1, W = W_1+ \hat W_2$  similar to 
  \[
       \tf{d}{dt} E_4 \leq- \lam E_4 +  C E_4^2 + \e    ,
    \]
where $-\lam E_4$ with $\lam >0$ reflects linear stability, $C E_4^2$ bounds nonlinear terms with $C = C(\bar \omega, \bar \theta, \psi, \varphi) $, and $\e$ bounds the weighted norm of profile residual. A sufficient condition for nonlinear stability  $E_4(t) < E^*$ for all $t \geq 0$ is $\e < \e^* = \lambda^2 / (4C)$,  which quantifies the required profile accuracy. Condition \eqref{eq:PDE_nondiag} in the stability Lemma \ref{lem:PDE_nonstab} provides similar pointwise constraints on error $\bar \e(x)$. Here, a key difference from the linear stability analysis (Sections \ref{sec:linf_decay}–\ref{sec:rank1}) is the small parameter $\e$. Together with the linear damping $-\lam E_4$, it allows us to tolerate large constants $C$ in the nonlinear estimates.  Thus, the nonlinear estimates are simpler and treated perturbatively.

We perform weighted estimates on $W_{1, i } \rho_i$ using \eqref{eq:lin_U}, which modifies decomposition  \eqref{eq:bous_decoup2}
\beq\label{eq:non_U}
\bal
 \pa_t( W_{1, i} \rho_i) + (\bar c_l x + \bar \uu^N + \UU) \cdot \na ( W_{1, i} \rho_i)
  & =  d^{\mf{num}}_{i, L}(\rho_i)  W_{1, i} \rho_i + B_{\mf{modi},i} \rho_i + \td \cN_i(\rho_i) + (\bar \cF_{loc, i} - \cR_{loc, i}) \rho_i . \\ 
\eal
\eeq

We put the perturbative velocity 
$\UU$ in the transport term so that the nonlinearity $\td \cN(\rho_i)$ 
does not contain derivatives $\na W_1$, and we have an \emph{extra nonlinear term} $\cT_{d, N}$ from the weight in  
\eqref{eq:tran_extra}. 
Here, we modify the linear damping terms \eqref{eq:dp} below and $\mf{num}$ is short for numerics
\beq\label{eq:dp_num}
\bal
 \cT_d^{\mf{num}}(\rho) & =  \rho^{-1} \big(  (\bar c_l x + \bar \uu^{N}  ) \cdot \na \rho \big) , &  \quad  d^{\mf{num}}_{1, L}(\rho)  &  = \cT^{\mf{num}}_d( \rho ) + \bar c^N_{\om} , \\
    d^{\mf{num}}_{2, L}(\rho) & = \cT^{\mf{num}}_d( \rho) + 2 \bar c^N_{\om} - \bar u^N_x,  &  \quad
 d^{\mf{num}}_{3, L }(\rho) &  = \cT^{\mf{num}}_d( \rho) + 2 \bar c^N_{\om} + \bar u^N_x .
\eal
\eeq
Here, the \emph{weighted} nonlinear term $\td \cN_i(\rho)$ is defined as follows and satisfies
\bseq\label{eq:tran_extra}
\beq
 \td \cN_i( \rho ) \teq ( \td \cN_i  + \UU \cdot \na W_{1, i} ) \cdot \rho + \cT_{d, N}(\rho) \cdot W_{1, i} \rho,
 \quad \td \cN_i( 1 ) \equiv \td \cN_i + \UU \cdot \na W_{1, i}, 
\eeq
where the \emph{unweighted} nonlinearity $\td \cN_i$ is defined in  \eqref{eq:bous_nonM},
The extra term arises from the weight: 
\beq
\UU \cdot \na ( W_{1, i} \rho ) =
\rho( \UU \cdot \na  W_{1, i})  + \rho^{-1} (\UU \cdot \na \rho ) \cdot (
W_{1, i} \rho ) , \quad \cT_{d, N}(\rho) \teq \rho^{-1} (\UU \cdot \na \rho ). 
\eeq
\eseq

\paragraph{\bf Decomposition of  $ \td \cN_i( \rho )$}

We decompose $\td \cN_i(\rho_i)$ into two types of terms to be estimated
separately: those involving $\UU$ and $W_1$, and those involving $\UU$ and
$\hat W_2$.

Using \eqref{eq:tran_extra} and $\UU=\td \UU+U_x(0)(x,-y)$ from \eqref{eq:U}, we first define the former by
\beq\label{def:N_W1}
\bal
N_{W_1, 1} & \teq U_x(0) \om_1 ,\\
N_{W_1, 2} & \teq  2 U_x(0) \eta_1 - U_x \eta_1 - V_x \xi_1
  =  U_x(0) \eta_1 - \td U_x \eta_1 - \td V_x \xi_1, \\
N_{W_1, 3} & \teq  2 U_x(0) \xi_1 - U_y \eta_1 - V_y \xi_1 
 = 3 U_x(0) \xi_1 - \td U_y \eta_1 - \td V_y \xi_1 ,
\eal
\eeq

Using $\td \cN_i$ \eqref{eq:bous_nonM}, $\cB_{op, i}(\UU, W_1)$  \eqref{eq:Blin}, 
identities in \eqref{def:N_W1}, and $W_1 + \hat W_2 = (\om,\eta,\xi)$, we expand
\beq\label{eq:non_dec1}
\bal 
\td \cN_1(\rho_1) = &  (\UU \cdot \na \rho_1) \cdot   \om_1  + U_x(0) \om_1 \rho_1  \qquad \qquad \qquad \quad  + \cB_{op, 1} ( \UU, \wh W_2) \rho_1 - U_x(0) \om_{xy}(0) f_{\chi, 1} \rho_1
, 
 \\
\td \cN_2(\rho_2) = &  ( \UU \cdot \na \rho_2) \cdot \eta_1  + ( U_x(0) \eta_1 
- \td U_x \eta_1 -\td  V_x \xi_1 ) \rho_2 +  \cB_{op,2}( \UU, \wh W_2) \rho_2
- U_x(0) \eta_{xy}(0) f_{\chi, 2} \rho_2
, \\
 \td \cN_3(\rho_3) = &  (\UU \cdot \na \rho_3) \cdot \xi_1  + ( 3 U_x(0) \xi_1 
- \td U_y \eta_1 - \td V_y \xi_1 ) \rho_3 + \cB_{op,3}( \UU,\wh W_2) \rho_3 
- U_x(0) \xi_{xx}(0) f_{\chi,3} \rho_3 .
\eal
\eeq

Recall from \eqref{eq:reg_W1} that $\om = \om_1 + \hat \om_2$ and $\om_1 = O(|x|^3)$ near $0$. 
We have $\om_{xy}(0) = \pa_{xy} \hat \om_{2}(0)$. Using \eqref{eq:Blin}, \eqref{eq:U}, 
we define the nonlinearity involving $\UU$ and $\hat W_2$ by combining the last two terms in $\td \cN_i(\rho_i)$ and isolating the $U_x(0)$-terms 
\bseq\label{eq:W2_M2}
\beq
\bal
  N_{\hat W_2 ,i} & \teq  \cB_{op, i}(  \UU, \hat W_2)  - U_x(0) D_i^2 \hat W_{2, i}(0) f_{\chi, i} \\
  & =  \cB_{op, i}( \td \UU, \hat W_2)
  + \cB_{op, i}( U_x(0)(x, -y), \hat W_2) -  U_x(0) D_i^2 \hat W_{2, i}(0) f_{\chi, i}  \\
  &  =   \cB_{op, i}( \td \UU, \hat W_2) + U_x(0) \hat W_{2,  i, M} .
\eal
\eeq
where $D^2 = (\pa_{xy}, \pa_{xy}, \pa_x^2)$, and $\hat W_{2,  \cdot, M }$ is the modification of $\hat W_2$ and is defined as 
\beq
\bal
& \hat W_{2,  \cdot, M } =(\hat \om_{2, M}, \hat \eta_{2, M}, \hat \xi_{2, M} ),  & \quad 
 \hat \om_{2, M} &  \teq \hat \om_2 - x \pa_x \hat \om_2 + y \pa_y \hat \om_2 - \hat \om_{2, xy}(0) f_{\chi, 1} , \\
 &  \hat \eta_{2, M}   \teq \hat \eta_2 - x \pa_x \hat \eta_2 + y \pa_y  \hat \eta_2 -  \hat \eta_{2, xy}(0) f_{\chi, 2} , \  
   & \hat \xi_{2, M} &  \teq 3 \hat \xi_2- x \pa_x \hat \xi_2 + y \pa_y \hat \xi_2 
- \hat \xi_{2, xx}(0) f_{\chi, 3}.
\eal
\eeq
\eseq
Note that in $\cB_{op, i}(\td \UU, \hat W_2)$ in \eqref{eq:Blin}, 
using \eqref{eq:u_tilde}, we have $\td U_x(0) =  \pa_x U(0) - U_x(0)=0$. Using the above derivations, we rewrite the weighted nonlinearity in \eqref{eq:tran_extra} as 
\beq\label{def:non_wg}
\td \cN_i(\rho) = \cT_{d, N}(\rho) \cdot W_{1,i} \rho 
+ N_{W_1, i} \cdot \rho + N_{\hat W_2, i} \cdot \rho.
\eeq

Below, we discuss the estimates of these three types of nonlinear terms.

\subsubsection{Main nonlinear term}\label{sec:non_mainT}

Among the estimates of nonlinear terms in \eqref{eq:non_dec1},\eqref{eq:W2_M2},  that of
\[
   U_x(0)\hat W_{2,i,M},
\]
involves the largest constants. Recall that the functional estimates
 \beq\label{eq:cw_poor1}
|c_{\om}(\om) |= |u_x(\om)(0)| \leq \mu_6 E_4, \quad  | c_{\om}(\om_1)| = |u_x( \om_1)(0)| \leq \mu_5 E_4,
\eeq
by the energy $E_4$ \eqref{energy4} carry large constants. Since $U_x(0) = u_x(\om)(0)+ u_x(\bar \e)(0)$ by \eqref{eq:U:a} and $c_{\om}(\om) = u_x(\om)(0)$, 
under the bootstrap assumption $E_4 < E_*$, using \eqref{eq:cw_poor1}, we bound
\beq\label{eq:U_bd}
 |U_x(0)| < \mu_6 E_4 + |u_x(\bar \e)(0)| < \mu_6 E_* + |u_x(\bar \e)(0)|.
\eeq
From \eqref{eq:W2_bi} and Section \ref{sec:W2}, we have a large constant $c_1 = \mu_5 $ in our estimate of $\hat \om_2$. Thus, in the estimate of $U_x(0) \widehat W_{2, i, M}$ in \eqref{eq:W2_M2}, we incur a large constant about $\mu_5 \mu_6$, with $\mu_5 \mu_6 \approx 4700$.\footnote{
We also have a large constant $\mu_6^2$ in the estimate of nonlinear terms for $\cB_{op, i}(\td \UU, \widehat W_2)$ from \eqref{eq:W2_bi} since $\td \UU$ \eqref{eq:U} involves 
$ \td \uu(-\D \hat \phi^N)$ via $\UU_{\FM}$. 
 Since $\td \UU = O(|x|^3), \widehat W_2 = O(|x|^2)$, these nonlinear terms have a higher vanishing  $O(|x|^4)$ near $x=0$. 
Due to this higher vanishing  and the decay of $\hat \om_2$, our estimates of these nonlinear term are smaller than $U_x(0) \widehat W_{2, i, M} $ \eqref{eq:W2_M2}.
}
By contrast, constants in the estimates for $W_1$-term are much smaller.  
We have $ | \om_1 \vp_1 | \leq E_4$ with constant $1$ by \eqref{energy4}.
Similarly, the estimate of $\uu_A(\om_1)$ \eqref{eq:vel_rem} has constants  of order $1$.  Thus, $U_x(0) \widehat W_{2, i, M} $ is the main nonlinear term in \eqref{eq:non_dec1}, \eqref{eq:W2_M2} with large constants in its estimate.

After isolating the main term, we separate the nonlinear terms according to whether they contain $W_1$: $N_{W_1,i}$ \eqref{def:N_W1} and $\cT_{d, N}$ \eqref{eq:tran_extra}, or contain \(\widehat W_2\): $N_{\hat W_2,i}$  \eqref{eq:W2_M2}. We estimate these terms using simple  $L^{\infty}, C^{1/2}$ product estimates  \eqref{eq:hol_prod}-\eqref{eq:hol_sq}, the functional inequalities \eqref{eq:u_linf_est12}, \eqref{eq:u_linf_est22} and $C^{1/2}$ estimates for $\uu_A(\om_1), (\na \uu)_A(\om_1)$ and corresponding estimates for the nonlocal term $\UU_A$, which combines the perturbation and error.  Terms $N_{W_1,i}$ are controlled using 
 weighted $L^{\infty}$, $C^{1/2}$ estimates for $W_1$ by energy $E_4$ \eqref{energy4} and the nonlocal estimates for $\UU$.  Terms $N_{\hat W_2, i}$ 
are more regular and are controlled using the rigorous precomputed \(C^3\) bounds from Section~\ref{sec:W2} 
(see e.g. \eqref{eq:W2_res1} and the nonlocal estimates for $\UU$. Below, we focus on some typical terms whose estimates require careful distribution of the singular weights.

\subsubsection{$L^{\inf}$ estimates}\label{sec:non_linf}

Using \eqref{eq:U}, \eqref{eq:u_tilde}, we decompose $\td \UU$ 
and the nonlinearity in $N_{\hat W_2,i}$ \eqref{eq:W2_M2} 
\beq\label{eq:U_dec1}
\bal
 \td \UU  & = \UU_A + \UU_{\FM}, \\ 
  \cB_{op, i}(\td \UU, \wh W_2) & = \cB_{op, i}( ( \UU_A, (\na \UU)_A ) , \wh W_2)  +  \cB_{op, i}( (\UU_{\FM}, (\na \UU)_{\FM} ), \wh W_2)  \teq I_1 + I_2.
 \eal
\eeq
We estimate $C^3$ bounds of the $\widehat W_2$ terms in \eqref{eq:W2_M2}, \eqref{eq:non_dec1} using \eqref{eq:W2_res1}, \eqref{eq:W2_res0} and following   Part II \cite[Section 3.6]{ChenHou2023b}.  
Then we can estimate $\hat W_{2, i, M }$ \eqref{eq:W2_M2} and apply the same estimate for the nonlocal terms $(\uu_A, \bar W)$ \eqref{eq:lin_main},\eqref{eq:lin} in Section \ref{sec:linf_decay} to $I_1$.
From the discussion below \eqref{eq:U}, we can estimate piecewise $C^3$ bounds for $\UU_{\FM}$. 
Then we obtain the estimate for $I_2$.

Recall  $W_1 =(\om_1, \eta_1, \xi_1)$. 
By definition of $E_1 $ \eqref{energy1}, $E_3$ \eqref{energy3}, we obtain
\[
 \nlinf{ W_{1, i} \vp_i } \les E_1, \nlinf{ \om_1 \vp_4} \les E_1,
 \quad  \nlinf{ W_{1, i} \vp_{g, i} } \les E_3. 
 \]
Since we have both piecewise  weighted and \emph{unweighted} $L^{\inf}$ bound for $ \na \UU,  \na \td \UU$ (see discussions around \eqref{eq:U}), 
the weighted $L^{\infty}$ estimate of $N_{W_1, i}$ in \eqref{def:N_W1} follows from trivial product rule.

Next, we consider the nonlinearity $\cT_{d, N}(\rho) W_{1, i} \rho$ in $\td \cN_i(\rho) \eqref{def:non_wg}$, \eqref{eq:non_U}.
For 
\[
(c, W_{1, i}, \rho) = (1, W_{1, i}, \vp_i),  
\qquad 
 (\mu_{g, i}, W_{1, i}, \vp_{g,i}) , \qquad 
( \sqrt{2} \tau_1^{-1}, W_{1, 1}, \vp_4) ,
\]
associated with the energy $E_1$ \eqref{energy1}, $E_3$ \eqref{energy3},
$L^{\infty}(\vp_4)$ in $E_1$ \eqref{energy1}, 
respectively, we estimate 
\begin{align}\label{eq:linf_non_tran}
c | \cT_{d, N}(\rho) W_{ 1, i} \rho |
  &  \leq ( |\f{U}{x} \cdot \f{x \pa_x \rho}{\rho} | + |\f{V}{y} \cdot \f{y \pa_y \rho}{\rho} |   )    || c W_{1, i} \rho||_{\inf} 
 \leq ( |\f{U}{x} \cdot  \f{x \pa_x \rho}{\rho} | + |\f{V}{y} \f{y \pa_y \rho }{\rho }   )   E_4 \notag \\
 & \leq \cT_\mw{wg}(\UU, \rho)(x) E_4, 
 \end{align}
and estimate piecewise $L^{\inf}$ bounds for $U/x, V /y, \pa_{x_j} \rho_i / \rho_i$.  Here, we denote 
\beq\label{eq:linf_non_tran2}
\cT_\mw{wg}(\UU, \rho)(x) \teq |\f{U}{x} |  \cdot |\f{x \pa_x \rho}{\rho} | + |\f{V}{y}| \cdot |\f{  y \pa_y \rho}{\rho} | ,
\eeq
Here \emph{wg} is short for weight,  and we use it to emphasize that the second component is a weight. Operator $\cT_\mw{wg}(\UU,\rho)$ is a
bookkeeping device for estimating  $\cT_{d, N}$ \eqref{eq:tran_extra} due to the weight.

\begin{remark}[\bf Two estimates]
In nonlinear estimates, we optimize two estimates of $\uu_A(\om_1)$, $(\na \uu)_A(\om_1)$ using functional inequalities based on  $|| \om_1 \vp_1||_{\inf}, 
[ \om_1 \psi_1 ]_{C_{g_1}^{1/2}}$ 
in Section \ref{sec:linf_decay} and $|| \om_1 \vp_{g, 1}||_{\inf}$, $[ \om_1 \psi_1 ]_{C_{g_1}^{1/2}}$  
in Section \ref{sec:linf_grow} so that $\na \uu_A \in L^{\inf}$. Similarly, we apply two functional inequalities to $\uu_A(\e_1), \e_1 = \bar \e_1, \hat \e_1 $ in Section \ref{sec:comb_vel_err} using  $|| \e_1 \vp_{elli}||_{\inf}, || \e_1 \vp_{g,1}||_{\inf}$ so that $ \na \UU \in L^{\inf} $.
\end{remark}

\subsubsection{H\"older estimate of typical terms}\label{sec:non_hol}

Recall the difference operator $\d$ from \eqref{eq:hol_dif}.

\vs{0.1in}
\paragraph{\bf{Nonlinearity $N_{W_1,i}$ involving $\UU, W_{1, i}$}}

We illustrate the $C^{1/2}$ estimate using the typical term $\td U_x \eta_1$ in $N_{W_1,2}$
\eqref{def:N_W1} with weight $\psi_2$.  Using the $C^{1/2}$ estimate of $\psi_1 \UU_A$ by the energy, $C^{1/2}$ estimate of $\psi_1^{-1}, \UU_{\FM}$ (see below \eqref{eq:U}), and \eqref{eq:hol_sq}, we have $C^{1/2}$ estimate $ \td \UU$. We then  estimate $\d( \td U_x \cdot \eta_1 \psi_2 ) $ using \eqref{eq:hol_loc}. For $x, z$ in the far-field, we need another decomposition and estimate since $\psi_1^{-1}$ \eqref{wg:hol} in this estimate is not bounded. 
We decompose $\td U_x = U_{x, \FM} + U_{x, A}$  by \eqref{eq:U}, and can still estimate $\d(  U_{x,\FM} \cdot \eta_1 \psi_2 ) $ using \eqref{eq:hol_loc}. 
For $ U_{x, A} \eta_1 \psi_2$, we use  \eqref{eq:hol_dif}, \eqref{eq:hol_prod} to get
\[
 \d( U_{x, A}  \eta_1 \psi_2)
  = \d( U_{x, A} \psi_1 \cdot \eta_1 \f{ \psi_2}{\psi_1}  )
 = \d(U_{x, A} \psi_1 ) \f{\eta_1 \psi_2}{ \psi_1}(x)
 + ( U_{x, A} \psi_1)(z) \d(\eta_1  \f{ \psi_2}{\psi_1} ) = I_{\mf{h}, 1} + I_{\mf{h},2}.
\] 
To bound $I_{\mf{h}, 1}$, using the energy $E_4$ \eqref{energy4}, we can bound $C^{1/2}$ of $U_{x, A} \psi_1$ and $\f{\eta_1 \psi_2}{ \psi_1}$ in $L^{\inf}$
\beq\label{eq:wg_linf_grow_hol}
 |  \f{ \eta_1 \psi_2 }{  \psi_1 }  | \leq  || \eta_1 \vp_{g, 2}||_{\inf} \f{ \psi_2 }{ \psi_1 \vp_{g, 2}  }. 
\eeq
To ensure $\psi_2 / ( \psi_1  \vp_{g, 2 }) \in L^{\inf}$, comparing far-field behaviors, $\psi_2 \sim |x|^{1/6}, \psi_1 \sim |x|^{-1/6}$, we need $\vp_{g, 2} \gtr |x|^{1/3}$. 
To achieve this condition, we choose $\al_{g, n} = \f{1}{3} + 10^{-8}$ in \eqref{eq:wg_linf_grow}. For $I_{\mf{h}, 2}$, we use
\[
|I_{\mf{h},2} | \leq |(U_{x, A} \psi_1)(z)|  \cdot ( \psi_1^{-1}(z) |\d( \eta_1 \psi_2)|
+ |\d( \psi_1^{-1})| \cdot |\eta_1 \psi_2|(x) ).
\]
Since $U_{x, A} \in L^{\inf}$, we can bound both terms using the energy $E_4$.

\vs{0.05in}
\paragraph{\bf Nonlinearity $\cT_{d,N}(\psi_i) \cdot W_{1,i}\psi_i $}
To estimate
$\cT_{d,N}(\psi_i)(W_{1,i}\psi_i)$ in \eqref{def:non_wg}, \eqref{eq:tran_extra}, the energy already 
provides the required $C^{1/2}$ control of $W_{1,i}\psi_i$  and we apply \eqref{eq:hol_loc}. We only need to control $\cT_{d, N}(\psi_i)$. Using \eqref{eq:U}, we perform the decomposition 
\beq\label{eq:non_tran_hol}
\cT_{d, N}(\psi_i) = \f{\UU \cdot \na \psi_i }{ \psi_i}
= U_x(0) \f{( x \pa_x \psi_i - y \pa_y \psi_i)}{ \psi_i }
 + \f{ \UU_A  \cdot \na \psi_i}{\psi_i}
 + \f{ \UU_{\FM}  \cdot \na \psi_i}{\psi_i}
\teq \cT_{ c_{\om}}( \psi_i ) + \cT_{uA}(\psi_i ) + \cT_{uR}(\psi_i).
\eeq

Recall $\psi_u$ from \eqref{wg:elli}. For $\cT_{uA}$, since we have piecewise $C^{1/2}$ estimates of $\UU_A \psi_u$ (generalizing \eqref{eq:holest:u}), $C^3$ estimate of 
$\UU_{\FM}$ with  $ \UU_{\FM}= O(|x|^3)$, we decompose it as follows 
\beq\label{eq:non_dec3}
\cT_{uA} = ( U_A  \psi_u ) \cdot \f{\pa_x \psi_i}{ \psi_i \psi_u}
+ ( V_A  \psi_u ) \cdot \f{\pa_y \psi_i}{ \psi_i \psi_u}, 
\quad \cT_{u R} =  \f{ \UU_{\FM} }{|x|^2}
 \cdot |x|^2 \f{\na \psi_i}{\psi_i} ,
\eeq
and then estimate each product $f_1 \cdot f_2 $ using \eqref{eq:hol_sq}, \eqref{eq:hol_sq2}. The explicit function $\cT_{c_{\om}}$ is not $C^{1/2}$ near $0$, but we can bound $\d(\cT_{c_{\om}}, x, z) g(x-z) |p|^{1/2}$  with 
the \emph{extra vanishing power} $|p|^{1/2}$ for $p = x, z$. See Remark \ref{rem:dp_near0} and the supplementary material I \cite[Section {\suppsecnondphol}]{ChenHou2023aSupp}.

\vs{0.1in}
\paragraph{\bf{Nonlinearity $N_{\hat W_2, i}$ involving $ \UU, \hat W_2$}}

Recall $N_{\hat W_2, i}$  from \eqref{eq:W2_M2}  and decomposition in \eqref{eq:U_dec1}
\beq\label{eq:N_W2:hol_dec}
U_x(0) \hat W_{2,  i, M},  \quad  \cB_{op, i}(\td \UU, \wh W_2)  = \cB_{op, i}( \UU_A, (\na \UU)_A, \wh W_2)  +  \cB_{op, i}( \UU_{\FM}, (\na \UU)_{\FM}, \wh W_2)  = I_1 + I_2.
\eeq
We estimate a typical term
$P_6 = \UU_A \cdot \na \wh W_{2, i} \psi_i $ from \eqref{eq:Blin},\eqref{eq:non_dec1} in $I_1$. Using \eqref{eq:hol_sq}-\eqref{eq:hol_sq2} gets
\[
|\d( P_6 )| \leq \d_{\sq}( U_A \psi_u, \pa_x \hat W_{2, i} \f{\psi_i}{\psi_u}, h )
+ \d_{\sq}( V_A \psi_u, \pa_y \hat W_{2, i} \f{\psi_i}{\psi_u}, h ) .
\]
Recall $\psi_u, \psi_2, \psi_3 \sim |x|^{-2.5}, \psi_1 \sim |x|^{-2}$ near $x=0$ from \eqref{wg:hol},\eqref{wg:elli}.
Near $0$, we have $\na \wh W_{2, i} \f{\psi_i}{ \psi_u} = O(|x|)$, and we can estimate its $C^1$ bound. We bound $\d_i(  U_A \psi_u)$ using the energy.

For $I_2$, since  $\wh W_{2, i} = O(|x|^2), \UU_{\FM} = O(|x|^3)$ 
admit piecewise $C^3$ bounds, from \eqref{eq:Blin_gen}, we get
\[
  I_2 = \cB_{op, i}( ( \UU_{\FM}, (\na \UU)_{\FM} ), \hat W_2) = O(|x|^4), 
  \quad \mbox{near} \  x=0 \, ,
\]
 and we can estimate the $C^1$ bound of $\cB_{op, i} \psi_i$ and then its $C^{1/2}$ bound using \eqref{eq:hol_sq}-\eqref{eq:hol_sq2}.  

The main terms 
$U_x(0) \hat W_{2,  i, M}$ in \eqref{eq:W2_M2} have a vanishing order $O(|x|^3)$, and we can estimate their $C^3$ bounds. Since $\psi_i \les |x|^{-5/2}$, we can estimate $C^{1/2}$ bound of $\psi_i \hat W_{2,  i, M}$. In the supplementary material I \cite[Section {\suppsecbdother}]{ChenHou2023aSupp} (attached to this paper), we discuss the piecewise $C^{1/2}$ estimates of $f(x) / |x|^{5/2}$ for $f \in C^3$ with $f = O(|x|^3)$ near $x=0$.

Other nonlinear terms are estimated similarly and directly using \eqref{eq:hol_prod}--\eqref{eq:hol_sq}. See supplementary material I \cite[Section {\suppsecnonest}]{ChenHou2023aSupp} 
for details.
For nonlinear stability, we can allow much larger constants in the estimates of terms other than the main terms \eqref{eq:W2_M2}, whose estimates is simple.

\vs{0.1in}
\paragraph{\bf{Other nonlinear estimates}}

We consider the nonlinear terms in the ODEs \eqref{eq:lin_main_cw}. 
We rewrite the nonlinear term $\cN_i$ using \eqref{eq:non_dec1}, \eqref{eq:W2_M2} with $\rho_i \equiv 1$. The term $U_x(0) \om_1, U_x(0) \eta_1$ in \eqref{eq:non_dec1} in the ODEs of $c_{\om}( f q )$
\eqref{eq:lin_main_cw} with cutoff $q = \chi_{\mf{ode}}$ \eqref{eq:cw_chi} and $q=1$ \eqref{eq:ode_cw_all} 
contribute to  $U_x(0) \la \om_1, f_* q  \ra , f = \om_1, \eta_1$, which can be bounded by the energy $E_4$ \eqref{eq:cw_est}, \eqref{energy4} directly
\beq\label{eq:cw_est_boot}
\bal
 &  \tf{4}{\pi} |U_x(0)| \cdot |\la \om_1, f_* q\ra| 
 = |U_x(0) c_{\om}(\om_1 q)|
\leq  \g_{1} |U_x(0) | E_4, \\
&  \tf{4}{\pi} 
 |U_x(0)| \cdot |\la \eta_1, f_* q \ra|
 = |U_x(0) c_{\om}(\eta_1 q) |
\leq \g_2 |U_x(0)|  E_4, 
\eal
\eeq
where $ (q, \g_1, \g_2) = (\chi_{\mf{ode}}, \mu_{51}, \mu_{52}), (1, \mu_5, \mu_{62})$. The estimates of other nonlinear terms in these ODEs follow Section \ref{sec:non_linf} and the argument in Section \ref{sec:rank1}, e.g. integration by parts.

For the energy estimates beyond our computational domain $[0, D]^2, D \geq 10^{15}$, we estimate the asymptotics of the profile \eqref{eq:ASS_decomp1} in Appendix {\secestapprfar} and the nonlocal terms in Section {\secholfar} in Part II \cite{ChenHou2023b}. Since the coefficients of the nonlocal terms decay, e.g. $\na \bar \om, \na \bar  \th$, the equations \eqref{eq:lin_main}, \eqref{eq:bous_decoup2} are essentially local in the far-field. We have much larger damping factors and can afford much larger constants in the estimate of nonlocal terms. 
Since far-field estimates are \emph{completely perturbative}, we refer them to supplementary material I \cite[Sections {\suppseclinffar}, {\suppsecholfar}]{ChenHou2023aSupp}.

\begin{figure}[h]
   \centering
      \includegraphics[width = \textwidth  ]{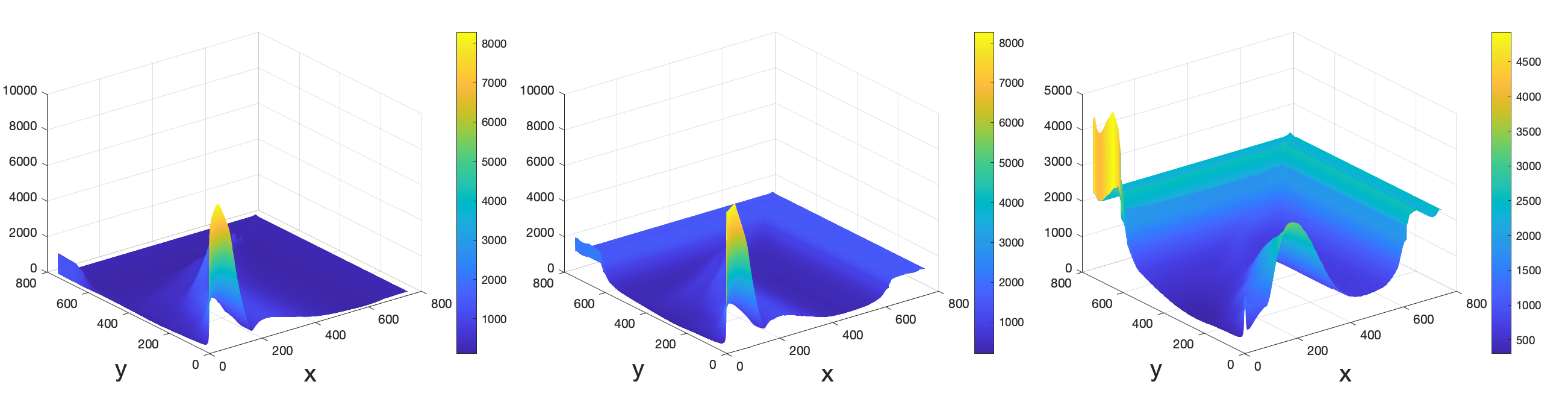}
       \caption{Weighted $L^{\inf}(\vp_i)$ estimate of the nonlinear terms.
       }
            \label{fig:N_linf_WG1}
 \end{figure}

In Figure \ref{fig:N_linf_WG1}, we plot the rigorous piecewise bounds $C_{i, kl}$ for the  
weighted nonlinear terms $\td \cN_i(\vp_i)$ defined in \eqref{eq:non_dec1} with $\rho_i = \vp_i$ under the bootstrap assumption $\sup_{s \leq t}E_4(s) < E_*$. 
On each mesh grid $[y_k, y_{k+1}] \times [y_l, y_{l+1}]$, whose union covers $[0, 10^{15}]^2$, we have 
$ | \td \cN_i(\rho_i)(t)| \leq C_{i, kl} E_*^2$, in the $\om, \eta, \xi$ equation, respectively. The largest terms for these three equations are bounded by $8300, 8300, 5000$, respectively. 
\footnote{
In Figure \ref{fig:N_linf_WG1}, for very small $x$ and very large $y$ (adaptive mesh $y_{2, n} > 10^6$ for $n\geq 600$), we observe a jump in the estimate, especially in the $\xi$-equation. This arises from the piecewise estimate of $\vp_3 / \vp_2$ (see \eqref{wg:linf_decay}) used in bounding the nonlinear term $|u_y \eta_1 \vp_3| \leq |u_y \f{\vp_3}{\vp_2}| || \eta_1 \vp_2||_{\inf}$. 
Both weights involve the $|x_1|^{-1/2}$ singularity along $x_1 = 0$. We can refine the estimates to get a smoother bound. Yet, since we have a large damping factor (bigger than $1$) in this far-field region, we can tolerate a constant up to $20$ times larger (e.g., $- E_* + 20 \cdot 5000 E_*^2 \leq -0.5 E_* $) and do not need the refinement.

}

\paragraph{\bf{Estimate of the residual error of the profile}}

For the residual error $ \bar \cF_{loc, i}$ in \eqref{eq:bous_errM} modified from \eqref{eq:bous_decoup2}, \eqref{eq:bous_err}, it is essentially local and its estimate follows standard numerical analysis. We estimate them following Part II \cite[Sections 3.6]{ChenHou2023b} with some details in Part II \cite[Appendix {\secresid}]{ChenHou2023b}. We have plotted rigorous piecewise bounds for the local part $\bar \cF_{loc, i} \vp_i, i= 1,2$ in Figure \ref{fig:prof_err}. Note that the weighted residual error away from the first few grids and in the bulk region is very small ($\leq 5 \cdot 10^{-8}$) relative to the bootstrap threshold $E_*$ \eqref{eq:Ec}.

\vs{0.05in}
\paragraph{\bf{Estimate in the region with small stability factor}}

In the linear weighted $L^{\inf}(\vp_i)$ estimate, we have a minimum stability factor about $0.04$. We have a small stability factor below $0.08$ only in the bulk region $D_{B} = [0, 1000]^2 \bsh [0, 1]^2$. Since it is away from $0$, the singular weight $\vp_i$ becomes much smaller and both the estimate of nonlinear terms and the residual error becomes much smaller in $D_B$. See Figures \ref{fig:prof_err}, \ref{fig:N_linf_WG1}. Similar discussion applies to the H\"older estimate.

\begin{remark}[\bf Local stability condition]

We establish nonlinear stability using the stability condition \eqref{eq:PDE_nondiag}. 
A key advantage of this condition  is that it depends on the estimate locally. 
Thus, we do not need to compare the minimum damping coefficients with the $L^{\inf}$ bound of the nonlinear terms and error terms, or verify
\[
 \min_{x, z} a_{ii}(x, z, t) E_*  - \max_{x ,z} 
 \tts{\sum}_{j \neq i} 
 ( |a_{ij}|  E_* + |a_{ij, 2}| E_*^2 
+ |a_{ij,3}(x, z, t) | ) >  \e_0 ,   \\
\]
for some $\e_0 > 0$, which is a much tighter constraint for stability
\footnote{
In comparison, if we use $L^2$-based stability estimates similar to those in \cite{chen2019finite,chen2021HL}, the estimates rely on a \textit{global} stability condition and require a much smaller  residual error to obtain nonlinear stability.
}.

\end{remark}

We remark that we have large damping factors in the far-field since the coefficients for the nonlocal terms in the linearized equations \eqref{eq:lin}, \eqref{eq:lin_main} decay. Thus, it is much easier to establish the local stability condition \eqref{eq:PDE_nondiag} for large $x$.

\subsubsection{Nonlinear stability and finite time blowup}\label{sec:non_stab}

To close the nonlinear estimates, for the bootstrap argument in Lemma \ref{lem:PDE_nonstab}, we choose the threshold 
\beq\label{eq:Ec}
E_* = 5\cdot 10^{-6}.
\eeq
We choose the bootstrap threshold guided by the quadratic inequality on $E_*$ \eqref{eq:PDE_nondiag} for $x$ (or $x,z$) in the region with small damping coefficients. Under this bootstrap threshold, the largest part of the nonlinear terms in the weighted $L^{\inf}(\vp_i)$ estimates are bounded by $(0.0415, 0.0415, 0.025) E_*$ (see Figure \ref{fig:N_linf_WG1}).
\footnote{
  One should compare $(0.0415, 0.0415, 0.025) E_*$ with remaining damping terms $d_i E_*$, where 
  $d_i$ is LHS of \eqref{eq:linf_diag_decay}.
}
Figures \ref{fig:stab_WG1} and \ref{fig:stab_WG1g} illustrate the nonlinear
stability conditions \eqref{eq:PDE_nondiag} corresponding to the full
$L^{\inf}(\vp_i)$ and $L^{\inf}(\vp_{g,i})$ estimates. The plotted quantities are 
piecewise lower bound of $E_*^{-1}\times \mathrm{LHS \ of \ }  \eqref{eq:PDE_nondiag}$ 
for $x \in [0, 10^{15}]^2$. Their positivity illustrates that, at the threshold $E_*$, the damping dominates
the nonlinear and error contributions.

The nonlinear stability conditions \eqref{eq:PDE_nondiag} associated with the
energy $E_4$ in \eqref{energy4} are stated in Appendix \ref{sec:ineq}. 
The above plots are included to \emph{only} illustrate the nonlinear stability estimates; see
Remark \ref{rem:plot}. For the
actual verification of the $L^\infty$ stability inequalities, 
such as \eqref{stab:linf_WG1},\eqref{stab:linf_WG4},\eqref{stab:linf_WG1g}, we discretize the
range of $x$, derive rigorous piecewise bounds on each subregion, and verify
the resulting inequalities using computer-assisted interval arithmetic. For
the $C^{1/2}$ stability inequalities \eqref{stab:hol}, we follow the methods in Section
\ref{sec:EE_hol_sum} and apply the same interval-arithmetic verification
procedure to the resulting piecewise bounds.

Since $\hat W_2(0)=0,  W(0) = W_1(0)$ initially
, if the initial perturbation  $W$ satisfies 
\bseq\label{eq:boot_stab_E4}
\beq\label{eq:boot_stab_E4:init}
E_4( W(0), W(0) ) < E_*, 
\eeq
using Lemma \ref{lem:PDE_nonstab}, we have 
\beq\label{eq:boot_stab_E4:dyn}
E_4( W_1(t), W(t)) < E_*,
\eeq
\eseq
for all time $t> 0$. With the estimates of $W_1$, we can control $\wh W_2$ using the estimates in Sections \ref{sec:rank1}, \ref{sec:W2}. 
In particular, using the energy $E_4$ and embedding inequalities, we verify
\beq\label{W_bound}
|| W_{1,i } + \wh W_{2,i}||_{\inf} < 200 E_*, \quad |c_{\om}(\om)| < 100 E_*.
\eeq
The bounds for $W_{1, i}, c_{\om}(\om)$ follow from the definition of the energy \eqref{energy1}, \eqref{energy3}, \eqref{energy4}. From the definitions of the weights $\vp_i, \vp_{g, i}, \mu_{g,i}$ \eqref{wg:linf_decay}, \eqref{wg:linf_grow}, \eqref{wg:EE}, we verify directly $|W_{1,i}| \leq \mu_{g,i}^{-1} \vp_{g,i}^{-1} E_4 < 100 E_* $. We verify $| \wh W_{2,i}| < 100 E_4 $ and collect this inequality in \eqref{W2_bound}. 
The normalization condition \eqref{eq:normal_pertb} and energy $E_4$ \eqref{energy4} imply $|u_x(0)| = |c_{\om}| < \mu_6 E_*< 100 E_*$.

Since $\hat W_2 |_{t =0} \equiv 0$ by \eqref{eq:bous_W2_appr}, we have $W_1 = (\om, \eta,\xi)
= (\om, \th_x, \th_y)$ initially.
Hence \eqref{eq:boot_stab_E4:init} 
is precisely a smallness condition on the \emph{full initial perturbation}  $(\om, \th_x, \th_y)$. Passing from stability estimates
\eqref{eq:boot_stab_E4},\eqref{W_bound} to finite time blowup follows the standard rescaling argument \cite{chen2019finite,chen2019finite2,chen2021HL}.

\subsection{Nonlinear stability and blowup theorem}\label{sec:EE_stab_thm}

Using the estimates and bootstrap argument in the preceding subsections, 
we prove the following nonlinear stability theorem, 
which strengthens Theorem~\ref{thm:main}. 
Before we state the result, we first recall the energies for $W_1, \hat W_2$
\[
E_1 \leq E_2 \leq E_3 \leq E_4
\]
defined in \eqref{energy1},  \eqref{energy2}, \eqref{energy3}, \eqref{energy4}.
 Let $\psi_i, \vp_i, \vp_{g, i}, g_i$ be the  weights defined in \eqref{wg:hol}, \eqref{wg:linf_decay}, \eqref{wg:linf_grow}, $ [  \cdot ]_{ C_{g_i}^{1/2}}$ be the H\"older seminorm \eqref{hol:g} in $\R_2^{++}$, and $\mu_{ij}, \tau_1, \tau_2 \in \R_+,  \mu_h, \mu_g \in \R_+^3$ be the parameters chosen in \eqref{wg:EE}. We define the initial energy $\Ein$ on three variables $ f_1, f_2, f_3 $ as 
 \begin{align}\label{energy}
 P_1(f)  & \teq \max(  \max_{i=1,2,3}  || f_i \vp_i ||_{\inf}, \tau_1^{-1}   \sqrt{2} || f_1  \vp_4 ||_{\inf}   ),  \notag \\
P_2(f) &\teq \max\nolim_{i=1,2,3}  ( \mu_{h, i} [ f_i \psi_i ]_{C_{g_i}^{1/2} } ) , \notag \\  
P_3(f) & \teq \max\nolim_{i=1,2,3} ( \mu_{g,i} \nlinf{ f_i \vp_{g, i}} )  , \notag \\
P_4(f) & \teq  \max\big(  \mu_{51}^{-1}  |c_{\om}( f_1 \chi_{\mf{ode}})|, 
\mu_{52}^{-1} | c_{\om}( f_2 \chi_{\mf{ode}} ) |,  \mu_{6}^{-1} | c_{\om}( f_1)|,  \mu_{62}^{-1} |c_{\om}( f_2)|,   \mu_7^{-1}  |f_{2,xy}(0)|, \mu_8^{-1} |f_{1, xy}(0)| \big), \notag \\
\Ein(f)   & \teq  \max\big( P_1(f), P_2(f), P_3(f), P_4(f) \big).
\end{align}
where $ c_{\om}(f) =  u_x(f)(0) = -\f{4}{\pi} \int_{\R_2^{++}} \f{y_1 y_2}{|y|^4} f(y) dy$, $\chi_{\mf{ode}}$ is defined in \eqref{eq:cw_chi}. 

With this notation, for $W = W_1 + \hat W_2$ with $\hat W_2= 0$, we have 
$
\Ein( W_1)  = E_4( W_1, W_1 ).$

We now state the full nonlinear stability and blowup result.

\begin{thm}
\label{thm:main_stab_full}
Let $(\bar{\th},\bar{\om}, \bar \uu,  \bar c_l, \bar c_{\om})$ be the approximate self-similar profile constructed in Section \ref{sec:ASS} and $E_* = 5 \cdot 10^{-6}$. 
Assume that the initial data satisfy the symmetry conditions \eqref{eq:Sym}, 
the regularity  $\om_0  \in  C^{k, \al}, \th_0  \in C^{k + 1, \al}$ 
 for some $k \geq 5, \al \in (0, 1)$, and 
\[
 \quad \Ein ( \om_0 - \bar \om,  \th_{0,x} - \bar \th_x ,  \th_{0, y} - \bar \th_y ) <  E_* .
\]

Let $W_1, \hat W_2$ be the solution to \eqref{eq:bous_decoup2} and $E_4$ be the energy \eqref{energy4}. Then for any $t \geq 0$, we have 
\[
\bga
E_4(W_1(t), W_1(t) + \hat W_2(t)) < E_* , \\
 || \om(t) - \bar \om ||_{L^{\inf}}, \  || \th_x(t) - \bar \th_x ||_{L^{\inf}} ,
\  || \th_y(t) - \bar \th_y ||_{\inf} < 200 E_* , \quad | u_x(t, 0 ) - \bar u_x( 0 )| , \quad  | \bar c_{\om} - c_{\om}(t)| < 100 E_* .
\ega
\]

Consequently, we can choose smooth initial data $\om_0, \th_0 \in C_c^{\inf}$ with finite energy
$\uu_0 \in L^2$ such that the corresponding solution of the physical equations \eqref{eq:bous_sys} blows up in finite time $T$.
\end{thm}

The regularity assumption \(k\ge 5\) is used only to establish the qualitative
regularity of \(W_1\) and \(\widehat W_2\) in Section~\ref{sec:reg_van}, which justifies that
the  weighted energy \(E_4\) is well defined; see Remark \ref{rem:E4_reg}.

\section{Finite time blowup of 3D axisymmetric Euler equations with boundary}\label{sec:euler}

In this section, we prove Theorems \ref{thm1b}, \ref{thm:euler} on finite time blowup of the axisymmetric Euler equations with smooth initial data and boundary. We follow the proof strategy from  our previous work \cite{chen2019finite2}. Section \ref{sec:dyn} introduces the dynamic rescaling formulation for 3D Euler.
In Sections \ref{sec:3Delli},  \ref{sec:main_vel_euler}, we establish the localized elliptic estimates. 
In Section \ref{sec:3Dnon}, we show that the differences between 
the nonlinear stability estimates for 3D Euler and 2D Boussinesq  are perturbative, using the smallness of the support size. In Section \ref{sec:euler_blowup_pf}, we close the bootstrap assumptions and prove finite-time blowup. Thus, the extension from 2D Boussinesq to 3D Euler is  obtained \emph{using a purely analytic argument, without any additional computer-assisted estimates}.

Throughout this section, since we mainly perform perturbative analysis, we treat all parameters and constants determined in the analysis of 2D Boussinesq as \textit{absolute} constants.

\vspace{0.1in}
\paragraph{\bf{Notations}} 
In this section, we use $x_1, x_2,  x_3$ to denote the Cartesian coordinates in $\R^3$, and 
\beq\label{eq:cylinder}
r= \sqrt{x_1^2 + x_2^2}, 
\quad z = x_3, \quad \vartheta = \arctan(x_2 /x_1) 
\eeq
 to denote the cylindrical coordinates. Let $\uu_{\mf{axi}}$ be the axi-symmetric velocity and $\vom = \na \times \uu_{\mf{axi}}$ be the vorticity vector. In the cylindrical coordinates, we have the following representation
\[
\uu_{\mf{axi}}(r, z) = u^r(r, z) \ee_r + u^{\th}(r, z ) \ee_{\th}  + u^z(r, z ) \ee_z, 
\quad 
\vom = \om^r(r, z) \ee_r + \om^{\th}(r, z ) \ee_{\th}  + \om^z(r, z ) \ee_z,
\]
where $\ee_r, \ee_{\th}, \ee_z$ denote the basis:
$
\ee_r = ( \f{x_1}{r}, \f{x_2}{r}, 0 ) , \ \ee_{\th} = (- \f{x_2}{r},  \f{x_1}{r}, 0  ), \  \ee_z = (0, 0, 1)$.

We study the 3D axisymmetric Euler equations in a cylinder $D = \{ (r,z) : r \in [0,1], z \in \BT \}, \BT = \R / ( 2 \BZ)$
that is $2$-periodic in $z$. 
The 3D axisymmetric Euler equations are given below:
\beq\label{eq:euler1}
\pa_t (ru^{\th}) + u^r (r u^{\th})_r + u^z (r u^{\th})_z = 0, \quad 
\pa_t (\f{\om^{\th}}{r}) + u^r ( \f{\om^{\th}}{r} )_r + u^z ( \f{\om^{\th}}{r})_z = \f{1}{r^4} \pa_z( (r u^{\th})^2 ).
\eeq
The radial and axial components of the velocity can be recovered from the Biot-Savart law
\beq\label{eq:euler2}
-(\pa_{rr} + \f{1}{r} \pa_{r} +\pa_{zz}) \td{\phi} + \f{1}{r^2} \td{\phi} = \om^{\th}, 
 \quad  u^r = -\td{\phi}_z, \quad u^z = \td{\phi}_r + \f{1}{r} \td{\phi}  
\eeq
with a periodic boundary condition in $z$ and a no-flow boundary condition on the boundary $r = 1$
\beq\label{eq:euler21}
\td{\phi}(1, z ) = 0 .
\eeq

We consider solution $\om^{\th}$ with odd symmetry in $z$, which is preserved by the equations dynamically. Then $\td \phi$ is also odd in $z$. Moreover, since $\td \phi$ is 2-periodic in $z$, we obtain 
\beq\label{eq:euler22}
\td \phi(r, 2k-1)  = 0 . \quad \textrm{for all \ } k \in \BZ .
\eeq
This setup of the problem is essentially the same as that in \cite{luo2014potentially,luo2013potentially-2}.

Due to the periodicity in $z$ direction, it suffices to consider the equations in the first period $D_1 =\{ (r,z) : r \in [0,1], |z| \leq 1 \} $. We have the following pointwise estimate on $\td \phi$ from \cite{chen2019finite2}, which will be used to estimate $\td \phi$ away from the $\supp(\om^{\th})$ in Section \ref{sec:3Delli}.

\begin{lem}\label{lem:biot1}
Let $\td{\phi}$ be a solution of \eqref{eq:euler2}-\eqref{eq:euler21}, and $\om^{\th} \in C^{\al}(D_1)$ for some $\al>0$ be odd in $z$ with $\supp(\om^{\th}) \cap D_1\subset \{ (r ,z) : (r-1)^2 + z^2 < 1/4  \}$. 
For $  \f{1}{4} < r \leq 1, |z|\leq 1$, we have
\[
| \td{\phi}(r, z) |  \les \int_{D_1} | \om^{\th}(r_1, z_1) | \B( 1 + |\log( (r-r_1)^2 + (z- z_1)^2) |  \B) r_1 d r_1 d z_1.
\] 
\end{lem}

If the initial data $u^{\th}$ of \eqref{eq:euler1}-\eqref{eq:euler21} is non-negative, 
$u^{\th}$ remains non-negative before the blowup, if it exists. Then, $u^{\th}$ can be uniquely determined by $(u^{\th} )^2$.
We introduce the following variables 
\beq\label{eq:omth}
\td{\th} \teq (r u^{\th})^2,  \quad \td{\om} = \om^{\th} / r.
\eeq

We reformulate \eqref{eq:euler1}-\eqref{eq:euler21} as 
 \beq\label{eq:euler31}
\bal
\pa_t \td{\th} + u^r \td{\th}_{ r} + u^z \td{\th}_{z} &= 0,  \quad  \pa_t \td{\om} + u^r \td{\om}_r + u^z \td{\om}_{ z} = \f{1}{r^4} \td{\th}_z  ,  \\
-( \pa^2_{r} +\f{1}{r} \pa_r + \pa_z^2 -\f{1}{r^2} )  \td{\phi} & = r \td{\om} , \quad \td{\phi}(1, z) = 0 ,  \quad u^r = -  \td{\phi}_{ z} , \quad 
u^z = \f{1}{r} \td{\phi }+  \td{\phi}_{ r} .
\eal
\eeq

\subsection{Dynamic rescaling formulation}\label{sec:dyn}
We will construct blowup solution near $r=1, z=0$. We introduce new coordinates $(x, y)$ centered at $r = 1, z= 0$ 
\beq\label{eq:euler_polar}
x =  C_l(\tau)^{-1} z, \quad  y = (1-r) C_l(\tau)^{-1}, 
\eeq
where $C_l(\tau)$ is defined below \eqref{eq:rescal42}. By definition, we have 
\beq\label{eq:label_rs}
z = C_l(\tau)x, \quad r = 1 - C_l(\tau) y . 
\eeq

We consider the following dynamic rescaling formulation centered at $r = 1, z= 0$
\beq\label{eq:rescal41}
\bal
\th(x, y, \tau) &= C_{\th}(\tau) \td{\th}( 1 - C_l(\tau) y,   C_l(\tau) x , t(\tau) ), \\
  \om(x, y, \tau) &=  C_{\om}(\tau) \td{\om}( 1 - C_l(\tau) y,  C_l(\tau) x  , t(\tau)) , \\
 \phi(x, y, \tau)  & =  C_{\om}(\tau) C_l(\tau)^{-2} \td{\phi} (1 - C_l(\tau) y,  C_l(\tau) x, t(\tau)),
\eal
\eeq
where $C_l(\tau), C_{\th}(\tau), C_{\om}(\tau), t(\tau)$ are given by 
\beq\label{eq:rescal42}
\bal
C_{\th}(\tau )  & = C^{-1}_l(0)  \exp\big( \int_0^{\tau} c_{\th} (s)  d s \big),\
C_{\om}(\tau)  =  \exp \big( \int_0^{\tau} c_{\om} (s)  d s \big), \
 C_l(\tau) = C_l(0) \exp\big( \int_0^{\tau} -c_l(s) ds \big) , 
\eal
\eeq
and $t(\tau) = \int_0^{\tau} C_{\om}(s ) d s$.
The rescaling parameters $c_l(\tau), c_{\th}(\tau), c_{\om}(\tau)$ satisfy $c_{\th}(\tau) = c_l(\tau ) + 2 c_{\om}(\tau)$ and $C_{\th}, C_l, C_{\om}$ satisfy $C_{\th} = C^2_{\om} C_l^{-1}$, by the same reasoning as in \eqref{eq:rescal3}. Since we rescale the cylinder $D_1 = \{ (r, z) : r \leq 1, |z|\leq 1  \}$, from \eqref{eq:euler_polar}, the domain for $(x, y)$ is 
\beq\label{eq:rescale_D}
\td D_1 \teq \{ (x, y) :  |x| \leq C_l^{-1}, y \in [0, C_l^{-1}] \}.
\eeq

The change of variable \eqref{eq:label_rs} leads to a  minus sign  between $\pa_r$ and $\pa_y$
\[
\pa_y \th = -C_{\th} C_l(\tau) \td{\th}_r , \quad  \pa_y \om = - C_{\om} C_l(\tau) \td{\om}_r,  \quad \pa_y \phi = - C_{\om} C_l(\tau)^{-1} \td{\phi}_r.
\]

Let $(\td{\th}, \td{\om})$ be a solution of \eqref{eq:euler31}. 
As in Section \ref{sec:dyn_res}, we use $t$ for the rescaled time below.  
Then it is easy to show that $\om, \th$ satisfy
\[
\th_t +  c_l \xx \cdot \na \th + (-u^r) \th_y +  u^z  \th_x  = c_{\th} \th , \quad \om_t + c_l \xx \cdot \na \om  + (-u^r)  \om_y +  u^z \om_x = c_{\om } \om + \f{1}{r^4} \th_x .
\]
The Biot-Savart law in \eqref{eq:euler31} depends on the rescaling parameter $C_l, \tau$
\[
-(\pa_{xx}  + \pa_{yy}) \phi + \f{1}{r} C_l \pa_y \phi  + \f{1}{r^2 } C^2_l \phi = r \om ,
\quad   u^r(x, y) = -  \phi_x(x, y),  \quad u^z(x, y) =  \f{1}{r}  C_l(\tau)  \phi(x, y) -  \phi_y(x ,y),  
\]
where $r = 1 - C_l(\tau) y$ \eqref{eq:label_rs}. We introduce 
\[
u(x, y) = u^z(x, y), \quad v(x, y) = - u^r(x, y) .
\]
Then, we can further simplify
\beq\label{eq:euler4}
\bal
& \th_t + (c_l \xx + \uu \cdot \na ) \th = c_{\th} \th , \quad  \om_t  + ( c_l \xx + \uu \cdot \na ) \om = \th_x + \f{1 - r^4}{r^4} \th_x, \\
& -(\pa_{xx}  + \pa_{yy}) \phi + \f{1}{r} C_l \pa_y \phi + \f{1}{r^2} C_l^2 \phi= r \om , \quad 
u =    - \phi_y + \f{1}{r} C_l \phi, \quad v  = \phi_x ,
\eal
\eeq
with boundary condition $\phi(x, 0 ) \equiv 0$. If $C_l\ll 1 $, \eqref{eq:euler4} is  the
perturbation of
 the dynamic rescaling equations \eqref{eq:bousdy1} for the Boussinesq system.
We seek solutions of \eqref{eq:euler4} with symmetries
\beq\label{eq:euler_sym}
\om(x, y)= -\om(-x, y), \quad \th(x, y) = \th(-x, y),
\quad \phi(x, y) = - \phi( - x, y),
\eeq
which are preserved by the equations. Hence, it suffices to solve \eqref{eq:euler4} on $x, y \geq 0$. 

\begin{definition}\label{def:supp}
We define the size of support of $(\th, \om)$ of \eqref{eq:euler4} 
\[
S(\tau) = \max( \mathring{S}(\tau), \mathring{S}(0), 1),
\quad \mathring{S}(\tau) 
\teq  \mathrm{essinf}  \{\rho :  \th(x, y,\tau) =0, \om(x, y, \tau ) = 0 \textrm{ for } x^2 + y^2 \geq \rho^2 \} .
\]
\end{definition}

We impose $S(\tau) \geq 1$ to simplify later estimates. 
We will construct initial data with compact support $S(0) < + \infty$.
After rescaling the spatial variable, the support of $(\td{\th}, \td{\om})$ of \eqref{eq:euler31} satisfies 
\[
\mathrm{supp}  \ \td{\th}(t(\tau)),  \ \mathrm{supp}  \ \td{\om}(t(\tau)) \subset \{ (r, z) :  ( (r-1)^2 + z^2)^{1/2} \leq  C_l(\tau) S(\tau)  \}.
\]

For a large $\nu>0$ to be chosen, we define the weights $\vp_{3,m}, \vp_{g,3,m}$ in \eqref{wg:xi_modi}  modifying $\vp_3, \vp_{g,3}$, and define the modified  energy 
$E_{\mf{in}, m}$ of $\Ein$ in \eqref{energy} 
with $\vp_3, \vp_{g,3}$ replaced by $\vp_{3,m}, \vp_{g,3,m}$:
\beq\label{energy:modi}
\bal
E_{\mf{in},m}(f) \teq \max & \big(   \max_{i=1,2}(  || f_i \vp_i ||_{\inf} ) , \max_{i=1,2}( \mu_{g, i} \nlinf{ f_i \vp_{g, i}}  ), \tau_1^{-1}   \sqrt{2} || f_1  \vp_4 ||_{\inf}, \\
& \quad  \nlinf{f_3 \vp_{3,m}} , \mu_{g, 3} \nlinf{f_3 \vp_{g,3,m} }   ,  P_2(f), P_4(f) \big).
\eal
\eeq
where $P_2(\cdot), P_4(\cdot)$ are defined in \eqref{energy}.
We use the same functional as in $P_4(f)$ in 2D Boussinesq energy \eqref{energy}: $c_{\om,2D}(f) = u_{x,2D}(f)(0) = -\f{4}{\pi} \int_{\R_2^{++}} \f{y_1 y_2}{|y|^4} f(y) dy$. 
\footnote{
For 3D Euler, \eqref{eq:E_normal2} for $c_{\om}(f)$ perturbation differs from its 2D counterpart \eqref{eq:normal_pertb} by $ C C_l S \| f \|_{X_1}$, which can be made arbitrary small. See 
\eqref{eq:elli_main3}, \eqref{eq:u_dif_3D2}.
}
 
 We now state a more precise version of Theorem \ref{thm1b}.

\begin{thm}\label{thm:euler}

There exists a large $\nu >0$ such that the following hold. Let $(\bar{\th}_0,\bar{\om}_0, \bar c_l)$ be the approximate self-similar 
profile constructed in Section \ref{subsec:non_appr}, which depends on $\nu$, and let $\Einm$ be defined in \eqref{energy:modi}, which depends on $\nu$ via the weights $\vp_{3,m}, \vp_{g,3,m}$ \eqref{wg:xi_modi}. Set $E_* = 5 \cdot 10^{-6}$. 
 Assume that 
the initial data $\th_0, \om_0$ of \eqref{eq:euler4} are compactly supported with size $S(0)$  defined in Definition \ref{def:supp}. 
 Assume also that  $\th_0$ is even in $x$, $\om_0$ is odd in $x$,  $ \om_0  \in  C^{\kin, \ain}$, $\th_0  \in C^{\kin + 1, \ain}$ for some $\kin \geq 5, \ain \in (0, 1)$ and
$
 \Einm( \om_0 - \bar \om_0,  \th_{0,x} - \bar \th_{0, x} ,  \th_{0, y} - \bar \th_{0, y} ) <  E_*.
 $

There exists a constant  $C( S(0))$ depending on $S(0)$ such that if the initial rescaling factor $C_l(0)$ \eqref{eq:rescal42} satisfies $C_l(0) < C( S(0))$, 
\footnote{
The rescaling factor $C_l(t) = C_l(0) e^{- \bar c_l t}$ 
(see \eqref{eq:E_normal:c} enters equation \eqref{eq:euler4} 
as a parameter.
}
we have 
\[
 || \om(t) - \bar \om_0 ||_{L^{\inf}}, \  || \th_x(t) - \bar \th_{0,x} ||_{L^{\inf}} ,
\  || \th_y(t) - \bar \th_{0, y} ||_{\inf} < 200 E_* , \  | u_x(t, 0) - \bar u_x(0)| , \  | \bar c_{\om} - c_{\om}(t)| < 100 E_*,
\]
for any $t \geq 0$, 
where $ \bar \uu,  \bar c_{\om}$ are the profile variables associated with  $\bar \th_0, \bar \om_0, C_l(0)$. In particular, there exists smooth initial data $\om_0, \th_0 \in C_c^{\inf}$ in this class with finite energy $||\uu_0||_{L^2} < +\infty$, such that the corresponding solution to the physical equations 
\eqref{eq:euler1}-\eqref{eq:euler21} 
 blows up in finite time $T$. 
\end{thm}

We need to choose a small rescaling factor $C_l(0)$ so that the solution in the physical space is confined in the cylinder, 
 $D_1 = \{ (r,z) : r \in [0,1], |z| \leq 1 \}$, which is not scaling invariant.

\subsection{The elliptic estimates} \label{sec:3Delli}

We first bound $\phi$ away from $\supp(\omega)$, then localize the equation and perform weighted $L^{\infty}$ and H\"older estimates. 
We will show that within the support of $\om, \th$, the elliptic estimates are the same as those in the 2D Boussinesq  up to lower order terms, which can be made arbitrarily small. Throughout this section, we assume that $\om(x, y)$ is odd in $x$.

\begin{remark}[\bf Small parameters]\label{rem:order1}
We treat $ C_l(\tau), C_l(\tau) S(\tau)$ as small parameters in the following estimates, 
and justify this by taking $C_l(0)$ sufficiently small at the 
end of the proof. One can essentially treat $C_l(\tau) \approx 0$. From \eqref{eq:label_rs}, $r = 1 - C_l y  \approx 1$ in the support of the solution.
\end{remark}

The term $\f{1}{r} C_l \pa_y \phi$  in \eqref{eq:euler4} leads to a few technical difficulties. 
To overcome these, we multiply \eqref{eq:euler4} by the integrating factor $r^{1/2}$. Using $\pa_y r^{1/2} = - C_l r^{-1/2} /2, \pa_{yy} r^{1/2} = - \f{1}{4} C_l^2 r^{-3/2}$, 
\[
\pa_{yy} ( \phi r^{1/2}) = r^{1/2} \pa_{yy} \phi + 2 \pa_y \phi \pa_y r^{1/2}
+ \phi \pa_{yy} r^{1/2}
=  r^{1/2} \pa_{yy} \phi  - \f{ C_l }{r^{1/2}} \pa_y \phi
- \f{1}{4} C_l^2 r^{-3/2} \phi,
\]
and denoting $\D = \pa_{xx} + \pa_{yy}$, we rewrite \eqref{eq:euler4} as follows 
\beq\label{eq:E_elli}
- \D \phi_\dag = \Om_\dag - \f{ a C_l^2}{r^2} \phi_\dag, \quad \Om_\dag = \om r^{3/2}, \quad \phi_\dag \teq \phi r^{1/2}, \quad a = \f{3}{4} .
\eeq
Since $r, r^{-1}$ are smooth within the support of $\om, \th$, with the estimate of $\phi r^{1/2}$, we can recover that of $\phi$. Our goal is to show that $\phi_\dag, \phi$ enjoy estimates similar to those for $ (-\D_{2D})^{-1} \om$.

\subsubsection{Estimate of $\phi$ away from the support}

To localize the elliptic equations, we first estimate $\phi$ away from the support of the solution. Based on Lemma \ref{lem:biot1}, we have the following estimate. 

\begin{lem}\label{lem:far}
Let $S(\tau)$ be the support size of $\om(\tau), \th(\tau)$. Assume the assumptions in Lemma \ref{lem:biot1} hold  and $C_l(\tau) S(\tau) < \f{1}{4}$. For any 
$\xx=(x,y) \in \td D_1$ \eqref{eq:rescale_D} with $|\xx| > 2 S$,
$r = 1- C_l(\tau) y > \f14$   and $\b \in [0, 1)$, the solution to \eqref{eq:euler4} satisfies 
\[
|\phi(\xx)|  \les || \om (1 + |\xx|^{\b} ) ||_{L^{\inf}} (1 + | \log( C_l |\xx|)|) S^{2-\b}. 
\]
\end{lem}

Since $x$ is away from the support, the proof follows from the rescaling relation \eqref{eq:rescal41} and the estimate in Lemma \ref{lem:biot1} by putting $ \om(y) (1 + |y|^{\b})$ in $L^{\inf}$. We defer the proof to Appendix \ref{app:euler_stream}.

\subsubsection{Localize the elliptic equation}

We will take advantage of the fact that $C_l(\tau) S(\tau)$ can be extremely small and localize the elliptic equation. Firstly, we assume that $C_l(\tau) S(\tau) < \f{1}{4}$. Recall the relation \eqref{eq:label_rs} about $r$. Within the support, we have $r = 1 - C_l y  \geq \f{3}{4}, r^{-1} \les 1$.

Let $\chi(\cdot) : \R_2^+ \to [0, 1]$ be a smooth cutoff function even in $x_1$, such that $\chi(x) = 1$ for $|x| \leq 1$, $\chi(x) = 0$ for $|x| \geq 2$. 
\footnote{
  The reader should not confuse $\chi$ with cutoff functions used in Section \ref{sec:finite_rank}, \ref{sec:EE} and in \eqref{eq:cutoff_near0_all}. They are different.
}
It is easy to verify that 
 \beq\label{eq:chi2}
 | \na^k \chi(x/ R)| \les  R^{-k} \one_{ R \leq |x| \leq 2 R} ,
 \eeq
for $1\leq k\leq 5$. 
Next, we choose several radii and define the related cutoff function
\beq\label{eq:chi_Ri}
R_i = 4^{-i} C_l^{-1} , \quad \chi_i(x) \teq \chi( x / R_i ),  \quad 1\leq i \leq 5.
\eeq
Clearly, $\chi_i = 1$ in the support of $\chi_{i+1}$. Multiplying \eqref{eq:E_elli} with $\chi_i$ 
yields the equation of $ \phi_\dag \chi_i$
\beq\label{eq:E_elli1}
- \D( \phi_\dag  \chi_i  ) 
= \Om_\dag \chi_i -  Z_{1, i} - Z_{2, i}, \quad 
Z_{1, i} \teq   2 \na \phi_\dag  \cdot  \na \chi_i + \phi_\dag \D \chi_i , \quad 
Z_{2, i} \teq  \f{a C_l^2}{r^2} \phi_\dag \chi_i,
\eeq
with boundary condition 
\[
 (\phi_\dag \chi_i)( 0, y) = 0, \quad (\phi_\dag \chi_i)(x, 0) = 0.
\]
After we localize the elliptic equation, \eqref{eq:E_elli1} can be seen as an elliptic equation in $\R_2^{+}$ with compactly supported source term. Since the solution $\phi_\dag \chi_i$ decays for large $|x|$, it agrees with the solution defined by the Green function $\log( |x- y|)$ in the upper half space:
\beq\label{eq:lap1}
( (-\D)^{-1} f)(x)  = -\f{1}{2\pi} \int_{\R_2^+} ( \log |x - y| 
 - \log|  (x_1 - y_1, x_2 + y_2) |)  f(y) dy 
= -\f{1}{2\pi} \int_{\R_2}  \log |x - y|  F(y) dy.
\eeq
where $F$ is the odd extension of $f$ from $\R_2^+$ to $\R_2$. Similar formula also holds for $ \na (-\D)^{-1} f$
\beq\label{eq:lap2}
\pa_i (-\D)^{-1} f = -\f{1}{2\pi} \int_{\R_2} \f{x_i - y_i}{ |x-y|^2} F(y) dy.
\eeq

\paragraph{\bf{Ideas of the estimates}}

Our goal is to apply the weighted nonlocal estimates for  
2D Boussinesq, e.g. \eqref{eq:holest},
to estimate $\phi_\dag \chi_i$ with source term $ \Om_\dag \chi_i - Z_{1, i} - Z_{2, i}$ in 
\eqref{eq:E_elli1}, and to show that such estimates are perturbations of $-\D \phi = \om$.
For this purpose, 
through several estimates, we show that $Z_{1, i}, Z_{2, i}$, with corrections, are small in the weighted norm $X_1$  related to energy $E_4$ \eqref{energy4m}.

We require that the support satisfies
\beq\label{eq:E_supp}
S(\tau) < R_i, \quad 
\iff  C_l S(\tau) < 4^{-i}, \quad 1\leq i \leq 5 ,
\eeq
so that $\Om_\dag \chi_i = \Om_\dag$. We will choose $C_l S $ to be sufficiently small.

We need to estimate the $L^{\inf}$ norm of $\na \phi_\dag$ and its H\"older norm. We will first estimate $\na \phi_\dag$ for $|x| \leq 2 R_2 = R_1/2$, and then $\na^2 \phi_\dag$ for $|x| \leq R_2/2$. Once we obtain the estimates of $\na \phi_\dag, \na^2 \phi_\dag$, due to the small parameters on the right hand side of \eqref{eq:E_elli1} and the decay of the solution, we establish the desired estimate. We need several weighted estimates of the Laplacian in $\R_2^{+}$.

\begin{lem}\label{lem:E_ker}
Suppose that $-\D \phi = \om$ in $\R_2^{+}$, $\om$ is odd in $x$, and $\phi$ satisfies the Dirichlet boundary condition. For $\al >0$ and $\b \in (0, 1)$, we have 
\bseq\label{eq:E_ker_est}
\beq\label{eq:E_ker_est:a}
 | \na \phi | \les_{\al, \b}  |x| \wedge |x|^{1-\b} || \om (|x|^{-\al} + |x|^{\b}) ||_{L^{\inf}}.
\eeq
For $ \al \in (0, 2)$ and $\om \in C^{1,\g}$ with some $\g > 0$, 
\footnote{
When $\al > 1$, if $\om_x(0) \neq 0$,  the upper bound is $\infty$ and the inequality is trivial. 
If $\om_x(0) = 0$, from the derivations at the beginning of Section \ref{sec:appr_vel_near0}, 
we have $\na^3 \phi(0) = 0$. Thus, there are no cubic-order terms in the expansion of 
$\phi - \phi_{xy}(0) x_1 x_2$ near $0$. The same reasoning applies to the 4th and 5th estimates in Lemma \ref{lem:E_ker}.
\label{fn:phi_expand}
}
we have 
\beq\label{eq:E_ker_est:b}
\bal
 | \na ( \phi - \phi_{xy}(0) x_1 x_2 ) |  & \les_{\al, \b}  |x|^{1 + \al} || \om ( |x|^{-\al} +|x|^{\b} )||_{L^{\inf}} , \\
  |  \phi - \phi_{xy}(0) x_1 x_2  | & \les_{\al, \b}  |x|^{2 + \al} || \om ( |x|^{-\al} +|x|^{\b} )||_{L^{\inf}}  .
  \eal
\eeq
For $\al \in (2, 3]$ and $\om \in C^{2, \g}$ with some $\g>0$, we have 
\beq\label{eq:E_ker_est:c}
\bal
 |\na (\phi - \phi_{xy}(0) x_1 x_2 ) | & \les_{\al, \b} |x|^{3} || \om ( |x|^{-\al} +|x|^{\b} )||_{L^{\inf}}, \\
   |  \phi - \phi_{xy}(0) x_1 x_2  | & \les_{\al, \b}  |x|^{4} || \om ( |x|^{-\al} +|x|^{\b} )||_{L^{\inf}}  .
   \eal
\eeq
\eseq
\end{lem}

We will mostly use $\al = 1, \, 2.9, \, \b = \f{1}{16}$ related to the weight $\vp_1, \vp_{g, 1}$ for $\om_1$ 
\eqref{wg:linf_decay}, \eqref{wg:linf_grow}.  We prove \eqref{eq:E_ker_est:a} below, and defer the proof of \eqref{eq:E_ker_est:b}, \eqref{eq:E_ker_est:c} to Appendix \ref{app:euler_stream}, which are similar.

\begin{proof}[Proof of \eqref{eq:E_ker_est:a}]
In the following proof, the implicit constant in $\les$ can depend on $\al, \b$. We drop it to simplify the notations. Denote by $M = || \om ( |x|^{-\al} + |x|^{\b})||_{\inf}$. Clearly, we have 
\beq\label{eq:E_ker1}
\bal
 |\om|  & \leq \min( |x|^{\al}, |x|^{-\b} ) M, \\
\pa_i  \phi & =  C_i \int K_i(x - y)  W(y) dy, \quad K_i(z)  = \f{z_i}{|z|^2},
\eal
\eeq
where $W$ is the odd extension of $\om$ from $\R_2^+$ to $\R_2$. Then $W$ is odd in both $x$ and $y$. Without loss of generality, we consider $i = 2$. For a fixed $x$, we partition the integral into three regions: 
\[
Q_1 = \{ y : |y| \geq 2|x|  \}, \quad Q_2 = \{ y: |y-x| \leq  |x| / 2 \},  \quad Q_3 = (Q_1 \cup Q_2)^c .
\]

In $Q_1$, symmetrizing the kernel, we need to estimate 
\[
\bal
I_1 \teq  \int_{Q_1} K_2(x-y) W(y)  dy
 = \int_{Q_1, y_1 \geq 0} ( K_2( x_1 - y_1, x_2 - y_2) - K_2( x_1 + y_1, x_2 - y_2) ) W(y) dy . 
\eal
\]
Since $K_2(z)$ is odd in $z_2$ and even in $z_1$, $|\na K_2(z)| \les |z|^{-2}$,  for $|y| \geq 2|x|$, we get 
\[
 |K_2( x_1 - y_1, x_2 - y_2) - K_2( x_1 + y_1, x_2 - y_2) |
 =  |K_2( y_1 - x_1, x_2 - y_2) - K_2( x_1 + y_1, x_2 - y_2) |
 \les  \f{ |x_1|}{|y|^2}.
\]
Using \eqref{eq:E_ker1}, we get 
\[
\bal
|I_1| & \les M |x_1|  \int_{|y|\geq 2|x|} |y|^{-2} \min( |y|^{\al}, |y|^{-\b}) dy 
 \les M |x_1| \min(1 , |x|^{-\b})
\les M \min( |x|, |x|^{1-\b}).
\eal
\]

In $Q_2$, since $|x-y| \leq |x|/2$, we have $|y| \asymp |x|$, $\min( |x|^{\al} , |x|^{-\b}) \asymp \min( |y|^{\al} , |y|^{-\b}) $, and obtain
\[
\B| \int_{Q_2} |K_2(x-y) W(y) d y \B| 
\leq M  \int_{ |x-y| \leq |x| /2} |x-y|^{-1} \min( |y|^{\al}, |y|^{-\b}) dy
\les M \min( |x|^{\al}, |x|^{-\b})  |x| .
\]

In $Q_3$, we have $|x|/2 \leq |x-y| \leq 3|x|, |y| \leq 2 |x|$. Using this estimate and \eqref{eq:E_ker1}, we obtain 
\[
\bal
|\int_{Q_3}  K_2(x- y) W(y) dy | 
&\les  M \int_{Q_3} |x-y|^{-1}  \min( |y|^{\al}, |y|^{-\b}) dy 
\les M |x|^{-1}  \int_{ |y| \leq 2 |x|}  \min( |y|^{\al}, |y|^{-\b}) dy \\
& \les M |x|^{-1}  \min( |x|^{2 + \al}, |x|^{ 2 -\b}) 
\les M  \min( |x|^{1 + \al}, |x|^{ 1 -\b}) .
\eal
\]
Combining the above estimates, we prove \eqref{eq:E_ker_est:a} in Lemma \ref{lem:E_ker}.
\end{proof}

\subsubsection{Estimate of $\na \phi_\dag$}

We have the following estimate of $\na \phi_\dag$ in $|x| \leq R_2$. 
\begin{prop}\label{prop:C1}
Let $\phi_\dag$ be the solution in \eqref{eq:E_elli} and $\al > 0, \b \in (0, 1)$. There exists some absolute constant $\nu_1(\al, \b) < \f{1}{32}$ such that if $C_l(\tau) (1 + S(\tau)) < \nu_1$, we have 
\[
\max_{|x| \leq 2 R_2} | \na \phi_\dag| ( |x|^{-1} + | x|^{-1 + \b})
\les  || \Om_\dag  (|x|^{-\al} + |x|^{\b}) ||_{\infty}.
\]

For $8 S \leq |x| \leq R_1 / 2 $ away from the support of $\om$, we have an improved 
decay estimate 
\beq\label{eq:prop_C1_impr}
 |\na \phi_\dag| \les  || \Om_\dag  (|x|^{-\al} + |x|^{\b}) ||_{\infty} 
\big( \f{S^{3-\b} }{ |x|^2}  + |x|^{1-\b} (C_l S)^{2- \b} \big) .
\eeq
For large $|x|$, compared to the first estimate,  \eqref{eq:prop_C1_impr} has a small factor 
$(C_l S)^{2- \b}$ for $|x|^{1-\b}$.

\end{prop}

We use the smallness related to the support size of $\om$ to obtain the improved estimates \eqref{eq:prop_C1_impr}.

In the following estimate, since $r$ is sufficiently close to $1$ within the support of $\om$, one may simply treat $\Om_\dag = \om r^{3/2}$ and $\om$ as the same.

\begin{proof}

We choose $i=1$ in \eqref{eq:E_elli1}. Denote 
\[
B_1 \teq \max_{ |x| \leq  2 R_2} |\na \phi_\dag ( |x|^{-1} + |x|^{-1 + \b})| , \quad 
M_1 \teq  || \Om_\dag  (|x|^{-\al} + |x|^{\b}) ||_{\infty}.
\]

Inverting $-\D$ and then apply $\na$ to \eqref{eq:E_elli1} with $i=1$, we obtain 
\beq\label{eq:prop_C10}
\bga
\na( \phi_\dag \chi_1) = \na (-\D)^{-1} ( \Om_\dag \chi_1) 
-  \na (-\D)^{-1} Z_{1,1} - \na (-\D)^{-1} Z_{2, 1} \teq I_1 + I_2 + I_3 , \\
Z_{1, 1} = 2 \na \phi_\dag  \cdot  \na \chi_1 + \phi_\dag \D \chi_1, \quad 
\quad Z_{2, 1} = \f{a C_l^2}{r^2} \phi_\dag \chi_1 .
\ega
\eeq

Our goal is to prove the following estimate 
\beq\label{eq:prop_C11}
B_1 \leq C_{\al, \b}( M_1 + (C_l S)^{2 -\b}  \cdot B_1  ).
\eeq
Then as long as $C_l S$ is small, we can obtain the bound for $B_1$.

For $I_1, I_3$, applying Lemma \ref{lem:E_ker}, we get
\[
\bal
|I_1| & \leq \min( |x|, |x|^{1-\b}) M_1, \quad  |I_3 | 
\leq \min( |x|, |x|^{1-\b}) || Z_{2, 1} ( |x|^{-1} + |x|^{\b}) ||_{\infty}, \quad 
Z_{2, 1}   =  \f{ a C_l^2 }{r^2} \phi_\dag \chi_1 .
\eal
\]

It suffices to bound the norm of $Z_{2, 1}$. For $|x| \leq 2 S < R_2$ \eqref{eq:E_supp}, using the definition of $B_1$, $\phi(0) = 0$, and integration, we get 
\[
\bal
 & |\phi_\dag|   \les B_1 \min( |x|^2 , |x|^{2 - \b}) ,  \\
\eal
\]
which along with $r^{-1} \les 1$ within the support of $\chi_1$ yields 
\[
 |Z_{2,1}|  ( |x|^{-1} + |x|^{\b} )
 \les B_1  C_l^2 \min( |x|^2 , |x|^{2 - \b}) ( |x|^{-1} + |x|^{\b} ) 
\les  B_1 C_l^2 ( |x|^2 + |x|) \les B_1  C_l^2 S^2,
\]

For $|x| \geq 2S > 2$, using Lemma \ref{lem:far} on $\phi_\dag$, we yield 
\[
|Z_{2, 1} (|x|^{-1} + |x|^{\b})|
\les  C_l^2 (1 +  |\log (C_l |x| ) | )  S^{2-\b} |x|^{\b} \one_{ |x| \leq  2 R_1}  
|| \om (1 + |x|^{\b})||_{\infty} .
\]

By definition, we have $C_l R_1  \in [0, 1/4]$ \eqref{eq:chi_Ri}. Within the support of $\om$, 
since $|\om r^{3/2}| = |\Om_\dag|$ and $r \asymp 1$, we have  $\Om_\dag \asymp \om$. Hence, we obtain 
$ (1 + | \log C_l x |) | C_l x|^{\b} \les 1$ and 
\beq\label{eq:prop_C1_Z21}
|Z_{2, 1} (|x|^{-1} + |x|^{\b})| \les (C_l S)^{2-\b} || \Om_\dag (1 + |x|^{\b}) ||_{\infty}.
\eeq

Next, we estimate $I_2$ defined in \eqref{eq:prop_C10}. Since $R_2 \leq \f{R_1}{4}$, for $|x| \leq 2 R_2$ and $y \in \supp(Z_{1,1})$, we have $2 |x| \leq |y|$. 
We estimate a typical term in $Z_{1, 1}$. To use the formula \eqref{eq:lap1}, \eqref{eq:lap2}, we extend $\phi_\dag$ naturally from $\R_2^+$ to $\R^2$ as an odd in $x_2$,
and $\chi_1$ as an even function in $x_2$. For $i, j \in \{ 1, 2 \}$, using integration by parts,
we get  
\[
\bal
J &\teq \pa_i( -\D)^{-1} ( \pa_j \phi_\dag \pa_j \chi_1 )
= C \int_{ \R_2}  \f{ x_i - y_i}{ |x-y|^2} ( \pa_j \phi_\dag  \pa_j \chi_1) 
= - C \int_{ \R_2} \pa_j (  \f{ x_i - y_i}{ |x-y|^2} \pa_j \chi_1 ) \phi_\dag \teq J_1 + J_2, \\
J_1 & \teq  - C \int_{\R_2} \pa_j  \f{ x_i - y_i}{ |x-y|^2}  \pa_j \chi_1 \phi_\dag, \quad 
J_2 \teq - C \int_{\R_2}  \f{ x_i - y_i}{ |x-y|^2}  \pa^2_j \chi_1 \phi_\dag .
\eal
\]
For $\pa_j$ with $j=2$, the boundary term on $x_2=0$ vanishes since 
$\phi_\dag(x_1,0) = 0,  \pa_2 \chi_1(x_1,0)=0$.

Since the singularity $x$ of the kernel is away from the support of the integrand, the integrand is smooth. 
We estimate $J_1$ with $i = 2, j= 1$. Estimates of other cases and the second term are similar. Denote $K(z) = \f{z_1 z_2}{|z|^4}$. Since $\pa_1 \chi_1 \phi_\dag$ is even in $y_1$, symmetrizing the kernel in $y_1$, we get 
\[
J_1 = C \int K(x-y) \pa_1 \chi_1 \phi_\dag(y) dy 
= C \int_{y_1 \geq 0} ( K(x_1 - y_1, x_2 - y_2) + K(x_1 +y_1, x_2 -y_2) ) \pa_1 \chi_1 \phi_\dag(y) dy .
\]

Since $K(z)$ is odd in $z_1$ and $|y| \geq 2|x|$ for $y$ in the support of the integrand, we get 
\[
|  K(x_1 - y_1, x_2 - y_2) + K(x_1 +y_1, x_2 -y_2) | 
= |K( y_1 + x_1, x_2 - y_2 ) - K(y_1 - x_1, x_2 - y_2)| 
\les \f{2 x_1}{|y|^3}. 
\]

Using Lemma \ref{lem:far}, \eqref{eq:chi2}, and \eqref{eq:chi_Ri}, we get 
\[
|\pa_1 \chi_1 \phi_\dag | \les \f{1}{R_1}  S^{2-\b} \one_{ R_1 \leq | y| \leq 2R_1} M_1.
\]
It follows 
\[
|J_1| \les  \f{ |x| S^{2-\b}}{R_1} M_1 \int   \one_{ R_1 \leq | y| \leq 2R_1} |y|^{-3}  dy 
\les \f{ |x| S^{2-\b}}{ R_1^2} M_1.
\]

Next, we estimate $\na (-\D)^{-1} (\phi_\dag \D \chi_1)$ in $I_2$ in \eqref{eq:prop_C10}. Using a similar symmetrization argument, the fact that $x$ is away from the singularity of the kernel when $|x| \leq  2 R_2 < R_1$ (as the integrand is supported in $|y| \geq R_1  \geq 2 |x|$), and $ R_1 \asymp  C_l^{-1}$ \eqref{eq:chi_Ri}, we obtain 
\beq\label{eq:prop_C1_Z11}
 |\na (-\D)^{-1} Z_{1, 1}| \les \f{ |x| S^{2-\b}}{R_1^2 } M_1
 \les |x| C_l^2 S^{2-\b} M_1
 \les \min( |x|, |x|^{1-\b}) (C_l S)^{2-\b} M_1.
 \eeq

Combining the above estimate and using $C_l S \leq 1$, we obtain 
\[
|\na \phi_\dag | \les_{\al, \b} \min( |x|, |x|^{1-\b}) (  M_1 + (C_l  S )^{2-\b} B_1  ).
\]

Using \eqref{eq:prop_C10}, and taking the maximum of $x$ over $|x| \leq  2 R_2$, we prove \eqref{eq:prop_C11}, which further implies the desired result.

\vs{0.1in}
\paragraph{\bf{Improved estimate}}
For $ 8 S \leq |x|  \leq R_1 / 2$, we refine the estimates of $I_1, I_3$. In $I_1$, 
for $y \in \supp(\Om_\dag)$, we have $|x| \geq 2 |y|$. For $K(z) = \f{z_i}{|z|^2}$, using the same symmetrization argument, we get
\[
\bal
 |\int_{\R_2} K(x - y) \Om_\dag dy |
&  = \B| \int_{y_1 \geq 0}  ( K(x_1 - y_1, x_2 - y_2) - K(x_1 + y_1, x_2 - y_2) ) \Om_\dag d y \B| \\
&\les \int_{y_1 \geq 0} \f{2y_1}{|x|^2} |\Om_\dag(y) | dy 
\les |x|^{-2} M_1 \int_{|y|\leq S} |y| \min( |y|^{\al}, |y|^{-\b}) dy 
\les M_1 |x|^{-2} S^{3-\b} .
\eal
\]

We have estimated $ \na (-\D)^{-1} Z_{1,1}$ in \eqref{eq:prop_C1_Z11}. For $Z_{2, 1}$
and $I_3$ \eqref{eq:prop_C10}, we improve the estimate of $|| Z_{2, 1} ( |x|^{-1} + |x|^{\b}) ||_{\inf}$. For $|x| \geq 2S$, we have the estimate \eqref{eq:prop_C1_Z21}. For $|x| \leq 2S$, using the first estimate in Proposition \ref{prop:C1} we just proved, $S \geq 1$, and using $\na \phi_\dag$ to bound $\phi_\dag$, we obtain
\[
|Z_{2, 1}| (|x|^{-1} + |x|^{\b})
\les C_l^2  \min( |x|^2, |x|^{2-\b})(|x|^{-1} + |x|^{\b})  M_1
\les   C_l^2 ( |x| + |x|^2) M_1 \les (C_l S)^{2-\b} M_1. 
\]
Combining the estimates in two cases, we get
\[
|| Z_{2, 1}| (|x|^{-1} + |x|^{\b}) ||_{\infty} \les  ( C_l S)^{2-\b} M_1.
\]
Applying Lemma \ref{lem:E_ker} again, we obtain 
\[
|\na^{\perp}(-\D)^{-1} Z_{2,1}| \les |x|^{1-\b}   ( C_l S)^{2-\b} M_1.
\]
Combining the above estimates and \eqref{eq:prop_C1_Z11}, we prove \eqref{eq:prop_C1_impr}.
\end{proof}

\subsubsection{Estimate of $\na^2 \phi_\dag$}

Using Proposition \ref{prop:C1} for $\na \phi_\dag$, we further estimate $\na^2 \phi_\dag$ for $|x| \leq R_2/2$.

\begin{prop}\label{prop:C2}
Let $\phi_\dag$ be the solution in \eqref{eq:E_elli}, $ \b \in (0, 1)$, and $\al \in (0, 1]$. There exists some absolute constant $\nu_2(\al, \b) < \f{1}{32}$ such that if $C_l(\tau) (1 + S(\tau)) < \nu_2$, for $|x| \leq R_2 / 2$, we have 
\[
|\na^2  ( \phi_\dag \chi_2) -   \na^2 (-\D)^{-1} \Om_\dag| 
\les_{\al, \b}  C_l^{\b} || \Om_\dag  (|x|^{-\al} + |x|^{\b}) ||_{\infty}.
\]
In particular, we have 
\beq\label{eq:psixy0}
\pa_{xy}\phi_\dag (0) =\pa_{xy} \phi(0), \quad  | \pa_{xy} \phi_\dag(0) - \pa_{xy} (-\D)^{-1}\Om_\dag(0) | 
\les_{\al, \b}  C_l^{\b} || \Om_\dag  (|x|^{-\al} + |x|^{\b}) ||_{\infty}.
\eeq

\end{prop}

For the perturbation $\Om_\dag$, we further bound $\na^2(-\D)^{-1}\Om_\dag$ using the Boussinesq energy \eqref{energy4}.

\begin{proof}
We consider \eqref{eq:E_elli1} with $i = 2$. Denote 
\[
M_1 = || \Om_\dag  (|x|^{-\al} + |x|^{\b}) ||_{\infty}.
\]

Using \eqref{eq:E_elli1} with $i=2$ and $Z_{1, 2}, Z_{2, 2}$ defined therein, we have
\[
\na^2( \phi_\dag \chi_2) = \na^2(-\D)^{-1} \Om_\dag - \na^2 (-\D)^{-1} Z_{1, 2} - \na^2 (-\D)^{-1} Z_{2, 2 } \teq I_1 + I_2 + I_3,
\]
where we have used $\Om_\dag \chi_i = \Om_\dag$ by requiring $C_l S$ small. We only need to estimate $I_2, I_3$. The estimate of $I_2$ is similar to that in the proof of Proposition \ref{prop:C1}. We consider the typical term 
\[
J = \pa_{12}(-\D)^{-1} ( \pa_1 \phi_\dag \pa_1 \chi_2).
\]

For $|x| \leq R_2/2$, it is away from $\supp(\pa_1 \phi_\dag \pa_1 \chi_2) \subset \{ |x| \in [R_2, 2 R_2] \}$. We have 
\[
J(x) = C  \int K(x - y) ( \pa_1 \phi_\dag \pa_1  \chi_2)(y) dy, 
\quad  K(z) \teq \f{z_1 z_2}{|z|^4} .
\]

Using Proposition \ref{prop:C1}  and \eqref{eq:chi2}, we get 
\beq\label{eq:psi_Z11}
\bal
|\pa_1 \phi_\dag \pa_1  \chi_2| 
&\les R_2^{-1} \one_{ R_2 \leq |y| \leq 2R_2} |\pa_1 \phi_\dag  |
\les M_1 R_2^{-1}  R_2^{1-\b}\one_{ R_2 \leq |y| \leq 2R_2}  \les M_1  R_2^{-\b} \one_{ R_2 \leq |y| \leq 2R_2} .
\eal
\eeq

Since $|x| \leq  R_2 / 2 \leq |y| / 2$ and $|K(x-y)| \les |y-x|^{-2} \les|y|^{-2}$, we get 
\[
\bal
|J| &\les  M_1   R_2^{-\b}  \int_{ R_2 \leq |y| \leq 2 R_2} |K(x- y)| dy
\les  M_1   R_2^{-\b}  \int_{ R_2 \leq |y| \leq 2 R_2} |y|^{-2} dy  
\les M_1 R_2^{-\b}.
\eal
\]

Using similar estimates for other terms in $I_2$ and \eqref{eq:chi_Ri}, we obtain
\[
|I_2| \les M_1   R_2^{-\b} \les M_1  C_l^{\b} ,
\]

For $I_3$, we estimate $\pa_{12}(-\D)^{-1} Z_{2,2}$. Recall $Z_{2, 2} = \f{a C_l^2}{r^2} \phi_\dag \chi_2$ from \eqref{eq:E_elli1}. Other derivatives are similar. By definition, we have 
\[
 \pa_{12}(-\D)^{-1} Z_{2,2} = C \int K( x - y) Z_{2, 2}(y) dy.
\]

For a fixed $x$, we partition the region of the integral into three parts
\[
Q_1 = \{ y : |y| \geq 2|x|  \}, \quad Q_2 = \{ y: |y-x| \leq  |x| / 2 \},  \quad Q_3 = (Q_1 \cup Q_2)^c .
\]

Applying Proposition \ref{prop:C1}, using \eqref{eq:chi2} and $|\pa_i r^{-1}| \les C_l, r \asymp 1$ when $|r| \geq \f{1}{2}$, for $ |x| \leq 2 R_2$, we have $C_l |x| \leq 1$ and 
\[
\bal
 |\pa_i Z_{2, 2}| 
& \les C_l^2 ( |\pa_i r^{-1}  \phi_\dag \chi_2| + |\pa_i \phi_\dag \chi_2 | 
+ | \phi_\dag  \pa_i \chi_2 |  ) \\
& \les  M_1 C_l^2 \one_{|x| \leq 2 R_2} (   C_l  |x|^2 \wedge |x|^{2-\b}+ |x| \wedge |x|^{1-\b}
+ (|x|^2 \wedge |x|^{2-\b} ) R_2^{-1}  )  \\
& \les M_1 C_l^2 R_2^{1-\b}  \les C_l^{1+\b} M_1.
\eal
\]

We also have the pointwise estimate 
\beq\label{eq:psi_Z21}
|Z_{2, 2}| \les C_l^2 |\phi_\dag \chi_2| \les C_l^2 \min( |x|^2, |x|^{2-\b}) \one_{|x| \les R_2} M_1
\les M_1 C_l^{\b}.
\eeq

Since $Q_1 = \{ y : |y| \geq 2|x|  \}$, using the above pointwise estimate,  
we obtain for $|x| \leq R_2/2$:
\[
 |\int_{Q_1} K(x- y) Z_{2, 2}(y) dy| \les M_1 \int_{ 2 |x| \les |y|\les R_2 } |y|^{-2} C_l^2 
 \min( | y|^{2}, |y|^{2-\b}) dy
 \les C_l^2 R_2^{2-\b} M_1 \les M_1 C_l^{\b}.
 \]

Using $ P.V. \int_{|x-y| \leq |x|/2} K(x-y) d y = 0$ and the above pointwise estimate, we obtain 
\[
 |\int_{|x-y| \leq |x|/2} K(x- y) Z_{2, 2}(y) dy| 
\les  \max_{ |x-y | \leq |x|/2} |\na Z_{2, 2}(y)|
\int_{|x-y| \leq |x|/2} |x-y|^{-1} dy 
\les C_l^{1+\b} M_1 |x| \les M_1  C_l^{\b}.
\]

Since $Q_3 \subset \{ y: |x| / 2 \leq |x-y| \leq 3|x| \}$, 
using \eqref{eq:psi_Z21}, we get 
\[
 |\int_{ Q_3} K(x- y) Z_{2, 2}(y) dy| 
 \les |x|^{-2} |\int_{Q_3} Z_{2, 2}(y) dy|
 \les |x|^{-2} M_1 C_l^{\b} \int_{ |y| \leq 4|x|} dy \les M_1 C_l^{\b}.
\]
Combining the above estimates, we prove the desired result.

To obtain \eqref{eq:psixy0}, we use $r = 1$ when $x=y =0$ and 
$
\pa_{xy}\phi_\dag(0) = \pa_{xy} ( \phi r^{1/2})(0) = \pa_{xy} \phi(0).
$
\end{proof}

\subsubsection{Weighted $L^{\inf}$ and H\"older estimate}

  Recall the weights $\vp_1, \vp_{g,1}, \psi_1$ and $g_1$ from \eqref{wg:hol}, \eqref{wg:linf_decay}, \eqref{wg:linf_grow}. 
The energy \eqref{energy4} for 2D Boussinesq in Section \ref{sec:EE} controls the following norm of $\om_1$:
\beq\label{norm:X}
|| \Om ||_{X_1} \teq \max( || \Om \vp_1 ||_{L^{\inf}},  || \Om \vp_{g,1} ||_{L^{\inf}} , 
[ \Om \psi_1 ]_{C_{g_1}^{1/2}}  ). 
\eeq
Recall $\vp_{g, 1} \gtr c \la x \ra^{1/16}$ 
from \eqref{wg:linf_grow}. 
For applications in Propositions \ref{prop:cor1}, \ref{prop:elli}, we take 
$\b = \tf{1}{16}$ .

We aim to apply the weighted nonlocal estimates 
in 2D Boussinesq, e.g.  \eqref{eq:holest}, to estimate $\phi_\dag \chi_i$ with source $ \Om_\dag \chi_i - Z_{1, i} - Z_{2, i}$ in \eqref{eq:E_elli1}. However, since $Z_{2, 4}$ vanishes only as $O(|x|^2)$ near $x=0$, $Z_{2, 4} \vp_1 \notin L^{\infty}$ and $Z_{2, 4} \notin X$. Instead, we show that $Z_{2, 4}$, with a rank one correction near $x=0$, and $Z_{1, 4}$ are small in $X_1$ using Propositions \ref{prop:C1} and \ref{prop:C2}.  Below, the sizes of $Z_{1, i }, Z_{2, i}$ are very small. The reader can mainly focus on the vanishing order of these terms near $x=0$.

\begin{prop}\label{prop:cor1}
Let $\phi_\dag$ be the solution to \eqref{eq:E_elli}. Suppose that $ \Om_\dag \in X_1 $ and $C_l ( 1 + S) < \nu_2$, where $\nu_2$ is the constant in Proposition \ref{prop:C2}. For $|x| \leq 2, \al, \b  >0, \al < 2$,  we have
\footnote{
For $\al > 1$, there are no cubic order terms in the expansion of $ \phi_\dag -   x_1 x_2  \pa_{xy} \phi_\dag(0) $ for reasons in Footnote \ref{fn:phi_expand}.
}
\[
\bal
  | \na( \phi_\dag -   x_1 x_2  \pa_{xy} \phi_\dag(0)) | & \les_{\al, \b}   |x|^{1 + \al } || \Om_\dag ( |x|^{-\al} + |x|^{\b} ) ||_{L^{\inf}} , \\
  |  \phi_\dag -   x_1 x_2  \pa_{xy} \phi_\dag(0)  |  & \les_{\al, \b}   |x|^{2 + \al } || \Om_\dag  ( |x|^{-\al} + |x|^{\b} ) ||_{L^{\inf}} .
  \eal
\]
\end{prop}

\begin{proof}
Using \eqref{eq:E_elli1} with $i=4$, we get 
\[
\phi_\dag \chi_4 = (-\D)^{-1} ( \Om_\dag -  Z_{1, 4} - Z_{2, 4}) . 
\]
Using Lemma \ref{lem:E_ker}, we only need to prove that 
\[
||  (\Om_\dag -  Z_{1, 4} - Z_{2, 4} )( |x|^{-\al} + |x|^{\b}) ||_{\inf}
\les ||  \Om_\dag  ( |x|^{-\al} + |x|^{\b}) ||_{\inf}.
\]
Since $\chi_1 = 1$ in the support of $\chi_4$, the estimates of $Z_{1, 4}, Z_{2, 4}$ follows from Proposition \ref{prop:C1} and its proof. We only consider a typical term. For $\pa_1 \phi_\dag \pa_1 \chi_4$ in $Z_{1, 4}$, using Proposition \ref{prop:C1}, we get 
\[
\bal
| \pa_1 \phi_\dag \pa_1 \chi_4| (|x|^{-\al} + |x|^{\b})
&\les ||  \Om_\dag  ( |x|^{-\al} + |x|^{\b}) ||_{\inf} \min( |x|, |x|^{1-\b}) R_4^{-1} \one_{ |x| \asymp R_4} (|x|^{-\al} + |x|^{\b})  \\
& \les  ||  \Om_\dag  ( |x|^{-\al} + |x|^{\b}) ||_{\inf} .
\eal
\]
We require $\al < 2$ since $Z_{2, i}$ in \eqref{eq:E_elli1} only vanishes to order $|x|^2$ near $x = 0$.
\end{proof}

We are in a position to show that the $Z$ term in \eqref{eq:E_elli1} with a correction is small in space $X_1$.

\begin{prop}\label{prop:elli}
Suppose that $ \Om_\dag \in X_1 $ and $C_l (1 + S) < \nu_2$, where $\nu_2$ is the constant in Proposition \ref{prop:C2}. 
For $Z_{1, 4}, Z_{2, 4}$ defined in \eqref{eq:E_elli1}, 
 a rank-one correction of $Z_{1, 4} + Z_{2, 4}$ satisfies 
 \beq\label{eq:prop_elli}
 \bga
  || Z ||_{X_1}  \les C_l S M , \quad || \na Z ||_{L^{\inf}} \les C_l M,  \quad 
Z \teq Z_{1 , 4} + Z_{2, 4} - \pa_{xy} Z_{2, 4}(0) (-\D \kp) , \\
  \ega
\eeq
where $\kp = -\f{x y^3}{ 6  } \chi, \chi$ is some cutoff function supported near $x=0$, e.g. \eqref{eq:chi2}, and $M$ denotes 
\beq
   M_{\eqref{def:M_Z}} \teq || \Om_\dag ( |x|^{-1} + |x|^{\b} ) ||_{\inf} + || \chi_3 \na^2 (-\D)^{-1} \Om_\dag||_{\infty}. 
   \label{def:M_Z}
\eeq
\end{prop}

We will apply Proposition \ref{prop:elli} to $\Om_\dag = \om r^{3/2}$ with $\om \in X_1$ or $\om = \bar \om \chi(x / \nu)$ with $1< \nu < R_4$, where $\bar \om$ is the approximate steady state for the 2D Boussinesq equations. In both cases, we can further bound the right hand side as follows
\begin{align}
& M_{ \eqref{def:M_Z} }
  \les_{\al} || \Om_\dag ( |x|^{-1} + |x|^{\b} ) ||_{\inf} 
 + || \la x \ra^{-p} \Om_\dag||_{C^{\al}}, \quad \forall  \ \al, \b > 0 , p \geq 0  \notag \\ 
& ||  \bar \om \chi_{\nu} r^{3/2} ( |x|^{-1} + |x|^{\b} ) ||_{\inf} + || \chi_3 \na^2 (-\D)^{-1} (\bar \om \chi_{\nu})||_{\infty}  \les 1.
\label{prop:elli_cor}
\end{align}
Estimates of $\na^2(-\D)^{-1} f, f = \Om_\dag, \bar \om \chi_{\nu}$ in \eqref{prop:elli_cor} follow from  estimate below \eqref{eq:u_est_direct}  and \eqref{solu:decay}.

\begin{proof}

Let $\phi_\dag$ be the solution of \eqref{eq:E_elli1}. Now, using Proposition \ref{prop:C2}, for $|x| \leq R_2/2$, we yield 
\[
 |\na^2 \phi_\dag(x)| \les M.
\]
where $M$ is defined in \eqref{eq:prop_elli}, and the norm of $\Om_\dag$ in Propositions \ref{prop:C1}, \ref{prop:C2} with $\b = \f{1}{16}$  can be bounded by $M$.

\vs{0.1in}
\paragraph{\bf{Estimate of $Z_{1, 4}$}}

Firstly, we estimate $Z_{1, 4}$ \eqref{eq:E_elli1}. Recall $Z_{1, 4}$ from \eqref{eq:E_elli1} 
\[
  Z_{1,4} = 2 \na \phi_\dag  \cdot  \na \chi_4 + \phi_\dag \D \chi_4 .
\]
From the definition of $\chi_4$ \eqref{eq:chi_Ri}, for $x \in \supp( Z_{1, 4} )$, we have $|x| \asymp R_4$. 
From \eqref{eq:chi_Ri}, we have $C_l R_4 \asymp 1$. 
Since $\phi_\dag = \phi r^{1/2}$, using Lemma \ref{lem:far} and \eqref{eq:prop_C1_impr} in Proposition \ref{prop:C1} for $\na \phi_\dag$, we obtain 
\[
| \phi_\dag(x)|  \les S^{2-\b} M, \quad  | \na \phi_\dag(x) | \les |x|^{1-\b}  M
(C_l S)^{2-\b} ,
\]
where we have used $S \les |x| , C_l |x| \asymp  C_l R_4 \asymp 1 $ to simplify the upper bound in \eqref{eq:prop_C1_impr}. Using the above estimate and $R_4^{-1} \asymp C_l$, we obtain the pointwise estimate 
\[
\bal
|Z_{1, 4}|  & \les  \one_{ |x| \asymp R_4 } ( R_4^{-1}  |\na \phi_\dag|  + R_4^{-2} |  \phi_\dag| )
\les M  \one_{ |x| \asymp R_4  }  ( R_4^{-1} |x|^{1-\b}  (C_l S)^{2-\b}
+ S^{2-\b} R_4^{-2}  ) \\
& \les  C_l^{\b} ( C_l S )^{2-\b} M  \one_{ |x| \asymp R_4 . }
\eal
\]

Using \eqref{eq:chi2} for $\na^k \chi_4$, Lemma \ref{lem:far}, and Propositions \ref{prop:C1}, \ref{prop:C2} for $ \na^k \phi_\dag$, we obtain 
\[
\bal
|\na Z_{1, 4}| & \les \one_{ |x| \asymp R_4 } ( 
|\na^2 \phi_\dag| \cdot  |\na \chi_4| 
+ | \na \phi_\dag | \cdot |\na^2 \chi_4| 
+ |\phi_\dag| \cdot |\na^3 \chi_4| ) \\
& \les \one_{ |x| \asymp R_4 } M
( R_4^{-3} |x|^{2-\b} + R_4^{-2} |x|^{1-\b}  + R_4^{-1} )
\les \one_{ |x| \asymp R_4 } M  R_4^{-1} 
\les \one_{ |x| \asymp R_4 } M  C_l.
\eal
\]

Note that the weights $\vp_1, \vp_{g, 1}$ \eqref{wg:linf_decay}, \eqref{wg:linf_grow} involve $|x_1|^{-1/2}$ which is singular along $x_1=0$. For any power $\g \in [-3, \b]$, we have 
\[
| Z_{1, 4} |x|^{\g} | \les C_l^{\b} |x|^{\g} (C_l S)^{2-\b} \one_{|x| \asymp R_4} M 
\les (C_l S)^{2-\b} \one_{|x| \asymp R_4} M .
\]

For any power $\g \in [- 5/2, 0)$ and $x_1 \leq 1$, since $ x \in \supp(Z_{1, 4})$ implies 
$|x| \les R_4$, we have 
\[
 | Z_{1, 4} |x|^{\g} |x_1|^{-1/2} |
 \les  \max_{ |z| \leq 1 } | \na Z_{1, 4}(z, x_2)| \cdot |x|^{\g}
 \les R_4^{-1 + \g} M
 \les C_l  M .
\]

To estimate the H\"older norm of $Z_{1, 4} \psi_1$, following similar estimates, we obtain 
\[
 |Z_{1,4} \psi_1| \les ( C_l S)^{2-\b}  M, \quad  |\na ( Z_{1,4} \psi_1)| \les   (C_l + (C_l S)^{2-\b} ) M.
\]

We can obtain better estimates due to the decay of $\psi_1(x)$ for large $x$, see \eqref{wg:hol}. But we do not need this extra smallness. Using the above estimates and embedding inequalities, we obtain 
\[
| Z_{1, 4} \vp_1| \les | Z_{1,4} \vp_{g,1} | \les  ( C_l + (C_l S)^{2-\b} ) M, 
\quad [ Z_{1,4} \psi_1 ]_{C_{g_1}^{1/2}} 
\les ( C_l + (C_l S)^{2-\b} ) M ,
\]
which along with the above estimate for $\na Z_{1, 4}$ establish the estimates for $Z_{1,4}$ in \eqref{eq:prop_elli}.

\vs{0.1in}
\paragraph{\bf{Estimate of $Z_{2, 4}$}}
Recall 
\[
Z_{2, 4} = \f{a C_l^2}{r^2} \phi_\dag \chi_4. 
\]

Clearly, we have $\pa_{xy} Z_{2,4}(0) = a C_l^2 \pa_{xy} \phi_\dag(0)$. We perform the following decomposition 
\[
Z_{2,4} - \pa_{xy} Z_{2, 4}(0) ( -\D \kp) 
= \f{a C_l^2}{r^2} (\phi_\dag \chi_4 - \pa_{xy} \phi_\dag(0)  (-\D \kp )  )
+ \pa_{xy} \phi_\dag(0) \f{a C_l^2}{r^2 } (1- r^2) (-\D \kp) \teq I + II.
\]

\paragraph{\bf Estimate $II$ in $Z_{2, 4}$}
From the definition of $\kp$ in Proposition \ref{prop:elli}, for $|x| \leq 1$, we have 
\[
-\D \kp = x_1 x_2, \quad 1 - r = C_l x_2 .
\]

Thus $II = O( x_1 x_2^2)$ near $x=0$. Using Proposition \ref{prop:C1} for 
\[
|\na^2 \phi_\dag(0) | \les \limsup_{x\to 0} | \na \phi_\dag |x|^{-1} |, 
\]
Proposition \ref{prop:C2} for $ \pa_{xy} \phi(0)$ and the fact that $II$ is supported near $x=0$, we get 
\[
|| II||_{X_1} \les C_l^3 M , \quad |\na II | \les C_l^3 M.
\]

\paragraph{\bf Estimate $I$ in $Z_{2, 4}$}
Applying Propositions \ref{prop:C1} and \ref{prop:cor1} with $\al = 1$,  we obtain 
\beq\label{eq:prop_elli_I1}
| \na I| \les C_l^2 \min( |x|^2, |x|^{1-\b} ) M,  \quad |I| \les C_l^2 x_1 \min( |x|^2, |x|^{1-\b} ) M,
\eeq
where to obtain the second bound $C_l^2 x_1 |x|^{1-\b}$, we have used Proposition \ref{prop:C1} and integrated the estimate for $\pa_1 \phi$  in $x_1$ to estimate $\phi_\dag$. Note that if the derivative acts on $r^{-1}$, we get $ |\na r^{-1}| \les C_l$ and then use $C_l |x| \les 1$ to remove a growing power $|x|$. Since $C_l |x| \les 1$ in the support of $I$, combining the estimate of $\na I, \na II$, we obtain the estimate of $\na (I + II)$ in \eqref{eq:prop_elli}.

For large $|x| \geq 8 S$, the correction vanishes $(-\D \kp) = 0 $. Using Lemma \ref{lem:far} and the improved estimate \eqref{eq:prop_C1_impr}, for $ 8 S \leq |x| \leq R_1 / 2 \les C_l^{-1}$, we have $|x|^{1-\b} (C_l S)^{2-\b} \les C_l S^{2-\b} \les S^{1-\b}$ and 
\[
\bal
|\phi_\dag |  & \les (1 +| \log( C_l S) | ) S^{2-\b}  M, \quad  |\na( \phi_\dag / r)| 
\les |\na \phi_\dag| + C_l |\phi_\dag| \les S^{1-\b}  M,  \\
\eal
\]
where we have used $ C_l S^{2-\b} (1 + | \log |C_l S| | ) \les S^{1-\b}$ to absorb the logarithm factor. It follows 
\beq\label{eq:prop_elli_I13}
\bal
  | I|  & \les  C_l^2 (1 +| \log( C_l S) | ) S^{2-\b} M
 \les  (C_l S)^{1-\b} (1 + |\log (C_l S )|)  C_l^{1 + \b} S M
 \les C_l^{1 + \b} S  M , \\
| \na I| & \les C_l^2  S^{1-\b}  M \les C_l M,
 \eal
\eeq
for $|x| \in [8 S, R_1/ 2]$. Therefore, for any $\g_1 \in [-3, \b]$, combining the above estimates and \eqref{eq:prop_elli_I1} and using $|x| \leq R_1 \les C_l^{-1}$ within the support of $I$, we have 
\beq\label{eq:prop_elli_I12}
 |x|^{\g_1} |I| \les C_l S  M.
\eeq

Next, we bound $|I| |x|^{\g_2} |x_1|^{-1/2}$ for $\g_2 \in [-5/2, 0]$. If $|x_1| \geq 1$, it follows from the above bound. If $|x_1| \leq 1$, using \eqref{eq:prop_elli_I1}, we obtain 
\[
|x|^{\g_2} |x_1|^{-1/2} |I|
 \les   |x|^{\g_2} C_l^2 |x_1|^{1/2} \min( |x|^2, |x|^{1-\b}) M
 \les C_l^2 \min( 1,  |x|^{1-\b} )M \les C_l^{1+\b} M.
\]

The above estimates imply 
\[
 || I \vp_1||_{\inf} \les || I \vp_{g, 1}||_{\inf} \les (C_l + C_l S)  M
\les C_l S || \Om_\dag||_{X_1}, \quad || I \psi_1 ||_{\inf} \les C_l S  M.
\]
Since $\psi_1 \les  |x|^{-2} + 1, |\na \psi_1| \les |x|^{-3} + |x|^{-1}$, for $I$, using \eqref{eq:prop_elli_I1},\eqref{eq:prop_elli_I12}, we get 
\[
\bal
| \na ( I \psi_1) | & \les |\na I| \psi_1 + |I \na \psi_1| 
\les C_l^2 \min( |x|^2, |x|^{1-\b}) ( |x|^{-2} + 1)  || \Om_\dag||_{X_1}+ C_l  S M \\
& \les  ( C_l^2 (1 + |x|^{1-\b} ) + C_l S ) M
\les  C_l S  M, \\
 |I \psi_1 | & \les C_l S M.
\eal
\]
In the last inequality for $ \na ( I \psi_1) $, we have used $ C_l^{1-\b} |x|^{1-\b} \les 1$ within the support of $I$ and $C_l \ll  1< S$. Using embedding inequalities, we obtain the H\"older estimate of $ I \psi_1$. We conclude 
\[
|| Z_{2, 4}  - \pa_{xy}Z_{2,4}(0) (-\D) \kp ||_{X_1} \les || I ||_{X_1} + || II||_{X_1} \les C_l S M.
\]
Since $C_l S \les 1, C_l \les 1$, combining the estimate of $Z_{1, 4}, Z_{2, 4}$, we complete the proof.
\end{proof}

\subsection{Main terms for the velocity}\label{sec:main_vel_euler}

In this section, we rewrite the elliptic equation in a form whose leading part is the
standard 2D Poisson equation, and apply estimates in Section \ref{sec:3Delli} to 
bound the $C_l$-dependent terms perturbatively. 
Below, we treat $C_l S, C_l$ as small parameters.

Based on Proposition \ref{prop:elli}, we rewrite \eqref{eq:E_elli1} with $i = 4$ as: 
\[
- \D( \phi_\dag \chi_4 + a C_l^2 \phi_{xy} (0) \kp  )  = \Om_\dag \chi_4 - Z_{1, 4} - ( Z_{2, 4} - \pa_{x y} Z_{2, 4}(0) (-\D \kp) ), 
\]
where we have used \eqref{eq:psixy0} and 
\[
\pa_{xy} Z_{2,4}(0) = a C_l^2 \pa_{xy} \phi_\dag(0) = a C_l^2 \pa_{xy} \phi(0). 
\]

Recall the definitions of $\phi_\dag, \Om_\dag$ from \eqref{eq:E_elli1}, \eqref{eq:E_elli}, and $\kp$ from Proposition \ref{prop:elli}. We introduce 
\beq\label{eq:elli_main}
\bal
\Psi &= \phi r^{1/2} \chi_4 + a C_l^2  \phi_{xy}(0) \kp , \quad 
\Psi_2  = -a C_l^2 \phi_{xy}(0) \kp , \\
 \Om &  = \om r^{3/2} \chi_4 - Z_{1, 4} - ( Z_{2, 4} - \pa_{x y} Z_{2, 4}(0) (-\D \kp) )  . \\
\eal
\eeq
Then we obtain 
\beq\label{eq:elli_main2}
- \D \Psi = \Om , \quad  \phi r^{1/2} \chi_4= \Psi + \Psi_2.
\eeq

Recall $r = 1 - C_l x_2$ and $\Om_\dag = \om r^{3/2}$ from \eqref{eq:E_elli}. 
By \eqref{eq:chi_Ri} and \eqref{eq:E_supp}, within the support of $\om$,   we have $\om \chi_4 = \om$,  $r \asymp 1, r \in C^{\infty}$ and $|r-1| \les C_l |x_2| \leq C_l S$. Using Proposition \ref{prop:elli}, \eqref{prop:elli_cor}, 
and simple $C^{1/2}$ product estimate,
for $i=0, 1, 2$,
we have
\bseq\label{wg:equi}
\begin{align}\label{wg:equi:a}
 & || \f{|x|^i}{ \la x \ra^i}  ( \Om - \om ) ||_{X_1}  \les  \| \f{|x|^i}{ \la x \ra^i}  \om (r^{3/2} - 1) \|_{X_1}
 +  C_l S ( \nlinf{ ( |x|^{-1} + |x|^{\f{1}{16}} ) \om } + [ \la x \ra^{- \f{1}{6} } \om ]_{C^{1/2}} ) \\
& \les 
C_l  \| \f{|x|^i}{ \la x \ra^i}   x_2 \om \|_{X_1} 
+  C_l S ( \nlinf{ ( |x|^{-1} + |x|^{ \f{1}{16} } ) \om } + [ \la x \ra^{- \f16} \om ]_{C^{1/2}} ) 
\les C_l S  (  \|( \f{|x|}{ \la x \ra} )^{\min(i,2)}  \om \|_{X_1}
+  \nlinf { |x|^{-1} \om } ). \notag
\end{align}
We require $i\leq 2$ so that $ \|( \f{|x|}{ \la x \ra} )^{\min(i,2)}  \om \|_{X_1}$ is stronger 
than the weighted $L^{\infty}$-norm in \eqref{wg:equi:a}.  For $\al \in (0,1]$, using Proposition \ref{prop:elli}, \eqref{prop:elli_cor}, we bound 
\[
||   \Om - \om  ||_{C^{0, \al}} \les 
\|   \om (r^{3/2} - 1) \|_{C^{0,\al}} 
+ C_l S ( \| \om \|_{C^{0,\al}}
+ \nlinf{ ( |x|^{-1} + |x|^{\f{1}{16}} ) \om } ) .
\]

In particular, for $i=0,1,2$, 
and $\al \in (0, 1]$, we obtain 
\beq
\|   \tf{|x|^i}{ \la x \ra^i} \Om \|_{X_1} \les  \|   \tf{|x|^i}{ \la x \ra^i} \om \|_{X_1} 
+ \nlinf { |x|^{-1} \om },
\quad  \| \Om \|_{C^{0,\al}}  \les \| \om \|_{C^{0,\al}}
+ \nlinf{ ( |x|^{-1} + |x|^{\f{1}{16}} ) \om } .
\eeq
\eseq

Thus, $\Om$ and $\om$ enjoy almost the same estimates up to an error of size
$C_lS$.

From Propositions \ref{prop:C2}, \ref{prop:elli}, 
using \eqref{wg:equi}, $\phi_{xy}(0)$ satisfies 
\beq\label{eq:elli_main3}
|\phi_{xy}(0)| \les || \Om_\dag||_{X_1} \les || \Om||_{X_1}
\les \| \om \|_{X_1},   \quad \Psi_{xy}(0) = \phi_{xy}(0).
\eeq

Therefore, the terms $a C_l^2 \phi_{xy}(0) \kp, \Psi_2$ are very small and vanish to order $|x|^4$ near $x=0$.

\subsubsection{Main terms for the velocity}

Recall $u, v$ from \eqref{eq:euler4}. Since we only need the estimate of the velocity within the support of the solution, where $\chi_i = 1$, in the following derivation, we drop the cutoff functions $\chi_i$ to simplify the notation. Firstly, from \eqref{eq:elli_main2}, we yield 
\[
 \phi r^{1/2} = \Psi + \Psi_2, 
 \quad \Psi_2  = -a C_l^2 \phi_{xy}(0) \kp ,\ \phi_{xy}(0) = \Psi_{xy}(0). 
\]

The term $\Psi_2$ is smooth with vanishing order $|x|^4$, compactly supported, and small. We treat it as a lower order term and do not expand its derivation below. The velocity depends on the derivatives of $\phi$. Using $\pa_x r = 0, \pa_y r = - C_l$ defined in \eqref{eq:label_rs} 
(here $r=1-C_l y$ as in \eqref{eq:label_rs}; it is \emph{not} the radius $r = \sqrt{x^2+y^2}$),  
we rewrite $\na \phi$ as follows 
\beq\label{eq:E_psi}
\bal
\phi_y & = ( r^{1/2} \phi r^{-1/2} )_y 
= (\Psi r^{-1/2})_y + (\Psi_2 r^{-1/2})_y 
= \Psi_y r^{-1/2} + \f{ C_l}{2 r^{3/2}} \Psi +  (\Psi_2 r^{-1/2})_y  ,  \\
\phi_x & = r^{-1/2} \Psi_x + r^{-1/2} \Psi_{2, x}.
\eal
\eeq

Then using \eqref{eq:euler4}, we can rewrite $u, v$ as follows 
\beq\label{eq:E_u}
u  = -\phi_y + \f{1}{r} C_l \phi 
= - \Psi_y r^{-1/2} -  \f{ C_l}{2 r^{3/2}} \Psi + \f{1}{r^{3/2}} C_l \Psi
 - (\Psi_2 r^{-1/2})_y + \f{1}{r^{3/2}} C_l \Psi_2  \teq u_M + u_R,
\eeq
where the main term $u_M$ and the remainder $u_R$ are given by 
\beq\label{eq:E_uR}
u_M = - \Psi_y, \quad  u_R = - \Psi_y( r^{-1/2} - 1) + \f{ C_l}{2 r^{3/2}} \Psi  - (\Psi_2 r^{-1/2})_y + \f{1}{r^{3/2}} C_l \Psi_2. \\
\eeq

Note that the first and the second terms in $u_R$ cancel each other near the origin. Since
\beq\label{eq:E_r1}
r^{-1/2} - 1 = \f{1 - r^{1/2}}{r^{1/2}} = \f{1 - r}{ r^{1/2} (1 + r^{1/2})  } 
= \f{ C_l y}{ r^{1/2} (1 + r^{1/2}) },
\eeq
we obtain
\begin{align}\label{eq:E_uR2}
&- \Psi_y (r^{-1/2} - 1) + \f{ C_l}{2 r^{3/2}} \Psi  \\
 = & - (\Psi_y -  \Psi_{xy}(0) x ) (r^{-1/2 } - 1) 
+ \f{C_l}{2 r^{3/2}} ( \Psi - \Psi_{xy}(0) x y)  
 +  \Psi_{xy}(0)  ( - x (r^{-1/2} - 1) + \f{C_l x y}{ 2 r^{3/2} }   ) , \notag  \\
 = & - (\Psi_y -  \Psi_{xy}(0) x ) (r^{-1/2 } - 1) 
+ \f{C_l}{2 r^{3/2}} ( \Psi - \Psi_{xy}(0) x y)   + \Psi_{xy}(0) \f{ C_l x y( 1 + r^{1/2 } -2r )}{ 2 r^{3/2} (1 + r^{1/2})}. \notag 
\end{align}
The last term vanishes to order $O(|x|^3)$ near $x=0$. We treat $u_R$ as the remainder since 
it vanishes to order $O(|x|^3)$ near $x=0$ and contain the small factor $C_l$. Within the support of the solution, we get $C_l |x| \leq C_l S$, which is small.

For $v$ in \eqref{eq:euler4}, using \eqref{eq:E_psi} we have 
\beq\label{eq:E_v}
\bal
& v = \phi_x =  \Psi_x  + ( r^{-1/2} - 1) \Psi_x + r^{-1/2} \Psi_{2, x}   
\teq v_M + v_R, \\
& v_M = \Psi_x, \quad v_R = ( r^{-1/2} - 1) \Psi_x + r^{-1/2} \Psi_{2, x}   .
\eal
\eeq 
Within the support, since $v_R$ has a small size of order $C_l S$, we treat $v_R$ as a remainder term. 

\begin{remark}[\bf Vanishing order of $v_R$]
Near $x =0$, the vanishing order of $v_R$ is $O(|x|^2)$, which is lower than that of $u_R$ ($O(|x|^3)$). 
However, in \eqref{eq:euler4}, the coefficients of $v$, e.g. $\th_y, \om_y$, have higher vanishing order near $x=0$ than those of $u$, e.g. $\th_x, \om_x$. The remainder $u_R,v_R$, together with their coefficients, have enough vanishing order near $0$ for our weighted estimates.

\end{remark}

Using \eqref{eq:E_u},\eqref{eq:E_v}, we can compare 3D velocity $\uu$, expressed in terms of $\Psi(\om)$ \eqref{eq:elli_main},\eqref{eq:elli_main2}, with 2D velocity $\uu_{2D} = \na^{\perp} (-\D)^{-1} \om$. We refer to Section \ref{sec:E_lin_M} to bound this error.

\subsubsection{Main terms for the velocity of the approximate steady state}\label{sec:E_ub}

  We will construct an approximate steady state with bounded support in Section \ref{subsec:non_appr}
  \beq\label{def:om_nu1}
\bar \om_{\nu} \teq \chi_{\nu}(x) \bar \om, 
\quad 
\chi_{\nu}(x) \teq \chi_\dap( x / \nu), 
\eeq
where $\chi$ is the cutoff function chosen in \eqref{eq:profile2}. 
For $\nu$ sufficiently large, we will show that the associated velocity \eqref{eq:euler4} is close to that in the 2D Boussinesq equation.

To avoid confusion, we use subscript $2D$ to indicate the variables 
in 2D Boussinesq and  denote 
\[
\bar \phi_{2D} = (-\D)^{-1} \bar \om.
\]
Using \eqref{eq:elli_main}, \eqref{eq:elli_main2} with $\om = \bar \om_{\nu}$ in \eqref{def:om_nu1}, and subtracting \eqref{eq:elli_main} by $-\D \bar \phi_{2D} = \bar \om$, we yield
\[
(-\D) (\phi r^{1/2} \chi_4 + a C_l^2  \phi_{xy}(0) \kp - \bar \phi_{2D} )
 = \bar \om  \chi_{\nu} r^{3/2} \chi_4 - \bar \om - Z_{1, 4} - ( Z_{2, 4} - \pa_{xy} Z_{2, 4}(0) (-\D) \kp).
\]
Applying Propositions \ref{prop:elli}, \eqref{prop:elli_cor} with $\Om_\dag =r^{3/2} \bar \om_{\nu}$, 
we estimate the $Z$ terms as 
\[
|| Z_{1, 4} + ( Z_{2, 4} - \pa_{xy} Z_{2, 4}(0) (-\D) \kp)||_{X_1} \les C_l S.
\]

Note that the source term of the elliptic equation only vanishes to order $|x|^2$ near $x=0$:
\[
\bar \om \chi_{\nu} r^{3/2} \chi_4 - \bar \om
 = \bar \om ( r^{3/2} - 1 ) = -\tf{3}{2} C_l \bar \om_x(0) x y + O(|x|^3) . 
\]

We remove the quadratic remainder by adding a correction $ \f{3}{2} C_l \bar \om_x(0) (-\D \kp)$ with $-\D \kp = xy + l.o.t.$ to the above elliptic equation 
\[
\bal
& (-\D) \big( \phi r^{1/2} \chi_4 + (a C_l^2 \phi_{xy}(0)  + \bar \om_x(0) \tf{3}{2} C_l  ) \kp 
- \bar \phi_{2D}
 \big)  = \Om_{\nu, R},
\eal
\]
where 
\beq\label{eq:elli_mainb}
\bal
\Om_{\nu, R} &\teq \bar \om \chi_{\nu} \chi_4 (r^{3/2}-1) + \tf{3}{2} C_l \bar \om_x(0) (-\D \kp)
+ \bar \om( \chi_4 \chi_{\nu} - 1)  \\
 &\quad - Z_{1, 4} - ( Z_{2, 4} - \pa_{xy} Z_{2, 4}(0) (-\D) \kp), \\
 \bar \Psi  & = \phi r^{1/2} \chi_4 + (a C_l^2 \phi_{xy}(0)  + \bar \om_x(0) \tf{3}{2} C_l  ) \kp , \quad 
\bar \Psi_2 = - (a C_l^2 \phi_{xy}(0)  + \bar \om_x(0) \tf{3}{2} C_l  ) \kp .
 \eal
\eeq
Then the source term $\Om_{\nu, R}$ vanishes near $x=0$ to the order $|x|^3$.  
\footnote{
From the definitions of $r = 1 - C_l y$ in \eqref{eq:euler_polar} and $\kp$ in Proposition \ref{prop:elli}, 
near $x= 0$, we have $ \chi_4 = \chi_{\nu}  = 1$ and $- \D \kp =  x y $. 
Near $x=0$, it follows $\bar \om (\chi_4 \chi_{\nu} - 1) = 0$ and 
\[
\bar \om \chi_{\nu} \chi_4 (r^{3/2} - 1)
+ \tf{3}{2} C_l \bar \om_x(0) (-\D \kp )
= -\tf{3}{2} C_l y \cdot \bar \om_x(0) x + 
\tf{3}{2} C_l \bar \om_x(0) \cdot x  y +  O(|x|^3) = O(|x|^3)  .
\] 
}
We yield 
\beq\label{eq:elli_mainb2}
- \D(  \bar \Psi - \bar \phi_{2D}) = \Om_{\nu, R}, \quad -\D \bar \Psi = \Om_{\nu, R} + \bar \om.
\eeq

Since $\bar \om \in C^1$ and $ |\bar \om| \les \min( |x_1| , |x|^{-1/6}), |\na \bar \om| \les \min(1, |x|^{-7/6})$, using Proposition \ref{prop:elli} and \eqref{prop:elli_cor}, we obtain 
\beq\label{eq:elli_mainb_err}
|| \Om_{\nu, R} ||_{X_1} \les || \bar \om ( \chi_4 \chi_{\nu} - 1) ||_{X_1} + C_l S,
\quad ||  \Om_{\nu, R} ||_{ C^1} \les ||  \bar \om ( \chi_4 \chi_{\nu} - 1) ||_{C^1} 
+ C_l S.
\eeq

Since $\bar \om$ decays for large $|x|$ and $ | \chi_4 \chi_{\nu} - 1| \leq \one_{|x| \geq \nu}$, 
by choosing $ \nu$ sufficiently large and $ C_l S$ sufficiently small, we can ensure that $\bar \Psi - \bar \phi_{2D}$ is very small. See Remarks \ref{rem:order1}, \ref{rem:order2} for the size of these parameters. Similar to \eqref{eq:E_uR} and \eqref{eq:E_v}, based on $\bar \Psi, \bar \Psi_2$ in \eqref{eq:elli_mainb} and 
\[
\phi r^{1/2} = \bar \Psi + \bar \Psi_2,
\]
we decompose the velocity in \eqref{eq:euler4} associated with $\bar \om_{\nu}$ as follows 
\beq\label{eq:E_ub}
\bar u = - \bar \Psi_y + \bar u_R, \quad \bar v = \bar \Psi_x + \bar v_R,
\eeq
for $|x| \leq R_4$. The formulas of $\bar u_R, \bar v_R$ are similar to those in \eqref{eq:E_uR}, \eqref{eq:E_v} with $\Psi,  \Psi_2$ replaced by $\bar \Psi, \bar \Psi_2$. The remaining terms vanish near $0$ with order $\bar u_R = O(|x|^3), \bar v_R = O(|x|^2)$.

Since $ || f||_{L^p} \les || f \vp_{g, 1} ||_{\inf} \les || f||_{X_1}$ for $p>100$, 
using embedding \eqref{eq:elli_embed} established below and \eqref{eq:elli_mainb2}, \eqref{eq:elli_mainb_err}, we get 
\beq\label{eq:elli_mainb_err2}
\bal
 &  || \na^2 (\bar \Psi - \bar \phi_{2D}) ||_{L^{\inf}} + || \na^2 (\bar \Psi - \bar \phi_{2D} )||_{C^{1/2}}
 \les || \Om_{\nu, R}||_{X_1} + || \Om_{\nu, R}||_{C^1} \\
 & \les ||  \bar \om ( \chi_4 \chi_{\nu} - 1) ||_{C^1}  + ||  \bar \om ( \chi_4 \chi_{\nu} - 1) ||_{X_1} + C_l S  .
\eal
\eeq

For the remaining terms $\bar \uu_R = (\bar u_R, \bar v_R)$, for $|x| \leq 2S$, using $|r-1| , C_l |x|\les C_l S$, the elliptic equation for $\bar \Psi$ \eqref{eq:elli_mainb2}, and embeding inequalities, we get 
\beq\label{eq:elli_mainb_err3}
|| \na^2 \bar \Psi ||_{L^{\inf}}
+ || \na^2 \bar \Psi ||_{C^{1/2}} \les 1, \ 
 \sup_{|x|\leq 2S}| \na \bar \uu_R(x) |   + 
  \sup_{ |x|, |z| \leq 2S} \f{| \na \bar \uu_R(x) - \na \bar \uu_R(z) |}{|x-z|^{1/2}} \les C_l S .
\eeq

\subsubsection{Estimate of the 2D velocity}

We need several weighted estimates of  $(-\D)^{-1} \Om$ for the main terms in the velocity. 
The implicit constants will be absorbed by the small parameters $\nu^{-1}, C_l(0)$.

\begin{lem}\label{lem:E_vel}
Suppose that $\Om \in X_1$ \eqref{norm:X} is odd and $-\D \Psi = \Om$. We have 
\beq\label{eq:E_vel1}
\B| \na^2 \B( \Psi - \Psi_{xy}(0) xy  -  \f{1}{6} \pa_{1112} \Psi(0)  (x_1^3 x_2 -x_1 x_2^3)   \B)   \B|
\les |x|^{5/2} || \Om||_{X_1} , 
\eeq
for $|x| \leq 1$ and 
\[
\bal
\B| \na^2( \Psi - \Psi_{xy}(0) xy   ) \B| \les \min( |x|^{2}, 1) || \Om||_{X_1} ,
\quad  |\Psi_{xy}(0)| , \ |\pa_{1112} \Psi(0)| \les || \Om \vp_{g, 1}||_{\inf}.
\eal
\]
\end{lem}

The formula of $\pa_{1112} \Psi(0)$ can be written as an integral of $\Om$; see \eqref{eq:u_appr_near0_K}. In  Part II \cite[Section 4]{ChenHou2023b}, we develop the sharp version of the above estimates with better constants.  
In Appendix \ref{app:E_vel}, we present the proof, which also helps to illustrate the ideas for  Part II \cite[Section 4]{ChenHou2023b}.

In Section \ref{sec:sharp} and Part II \cite[Section 4]{ChenHou2023b}, for $ \uu = \na^{\perp}(-\D)^{-1} \Om$ with $\Om \in X_1$,  we develop weighted estimate for $ \na \uu -\wh{\na \uu} $ with 
approximations $\wh{\na \uu}$ constructed in Section \ref{sec:appr_vel}. In particular, 
 $ f = u_x(\Om) ,u_y(\Om), v_x(\Om)$, we obtain
\beq\label{eq:E_vel_hol}
|x-z|^{-1/2}  | ( f- \hat f) \psi_1(x) - ( f- \hat f) \psi_1(z)  |  \les || \Om \psi_1 ||_{C_{g_1}^{1/2}} + || \Om \vp_1||_{\inf} \les || \Om ||_{X_1}, 
\eeq
for $x, z$ with $x_1 = z_1$ or $x_2 = z_2$ and $|x| \leq |z| \leq (1 + \mu) |x|$ for some $\mu \in (0, 1)$.  Up to absolute constants, 
this estimate can be established by following the decomposition and argument in Part II \cite[Section 4]{ChenHou2023b} and using the asymptotics of the weights $\vp_1, \psi_1$ \eqref{wg:hol}, \eqref{wg:linf_decay}. When $ |z| > (1 + \mu) |x|$, the estimate follows from the above Lemma, the triangle inequality, 
$|x|^{-1/2}  \psi_1 \les \vp_1, || \Om \vp_1 ||_{\inf} \les || \Om ||_{X_1}$ 
by \eqref{wg:linf_decay}, \eqref{wg:hol}, and \eqref{norm:X}. Since near $x=0$, $f - \hat f$ agree with the left hand side of \eqref{eq:E_vel1} (see \eqref{eq:u_appr_near0_coe} and \eqref{eq:u_appr_1st}), we obtain
\[
||  (f - \hat f)  \psi_1 ||_{ C^{1/2}}  \les || \Om ||_{X_1}. 
\]

Note that the approximations in \eqref{eq:u_appr_1st}, \eqref{eq:u_appr_2nd} except for 
\[
C_{f0}(x)  ( -\pa_{12}( -\D)^{-1} \Om )(0) =C_{f0}(x)  u_x(\Om)(0)  , \quad C_f(x)\cK_{00}(\Om) \chi_0
\]
are supported away from $0$ with smooth coefficients. Moreover, the functionals in \eqref{eq:u_appr_1st}, \eqref{eq:u_appr_2nd}, e.g. $\fNS(\Om)(x_i, 0), \cK_{00}(\Om)$, are bounded by $|| \Om||_{X_1}$. 
See \eqref{eq:functional_w_est}. Using triangle inequality, we yield 
\[
|| \big(f - C_{f0}(x) (-\pa_{12}(-\D)^{-1} \Om)(0) - C_f(x) \cK_{00}(\Om) \chi_0 \big) \cdot \psi_1 ||_{C^{1/2}} \les || \Om||_{X_1}, 
\] 
where $\chi_0$ is defined in Section \ref{sec:u_appr_1st}, $\chi_0 = 1$ near $0$, and supported near $0$. In summary, we have 
\begin{lem}\label{lem:E_vel2}
Suppose that $\Om \in X_1$ \eqref{norm:X} is odd and $-\D \Psi = \Om$. We have 
\[
|| \psi_1 \B( \pa_{ij} ( \Psi - \Psi_{xy}(0) x_1 x_2 ) - \chi_0  \pa_{1112} \Psi(0) \pa_{ij} G(x)  \B)||_{C^{1/2}} \les || \Om||_{X_1}, 
\quad G(x) \teq \f{1}{6}( x_1^3 x_2 - x_1 x_2^3).
\]
\end{lem}

Using the above estimates for $\na^2 (-\D)^{-1}\Om$, we obtain the estimate for $\na (-\D)^{-1} \Om(x)$ by integration from $0$ to $x$, which is more regular.

\vs{0.05in}
\paragraph{\bf Classical embedding of $(-\D)^{-1}$}
Let $f(x)$ be odd in $x_1$, and $\phi_{2D}(f) = (-\D)^{-1} f, \uu_{2D}(f) = \na^{\perp} \phi_{2D}(f)$. By symmetrizing and decomposing  the singular integral, we obtain
\footnote{
\label{foot:u_direct}
We extend $f$ from $\R^2_+$ to $\R^2$  as an odd function $F$ in $x_2$. Then  $\na \uu_{2D}(f) = C_1 \int_{\R^2} K(x-y) F(y) dy + C f(x) $
for some constant $C$ and $K(z) = \f{z_1 z_2}{|z|^4}, \f{z_1^2  -z_2^2}{|z|^4}$.
 Since $f \in C_{x_1}^{\al}(\R^2_+)$, we obtain $F \in C_{x_1}^{\al}$.
 We decompose $\R^2$ into the singular region $S = \{ y: |x-y| \leq \la x \ra^{-100 - p} \}$ and regular region $\R^2 \bsh S$. Let $I_S, I_R$ be the integral in $S$ and out of $S$. Then we bound
 \[
| I_S | \les_{\b, p}   \la x \ra^{-1} ( [ F \la x\ra^{-p}]_{C_{x_1}^{\al}} + \nlinf{ \la x \ra^{\b} F } ) , \quad I_R \les_{\b,  p} \la x \ra^{-\b} \log ( | x |+ 2) \nlinf{ \la x \ra^{\b} F }
\les \la x \ra^{-\b/2 } \nlinf{ \la x \ra^{\b} F }.
 \]
}
\bseq\label{eq:u_est_direct}
\beq
   |\na \uu_{2D}(f) | 
   \les_{\al, \b, p}  \la x \ra^{- \b/2 }  ( \|  \la x \ra^{\b}   f \|_{L^{\infty}}
   + \| \la x \ra^{- p }  f \|_{C^{\al}(\R^2_+)} ), \quad   \b \in [0, \tf14 ], \  p \geq 0,  \
   \al > 0.
\eeq

Since $u_{2D}(0, x_2) = 0, v_{2D}(x_1,0)=0$, 
$(-\D)^{-1} f = 0$ on $x_1=0, x_2=0$, 
for  $  \b  \in (0, \tf14 ], \, p \geq 0, \al > 0$, integrating the above estimates, we obtain 
\beq
  |  \f{ \phi_{2D}(f) }{x_1 x_2 }| +  | \f{ u_{2D}(f) }{x_1} | +      | \f{ v_{2D}(f) }{x_2}  |  
   \les_{\al, \b ,   p}  \la x \ra^{- \b/2 }  ( \| \la x \ra^{\b}  f \|_{L^{\infty}}
   + \| \la x \ra^{- p}  f \|_{C^{\al}(\R^2_+)} ) . 
\eeq
\eseq

Similarly, we have the classical embedding inequalities 
\beq\label{eq:elli_embed}
||\na^2(-\D)^{-1} f ||_{L^{\inf}} \les_{ p} || f ||_{C^{1/2}} + || f ||_{L^p}, \ 
p \in (1, \infty),
\quad [\na^2(-\D)^{-1} f ]_{C^{1/2}} 
\les [f]_{C^{1/2}}.
\eeq

\subsubsection{Estimate of 3D velocity}

Below, we  consider the elliptic equation  $-\D \Psi =\Om(f)$ \eqref{eq:elli_main2} with source $\om = f$ in \eqref{eq:elli_main}. 
Let $\uu(f)$ be the associated 3D velocity.  
Using \eqref{eq:u_est_direct} with $\al = \f12, p=0, \b = \f{1}{16}$, 
 \eqref{wg:equi} with $i=2$, and $X_1$-norm \eqref{norm:X},  we bound 
\bseq\label{eq:xi_Psi_decay0}
\beq
| \f{1}{x_1 x_2} \Psi | +
 | \na^2 \Psi | + | \f{1}{x_1} \pa_2 \Psi | 
 + | \f{1}{x_2} \pa_1 \Psi| 
 \les \la x \ra^{-1/32} \big\| \f{|x|^2}{\la x \ra^2 } \Om(f  ) \big\|_{X_1} 
 \les \la x \ra^{-1/32} \big\| \f{|x|^2}{\la x \ra^{2} } f \big\|_{X_1}. 
\eeq

Within the support $|x| \leq S$, using $C_l S < 1, r \asymp 1$, 
and formulas of $\uu$ \eqref{eq:E_u}, \eqref{eq:E_v},  we bound 
\beq
   |\na \uu(f) | 
   +     | \tf{1}{x_1}  u(f) | +      | \tf{1}{x_2}  v(f)   |  
 \les  \la x \ra^{-1/32} \big\| \tf{|x|^2}{\la x \ra^{2} } f \big\|_{X_1}. 
\eeq
\eseq

Since the cutoff function $\chi_0$ in Lemma \ref{lem:E_vel2} is supported
in $|x| \leq 1$ by Appendix \ref{app:appr_vel_para},  for $\nu > 4$
and $|x|, |z| \geq \nu/4$, using Lemma \ref{lem:E_vel2} with $i+ j=2$, we bound 
$\uu_{2D} = \na^{\perp}(-\D)^{-1} f$:
\bseq\label{eq:E_du_hol}
\beq\label{eq:E_du_hol:2D}
|\d(  \psi_1 ( \na \td \uu_{2D}(f) ), x, z  ) |
  \cdot |x-z|^{-1/2}
 \les  \| f \|_{X_1}, \quad \d(F, x, z) \teq F(x) - F(z) .
\eeq
where as in \eqref{eq:u_tilde}, we define 
\[
 \td \uu_{2D} =   ( u_{2D} -u_{2D, x}(0) x , v_{2D} + u_{2D, x}(0) y ), \quad  u_{2D, x}(0) = -\pa_{12}(-\D)^{-1} f(0).
\]

 Applying \eqref{eq:E_du_hol:2D} to  $-\D \Psi = \Om(f)$ in \eqref{eq:elli_main2}, 
 \eqref{eq:xi_Psi_decay0} for the lower order terms in 
the formulas  \eqref{eq:E_u}, \eqref{eq:E_v} for $\uu$, and estimate \eqref{wg:equi},  for $|x|, |z| \in [ \nu / 4, |S| ]$, we bound
\beq\label{eq:E_du_hol:X1}
   ( | \d( \psi_1 ( \na  \uu ) ( f ) ,x, z)| 
   +  | \d( \psi_1 ( \na \td \uu ) ( f ) ,x, z)| ) \cdot |x-z|^{-1/2} \les 
   \| \Om(f ) \|_{X_1} \les \| f \|_{X_1} \, , 
\eeq
where we  used $ \psi_1 \in C^{1/2}(B_1^c)$
and $ |u_{2D, x}(f)(0 )|, |u_{x}(f )(0)| \les \xxa{f }$ by 
\eqref{eq:elli_main3} and \eqref{eq:u_est_direct}.

Recall the support size $S$ from Definition \ref{def:supp}. Applying  \eqref{eq:elli_embed} for the top order $\Psi$-term in \eqref{eq:E_u}, \eqref{eq:E_v} for $\uu$, estimate \eqref{eq:xi_Psi_decay0} for lower order terms, and \eqref{wg:equi}, for $|x|, |z| \leq |S|$, we bound
\beq\label{eq:E_du_hol:reg}
\| \na \uu(f) \|_{C^{1/2}( \{ |x| \leq |S| \}) }
\les \| \Om( f) \|_{C^{1/2}} + \| \tf{|x|^2}{ \la x \ra^2} f\|_{X_1} 
\les \| f \|_{C^{1/2}} + \| \tf{|x|^2}{ \la x \ra^2} f\|_{X_1} .
\eeq

\eseq

\subsection{Nonlinear stability}\label{sec:3Dnon}

In this section, we extend the nonlinear stability estimates
for 2D Boussinesq to 3D Euler \eqref{eq:euler4}.
 In Section \ref{subsec:non_boot}, we impose the bootstrap assumption on the support size. In Section \ref{subsec:non_appr}, we construct the approximate steady state and impose the normalization conditions. 
In the remaining parts,
we show that the differences between the two systems are perturbative,
and can be controlled  by the small parameters $C_lS$, $C_l$, and 
$\nu^{-\gamma}$.

\subsubsection{Bootstrap assumption on the support size}\label{subsec:non_boot}

We will use Propositions \ref{prop:C1}, \ref{prop:C2}, \ref{prop:elli} 
with $\al = 1$ or $\al = 2.9$ and  $\b = \f{1}{16} $, which relate to the singular weights 
$\vp_1, \vp_{g, 1}$ in \eqref{wg:linf_decay}, \eqref{wg:linf_grow}.
Then the constants $\nu_1, \nu_2$ in these propositions are determined. 

\begin{assumption}[\bf Support size]\label{boot_ass:1} We impose the first bootstrap assumption: 
\beq\label{eq:boot5}
\bal
C_l( t) (1 + \max( S(t ), S(0)) ) < \min( \nu_2, 4^{-6}), 
 \quad \forall \, t \geq 0.
\eal
\eeq
\end{assumption}
Under \eqref{eq:boot5}, the support of $\om, \th$ in $D_1$ does not touch the symmetry axis and $z = \pm 1$, the cutoff functions \eqref{eq:chi_Ri}  satisfy $\chi_i = 1 , i \leq 5$  
for $x$ in the support, 
and assumptions in Propositions \ref{prop:C1}, \ref{prop:C2}, \ref{prop:elli} hold.  We will choose $C_l(0)$ at the final step to ensure the smallness in \eqref{eq:boot5}.

\subsubsection{Approximate steady state and normalization conditions}\label{subsec:non_appr}
Since the rescaled domain $\td D_1$ \eqref{eq:rescale_D} is bounded, we localize the approximate steady state $\bar{\om}, \bar{\th}$ for the 2D Boussinesq constructed in Section \ref{sec:ASS} to construct the approximate steady state for \eqref{eq:euler4} with bounded support 
\beq\label{eq:profile2}
\bar{\om}_0 \teq \chi_{\nu} \bar{\om}, \quad \bar{\th}_0 \teq \chi_{\nu} ( m_{\nu} + \bar{\th}),
\quad  \chi_{\nu}(x)= \chi_\dap( |x| / \nu), \quad m_{\nu} = 2 \max\nolim_{|x| \leq 2 \nu} |\bar \th(x) |,
\eeq
with $\nu \geq 1$, Here
$\chi_\dap(y) = \td \chi_\dap(|y|)^2$, $\td \chi_\dap: \R \to [0,1]$ is 
a smooth, even cutoff function, and satisfies
\beq\label{def:chi_dap}
\chi_\dap(y) = \td \chi_\dap(|y|)^2, \quad 
 \td \chi_\dap(z)  \ \mbox{is decreasing in } |z|  ,  \quad 
\td \chi_\dap(z) = 1 , |z| \leq 1, 
\quad \td \chi_\dap(z) = 0, z \geq 2. 
\eeq

Clearly, from Definition \ref{def:supp}, the support size of $\bar{\om}_0, \bar{\th}_0$ is $2 \nu$. We add $m_{\nu}$ and truncate $m_{\nu} + \bar \th$ rather than $\bar \th$ so that $2 m_{\nu} \geq m_{\nu} + \bar \th  \geq m_{\nu}/2$ for $|x| \leq 2 \nu$ and $(m_{\nu} + \bar \th)^{1/2}$ has the same regularity as $\bar \th$. This idea follows \cite{chen2019finite2}. By definition, we obtain
\beq\label{eq:profi_agree}
 \bar \om_0(x) = \bar \om(x), \quad \na \bar \th_0(x) = \na \bar \th(x),  \quad \forall |x| \leq \nu. 
\eeq
We will choose initial data with support size $S(0) \geq \nu$. By Definition \ref{def:supp}, we obtain 
\beq\label{eq:supp_para}
\nu \leq S(0) \leq S(\tau) , \quad \tau \geq 0. 
\eeq

Consider the polar coordinates for $(x, y)$: $\sig = \arctan \f{y}{x}, \rho = \sqrt{x^2 + y^2}$. 
Using the formula 
\[
\pa_x g = (\cos \sig \pa_\rho - \rho^{-1} \sin \sig   \pa_{\sig} ) g, \quad  \pa_y  g= (\sin \sig \pa_\rho + \rho^{-1}  \cos \sig  \pa_{\sig} ) g ,
\]
we derive
\[
\bal
\pa_x \bar \th_{0 } = \chi_{\nu} \bar \th_x + \nu^{-1} \cos \sig \chi_\dap^{\pr}( |x| / \nu ) ( m_{\nu} + \bar \th)  , 
\quad \pa_y \bar \th_{ 0 } = \chi_{\nu} \bar \th_y + \nu^{-1} \sin \sig \chi_\dap^{\pr}( |x| / \nu ) (m_{\nu} + \bar \th)  . 
 \eal
 \]

To show that $\bar \om - \bar \om_0, \bar \th- \bar \th_0$ is small, 
we will use the decay estimates \eqref{solu:decay} and choose large $\nu$. 
To distinguish the notations between 3D Euler and  2D Boussinesq equations, we write 
\beq\label{eq:u_2D}
\bar \phi_{2D} = (-\D)^{-1} \bar \om, \qquad \bar \uu_{2D} = \na^{\perp} \bar \phi_{2D}
\eeq
for 2D Boussinesq. Let $\bar{\phi}$ and $\bar \uu$ be the stream function and velocity in \eqref{eq:euler4} associated with $\bar \om_{0}$.  
The leading-order decomposition of $\bar\uu$ is given in \eqref{eq:E_ub}.
See more discussions in Section \ref{sec:E_ub}.

We need to adjust the time-dependent normalization condition for $c_{\om}(t), c_l(t)$. We impose 
\footnote{
  We remark that $\bar{c}_{\om}(t) $ is \emph{time-dependent} since it depends on $ \bar u_x(0)$ and the elliptic equation in \eqref{eq:euler4} depends on the rescaling factor $C_l(t)$.
  Since $C_l(t)$ given in  is \emph{solution-independent}, $\bar{c}_{\om}(t)$ 
  is also \emph{solution-independent}.
}
\bseq\label{eq:E_normal}
\beq\label{eq:E_normal1}
\bar c_l = 2 \f{ \bar \th_{xx}(0) } { \bar \om_x(0)},  \quad \bar c_{\om}(t) = \f{1}{2} \bar c_l + \bar u_x(0), \quad \bar c_{\th}(t) = \bar c_l + 2 \bar c_{\om}(t)
\eeq
for the approximate steady state $ \bar \om_0, \bar \th_0$, and 
\beq\label{eq:E_normal2}
c_l(t) = 0, \quad c_{\om}(t) = u_x( t, 0)
\eeq
for the perturbations, 
where $u_x(t, 0)$ is obtained from velocity $u$ in \eqref{eq:euler4},  
and \emph{differs from} $-\pa_{xy} (-\D)^{-1} \om (0)$.  The above conditions are the same as \eqref{eq:normal} and \eqref{eq:normal_pertb}, and play the same role of enforcing \eqref{eq:normal1}. As a result,  perturbation $ \om, \na \th$ satisfies the vanishing condition \eqref{eq:normal_vanish}
\[
\om = O(|x|^2), \quad \na \th = O(|x|^2) ,
\quad \mbox{ near} \  x = 0.
\]
Since $ \na \bar \th_0 = \na\bar \th, \bar \om_0 = \bar \om $ near $x=0$, $\bar c_l$
in \eqref{eq:E_normal1} is the same as that for 2D Boussinesq \eqref{eq:normal}. 

From \eqref{eq:E_normal} and \eqref{eq:rescal2}, the scaling factor $C_l(\tau)$ 
is \emph{solution-independent} and reduces to
\beq\label{eq:E_normal:c}
\bar c_l + c_l(\tau) = \bar c_l ,\quad C_l(\tau)  = C_l(0) \exp( -\bar c_l \tau) .
\eeq
Thus, $C_l(t)$ in the elliptic equation \eqref{eq:euler4}
and \eqref{eq:boot5}
is determined \emph{a priori}, up to choosing $C_l(0)$.
\eseq

 From the estimate in Proposition \ref{prop:C2}, $\bar u_x(0)$ is very close to $-\pa_{xy}(-\D)^{-1} \bar \om_{0}$. For $\nu$ sufficiently large, comparing the above conditions and \eqref{eq:normal}, $\bar c_{\om}$ is very close to $ \bar c_{\om, 2D}$ \eqref{eq:bous_expo} used for the 2D Boussinesq  in Section \ref{sec:lin}. From \eqref{eq:E_normal1} and \eqref{eq:normal}, we yield 
\beq\label{eq:cw_dif}
\bar c_{\om} - \bar c_{\om, 2D} = \bar u_x(0) - \bar u_{x, 2D}(0).
\eeq

\begin{assumption}[Outgoing]\label{boot_ass:2}

We impose a bootstrap assumption on an outgoing condition 
\beq\label{eq:boot_ass:2}
( \bar c_l \xx + \uu_{\tot}( t, \xx) ) \cdot \xx |_{|\xx| = \nu/2} > 0, \, \forall  \, t \geq 0.
\eeq
\end{assumption}

Using normalizations \eqref{eq:E_normal1}, \eqref{eq:E_normal2} and restricting the $\th$-equation in \eqref{eq:euler4} on $x=0$ and taking $\pa_y$-derivative, we obtain 
the equation for the full solution $\th_{\tot} = \bar \th_0 + \th, v_{\tot} = \bar v + v$:
\[
 \pa_t \pa_y \th_{\tot} + ( \bar c_l y + v_{\tot}(0, y)) \pa_y ( \pa_y \th_{\tot}) = ( c_{\th, \tot} -\bar c_l - \pa_y v_{\tot}(0, y) )
 \cdot \pa_y \th_{\tot}.
\]

If $ \pa_y \th_{\tot}(0, y) = 0$ for $y \in [0, \nu / 2]$ initially, under Assumption \ref{boot_ass:2} 
and by  \eqref{eq:profile2}, 
we obtain 
\beq\label{eq:3D_thy_van}
 \pa_y \th_{\tot}(0, y, t) =0 , \  \pa_y \bar \th_0(0, y) = 0,
\ \Longrightarrow \quad  \pa_y \th(0, y, t) = 0, \quad \forall \, y \in [0, \nu/2].
\eeq

Unlike 2D Boussinesq, due to modification of the $ \th$-profile \eqref{eq:profile2}, 
the profile and perturbation do not satisfy $\pa_y \th_{\tot}(0,y)=0, \pa_y \th(0,y)$ \emph{for all} $y$. In \eqref{wg:xi_modi}, we will also modify the singular weight $\vp_{3,m}, \vp_{g, 3, m}$ for $\xi = \th_y$ away from $|(x, y)| \geq \nu/2$ so that the energy is well-defined.

\begin{remark}[\bf Size of parameters]\label{rem:order2}
We will first choose $\nu \gg 1$, then $\nu \leq S(0)$, and finally take $C_l$ sufficiently small so that $C_l \nu, C_l S \ll 1$. Due to these small factors and using \eqref{eq:elli_mainb}, \eqref{eq:E_ub}, we can treat 
$\bar \om_0 \approx \bar \om, \bar \th_0 \approx \bar \th$, $\bar \uu \approx \bar \uu_{2D}$. From Remark \ref{rem:order1} and the bootstrap assumption \eqref{eq:boot5}, we also have $C_l \approx 0, C_l S \approx 0, r \approx 1$.  We treat the errors in these approximations perturbatively.
\end{remark}

\subsubsection{Linearized equations}

The equations \eqref{eq:euler4} are slightly different from \eqref{eq:bousdy1} for the Boussinesq equations. Denote 
\[
\eta = \th_x, \quad \xi = \th_y. 
\]

Linearizing \eqref{eq:euler4} around the approximate steady state $(\bar \om_0, \bar \th_0, \bar c_l, \bar c_{\om})$ \eqref{eq:profile2}, \eqref{eq:E_normal1}, we obtain the equations for the perturbation $W =  (\om, \eta, \xi)$, which are similar to \eqref{eq:bous_lin11}, \eqref{eq:bous_lin12}
\beq\label{eq:Euler_lin}
\bal
\pa_t \om  &= - (\bar c_l x +\bar \uu) \cdot \na \om + \f{1}{r^4} \eta + \bar c_{\om} \om 
 -  \uu \cdot \na \bar \om_0 + c_{\om} \bar  \om_0  + \cN_1 + \bar \cF_1 
\teq \cL_1 W + \cN_1  + \bar \cF_1 , \\
 \pa_t \eta & =  - (\bar c_l x +\bar \uu) \cdot \na \eta +
(2 \bar c_{\om} - \bar u_x) \eta - \bar v_x \xi - \pa_x ( \uu \cdot \na \bar \th_0 )
 + 2 c_{\om} \bar \th_{0,x} 
+  \cN_2 + \bar \cF_2  
 \teq  \cL_2 W   + \cN_2 + \bar \cF_2
 , \\
  \pa_t \xi & = -  (\bar c_l x +\bar \uu) \cdot \na \xi +
(2 \bar c_{\om} - \bar v_y) \xi  - \bar u_y \eta - 
\pa_y( \uu  \cdot \na \bar \th_0 )
 + 2 c_{\om} \bar \th_{0, y}
+ \cN_3 + \bar \cF_3 
\teq \cL_3 W + \cN_3 + \bar \cF_3 ,
\eal
\eeq
where 
\[
\bar \cF_1
= - (\bar c_l x + \bar \uu ) \cdot \na \bar \om_0 + \tf{1}{r^4} \bar \th_{0, x} + \bar c_{\om} \bar\om_0 ,
\]
and we adopt similar notations $\cN_i, \bar \cF_i$ for other nonlinear terms and the error terms from \eqref{eq:bous_non}, \eqref{eq:bous_err}. The $\om$ equation is different from the corresponding equation in \eqref{eq:bous_lin11} since we have $\f{1}{r^4} \th_x$ in \eqref{eq:euler4}. The $\xi$ equation is also different from the corresponding equation in \eqref{eq:bous_lin12} since we do not have the same incompressible conditions $\bar u_x + \bar v_y = 0, \ u_x + v_y = 0$. We remark that the velocity $\uu, \na \uu$ in the above system are determined by the elliptic equation in \eqref{eq:euler4}.

To generalize the 2D Boussinesq analysis to 3D Euler, we show that the  additional terms 
can be made \emph{arbitrarily small}. 
In the following derivations, {\bf we use $f_{2D}$
to denote the quantity $f$ from the 2D Boussinesq}. 
For example, $\bar \uu_{2D}= \na^{\perp}  (-\D)^{-1}\bar \om$ 
\eqref{eq:u_2D} denotes the approximate steady state for the velocity for 2D Boussinesq. 

\vs{0.05in}
\paragraph{\bf Modified weights for $\xi$ and  energy}
Recall the weights $p_{6,5} |x|^{1/7}$ in $\vp_3$ and  $p_{9, 2} |x|^{\al_{g, n}}$ in $\vp_{g, 3}$ 
from \eqref{wg:linf_decay}, \eqref{wg:linf_grow}.  Since $\xi$ for $\th_y$ do not vanish on $x=0$ for any $y$ (see \eqref{eq:3D_thy_van} and discussions therein), we modify  $\vp_3, \vp_{g, 3}$ 
by replacing the weight beyond $|x| \geq \nu/2$ :
\begin{align}\label{wg:xi_modi}
\vp_{3, m}(x) = \one_{|x| < \f{\nu}{2}} \vp_3 + \one_{|x| \geq \f{\nu}{2}  } 
p_{6, 5} |x|^{1/ 7},  \ 
      \vp_{g, 3, m}(x)  \teq \one_{|x| <  \f{\nu}{2}}  \vp_{g, 3}
      +  \one_{|x| \geq \f{\nu}{2}  }  p_{9, 2} |x|^{\al_{g, n}} .
\end{align}

Recall $E_1$ \eqref{energy1}, $E_2,E_{2,0}$ \eqref{energy2}, $E_3$ \eqref{energy3}, and $ E_4,E_{4, 0}$ \eqref{energy4}. We modify energy $E_4$ by replacing the weight $\vp_{3}$ 
in $E_1$ \eqref{energy1} and $\vp_{g, 3}$ in $E_3$ \eqref{energy3} by $\vp_{3,m}, \vp_{g, 3, m}$ :
\beq\label{energy4m}
\bal
E_{4, m}(W_1, W) \teq  &
 \max\big(   \max_{i=1,2,}  || W_{1,i} \vp_i ||_{\inf}, 
\nlinf{ W_{1, 3} \vp_{3, m} }, \tau_1^{-1}  \sqrt{2} || \om_1 \vp_4 ||_{\inf} , \\
 &  E_{2, 0}(W_1) , \, \max_{i\leq 2} \mu_{g, i} || W_{1, i} \vp_{g, i}  ||_{\inf},  
 \mu_{g, 3} \nlinf{ W_{1, 3} \vp_{g, 3, m} }, \, E_{4, 0}(W_1, W)  \big) .
\eal
\eeq

We simplify $E_{4, m}(W_1(t), W(t))$ as $E_{4, m}(t)$.  The energy $E_4$ in \eqref{energy4} and $E_{4,m}$ differ only by two weighted $L^{\infty}$ norms on $W_{1,3} = \xi_1$.  We introduce the norm $X_i$ related to the energy $E_{4, m}$ 
\beq\label{norm:Xi}
\bal
|| f||_{X_i}   \teq || f \vp_i||_{\infty} + || f \vp_{g, i}||_{\infty} + [ f \psi_i ]_{C_{g_i}^{1/2}} ,
\  i = 1,2,  \
 || f ||_{X_3}   \teq || f \vp_{3, m}||_{\infty} + || f \vp_{g, 3, m}||_{\infty} + [ f \psi_3 ]_{C_{g_3}^{1/2}} .
\eal
\eeq

The norms for $f_1, f_2$ are the same as those for $\om_1, \eta_1$ in $E_4$ \eqref{energy4}.

\begin{remark}[\bf Modified weight]
Equation \eqref{wg:xi_modi} is the only modification of the 2D weights needed for the smooth 3D swirl construction. 
The modified profile in \eqref{eq:profile2} preserves the smoothness
of $\bar \th$, but it destroys vanishing of \(\theta_y(0, y)\) in the far-field with $y > \nu$.  The modified weights $\vp_{3,m}, \vp_{g,3,m}$ are used only 
for $|(x, y)| \geq  \nu / 2$. With this modification, the weighted perturbation are bounded :
$ \vp_{3, m} \xi_1,   \vp_{g, 3, m} \xi_1 \in L^{\infty}$.  
By choosing $\nu$ sufficiently large, the 
associated stability estimates in the region $|(x, y)| \geq  \nu/2$ is \emph{perturbative}. 
See more details in Sections \ref{sec:euler_other_xi}, \ref{sec:E_xi_lower}.  For $|(x, y)| <  \nu / 2$, the original Boussinesq weights and estimates generalize directly.
\end{remark}

\paragraph{\bf{Lower order terms in the linearized and nonlinear operator}}

Using \eqref{eq:E_uR}, we decompose the linear operator $\cL_i$ into the main part $\cL_{M,i}$
and the remainder part $\cL_{R,i}$
\beq\label{eq:E_lin_MR}
\bal
\cL_{i } W  & = \cL_{M, i } W + \cL_{R, i }( \om ,\eta), \qquad \qquad W = (\om,\eta,\xi),\\
\cL_{M, 1} &= - ( \bar c_l x + \bar \uu) \cdot \na \om +  \eta - \na^{\perp} \Psi \cdot \na \bar \om_0 - \Psi_{xy}(0) \bar \om_0  + \bar c_{\om} \om, \\ 
\cL_{M, 2} & =   - (\bar c_l x +\bar \uu) \cdot \na \eta +
(2 \bar c_{\om} - \bar u_x) \eta - \bar v_x \xi - 
\pa_x ( \na^{\perp} \Psi \cdot \na \bar \th_0) - 2 \Psi_{xy}(0) \bar \th_{0,x} , \\
\cL_{M, 3} & =   - (\bar c_l x +\bar \uu) \cdot \na \xi +
(2 \bar c_{\om} - \bar v_y) \xi - \bar u_y \eta - 
\pa_y ( \na^{\perp} \Psi \cdot \na \bar \th_0)- 2 \Psi_{xy}(0) \bar \th_{0,y} , \\
\cL_{R,1}  & = (\f{1}{r^4} - 1) \eta -  \uu_R \cdot \na \bar \om_0 ,\quad 
\cL_{R, 2}  = - \pa_x( \uu_R \cdot \na \bar \th_0)  ,\quad 
\cL_{R, 3}  = - \pa_y( \uu_R \cdot \na \bar \th_0) . 
\eal
\eeq
Note that $\cL_{R,i}$ \emph{only} depends on $\om, \eta$.
From \eqref{eq:elli_main3}, \eqref{eq:E_uR}, and \eqref{eq:E_v}, we have 
\beq\label{eq:E_uR0}
u_R = O(|x|^3), \quad u_{R, x}(0) = 0, \quad 
v_R = O(|x|^2), \quad 
 - \Psi_{xy}(0) = - \phi_{xy}(0) = u_{x}(0) = c_{\om}.
\eeq
See also the remark of $v_R$ below \eqref{eq:E_v}. We will estimate  $\cL_{R, i}$ in Section \ref{sec:E_lin_lower} and  bound it by $C_l S  E_{4, m}(t)$, where $E_{4, m}$ is defined in \eqref{energy4m}.

We decompose nonlinear terms $\cN_i$ \eqref{eq:bous_non} 
by decomposing velocity $\uu = \uu_M + \uu_R$  \eqref{eq:E_uR},\eqref{eq:E_v}:
\beq\label{eq:E_non_MR}
\bal
 \cN_1 & =  - \uu \cdot \na \om  + u_{M, x}(0) \om \teq  \cN_{M, 1}, \quad \cN_{R, 1} \equiv 0 , \\
\cN_2  &=  \big( -\uu \cdot \na \eta - u_{M, x} \eta - v_{M, x} \xi + 2 u_{M, x}(0) \eta \big)
+ \big( - u_{R, x} \eta - v_{R, x} \xi \big) \teq \cN_{M, 2} + \cN_{R, 2},  \\
\cN_3  &=  \big( -\uu \cdot \na \xi - u_{M, y} \eta - v_{M, y} \xi + 2 u_{M, x}(0) \xi \big)
+ \big( - u_{R, y} \eta - v_{R, y} \xi \big) \teq \cN_{M, 3} + \cN_{R, 3},  \\
\eal
\eeq
where we used  $u_{R, x}(0) = 0, u_{M, x}(0) =  u_{x}(0) = c_{\om}$ \eqref{eq:E_uR0}. 
We keep the \emph{full} velocity $\uu$ in the transport term. The lower order term $\cN_{R, 2}$ vanishes $O(|x|^4)$ near $x =0$ and is estimated by  the estimates of $\uu_{R, x}$ in Section \ref{sec:E_lin_lower} and the nonlinear terms in Section \ref{sec:non}.  In the weighted estimate, the transport term leads to the nonlinear term $\cT_{d, N}(\rho) W_{1, i} \rho$ in \eqref{eq:tran_extra}.  We decompose it into the main and lower order terms  similarly,
and estimate  the lower order terms as $\cN_{R, 2}$.

\vs{0.1in}
\paragraph{\bf{Lower order terms in the residual error}}

We introduce the difference variables
\beq\label{eq:dif_op}
\d f = f - f_{2D} ,\quad \d \bar \uu =  \bar \uu - \bar \uu_{2D},  \quad \d \bar c_{\om} = \bar c_{\om} - \bar c_{\om, 2D}. 
\eeq

For the residual error, using \eqref{eq:E_ub} and \eqref{eq:bous_err}, we obtain 
\begin{align}\label{eq:E_err_MR}
\bar \cF_i &= \bar \cF_{M, i} +  \bar \cF_{R, i},  \\
\bar \cF_{M, 1} & =  - (\bar c_l x + \bar \uu_{2D} )  \cdot \na \bar \om_0 +
\bar \th_{0, x} + 
 \bar c_{\om, 2D} \bar \om_0 , 
\quad 
\bar \cF_{R, 1} =  - \d \bar \uu \cdot \na \bar \om_0 + \d \bar c_{\om}  \cdot \bar \om_0 + \f{1 - r^4}{r^4} \bar \th_{0, x}, \notag 
\\
\bar \cF_{M,i+1} & = \pa_i \big( - ( \bar c_l x + \bar \uu_{2D} ) \cdot \na \bar \th_0 + (\bar c_l + 2 \bar c_{\om, 2D}) \bar \th_{0}\big) , \quad \bar \cF_{R, i+1} = \pa_i\big(  - \d \bar \uu \cdot \na \bar \th_0 
+ 2 \d \bar c_{\om} \cdot  \bar \th_0 \big), \quad i=1,2. \notag 
\end{align}

Using the decay estimates for $\bar \om, \na \bar \th$ \eqref{solu:decay}
and the weights \eqref{wg:linf_grow}, \eqref{wg:linf_decay}, we obtain 
\[
| |x|^i | \na^i \bar \om|  \vp_{g,1}| \les |x|^{-1/16}, \quad | |x|^i |\na^{i+1} \bar \th| \vp_{g, 2} | \les |x|^{-1/16},
\quad   |x| \geq 1 ,  \quad i \leq 3.
\]
By \eqref{eq:profi_agree}, 
for $|x| \leq \nu$, 
$\bar \cF_{M,i}= \bar \cF_i$ for the 2D Boussinesq error \eqref{eq:bous_err}, where $\nu$ is the size of the cutoff function in \eqref{eq:profile2}. Using \eqref{wg:xi_modi} and \eqref{norm:Xi}, we have 
\beq\label{eq:E_err_M}
 || \bar \cF_{M,i} - \bar \cF_{ 2D, i} ||_{X_i} \les \nu^{-\g}  \, ,
\eeq
for some $\g > 0$, e.g. $\g = \f{1}{8}$. The H\"older estimate of the tail is even smaller since $\psi_1$ and $ | |x|^k \na^k \bar \cF_{ 2D, i} |, | |x|^k \na^k \bar \cF_{M, i} |, k=0,1$ decay, which can be derived using the regularity and asymptotics of the profile (see Section \ref{sec:property}).

The error term $\bar \cF_{R, i}$ does not vanish to order $|x|^{3}$ near $0$. 
We use correction $D_i^2 \bar \cF_{R, i}(0) \cdot f_{\chi, i}, D^2 = (\pa_{xy}, \pa_{xy}, \pa_{xx})$ similar to \eqref{eq:appr_near0} and then incorporate it in the 
decomposition \eqref{eq:euler_decoup}.

\vs{0.1in}
\paragraph{\bf Estimates of difference}

Below, we estimate several difference between  estimates for 2D Boussinesq and 3D Euler. 
In Section \ref{sec:E_lin_lower}, we estimate the lower order linear operator $\cL_{R, i}$ \eqref{eq:E_lin_MR}.
In Section \ref{sec:E_err_lower}, we estimate the lower order terms in the residual error 
$\bar \cF_{R, i}$ \eqref{eq:E_err_MR}. 
In Section \ref{sec:3D_finite_rank}, we modify the finite rank decomposition 
in Section \ref{sec:finite_rank}. 
In Section \ref{sec:E_xi_lower}, we perform the weighted $L^{\infty}$ estimates 
for $\xi_1$ for large $|x|$ due to modification of weights in \eqref{wg:xi_modi}.
In Section \ref{sec:euler_other_xi}, we estimate the difference in other estimates 
involving $\xi_1$ due to change of weights in \eqref{wg:xi_modi}, \eqref{energy4m}. 
In Section \ref{sec:E_lin_M}, we summarize the differences in the energy estimates.

\subsubsection{Estimate the lower order linearized operator $\cL_{R,i}$}\label{sec:E_lin_lower}

In this section, for \emph{arbitrary} $\om \in X_1, \eta \in X_2$ 
\footnote{ 
We apply \eqref{eq:E_lin_lower} only to $\om_1$ and $\eta_1$, which solve the $W_1$-equation in \eqref{eq:euler_decoup}; it is not applied to $\hat\om_2$ and $\hat\eta_2$, which solve the $\hat W_2$-equation in \eqref{eq:euler_decoup}, nor to $(\om,\eta)$, which solves \eqref{eq:Euler_lin}.
}
with support size $S$ satisfying  $C_l( t) (1 + S ) < \min( \nu_2, 4^{-6})$
similar to \eqref{eq:boot5}, we estimate 
\beq\label{eq:E_lin_lower}
 ||\cL_{R, i} \vp_{i} ||_{\infty} +
||\cL_{R, i} \vp_{g, i} ||_{\infty} + [ \cL_{R, i} \psi_{ i} ]_{C_{g_i}^{1/2}} 
+ || \cL_{R, 1} |x_1|^{-1/2} \psi_1 ||_{\inf}    \les 
 C_l S   (  \xxa{\om} + \xxb{\eta} ),
\eeq
 for $\cL_{R, i} = \cL_{R,i}(\om , \eta)$ and $i\leq3$. 
For convenience, we estimate weighted $L^{\infty}$ norms with
weights $\vp_3,\vp_{g,3}$ stronger than $ \vp_{3, m}, \vp_{g, 3, m}$ \eqref{wg:xi_modi}. Given $\om $, we define $\Om$ via \eqref{eq:elli_main}. 
By \eqref{wg:equi}, we have 
\beq\label{norm:X_E4}
|| \Om||_{X_1} \les  || \om ||_{X_1}.
\eeq

The estimate of $(r^{-4}-1) \eta$ 
in $\cL_{R, 1}$ \eqref{eq:E_lin_MR} follows from $ |r^{-4}-1| \les C_l x_2 \les C_l S$ 
for $x \in \supp(\eta)$ and the norm
$\xxb{\eta}$ for $\eta$.  Other terms in $\cL_{R, i}$ are nonlocal, involving $u_R, v_R$. We estimate a typical term $\pa_x u_R \bar \th_{0, x}. $

\vs{0.1in}
\paragraph{\bf{Estimate of $\Psi_2$}}

Recall the formulas of $u_R , v_R$ from \eqref{eq:E_uR}, \eqref{eq:E_v}  and $\Psi, \Psi_2$ from \eqref{eq:elli_main}. The estimates of the terms involving $\Psi_2$ are simple since 
\[
\Psi_2 = - a \, C_l^2 \phi_{xy}(0) \kp . 
\]

Recall the definition of $\kp$ from Proposition \ref{prop:elli}. Using \eqref{eq:elli_main3}, we yield 
\beq\label{eq:E_psi2}
|\na^k \Psi_2| \les C_l^2 |x|^{4-k} \xxa{\om}, \quad k \leq 4.
\eeq

For $|x| \leq 2$, we have 
\[
 |\na \bar \th_0| \les |x_1|, \quad | \na^2 \bar \th_0| \les 1.
\]

Consider a typical term related to $\Psi_2$  in  $\pa_x u_R \bar \th_{0, x}$, e.g. $\pa_{xy} \Psi_2 r^{-1/2} \bar \th_{0, x}$ \eqref{eq:E_uR}. We estimate
\[
\bal
 |\pa_{xy} \Psi_2 r^{-1/2} \bar \th_{0, x} \vp_2| & \les 
 C_l^2 |x|^2 |x_1|  \cdot  |x|^{-{5/2}} |x_1|^{-1/2} \one_{|x| \leq 2}  \xxa{\om}  \les C_l^2 \xxa{\om}. 
 \eal
\]
Note that for $|x| \leq 1$, we have $\chi = 1, \kp = - \f{x y^3}{6}$ (see Proposition \ref{prop:elli}), 
\[
\pa_{xy} \Psi_2 r^{-1/2} \bar \th_{0, x} \psi_2  = C  C_l^2 r^{-1/2} y^2 \bar \th_{0, x} \psi_2,
\]
for some absolute constant $C$. Since $\psi_2 \asymp c |x|^{-5/2}$ near $x=0$,  $ f= r^{-1/2} y^2 \bar \th_{0, x} \psi_2$ vanishes to order of $|x|^{1/2}$ near $x=0$. For $|x| > 1/2$, $f$ is smooth and is supported near $x=0$. Hence, we obtain that $f$ is in $C^{1/2}$ and 
\[
|| \pa_{xy} \Psi_2 r^{-1/2} \bar \th_{0, x} \psi_2||_{C_{g_2}^{1/2}} \les C_l^2 \xxa{\om}.
\]

The estimates of other terms related to $\Psi_2$ or $\bar \Psi_2$ \eqref{eq:E_ub}  in the residual error, nonlinear terms, or linear parts related to $\Psi_2$ follow similar estimates since $\Psi_2, \bar \Psi_2 = O(|x|^4)$ near $x=0$ and contain the small factor $C_l^2$. We treat them as lower order terms.

\vs{0.1in}
\paragraph{\bf{Estimate of $\Psi - \Psi_{xy}(0) x y$}}

Next, we estimate other terms in $\uu_R$ related to $\Psi$. Recall the decomposition \eqref{eq:E_uR2}. The third term in \eqref{eq:E_uR2} follows an estimate similar to that of $\Psi_2$ performed above. The first two terms vanish to a higher order near $x=0$. We estimate a typical term 
$\Psi_y( r^{-1/2} - 1 )$ in $\pa_x u_R \bar \th_{0, x}$ related to $\Psi$:
\beq\label{eq:E_lin_lower_I}
 I \teq \pa_x (\Psi_y -  \Psi_{xy}(0) x ) (r^{-1/2 } - 1) \bar \th_{0, x}
 = ( \Psi_{xy} - \Psi_{xy}(0) ) (r^{-1/2 } - 1) \bar \th_{0, x}.
\eeq

Note that $|r^{-1/2} - 1| \les C_l |x_2| \les C_l S$ within the support of the solution. The weighted $L^{\inf}$ estimate is simple and follows from Lemma \ref{lem:E_vel}. For example, 
using $\bar \th_{0, x} \les \min( |x_1|,  |x|^{-3/5})$ from \eqref{solu:decay},
and the definition of $\vp_{g, 2}$ from \eqref{wg:linf_grow}, we have 
\[
\bal
| I  \vp_2 | \les | I \vp_{g, 2}| \les \one_{|x| \leq S} C_l |x_2|  \min(|x|^2, 1) 
\min( |x_1|,  |x|^{-3/5} ) \vp_{g,2} || \Om||_{X_1} \les  C_l S \xxa{\om} ,
\eal
\]
where the weights $\vp_2, \vp_{g, 2}$ are defined in \eqref{wg:linf_grow}, \eqref{wg:linf_decay} with  $\vp_{g, 2}\les |x|^{ \f{1}{2}} + |x_1|^{- \f{1}{2}} ( |x|^{-1/6} + |x|^{-5/2} )$.

For the H\"older estimate, we use Lemma \ref{lem:E_vel2}. Recall $G(x)$ from Lemma \ref{lem:E_vel2} and $\chi_0$ from \eqref{eq:vel_cutoff_bd}, which satisfies $
\chi_0(x, y) = 1$ near $0$. Firstly, we rewrite $I$ as follows 
\[
I = (\Psi_{xy} - \Psi_{xy}(0) - \chi_0 \pa_{xy} G(x) \Psi_{xxxy}(0)) (r^{-1/2} - 1) \bar \th_{0, x}
+ \chi_0 \Psi_{xxxy}(0) \pa_{xy} G(x) (r^{-1/2} - 1) \bar \th_{0, x} \teq I_1 + I_2.
\]

Denote 
\[
\Psi_{xy, A} \teq \Psi_{xy} - \Psi_{xy}(0) - \chi_0 \pa_{xy} G(x) \Psi_{xxxy}(0).
\]

The estimate of $I_2$ is simple since the coefficient vanishes to $O(|x|^4)$ and we obtain a small factor: $|r^{-1/2} - 1| \les C_l y \les C_l S$ within the support of the solution. In particular, we have 
\[
|| I_2 \psi_2 ||_{C^{1/2}} \les  C_l S  | \Psi_{xxxy}(0)| \les C_l S \cdot \xxa{\Om}
\les C_l S \xxa{\om}.
\]

For $I_1$, we first rewrite $I_1 \psi_2$ as follows 
\[
I_1 \psi_2 = (\Psi_{xy, A} \psi_1)  \f{ \psi_2}{\psi_1} (r^{-1/2} - 1) \bar \th_{0,x}. 
\]
From Lemmas \ref{lem:E_vel}, \ref{lem:E_vel2}, we have $\Psi_{xy, A} \psi_1 \in L^{\inf} \cap C^{1/2}$. The coefficient $ \f{ \psi_2}{\psi_1} (r^{-1/2} - 1) \bar \th_{0,x}$ vanishes 
$O(|x|^{3/2})$ near $x=0$. 
Within $\supp( \bar \th_{0, x})$, since $|x| \leq 2 \nu  \leq 2 S$ by \eqref{eq:supp_para}, applying
\[
| \f{ \psi_2}{\psi_1} (r^{-1/2} - 1) \bar \th_{0,x} | \les C_l S,
\quad | \na (\f{ \psi_2}{\psi_1} (r^{-1/2} - 1) \bar \th_{0,x} ) | \les C_l S ,
\]
yields H\"older estimates for $I_1 \psi_2, I \psi_2$. We estimate other terms in $\uu_R$ similarly. We prove \eqref{eq:E_lin_lower}.

\subsubsection{Estimates of the lower order terms in the residual error}\label{sec:E_err_lower}

In this section, we estimate the lower order terms in the residual error \eqref{eq:E_err_MR} and show that 
\beq\label{eq:E_err_lower}
\bga
|| \bar \cF_{R, i} -  c_i  f_{\chi, i}  ||_{X_i}  \les C_l S + 
|| \Om_{\nu, R}||_{X_1}, \\
(c_1, c_2 , c_3 ) = ( \pa_{xy} \bar \cF_{R, 1}(0),  \, \pa_{x y} \bar \cF_{R, 2}(0),   \, 
\pa_{xx} \bar \cF_{R,3}(0)  ),
\ega
\eeq
where the norm $X_i$, $f_{\chi,i}$, $\Om_{\nu, R}$ are defined in \eqref{norm:Xi},  \eqref{eq:appr_near0}, \eqref{eq:elli_mainb}, respectively. 
We will further bound $ || \Om_{\nu, R}||_{X_1}$  and show that it is small 
using \eqref{eq:elli_mainb_err} and \eqref{eq:wb_decay2}.

We focus on the case $i=2$, i.e. the estimate of $\bar \cF_{R, 2} - c_2 f_{\chi, 2}$. Firstly, from \eqref{eq:cw_dif}, we have 
\[
\d \bar c_{\om} = \bar c_{\om} - \bar c_{\om, 2D} = \bar u_x(0) - \bar u_{2D, x}(0) = \d \bar u_x(0). 
\]

Note that $\bar u_x(0) + \bar v_y(0) = 0$ for the velocity \eqref{eq:E_ub}. A direct computation yields 
\[
c_1 = \pa_{xy}  \bar \cF_{R,1}(0) = \d \bar c_{\om} \bar \om_{xy}(0) + 4 C_l \bar \th_{xx}(0), 
\quad c_2 =  \d \bar c_{\om} \bar \th_{xxy}(0),  \quad c_3 = \d \bar c_{\om} \bar \th_{xxy}(0). 
\]

We can rewrite $\cF_{R, 2} - c_2 f_{\chi, 2}$ \eqref{eq:E_err_MR} as follows 
\[
\bal
I & = \pa_x( - \d \bar \uu \cdot \na \bar \th_0 + 2 \d \bar u_x(0) \bar \th_0 ) 
- c_2 f_{\chi, 2} \\
& = - (\d \bar u - \d \bar u_x(0) x) \bar \th_{0, xx} 
-  \pa_x (\d \bar u - \d \bar u_x(0) x) \bar \th_{0, x} 
- (\d \bar v - \d \bar v_y(0) y) \bar \th_{0, xy}  \\
& \quad -  \pa_x (\d \bar v - \d \bar v_y(0) y) \bar \th_{0, y}  
 + \d \bar u_x(0) ( \bar \th_{0, x} - x \bar \th_{0, xx} + y \bar \th_{0, x y} - \bar \th_{xxy}(0) f_{\chi, 2}  ) \\
 &  \teq I_1 + I_2 + I_3 + I_4 + I_5.
\eal
\]

\paragraph{\bf Estimate $I_5$}
For $I_5$, the coefficient is odd and $C^2$, has sufficiently fast decay, 
and vanishes $O(|\xx|^3)$ near $\xx=0$. Moreover, using \eqref{eq:elli_mainb}, \eqref{eq:elli_mainb2}, \eqref{eq:E_ub}, and Proposition \ref{prop:C2}, we have 
\[
 | \d \bar c_{\om} | = | \d \bar u_x(0)| \les || \Om_{\nu, R} ||_{X_1}.
\]
Thus, we can obtain
\[
|| I_5 ||_{X_2} \les || \Om_{\nu, R} ||_{X_1},
\]
where $X_2$ is the norm associated with the energy for $\eta$ and is defined in \eqref{norm:Xi}.

\paragraph{\bf Estimate $I_2$}
The estimates of $I_j , 1\leq j \leq 4$ are similar. We focus on the typical term in $I_2$
\beq\label{eq:E_err_lower_I2}
I_2 = -  \pa_x (\d \bar u - \d \bar u_x(0) x) \bar \th_{0, x} .
\eeq

 Recall $\bar \Psi$ from \eqref{eq:elli_mainb} and $\bar \phi_{2D} = (-\D)^{-1} \bar \om$.  Denote 
\[
\td \Psi = \bar \Psi -  \bar \phi_{2D}. 
\]

Recall the formula $\bar \uu = \na^{\perp} \bar \Psi +  \bar \uu_R$ from \eqref{eq:E_ub}. We have 
\[
\d \bar u - \d \bar u_x(0) x 
= - \B( \pa_y (  \bar \Psi -  \bar \phi_{2D}) -  \pa_{xy} (  \bar \Psi -  \bar \phi_{2D})(0) x\B) + \bar u_R
= - ( \pa_y \td \Psi  - \pa_{xy} \td \Psi(0) x) + \bar u_R,
\]
where the remainder $\bar u_R$ is defined in \eqref{eq:E_uR} with $\Psi, \Psi_2$ replaced by $\bar \Psi, \bar \Psi_2$. From \eqref{eq:elli_mainb2}, we have 
\[
-\D \td \Psi = \Om_{\nu, R}, \quad \Om_{\nu, R} \in X_1.
\]

Then the estimate of 
\[
\pa_x ( \pa_y \td \Psi  - \pa_{xy} \td \Psi(0) x) \bar \th_{0, x}
\]
in $I_2$ follows from the estimate of $\Psi -\Psi_{xy}(0) xy$ at the end of Section \ref{sec:E_lin_lower}, and we  obtain 
\[
|| \pa_x ( \pa_y \td \Psi  - \pa_{xy} \td \Psi(0) x) \bar \th_{0, x} ||_{X_2}
\les || \Om_{\nu, R}||_{X_1}.
\]

Other terms in $I_j, 1\leq j\leq 4$ related to $\td \Psi$ can be estimated similarly. 

For the term in $I_2$ \eqref{eq:E_err_lower_I2} related to $\bar u_R$, we have several terms due to the formula \eqref{eq:E_uR}, \eqref{eq:E_uR2}. The term involving $\bar \Psi_2$ is simple and its estimate follows from the estimate of $\Psi_2$ in Section \ref{sec:E_lin_lower}. For other terms, we estimate a typical term 
\[
 J = \pa_x( \bar \Psi_y - \bar \Psi_{xy}(0) x ) \cdot ( r^{-1/2} - 1)  \bar \th_{0, x}.
 \] 
Since $\bar \Psi$ is close to $\bar \phi_{2D}$, we use the decomposition $ \bar \Psi = \td \Psi + \bar \phi_{2D}$ and 
\[
J = \pa_x( \td \Psi_y - \td \Psi_{xy}(0) x ) \cdot ( r^{-1/2} - 1)  \bar \th_{0, x}
+ \pa_x( \bar \phi_{2D, y} - \bar \phi_{2D, xy}(0) x ) \cdot ( r^{-1/2} - 1)  \bar \th_{0, x}
\teq J_1 + J_2.
\]
The term $J_1$ follows from the above estimate. For $J_2$, note that $\bar \phi_{2D}$ satisfies the elliptic equation $-\D \bar \phi_{2D} = \bar \om$. 
Recall $\bar \om \in C^2$ with decay estimates \eqref{solu:decay}.

From $- \D \bar \phi_{2D, x} = \bar \om_x, -\D \bar \phi_{2D} = \bar \om$, 
using embedding \eqref{eq:elli_embed} and estimates \eqref{solu:decay}, we obtain 
\[
|\na^2 \bar \phi_{2D,x}| \les  1, \quad  | \na^2 \bar \phi_{2D}| \les 1.
\]

Using $\bar \phi_{2D,yyy} = -\bar \om_y - \bar \phi_{2D, xxy}$ and the above estimate, we yield $|\bar \phi_{2D, yyy}| \les 1$, and thus $|\na^3 \bar \phi_{2D}| \les 1$. Now, using the above estimate of $\bar \phi_{2D}$ 
and $|\bar \th_{0, x}| \les \min( |x_1|, |x|^{-3/5} )$ by \eqref{solu:decay}, and the smallness of $|r^{-1/2} - 1 | $ \eqref{eq:E_r1} within the support of the solution, we yield 
\[
\bal
 |J_2| & \les C_l  |x|\min( |x|, 1 ) \min( |x|^{-3/5}, |x_1|) \one_{|x| \leq S}, \\
 |\na J_2| & \les C_l  \min( |x|, 1 ) \min( |x|^{-3/5}, |x_1|) \one_{|x| \leq S} .
 \eal
\]
The term $J_2$ vanishes to order $|x_1| \cdot |x|^{2}$ near $x=0$. It follows the weighted $L^{\inf}$ estimate
\[
| J_2 \vp_2| \les C_l S,  \quad   | J_2 \vp_{g,2}| \les C_l S, \quad |J_2 \psi_2| \les C_l S \min(|x|^{1/2}, 1).
\]

Recall the weight $\psi_2 \asymp |x|^{-5/2} + |x|^{1/6}$ from \eqref{wg:hol}. We have
\[
|\na (J_2 \psi_2)| \les C_l S ( |x|^{-1/2} + 1 ).
\]

Combining the $L^{\inf}$ and $C^1$ estimates of $J_2 \psi_2$, we obtain the $C^{1/2}$ estimate of $J_2 \psi_2$ 
\footnote{ 
Suppose $|x| \leq |z|$. To estimate $|A(x) - A(z)| / |x-z|^{1/2}$ using $L^{\infty}, C^1$ bounds,
one can discuss two cases 
(a) $|x-z| \leq |x|/2$ and (b) $|x-z| \geq |x|/2$. Then apply $C^1$ bound in case (a) and $L^{\infty}$ bound in case (b).
}
\[
  || J_2 \psi_2 ||_{C^{1/2}} \les C_l S.
\]

Other terms follow similar estimates. We prove \eqref{eq:E_err_lower}.

\subsubsection{Modified finite-rank perturbation}\label{sec:3D_finite_rank}

Due to the difference of the operators between the 3D Euler \eqref{eq:Euler_lin} and the 2D Boussinesq \eqref{eq:lin}, we modify the decomposition \eqref{eq:bous_decoup2} 
for the perturbation $W= (\om, \eta, \xi) $ into $W_1 + \hat W_2$ and nonlinear perturbation $\cNF_i$ as follows 
\beq\label{eq:euler_decoup}
\bal
\pa_t W_{1, i}  & = (\cL_i -\cK_{1i} - \cK_{2 i}) (W_1) 
+ (\cL_i - \cL_{2D, i}) \hat W_2
+ \cN_i(W_1 + \hat  W_2) + \bar \cF_i  \\ 
& \quad - \cNF_i( W_1,  \hat W_2) - \cR_i(W_1, \hat W_2) ,\\
\pa_t \hat W_{2, i} & =  \cL_{2D,i }  \hat W_2 + \cK_{1i }(W_1) + \cK_{2i}(W_1) + \cNF_i(W_1,  \hat W_2) 
+ \cR_i( W_1, \hat W_2) , \\
\cNF_{i} &=  ( c_{\om}(\om) D_i^2 (W_1 + \hat W_2)(0)  + D_i^2 \big( \bar \cF_i   +  (\cL_i - \cL_{2D, i}  ) \hat W_2  \big)(0)  ) f_{\chi, i}
\eal
\eeq
where $D^2 = (\pa_{xy}, \pa_{xy}, \pa_{xx})$. Since the stream function $\Psi$ in \eqref{eq:E_lin_MR} is obtained from a modified source term $\Om$ \eqref{eq:elli_main}, 
we also modify the finite-rank operator $\cK_{2 i }$ \eqref{eq:appr_vel}, \eqref{eq:u_appr_1st}, \eqref{eq:u_appr_2nd} 
\footnote{
In $\cL_i(W_1), \cK_{1i}(W_1), \cK_{2i}(W_1)$ \eqref{eq:euler_decoup}, 
for a fixed $C_l$, 
we define $\Om(\om_1)$ and $\Psi(\om_1)$ linearly on $\om_1$ 
via the relation \eqref{eq:elli_main} and \eqref{eq:elli_main2} 
with $\om$ replaced by input $\om_1$. We first obtain $\phi_\dag(\om_1), \phi(\om_1)$ from Poisson equations \eqref{eq:euler4},\eqref{eq:E_elli},\eqref{eq:E_elli1} 
with source term $\om_1 = W_{1, 1}$ instead of $\om$.
We then obtain $Z_{i, j}(\om_1)$ via \eqref{eq:E_elli1}, and further construct $\Om(\om_1)$ from \eqref{eq:elli_main} with $\om$ replaced by $\om_1$ and these $Z_{i,j}(\om_1)$.
}
\[
\cK_{2i} (W_1) \teq \cK_{2D, 2i}( \Om(\om_1)).
\]
Note that we can still represent $\cK_{2 j}(W_1)$ as follows 
\beq\label{eq:euler_cK}
  \cK_{2 j}( W_1(t)) = \tts{ \sum_i } a_i(\Om(\om_1(t)) ) \bar f_{ij}, \quad j = 1,2,3 \, ,
\eeq
for some 
\emph{time-independent} functions $\bar f_{ij}$, and \emph{spatial-independent} factor $a_i(\Om( \om_1(t)))$ similar to \eqref{eq:bous_finite_rank1}. Thus, we can apply the same constructions of $\hat W_2$ and $\cR$ in Section \ref{sec:decoup_modi}, and use the same approximate space-time solution $\hat F_i, \hat F_{\chi, i}$ in \eqref{eq:bous_W2_appr}, \eqref{eq:bous_err_op}. Due to \eqref{wg:equi}, the linear modes $a_i(\Om(\om_1)(t))$ and $a_i(\om_1(t))$ satisfy almost the same estimate up to 
an error of size $C_l S || \om_1 ||_{X_1}$. Since $\hat W_2 = W - W_1$, we can rewrite  $  (\cL_i - \cL_{2D, i}  ) \hat W_2  \big)(0)$  as a linear functional on $W_1$ and $W$. 
\footnote{
Then the forcing term  
$\cK_{1i }(W_1) + \cK_{2i}(W_1) + \cNF_i(W_1,  \hat W_2)$ 
is \emph{linear} in $W_1$ and depends on $W = W_1 + \hat W_2$ quadratically.
Using this structure and the argument in Section \ref{sec:reg_W12}, 
we rewrite the $W_1$-equation in \eqref{eq:euler_decoup} 
as a \emph{linear} equation in $W_1$ with a forcing depending 
on the full solution $W$. 
This structure allows us to obtain local well-posedness 
following Section \ref{sec:reg_W12}.
}

Recall $E_*$ from Theorem \ref{thm:main_stab_full}. To control the nonlinear mode $a_{\mf{nl},i}$ in \eqref{eq:bous_W2_appr}, \eqref{eq:bous_err_op}, we modify the bootstrap condition \eqref{eq:W2_non_boot} 
\beq\label{eq:W2_non_boot3}
\bal
  |c_{\om } \om_{xy}(0) + \pa_{xy} \bar \cF_1(0)
+ \pa_{xy} ( \cL_1 - \cL_{2 D, 1} ) \hat W_2(0) | 
 & 
 < 5 \mu_6 E_* \, ,  \\
   |c_{\om } \th_{xxy}(0) + \pa_{xy} \bar \cF_2(0)
+ \pa_{xy} ( \cL_2 - \cL_{2 D, 2 } ) \hat W_2(0) |  & 
 < 10 \mu_6 E_* \, ,
\eal
\eeq

\begin{assumption}\label{boot_ass:E} 
We impose the bootstrap assumption on energy $E_{4,m}$ \eqref{energy4m}
\[
E_{4, m}(t) < E_*,  \ \forall \, t \geq 0.
\]
\end{assumption}

Since estimates in Section \ref{sec:rank1} for functionals 
$a_l(\cdot), a_{\mf{nl}, i}(\cdot)$ in  $\hat W_2$ 
\eqref{eq:bous_W2_appr} do not involve $\xi$-norm, 
under Assumption \ref{boot_ass:E}, using \eqref{eq:W2_res1},  
the support and vanishing properties in \eqref{eq:W2_reg}, we bound \
\beq\label{eq:euler_W2}
  \| \hat W_2 \|_{C^{4,1} } \les 1, 
  \quad   | \hat W_2(t, x) | \les |x|^2, 
  \quad  \supp(\hat W_2) 
  \subset [-L_2, L_2] \times [0, L_2].
\eeq

In Sections \ref{sec:E_xi_lower}, \ref{sec:euler_other_xi}, we estimate 
the difference in the energy estimates 
involving $\xi_1$ due to the change of weights and energies in \eqref{wg:xi_modi} and \eqref{energy4m}.

\subsubsection{$L^{\infty}$ estimate of $\xi_1$ for $|x| \geq \f{ \nu}{2}$}\label{sec:E_xi_lower}

We close the far-field weighted $L^\infty$ estimate of \(\xi_1\) caused
by modifying the weights in \eqref{wg:xi_modi}.  We use the full $\theta$-equation in \eqref{eq:euler4} rather than the linearized \(W_1\)-equation \eqref{eq:euler_decoup} to simplify presentation. We establish the following weighted $L^{\infty}$ estimates.

\begin{prop}\label{prop:xi_far_linfty}
Under Assumptions \ref{boot_ass:1}, \ref{boot_ass:2}, \ref{boot_ass:E}, 
for $\rho = \vp_{3,m}, \vp_{g, 3, m}$, 
there exists a large $\nu_{*,1}>0$, such that for any $\nu \geq \nu_{*,1}$
and $C_l S < \nu_{*,1}^{-1}$, we have the following stability estimates 
\beq
  \| \one_{|x |\geq \f{\nu}{2} }  \xi_1(t) \rho \|_{ L^{\infty} } 
  \leq \max( \sup_{s \leq t, |x| = \f{\nu}{2} }  | \xi_1(s, x) \rho(x) |,  \,
    \| \one_{|x |\geq \f{\nu}{2} }  \xi_1(0,\cdot) \rho(\cdot) \|_{ L^{\infty} },  
    C_{\eqref{eq:xi_far_linfty}} \nu^{-1/32}  )  .
    \label{eq:xi_far_linfty}
\eeq
for some absolute constant $C_{\eqref{eq:xi_far_linfty}}$.

\end{prop}

Using the weighted $L^{\infty}$ stability estimates of $\xi_1(t)$ in $|x| \leq \nu / 2$, 
which controls the boundary term $\sup_{s \leq t, |x| = \f{\nu}{2} }  | \xi_1(s, x) \rho(x) |$ and choosing small initial perturbation and large $\nu$, we control the far-field perturbation $\xi_1$ dynamically. 
See Section \ref{sec:euler_blowup_pf}.

 \begin{proof}
 Using \eqref{eq:euler4}, $c_{l,\tot} \equiv \bar c_l$ by \eqref{eq:E_normal}, 
and $c_{\th,\tot} = \bar c_l + 2 c_{\om,\tot}$, we derive 
\[
 \pa_t \pa_y \th_{\tot} + ( \bar c_l \xx + \uu_{\tot}) \cdot \na   \pa_y \th_{\tot} 
= ( c_{\th, \tot} - \bar c_l ) \pa_y \th_{\tot}
- \pa_y \uu_{\tot}  \cdot \na \th_{\tot}
= 2 c_{\om, \tot} \pa_y \th_{\tot} - \pa_y \uu_{\tot}  \cdot \na \th_{\tot}.
\]

By \eqref{eq:3D_thy_van}, we obtain $\pa_y \th_{\tot}(0, y) =0 $ for $y \leq \nu/2$. 
For $|x | \geq \nu / 2$, decomposing $\xi_\tot = \pa_y \bar \th_0 + \xi$, using $\hat \xi_2 =0$ ,and multiplying by a weight $\rho$, we obtain 
\beq\label{eq:xi_far}
 \pa_t ( \xi_1  \rho) +  ( \bar c_l \xx + \uu_{\tot}) \cdot \na ( \xi_1 \rho)
= D(\rho) \cdot ( \xi_1 \rho)
+ \rho G(x) \, ,
\eeq
where $D(\rho), G$ denote the damping term and the force
\beq\label{eq:xi_far_dp}
\quad D(\rho) =  2 c_{\om, \tot} + \rho^{-1} ( \bar c_l \xx + \uu_{\tot}) \cdot \na \rho ,
\quad  G = 
2 c_{\om, \tot} \pa_y \bar \th_0 - (\bar c_l \xx + \uu_{\tot} ) \cdot \na \pa_y \bar \th_0
-   \pa_y \uu_{\tot}  \cdot \na \th_{\tot} , 
\eeq
 and $\uu_{\tot}= \uu(\om_{\tot})$ is the velocity associated  with the full solution $\om_{\tot} = \bar \om_0 + \om_1 + \hat \om_2$.

\vs{0.1in}
\paragraph{\bf Damping coefficients $D(\rho)$}

Due to symmetry in $x_1$, we consider $x_1 , x_2 \geq 0$. 
For $|x| \geq \nu/2$, by definition \eqref{wg:xi_modi}, we have 
$\rho = b |x|^{\ell}$ for $(\rho,  \ell) = (\vp_{3, m}, \f17)$ or $( \vp_{g,3,m}, \al_{g, n})$. 
We compute
\beq\label{eq:xi_dp_far1}
 \f{ x \cdot \na \rho}{\rho} = \ell, \quad  |\f{x_i \pa_i \rho}{\rho}| = |\f{x_i^2}{|x|^2} \ell  | \les 1.
\eeq

Since $\om_{\tot} = \bar \om_0 + \om_1 + \hat \om_2$, 
and $\bar \om_0, \hat \om_2$ vanishes $O(x_1)$ near $x=0$, under Assumption \ref{boot_ass:E}, using  estimates \eqref{eq:euler_W2} for $\hat \om_2$,
energy $E_{4,m}$ \eqref{energy4m} for $\om_1$, decay estimates \eqref{solu:decay} for $\bar \om_0, 
\bar \phi^N$, and estimates \eqref{eq:xi_Psi_decay0}  for $f = \bar \om_0, \om_1, \hat \om_2, \om_{\tot},
-\D \bar \phi^N$, 
for $|x| \leq S$, we bound
\beq\label{eq:xi_Psi_decay}
   |\na \uu(f) | 
   +     | \tf{1}{x_1}  u(f) | +      | \tf{1}{x_2}  v(f)   |  
 \les  \la x \ra^{-1/32} \big\| \tf{|x|^2}{\la x \ra^{2} } f \big\|_{X_1}
 \les \la x \ra^{-1/32} ( 1  + \xxa{\om_1} ) \les \la x \ra^{-1/32}.
\eeq

Using \eqref{eq:xi_dp_far1} and \eqref{eq:xi_Psi_decay},
for $\f{\nu}{2} \leq |x| \leq S$ with $S$ being the support of the solution, we obtain 
 \[
\B|\f{  \uu(\om_{\tot}) \cdot \na \rho}{\rho} \B| 
  \leq \la x \ra^{-1/32}
  ( | \f{x_1 \pa_1 \rho}{\rho} | + | \f{x_2 \pa_2 \rho}{\rho} | ) \les  \nu^{-1/32},
  \quad - \f{ ( \bar c_l \xx + \uu_{\tot} ) \cdot \na \rho }{ \rho}
 = - \bar c_l \ell +  O( \nu^{-1/32}).
  \]

Using definition of $D(\rho)$ in \eqref{eq:xi_far}, we prove 
\beq\label{eq:xi_damp_err}
\bal
- D(\rho) \geq -\bar c_l  \ell  -2 c_{\om, \tot} - C \nu^{-1/32}  ,
\quad | D(\rho) | \les 1.
 \eal
\eeq

By \eqref{eq:E_normal} and \eqref{eq:E_u}, we obtain $c_{\om,\tot} = \f{1}{2} \bar c_l + \uu_x(\om + \bar \om_0)(0)= \f{1}{2} \bar c_l - \Psi_{xy}( \om + \bar \om_0)(0)$. 
Denote $\phi_{2D}(f) =(-\D)^{-1}(f)$.
Recall the 2D factor $\bar c_{\om, 2D}(\bar \om) + c_{\om, 2D}(\om) = \f12 \bar c_l - \phi_{2D, xy}(\bar \om + \om)$ from \eqref{eq:normal}. 
Let $f = \om + \bar \om_0$. Using estimate \eqref{eq:u_est_direct} and \eqref{wg:equi}, we bound 
\[
|c_{\om, \tot} -\bar c_{\om, 2D}(\bar \om) - c_{\om, 2D}(\om)|
\leq | \Psi_{xy}(f) -\phi_{2D, xy}( f)| 
+  |\phi_{2D, xy}(\bar \om_0 - \bar \om)(0)| \teq I_1 + I_2
\]

Since $-\D \Psi(f) = \Om(f)$ with \eqref{eq:elli_main}, \eqref{eq:elli_main2} 
and $\om=f$, using \eqref{wg:equi}, $X_1$ norm \eqref{norm:X} for $\om_1$,
\eqref{eq:euler_W2} for $\hat \om_2$, \eqref{solu:decay}, \eqref{eq:profile2} for $\bar \om_0, \bar \om$, and $\phi_{2D, xy}(f)(0) = C \int \f{y_1 y_2}{|y|^4}f(y) dy $ , we bound 
\[
I_1 \les \| \tf{|x|^2}{\la x \ra^2} (\Om(f) - f) \|_{X_1}
\les C_l  S \xxa{\tf{|x|^2}{\la x \ra^2} f } \les C_l  S ,
\quad | I_2 | \les \nu^{-1/32}.
\]

Since $c_{\om}(\bar \om + \om)
= c_{\om}(-\D \bar \phi^N) + c_{\om}( \bar \e + \om)
= \bar c_{\om,2D}^N + U_{2D, x}(0)$ by \eqref{eq:wb_vel_err}, \eqref{eq:U}, 
under Assumption \ref{boot_ass:E}, using  the above estimates, 
estimate \eqref{eq:U_bd} for $U_x(0)$, and \eqref{stab:infty}, we bound 
\[
\bal
- D(\rho) & \geq - \bar c_l \ell 
  - 2 \bar c_{\om,2D}^N - 2  U_{2D, x}(0)
-   C( C_l S + \nu^{-1/32}) \\
& \geq - \bar c_l \ell 
  - 2 \bar c_{\om,2D}^N - 2 ( \mu_6 E_* + | u_{x, 2D}(\bar \e) | ) -   C( C_l S + \nu^{-1/32}) 
  \geq \lam -   C( C_l S + \nu^{-1/32}).
  \eal
\]
where $\lam>0$ is an absolute constant independent of $\nu, C_l$. By requiring $\nu^{-1}, C_l S$ small, we prove 
\beq
- D(\rho) > C_{ \eqref{eq:xi_far_damp} } > 0 , \quad \forall \,  t \geq 0. 
\label{eq:xi_far_damp}
\eeq
with constant $C_{ \eqref{eq:xi_far_damp} } $ independent of $\nu, C_l, S(0)$.

\vs{0.1in}

\paragraph{\bf Forcing in $|x| \geq \nu/2$}

By \eqref{wg:xi_modi} and \eqref{wg:linf_decay} for $\vp_2, \vp_3$,
\eqref{wg:linf_grow} for $\vp_{g, 2}, \vp_{g, 3}$, we obtain $\vp_{3, m} \les \vp_3 \asymp \vp_2$, $\vp_{g,3,m} \asymp \vp_{g, 3} \asymp \vp_{g, 2}$. 
Note that $\na \th_{\tot} =  \na \bar \th_0 + (\eta_1, \xi_1)$. 
Due to the vanishing of $\pa_y \bar \th_0(0, y)$ on $y \in [0, \nu]$ \eqref{eq:3D_thy_van},
the odd symmetry of $\pa_x \bar \th_0$, \eqref{solu:decay}, and \eqref{eq:profile2}
for $|x| \geq \nu/2$,  we obtain 
\beq\label{eq:euler_thy_0}
\bal
|\na^{1+i} \bar \th_0 |  ( \vp_{3, m}(x) + \vp_{g, 3, m} )
 & \les \la x \ra^{2 \bar \al_1 - i + \al_{g, n}}    \les \nu^{-1/10-i}, \quad i = 0, 1. \\
 \eal
\eeq

Using  energy $E_{4, m}$ \eqref{energy4m}, norms \eqref{norm:Xi}, 
the estimate of $\uu$ in \eqref{eq:xi_Psi_decay}, under Assumptions \ref{boot_ass:1}-\eqref{boot_ass:E}, 
for $\rho = \vp_{3,m}, \vp_{g, 3, m}$, 
for $|x| \geq \nu/2$,  we bound $G$ in \eqref{eq:xi_far_dp}
\beq\label{eq:xi_far_force}
\rho |G| 
\les \rho ( |\xx| \cdot | \na^2 \bar \th_0| + |\na \bar \th_0| ) 
+ |\na \uu_{\tot}(x)| ( \nlinf{ \eta_1 \rho} + \nlinf{ \xi_1 \rho} + \rho |\na \bar \th_0 |)
   \les \la x \ra^{-1/32}
   \les \nu^{-1/32}.
\eeq

Estimating $\xi_1 \rho$ in \eqref{eq:xi_far} along the trajectory $X(t, x)$ associated with $\bar c_l \xx + \uu_{\tot}$, 
we obtain 
\[
 \tf{d}{dt} ( \xi_1 \rho) \circ X(t, x) =  ( D(\rho)   \xi_1 \rho)  \circ X(t, x) 
 +   ( \rho  G ) \circ X(t, x). 
\]

Fixed $X_0$ with $ |X_0 | \geq \nu / 2$ and $t \geq 0$. Let $X(t, x_0) = X_0$ and $t_0 \in [0, t]$ be the first time where the trajectory reaches $ |X(t_0, x_0)| = \nu/2$, if exists, and  $t_0 = 0$ otherwise. For $\rho \in \{ \vp_{3,m}, \vp_{g, 3, m} \} $, using estimate \eqref{eq:xi_far_force}, and \eqref{eq:xi_far_damp}, and Gronwall's inequality, we obtain
\[
\bal
|(\xi_1 \rho)(t,  X_0)| & \leq e^{- C_{ \eqref{eq:xi_far_damp} }  (t- t_0)}  
| (\xi_1 \rho)(X(t_0, x_0), t_0)| + C \int_{t_0}^t e^{- C_{ \eqref{eq:xi_far_damp} }  (t- s)} \nu^{-1/32}  d s  \\
& \leq \max(  | (\xi_1 \rho)(X(t_0, x_0), t_0)|, C \nu^{-1/32}).
\eal 
\]

The term $ | (\xi_1 \rho)(X(t_0, x_0), t_0)|$ is bounded by the boundary value on $|x| = \nu/2$ 
or the initial data, as in \eqref{eq:xi_far_linfty}. Taking supremum over $|x| \geq \nu/2$, 
 we prove Proposition \ref{prop:xi_far_linfty}.
 \end{proof}

\subsubsection{Other estimates involving $\xi_1$}\label{sec:euler_other_xi}

We perform weighted $L^{\infty}$ estimates on $\om_1, \eta_1$ 
and $C^{1/2}$ estimates  as in the 2D Boussinesq case. The change of weights for $\xi$ 
and $ \pa_y \bar \th_0$ also affects the $C^{1/2}$ estimates for $\xi_1, \eta_1$ 
and the weighted $L^{\infty}$ estimates for $\eta$ in the region $|x| \geq \nu/2$. We show that these changes introduce only small perturbative errors.

\vs{0.05in}
\paragraph{\bf Difference in weighted $L^{\infty}$ estimates for $\eta$}

Variables  $\xi, \pa_y \bar \th_0$ enters the $\eta$-equation via the 
terms $- \bar v_x \xi, - v_x \pa_y \bar \th_0 $ in \eqref{eq:Euler_lin}, $ - v_x \xi$ in the nonlinear term $\cN_2$ \eqref{eq:bous_non},  \eqref{eq:E_non_MR}, and $\bar v_x \pa_y \bar  \th_0$ in 
$\bar \cF_i$ \eqref{eq:bous_err}, \eqref{eq:Euler_lin}, where $\bar v_x, v_x$ are the 3D velocity associated with $\bar \om_0, \om$.

Consider $(\rho_2, \rho_3,  \rho_{3, m})  = (\vp_2, \vp_3, \vp_{3,m})$ or $( \vp_{g, 2}, \vp_{g, 3} , \vp_{g, 3, m})$. 
Since $\rho_3 = \rho_{3, m}$ for $|x| < \nu/2$ by \eqref{wg:xi_modi}, the estimate of these two terms are the same for $|x| \leq \nu/2$.

For $|x| \geq |\nu| /2 > 1$, by definition of weights in \eqref{wg:linf_decay}, \eqref{wg:linf_grow}, 
and $\psi_1$ in \eqref{wg:hol}, we obtain $ \rho_2 \asymp \rho_3$,
\[
   \rho_2 \les \rho_{3,m} + |x_1|^{-1/2} |x|^{-1/6}
   \les \rho_{3,m} + |x_1|^{-1/2} \psi_1, \quad \rho_{3, m} \gtr |x|^{1/4}. 
\]

For $|x| \geq \nu/2$, $|x_1| \leq 1$, since $|(0, x_2) | \geq \nu / 4$ 
and $v_x(0, x_2 )= 0$,  applying \eqref{eq:E_du_hol:X1}, we bound
\[
| |x_1^{-1/2} \psi_1  v_x |
\les | \d( \psi_1 v_x, x, (0, x_2) ) | \cdot |x_1|^{-1/2}
  \les \| \om_1 \|_{X_1} .
\] 

Applying estimates \eqref{eq:E_du_hol:reg} with $f = \bar \om_0$ \eqref{eq:profile2}, $-\D \bar \phi^N, \hat \om_2$,  and using \eqref{solu:decay} for $\bar \om , \bar \phi^N$, and \eqref{eq:euler_W2} for $\hat W_2$, for $|x|, |z| \leq |S|$, we bound 
\beq\label{eq:E_bar_hol}
 \| \na \uu( f) \|_{C^{1/2} \{ |x| \leq |S| \} }  \les 1, \quad  
f = \bar \om_0,  \ \bar \om, \ -\D \bar \phi^N, \ \hat \om_2.
\eeq

Using \eqref{eq:E_bar_hol} with $f = \bar \om_0$, we get $| x_1^{-1/2} \psi_1  \bar v_x | \les  1$ for $ |x| \in [\nu/2, |S| ]$. 
Thus, for $|x| \geq \nu / 2$, Under Assumption \ref{boot_ass:E}, using above estimates
and  \eqref{eq:xi_Psi_decay} with $f = \om_1, \hat \om_2, \bar \om_0$, we bound 
\[
\bal
   \rho_2 (|v_x| + |\bar v_x|) ( |\xi| & + |\pa_y \bar \th_0| ) 
  \les (|v_x| + |\bar v_x|) \cdot \rho_{3, m}( |\xi| + |\pa_y \bar \th_0| ) 
 +  |x_1|^{-1/2} \psi_1 (  |v_x|  + |\bar v_x| )  \cdot  ( |\xi| + |\pa_y \bar \th_0| )  \\
 &  \les ( |v_x | + |\bar v_x| ) \cdot ( \| \xi \|_{X_3} + 1 )
  + 
   ( |x|^{-1/4}  \nlinf{ \xi \rho_{3,m} } 
  +  \nu^{-1/10} ) \les \la x \ra^{-1/32} \les \nu^{-1/32}.
\eal
\]

\vs{0.05in}
\paragraph{\bf ODE estimates for $\eta$}

In the ODEs of $\eta$ (3D Euler version of \eqref{eq:lin_main_cw}), 
we estimate the $\xi_1$-terms by their weighted $L^{\infty}$-norm. 
Similar to the above estimates, these estimates differ by $O(\nu^{-1/32})$.

\vs{0.05in}
\paragraph{\bf Difference in weighted $C^{1/2}$ estimates for $\xi_1, \eta_1$}

Following Section \ref{sec:comb_vel_err} and the beginning of Section \ref{sec:non}, we derive 
the weighted modified decomposition of the $\xi_1$-equation 
\[
\bal
 \pa_t (\xi_1 \psi_3) & + ( \bar c_l \xx + \uu(\bar \om_0 + \om) ) 
 \cdot \na (\xi_1 \psi_3)   =   d_{3, L}^{\numm}(\psi_3)  \cdot (\xi_1 \psi_3)
 + B_{\mf{modi}, 3} \psi_3  \\
 & \quad +  \td \cN_3(\psi_3) +  ( \bar \cF_{loc, 3} - \cR_{loc, 3} ) \cdot \psi_3
+ ( (\cL_{3} - \cL_{2D, 3}) \hat W_2 - \pa_x^2 ( (\cL_{3} - \cL_{2D, 3}) \hat W_2 )(0) f_{\chi, 3} ) 
\cdot \psi_3 .
\eal 
\]
Let $\bar \phi^N$ be the approximate stream function constructed  for 2D Boussinesq.
Here,  we define 
\[
\bar \uu^N =  \uu( (-\D) \bar \phi^N), \quad \UU = \uu(\om + \bar \e ),
\quad  \bar \e = \bar \om_0 -(-\D) \bar \phi^N  .
\]
similar to \eqref{eq:non_U},\eqref{eq:wb_vel_err}.  Note that the velocity $\uu(f)$ is \emph{the 3D velocity} for $f$ defined via \eqref{eq:euler4}.  As in \eqref{eq:tran_extra},\eqref{eq:dp_num}, we define  $\cT_{d, N}$ and damping coefficient
$d_{3, L}^{\numm}$:
\beq\label{def:xi_dp_far}
d_{3,L}^{\numm}(\psi_3)   =  2 \bar c_{\om}^N + \bar u_x^N + 
\cT_d^{\numm}, \  \cT_d^{\numm} \teq  \psi_3^{-1} ( (\bar c_l \xx + \bar \uu^N) \cdot \na \psi_3 ), \  
 \cT_{d, N}( \psi_3) =  \psi_3^{-1} ( \UU \cdot \na \psi_3 ),
\eeq

We define $\td \cN_3(\psi_3), B_{\mf{modi},3}, \bar \cF_{loc, 3}, \cR_{loc, 3}$ following Section \ref{sec:comb_vel_err}. 
The term $(\cL_{3}-\cL_{2D,3})\hat W_2$ comes from the same term in the $W_1$-equation in \eqref{eq:euler_decoup}, while $f_{\chi,3}$-term comes from $\cNF_3$. We show that the $\hat W_2$-term is small in Section \ref{sec:E_lin_M}.

 Since the weight $ \vp_3(y) = \vp_{3,m}(y), \vp_{g,3}(y) = \vp_{g, 3,m}(y)$ for $|y| < |\nu|/2$   by \eqref{wg:xi_modi}, in the $C^{1/2}$ estimates for $x, z \in \R^2_{++} \cap \supp(\xi)$, 
the changes of weights in \eqref{wg:xi_modi} only affects the estimates 
of $\xi$-terms for $\max( |x|, |z| ) \geq \nu / 2$. 
These $\xi$-terms include :
\bseq\label{eq:E_hol_xi_damp}
\beq
d_{3, L}^{\numm}(\psi_3)  \cdot (\xi_1 \psi_3), 
 \quad    \td \cN_3(\psi_3)  \eqref{eq:non_dec1} : ( 3 U_x(0) - \td V_y) \xi_1 
 = ( 2 U_x(0) - V_y ) \xi_1 ,
 \   \cT_{d, N}(\psi_3) \xi_1.
\eeq

Below, we assume $|x|\leq |z|$ and $|z| \geq \nu / 2$, and simplify  $ d_{3, L}^{\numm} =  d_{3, L}^{\numm}(\psi_3)$. In the weighted $C^{1/2}$ estimates for $\eta_1$, we need to estimate 
\beq
  B_{\mf{modi,2}} \, \eqref{eq:bad_lin}:  - \bar v_x^N  \xi_1,  \quad N_{W_1, 2}, \td \cN_2(\psi_2) 
\eqref{def:N_W1}, \eqref{eq:non_dec1} :  - V_x \xi_1 
\eeq
\eseq

Below, under Assumptions \ref{boot_ass:1}, \ref{boot_ass:2}, \ref{boot_ass:E}, we show that the 
difference of the estimates for the above terms are perturbative.

\vs{0.05in}
\paragraph{\bf Damping coefficients}

Following the derivations in \eqref{eq:hol_est1}, \eqref{eq:hol_nota}, 
we have two  damping terms 
\[
I_1 \teq \d( d_{3, L}^{\numm}(\psi_3)  J_3 , x, z) g_3(x-z),  \quad I_2 \teq \f{ ( b_{\UU}(x) - b_{\UU}(z) ) \cdot \na g_3(x- z)}{ g_3(x-z)} \cdot \d( J_3 , x, z) g_3(x-z), 
\quad J_3 \teq \xi_1 \psi_3, 
\]
where $\d(f, x, z) \teq f(x) - f(z)$. The reader should not confuse it with $\d f$ in \eqref{eq:dif_op}. The $I_2$-term is a damping term appearing in \eqref{stab:hol} and its estimate \emph{does not involve} 
 $L^{\infty}$-norm of $\xi_1$.

For $(p, q) = (x, z)$ or $(z, x)$, we decompose $I_1$ as
\[
  I_1 =  d_{3, L}^{\numm}(\psi_3) (q) \cdot \d(J_3, x, z) g_3 + \d(d_{3, L}^{\numm}(\psi_3) , x, z) 
\cdot  J_3(p) g_3  \teq I_{1,1} + I_{1, 2}
\]
where  $I_{1,1}$ is a damping term, and we bound $|I_{1,2}|$ perturbatively by the energy in
Section \ref{sec:EE_hol_bad}.

Using $ | x_i \pa_i  \psi_3 | \les \psi_3$ by definition \eqref{wg:hol} and using \eqref{eq:xi_Psi_decay} 
with $f= \om_1, \hat \om_2, \bar \om_0, -\D \bar \phi^N$,  we bound 
\beq\label{eq:xi_U_decay}
   |d_{3,L}^{\numm}|\les 1 + |\tf{\bar u^N}{x_1}| + |\tf{\bar v^N}{x_2}| \les 1,  \quad 
  | \na V| + | \na \bar \uu^N|  \les \la x \ra^{- 1/32},
 \quad | \cT_{d, N} |  \les |\tf{U}{x_1}| + |\tf{ V }{x_2}| \les 1. 
\eeq

If $|x| \leq |z|/2$, using \eqref{eq:xi_U_decay} , 
$|x-z| \gtr |z| \gtr \nu$, and  \eqref{eq:xi_Psi_decay}: 
we bound 
\[
\bal
  | \d(  d_{3, L}^{\numm}  ) | g_3(|x-z|) \les |x-z|^{-1/2} \les \nu^{-1/2}.
\eal
\]
Since $\psi_3 \les \vp_{g,3,m}$ by \eqref{wg:hol}, \eqref{wg:xi_modi}, we bound $|I_{1,2}| \les \nu^{-1/2} \nlinf{ \xi \vp_{g, 3, m} }$ . 

Next, we consider $  |z|/ 2 \leq  |x|,|z|\leq S$. By \eqref{wg:hol}, we expand $\psi_3(y) = \sum a_i |y|^{\ell_i}$ with $\vec{\ell} = (-\f52, -1,-\f12, \f16)$.
Using estimate \eqref{eq:xi_Psi_decay} for $\bar \uu^N$,
and \eqref{def:xi_dp_far},
for $|y| \gtr \nu > 1$, we get 
\[
\f{y \cdot \na \psi_3}{\psi_3}
 = \f16 + \f{\sum_{\ell_i< 1/6} a_i (\ell_i - \f16 )|y|^{\ell_i} }{ \sum a_i |y|^{\ell_i} }, 
\quad 
| \na \cT_{d}^{\numm}(y) | 
\les |y|^{-1} |\na  \bar \uu^N | + |y|^{-2} |\bar \uu^N |
+ |y|^{-1 -\f12 - \f16} \les |y|^{-1 - \f{1}{32}}.
\]

Using \eqref{eq:E_bar_hol} with $f =-\D \bar \phi^N$ for $\bar \uu^N$, 
for $|x|, |z| \in [\nu/4, |S|]$, we obtain
\[
 |\d( d_{3, L}^{\numm},x, z)| g_3(x-z)
 \les [ \na \bar \uu^N ]_{C^{1/2}}
 +   |x-z|^{ \f12 }  \sup_{ |y| \in [ \f14 |z|, 2 |z|]} | \na \cT_{d}^{\numm}(y) |
 \les 1 +  \sup_{ |y| \in [  |z|/4, 2 |z|]}  |y|^{ \f12- 1- \f{1}{32} } \les 1.
\]

Since $\psi_3 \les |x|^{-1/6} \psi_{g,3,m}$ by \eqref{wg:hol}, \eqref{wg:linf_grow}, 
and $|x|, |z| \geq \nu/4$, we bound 
\[
|I_{1, 2} | \les  \max(| \xi_1 \psi_3|(x), | \xi_1 \psi_3|(z) ) \les \nu^{-1/6} E_{4, m}
\les \nu^{-1/6}.
\]

Thus, in both cases, we get $ |I_{1,2} | \les \nu^{-1/6}$. Moreover, we have the maximal damping
\[
 II \cdot   \d(J_3, x, z) g_3, \quad  II  = \min( d_{3, L}^{\numm}(x), \, d_{3, L}^{\numm}(z) ) ,
\]
with $- II \geq  -d_{3 ,L}^{\numm}(p_{x, z})$ in  \eqref{stab:hol} up to the difference between $\bar \uu^N$ and $\bar \uu_{2D}^N$.

\vs{0.05in}
\paragraph{\bf Estimate of forcing}

For the $C^{1/2}$ estimate of $U_x(0)$-term \eqref{eq:E_hol_xi_damp} in $\xi_1$-equation, we get 
\[
\mu_{h, 3} |\d( 2 U_x(0) \xi_1 \psi_3,x, z ) g_3| \leq  2 |U_x(0) | \mu_{h, 3} [ \xi_1 \psi_3 ]_{C_{g_3}^{1/2}} 
\leq 2 |U_x(0)| E_{4, m}, 
\]
which \emph{does not}  involve $L^{\infty}$-norm and \emph{improves} the estimate $3 |U_x(0)| E_{4, m}$ in $|N_{W_1, 3 ,\mf{hol}}|$ in \eqref{hol:nlin:b2}.

Below, we simplify $\d(f, x, z) $ as $\d(f)$.  If $|x| \leq |z|/2$, we obtain $|x-z| \geq |z|/2 \gtr \nu$. Using \eqref{eq:xi_U_decay},  $\psi_2 = \psi_3 \les \la x \ra^{-1/6} \psi_{g, 3, m}, g_2 = g_3$ by \eqref{wg:hol}, \eqref{wg:linf_grow}, 
we bound 
\[
 (| \d(   V_y  \xi_1 \psi_2)| + 
  | \d(  \bar v_x^N \xi_1 \psi_2)|
  +   | \d(  V_x \xi_1 \psi_2)| ) g_3(x-z) 
  \les \nlinf{ \xi_1 \vp_{g,3,m} } |x-z|^{-1/2} \les \nu^{-1/2}.
\]

If $ |z|/2 \leq |x| \leq |z| \leq S$, $\cT_{d, N}(\psi_3)$ in \eqref{eq:E_hol_xi_damp} is not singular. We focus on weighted $C^{1/2}$ estimate for a typical term $-\td V_y \xi$ in \eqref{eq:E_hol_xi_damp} :
\[
 |\d(  V_y \psi_1 \cdot \xi \psi_2 \cdot \psi_1^{-1} ) |
 \les  |\d(  V_y \psi_1 )|  \cdot |\xi \psi_2 \cdot \psi_1^{-1}(x) |
 +|  V_y \psi_1(z) | \cdot ( | \d( \xi \psi_2) |  \cdot  \psi_1^{-1}(x) 
+ |\xi \psi_2(z) | \cdot | \d( \psi_1^{-1} ) | ) .
\]
 Using  \eqref{eq:E_du_hol:X1} with $f =  \om_1$ and \eqref{eq:E_bar_hol} 
 for  $C^{1/2}$ estimates,
 \[
\psi_1(x) \asymp \psi_1(z) , 
 \quad | \d(\psi_1^{-1},x,z) | \les |x-z| \cdot  |x|^{-1} \psi_1^{-1}(x) , 
 \]
by $|x| \in [|z|/2, |z|]$,  and  using $\| \xi_1 \|_{X_3} \les 1$ by Assumption \ref{boot_ass:E}, we bound 
 \[
 | \d(  V_y \xi \psi_2 ) |
 \les\| \xi \|_{X_3} ( |x-z|^{1/2}  \f{\psi_2 }{\psi_1 \vp_{g, 3, m}}
 + |x-z|^{1/2}   |V_y| 
 + |V_y| \cdot |x-z| \cdot |x|^{-1} \f{\psi_2}{ \vp_{g, 3, m}} ) .
 \]
 Since  $ \psi_2 \psi_1^{-1} \les \la x \ra ^{1/3 - \al_{g, n} } \vp_{g, 3,m}$, 
   $\psi_2 \les \vp_{g,3,m}$ by  \eqref{wg:hol}, \eqref{wg:linf_grow},
   using  \eqref{eq:xi_Psi_decay} and $|x|\asymp |z| \gtr \nu$ : 
   \beq
 | \d(  V_y \xi \psi_2 ) | \cdot |x-z|^{-1/2}
 \les \la x \ra^{1/3 - \al_{g, n}} + \la x \ra^{-1/32} \les \nu^{ -\g_{ \eqref{eq:xi_force_lot} } },  \quad \g_{ \eqref{eq:xi_force_lot} }= -1/3 + \al_{g, n}  >0.
 \label{eq:xi_force_lot}
   \eeq

Similarly, we prove that other forcing terms are bounded by $C \nu^{ -\g_{ \eqref{eq:xi_force_lot}  } }$.

\subsubsection{Comparison between the estimates}\label{sec:E_lin_M}

In this section,  under Assumption \ref{boot_ass:1}, \ref{boot_ass:2}, \ref{boot_ass:E}, we show that the differences between the energy estimates for 3D Euler  \eqref{eq:E_lin_MR}, \eqref{eq:E_non_MR}, \eqref{eq:euler_decoup}  and 2D 
Boussinesq  \eqref{eq:lin}, \eqref{eq:bous_decoup2}  are perturbative. We have estimated the lower order operators in Section \ref{sec:E_lin_lower}, \ref{sec:E_err_lower}. Here, we only focus on the main terms. 
By \eqref{eq:chi_Ri} and Assumption \ref{boot_ass:1}, we get $\chi_4 \chi_{\nu} = \chi_{\nu}$. Recall that we perform energy estimate on $W_1$ in \eqref{eq:euler_decoup} with energy $E_{4,m}$ \eqref{energy4m}.

There are three differences between $\cL_{M, i}$ in \eqref{eq:E_lin_MR}, \eqref{eq:euler_decoup} and $\cL_{2D, i}$ in \eqref{eq:lin},\eqref{eq:bous_decoup2}.

First, we use $\bar \uu$ in the transport term instead of $\bar \uu_{2D} = \na^{\perp} \bar \phi_{2D}$. We estimate the difference $\bar \uu - \bar \uu_{2D}$ using \eqref{eq:elli_mainb_err2},\eqref{eq:elli_mainb_err3},\eqref{eq:elli_mainb2},\eqref{eq:E_ub}, and bound $\bar \om (\chi_4 \chi_{\nu} - 1)=\bar \om  (\chi_{\nu} - 1)$ using decay \eqref{solu:decay}
\beq\label{eq:wb_decay2}
|| \bar \om (\chi_4 \chi_{\nu} - 1) ||_{X_1} +  || \bar \om (\chi_4 \chi_{\nu} - 1) ||_{C^1} \les \nu^{-1/32} .
\eeq
Thus, this difference in the linear stability estimate is bounded by 
\[
C ( C_l S + \nu^{-1/32} ) E_{4,m},
\]
where $E_{4,m}$ is the energy \eqref{energy4m}. 
Hereafter, $C$ denotes absolute constants independent of $C_l, S , \nu$.

Second, we use the truncated profile $\bar \th_0, \bar \om_0$ in \eqref{eq:E_lin_MR} rather than the $(\bar \th, \bar \om)$ in \eqref{eq:lin}. We estimate it using the decay of the profiles \eqref{solu:decay},
the asymptotics of the weights, and the elliptic estimates in Lemmas \ref{lem:E_vel}, \ref{lem:E_vel2}
and \eqref{eq:xi_Psi_decay0} and \eqref{eq:E_du_hol}. For example, in $\cL_{M,1} - \cL_{2D, 1}$,  
since $\bar \om - \bar \om_0$ is supported away in $|x| \geq \nu$, 
using  \eqref{eq:xi_Psi_decay0}, we have 
\beq\label{eq:err_decay1}
| u(\om_1) \cdot \pa_x (\bar \om - \bar \om_0) \vp_{g, 1} |
\les |x_1| \cdot |x|^{-5/4} (|x|^{1/16} + |x_1|^{-1/2}) \one_{|x| \geq \nu}
\| \om_1 \|_{X_1} 
\les \nu^{-1/ 32}  E_{4,m}. 
\eeq
This difference in the linear stability analysis is bounded by 
\[
C \nu^{-1/32} E_{4,m} .
\]

Third, 
the main term in the velocity $\uu_M(f) = \na^{\perp} \Psi(f)$ depends on the modified stream function $\Psi(f) = (-\D)^{-1} \Om(f)$  \eqref{eq:elli_main}, \eqref{eq:elli_main2} 
with $\om = f$ rather than 
$\uu_{2D}(f)\teq \na^{\perp}(-\D)^{-1} f$. 
The same argument applies to $\cK_{2i}$ \eqref{eq:euler_cK}.  Using estimates \eqref{wg:equi}, we decompose 
\beq\label{eq:u_dif_3D2}
 \uu_M( f ) = \uu_{2D}(f ) + \na^{\perp}(-\D)^{-1} f_{\mf{err}},
 \quad f_{\mf{err}} \teq \Om(f) - f,
 \quad  \| f_{\mf{err}} \|_{X_1} \les  C_l S \| \tf{|x|}{ \la x \ra } f  \|_{X_1}.
\eeq

For $ f = \hat \om_2$ and nonlocal error $\bar \e, \hat \e$ defined in \eqref{eq:wb_vel_err}, \eqref{eq:lin_evo_main1}, 
due to the extra vanishing power, we obtain  $\f{|x|}{\la x \ra} f \in X_1$. Applying the above estimates with $ f = \om_1$, this leads to a difference in the linear stability analysis bounded by 
\[
C C_l S E_{4,m}.
\]
We also refer to Section \ref{sec:E_lin_lower} for the estimate of the lower order part $\uu_R$, which is small.

The nonlinear terms in \eqref{eq:bous_non}, \eqref{eq:E_non_MR} 
with modification as in \eqref{eq:bous_nonM} and  \eqref{eq:non_dec1} in Section \ref{sec:non} all involve the nonlocal terms determined by $\UU= \uu(\om + \bar \e)$ defined in \eqref{eq:U}.  Recall that in the energy estimate of the 2D Boussinesq, we treat the nonlocal term $\uu_{2D}$ perturbatively. 
Using the estimate of the lower order part $\uu_R$ in Section \ref{sec:E_lin_lower} and the above argument to estimate $\uu_M(f) - \uu_{2D}(f) $ for $f = \hat \om_2, \bar \e, \hat \e$, we have a difference in the nonlinear stability bounded by 
\[
\sup\nolim_{s \leq t} C  C_l S     (E_{4,m}(s) + 1)^2 . 
\]

\vs{0.05in}
\paragraph{\bf Difference in the estimates for $\xi_1$ }

The change of weights from $(\vp_3, \vp_{g,3})$ for 2D Boussinesq 
to $(\vp_{3,m}, \vp_{g, 3,m} )$ defined in \eqref{wg:xi_modi} 
for 3D Euler only affect: (A) $L^{\infty}$ estimates in the region $|x| \geq \nu/2$, 
and (B) $C^{1/2}$ estimates with $\max( |x|, |z| ) \geq \nu/2$. 
In this region, under Assumption \ref{boot_ass:1}, \ref{boot_ass:2}, \ref{boot_ass:E} and 
for $C^{1/2}$ estimates, as estimated in Section \ref{sec:euler_other_xi}, 
the damping coefficients do not decrease; and the $\xi_1$-forcing terms 
are estimated  perturbatively and are bounded by 
\beq
 C ( \nu^{ -\g_{ \eqref{eq:xi_force_lot}  } } + \nu^{-1/32} ).
 \label{eq:xi_diff_1}
\eeq
Thus, the increase of the estimates for these terms compared to those  for 2D Boussinesq
are bounded by \eqref{eq:xi_diff_1}. 
In $L^{\infty}(\vp_2), L^{\infty}(\vp_{g, 2})$ estimates for $\eta_1$, 
as proved in Section \ref{sec:euler_other_xi}, the increase of the estimates of forcing terms involving $\xi_1$ are also bounded by \eqref{eq:xi_diff_1}.

We split weighted $L^\infty$ estimates for $\xi$ into $|x|\le \nu/2$ and $|x|\ge \nu/2$, rather than all of $\mathbb{R}^2_+$ as in 2D Boussinesq case, \emph{solely} to simplify the presentation. 
We close both estimates by bootstrap.

\vs{0.1in}
\paragraph{\bf{Difference between operators for $\widehat W_2$}}

Comparing \eqref{eq:bous_decoup2} and \eqref{eq:euler_decoup}, we have extra terms 
\[
I = (\cL_i - \cL_{2D, i}) \hat W_2 - D_i^2 (\cL_i - \cL_{2D, i})\hat W_2(0) f_{\chi, i} .
\]
Due to correction, $I$ vanishes $O(|x|^3)$ near $0$. 
Recall $\hat W_2 \in C_c^{4, 1}$ from  \eqref{eq:euler_W2}. Under Assumption \ref{boot_ass:E}, we obtain $|| \hat W_2(t) ||_{C^{4,1}} \les 1$ 
by \eqref{eq:euler_W2}. For both local and nonlocal terms in $I$, e.g. 
\[
\bal
& J_{ \bullet } - D_i^2 J_{ \bullet }(0) f_{\chi, i}, \quad J_{ \mf{loc} } =  \bar \uu \cdot \na \hat W_{2, i} - \bar \uu_{2D} \cdot \na \hat W_{2, i} = \d \bar \uu \cdot \na \widehat W_{2, i},  \\
& J_{ \nlocc } = \uu( \hat \om_2) \cdot \na \bar \om_0 - \uu_{2D}( \hat \om_2) \cdot \na \bar \om
= \d \uu(\hat \om_2) \cdot \na \bar \om_0 - \uu_{2D}(\hat \om_2) \cdot \na (\bar \om - \bar \om_0) ,
\eal
\]
for $\bullet = \mf{loc}, \nlocc$, we apply the same estimates as those of the lower order part of residual error in Section \ref{sec:E_err_lower} by replacing $(\bar \om_0, \bar \th_{0,x}, \bar \th_{0, y})$ by $\widehat W_{2}$. To estimate $\d \uu( \hat \om_2)$, we apply \eqref{eq:elli_mainb2}, \eqref{eq:elli_mainb_err}, \eqref{eq:E_ub}. Since $\hat \om_2$ has compact support in $[-L_2, L_2]^2$ 
\eqref{eq:euler_W2}, and $\nu, R_4$ will be chosen to sufficiently large, instead of \eqref{eq:elli_mainb_err}, we have 
$ \hat \om_2 (\chi_4 \chi_{\nu}-1) =0$ and 
\[
|| \Om(\hat \om_2)_{\nu, R}(t) ||_{X_1}  + || \Om(\hat \om_2)_{\nu, R}(t) ||_{C^1} \les C_l S , 
\] 
where we define $\Om(\hat \om_2)$ via \eqref{eq:elli_main}, \eqref{eq:elli_main2}. We estimate the cutoff error  $\bar \om_0 - \bar \om$ using \eqref{solu:decay}.  

In summary,  we have smallness from the difference between two nonlocal operators $\uu - \uu_{2D}$ or the decay of the profile in the estimate of this difference, and can bound it at time $t$ by 
\[
C( C_l S + \nu^{-1/32} ).
\]

\subsection{Nonlinear stability and finite time blowup}\label{sec:euler_blowup_pf}

Let the initial perturbation $(\om,\eta,\xi)$ belong to the energy class
$\Einm$ \eqref{energy:modi}, satisfy $E_0 \teq \Einm(\om,\eta,\xi) < E_* $, 
\footnote{
Although the weights $\vp_{3,m}, \vp_{g,3,m}$ \eqref{wg:xi_modi} jump across $|x| = \nu/2$, 
due to regularity of $\xi|_{t=0}$, 
condition $E_0 = \Einm(\om,\eta,\xi) < E_*$ implies $ |\xi \vp_{ 3}|,   
\mu_{g, 3} |\xi \vp_{g, 3}| \leq E_0 < E_* $ on $|x| = \nu/2$.
}
and obey the symmetry condition \eqref{eq:euler_sym}. 
\footnote{
Recall energy $E_{4,m}$ \eqref{energy4m} 
for $W_1, \hat W_2$. 
Initially, we have $\hat W_2=0$, and  $E_{4,m}$ and $\Einm$ are the same.
}
We can choose small initial perturbation $\om, \th= \chi_\dap g$ with
\[
 \om_{\tot} = \om + \bar \om_0 \in C_c^{\infty}, \quad  \th_{\tot} =\th + \bar \th_0 = \chi_\dap ( \bar g + \bar \th + M) ,
 \quad \bar g + \bar \th  \in C^{\infty}, \quad  \bar g + \bar \th  + M \geq M / 10,
\]
which gives smooth initial data $\om_0^{\th}, u_0^{\th}$ defined via \eqref{eq:omth}.

Under the bootstrap Assumptions \ref{boot_ass:1}, \ref{boot_ass:2}, \ref{boot_ass:E}, 
we perform nonlinear energy estimates as in the 2D Boussinesq case in Section \ref{sec:EE}. Combining \eqref{eq:E_err_M}, the estimates in Section \ref{sec:E_lin_lower}, \ref{sec:E_err_lower}, 
\ref{sec:E_xi_lower},\ref{sec:euler_other_xi}, and the comparison in Section \ref{sec:E_lin_M}, we  bound the additional terms  from  the differences between the 3D Euler and 2D Boussinesq estimates, including the weighted $L^{\inf}$ and $C^{1/2}$ estimates, and the nonlinear modes \eqref{eq:W2_non_boot},\eqref{eq:W2_non_boot3} (the coefficients of $\cNF_i$ in \eqref{eq:appr_near0}, \eqref{eq:euler_decoup}) by 
\[
\bal
  |J | &\leq C (C_l S +   || \Om_{\nu, R} ||_{X_1} + \nu^{-1/32} 
+ \nu^{ -\g_{ \eqref{eq:xi_force_lot}  } } ) ( 1 + \sup\nolim_{s \leq t} E_{4,m}(s) + E_{4,m}(t)^2) \\
& \leq C_{1, *} (C_l S + \nu^{-\g_{ \eqref{eq:xi_force_lot}  } }    ) ( 1 + E_* + E_*^2) ,
\eal 
\]
for $C_{1,*}$ independent of $C_l, S, \nu$, where  we further bound $|| \Om_{\nu, R}||_{X_1}$ using \eqref{eq:elli_mainb_err}, \eqref{eq:wb_decay2}. Recall $c_{\om}, \bar c_{\om}$  from \eqref{eq:E_normal}. From the energy estimate and the definition of $E_{4,m}$ \eqref{energy4m}, we also have 
\beq
|c_{\om} + \bar c_{\om} - \bar c_{\om, 2D}| \leq 100 E_{4,m}(t) + C_{2, *} ( \nu^{-
\g_{ \eqref{eq:xi_force_lot}  } } + C_l S  ).
\label{eq:boot_cw}
\eeq

The energy estimates for 2D Boussinesq satisfy the nonlinear stability conditions \eqref{eq:PDE_nondiag} with  $\e_0 > 0$, 
and the second inequalities in \eqref{eq:W2_non_boot2} are strict with a gap $\e_1>0$. 
We  choose $\nu > \nu_* > \nu_{*, 1}$  with $\nu_*$ large enough  and $\nu_{*,1}$ chosen in Proposition \ref{prop:xi_far_linfty} and choose a small $\d$ such that 
\beq
( C_{1,* } + C_{2, *} + C_{ \eqref{eq:xi_far_linfty} } ) ( \nu_*^{-1/32}+ \nu_*^{- \g_{ \eqref{eq:xi_force_lot}  } } + \d + \d^{\b}) 
 < \tf{1}{100} \min(  \e_0 , \, \e_1, 
\, 10^{-4}, E_* ).
\label{eq:E_non}
\eeq

\begingroup
\renewcommand{\theassumption}{1$'$}
\begin{assumption}\label{boot_ass:1p}
We impose $S(0) \geq \nu$ and a stronger bootstrap assumption than 
Assumption \ref{boot_ass:1}
\beq\label{eq:boot6}
 C_l(t) (1 + S(t) ) < \min( \d , \nu_2, 4^{-6}).
\eeq
\end{assumption}
\endgroup

Under the Assumptions \ref{boot_ass:1p}, \ref{boot_ass:2},  \eqref{eq:boot6}, \eqref{eq:W2_non_boot3}, and the energy assumption \ref{boot_ass:E}: 
\beq\label{eq:E_boot}
E_{4, m}(t) < E_*,
\eeq
using the nonlinear stability estimate for the 2D Boussinesq equations, 
Lemma \ref{lem:PDE_nonstab}, \eqref{eq:W2_non_boot2}, and the smallness condition \eqref{eq:E_non}, we improve  \eqref{eq:W2_non_boot3}  and the estimates of the energies in $E_{4,m}$ except for $\nlinf{ \xi_1 \vp_{3,m}}, \nlinf{\xi_1 \vp_{g,3,m}}$ in \eqref{energy4m} by $E_{ \mf{sm} } \teq \max( E_0 , (1 - \d_{\eqref{eq:PDE_nonboot_impr} }   ) E_*)$.

Due to the outgoing Assumption \ref{boot_ass:2} of the flow,  solution $\xi$ in $|x| \leq \nu/2$ 
does not depend on $\xi(t, x)$ for $|x| > \nu/2$. Since $\xi(0, y), \pa_y \bar \th_0(0,y) =0$ for $|y|\leq \nu$ by \eqref{eq:3D_thy_van}, as in 2D Boussinesq, applying  $L^{\infty}(\vp_3), L^{\infty}(\vp_{3,m})$ stability estimates 
 on $\xi_1$ and using the smallness of the difference, we obtain
$\ \| \xi_1(s) \vp_{3} \|_{L^{\infty}(B_{\nu/2})} , \mu_{g, 3} \| \xi_1(s) \vp_{g, 3} \|_{L^{\infty}(B_{\nu/2})} 
   \leq E_{ \mf{sm} }$ for $s \leq t$. 
   Since $\vp_{3,m} \leq \vp_{3}, \vp_{g,3,m} \leq \vp_{g, 3}$,
   this estimate controls the boundary value in Proposition \ref{prop:xi_far_linfty}.
   Thus, combining this estimates, Proposition \ref{prop:xi_far_linfty},
    $E_0 < E_*$, and the smallness in \eqref{eq:E_non}, 
   we prove 
   \[
 \| \xi_1(s) \vp_{3,m} \|_{L^{\infty}(B_{\nu/2}^c)} , \ \mu_{g, 3} \| \xi_1(s) \vp_{g, 3,m} \|_{L^{\infty}(B_{\nu/2}^c)}  \leq E_{\mf{sm}}.
   \]
Combining the above estimates, we improve the  bootstrap Assumption \ref{boot_ass:E}, \eqref{eq:E_boot}
on energy $E_{4, m}$.

Moreover, using $ | \om| \les \la x \ra^{-1/32} E_{4,m}$ from the $L^{\inf}(\vp_{g,1})$ estimate and \eqref{eq:boot_cw}, we have 
\beq\label{eq:E_u_sublinear}
| \uu + \bar \uu| \les  |x|^{ 1-1/64}, \quad \bar c_l + c_l = \bar c_l  > 2, \quad \bar c_{\om} + c_{\om} <-\tf{1}{2}.
\eeq
Thus, the whole velocity \emph{grows sublinearly} and the blowup is focusing ($\bar c_l > 2$). Following the argument in \cite[Section 9.3.5]{chen2019finite2}, under bootstrap assumptions, we control the support 
\[
C_l(t) (1 + S(t)) \leq C( S(0)) C_l(0),
\] 
uniformly in $t$, for some constant $C$ depending on $S(0)$. 

Upon further requiring $ \nu$ large, 
using \eqref{eq:E_u_sublinear}, we obtain
\[
  ( \bar c_l \xx +  \uu + \bar \uu ) \cdot \xx \geq \bar c_l |\xx|^2 - C |\xx|^{2-1/64} \geq \tf12 \bar c_l |\xx|^2 , \quad |\xx| = \tf12 \nu ,
\]
which improves Assumption \ref{boot_ass:2}. We then fix $\nu$ for the profile \eqref{eq:profile2} 
and for the weight $\vp_{3, m}, \vp_{g, 3, m}$.

For any $S(0) \in [\nu, \infty)$, 
by choosing $C_l(0)$ sufficiently small and using the above estimate, we improve the stronger Assumption \ref{boot_ass:1p} and hence also  Assumption \ref{boot_ass:1}. Thus, we establish nonlinear stability estimates \eqref{eq:E_boot},\eqref{eq:W2_non_boot3},\eqref{eq:boot6},\eqref{eq:boot5} for any $t \geq 0$.

Passing from nonlinear stability to finite time blowup with smooth data $\om^{\th}, u^{\th}$ compactly supported near $(r, z) = (1, 0)$ follows the argument in \cite{chen2019finite2}. 
This proves Theorems \ref{thm:euler} and \ref{thm1b}.

\section{Construction of an approximate steady state}\label{sec:ASS}

Following our previous works with Huang on the De Gregorio model \cite{chen2019finite} and the Hou-Luo model \cite{chen2021HL}, we construct the approximate steady state to the dynamic rescaling equations \eqref{eq:bousdy1} with the normalization conditions \eqref{eq:normal} by solving \eqref{eq:bousdy1} numerically for a long enough time. The residual error is estimated {\it a-posteriori} and incorporated in the energy estimate as a small error term. It is extremely challenging to obtain an approximate steady state with a {\it sufficiently} small residual error in the weighted energy space \eqref{energy4}, e.g. of order $10^{-7}$, since the weight is singular of order $|x|^{-\b}, \b \geq 2.9$ near $0$ and the solution is supported on the whole $\R^2_+$  with a slowly decaying tail in the far-field, e.g., $\om(t, x) \sim |x|^{-c}, c \approx \f{1}{3}$ for large $x$. See \eqref{eq:ASS_asym}.

To overcome these difficulties and control the round-off error, the construction has three main components. First, Section~\ref{sec:ASS_far_field} derives the far-field asymptotics and motivates the decomposition into a semi-analytic tail and a faster-decaying numerical remainder to reduce 
far-field round-off error.  Second, Sections~\ref{sec:fit_profile}--\ref{sec:ASS_refine} describe how the angular profiles and B-spline remainder are constructed. Third, Section~\ref{sec:ASS_multi} explains the multi-level representation and near-origin corrections used to reduce round-off error in the high-order derivative estimates.
Without these methods, the computation would suffer from relatively large round-off errors. 

\subsection{Far-field asymptotics}\label{sec:ASS_far_field}

The far-field asymptotics of the profiles are determined by the leading-order balance between the scaling fields \(c_l \xx \cdot\nabla\) and the terms \(c_\omega\omega, c_\theta\theta\) in \eqref{eq:bousdy1}, and \emph{ is not} an independent numerical ansatz.
Consider the polar coordinates in $\R_2^+$: $ r = (x^2 + y^2)^{1/2},  \b = \arctan(y / x)$.
As in \cite{chen2021HL}, we consider the approximate steady state with the  far-field asymptotics
\beq\label{eq:ASS_asym}
 \om(r, \b) \sim  r^{ \al} g_1(\b) , \quad  \th(r, \b) \sim x \cdot r^{ 2 \al} g_2(\b),  \quad \al = \f{c_{\om}}{c_l}< 0 , \quad \al \approx -\f{1}{3},
\eeq
for some angular profiles $g_1(\b), g_2(\b)$, under the mild assumption that $\om$ decays for large $|x|$, $c_l >0$, and $c_{\om } < 0$. 
In fact, if $\om$ decays for large $|x|$, the velocity $\uu = \na^{\perp}(-\D)^{-1} \om$ has a sublinear growth: $\f{ | \uu(x) |  }{r} \to 0 $ as $r \to \infty$. 
Since $ x \cdot \na = r \pa_r $, passing to the polar coordinate and dropping the lower order terms and $\pa_t \om, \pa_t \th$ in \eqref{eq:bousdy1}, we yield 
\[
c_l r \pa_r \om(r, \b) = c_{\om} \om + \th_x + l.o.t., \quad c_l r \pa_r \th(r, \b) = (2 c_{\om} + c_l) \th + l.o.t.. 
\]

Under the ansatz $\om(r, \b) = r^{k} g_1(\b), \th(r, \b) = x  r^{\ell} g_2(\b)$,  
  matching powers of $r$  for $\th$-equation and then for $\om$-equation
   yields $\ell=2\al, k=\al$,
which gives \eqref{eq:ASS_asym}. Thus, we represent $\bar \om, \bar \th$ as
\beq\label{eq:ASS_decomp1}
\bar \om = \bar \om_1 + \bar \om_2, \quad \bar \th = \bar \th_1 + \bar \th_2 , \quad \bar \om_1 = \chi(r) r^{ \al} \bar g_1(\b), \quad 
\bar \th_1 = \chi(r) x \cdot  r^{  2 \al}  \bar g_2(\b),
\eeq
where $\chi(r)$ is the radial cut-off function defined in \eqref{eq:cut_radial} 
and \emph{differs} from $\chi_{rati}$ in \eqref{eq:cut_radial}. The crucial semi-analytic parts $(\bar \om_1, \bar \th_1)$ capture the far-field asymptotic behavior \eqref{eq:ASS_asym}. 
The second remainder part has a much faster decaying rate, and we construct it using numerical computation with  sixth order B-splines  on the adaptive mesh for a large bounded domain. This decomposition is essential both for reducing round-off error and obtaining small residual error.

\subsection{Angular profiles and the representation}\label{sec:ASS_rep}

By symmetry in $x$, we compute \eqref{eq:bousdy1} on $[0, L]^2$ with $L \approx 10^{13}$. 
In $\R^2_{++}$, the stream-function components $\bar \phi_2 , \bar \phi_3, \bar \phi_{\mf{cor}}$ in \eqref{def:phi_N} are supported on a larger domain $D_{\mf{lg}} =[0, L_2]^2$ with $L_2 \approx 10^{15}$. We partition $[0, L_2]$ using a non-uniform mesh 
\[
0 = y_0 < y_1 < \cdots < y_{N-1} = L_2.
\]
Since $\th(t, 0, y) =0 $ and $\th(t, x, y)$ is even in $x$, $\th(t, x, y)$ vanishes quadratically on $x= 0 $. Thus, instead of using $\th$ in our computation, we consider $\zeta(t, x, y) = \f{1}{x} \th(t, x, y) $. Then $\zeta$ is odd in $x$.

Using $c_{\th} =c_l + 2 c_{\om}$ and normalization conditions  \eqref{eq:normal}, we rewrite 
equations \eqref{eq:bousdy1} equivalently: 
\beq\label{eq:SS_force}
\bal
\om_t & = F_{\om}(\om, \zeta) \teq  - (c_l \xx + \uu) \cdot \na \om + x \pa_x \zeta + \zeta + c_{\om} \om , \\
 \zeta_t & = F_{\zeta}(\om, \zeta) \teq -   (c_l \xx + \uu)\cdot \na  \zeta +  ( 2 c_{\om} - \tf{u}{x} )  \zeta , \quad  \uu = (u, v ) = \na^{\perp} \phi, \quad -\D \phi = \om .
 \eal
\eeq

Without the semi-analytic part, we represent $\om , \zeta$ using $6$th order B-splines in $x$ and $y$, e.g.
\beq\label{eq:spline_solu1}
\om(t, x, y) = \tts{ \sum_{i, j} }  \, a_{i,j} B_{1,i}(x) B_{2,j}(y) \, , 
\eeq
where $B_{1,i}(\cdot), B_{2,j}(\cdot)$ are the B-spline basis and are piecewise polynomials.
We choose $B_{1,i}(x)$ odd in $x$ so that $\om, \zeta$ are odd in $x$. The subscripts $1$ and $2$ in $B_{1,i}$ and $B_{2,j}$ indicate the $x$- and $y$-directions, respectively. We represent $\phi$ using B-splines with an extra weight $\rho_p(y)$ vanishing on  $y = 0$ to enforce the boundary condition $\phi(x, 0 ) = 0$. Given the grid point values of $\om$ using suitable extrapolation of $\om$ in the far-field, we obtain the coefficients $a_{i,j}$ of $\om$ by interpolation  and solving the linear equations \eqref{eq:spline_solu1} for $a_{i,j}$. 
We solve the Poisson equations 
\beq\label{eq:ASS_pois}
-\D \phi = \om
\eeq
using a B-spline based finite element method and obtain the B-spline coefficients 
of $\phi$. After obtaining the coefficients $a_{i, j}$, we compute the derivatives of $\om$ using the basis functions 
\[
\pa_x^k \pa_y^\ell \om(t, x, y) = \tts{ \sum_{i, j} } \, a_{i,j} \pa_x^k B_{1,i}(x) \pa_y^\ell B_{2,j}(y) .
\]
Similar considerations apply to $\zeta, \th, \phi$. 
See Part II \cite[Appendix {\appsolurep}]{ChenHou2023b} for more details on mesh $y_i$,  B-spline basis, and representations. 
Similar  B-spline  representations  have been used in \cite{luo2013potentially-2}.

To discretize the temporal variable, we use a second order Runge-Kutta method. 
Using the above methods, given $(\om, \zeta, c_l, c_{\om})$, we derive the forcing $F_{\om}, F_{\zeta}$ in \eqref{eq:SS_force} and  update the PDE.

To construct the decomposition in \eqref{eq:ASS_decomp1}, we proceed in two
stages. First, we determine the exponent $\bar\al_1$ and the angular profiles
$g_i(\b)$, which determine the shape of the semi-analytic far-field parts $\bar \om_{10}, \bar \th_{10}, \bar \zeta_{10} = \f{1}{x} \bar \th_{10}$. Then, during the computation, we use $\bar \om_1 = c_1 \bar \om_{10}, \bar \th_1 = c_2 \bar \th_{10}$ and the ansatz \eqref{eq:ASS_decomp1}, adjust the prefactors  $c_1, c_2$, and refine the construction of $\bar\om_2,\bar\th_2$ in \eqref{eq:ASS_decomp1}.

\subsubsection{Fitting the angular profile and the exponent}\label{sec:fit_profile}
This subsection determines the fixed angular profiles in the far-field part in \eqref{eq:ASS_decomp1}.
First, we solve \eqref{eq:bousdy1} numerically using the above method without the semi-analytic part, i.e. $\bar \om_1 = 0, \bar \th_1 = 0$, to obtain an approximate steady state $ (\bar \om, \bar \zeta  ), \bar \th = x \bar \zeta$. 
Using ansatz \eqref{eq:ASS_decomp1} and fitting the angular part of the far-field of $r^{-\bar \al_1} \bar \om, x^{-1} r^{- 2 \bar \al_1} \bar \th  = r^{- 2 \bar \al_1}   \cdot \bar \zeta$  with exponent $\bar \al_1 = \f{ \bar c_{\om, \mf{pre}}}{\bar c_{l, \mf{pre}} }$ \eqref{eq:ASS_asym}
\footnote{
We use preliminary values of $\bar c_{\om, \mf{pre}} , \bar c_{l,\mf{pre}}$ to determine $\bar \al_1 =  \bar c_{\om, \mf{pre}} / \bar c_{l,\mf{pre}}$.
These values \emph{differ from} $\bar c_{\om}, \bar c_l$ in the final profile discussed in Section \ref{sec:property}, but they are very close. This exponent $\bar \al_1$ is the \emph{same} as that 
in \eqref{solu:decay}.
}
, we find the following approximate profiles 
\footnote{ 
We \emph{intentionally} use $\td \beta$ for $g_{10}$ and $\cos \beta$ for $g_{20}$. These  two different  forms are \textit{not} typos. 
}
\[
g_{10}(\b) = \f{ a_{11} \td \b ( 1 + a_{15}  \td \b^2   )  }{ 
( \td \b^2 + a_{12})^{2/3} + a_{13} \td \b^2 + a_{14} }, 
   \ \td \b = \f{\pi}{2} - \b, 
   \quad g_{20}(\b) =  \f{  a_{21} \cos \b ( 1 + a_{25} \sin \b)  }{  
( \cos^2 \b + a_{22})^{2/3} + a_{23} + a_{24} \cos^2 \b } \, ,
\]
for some parameters $a_{ij}$. We use the factor $\f{\pi}{2} - \b$ since $\om$ is odd in $x$ and $g_{10}(\b)$ is odd with respect to $\b = \pi/2$. Similarly,  since $\theta(0, y) = 0$ and $\th_x, \zeta = \f{1}{x} \theta$ are odd in $x$,
 $g_{20}(\b)$ is odd with respect to $\b = \pi/2$. 
 Hence, we add the factor $\cos \b$ in  $g_{20}(\b)$.   After obtaining the above analytic formulas, we further approximate the above profiles $g_{i0}(\b)$ by $g_i(\b)$, which are 
linear combinations of $8$-th order B-splines  $B_{ i}^{(8)}$  ( see Appendix C.1 of Part II \cite{ChenHou2023b})
\[
g_1( \f{\pi}{2} - s) = \sum\nolim_{ 1 \leq i \leq n} b_{1i} B_{ i}^{(8)}( s ), 
\quad  g_2( \f{\pi}{2}- s ) =  \sum\nolim_{ 1 \leq i \leq n} b_{2i} B^{(8)}_{ i}( s ), 
\quad s \in [-\f{\pi}{2} , \, \f{\pi}{2}],
 \]
for some coefficients $b_{ji}$. 
Here $s$ is a dummy variable, and $B_i^{(8)}(s)$ are chosen to be \emph{odd} Bspline basis functions 
in $s$. Thus, the resulting
profiles $g_i( \b)$ are odd with respect to the physical angle $ \b=\pi/2$. The explicit formulas for \(g_{10}\) and \(g_{20}\) are used only to obtain an initial fit of the far-field shape.  We use the B-spline representation of \(g_1\) and \(g_2\)  in the final construction.

We further use the B-spline representation $g_i(\b)$ rather than $g_{i0}(\b)$ for the angular profiles since the high-order derivatives 
of the B-spline can be bounded systematically, while the derivatives of 
$g_{i0}(\b)$ are  complicated to estimate. These derivative estimates are used to estimate the profile.

Once we obtain the \emph{approximate angular profile} $g_i(\b)$, we construct the semi-analytic part 
\beq\label{eq:ASS_semi1}
 \quad \bar \om_{10} = \chi(r) r^{ \bar\al_1} g_1(\b), \quad 
 \bar \zeta_{10} = \chi(r) r^{2 \bar \al_1} g_2(\b),
 \quad  \bar \th_{10} = x \cdot \bar \zeta_{10}.
\eeq

Given the asymptotic behavior of $\bar \om_{10}$ in \eqref{eq:ASS_semi1}, 
for $  (-\D)^{-1} \bar \om_{10}$, we consider the far-field asymptotic ansatz   $r^{2 + \bar \al_1} f(\b)$ with some profile $f(\b)$. We construct $f(\b)$ by solving
 \[
  - \D (  r^{2 + \bar \al_1} f(\b)) = r^{\bar \al_1} g_1(\b)
 \]
 with boundary condition $f(0) = f(\pi/2) = 0$ due to the Dirichlet boundary condition $\phi(x,0)=0$ and the odd symmetry of $\om$ in $x$. In the polar coordinate, the above equation is equivalent to 
\beq\label{eq:ASS_pois_1D}
 (-\pa_{\b}^2 - (2 + \bar \al_1)^2 ) f(\b) = g_1(\b) , \quad f(0) = f(\pi/2) = 0.
\eeq
We represent $f(\b)$ using a weighted $8$th order B-spline and solve the above elliptic equations using the finite element method. Then, we construct the semi-analytic part for $\phi$ as follows 
\beq\label{eq:ASS_semi2}
\bar \phi_{10} = \chi(r) r^{2 + \bar \al_1} f(\b). 
\eeq

\subsubsection{Refinement}\label{sec:ASS_refine}

We now refine the approximate profile while keeping the exponent \(\bar\alpha_1\) and the angular profiles \(g_1,g_2\) fixed. The purpose of this step is to determine the scalar amplitudes of the semi-analytic part and to represent the remaining faster-decaying part by B-splines.

Given the grid values of $\om(t,x,y)$ at $(y_i,y_j) \in [0,L]^2$, we first determine $c_1(t)$ by fitting $c_1(t)\bar\om_{10}(x,0)$ to $\om(t,x,0)$ over far-field grid points away from the computational boundary, for example, $x=y_i$ with $10^{-4}L<|y_i|<L/4$. We then define $\om_2(t,x,y)=\om(t,x,y)-c_1(t)\bar\om_{10}(x,y)$ and represent $\om_2$ by the B-spline interpolation \eqref{eq:spline_solu1}. Thus, $\om$ is represented as $c_1(t)\bar\om_{10}+\om_2$, where $\om_2$ is piecewise polynomial on the computational domain. The same procedure is applied to $\zeta$ 
with far-field components $c_2(t) \bar \zeta_{10}$. Updating $c_1(t)$ primarily adjusts the far-field behavior of the representation, while the corresponding update of $\om_2$ preserves the prescribed grid values of $\om$ on $[0,L]^2$ up to the B-spline interpolation error. To update the stream functions $\phi$, we use $c_1(t)  \bar \phi_{10}$ to capture the far-field of $\phi$ and then construct the near-field part by solving 
\beq\label{eq:ASS_pois2}
 -\D ( \phi_2 + c_1(t) \bar \phi_{10}   ) = \om_2 + c_1(t) \bar \om_{10}, \quad  \mathrm{or \quad } -\D \phi_2
 = \om_2 + c_1(t) (\bar \om_{10} + \D \bar \phi_{10}).
\eeq
Then the stream function is represented as $\phi_2 + c_1(t) \bar \phi_{10}$.

The advantage of \eqref{eq:ASS_pois_1D}--\eqref{eq:ASS_pois2} over directly
solving \eqref{eq:ASS_pois} is the reduced far-field round-off error.
The large computational domain and the adaptive mesh leads to a
poorly conditioned stiffness matrix. Moreover, the source $\om$ decays
slowly, like $r^{\bar \al_1}\approx r^{-1/3}$. Hence, directly solving
\eqref{eq:ASS_pois} produces significant round-off error. In
\eqref{eq:ASS_pois2}, the semi-analytic part $c_1(t)\bar\om_{10}$ captures the
far-field asymptotics of $\om$, making $\om_2$ much smaller than $\om$ in the far field.
Moreover, by \eqref{eq:ASS_semi1}--\eqref{eq:ASS_semi2}, the far-field part of
$\bar\om_{10}+\D\bar\phi_{10}$ is only of size $\epsilon r^{-1/3}$ with a
small constant $\epsilon$. Thus the source in \eqref{eq:ASS_pois2} is much
smaller  than  $\om$ in the far field, significantly reducing the round-off error. 
\footnote{
A similar technique has been used in the Hou–Luo model \cite{chen2021HL} to overcome the  significant round-off errors.
}

The terms $c_1(t) \bar \om_{10}, c_2(t) \bar \zeta_{10}$ capture the far-field components of $(\om, \zeta)$, respectively, out of the computational domain $[0, L]^2$. 
These far-field components  affect  the equations \eqref{eq:SS_force} in $[0, L]^2$ through the nonlocal velocity $\uu( \om )$, which is determined by the stream function and \eqref{eq:ASS_pois2}. 
\footnote{
For $(\om, \theta)$ close to $(\bar \om, \bar \th)$, 
 we expect that the outgoing condition \eqref{eq:outgoing} 
 holds for the associated $(c_l, \uu)$. 
 This condition implies that information propagates outward near the boundary of the computational domain $[0, L]^2$.
 }  
 After obtaining the stream function, we derive the forcing $F_{\om}, F_{\zeta}$ and update the PDE \eqref{eq:SS_force} using the second order Runge-Kutta method. 
We stop the computation at time $t_*$ if the residual error on the grid points is about the round-off error. Then we finalize the semi-analytic part in \eqref{eq:ASS_decomp1} as
\beq\label{eq:ASS_semi3}
\bal
& \bar \om_1 = \bar c_1 \bar \om_{10} 
= \chi(r) r^{\bar \al_1 }( \bar c_1 g_1(\b)) , \qquad \bar \th_1 = \bar c_2 \bar \th_{10} 
= \chi(r) x \cdot r^{ 2 \bar \al_1 } ( \bar c_2 g_2(\b)  ), \\
& \bar \phi_1 = \bar c_1 \bar \phi_{10} = \chi(r) r^{2 + \bar \al_1} ( \bar c_1 f(\b) )
\teq \chi(r) r^{2 + \bar \al_1} \bar f(\b),  
\eal
\eeq
where $\bar c_1 \bar \om_{10},\bar c_2 \bar \zeta_{10}$ best approximate $ \om(t_*, x, y), \zeta(t_*, x, y)$ in the far-field, respectively. We construct $\bar \om_2, \bar \th_2  = x \bar \zeta_2$ in \eqref{eq:ASS_decomp1} by interpolating the grid point values of $\om - \bar \om_1, \th - \bar \th_1$ and applying a low-pass filter to the solution to reduce the round-off error.

In Part II \cite[Appendix C]{ChenHou2023b}, we estimate the derivatives of the approximate steady state rigorously, which will be used to verify the residual error.

\subsubsection{A multi-level representation and round-off error}\label{sec:ASS_multi}

The previous subsections address the far-field tail. We now address another numerical difficulty: the rigorous weighted residual estimates require high-order derivative bounds near the origin, which may not be small due to the round-off error.  Below, we explain this difficulty  and 
design the multi-level representation to reduce this error without needing to use high precision computation, which could be significantly slower.

Near the origin, the mesh size $y_{i+1} - y_i$ (see $y_i$ in Section \ref{sec:ASS_rep}) is small, e.g. $h \leq \f{1}{256}$ in our computation. Hence, when computing high order derivatives, e.g. $\na^4 \bar \om$, the round-off error may not be relatively small. To construct the approximate steady state and the approximate space-time solution in  Part II \cite[Section 3]{ChenHou2023b}, we only use lower order derivatives $ \na \om, \na \eta, \na \xi, \na^2 \phi , \na \th $, and the round-off error is negligible. 
However, verifying the smallness of weighted norm of the error, e.g. $\bar \cF_i$ in \eqref{eq:bous_err}, requires piecewise $C^3$ bounds and estimates \eqref{eq:key_weight_est} since the weight $\vp$ is singular like $|x|^{-3}$ near the origin. Thus, we evaluate $ \na^3 \bar \cF_i, \na^4  \om, \na^4  \eta, \na^4 \xi, \na^5 \th$ on some grid points based on the estimates in  supplementary material II of Part II \cite[Section {\secCIIItonorm}]{ChenHou2023b}.

To obtain rigorous bounds, we use interval arithmetic. For each operation, e.g. 
$a \in [a_l, a_u], b \in [b_l, b_u]$, the interval bound for $ a b $ is obtained by considering the worst case. If we use interval arithmetic with lower precision, e.g. the double precision with machine error about $10^{-16}$, 
the interval bound for $\na^3\bar\cF_i$ may greatly overestimate the actual round-off error. One may overcome this problem  using higher  precision, e.g. interval arithmetic with quadruple precision.

To save the computational cost, we instead refine the B-spline representation of the solution $f$ so that $\na^k f$ has much smaller round-off error. The round-off error of $\na^k f$ is about $C \e h^{-k}$, where $C$ is the size of the B-spline coefficient for $f$, $\e$ is the machine precision, and $h$ is the mesh near $0$.  
Thus the round-off error can be reduced either by increasing $h$ or by reducing the size $C$ of the B-spline coefficients. 
We use a multi-level B-spline representation 
\[
f  = f_1 + f_2 + ...+ f_n, \quad  f_k  \  \mbox{uses supporting points }  M^{(k)} = \{ (y^{(k)}_i, y^{(k)}_j ) \}   , \quad 
M^{(k)} \subset M^{(k+1)}, 
\]
to achieve both effects, where $M^{(k+1)}$ is a refined mesh of $M^{(k)}, k \leq n-1$. 
\footnote{
Fix $m_k\geq 2$. Each interval $[y_i^{(k)},y_{i+1}^{(k)}]$ of
$M^{(k)}$ is divided into $m_k$ equal subintervals in $M^{(k+1)}$. 
We have
$
  y_i^{(k)}=y_j^{(k+1)}<\cdots<y_{j+m_k}^{(k+1)}=y_{i+1}^{(k)},
$
for some $j$ with $
  y_{\ell+1}^{(k+1)}-y_\ell^{(k+1)}
  ={1\over m_k} (y_{i+1}^{(k)}-y_i^{(k)} ),  j\leq \ell\leq j+ m_k-1 $.
}
Since the profile $\bar f$ is quite smooth, we use the first level representation 
$f_1$  with  mesh $M^{(1)}$ to interpolate $\bar f$ and the round-off error for $\na^k f_1$ is very small since  its mesh size is much larger. At the next level, we use 
mesh $M^{(2)}$ with smaller mesh size,  and use $f_2$ to interpolate $\bar f - \bar f_1 $. Since $\bar f - \bar f_1$ is much smaller than $\bar f$, the coefficients for the B-spline $f_2$ are small and the round-off error is small. The same procedure applies to other levels. We use the nested mesh structure $M^{(k)} \subset M^{(k+1)}$ so that $f=f_1+\cdots+f_n$ is still a piecewise polynomial on the finest mesh $M^{(n)}$ for $f_n$.
Then we estimate the piecewise derivatives of $f$ using the method in Part II \cite[Appendix {\secpiecepol}]{ChenHou2023b}.

To obtain a \emph{very small} residual error in solving the stream function $\phi_2$ \eqref{eq:ASS_pois2} near the origin, we use another round-off error control method. The leading \(xy\)-mode of \(\bar\phi_2\) can contribute relatively large B-spline coefficients on the finest mesh. Hence, we remove this mode analytically before representing the remaining part by the multi-level B-spline expansion. Near $x=0$, since $\bar \phi_2 = \pa_{xy} \bar \phi_2(0) xy + O(|x|^3)$, we approximate it using an analytic profile 
\beq\label{eq:psi_near0}
  \bar \phi_3 = a \cdot  \chi_{\phi, 2D} , \quad \chi_{\phi, 2D} = - x y \chi_{\phi}(x) \chi_{\phi}(y) ,
\eeq
where $a$ is chosen to approximate $-\pa_{xy} \bar \phi^N(0)$, with $\bar \phi^N$ denoting the whole stream function, and $\chi_{\phi}$ is some cutoff function satisfying $\chi_{\phi}(x) = 1 + O(|x|^2)$ near $0$ and is constructed in \eqref{eq:cutoff_psi_near0}. This is a second round-off error control mechanism, specifically for the leading $xy$-mode of the stream function near the origin. We add the negative sign to normalize 
\[
u_x( (-\D)\chi_{\phi, 2D})(0) = - \pa_{xy}(-\D)^{-1}(-\D) \chi_{\phi, 2D}(0) = -\pa_{xy} \chi_{\phi,2D}(0)=1 .
\]

In solving the approximate steady state, at the $n-th$ time step, we choose $a = -\pa_{xy} \phi_{n-1}(0)$, where $\phi_{n-1}$ is the stream function obtained at the $(n-1)$-th time step.
Then we use the multi-level B-spline representation $\phi_2 =\phi_{2, 1} + ... + \phi_{2, n}$ by solving a modification of \eqref{eq:ASS_pois2}
\[
 - \D \phi_2 =  \om_2 + c_1(t) (\bar \om_{10} + \D \bar \phi_{10}) + \D \phi_3 \teq S.
\]
The approximation term $\phi_3$ allows us to obtain smaller spline coefficients for $\phi_2$ and 
significantly reduce the round-off error near $0$, where a very small error is required.
To obtain the top-level B-spline $\phi_{2, 1}$ on the coarse mesh $M^{(1)}$, we first restrict $S$ on the mesh $M^{(1)}$ and interpolate $S( y_i^{(1)}, y_j^{(1)})$ using the single level B-spline $S^{(1)}$ with supporting points on $y^{(1)}$. Then we use the B-spline based finite element method to solve $ - \D \phi_{2, 1} = S^{(1)}$. We evaluate $ \D \phi_{2, 1}$ on the fine mesh and further solve $\phi_{2, 2}, \phi_{2, 3}...$ recursively from the remaining source part $S + \D \phi_{2, 1}$.

After we obtain the above stream function, we further add a rank-one correction $\bar \phi_{\mf{cor}}$ near $0$ 
and construct the final approximate stream function $\bar \phi^N$   
\beq\label{def:phi_N}
 \bar \phi^N = \bar \phi_1 + \bar \phi_2 + \bar \phi_3 +  \bar \phi_{ \mf{cor}},
 \quad \bar \phi_{\mf{cor}} = - c \cdot \f{x y^2}{2} \kp_*(x) \kp_*(y) , \quad c  = \pa_x( \bar \om + \D (\bar \phi_1 + \bar \phi_2 + \bar \phi_3) )(0),
\quad 
\eeq
where $\kp_*(x) = 1 + O(|x|^4)$ is defined in \eqref{eq:cutoff_near0_chi3} and  \(N\) stands for numerics. Here, we define $c$ \emph{analytically} instead of \emph{numerical evaluation},
so that  $ \pa_x (-\D) \bar \phi_{\mf{cor}}(0) = c $ 
\emph{exactly}. Then the error of solving the Poisson equation satisfies the \emph{exact quadratic vanishing} near $x=0$
\beq\label{eq:poisson_err}
\bar \e =  \bar \om + \D \bar \phi^N = O(|x|^2) .
\eeq
For our solution, $|c| < 10^{-10}$ is very small. Since the exact stream function \((-\Delta)^{-1}\bar\omega\) depends on \(\bar\omega\) nonlocally, we cannot construct it exactly. Instead, we use
\(\bar\phi^N\) as a numerical approximation. The nonlocal error $\bar \e$ is decomposed and estimated in Section \ref{sec:comb_vel_err}.

\vs{0.05in}

\paragraph{\bf Difference between $\bar{\e}$ and $\bar{\cF}_i$}

We distinguish the Poisson residual $\bar{\e}$ \eqref{eq:poisson_err}
from the profile residual errors $\bar{\cF}_i$ in \eqref{eq:bous_err}.
The error $\bar{\e}$ enters $\bar{\cF}_i$ only through the nonlocal velocity
$\uu(\bar{\e})$:
\[
\bal
\bar{\cF}_1^{\nlocc}
&= -\uu(\bar{\e})\cdot \na\bar{\om}
= -|\xx|^{-1}\uu(\bar{\e})\cdot |\xx|\na\bar{\om},\\
\bar{\cF}_{i+1}^{\nlocc}
&= -\pa_i\uu(\bar{\e})\cdot \na\bar{\th}
-\uu(\bar{\e})\cdot \na\pa_i\bar{\th}  = -\pa_i\uu(\bar{\e})\cdot \na\bar{\th}
-|\xx|^{-1}\uu(\bar{\e})\cdot |\xx|\na\pa_i\bar{\th},
\qquad i=1,2 .
\eal
\]
By estimates similar to \eqref{eq:elli_embed} and estimates \eqref{solu:decay}, $\na\uu(\bar{\e})$ and
$|\xx|^{-1}\uu(\bar{\e})$ have at least the same far-field decay rate as $\bar{\e}$.
Moreover, the profile decay \eqref{solu:decay} gives additional decaying factors
\[
|\xx|\,|\na\bar{\om}|\les \jan{\xx}^{\bar{\al}_1},
\qquad
|\na\bar{\th}|+|\xx|\,|\na\pa_i\bar{\th}|
\les \jan{\xx}^{2\bar{\al}_1},
\qquad \bar{\al}_1\approx -1/3.
\]
Due to these decaying factors, for very large $|x|$, the nonlocal error $\bar{\cF}_i^{\nlocc}$, and consequently the profile residual $\bar{\cF}_i$, can be much smaller than the Poisson residual $\bar{\epsilon}$. 
See \cite{ChenHou2026noteErrors} for more discussion.

\appendix

\section{Lemmas for stability estimates}\label{sec:lem_stability}

We use Lemma \ref{lem:PDE_stab} for the linear stability estimates in 
the weighted $L^\infty$ and H\"older energy estimate, and Lemma \ref{lem:PDE_nonstab} 
for nonlinear stability estimates.

\begin{lem}\label{lem:PDE_stab}
Suppose that $f_i(x, z, t) : \R^2_{++} \times \R^2_{++} \times [0, T] \to \R, 1\leq i \leq n$, satisfies 
\beq\label{eq:PDE_stab_1}
\pa_t f_i + v_i(x, z, t) \cdot \na_{x, z} f_i = -a_{ii}(x, z, t)  f_i+ B_i(x, z, t),
\eeq
where $v_i(x, z, t) \in \R^4$ are some  Lipschitz  vector fields in $x, z$ with 
$v_i \cdot \nn = 0$ on $x_j=0$ or $z_j=0$ with $j =1,2$, where $\nn$ is the associated unit outer normal,  
and $B_i$ satisfies the following estimate 
\beq\label{eq:PDE_stab_2}
|B_i(x, z, t) |\leq  \sum\nolim_{j \neq i} |a_{ij}(x, z, t)| \cdot || f_j||_{L^{\inf}}.
\eeq
If there exists some constants $M ,\lam, \mu_i >0 $ such that for all $(x, z)$ and $t\in[0, T]$, we have
\beq\label{eq:PDE_diag}
 a_{ii}(x,z, t)  - \sum\nolim_{j \neq i} |a_{ij}(x, z, t)| \mu_i \mu_j^{-1} \geq \lam , \quad 
\sum\nolim_{j \neq i} \mu_i \mu_j^{-1} |a_{ij}(x, z, t)| \leq M .
 \eeq
 Then for $E(t) = \max_i( \mu_i || f_i(t)||_{\infty})$, which is Lipschitz,  and $  0 \leq t_0 < t \leq T$ , we have 
 \footnote{
 The energy $E(t) $ implies $ || f_i(t)||_{\infty} \leq \mu_i^{-1} E(t)$. 
 The parameters $\mu_i$ capture the relative size of $f_i(t)$. 
 }
\[
E(t) \leq e^{- \lam(t - t_0)} E_0 ,\quad E_0 = E(t_0).
\]

\end{lem}

The condition \eqref{eq:PDE_diag} means that the damping term dominates the bad terms, leading to the stability. We apply $f_i(x, z, t) = ( ( S_i \psi_i)(x) - (S_i \psi_i)(z)) g_i(x, z) $ in the weighted H\"older estimate, and $f_i(x, z, t) = (S_i \vp_i)(x)$ in the weighted $L^{\inf}$ estimate, $S_1 = \om, S_2  =\eta, S_3 = \xi$. 
In the weighted $L^{\inf}$ estimate, we do not need the extra variable $z$, and $f_i$ is constant in $z$. 
We also perform energy estimates on some scalars $a_i(t)$ and choose $f_i(x, z, t) = a_i(t)$. In this case, the advection term vanishes, and $a_{ii}, a_{ij}, B_i$  depend only on $t$.
For the Boussinesq equations \eqref{eq:lin_main}, we choose 
\[
b(x, t) =  \bar c_l x + \bar \uu(x) + \uu(\om)(x, t), \quad 
v_i(x, t) = b(x, t), \quad \mathrm{or} \quad v_i(x, z, t) = (b(x, t), b(z, t)).
\]

\begin{proof}
For simplicity, we assume that the condition \eqref{eq:PDE_diag} holds for $\mu_i = 1$. Otherwise, we can estimate the renormalized variables $\mu_i f_i$ and introduce $\td a_{ij} = a_{ij} \mu_i \mu_j^{-1}$.  Then the equations and estimates \eqref{eq:PDE_stab_1}, \eqref{eq:PDE_stab_2} become 
\[
\bal
 & \pa_t (\mu_i f_i ) + v_i(x, z, t) \cdot \na_{x, z} (\mu_i f_i) = -a_{ii} (\mu_i f_i )+ \mu_i B_i(x, z, t)   \\
 & \mu_i |B_i(x, z, t)|  \leq  \tts{ \sum_{j \neq i} }  \mu_i \mu_j^{-1} | a_{ij}(x, z, t) | \cdot \mu_j || f_j ||_{L^{\inf}} = \tts{ \sum_{j \neq i} } | \td a_{ij}(x, z, t) | \cdot \mu_j || f_j ||_{L^{\inf}} .
\eal
\] 
The condition \eqref{eq:PDE_diag} for $a_{ij}$ becomes the condition for $\td a_{ij}$ with equal weights. Thus, it suffices to consider the case $\mu_i = 1$, $1\le i\le n$.

Formally, we can perform $L^{\inf}$ estimate on  \eqref{eq:PDE_stab_1} and then evaluate \eqref{eq:PDE_stab_1} at the maximizer to obtain the desired result. To justify it rigorously, we use the characteristics, Duhamel's principle, and a bootstrap argument. 
We define the characteristics associated with $v_i$
\beq\label{eq:PDE_stab_char1}
 \tf{d}{dt}(X_i(t), Z_i(t)) = v_i( X_i, Z_i, t), \quad X_i( t_0) = x_0, \quad Z_i(t_0) = z_0, \quad 
 F_i(t) = f_i( X_i(t), Z_i(t), t)
\eeq
To simplify the notation, we drop $x_0, z_0$. Denote 
\beq\label{eq:PDE_stab_char2}
A_i(t) = a_{ii}( X_i(t), Z_i(t), t) , \quad C_i(t) \teq \tts{ \sum_{j \neq i} } |a_{ij}(X_i(t), Z_i(t), t)| .
\eeq

It suffices to prove that for small $\e >0$, we have
\beq\label{eq:PDE_stab_boot1}
E(t) \leq (1 + c(M) \e) e^{ -\lam_{\e}  (t-t_0) } E_0, \quad 
\lam_{\e} = \lam - \e, \quad 
c(M) = \tf{2}{M}, 
\eeq
where $M$ is the upper bound in \eqref{eq:PDE_diag}. Then taking $\e \to 0$ completes the proof.

We use a bootstrap argument to prove \eqref{eq:PDE_stab_boot1}. First, since $E( t_0) = E_0$ and $E(t)$ is Lipschitz,  condition \eqref{eq:PDE_stab_boot1} holds for $t \in [t_0, t_0 + T_1]$ with some $T_1>0$. 
Below, under \eqref{eq:PDE_stab_boot1}, we  improve
\beq\label{eq:PDE_stab_boot2}
E(t) \leq (1 + c(M) \e/ 2 ) e^{ -\lam_\e  (t-t_0) } E_0. 
\eeq

By definition, along the characteristics, we get 
\[
 \tf{d}{dt} F_i(t) = -A_i(t) F_i(t) + B_i(X_i(t), Z_i(t), t), 
 \]

Using the estimates \eqref{eq:PDE_stab_2} and definition \eqref{eq:PDE_stab_char2}, we yield 
\[
 |B_i(X_i(t), Z_i(t), t)| \leq \tts{ \sum_{ j \neq i} }  |a_{ij}| E(t) \leq C_i(t) E(t).
\]

Using Duhamel's principle and the above estimate, we obtain 
\beq\label{eq:PDE_duham}
\bal
F_i(t)  & = e^{ \int_{t_0}^{t}  -A_i(s) ds } F_i(t_0) 
 + \int_{t_0}^t  B_i(X_i(s), Z_i(s), s) e^{ \int_{s}^t - A_i(\tau) d\tau  } ds ,  \\
 |F_i(t)| & \leq e^{ \int_{t_0}^{t}  -A_i(s) ds } |F_i(t_0) |
 + \int_{t_0}^t C_i(s) E(s) e^{ \int_{s}^t - A_i(\tau) d\tau  } d \, s \teq I + II.
 \eal
\eeq

For the second term, using the bootstrap assumptions \eqref{eq:PDE_stab_boot1}, we yield 
\[
|II| \leq  \int_{t_0}^t
C_i(s) (1 + c(M) \e) e^{ -\lam_\e (s - t_0)} E_0 e^{ - \lam_{\e}(t-s) -  \int_s^t (A_i( \tau ) - \lam_{\e}) d \tau }  .
\]

Using $C_i(s) \leq M, C_i(s) \leq A_i(s) - \lam$ \eqref{eq:PDE_diag} and the choice of $c(M) = \f{2}{M}$ in \eqref{eq:PDE_stab_boot1}, we get 
\[
\f{C_i(s) }{  C_i(s) +\e }
\leq \f{ M }{ M + \e} \leq  \f{ 1 +  \e / M}{ 1 + 2 \e / M} = \f{ 1 + c(M) \e / 2 }{ 1 + c(M) \e},
  \quad C_i(s) + \e \leq A_i(s) - \lam + \e = A_i(s) - \lam_{\e},
\]
which implies 
\[
C_i(s) (1 + c(M) \e) \leq  (C_i(s) + \e ) (1 + c(M) \e / 2) \leq (A_i(s) - \lam_{\e}) (1 + c(M) \e / 2).
\]
Hence, we can simplify the bound of $II$ as follows 
\[
\bal
|II| &\leq  (1 + c(M) \e / 2) e^{-\lam_{\e}(t-t_0)} E_0 \int_{t_0}^t  (A_i(s) - \lam_{\e}) 
e^{ - \int_s^t (A_i(\tau) - \lam_{\e}) d \tau} d s \\
&=  (1 + c(M) \e / 2) e^{-\lam_{\e}(t-t_0)} E_0
( 1 - e^{ - \int_{t_0}^t (A_i(\tau) - \lam_{\e}) d \tau} ) .
\eal
\]

The estimate of $I$ is straightforward. Since $|F_i(t_0)| \leq E_0$, we have 
\[
|I | = 
|F_i(t_0)| e^{-\lam_{\e}(t-t_0)} e^{ - \int_{t_0}^t (A_i(\tau) - \lam_{\e}) d \tau} 
\leq E_0 e^{-\lam_{\e}(t-t_0)} e^{ - \int_{t_0}^t (A_i(\tau) - \lam_{\e}) d \tau} ,
\]
which along with the estimate of $II$ yields
\[
 |F_i(t)| \leq |I| + |II| \leq  (1 + c(M) \e / 2) e^{-\lam_{\e}(t-t_0)} E_0.
\] 

Since the above estimate holds for any initial data $x_0, z_0$ and $i$, taking the supremum, we prove \eqref{eq:PDE_stab_boot2}. Then the standard bootstrap argument implies the desired estimate \eqref{eq:PDE_stab_boot1}.
\end{proof}

The next lemma is the nonlinear version of Lemma~\ref{lem:PDE_stab}. It is designed for the
bootstrap argument in Sections~\ref{sec:EE}, \ref{sec:euler}: 
\(B_i\) are linear  terms, \(N_i\) are nonlinear terms, and \(\bar\varepsilon_i\) are residual errors.

\begin{lem}\label{lem:PDE_nonstab}
Suppose that $f_i(x, z, t) : \R^2_{++} \times \R^2_{++} \times [0, T] \to \R, 1\leq i \leq n$, satisfies 
\beq\label{eq:PDE_nonstab_1}
\pa_t f_i + v_i(x, z, t) \cdot \na_{x, z} f_i = -a_{ii}(x, z, t)  f_i+ B_i(x, z, t) + N_i(x, z, t) + \bar \e_i,
\eeq
where $v_i(x, z, t)$ are some vector fields Lipschitz in $x, z$  and satisfy $v_i\cdot\nn=0$ on $x_j=0$ or $z_j=0$, $j=1,2$, where $\nn$ is the corresponding outer normal.
 For some $\mu_i > 0$, we define the energy 
\[
E(t) = \max_{1 \leq i \leq n} (\mu_i || f_i||_{L^{\inf}}),
\quad E_{\mf{max} }(t) \teq \sup_{s \leq t} E(s).
\]
Suppose that $B_i, N_i$ and $\bar \e_i$ satisfy the following estimate for any $1\leq i \leq n$
\beq\label{eq:PDE_nonstab_2}
\bal
\mu_i (|B_i(x, z, t) | + |N_i(x, z, t)| +| \bar \e_i|) \leq  \sum_{j \neq i}  ( |a_{ij}(x, z, t)| E_{\mf{max}}(t) 
+ | a_{ij, 2}(x, z, t)| E_{\mf{max}}^2(t) + | a_{ij, 3}(x, z, t)|  ).
\eal
\eeq

If there exists some $E_*, \e_0, M > 0$ such that 
\beq\label{eq:PDE_nondiag}
\bal
  a_{ii}(x, z, t) E_*  - 
\tts{\sum_{j \neq i}}
  ( |a_{ij}(x, z, t)|  E_* + |a_{ij, 2}(x, z, t)| E_*^2 
+ |a_{ij,3}(x, z, t) | ) & >  \e_0 ,   \\
  \tts{ \sum_{j \neq i} } ( |a_{ij}(x, z, t)|  E_* + |a_{ij, 2}(x, z, t)| E_*^2 
+ |a_{ij,3}(x, z, t) | ) & < M ,
\eal
\eeq
for all $x, z$, $t \in [0, T]$, and $i \leq n$. Then for $E(0)< E_*$, we have 
\beq
E(t) \leq \max( E(0), (1-\d_{\eqref{eq:PDE_nonboot_impr} } ) E_* ), \quad \d_{\eqref{eq:PDE_nonboot_impr}} = \f{\e_0}{ M + \e_0}, \quad  \forall  t \in [0, T] .
\label{eq:PDE_nonboot_impr}
\eeq
\end{lem}

The second inequality in \eqref{eq:PDE_nondiag} is only qualitative. Note that the factor $a_{ij}$ \eqref{eq:PDE_nonstab_2} for linear terms is different from that in \eqref{eq:PDE_stab_2}. We have combined the weight $\mu_j$ with $a_{ij}, j \neq i$ in \eqref{eq:PDE_nonstab_2}.

\begin{proof}
The proof is very similar to that of Lemma \ref{lem:PDE_stab}. We fix $E(0)$. Without loss of generality, we assume $\mu_i = 1$. Otherwise, we rewrite \eqref{eq:PDE_nonstab_1} in terms of $\mu_i f_i$. It suffices to prove that under the bootstrap assumption 
\beq\label{eq:PDE_nonboot}
E(t) < E_*,
\eeq
on $[0, T_1]$, for $\d$ chosen in \eqref{eq:PDE_nonboot_impr}, we can improve 
\beq\label{eq:PDE_nonboot2}
E(t) \leq \max( E(0), (1 - \d ) E_*), 
\quad  t \in [0, T_1] .
\eeq

Since $E(0) < E_*$ and $E(t)$ is Lipschitz, the bootstrap assumption 
\eqref{eq:PDE_nonboot} holds for some  $T_1 > 0$. 

We adopt most notations from the proof of Lemma \ref{lem:PDE_stab} but use
\[
C_i(t) \teq \tts{ \sum_{j \neq i} } ( |a_{ij}(X_i(t), Z_i(t), t)| E_{\mf{max}}(t) 
+ | a_{ij, 2}( X_i(t), Z_i(t), t)| E_{\mf{max}}^2(t) + | a_{ij, 3}( X_i(t), Z_i(t), t)|  ).
\]

Using these notations 
and estimates similar to those in the proof of Lemma \ref{lem:PDE_stab}, we obtain 
\[
|F_i(t)| 
\leq e^{ -\int_0^t A_i(s) ds} |F_i(0)|
+ \int_0^t C_i(s) e^{ - \int_s^t A_i(\tau) d \tau } ds.
\]

Using the bootstrap assumption and \eqref{eq:PDE_nondiag}, we obtain 
\[
C_i(s) < \min( M, A_i(s) E_* - \e_0 )  
\leq \f{\e_0}{M+ \e_0} \cdot M + \f{M}{ M + \e_0} \cdot (A_i(s) E_* - \e_0 )
= (1-\d) A_i(s)E_*, 
\]
where $\d = \f{\e_0}{M + \e_0}$.
Now, we obtain 
\[
\bal
|F_i(t) |  & \leq e^{ -\int_0^t A_i(s) ds } |F_i(0)|
+ \int_0^t ( A_i(s) E_* - \e_0 ) e^{ - \int_s^t A_i(\tau) d \tau } ds \\
& \leq 
e^{ -\int_0^t A_i(s) ds} | F_i(0)|
+ (1 -\d) \int_0^t A_i(s) E_* e^{ - \int_s^t A_i(\tau) d \tau } ds \\
&= e^{ -\int_0^t A_i(s) ds } | F_i(0) | + (1 - \d) E_* ( 1 -e^{ -\int_0^t A_i(s) ds} ) \\
& \leq \max( | F_i(0)| , (1 - \d) E_* )
\leq \max(  E(0)  , (1 - \d) E_* ).
\eal
\]

Taking the supremum over the initial data of the trajectory and $i$, we prove
\[
E(t) \leq \max(  E(0) , (1 - \d) E_* ) < E_* ,
\]
which is \eqref{eq:PDE_nonboot2}. The above estimate strictly improve \eqref{eq:PDE_nonboot}.
Using the bootstrap argument, 
we prove \eqref{eq:PDE_nonboot}. 
Then the above argument further implies  \eqref{eq:PDE_nonboot2}. We complete the proof.
\end{proof}

\subsection{Proof of Lemma \ref{lem:hol_comp}}\label{app:lem:hol_comp}

The proof is a direct commutator computation for the weighted difference
\(\delta(f\vp)(x,z)\), followed by multiplication by the H\"older weight
\(g(x-z)\).

\begin{proof}
 Using \eqref{eq:hol_comp0}, we first derive the equation for $f \vp$
\[
\pa_t ( f \vp) + b(x) \cdot \na (f \vp) = c(x) f \vp + 
( b \cdot \na \vp ) f  + 
 \cR \vp = d(x) f \vp + \cR \vp,
\]
where $d(x) = c(x) + \f{b \cdot \na \vp}{\vp}$ is defined in Lemma \ref{lem:hol_comp}. For $x, z \in \R^2_+$, we derive the equation of $\d( f \vp)(x, z) = f\vp(x) - f\vp(z)$:
\[
\pa_t \d(f \vp) (x, z, t) + b(x) \cdot \na_x ( f \vp)(x) - b(z) \cdot \na_z( f \vp)(z) =  (d \, f \vp)(x) -  ( d \, f \vp)(z) + \d (\cR \vp) .
\]
Since 
\[
\bal
\na_x (f \vp)(x)  & = \na_x(  (f \vp)(x) - (f\vp)(z) ) = \na_x \d (f \vp), \quad
\na_z (f \vp)(z) =  -\na_z( \d (f\vp)) , \\
 d f \vp(x) - d f \vp(z) & = d(x) ( f \vp(x) - f\vp(z)) + (d(x ) - d(z) ) f \vp(z) \\ 
 & = d(x) \d( f \vp)(x, z)  + (d(x ) - d(z) ) f \vp(z),
\eal
\]
 we obtain
\beq\label{lem:hol_comp_pf1}
\pa_t \d(f\vp)+ ( b(x) \cdot \na_x + b(z) \cdot \na_z) \d( f\vp)
 = d(x) \d(f\vp)(x,z) + (d (x) - d(z) ) ( f\vp)(z) + \d( \cR \vp).
\eeq
Since $g(h)$ is even in $h_1, h_2$, $\pa_i g$ is odd in $h_i$ and we have
\[
\bal
&( b(x) \cdot \na_x + b(z) \cdot \na_z) ( \d( f \vp) g(x-z) ) \\
= & g(x-z)  \cdot (  b(x) \cdot \na_x + b(z) \cdot \na_z) \d(f\vp) 
+ \d(f\vp) \cdot  ( b(x) \cdot \na_x + b(z) \cdot \na_z)  g(x-z)\\
=& g(x-z)  \cdot (  b(x) \cdot \na_x + b(z) \cdot \na_z) \d(f\vp)  + \d(f\vp) \cdot  ( b(x) -b(z)) \cdot (\na g)(x-z).
\eal
\]

We further multiply both sides of \eqref{lem:hol_comp_pf1} by $g(x-z)$ and use $F(x, z, t) = \d(f \vp)(x, z) g(x-z)$ and the above identity to yield
\[
\pa_t F +  ( b(x) \cdot \na_x + b(z) \cdot \na_z) F
= (d(x) + \f{ ( b(x) -b(z)) \cdot (\na g)(x-z)}{g(x-z)}  ) F 
+ ( (d (x) - d(z) ) ( f\vp)(z) + \d( \cR \vp) ) g(x-z),
\]
which concludes the proof of \eqref{eq:hol_comp}.
\end{proof}

\section{Additional proofs of sharp H\"older estimates}\label{app:sharp}

In this Appendix, we prove additional sharp H\"older estimates in Section \ref{sec:sharp}.
We have proved Lemmas \ref{lem:holx_ux}, \ref{lem:holy_ux} in Section \ref{sec:sharp} 
for sharp $C^{1/2}$ estimates of $u_x$. In Appendix~\ref{sec:vx_holx}, we prove the \(C^{1/2}_x\) estimates
for \(v_x, u_y\). In Appendix~\ref{hol:vx_cy}, we prove the 
\(C^{1/2}_y\) estimates for \(v_x, u_y\), where the boundary discontinuity requires a different
decomposition. In Appendix~\ref{app:OT_map}, we record the sign functions and transportation maps in these estimates. In  supplementary material II in Part II \cite[Section {\secholconst}]{ChenHou2023bSupp}, we estimate these explicit integrals using some integral formulas and numerical quadrature with computer assistance, and obtain rigorous upper bounds. The codes can be found in \cite{ChenHou2023code}.

\paragraph{\bf Transportation map}
We use $T$ to denote the transport map in $C_x^{1/2}$ estimates of $v_x, u_y$, and $\TTy$, written in \emph{sans-serif} font, for the map in $C_y^{1/2}$ estimates of $v_x, u_y$. 

\subsection{$C_x^{1/2}$ estimates of $v_x$ and $u_y$ }\label{sec:vx_holx}

We follow the ideas and argument in Section \ref{sec:idea_opt} to estimate the H\"older seminorm of $u_y, v_x$. Recall the kernel  for $u_y, v_x$:
\[
K_2(y) = \f{1}{2} \f{y_1^2 - y_2^2}{|y|^4} .
\]
First, we prove an identity for the principal value of the $K_2$-integral. By removing a thin strip in the \(y_1\)- or \(y_2\)-direction, 
we obtain two formulas with different local terms. 
In later optimal transport estimates for \(v_x\) and \(u_y\), 
we choose the formula that exploits more cancellations.

\begin{lem}\label{lem:pv}
Suppose that $f \in L^{\inf}$ and is H\"older continuous near $0$. For $0 < a, b < \infty$ and $Q = [0, a] \times [0, b], [0, a] \times [-b, 0], [-a, 0] \times [0, b]$, or $[-a, 0] \times [-b, 0]$, we have
\[
\bal
P.V. \int_Q K_2(y) f(y) dy 
= \lim_{\e \to 0} \int_{ Q \cap |y_1| \geq \e} K_2(y) f(y) dy 
- \f{\pi}{8} f(0) 
= \lim_{\e \to 0} \int_{ Q \cap |y_2| \geq \e} K_2(y) f(y) dy  + \f{\pi}{8} f(0).
\eal
\]

\end{lem}

In the strip $|y_1| \leq \e$, $K_2(y ) < 0 $ if $|y_2| > \e $. It contributes to $-\f{\pi}{8} f(0)$ in the first identity. In the strip $|y_2|\leq \e$, $K_2(y) >0$ if $|y_1| > \e$.
It contributes to $\f{\pi}{8} f(0)$ in the second identity.

\begin{proof}
As $K_2$ is even in $y_1, y_2$, we consider $Q=[0, a] \times [0, b]$ without loss of generality.  
We have 
\[
 P.V. \int_Q K_2(y) f(y) dy 
= \lim_{\e \to 0} ( \int_{ \e}^a \int_0^b + \int_0^{\e} \int_{\e}^b ) K_2(y) f(y) dy 
\teq \lim_{\e \to 0} (I_{\e} + II_{\e})
\]
We just need to compute $II_{\e}$. Since $f$ is H\"older continuous near $0$, we get 
\[
\lim_{\e \to 0} \int_0^{\e} \int_{\e}^b K_2(y) ( f(y) - f(0)) d y = 0.
\]
The first identity follows from 
\[
\bal
\lim_{\e \to 0} \int_0^{\e} \int_{\e}^b K_2(y) f(0) dy 
&=  f(0) \lim_{ \e \to 0 } \int_0^{\e} \f{1}{2 } \f{y_2}{y_2^2 + y_1^2} \B|_{\e}^b d y_1 \\
& 
= \f{f(0) }{2}\lim_{ \e \to 0 } \int_0^{\e}  \big( \f{b}{b^2 + y_1^2}  - \f{\e}{\e^2 + y_1^2} \big) dy_1
 = -  \f{\pi}{8} f(0).
 \eal
\]
The second identity follows from the same argument.
\end{proof}

Next, we perform the sharp H\"older estimates for $u_y, v_x$.
Without loss of generality, we assume $z_1 = -1, x_1 = 1$ and $z_2 = x_2 > 0$. Due to the boundary, we do not have translation symmetry of the kernel $K_2(y)$ in $y_2$ and cannot assume $x_2 = 0$.
We are going to estimate 
\beq\label{eq:vx_holest0}
v_x(z) - v_x(x) = \f{1}{\pi}  \int \big( K_{2,B}( y_1 + 1, y_2) - K_{2, B}(y_1-1, y_2)  \big) W(y_1, x_2 - y_2) dy -
\f{1}{2}  ( \om(z) - \om(x) ) ,
\eeq
where  $K_{2,B}(y) = K_2(y) {\bf1}_{|y_1|\leq B, |y_2|\leq B}$ is the localized version of $K_2$ over $[-B,B]^2$, and $W$ is the odd extension of $\om$ from $\R^2_+$ to $\R^2$ \eqref{eq:ext_w_odd}. Denote 
 \beq\label{eq:K2_B}
A = \min(B, x_2),  \quad K^+ \teq K_{2, B}(y_1 + 1, y_2), \quad K^-  \teq K_{2, B}(y_1 - 1, y_2),
 \quad \D(y) = K^+ - K^-.
 \eeq

We focus on $ B \geq 2 |x-z| =4 \geq 3$. It is easy to see that $\D$ is odd in $y_1$. Since the transportation cost in the $y$ direction is cheaper (we will choose $\tau < 1$ in Lemma \ref{lem:holx_uy} to capture the property that $[\om]_{C_y^{1/2}}$ enjoys better energy estimate than 
$[\om]_{C_x^{1/2}}$), we shall use the $Y$-transportation as much as possible to obtain a sharp estimate. Due to the presence of the boundary and the discontinuity of $W$ across the boundary, we partition the domain into the inner part and the outer part 
\beq\label{def:Om_in}
\Om_{in} \teq \{ y_2 \in [ - A,   A] \}, \quad \Om_{out} \teq \{ y_2 \notin [ - A,  A]\}.
\eeq
Then $\om (\cdot, x_2 - \cdot) \in C^{1/2}( \Om_{in})$. We add  parameter $A$ in these domains due to the localized kernel.

Note that for a fixed $y_1$, $\D(y_1, y_2)$ may not have a fixed sign in $y_2$. We define 
\begin{equation}
\label{eq:Del1D}
\D_{1D}(y_1) = \int_{  - A}^{  A} \D(y_1, y_2) d y_2.
\end{equation}

Denote the vertical line (\tbf{vl}) and the horizontal line (\tbf{hl})
\[
\mbf{vl}_{y_1} \teq \{ (y_1, y_2 ): y_2 \in \R \} , \quad \mbf{hl}_{x_2} \teq \{ (y_1, x_2) : y_1 \in \R \} .
\]

 The estimates consist of four steps. In Step 0, we identify the sign of the kernel $\D, \D_{1D}$.

In the first two steps, we estimate the integral in $\Om_{in}$. In the first step, we fix $y_1$ and consider the 1D transportation problem on the vertical line $\mbf{vl}_{y_1}$ by moving the positive part of $\D$ to its negative part. If $|y_1| \leq 9$, we move the remaining part with total mass $\D_{1D}$ to the horizontal line $\mbf{hl}_{x_2}$. In this step, the estimate is bounded by $ C [\om]_{C_y^{1/2}}$. See the blue arrows and the blue line in the left figure in Figure \ref{fig:vx_OT} for an illustration of the moving direction on $\mbf{vl}_{y_1}$.  

In the second step, we study the transportation problem on $\mbf{hl}_{x_2}$, and move the remaining part with total mass $\D_{1D}(y_1)$ for $|y_1| \leq 9$ in the first step horizontally. The estimate will be bounded by $C [\om]_{C_x^{1/2}}$ for some constant $C$. In the third step, we estimate the integral in the outer domain $ \Om_{out}$. The estimate will be bounded by $ C [\om]_{C_y^{1/2}}$ for some constant $C$.  Lastly, we make a small modification near the singularity to simplify the computer-assisted integral bounds.

We focus on $|y_2 | \leq B$ since otherwise $\D(y) = 0$. We assume $x_2 > 0$. The case $x_2 =0$ can be obtained by taking limit $x_2 \to 0$.

\subsubsection{Sign of $\D$ and $\D_{1D}$}\label{hol:vx_step0}

\begin{figure}[t]
\centering
\begin{subfigure}{.4\textwidth}
  \centering
  \includegraphics[width=0.8\linewidth]{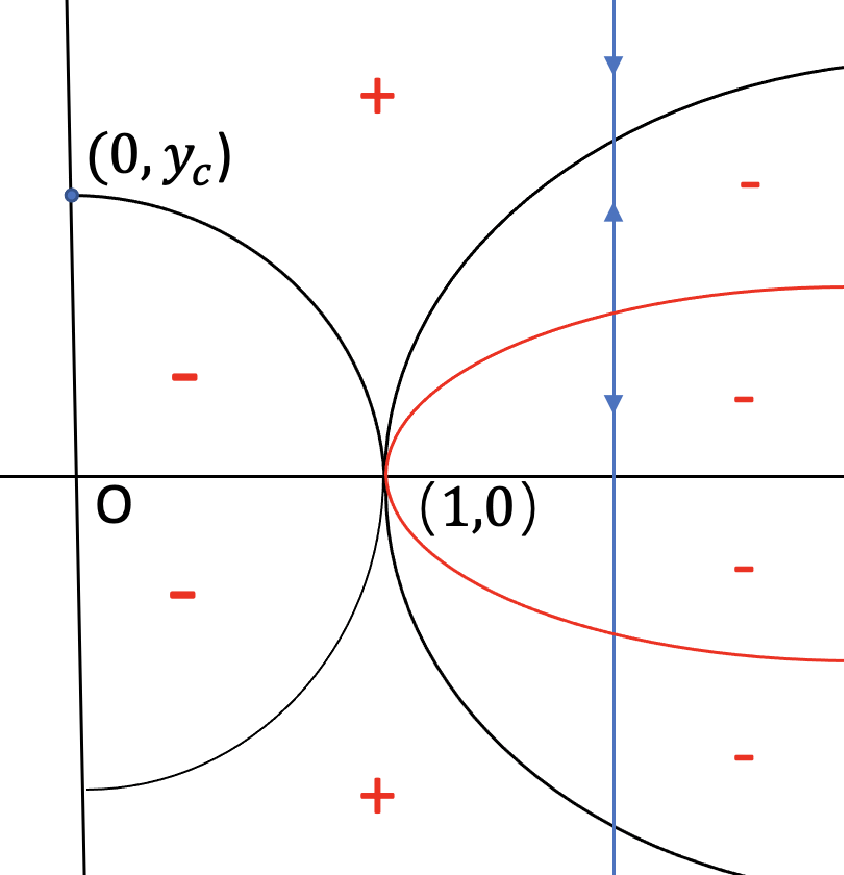}
    \end{subfigure}
    \begin{subfigure}{.59\textwidth}
  \centering
  \includegraphics[width=0.8\linewidth]{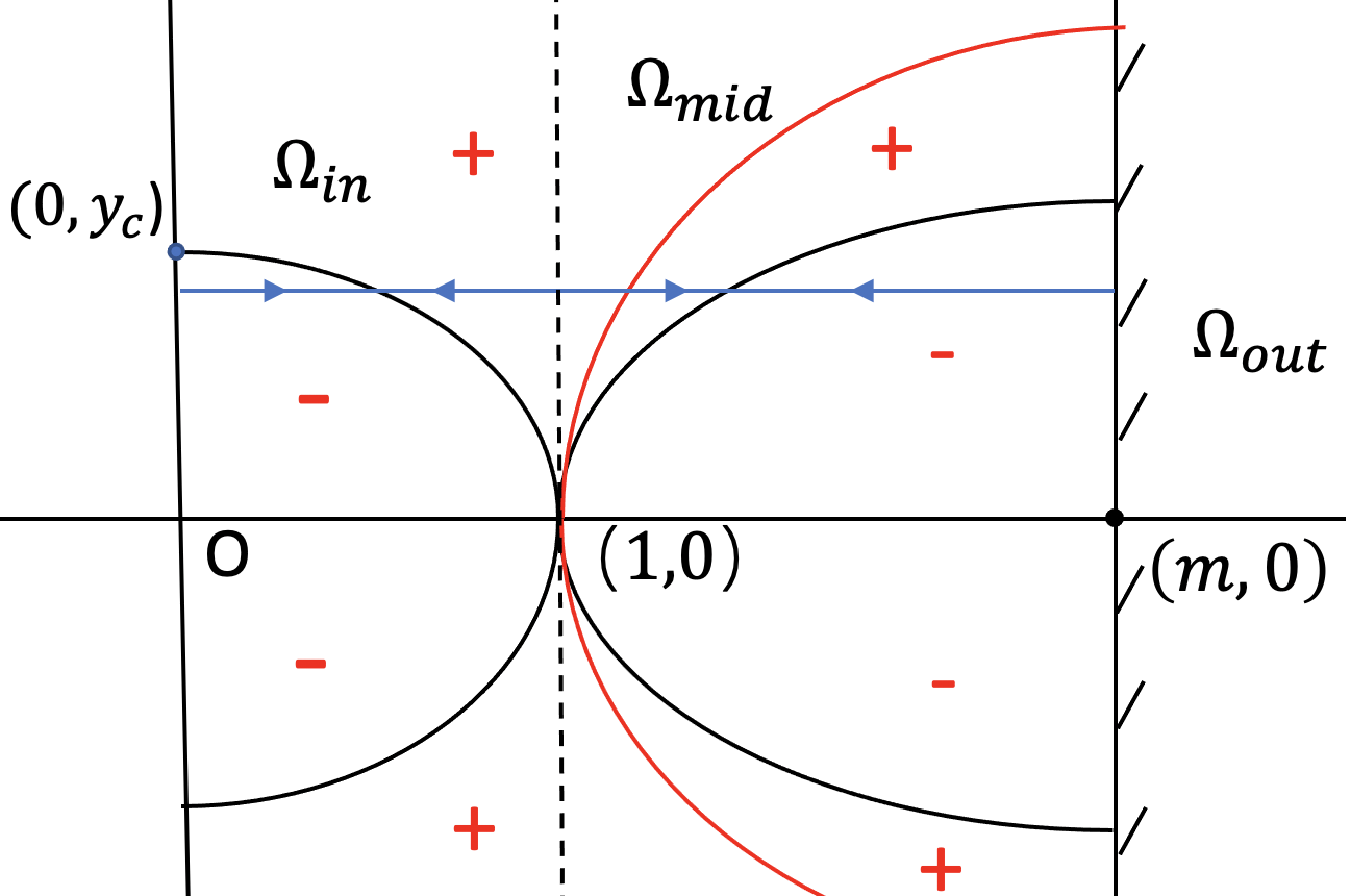}
    \end{subfigure}
  \captionsetup{width=0.9\textwidth} 
\caption{
Illustration of the sign of the kernel $\D(y)$ and transportation plan. The sign of $\D(y)$ in different regions is indicated by $\pm$. 
The blue arrows indicate the direction of 1D transportation plan. Left for $C_x^{1/2}$ estimate: The black curve and the red curve represent $ y_2 = \pm  s_c(y_1)$, $y_2 = \pm T(y_1, A)$ for $y_1 \geq 0$, respectively. Right for $C_y^{1/2}$ estimate: The black curve is for $y_2 = \pm s_c(y_1)$, or equivalently $y_1 = h_c^-(|y_2|)$ (two left black curves) and $y_1 = h_c^+(|y_2|) $ (two right black curves). The red curve represents $y_1 = \TTy(m, |y_2|)$. Note that these curves \textit{do not} agree with the actual functions.
 }
\label{fig:vx_OT}
\end{figure}

Due to the odd symmetry of $\D(y_1, y_2)$ in $y_1$, we focus on $y_1 \geq 0$. Solving $ K_2(y_1 + 1, y_2) - K_2(y_1 - 1, y_2) = 0$, we get $y_2 = s_c(y_1)$ \eqref{eq:vx_thres}. 
For $y_1 \neq 0$, we get
\beq\label{eq:vx_sign1}
\bal
K_2(y_1 + 1, y_2) - K_2(y_1 - 1, y_2) > 0 , \ | y_2 | > s_c(y_1),  \\
K_2(y_1 + 1, y_2) - K_2(y_1 - 1, y_2)  < 0 , \ | y_2 | < s_c(y_1) .
\eal
\eeq
See the left subplot in Figure \ref{fig:vx_OT} for an illustration of sign of $\D(y)$ in different regions. The black curve represents $ y_2 = s_c(y_1)$.
For $|y_1| \leq B-1$, we get 
\beq\label{eq:vx_del1}
\D(y) = K_2(y_1 + 1, y_2) - K_2(y_1 - 1, y_2) \teq \D_{all}(y).
\eeq
The sign of $\D(y_1, \cdot)$ is given above. For $ y_1 \in (B-1, B+1]$, we have 
\beq\label{eq:vx_del2}
\bal
\D(y) = - K_2(y_1 - 1, y_2) = - \f{1}{2} \f{ (y_1-1)^2 - y_2^2}{  ( (y_1-1)^2 + y_2^2)^2}. \\
\eal
\eeq
Since $B -1\geq 1$, for $y_1  \in (B-1, B+1]$, it satisfies 
\beq\label{eq:vx_sign2}
\D(y) > 0 , \ |y_2| > s_c(y_1) \teq y_1 - 1,
\quad \D(y) < 0,\  |y_2| < s_c(y_1) =  y_1 - 1 . 
\eeq
Similarly, we derive $s_c(y_1) = -y_1 - 1 $ for $y_1 \in [-B-1, 1-B)$. In Section \ref{sec:vx_holx}, we use  
\beq
s_c(y_1) \teq s_{c,  \eqref{eq:vx_thres}}(y_1), \ \mbox{for} \, |y_1| \leq B - 1,\quad  s_c(y_1) = |y_1|-1 , \ \mbox{for} \, |y_1 | \in (B-1, B+1] .
\label{eq:vx_sc_local}
\eeq

For $y_1 \geq B+1$, we have $\D(y) = 0$. Next, we compute $\D_{1D}$ defined in \eqref{eq:Del1D}. Since $\D$ is singular at $y = (\pm 1, 0)$ and $B> 2$, the singularity is in $J_1$ 
\beq\label{eq:vx_J1}
J_1 \teq [-9, 9] , \quad J_1^+ \teq J_1 \cap \R_+.
\eeq
In the inner part $\Om_{in}$ \eqref{def:Om_in}, for $y = (y_1, y_2)$, we have 
\beq\label{eq:vx_holest1}
S_{in} \teq \int_{ \Om_{in} }  \D(y) W(y_1, x_2 - y_2) dy  
= (\int_{ y_1 \in J_1} + \int_{ y_1 \notin J_1} ) \int_{-A}^A \D(y) W(y_1, x_2 - y_2) dy  
\teq S_1 + S_2 .
\eeq
By definition, we yield 
\begin{align}\label{eq:vx_holest2}
S_1 
& = \int_{ y_1 \in J_1} \int_{-A}^A  \D(y_1,  y_2) ( W(y_1, x_2 - y_2) - W(y_1, x_2) ) dy 
+   \int_{ y_1 \in J_1} \int_{-A}^A \D(y_1,  y_2) W(y_1, x_2) d y \notag \\
& \teq S_{11} + S_{12}.
\end{align}

For $S_{11}$, since $ |W(y_1, x_2 - y_2) - W(y_1, x_2) | \les |y_2|^{1/2}$, the integrand is locally integrable. We will estimate $S_{11}$ and $S_2$ in Section \ref{hol:vx_step1}. 

We should pay attention to the principal value in the singular integral in $S_{12}$ near the singularity $ (\pm 1, 0)$. 
Since $\D(y) = K_2(y_1 + 1, y_2) - K_2(y_1 - 1, y_2)$ near $y_1 = 1$, applying Lemma \ref{lem:pv} four times to $-K_2(y_1-1, y_2)$, which leads to $4 \cdot (-1) \cdot (-\f{\pi}{8} ) W(1, x_2) = \f{\pi}{2} W(1, x_2)$, we yield 
\beq\label{eq:vx_holest21}
\bal
 S_{12}^+  &\teq \int_{ y_1 \in J_1^+, |y_2|\leq A}   \D(y_1,  y_2) W(y_1, x_2) d y  \\
&=  \f{\pi}{2} W(1, x_2) + \lim_{\e \to 0} \int_{ y_1 \in J_1^+ \backslash [1-\e, 1+ \e] , |y_2|\leq A}   \D(y_1, y_2) W(y_1, x_2) d y .
\eal
\eeq

Recall the definition of $\D$ from \eqref{eq:kernel_du}, \eqref{eq:vx_del1}, \eqref{eq:vx_del2}. 
Denote 
\beq\label{eq:vx_del1D}
g_b(y) = \f{b}{y^2 + b^2}, \quad 
\D_{1D}(y_1) =  \f{A}{ (y_1+1 )^2 + A^2} \one_{ |y_1 + 1| \leq B} 
-\f{A}{ (y_1-1 )^2 + A^2} \one_{|y_1-1| \leq B} .
\eeq

Recall $\D$ and  $A \leq B$ from \eqref{eq:K2_B}. For $y_1 \notin [1-\e, 1+\e], |y_2| \leq A$, we have $K_{2,B}(y_1, y_2) = K_{2}(y) \one_{|y_1|\leq B}, K_2(y)  = \f{1}{2}\pa_{y_2} \f{y_2}{|y|^2} $ and 
\[
\bal
\int_{-A}^A \D(y_1, y_2) d y_2 
& = \f{1}{2} \B( \f{y_2 }{ (y_1+1)^2 + y_2^2} \one_{|y_1+1| \leq B} - \f{y_2 }{ (y_1-1)^2 + y_2^2} \one_{|y_1-1| \leq B} \B) \B|_{-A}^A  \\
& = g_A(y_1 + 1) \one_{|y_1+1| \leq B} - g_A(y_1 - 1) \one_{|y_1-1| \leq B} = \D_{1D}(y_1). 
\eal
\]

Plugging the above computation to the $P.V.$ integral yields 
\[
S_{12}^+ = \int_{J_1^+} \D_{1D}(y_1) W(y_1, x_2)  dy_1 + \f{\pi}{2} \om(1, x_2).
\]

The computation of the integral over $\R_-$ is similar due to symmetry. We yield a 1D-integral
\beq\label{eq:vx_holest_X1}
S_{12} = \int_{J_1 } \D_{1D}(y_1) W(y_1, x_2)  dy_1 + \f{\pi}{2} (\om(1, x_2) - \om(-1, x_2) ).
\eeq

\subsubsection{First step}\label{hol:vx_step1}

Next, we estimate $S_2$ \eqref{eq:vx_holest1} and $S_{11}$ \eqref{eq:vx_holest2}. 
Recall the sign of $\D$ from \eqref{eq:vx_sign1}, \eqref{eq:vx_sign2}
\beq\label{eq:D_vx_recall}
\D(y) > 0 , \ |y_2| > s_c(y_1), \quad \D(y) < 0, \ |y_2| < s_c(y_1),
\quad y_1 \neq 0.
\eeq
Since $\D$ is even in $y_2$ in $\Om_{in}$ and odd in $y_1$, we focus on the first quadrant. 

In this step, the variable \(y_1\) is fixed and the cancellation is performed
only in the vertical \(y_2\)-direction.  This is the direction in which we want
to use the ``cheaper" \(C^{1/2}_y\) control.

For a fixed $y_1 \geq 0$, we transport the positive part of $\D$ to its negative part on the line $\mbf{vl}_{y_1}$ in the first quadrant. We construct the transportation map $T(y)  > 0$ by solving 
\beq\label{eq:map_Cx_uy}
\int_{T(y)}^{ y_2} \D(y_1, s_2) d s_2 = 0. 
\eeq
By \eqref{eq:D_vx_recall} and \eqref{eq:map_Cx_uy}, 
for $y_1 \neq 0$, we get 
\beq\label{eq:Cx_uy_sign_order}
 A  \leq s_c(y_1) \leq T(y_1, A), 
 \  \mbox{for} \, A  \leq s_c(y_1), \qquad
   T(y_1, A)  \leq s_c(y_1) \leq A,  \ \mbox{for} \,  s_c(y_1) \leq A.
\eeq

For $y_1 \leq B-1$, $\D = K_2(y_1 + 1, y_2 ) - K_2(y_1 - 1, y_2)$. The map $T$ can be obtained from the cubic equation \eqref{eq:map_vx}. For $y_1 \in [ B-1, B+1]$, $\D = - K_2(y_1 - 1, y_2)$ and we get 
\beq\label{eq:map_K2_cy}
 0 = \int_{T(y)}^{ y_2} K_2(y_1-1, s_2) d s_2 
= \f{1}{2} \f{s_2}{ (y_1-1)^2 + s_2^2} \B|_{T}^{y_2}, \quad  T(y) = \f{ (y_1-1)^2 }{y_2}.
\eeq

Denote 
\[
\td W(y) = W(y_1, x_2 - y_2) - W(y_1, x_2) . 
\]

Using the above map, the estimates below,
\begin{align}\label{eq:vx_holest_cy1}
 |\td W(y_1, y_2) - \td W(y_1, T(y))| 
& = | W(y_1, x_2 - y_2) -  W(y_1, x_2 -T(y))|  \leq |y_2 - T(y)|^{1/2} [\om]_{C_y^{1/2}},  \notag \\
  |\td W(y_1, y_2)| & \leq |y_2|^{1/2} [\om]_{C_y^{1/2}}, 
\end{align}
and applying Lemma \ref{lem:trans} to the integral on $[T(y_1, A), A]$, we yield 
\beq\label{eq:vx_holest_cy2}
\bal
& \B|\int_0^A \D(y) \td W(y) d y_2 \B|
=\B| (\int_{T(y_1, A)}^A + \int_0^{T(y_1, A)} ) \D(y)  \td W(y) dy_2  \B|  \\
\leq & 
 \B( 
 \int_{ 0 }^A  \one_{y_2 \geq s_c(y_1)} |\D(y)|  |y_2 - T(y)|^{1/2}  dy_2
+ \int_0^A \one_{y_2 \leq T(y_1, A)} |\D(y)| |y_2|^{1/2} d y_2 \B) [\om]_{C_y^{1/2}}. 
\eal
\eeq
If $A\geq s_c(y_1)$, by $ T(y_1, A) \leq  s_c(y_1) \leq A $ \eqref{eq:Cx_uy_sign_order},
the ranges in the bounds reduces to $y_2 \in [s_c(y_1), A]$
and $y_2 \in [0, T(y_1, A)]$. If $A < s_c(y_1)$, 
we get  $ A\leq s_c(y_1) \leq T(y_1, A)$ by \eqref{eq:Cx_uy_sign_order}.
The first term in the upper bound vanishes. 
Using \eqref{eq:vx_holest_cy1}, we obtain \eqref{eq:vx_holest_cy2}. See the blue arrows in the left subplot in Figure \ref{fig:vx_OT} for an illustration of this transportation plan.

Due to the symmetry of $\D$ in $y_1, y_2$, we can estimate $S_{11}$ \eqref{eq:vx_holest2} as follows 
\beq\label{eq:vx_holest22}
|S_{11} | \leq  4 \int_{y_1 \in J_1^+} \B(  \int_{ 0 }^A 
\one_{y_2 \geq s_c(y_1)} 
|y_2 - T(y)|^{1/2} |\D(y) | dy_2  
+ 
 \int_0^A \one_{y_2 \leq T(y_1, A)}
  |y_2|^{1/2} |\D(y)|  dy_2 \B)   dy_1 [\om]_{C_y^{1/2}},
\eeq
where $J_1$ is defined in \eqref{eq:vx_J1} and the factor $4$ is due to the fact that we have $4$ quadrants.

The estimate of $S_2$ \eqref{eq:vx_holest1} is similar except that we \textit{do not} further transport the remaining negative part of $\D$ to the location $(y_1, x_2 )$. 
See Remark \ref{rem:restrict}. Consider the decomposition
\beq\label{eq:vx_holest23}
S_2 
=  \int_{  y_1 \notin J_1} ( \int_{ T(y_1, A) \leq  |y_2|  \leq A}
+  \int_{ |y_2| \leq  \min( T(y_1, A), A) } )  \D(y) W(y_1, x_2 - y_2) d y
\teq I + II. 
\eeq

For $I$, using symmetry of $\D(y)$ in $y_1, y_2$, we obtain 
\beq\label{eq:vx_holest3_I}
|I| \leq 2 \int_{ y_1 \notin J_1} \int_{0 }^A 
\one_{y_2 \geq s_c(y_1)}
|T(y ) - y_2|^{\f12 } |\D(y)| dy [\om]_{C_y^{\f12}}
= 4\int_{ y_1 \notin J_1, y_1 \geq 0} \int_{ 0 }^A 
\one_{y_2 \geq s_c(y_1)} |T(y ) - y_2|^{\f12} |\D(y)|  dy  [\om]_{C_y^{ \f12}} .
\eeq
If $A \leq s_c(y_1)$, we get $A \leq s_c(y_1) \leq T(y_1, A)$ by \eqref{eq:D_vx_recall}, \eqref{eq:map_Cx_uy}; both $y_2$-integrals in $I$ in \eqref{eq:vx_holest23}, \eqref{eq:vx_holest3_I} vanish. 
For $II$, we use the odd symmetry of $\D(y_1, y_2)$ in $y_1$ to get
\begin{align}\label{eq:vx_holest3}
|II|  &\leq \B| \int_{y_1 \notin J_1, y_1 \geq 0} \int_{ |y_2| \leq \min( T(y_1, A), A) } \D(y) (W(y_1, x_2 -y_2) - W(- y_1, x_2 -  y_2) ) \B| \\
&\leq \int_{y_1 \notin J_1, y_1 \geq 0} \int_{ |y_2| \leq \min( T(y_1, A) , A)}  \sqrt{2 y_1} |\D(y)| dy  
[\om]_{C_x^{1/2}}
 = \int_{y_1 \notin J_1, y_1 \geq 0} \sqrt{2 y_1} | \D_{1D}(y_1) | dy_1  [\om]_{C_x^{1/2}},  \notag
 \end{align}
where we used $\int_{T(y_1, A)}^A \D(y) dy_2 = 0$, $\D(y) \leq 0$ for $|y_2| \leq  \min( T(y_1, A),A) \leq  s_c(y_1)$ by \eqref{eq:Cx_uy_sign_order} and  
the following identity 
{\small
\beq\label{eq:holy_K2_rem}
\int_{ |y_2| \leq \min(T(y_1, A), A) } |\D(y) | dy_2 
=\B| \int_{|y_2| \leq \min(T(y_1, A), A) } \D(y) d y_2 \B|
= \B| \int_{ |y_2| \leq A}   \D(y)  dy_2\B| =  | \D_{1D}(y_1)|.
\eeq
}

\subsubsection{Second step: Estimate $S_{12}$ \eqref{eq:vx_holest_X1} }\label{hol:vx_step2}

The term \(S_{12}\) contains  1D residual mass left after the
vertical transportation in the inner region. We combine 
$S_{12}$ \eqref{eq:vx_holest_X1} and the local part of $v_x$ in the principal-value integral,  $- \f{\pi}{2}(\om(z) - \om(x) ) $  \eqref{eq:vx_holest0},
and perform horizontal transportation estimate.

For $v_x$, since $\om(z) - \om(x) = -\om(1, x_2) + \om(-1, x_2)$, we obtain 
\beq\label{eq:vx_holest_X2}
I \teq  S_{12} - \f{\pi}{2} ( \om(z) - \om(x) )
= \int_{J_1} \D_{1D}(y_1) W(y_1, x_2) d y_1 
+  \pi ( \om(1, x_2) - \om(-1, x_2)).
\eeq

Recall the definition of $\D_{1D}(y_1)$ \eqref{eq:vx_del1D}. Clearly, $\D_{1D}$ is odd and $\D_{1D} < 0$ for $y_1 > 0$. Note that for $k \in [0, 9]$, we have
\beq\label{eq:vx_holest_P}
P(k) \teq \int_{ k}^9 \D_{1D} d y_1\geq \int_{\R_+} \D_{1D}  d y_1
 = - \int_{-1}^1 \f{A}{y_1^2 + A^2} d y_1 = - 2 \arctan( \f{1}{A} ) \geq - \pi.
\eeq

We transport all the negative part of $\D_{1D}$ on $[1/9, 9]$ to $1$, 
and  transport all the positive part of $\D_{1D}$ on $(-9, -1/9]$ to $(-1)$. For $y_1 \in [-1/9, 0] \cup [0, 1/9]$, we move $y_1$ to $-y_1$. We do not move these parts to $\pm 1$ since  $ \sqrt{2 |y_1|} + \sqrt{2} \leq 2 |y_1-1|^{1/2}$ for $y_1 \leq 1/9$. See Remark \ref{rem:restrict}. Denote 
\[
J_2 = [1/9, 9] \subset J_1
\]
We derive the following estimate 
\begin{align}\label{eq:vx_holest_X3}
|I |
&=  \B|\int_{J_2 \cap \R_+} \D_{1D}(y_1) ( W(y_1, x_2)- W(1, x_2) ) d y_1  + \int_{ (- J_2) \cap \R_-} \D_{1D}(y_1) ( W(y_1, x_2)- W(-1, x_2) ) d y_1  \notag  \\
&\quad + \int_0^{1/9} \D_{1D}(y_1) (W(y_1, x_2) - W(-y_1, x_2) ) d y_1 
+ (\pi + P( \f{1}{9})) ( W(1, x_2)  - W(-1, x_2) )\B|  \notag \\
& \leq \B( 2 \int_{1/9}^9 | \D_{1D} |y_1-1|^{1/2} d y_1
+ \int_0^{1/9} |\D_{1D}| \sqrt{2 y_1} d y_1
+ (\pi + P( \f{1}{9})) \sqrt{2} \B) [\om]_{C_x^{1/2}} ,
\end{align}
where we have used the symmetry of $\D_{1D}$ to get the factor $2$.

\begin{remark}\label{rem:restrict}
We do not transport the negative part in $II$ in $S_{2}$ \eqref{eq:vx_holest23} to $(\pm y_1, x_2)$
(correspond to $(\pm z_1, 0 )$ below) for the following reasons. The remaining integrals in $S_2, S_{12}$ are similar to  
\[
 M = \int ( - \d_z + \d_{(1, 0 )} + \d_{ (-z_1, z_2)} - \d_{(-1, 0)} ) f(y) dy , \ f(y) = W( y_1, x_2 -y_2)
 \]
 for some $z_1 \geq 9, |z_2| \leq T(z_1, A) $. If we first move $(\pm z_1, z_2)$ to $( \pm z_1, 0)$ and then $(\pm z_1, 0)$ to $(\pm 1, 0)$, we obtain: 
$
 M     \leq  2|z_2 |^{1/2} [f]_{C_y^{1/2}} + 2 |z_1-1|^{1/2} [f]_{C_x^{1/2}} .
$
Factor $2$ comes from estimating both the $+$ case and the $-$ case.

We have another simple estimate:  $M \leq [f]_{C_x^{1/2}} (|2 z_1|^{1/2} + \sqrt{2})$. For $z_1 \geq 9$ or $z_1 \in [0,1/9]$, we get $2 |z_1 - 1|^{1/2}  \geq  \sqrt{2} + \sqrt{2 z_1}$, with equality for $z_1 = 9, 1/9$. 
Thus, the second estimate is better, and we use it in the estimates of $S_{12}, S_2$. This estimate also motivates the choice of $J_1$ \eqref{eq:vx_holest_X1}.

\end{remark}

\subsubsection{Third step}\label{hol:vx_step3}

It remains to estimate the integral in the outer part.  
If $B < x_2$, since $\D(y_1, y_2)$ is localized to $|y_2| \leq B = A $, the contribution from outer part $|y_2| > B$ is $0$. 

If $B > x_2 = A$, using the odd symmetry of $W(y)$ in $y_2$ \eqref{eq:ext_w_odd}: $W(y_1, x_2-y_2) = - W(y_1, y_2 - x_2)$  and the even symmetry of $\D$ in $y_2$, we yield 
\[
\bal
S_{out}  & = \int_{ \Om_{out}} \D( y) W(y_1, x_2 - y_2) d y
=\int_{\R} dy_1  ( \int_{x_2}^B   + \int_{-B}^{-x_2} ) \D(y) W(y_1, x_2 - y_2)  dy_2 \\
&= -\int_{\R}  dy_1 \int_{x_2}^B \D(y) \om (y_1, y_2 - x_2) dy_2 
+ \int_{\R} dy_1 \int_{x_2}^B \D(y) \om(y_1, x_2 + y_2) dy_2 .
\eal
\]

It follows 
\beq\label{eq:vx_holest4}
|S_{out}| \leq \one_{x_2 < B} \int_{\R} \int_{x_2}^B  |\D(y)| \cdot | \om(y_1, x_2 + y_2 ) -  \om(y_1, y_2 - x_2) \B| d y
\leq \one_{x_2 < B} \sqrt{2 x_2} \int_{\R} \int_{ x_2  }^B  |\D(y) | dy [\om]_{C_y^{1/2}}.
\eeq

\subsubsection{$C_x^{1/2}$ estimate of $u_y$ }\label{hol:uy}

Estimates of $u_y$ in Step 1-2 follow those for $v_x$ except that we \textit{do not} transport the remaining negative part of $\D(y)$ with $|y_2| \leq \min(T(y_1, A), A) $ to $(y_1, x_2)$ for \textit{any} $y_1 > 0$.
\footnote{
  We do not do so, since the remaining integrals in $S_{2}, S_{in}$ are similar to
$
\int ( \d_y - \d_{(-y_1, y_2)} ) W(z_1, x_2 - z_2)  d y$,  for $y$ in $\{ |y_2| \leq \min( T( |y_1|, A), A) \}$, and have a sharper 
bound $ |2 y_1|^{1/2} [\om]_{C_x^{1/2}}$. We apply this estimate in \eqref{eq:uy_holest4}.
}
The outer-part estimate in the third step is the same as that of $v_x$ in Section \ref{hol:vx_step3}.

Denote $J_{\e} =  [-1-\e, -1 + \e] \cup [1-\e, 1 + \e] $. Note that $\D = K_2(y_1 + 1, y_2) - K_2(y_1-1, y_2)$ has singularities at $(\pm 1, 0)$. Applying Lemma \ref{lem:pv} four times to $K_2(y_1+1, y_2)$ and $-K_2(y_1-1, y_2)$, respectively, we can rewrite $S_{in}$ \eqref{eq:vx_holest1} as follows 
\beq\label{eq:uy_holest1}
\bal
& S_{in}  \teq \int_{ \R } \int_{-A}^A \D(y) W(y_1, x_2 - y_2) dy  \\
&= 
 \f{\pi}{2} ( \om(1, x_2) - \om(-1, x_2))
 + \lim_{\e \to 0} 
 \int_{ J_{\e}^{c}} d y_1 \int_{-A}^A  
\D(y) W(y_1, x_2 - y_2) d y_2
\teq I + \lim_{\e \to 0}  II_{\e}.
\eal
\eeq

For $y_1 \notin J_{\e}$, we perform a decomposition 
\begin{align}\label{eq:uy_holest_dec1}
f(y_1) & \teq \int_{-A}^A \D(y) W(y_1, x_2 - y_2) d y_2  
= ( \int_{ T(y_1, A) < |y_2| \leq A}  + \int_{ |y_2| \leq \min( T(y_1, A), A) }  )\D(y) W(y_1, x_2 - y_2) d y_2  \notag \\
& \teq f_1(y_1) + f_2(y_1).
\end{align}
If $ A \leq s_c(y_1)$, by \eqref{eq:Cx_uy_sign_order},
we get $f_1(y_1) = 0$.  Otherwise, we estimate it using Lemma \ref{lem:trans} 
\beq\label{eq:uy_holest3}
|f_1(y_1)| \leq 2\int_{ 0 }^A \one_{y_2 \geq s_c(y_1) } | \D(y)| \cdot  |y_2 - T(y)|^{1/2} d y_2 \cdot  [\om]_{C_y^{1/2}}
 \eeq
Term $f_2(y_1)$ denotes purely negative part. 
Since $\D$ is odd in $y_1$ by \eqref{eq:vx_del1}, using  \eqref{eq:holy_K2_rem}, we get
 \beq\label{eq:uy_holest4}
 \bal
|f_2(y_1) + f_2(-y_1)|
&= \B| \int_{|y_2| \leq \min( T(y_1, A), A) } \D(y) (W(y_1, x_2 - y_2) -  W(-y_1, x_2 - y_2)) d y_2 \B| \\
&\leq \int_{|y_2| \leq \min(T(y_1, A), A)} \sqrt{2 y_1} |\D(y)| dy_2 [\om]_{C_x^{1/2}}
= |\D_{1D}(y_1)| \sqrt{2 y_1} [\om]_{C_x^{1/2}}.
\eal
 \eeq

Integrating \eqref{eq:uy_holest3}, \eqref{eq:uy_holest4} over $y_1 \notin J_{\e}$, we establish 
\[
|II_{\e}| \leq 4 \int_{ J_{\e}^c \cap \R_+} \int_{0 }^A
\one_{y_2 \geq s_c(y_1) }
 |\D(y)| \cdot |y_2 - T|^{1/2} d y 
[\om]_{C_y^{1/2}} 
+ \int_{ J_{\e}^c \cap \R_+} |\D_{1D}(y_1)| \sqrt{2 y_1} [\om]_{C_x^{1/2}}.
\]

Near the singularity $(1,0)$ of $\D(y)$, $|y_2 - T|^{1/2} \les |y_2|^{1/2} + |y_1- 1|^{1/2}$. Thus, the integrand in the first integral is locally integrable. Plugging the above estimate in \eqref{eq:uy_holest1}, we derive 
\beq\label{eq:uy_holest5}
\bal
& |S_{in} - \f{\pi}{2}( \om(1, x_2) - \om(-1, x_2)) |  \leq \lim\sup_{\e \to 0} \big( |I 
- \f{\pi}{2}( \om(1, x_2) - \om(-1, x_2) )|  
+ |II_{\e} | \big) \\
\leq  &
4 \int_{  \R_+} \int_{ 0}^A \one_{y_2 \geq s_c(y_1) } |\D(y)| |y_2 - T(y)|^{1/2} d y 
[\om]_{C_y^{1/2}} 
+ \int_{  \R_+} |\D_{1D}(y_1)| \sqrt{2 y_1} [\om]_{C_x^{1/2}}.
\eal
\eeq

Recall the definition of localized $u_y$ \eqref{eq:ux_local}. The term $ \f{\pi}{2} ( \om(1, x_2)  -\om(-1, x_2)  )$ in $S_{in}$ cancel the local term $ \f{\pi}{2} (\om(z) - \om(x))$ in $u_y(z) - u_y(x)$. Combining the above estimate and the estimate of $S_{out}$ in \eqref{eq:vx_holest4}, we prove the estimate of $u_y$.

\subsubsection{Modification near the singularity}\label{sec:holx_singu_md}

The modification introduced in this subsection is not needed for the analytic structure of
the proof; it is introduced to make the computer-assisted evaluation of the
explicit singular integrals more stable and easier to derive.

Near the singularity $s_* = (1, 0)$,  the integrand in the estimate in $S_{11}$ \eqref{eq:vx_holest2}, \eqref{eq:vx_holest22} is singular of order $|x|^{-3/2}$ and quite complicated. To ease our computation of the integral in supplementary material II of 
Part II \cite{ChenHou2023bSupp}, we use a simpler estimate 
in  $y \in [1-\d, 1 + \d] \times [0, A], \d < 1$ close to $s_*$.
Since $B \geq 3$, we have $y_1 \leq B-1, y_1 \leq 9$, $\D(y ) = \D_{all}(  y)$ \eqref{eq:vx_del1}, and $-K_2(y_1 -1, y_2)$ is the main term in $ \D(y)$. Instead of using \eqref{eq:vx_holest_cy2}, we separate two kernels and estimate 
\[
S^+_{11, \d } = I_+( \d) - I_-(\d),
\quad  I_{\pm} \teq \int_{1-\d}^{1 + \d} \int_0^{A} K_2(y_1 \pm 1, y_2) \td W(y) d y. 
\]

For $I_+$, the integrand is away from the singularity. Using \eqref{eq:vx_holest_cy1}, we get
\beq\label{eq:vx_sing_md1}
|I_+(\d)| \leq \int_{1-\d}^{1+\d } \int_0^A |K_2(y_1 + 1, y_2  )| y_2^{1/2} dy [\om]_{C_y^{1/2}}
\teq S_{in, y, \d, +} \cdot [\om]_{C_y^{1/2}}.
\eeq

In supplementary material of Part II \cite{ChenHou2023bSupp}, using an estimate similar to \eqref{eq:vx_holest_cy2}, we obtain
\beq\label{eq:vx_sing_md2}
|I_-(\d)| 
\leq  \B(  2 \int_{[0, \d] \times [0, A] \bsh Q_{a_{\d}}^+} |K_2(y)| y_2^{1/2} dy
+ 2 a_{\d}^{1/2} C_{K_2}  \B) [\om]_{C_y^{1/2}} 
\teq S_{in, y, \d, -} \cdot [\om]_{C_y^{1/2}},
\eeq
where 
\beq\label{eq:hol_singu_nota1}
\bal
 & a_{\d} = \min(A, \d), \quad Q_{a} \teq [-a, a] \times [0, a], \quad Q_a^+ \teq [0, a]^2,  \quad  C_{K_2}   \teq  C_{K_2, up} + C_{K_2, low}, \\
& C_{ K_2, up}   \teq \int_0^1 \int_{y_1}^1 |K_2(y)| \, | \f{y_1^2}{y_2} - y_2 |^{1/2} dy , \quad 
C_{K_2, low} \teq \int_0^1 \int_{0}^{y_1^2} |K_2(y)| \, |y_2|^{1/2} d y_2. \\
\eal
\eeq

We apply a similar modification in the estimate of $[u_y]_{C_x^{1/2}}$ \eqref{eq:uy_holest3}-\eqref{eq:uy_holest5} in the region $y \in [1-\d, 1 + \d] \times [-A, A]$. Using $\int_{T(y_1, A)}^A \D(y) d y_2 = 0$, we modify the decomposition \eqref{eq:uy_holest_dec1}
\[
\bal
& \int_{-A}^A \D(y) W(y_1 , x_2 - y_2) d y_2 
= \int_{-A}^A \B( \D(y) (W(y_1 , x_2 - y_2) - W(y_1, x_2) ) \B) d y_2  \\
& + \int_{|y_2| \leq \min(T(y_1, A), A)}  \D(y) W(y_1, x_2) d y_2  \teq \td f_1(y_1) + \td f_2(y_1) , 
\quad  \qquad
II_{i ,\d} \teq \int_{1-\d}^{1+\d} \td f_i(y_1)  d y_1 .
\eal
\]
We apply the above estimates of $S_{11,\d}^+$ to $II_{1, \d}$. For $\td f_2(y_1)$, using $|W(y_1, x_2) - W(-y_1, x_2)| \leq \sqrt{2 y_1} [\om]_{C_x^{1/2}}$, we obtain the same estimate as in \eqref{eq:uy_holest4}.

The above modification changes the estimate only by order $\d^{1/2}$, and we choose $\d=2^{-14}$.
We refer this estimate to supplementary material II \cite[Section {\secholconst}]{ChenHou2023bSupp} (attached to 
Part II \cite{ChenHou2023b}).

\vspace{0.1in}
\paragraph{\bf{Summary of the estimates of $[v_x]_{C_x^{1/2}}, [u_y]_{C_x^{1/2}}$ }} 
 We now collect the estimates obtained in the three regions into the final
computable upper bounds used in the stability estimates. For $[v_x]_{C_x^{1/2}}$, combining \eqref{eq:vx_holest1}, \eqref{eq:vx_holest22}, \eqref{eq:vx_holest3_I}, \eqref{eq:vx_holest3}, \eqref{eq:vx_holest4}, \eqref{eq:vx_holest_X3}, 
\eqref{eq:vx_sing_md1}, \eqref{eq:vx_sing_md2}, we establish
{\small
\begin{align}\label{eq:holx_vx}
& \f{1}{\sqrt{2}}| v_x(-1, x_2) - v_x(1, x_2)|  \leq  \f{1}{\pi \sqrt{2}}  \big( (S_{in, x} + S_{1D} ) [\om]_{C_x^{1/2}} + ( S_{in, y} + S_{out} ) [\om]_{C_y^{1/2}} \big) ,\\
& S_{in, x}  =   \int_{y_1 \notin J_1, y_1 \geq 0} \sqrt{2 y_1} | \D_{1D}(y_1) | dy_1 , 
\qquad S_{out}   = \one_{x_2 < B} \sqrt{2 x_2} \int_{\R} \int_{ x_2  }^B  |\D(y) | dy , \notag \\
& S_{in, y} =  4 \int_{J_1^+ \bsh [1-\d, 1+ \d] } \B(  \int_{ 0 }^A 
\one_{y_2 \geq s_c(y_1)}  |y_2 - T(y)|^{ \f{1}{2} }  |\D(y)| dy_2  +  \int_0^{ A } \one_{y_2 \leq T(y_1, A)}
 |y_2|^{ \f{1}{2} }  |\D(y)| dy_2 \B)    dy_1  \notag  \\
& \qquad \qquad + 4 ( S_{in, y, \d, +} + S_{in, y, \d, -}) + 4\int_{ y_1 \notin J_1, y_1 \geq 0} \int_{ 0}^A 
\one_{y_2 \geq s_c(y_1)}  |T(y ) - y_2|^{ \f{1}{2} } 
|\D(y)| dy   \notag \\
&\qquad = 4 \int_{\R^+ \bsh [1-\d, 1+ \d] }   \int_{ 0 }^A 
\one_{y_2 \geq s_c(y_1)}  |y_2 - T(y)|^{ \f{1}{2} }  |\D(y)| dy_2 d y_1  + 
4 \int_{J_1^+ \bsh [1-\d, 1+ \d] } \int_0^{ A } \one_{y_2 \leq T(y_1, A)}
 |y_2|^{ \f{1}{2} }  |\D(y)| dy_2     dy_1  \notag  \\
& \qquad \qquad + 4 ( S_{in, y, \d, +} + S_{in, y, \d, -}) , \notag \\ 
&
 S_{1D}  =    2 \int_{ \f{1}{9}}^9 | \D_{1D} |y_1-1|^{1/2} d y_1
+ \int_0^{ \f{1}{9}} |\D_{1D}| \sqrt{2 y_1} d y_1
+ (\pi + P( \f{1}{9})) \sqrt{2} , \notag
\end{align}
}
where 
$[1-\d, 1+\d] \subset J_1^+ $ \eqref{eq:vx_J1} for $\d<1$ and  $P(\cdot )$ is from \eqref{eq:vx_holest_P}. The factor $\f{1}{\sqrt{2}}$ is from $\f{1}{|x-z|^{1/2}} = \f{1}{\sqrt{2}}$. 
For $u_y$, combining \eqref{eq:uy_holest5} and \eqref{eq:vx_holest4}, 
and using $S_{out}$ defined above, we prove 
\begin{align}\label{eq:holx_uy}
&  \f{1}{\sqrt{2}} | u_y(-1, x_2) - u_y(1, x_2)|  \leq  \f{1}{\pi \sqrt{2}}  \big(  \td S_{1D}  [\om]_{C_x^{1/2}} +  (  S_{up} +  S_{out} ) [\om]_{C_y^{1/2}} \big) , \\
& S_{up} = 4 \int_{ y_1 \in \R_+ \bsh [1-\d, 1 + \d] }  \int_{ 0}^A 
\one_{y_2 \geq s_c(y_1) }
| T(y) - y_2|^{ \f{1}{2} }  |\D(y)| dy_2 dy_1 +  4 ( S_{in, y, \d, +} + S_{in, y, \d, -})  ,  \notag \\
&  \td S_{1D} = \int_{\R_+} |\D_{1D}(y_1)| \sqrt{2 y_1}  d y_1 . \notag 
\end{align}

Since $A =\min(x_2, B)  $,
these bounds read $C_1(x_2, B) [\om]_{C_x^{1/2}} + C_2(x_2, B) [\om]_{C_y^{1/2}}$. 
We  bound it:
\[
C_1(x_2, B) [\om]_{C_x^{1/2}} + C_2(x_2, B) [\om]_{C_y^{1/2}} 
\leq (C_1(x_2, B) + C_2(x_2, B) \tau) \max( [\om]_{C_x^{1/2}}, \tau^{-1} [\om]_{C_y^{1/2}} ),
\]
for any $\tau >0$, 
We partition the  parameter ranges and use monotonicity of the integrals in $A, B$ to obtain the uniform bound. See supplementary material of Part II \cite[Section {\secholconst}]{ChenHou2023bSupp} for details.

\subsection{$C_y^{1/2}$ estimate of $v_x$ and $u_y$}\label{hol:vx_cy}

Since $W$ is discontinuous across $y_2=0$, the localized $v_x$ or $u_y$ is not in $C_y^{1/2}$. Thus, we study the corresponding estimate without localization. 
Without loss of generality, we assume $z_2 = m + 1, x_2 = m - 1, x_1 = z_1 =0 $ with $m > 1$. The case $m=1$ can be obtained by taking a limit. The difference $v_x(z) - v_x(x)$ or $u_y(z) - u_y(x)$ is given by
\[
\bal
I\teq & \f{1}{\pi} \int \big( K_2( y_1,  1 + y_2 ) - K_2(y_1,   y_2 - 1 ) \big)  W(y_1, m - y_2) d y + s ( \om(z) - \om(x) ) , 
\eal
\]
where $s = \f{ 1}{2}$ for $u_y$ and $s = -\f{ 1}{2}$ for $v_x$. Denote 
\beq\label{eq:vx_eta}
\eta( y_1, y_2) = W( y_2, y_1), \quad \eta_m(y_1, y_2 ) = \eta(m-y_1, y_2).
\eeq

By definition, $\eta$ is odd in $y_1$,
$\eta_m$ is discontinuous across $y_1 = m$, and satisfies 
\beq\label{eq:vx_rota}
[\eta]_{C_x^{1/2} (\Om)} =  [\om]_{C_y^{1/2}} , \quad [\eta]_{C_y^{1/2} (\Om)} =  [\om]_{C_x^{1/2}},
\eeq
for $\Om = \{  y: y_1 \geq 0 \}$ or $\Om = \{  y: y_1 \leq 0 \}$ .

Swapping the dummy variables $y_1, y_2$ and then using $K_2(y_1, y_2) = - K_2(y_2, y_1)$, $\eta_m(y) = \eta(m-y_1, y_2) = W(y_2, m - y_1)$, we yield 
\beq\label{eq:vx_cy_holest1}
\bal
I  &= -\f{1}{\pi} \int ( K_2(  y_1 + 1, y_2 ) - K_2( y_1-1, y_2) ) \eta( m - y_1, y_2) dy
 + s ( \om(z) - \om(x) )  \\
 & = -\f{1}{\pi} \int \D(y) \eta_m(  y) dy + s ( \eta_m(-1, 0) - \eta_m(1, 0) ) , \quad s = -\f{ 1}{2} \mathrm{ \ for \ } v_x, \ s = \f{1}{2} \mathrm{ \ for  \ } u_y.
 \eal
\eeq
We perform the above reformulation so that we can adopt the analysis of 
\beq\label{eq:K2_del}
\D(y) =  K_2(  y_1 + 1, y_2 ) - K_2( y_1-1, y_2)
\eeq
in \eqref{eq:vx_sign1} and Section \ref{sec:vx_holx}. Since $\eta$ is discontinuous across $y_1 = 0$, which relates to $y_1 = m$ in the integral in \eqref{eq:vx_cy_holest1}, and the singularity of $\D$ is at $y = (\pm 1, 0)$, we decompose the integral into the inner region, the middle region, and the outer region 
\beq\label{eq:vx_cy_region}
\bal
& \Om_{in} \teq \{ y: |y_1| \leq 1 \},  \quad \Om_{mid} \teq 
\{ y: |y_1| \in [1, m] \}, \quad \Om_{out} \teq \{ y: |y_1| > m \}, \\
& S_{\al } \teq \int_{ \Om_{\al}} \D(y) \eta_m(y) dy , \quad \al \in \{ in, \, mid, \, out \}.
\eal
\eeq
See the right panel of Figure \ref{fig:vx_OT} for different regions in $\{ y_1 \geq 0 \}$. In each region, $\eta_m$ is H\"older continuous. Since we can obtain a smaller factor from $[\om]_{C_y^{1/2}}$ than $[\om]_{C_x^{1/2}}$, 
and in view of \eqref{eq:vx_rota},  we should use the $x$-transportation as much as possible to obtain a sharp estimate of \eqref{eq:vx_cy_holest1}.

\paragraph{\bf Ideas} 
The three regions are treated differently.  In \(\Omega_{\rm in}\), we evaluate the principal-value integral at $(1,0)$ and apply an \(x\)-direction transportation estimate.  In \(\Omega_{\rm mid}\), we use two
different estimates: one adapted to \(u_y\) and one adapted to the cancellation
in \(v_x\).  
In \(\Omega_{\rm out}\), the kernel is nonsingular, and the
estimate follows from the same tail bound used in the second estimate 
in \(\Omega_{\rm mid}\) for $v_x$.

First, we analyze the sign of $\D(y)$. Since $\D$ is odd in $y_1$ and even in $y_2$, we can focus on $y_1 , y_2 \geq 0$. For a fixed $y_2$, we have  
\beq\label{eq:vx_cy_sign}
\bal
&\D(y) < 0 , \  y_1 < h_c^-(y_2) \leq 1,  \quad
\D(y) > 0, \  y_1 \in ( h_c^-(y_2) , 1) , \quad y_2 \leq y_c \teq  3^{-1/2},  \\
& \D(y) < 0, \  y_1 > h_c^+(y_2) \geq 1, \quad \D(y) > 0, \   y_1 \in (1,  h_c^+(y_2)) ,
\eal
\eeq
where $h_c^{\pm}(y_2)$ solves $\D( h_c^{\pm}(y_2), y_2) = 0$ and is given explicitly in \eqref{eq:vx_thres}. The factor $y_c = 3^{-1/2}$ comes from solving 
$h_c^-(y_c) =0, y_c >0$. See the right figure in Figure \ref{fig:vx_OT} for $\sgn(\D(y))$ in different regions. Denote $Q_{\e} = [1-\e, 1 + \e] \times [-\e, \e]$ and $Q_i$ is the four quadrants with center at $(1, 0)$, e.g. $Q_1 = \{ y_1 \geq 1, y_2 \geq 0\}$. For the P.V. integral, since the kernel $ \f{y_1^2 -y_2^2}{|y|^4}$ has mean $0$ in each quadrant $\R^{\pm} \times \R^{\pm}$, it is not difficult to show that 
\[
\lim_{\e \to 0} \int_{  Q_{\e}^c} \D(y) \eta(m- y_1, y_2) \one_{y_1 \geq 0} dy 
= \sum\nolim_{ 1\leq i \leq 4 }  \lim_{\e_i \to 0} \int_{  Q_{\e_i}^c \cap Q_i } \D(y) \eta(m- y_1, y_2) \one_{y_1 \geq 0} dy .
\]
Thus, we can estimate the P.V. integral separately in each $Q_i$.

\subsubsection{Inner region $\Om_{in}$ \eqref{eq:vx_cy_region} }\label{hol:vx_cy_in}

This region contains the singularity \((\pm 1,0)\). We first evaluate the principal-value integral and then 
apply $x$-transportation estimates.

In $\Om_{in}$, we have $|y_1| \leq 1$. Denote $y_c = 3^{-1/2}$. Note that $\D = K_2(y_1 + 1, y_2) - K_2(y_1-1,y_2)$ is singular at $(1, 0)$. Applying Lemma \ref{lem:pv} to $-K_2(y_1-1, y_2)$ yields 
\beq\label{eq:vx_cy_holest2}
\bal
 \int_0^1 \int_{ 0}^{y_c} \D(y) \eta_m(y) dy 
&= \lim_{ \e \to 0} \int_0^1 \int_{\e}^{y_c} \D(y) \eta_m(y) dy  
- \f{\pi}{8} \eta_m(1, 0) 
 \teq \lim_{\e \to 0} I_{\e} - \f{\pi}{8} \eta_m(1, 0) . \\
\eal
\eeq

Let $\TTy(y) \geq 0$ be the map that solves 
\beq\label{eq:map_Cy_uy}
\int_{y_1}^{\TTy(y)} \D(s, y_2) ds = 0, \quad y_1 \in [0,\ h_c^-(y_2)],
\eeq
with formula in \eqref{eq:map_vx_cy}. Applying the sign inequality \eqref{eq:vx_cy_sign} and Lemma \ref{lem:trans} in $y_1$ direction  yields 
\[
\B|\int_{\e}^{y_c} d y_2 \int_0^{ \TTy( 0, y_2) }  \D(y) \eta_m(y) d y_1 \B| 
\leq 
\int_{\e}^{y_c} d y_2 \int_0^{ h_c^-(y_2)}  |\D(y)|  |y_1 - \TTy(y)|^{1/2} d y_1  
[\eta]_{C_x^{1/2}} \teq I_{in, ++}.
\]

See the blue arrows in $\Om_{in}$ in the right figure of Figure \ref{fig:vx_OT} for an illustration of this transportation estimate. Using the symmetry of $\D$ in $y_1, y_2$, we generalize the above estimate of the integral
in the region $ \{ y: \e \leq |y_2| \leq y_c, |y_1| \leq \TTy(0, |y_2|) \}$, which is bounded by $4 I_{in,++}$.

The remaining part of the integral in $\Om_{in}$ is in the following region 
\beq\label{eq:vx_cy_region1}
\bal
&R_{in, \e} \teq \{  |y_1|\leq 1, |y_2| \geq y_c\} \cup \{ \TTy(0, |y_2|) \leq |y_1| \leq 1 , \e \leq |y_2| \leq y_c \}, \\
& R^+_{in, \e} \teq R_{in, \e} \cap ( [0, 1] \times \R ) , \qquad \quad 
R^{++}_{in, \e} \teq R_{in, \e} \cap ( [0, 1] \times \R^+ ) \subset \R^2_{++} .
\eal
\eeq

Since $\D > 0$ in $R^+_{in, \e}$, we use the odd symmetry of $\D$ in $y_1$ and even symmetry in $y_2$ to obtain 
\[
\bal
&\B|\int_{R_{in, \e}} \D(y) \eta_m(y) dy\B|
=\B| \int_{R^+_{in, \e}} \D(y)  (\eta_m(y_1, y_2) - \eta_m( -y_1, y_2)) dy \B| 
\leq 
2 \int_{R^{++}_{in, \e}} |\D(y) | \sqrt{2 y_1} dy [\eta]_{C_x^{1/2}},
\eal 
\]
where we have the factor $2$ since the estimates in $y_2 \geq 0$ and $y_2 \leq 0$ are the same.

Plugging the above estimate in \eqref{eq:vx_cy_holest2} and using the symmetry of $\D$ in $y_1, y_2$, we derive 
\begin{align}\label{eq:vx_cy_in}
&\B|\int_{\Om_{in}} \D(y) \eta_m(y) + \f{\pi}{4} ( \eta_m(1, 0) - \eta_m(-1, 0))\B|  \\
\leq & \B( 4 
\int_{0}^{y_c} d y_2 \int_0^{ h_c^-(y_2)}  |\D(y)|  |y_1 - \TTy(y)|^{1/2} d y_1  
+  2 \int_{R^{++}_{in, 0}} |\D(y) | \sqrt{2 y_1} dy \B) [\eta]_{C_x^{1/2}} \teq C_{in} [\eta]_{C_x^{1/2}}. \notag 
\end{align}

\subsubsection{ Estimate in $\Om_{mid}$ \eqref{eq:vx_cy_region}}\label{hol:vx_cy_mid}

We develop two estimates for the integral \eqref{eq:vx_cy_region} 
\[
 S_{mid} \teq \int_{\Om_{mid}} \D(y) \eta_m(y) dy,
\]
since $u_y, v_x$ benefit from different cancellations. The first estimate transports in the \(y_1\)-direction and is
well suited to the \(u_y\) estimate.  The second estimate transports in the
\(y_2\)-direction and is better adapted to the cancellation of the local term in
the \(v_x\) bound.

\vs{0.05in}
\paragraph{\bf{First estimate for $u_y$ }}
The first estimate is similar to that in Section \ref{hol:vx_cy_in} and 
is useful for the estimate of $u_y$. Notice that the singularities of $\D$ \eqref{eq:K2_del} are $(\pm 1, 0)$. We first rewrite $S_{mid}$ as follows using Lemma \ref{lem:pv} twice with $Q_{1,\pm} = (\pm [1, m]) \times [0, 1]$ and $Q_{2,\pm} = (\pm [1, m]) \times [-1, 0]$
\beq\label{eq:vx_cy_holest3}
S_{mid} = \lim_{\e \to 0} \int_{ \Om_{mid} \cap |y_2| \geq \e  } \D(y) \eta_m(y) dy 
- \f{\pi}{4} (\eta_m(1, 0) - \eta_m(-1, 0)).
\eeq
To estimate the integral, we first study the sign of $\D(y)$. 
Recall $s_c$ from 
\eqref{eq:vx_thres}. For $y_1 \in [1, m]$,  we have 
\[
\D(y) > 0, |y_2| > s_c( m ),
\]

For $|y_2| < s_c( m )$, the sign of $\D(y)$ is given in \eqref{eq:vx_cy_sign}. We introduce the remainder region
\beq\label{eq:vx_cy_region2}
\bal
& R_{mid} \teq \{ |y_1| \in [1, m], |y_2| \geq s_c(m) \} \cup \{ |y_2| < s_c(m), 1\leq |y_1| \leq \TTy(m, |y_2|) \},  \\
& R_{mid}^+ \teq  R_{mid}  \cap \{ y_1 \geq 0 \}, \quad R_{mid}^{++} \teq  R_{mid}^+ \cap 
\{ y_2 \geq 0 \} ,
\quad R_{mid}^{+-} \teq  R_{mid}^+ \cap 
\{ y_2  < 0 \} . 
\eal
\eeq

In $\Om_{mid} \bsh R_{mid} $, using $\sgn(\D)$ \eqref{eq:vx_cy_sign} and applying Lemma \ref{lem:trans} in the $y_1$ direction, we yield 
\[
\B|\int_{\e}^{ s_c(m)} d y_2 \int_{ \TTy(m, y_2)}^m \D(y) \eta_m(y) dy_1 \B|
\leq \int_{\e}^{s_c(m)} d y_2 \int_{ h_c^+(y_2)}^m |\D(y)| |y_1 - \TTy(y)|^{1/2} dy_1 \teq I_{mid, ++}
\]
where $\TTy$ is given in \eqref{eq:map_vx_cy}. See the blue arrows in $\Om_{mid}$ in Figure \ref{fig:vx_OT} for an illustration of this transportation estimate. We generalize the above estimate of integral in the region $ \{ y: \e \leq |y_2| \leq s_c(m), \TTy(m, |y_2|) \leq |y_1| \leq m  \}$ using symmetry of $\D$, which is bounded by $4 I_{mid, + + }$. 

In the remainder region  $R_{mid}$,   $\D(y)$ is positive if $y_1 > 0$.  Using
\[
 |\eta_m(y_1, y_2) - \eta_m(-y_1, y_2)| \leq \sqrt{2 y_1} [\eta]_{C_x^{1/2}},
 \]
the odd symmetry of $\D$ in $y_1$, we obtain
\[
|\int_{R_{mid} \cap |y_2| \geq \e} \D(y) \eta_m(y)  d y|
= |\int_{R_{mid}^+ \cap |y_2| \geq \e} \D(y) (\eta_m(y ) - \eta_m(-y_1, y_2)  ) d y|
\leq \int_{ R_{mid}^+ } |\D(y)| \cdot \sqrt{2 y_1} [\eta]_{C_x^{1/2}} d y.
\]
Since the integral in 
$R_{mid}^{++}$ and $R_{mid}^{+-}$
are the same,  combining the above two integral estimates,  we obtain an estimate analogous to \eqref{eq:vx_cy_in}
\beq\label{eq:vx_cy_mid1}
\bal
| S_{mid} + \f{\pi}{4} ( \eta_m(1, 0 ) - \eta_m(-1, 0)) | & \leq 
 \B( 4 
\int_{0}^{ s_c(m)} d y_2 \int_{ h_c^+(y_2)}^m  |\D(y)|  |y_1 - \TTy(y)|^{1/2} d y_1   \\
&  +  2 \int_{R^{++}_{mid}} |\D(y) | \sqrt{2 y_1} dy \B) [\eta]_{C_x^{1/2}} \teq C_{mid, 1}(m) [\eta]_{C_x^{1/2}}.
\eal
\eeq

\paragraph{\bf{Second estimate for $v_x$}}
In the second estimate, instead of using transportation in the $y_1$ direction, we use transportation in the $y_2$ direction. This estimate will be very useful for $v_x$. We also combine the estimate of $S_{mid}, S_{out}$ \eqref{eq:vx_cy_region}. Recall $\Om_{mid} \cup \Om_{out} = \{ |y_1| \geq 1\}$
from \eqref{eq:vx_cy_region}. First, applying 
Lemma \ref{lem:pv} twice to $K_2(y_1+1, y_2)$ and $-K_2(y_1-1, y_2)$, respectively, we yield 
\beq\label{eq:vx_cy_mid20}
S_{mid} + S_{out} =  \f{\pi}{4} (\eta_m(1, 0 ) - \eta_m(-1, 0)) + \lim_{\e \to 0} \int_{ |y_1| \geq 1+\e }  \D(y) \eta_m(y) dy 
\teq  I + \lim_{\e \to 0} II_{\e}. 
\eeq

We remark that the above decomposition and the sign of $ (\eta_m(1, 0) - \eta_m(-1,0))$ are different from those in \eqref{eq:vx_cy_holest3} since we take the limit in different variables.
The above integral is similar to \eqref{eq:vx_holest1} and is simpler since we do not localize the kernel $\D$. We apply the same argument as that in Section \ref{hol:vx_step0}, \ref{hol:vx_step1}. Recall that $\D$ satisfies the sign condition \eqref{eq:vx_del1}. Note that
\beq\label{eq:vx_cy_holest4}
\int_0^{\inf} \D(y_1, y_2) dy_2 =\f12 \B( \f{ y_2}{y_2^2 + (y_1+1)^2} - \f{y_2}{y_2^2 + (y_1-1)^2} 
\B) \B|_0^{\inf} = 0,
\eeq
and $\D$ is even in $y_2$ and odd in $y_1$. 
Although \(\eta_m(y)\) in \eqref{eq:vx_eta} has a jump across \(y_1=m\) due to the boundary, it is \(C^{1/2}\) on each side of $y_1 = m$ : $\eta_m \in C^{1/2}(\Om_{in}), \, C^{1/2}(\Om_{out})$. In particular, for a fixed $y_1$, $\eta_m(y_1, \cdot) \in C_{y_2}^{1/2}$.  For $\d > 0$, applying Lemma \ref{lem:trans} in the $y_2$ direction, we obtain 
\beq\label{eq:vx_cy_mid21}
\bal
 |II_{\d}| \leq & 4 \int_{1 + \d}^{\inf} d y_1\int_{ s_c(y_1)}^{\inf} |\D(y)| \cdot |y_2 - T(y)|^{1/2} d y_2  [\eta]_{C_y^{1/2}} 
\teq C_{mid, 2}( \d) [\eta]_{C_y^{1/2}}, 
\eal
\eeq
where 
$s_c(y_1)$ and $T$ are \emph{the same as} those in Section \ref{hol:vx_step0},  \ref{hol:vx_step1} for $C_x^{1/2}$ estimate without localization, i.e. $B=\infty$, and are given in \eqref{eq:vx_thres}, \eqref{eq:map_vx}.  
The map $T$ \emph{differs} from $\TTy$, which is determined in \eqref{eq:map_Cy_uy}, for $C_y^{1/2}$ estimate. We have a factor $4$ since the same estimate applies to integral in each quadrant, $\pm  [1+ \d, \inf] \times \R_{\pm} $.

Recall the discussion in Section \ref{sec:holx_singu_md}. For $m>1$ close to $1$, to ease the computation, we use a simpler estimate. Using $y_2$-transport map $T(y) = \f{y_1^2}{y_2}$ \eqref{eq:map_K2_cy} for $K_2(y)$ 
and Lemma \ref{lem:trans},  we have 
\beq\label{eq:hol_strip_spec}
\bal
  \B| \lim_{\e \to 0}  \int_{ [a, b]\bsh [-\e, \e]}  d y_1  \int_0^{\inf} K_2(y) f(y) d y_2 \B|
   \leq \int_a^b \int_{|y_1| }^{\infty} |K_2(y)| \, | \f{y_1^2}{y_2} - y_2 |^{1/2} dy [f]_{C_y^{1/2}} 
& \teq I_{K_2, \inf}(a, b)  [f]_{C_y^{1/2}} ,
\eal
\eeq
for any $ a<b$. Clearly,  $I_{K_2,\inf}(a, b) = I_{K_2,\inf}(-b, -a)$ for  $0<a< b$.  Recall $\D(y) =K_2(y_1 + 1 , y_2) - K_2(y_1 -1, y_2)$ \eqref{eq:vx_del1}. Using a change of variable $(y_1 \pm 1, y_2) \to  y$ and \eqref{eq:hol_strip_spec}, we get  
\[
 |\lim_{ \e \to 0} \int_{ [a, b]\bsh [-\e, \e]} \int_{\R} 
 K_2(y_1 \pm 1 , y_2)    \eta_m(y) dy| 
\leq 2  I_{K_2, \inf}(a \pm 1, b \pm 1) [\eta]_{C_y^{1/2}}. 
\]
For $\e \to 0$, choosing $[a, b] = [\e+1, \d + 1], [-1-\d, -1-\e]$, we estimate the remaining part of the integral in $y_1 \in [a, b]$ as follows 
\beq\label{eq:vx_cy_mid2}
\bal
 & |S_{mid } + S_{out}  - \f{\pi}{4} (\eta_m(1, 0) - \eta_m(-1, 0))   | \leq \lim_{\e \to 0} |II_{\e}| 
 \leq C_{out}(1, \d) [\eta]_{C_y^{1/2} }, \\
 & C_{out}(1, \d) \teq 4 (I_{K_2, \inf}( 0, \d )  + I_{K_2, \inf}( 2, 2 + \d ) ) 
 + C_{mid, 2}(\d) . 
 \eal
 \eeq

When we combine different estimates to estimate $[v_x]_{C_y^{1/2}}$, the factor $- \f{\pi}{4} (\eta_m(1) - \eta_m(-1))$ 
allows us to partially cancel the local term $-\f{\pi}{2} \om$ in $v_x$ \eqref{eq:ux_local}. The factor $ \f{\pi}{4} (\eta_m(1) - \eta_m(-1))$ in \eqref{eq:vx_cy_mid1} has the opposite sign and does not offer such cancellation 
in the estimate of $v_x$.

\subsubsection{Estimate in the outer region}

Recall the integral \eqref{eq:vx_cy_region} in the outer region
\[
 S_{out} \teq \int_{\Om_{out}} \D(y) \eta_m(y) dy.
\]
In $\Om_{out}$, we have $\eta_m \in C_y^{1/2} ( [m ,\infty) \times \R  )$ and $\eta_m \in C_y^{1/2} ( (-\infty, - m] \times \R  )$. For $\d >0$ and $\d_m = \max(m, \d+1)$, following the second estimate of $S_{mid}$ in Section \ref{hol:vx_cy_mid}, applying estimates \eqref{eq:vx_cy_mid21} to $|y_1| \geq \d_m$, and \eqref{eq:hol_strip_spec} to the region $|y_1| \in [m \pm 1 ,  \d_m \pm 1 ]$, we yield 
\beq\label{eq:vx_cy_out1}
\bal
 |S_{out}| & \leq \big(  \, C_{mid, 2}(\d_m-1) + 4( I_{K_2,\inf}( m+1,  \d_{ m} + 1 ) + I_{K_2, \inf}(m-1, \d_{ m}-1) ) \,  \big) [ \eta ]_{C_y^{1/2}} \\
& \teq C_{out}(m, \d) [ \eta ]_{C_y^{1/2}} ,  
 \quad   \d_{ m} = \max(m, \d + 1).
\eal
\eeq
We have a factor $4$ since the same estimate applies to the region $\pm [m ,\infty] \times \R_{\pm}$. The above notation $C_{out}(m, \d)$ is consistent with $C_{out}(1, \d)$ in \eqref{eq:vx_cy_mid2}, where $\d_m = 1 + \d$.

\subsubsection{Modification near the singularity}\label{sec:holy_singu_md}

As in Section~\ref{sec:holx_singu_md}, this modification is introduced only to simplify the
computer-assisted evaluation near the singularity; it does not change the main
transportation structure of the estimate. Similar to Section \ref{sec:holx_singu_md},  
near the singularity $(1, 0)$, for $|y_2| \leq \d$, we modify the estimate of \eqref{eq:vx_cy_holest2}, \eqref{eq:vx_cy_mid1} by separating two kernels \eqref{eq:vx_del1}
\[
S_{in, \d}  \teq \lim_{\e \to 0}  \int_{ \e \leq |y_2| \leq \d ,|y_1|\leq 1} \D(y) \eta_m(y) dy , \quad 
S_{mid, 1, \d}  \teq  \lim_{\e \to 0} \int_{ \e \leq |y_2| \leq \d, 1 \leq |y_1| \leq m } \D(y) \eta_m(y) dy.
\]

In supplementary material II \cite[Section 5.5]{ChenHou2023bSupp} (attached to 
Part II \cite{ChenHou2023b}), we establish 
{\small 
\[
\bal
 |S_{in,  \d}|  & \leq 
  \B\{ \int_0^1 \f{  \sqrt{2}   y_1^{1/2} \d }{ (y_1+1)^2 + \d^2}  d y_1 
  + 2 \B( \sqrt{2} \f{1}{2} \arctan( \d )
+  2 \d^{1/2}(  f_s( \sqrt{1/\d}) - f_s(1) + C_{K_2, up}  ) \B)  \B\} [\eta]_{C_x^{1/2} }  \\
&  \teq C_{in, \d} [\eta]_{C_x^{1/2} } , \\
 | S_{mid, 1, \d } | & \leq  \B\{  \int_{|y_2| \leq \d} \int_1^m K_2(y_1 + 1, y_2) |2 y_1|^{1/2} d y_1  \\
 & \quad +  2 \B( \sqrt{2m } \f{1}{2} \arctan( \f{\d}{m-1} ) +  2 \d^{1/2}(  f_s( \sqrt{(m-1)/\d}) - f_s(1) + C_{K_2, up}  ) \B)   \B\}  [\eta]_{C_x^{1/2} }  \\
 & \teq C_{mid, 1, \d}(m) [\eta]_{C_x^{1/2} } ,
\eal
\]
}
where $C_{K_2, up}$ is defined in \eqref{eq:hol_singu_nota1}, and $f_s$ is given by 
\[
f_s(t) \teq  \int_0^t \f{ s^2}{1 + s^4} ds . 
\]

We obtain modified estimates of \eqref{eq:vx_cy_in}, \eqref{eq:vx_cy_mid1} with bounds $C_{\al}$ replaced by $\td C_{\al}, \al = in, mid$ 
\beq\label{eq:vx_cy_md}
\bal
 \td C_{in}(\d) & \teq 4 
\int_{\d}^{y_c} d y_2 \int_0^{ h_c^-(y_2)}  |\D(y)|  |y_1 - \TTy(y)|^{1/2} d y_1  
+  2 \int_{R^{++}_{in, \d}} |\D(y) | \sqrt{2 y_1} dy 
+ C_{in, \d}  , \\
 \td C_{mid, 1}(m, \d) & \teq 
  4 
\int_{0}^{ s_c(m)} \one_{y_2\geq \d} d y_2 \int_{ h_c^+(y_2)}^m  |\D(y)|  |y_1 - \TTy(y)|^{1/2} d y_1     +  2 \int_{R^{++}_{mid}, y_2  \geq \d} |\D(y) | \sqrt{2 y_1} d y  + C_{mid, 1, \d}(m) . 
\eal
\eeq

For $m> 1$ very close to $1$, we have an additional estimate for $S_{mid}$ \eqref{eq:vx_cy_mid1}
\begin{align}\label{eq:vx_cy_mid3}
& | S_{mid} + \f{\pi}{4} ( \eta_m(1, 0 ) - \eta_m(-1, 0)) |  
= |\lim_{\e \to 0} \int_{|y_2| \geq \e, |y_1| \in [1, m]} \D(y) \eta_m(y) dy |
\leq  C_{mid, 3}(m) [\eta]_{C_x^{ \f{1}{2}} }  ,  \notag  \\
&  C_{mid, 3}(m) \teq  \f{\pi}{4}\sqrt{2 m} + \f{1}{4}(m-1)^2 \sqrt{2 m}  + 
2  \B( \sqrt{2m} \f{\pi}{8} + 2 (m-1)^{1/2} C_{K_2, up} \B) . 
\end{align}

The above modifications and the refinements in \eqref{eq:vx_cy_mid2}, \eqref{eq:vx_cy_out1} are very tiny and of order $\d^{1/2}, |m-1|^{1/2}$. We choose a small $\d$ with $\d \leq y_c=3^{-1/2}$ \eqref{eq:vx_cy_sign}, e.g. $\d = 2^{-14}$, and use them to ease the computation of the integral near $(1, 0)$. If $\d = 0$, we recover the previous estimates.

\subsubsection{Summarize the estimates}

We now combine the estimates in the three regions and translate the bounds
from the rotated function \(\eta\) back to the original vorticity \(\omega\). Recall from \eqref{eq:vx_cy_holest1}
\[
\bal
v_x(z) - v_x(x)  &= - \f{1}{\pi} \int \D(y) \eta_m dy - \f{1}{2} (\eta_m(-1, 0) - \eta_m( 1, 0)) \\
&= - \f{1}{\pi} \B( S_{in} + S_{mid} + S_{out}  - \f{\pi}{2} (\eta_m(1, 0) - \eta_m( -1, 0))
 \B) , \\
 u_y(z) - u_y(x) & =- \f{1}{\pi} \B( S_{in} + S_{mid} + S_{out}  + \f{\pi}{2} (\eta_m(1, 0) - \eta_m( -1, 0))  \B).
 \eal
\]

Note that $| (\eta_m(1, 0) - \eta_m( -1, 0))| \leq \sqrt{2} [\eta]_{C_x^{1/2}}$. Using relation \eqref{eq:vx_rota}: $[\eta]_{C_{x}^{1/2}} = [\om]_{C_{y}^{1/2}}, 
[\eta]_{C_{y}^{1/2}} = [\om]_{C_{x}^{1/2}}$, and \eqref{eq:vx_cy_in}, the first estimate in $\Om_{mid}$ \eqref{eq:vx_cy_mid1}, \eqref{eq:vx_cy_out1},\eqref{eq:vx_cy_md},\eqref{eq:vx_cy_mid3}, for $u_y$, we prove
\[
\bal
 \f{| u_y(z) - u_y(x) | }{\sqrt{2}} & \leq \f{1}{\pi \sqrt{2}}
 \B( ( \td C_{in}(\d) + \min( \td C_{mid, 1}(m , \d), C_{mid, 3}(m) ) ) [\eta]_{C_x^{1/2}}
+ C_{ out}(m, \d)  [\eta]_{C_y^{1/2}}  \B) \\
& \leq 
\f{1}{\pi \sqrt{2}}
 \B( ( \td C_{in}(\d) + \min( \td C_{mid, 1}(m, \d), C_{mid, 3}(m) ) ) [\om]_{C_y^{1/2}}
+ C_{ out}(m, \d)  [\om]_{C_x^{1/2}}  \B) ,\\
\eal
\]
for $ 0< \d \leq y_c$ \eqref{eq:vx_cy_sign}. Applying \eqref{eq:vx_cy_in}, the second estimate in $\Om_{mid}$ \eqref{eq:vx_cy_mid2} to $v_x$, we prove 
\[
\bal
\f{ |v_x(z) - v_x(x) |}{\sqrt{2}} & \leq 
\f{1}{\pi \sqrt{2}} \B( \td C_{ in}(\d) [\eta]_{C_x^{1/2}} 
+  C_{out}(1, \d)  [\eta]_{C_y^{1/2}} + \f{\pi}{2} \sqrt{2} [\eta]_{C_x^{1/2}}   \B) \\
 &\leq  \f{1}{\pi \sqrt{2}} \B( \td C_{ in}(\d) [\om]_{C_y^{1/2}} 
+  C_{out}(1, \d)   [\om]_{C_x^{1/2}} + \f{\pi}{2} \sqrt{2} [\om]_{C_y^{1/2}}   \B) ,
\eal
\]
for $0< \d \leq y_c$ \eqref{eq:vx_cy_sign}, where 
$ \f{1}{\sqrt{2}}$ is from $|x-z|^{-1/2} =  \f{1}{\sqrt{2}}$, 
and $C_{out}(m, \d)$ is defined in \eqref{eq:vx_cy_out1}.

\subsection{Functions and transportation maps}\label{app:OT_map}

We record the formulas of the transportation maps and the functions for the sign of the kernels in the sharp H\"older estimate.  These formulas are used in the computer-assisted evaluation of integrals 
in the upper bounds. Recall 
{\small 
\[
K_1 = \f{y_1 y_2}{|y|^4}, \quad K_2 = \f{1}{2}\f{y_1^2 - y_2^2}{ |y|^4}.
\]
}

\subsubsection{Sign functions}

Solving $K_1(y_1 + 1/2, y_2 ) - K_1(y_1 - 1/2, y_2) = 0$ for $y_1 \geq 0$, we yield 
{\small 
\[
y_1 =  \B( \f{1/2 - 2 y_2^2 + \sqrt{ 16 y_2^4 + 4 y_2^2 + 1} }{6} \B)^{1/2}.
\]
}
See also \eqref{eq:ux_thres}. Solving $K_2(y_1 + 1, y_2 ) - K_2(y_1 -1, y_2) = 0$ for $y_2 \geq 0$, we yield 
{\small 
\beq\label{eq:vx_thres}
\bal
y_1 & = h_c^{\pm}(y_2) \teq \big( y_2^2 + 1 \pm  2 y_2 \sqrt{ y_2^2 + 1}  \big)^{1/2} , \\ 
  y_2 & = s_c(y_1) \teq  \big( \f{ -(y_1^2 + 1) + 2 (y_1^4 - y_1^2 + 1)^{1/2} }{3}   \big)^{1/2} .
\eal
\eeq
}

\subsubsection{Transportation maps}

Below, we record the equation for  transportation maps $T$ and $\TTy$.

\paragraph{\bf{Map for $u_x$}}
We derived the $T$-equation  in \eqref{eq:cubic_T0}, \eqref{eq:holx_ux_ker}. For a fixed $s_2 \neq 0$ and $s_1 > 0$, solving 
{\small 
\[
\int_{T(s)}^{s_1} ( K_1( t + \f12 , s_2) - K_1( t - \f12 , s_2) ) d t = 0,
\]
}
yields the equation of the transportation map in $x$ direction
\beq\label{eq:map_ux}
 T^3 + T^2 s_1 + T ( s_1^2 - \f{1}{2} + 2 s_2^2 ) - \f{ (4 s_2^2 + 1)^2}{ 16 s_1} = 0.
 \eeq
We rewrite the above equation as an equation for $Z = T + \f{s_1}{3}$
{\small 
\[
0 = 
 Z^3 + Z ( \f{2}{3} s_1^2 - \f{1}{2} + 2 s_2^2) 
- \B( \f{ (4 s_2^2 + 1)^2}{16 s_1}  + \f{7}{27} s_1^3 + \f{s_1}{3} ( 2 s_2^2 - \f{1}{2})  \B) 
\teq Z^3 + p(s_1, s_2) Z + q(s_1, s_2) .
\]
}
The discriminant is given by 
\beq\label{eq:cubic_discr}
\D_Z(s_1, s_2) = -( 27 q(s_1, s_2)^2 + 4 p(s_1, s_2)^3).  
\eeq

Note that
{\small 
\[
-q \geq \f{7}{27} s_1^3 - \f{s_1}{6} + \f{1}{16 s_1} \geq ( 2 \sqrt{ \f{7}{27} \cdot \f{1}{16}} - \f{1}{6}) s_1 \geq 0 ,
\]
}
and $-q, p $ are increasing in $|s_2|$. We yield 
{\small 
\[
-\D_Z(s_1, s_2) \geq -\D_Z(s_1, 0) =  \f{ (1-4 s_1^2)^2( 27 - 56 s_1^2 + 48 s_1^4)}{ 256 s_1^2}
 \geq 0.
\]
}

When $s \neq ( \f{1}{2}, 0)$, the above inequality is strict, and we have a unique real root. 
Using the solution formula for a cubic equation, we obtain the formula for the real root 
{\small
\beq\label{eq:cubic_T3}
Z = r_1  - \f{p}{3r_1} , \quad  
r_1 = \B( \f{1}{2} \big( -q + \sqrt{ q^2 + \f{4}{27} p^3 } \, \big)    \B)^{1/3}, 
\quad T = Z - \f{s_1}{3}.
\eeq
}

\paragraph{\bf{Map for $[ u_y]_{ C_x^{1/2} }, [ v_x]_{C_x^{1/2}} $}}

 We derived $T$-equation in \eqref{eq:map_Cx_uy}, \eqref{eq:vx_del1}. 
 Fix $y_1 \geq 0, y_1 \neq \pm 1 $, $y_2 \neq 0$. 
 The measure-zero set $y_1 =0, \pm 1, y_2 =0$ may simply be omitted in 2D integrals. Solving 
 {\small 
\[
\int_{T(y)}^{y_2}(  K_2(y_1 +1,  s_2 ) - K_2(y_1-1, s_2 ) ) d s_2 = 0,
\]
}
yields the equation of the transportation map in $y$ direction
\beq\label{eq:map_vx}
T^3 + T^2 y_2 +  T ( y_2^2 + 2 + 2 y_1^2) - \f{ (y_1^2 - 1)^2}{y_2} = 0.
\eeq
We rewrite the above equation as an equation for $W = T + \f{y_2}{3}$ 
{\small
\[
0 = W^3 + W( \f{2 y_2^2}{3} + 2 + 2 y_1^2)
 - \B( \f{ (y_1^2 - 1)^2}{ y_2} + \f{y_2^3}{27} + \f{y_2}{3} ( \f{2y_2^2}{3} +2 + 2 y_1^2 )  \B) 
 \teq W^3 + p_2(y) W + q_2(y) .
\]
}
 Since $p_2  > 0$, using the discriminant \eqref{eq:cubic_discr} yields $-\D_W(y) > 0$. Thus the cubic equation of $T$ or $W$ has a unique real root with formula similar to \eqref{eq:cubic_T3}.

\vspace{0.1in}
\paragraph{\bf{Map for $[ v_x ]_{ C_y^{1/2} },  [ u_y]_{ C_y^{1/2} } $}}

 We derived $\TTy$-equation in \eqref{eq:map_Cy_uy}, \eqref{eq:K2_del}. 
 Fix $y_1, y_2 \geq 0$. Solving 
 {\small 
\[
\int_{y_1}^{\TTy(y)} \D(s, y_2) ds = 0  ,
\]
}
with $\TTy(y)  \geq 0$  yields the equation of the transportation map in $x$ direction.  
In the following two formulas only, we write $T=\TTy(y)$; this $T$ is not the map in \eqref{eq:map_ux} or \eqref{eq:map_vx}. Then
\[
1 - T^2 - y_1^2 + T^2 y_1^2 - 2 y_2^2 - T^2 y_2^2 - y_1^2 y_2^2 - 3 y_2^4 = 0,
\] 
The above equation is  equivalent to
\beq\label{eq:map_vx_cy}
T^2 =  \f{ y_1^2 + 2 y_2^2 + y_1^2 y_2^2 + 3 y_2^4 - 1}{ y_1^2 - y_2^2 - 1}. 
\eeq

We apply the above map to the following two regions separately
\[
y_1 \in [0,1] , \ y_1 \leq h_{c}^-(y_2), \quad y_1 \in [1, \infty] ,\  y_1 \geq h_{c}^+(y_2).
\]

\section{Weights, parameters, and auxiliary stream function estimates}

In this appendix, we collect parameters and cutoff functions, and prove auxiliary estimates.
Appendix ~\ref{app:para_wg} records the weights,
cutoff functions, and coupling parameters used in Section~\ref{sec:EE}.  
Appendix~\ref{app:appr_vel_para}
records the parameters in the finite-rank approximation of the velocity.
In Appendix~\ref{app:euler_stream}, we prove several stream-function estimates in Section~\ref{sec:euler} for 3D Euler. In Appendix~\ref{app:idea_num_est}, we explain representative numerical-analysis ideas in   \cite{ChenHou2023b} to obtain rigorous weighted bounds.

\subsection{Parameters for the weights}\label{app:para_wg}

In the energy estimates, we use several weights. We group them by their roles:  \(\psi_i, g_i\)
 for weighted \(C^{1/2}\) estimates, $\vp_i$ and $\vp_{g, i}$ for 
decaying and growing weighted \(L^\infty\) estimates, respectively,
and $\rho_{ij}$ for weighted $L^{\infty}$ estimates of nonlocal terms.

Let  $x = (x_1, x_2 )\in \R^2$. For the H\"older estimates, we use the following weights 
\footnote{
In practice, we use the double-precision floating point values of these parameters which can differ from the values below by the machine precision, e.g. $10^{-16}$. 
}
\bseq\label{wg:hol}
\beq 
\bal
&\psi_1   = |x|^{-2} + 0.5 | x|^{-1} + 0.2 |x|^{-1/6},  \\
&\psi_2  = p_{2, 1} | x|^{-5/2} + p_{2,2} |x|^{-1} + p_{2,3} |x|^{-1/2} + p_{2,4} |x|^{1/6}, \\ 
\quad 
& \psi_3 \equiv \psi_2, \qquad  \qquad \vec{p}_{2, \cdot} = (0.46, 0.245, 0.3, 0.112), \\
\eal 
\eeq 
and 
\begin{align}
   & g_{i}(h)   = g_{i0}(h) g_{i0}(1, 0)^{-1}, \quad 
 g_{i0}(h) =  ( \sqrt{ |h_1| + q_{i1} |h_2| } + q_{i3} \sqrt{ |h_2| + q_{i2} |h_1| }   )^{-1}, 
 \  g_3(h) \equiv g_2(h),  \notag \\
 & \vec q_{1, \cdot }  = ( 0.12, 0.01, 0.25 ) , \quad  \vec q_{2,  \cdot } = ( 0.14, 0.005, 0.27 ), \quad 
\vec q_{3, \cdot } = \vec q_{2, \cdot }.
\end{align}
\eseq 
To estimate the weighted $L^{\inf}(\vp_{\mathrm{elli}})$ norm of the Poisson error, the weighted $L^{\inf}$ norm of $ \uu_A, \na \uu_A$, and the H\"older estimate of $\psi_{u}\uu_A, \psi_1 (\na \uu)_A$, 
for $x = (x_1, x_2) \in \R^2$, we use 
\beq\label{wg:elli}
\bal
 \vp_{elli} & = |x_1|^{-1/2} ( |x|^{-2} + 0.6 |x|^{-1/2}) + 0.3 |x|^{-1/6},  \quad
   \psi_{u} = |x|^{-5/2} + 0.2 |x|^{-7/6}, \\
  \rho_{10} & = \rho_{01} = |x|^{-3} + |x|^{-7/6}, \quad \rho_{ij} = \psi_1, \ i + j =2 , \\ 
\rho_3 & = |x|^{-1} + |x|^{-1/6}, \quad   \rho_4 =  |x_1|^{-1/2} ( |x|^{-2.5} + 0.6 |x|^{-1/2}) + 0.3 |x|^{-1/6},
\eal
\eeq
where $\psi_1$ is defined in \eqref{wg:hol}. The weight $\vp_{elli}$ is similar to $\vp_1$ except that we choose a less singular power for the first term. We use $\rho_3$ to capture the vanishing order of $ \pa_i ( \na \UU)_{\FM} \sim |x|$ near $x=0$ \eqref{eq:U} and estimate $ \rho_3 \pa_i (\na \UU)_{\FM}$. The weight $\rho_4$, which is singular along $x_1 = 0$, is used for another estimate of $ \rho_4 u_A$ using energy $|| \om \vp_1||_{\infty} $. See Appendix {\secudx} in Part II \cite{ChenHou2023b}.

For the weighted $L^{\inf}$ estimates with decaying weights, we use the following weights 
\footnote{
We write $p_{6, i}$ in the form $a_i p_{5, i}$ since we first determine $p_{5, i }$ for the weight $\vp_2$ of $\eta$ and then determine the parameter $p_{6, i}$ for weights $\vp_3$ relative to $\vp_2$. 
}
\beq\label{wg:linf_decay}
\bal
\vp_1 & = |x_1|^{-1/2} ( |x|^{-2.4} + 0.6 |x|^{-1/2} ) + 0.3 |x|^{-1/6},  \\
\vp_2 & = |x_1|^{-1/2} ( p_{5, 1} |x|^{-5/2} + p_{5, 2} |x|^{-3/2} + p_{5,3} |x|^{-1/6} ) + p_{5,4} |x|^{-1/4} +  p_{5, 5} |x|^{1/7}, \\
\vp_3 & = |x_1|^{-1/2} ( p_{6, 1} |x|^{-5/2} + p_{6, 2} |x|^{-3/2} + p_{6,3} |x|^{-1/6} ) + p_{6,4} |x|^{-1/4} +  p_{6, 5} |x|^{1/7}, \\
p_{5, \cdot} &= (0.42,  \ 0.144, \  0.198, \  0.172, \ 0.0383) \cdot \mu_0, \quad \mu_0 = 0.917, \\ 
p_{6, \cdot} & = ( 2.5 \cdot p_{5,1} , \ 2 \cdot  p_{5, 2},  \ 3.8 \cdot p_{5,3} ,  \ 1.71 \cdot p_{5,4},  \ 2.39 \cdot p_{5,5}  ). \\
\eal
\eeq

For the weighted $L^{\inf}$ estimates with growing weights and  $x = (x_1, x_2) \in \R^2$,  we use 
\footnote{
We choose $p_{9, 2}$ so that $ ( p_{9, 1 } , p_{9, 2 })$ is proportional to $( p_{8, 1},p_{8, 2})$ .
}
\beq\label{wg:linf_grow}
\bal
\vp_{g, 1} &= \vp_1 +  p_{7, 1} |x|^{1/16}, \   \vp_{g, 2} = \vp_2 + p_{8, 1} |x|^{1/4} + p_{8, 2} |x|^{\al_{g,n} } , \  \vp_{g,3} = \vp_3 +  p_{9, 1} |x|^{1/4}  + p_{9, 2}  |x|^{ \al_{g,n}} ,  \\
 p_{7, 1} &= 1,  \  p_{8, 1}  = 0.07, \  p_{8, 2} = 0.002, \quad p_{9,1} = 0.154, \quad p_{9,2} = p_{8, 2} p_{9,1} / p_{8, 1}, \quad \al_{g, n} = 1/3 + 10^{-8} ,
\eal
\eeq

\paragraph{\bf{Parameters in energies}}
We use parameters below for coupling weights in the energies \eqref{energy2},\eqref{energy3},\eqref{energy4}
\beq\label{wg:EE}
\bal
&\tau_1 = 5, \quad \mu_1 = 0.668, \quad \mu_2   = 2 \mu_1 = 1.336, 
\quad \mu_4 = 0.065, \quad \tau_2 = 0.23 , \\
& \mu_5 = 76 ,  \
\mu_{51} = 61, \  \mu_{52} = 15 , \quad  \mu_{6} = 61,  \  \mu_{62} = 35.88, \quad \mu_7 =  9.5, \quad \mu_8 = 4.5 , \\
& \mu_{h} = \tau_1^{-1}(1,  \mu_1, \mu_2 ) , \quad  \mu_g  = (\tau_2 \mu_4, \tau_2, \tau_2) .
\eal
\eeq

\subsubsection{Parameters for the cutoff functions}\label{app:cutoff}

We record  the cutoff functions 
for near-origin analytic
corrections, far-field profile construction, and residual/error
decompositions. The piecewise bounds of these cutoff functions are estimated in Part II \cite[Appendix {\seccutoffnear}]{ChenHou2023b}. First, we introduce
\beq\label{eq:cutoff_exp} 
\chi_e(x) = \big( 1 + \exp( \f{1}{x} + \f{1}{x-1} ) \big)^{-1}, \quad 
\kp(x; \nu_1, \nu_2 ) = \kp_1(\f{x}{ \nu_1} ) (1 - \chi_e( \f{x}{ \nu_2} )  ) ,  \quad \kp_1(x) = \f{1}{1 + x^4}, 
\eeq
 where $e$ denotes \textit{exponential}. 
We extend $\chi_e(x) = 0$ for $x \leq 0 $ and $\chi_e(x) = 1$ for $x \geq 1$. 
We mainly use
\beq\label{eq:cutoff_near0_chi3} 
 \kp_*(x) = \kp(x ; \, \tf13, \, \tf32). 
\eeq
 We construct 1D radial cutoff function for the semi-analytic parts in \eqref{eq:ASS_decomp1}
\, :
\beq\label{eq:cut_radial}
\bal
 \chi(x)  & = \chi_1(x) (1 - \chi_2(x) ) + \chi_2(x) ,  \quad a_1  = 10, \quad l_1= 50000, \quad a_2 = 10^5, \\
   \chi_1( x) & = \chi_{rati}( \f{x -a_1}{l_1^{1/2}} ), \quad 
\chi_{rati}(x) = \f{ x^7}{  (1 + x^2)^{7/2}} \one_{x \geq 0},  \quad 
 \chi_2(x) = \chi_{e}(  \f{x-a_2}{9a_2} ) , \quad x \in \R .
\eal
\eeq

Below, we consider $(x, y) \in \R^2$. We choose the cutoff functions $\chi_{\bar \e}, \chi_{\hat \e}$ in \eqref{eq:uerr_dec1}
\bseq\label{eq:cutoff_near0_all}
\beq
\bal
 \chi_{\bar \e}(x, y)  &= \kp(x; \nu_{\bar \e, 1}, \nu_{\bar \e, 2}) 
\cdot \kp(y; \nu_{\bar \e, 1}, \nu_{\bar \e, 2}) ,  \quad \nu_{\bar \e, 1 } = 1/192, \quad \nu_{\bar \e,  2} = 3/2 , \\ 
 \chi_{ \hat \e}(x, y) & = \kp_*(x) \kp_*(y) , 
\eal
\eeq
$f_{\chi, i}, \chi_{\cNF}$ in \eqref{eq:appr_near0} 
\beq
\bga
  \chi_{\cNF}(x, y) = \kp(x; 2,10)  \kp(y ; 2, 10 ),  \\
    f_{\chi, 1}  = \D( \f{x y^3}{6} \chi_{\cNF}(x, y)  ) ,
\quad f_{\chi, 2} = x y \cdot  \chi_{\cNF}(x, y), \quad f_{\chi, 3} = \f{x^2}{2} \cdot \chi_{\cNF}(x, y), \\
\ega
\eeq
$ \chi_{j, 2}$ with $j \leq 3$  in \eqref{eq:W2_err_est1} 
\beq
 \chi_{12}   = - \D \phi_{2}, \quad \phi_2 = -\f{x y^3}{ 6} \kp_*(x) \kp_*(y) , 
 \quad    \chi_{22} = xy \kp_*(x ) \kp_*(y), \quad \chi_{32} = \f{x^2}{2}  \kp_*(x) \kp_*(y),  
\eeq
and $\chi_{\mf{ode}}$ in \eqref{eq:cw_chi}
\beq
\chi_{\mf{ode}}(x, y) =  1 - \chi_e(  (x- \nu_{31} ) / \nu_{31} ) \cdot \chi_e( (y - \nu_{32}) / \nu_{32} ), \quad \nu_{31} = 80, \  \nu_{32} =  1200.
\eeq
\eseq

For the cutoff function in \eqref{eq:psi_near0}, we choose 
\beq\label{eq:cutoff_psi_near0}
 \chi_{\phi }(x) = \kp_2( \f{x} { \nu_{4,1}} ) (1 - \chi_e( \f{x}{ \nu_{4, 2} } )  ), \quad 
\kp_2(x) = \f{1}{1 + x^2}, \quad \nu_{4, 1} = 2,  \quad \nu_{4, 2} = 128.
\eeq

\subsection{Parameters for approximating the velocity}\label{app:appr_vel_para}

We collect the parameters used in the finite-rank approximations in Section~\ref{sec:appr_vel}. We choose the following interpolation points $x_i$ and width $t_i$ of the singular region  in the first approximation of velocity $u, u_x$ in \eqref{eq:vel_cutoff_bd}, \eqref{eq:u_appr_1st},  \eqref{eq:u_appr} 
\[
\bal
 u_x : \;\; & x = ( 1, 2,3, 4, 6, 8, 11, 16, 22, 32, 48 )\cdot 64 h_x,  \quad t = (16, 16, 20, 24, 32, 40, 56, 72, 96, 128, 256) h, \\
 u :\;\;&  x = ( 1, 2 , 4, 8, 12, 16, 22, 32, 64  ) \cdot 64 h_x, 
 \qquad \quad t = (8, 8, 24, 40, 56, 72, 96, 128, 256) h,  \\
    h_x & = 13 \cdot 2^{-12}, \quad h = 13 \cdot 2^{-11} . 
\eal
\]

For the parameters $y_0, x_0$ in $\cK_{00}\chi_0$ \eqref{eq:u_appr_1st}, we choose 
\[
y_0 = 256 h_x, \quad u_x: x_0 = 32 h_x,  \quad u : x_0 = 16 h_x , \quad  u_y, v_x, v: x_0 = 128 h_x,
\]
By the divergence-free condition, $v_y=-u_x$, so no additional approximation
for $v_y$ is needed.
Note that we choose the same $y_0$ for all cases in the cutoff function $\td \chi( \f{y - y_0}{y_0} )$ \eqref{eq:vel_cutoff_bd}.
For $u, v, u_x$, we choose the following parameters $R_i$ in the second far-field approximation \eqref{eq:u_appr_2nd}
\beq\label{para:appr_2nd}
R = (8, 16, 32, 64, 128, 256, 512, 1024) \cdot 64 h_x.
\eeq

\subsubsection{Initial conditions for the linearized equations}\label{app:init_data}

This subsection fixes the indexing convention for the rank-one modes used in
the Duhamel representation of \(\widehat{W}_{2}\) in Section \ref{sec:decoup_modi}. The notation below is purely
bookkeeping, but it is important for matching the finite-rank operators in
Section~\ref{sec:finite_pertb} with the approximate space-time solutions used in Section~\ref{sec:W2}. 

Recall the formula \eqref{eq:bous_W2_appr} and the approximation terms $\cK_{1i}$ \eqref{eq:appr_near0} near the origin, 
$\cK_{2i}$ \eqref{eq:appr_vel} for the finite-rank perturbation of velocity. We use the formulas \eqref{eq:u_appr_1st_sum}, \eqref{eq:u_appr_2nd_sum} to approximate the velocity. Denote by $n_{u_x}, n_u$ the number of the terms in \eqref{eq:u_appr_1st_sum} for $ f = u_x$ and $f = u$, respectively, and $n_R = m$ the number of terms in \eqref{eq:u_appr_2nd_sum}. See Appendix \ref{app:appr_vel_para}. We label the rank-one terms : 
\[
\bal
a_1(W_1) \bar F_1(0) = c_{\om}(W_1 ) \bar f_{c_{\om}}, \quad 
a_{ n_\mw{lin} + i}( W )  \bar F_{n_\mw{lin} + i} = a_{ \mf{nl}, i} ( W ) \bar F_{\chi, i} , \quad 
n_\mw{lin} = n_R + n_u + n_{u_x} + 2,
  \eal
  \]
for $i = 1,2, 3$, where $W_{1, 1} = \om_1$, $n_\mw{lin}$ is the total rank for approximating the linearized operator, and \textit{is not} short for $\textit{nonlinear}$. The term $a_{ \mf{nl}, i} ( W ) \bar F_{\chi, i}$ approximates the nonlinear and error terms \eqref{eq:appr_near02}. We denote by $a_i(W_1) \bar F_i(0) , 2 \leq i \leq n_\mw{lin}$ the rank-one terms generated by: 
\begin{itemize}
\item the approximation \eqref{eq:u_appr_2nd_sum} for $i =2, 3, .., n_R+1$ (total $n_R$ terms).

\item the approximation term $\cK_{00}$ of $v, v_x, u_y$ for $i = n_R+2$ (see the discussion above \eqref{eq:u_appr_1st}).

\item \eqref{eq:u_appr_1st_sum} with $ f =u_x$ for $i = n_R + 3 ,.., n_R + n_{u_x} + 2$ (total $n_{u_x}$ terms).

\item  \eqref{eq:u_appr_1st_sum} with $f = u$ for $i = n_R + n_{u_x} + 3, .., n_R + n_{u_x} + n_u + 2 = n_\mw{lin}$ (total $n_u$ terms).

\end{itemize}

For example, for $i=2$, the functional $a_2(W_1)$ and the initial condition $\bar F_2 |_{t=0}$ are given by 
\[
\bga
 a_2(W_1) = u_x(W_{1, 1})(0) - I_1(W_{1,1}) ,
 \quad 
 \bar F_2 |_{t=0} = ( \appo{ \uu} \cdot \na \bar \om , \  \appow{ \uu_x} \cdot \na \bar \th + \appo{ \uu} \cdot   \na  \bar \th_x, \  \appow { \uu_y} \cdot \na \bar \th + \appo{ \uu} \cdot \na \bar  \th_y)
\ega
\]
with 
\footnote{
We abuse notation $\appo{f}$ in  \eqref{eq:eg:mode_appo_f}, which denotes \emph{one mode} in \eqref{eq:u_appr}, rather than \emph{full
approximation} $\appo{f}$ in \eqref{eq:u_appr}.
}
\beq\label{eq:eg:mode_appo_f}
  \appo{ f} =  C_{f0} (1 -\chi_{tot}) S_1^R, \quad f = u, v, u_x, u_y, v_x ,
\eeq
where $u_x(0) - I_1$ is the integral in \eqref{eq:u_appr_2nd_sum} and $C_{f0}$ is defined in \eqref{eq:u_appr_near0}. 
Here, we use the opposite sign for $\bar F_2|_{t=0}$ from the convention in \eqref{eq:appr_vel}, while $a_2(W_1)\bar F_2$ has the same sign.

Although the nonlinear modes 
 $a_{ \mf{nl}, i} ( W ) \bar F_{\chi, i}, i=2,3$ are initially labeled separately, they can be
combined. From \eqref{eq:eta_xi}, \eqref{eq:bous_err}, and
\[
\eta_{xy}(0) = \xi_{xx}(0) = \th_{xxy}(0) , \quad 
\pa_{xy} \bar \cF_2(0) = \pa_{xxy} \bar F_{  \th}(0) = \pa_{xx} \bar \cF_3(0) ,
\]
we obtain $a_{ \mf{nl},2}(W) = a_{ \mf{nl}, 3}(W)$, 
where $\bar \cF_i$ in  \eqref{eq:bous_err} is the residual error and differs from $\bar F_{\chi, i}, 
f_{\chi, i}$ defined in \eqref{eq:bous_err_op}.
Recall $\bar F_{\chi,i } = f_{\chi, i} \ee_i$ from \eqref{eq:bous_err_op}. Thus, we can rewrite $\cNF(W)$ in \eqref{eq:appr_near0} equivalently as a rank-two operator 
instead of a rank-three operator:
\[
  \cNF(W) = a_{ \mf{nl},1}(W) ( f_{\chi, 1}, 0, 0) + a_{ \mf{nl},2}(W) ( 0, f_{\chi,2}, f_{\chi, 3} )
  = a_{ \mf{nl}, 1}(W) \bar F_{\chi, 1} + a_{ \mf{nl},2}(W) (  \bar F_{\chi, 2} + \bar F_{\chi, 3} ). 
\]
Hence, for the combined mode $a_{\mf{nl},2}(W)(\bar F_{\chi,2}+\bar F_{\chi,3})$, 
we construct one approximate space-time solution with initial data $(0,f_{\chi,2},f_{\chi,3})$, instead of constructing two separate solutions with initial data $(0,f_{\chi,2},0)$ and $(0,0,f_{\chi,3})$. We also combine
the terms $a_{\mf{nl},2}(W), a_{\mf{nl},3}(W)$ in
\eqref{eq:bous_W2_appr} and \eqref{eq:bous_err_op}.

\subsection{Estimate of the stream function in 3D Euler}\label{app:euler_stream}

In this section, we prove several stream-function estimates used in  Section \ref{sec:euler} for 3D Euler.

\subsubsection{Proof of Lemma \ref{lem:far} }\label{app:euler_stream_sym}

Recall $r = 1 -C_l y $ in \eqref{eq:label_rs} and $\td{\om }(r, z) = \om^{\th}(r, z) / r$ in \eqref{eq:omth}. Since the support size satisfies $C_l(\tau) S(\tau) < 1/4$, within the support of $\om^{\th}(r, z)$, we have $r \geq 1/2$. 
Hence, $| \td{\om} | \asymp  | \om^{\th} | $ and $|r \om^{\th}(r,z)| \les | \td \om(r,z)| $. 
We can apply Lemma \ref{lem:biot1} and \eqref{eq:rescal41} to get
\[
\bal
|\phi(x, y) | & \les C_{\om} C_l^{-2} \int   |\td{\om}(r_1, z_1)| \B( 1+ |\log( (r_1- (1 - C_l y) )^2 + (z_1 - C_l x)^2) | \B)  d r_1 d z_1  \\
&=  C_{\om}   \int  |\td{\om}(1 - C_l y_1, C_l x_1) | \B( 1 + | \log \big( C_l^2 ( ( y_1 -y)^2 + (x_1- x)^2 )  \, \big) | \B) d y_1 d x_1  ,
\eal
\]
where we used Lemma \ref{lem:biot1} and $r_1 \leq 1 $ in the first inequality, and used change of variables $r_1 = 1 - C_l y_1, z_1 = C_l x_1$ in the second identity. Since
$ C_{\om}  |\td{\om}(1 - C_l y_1, C_l x_1) |  = |\om( x_1, y_1)|$ and $|(x, y)| > 2S$, within the support of $\om$, we get $|(x, y) - (x_1, y_1)|^2 \asymp |(x, y)|^2$. It follows 
\[
\bal
|\phi(x, y)|  & \les \big(1 +| \log( C_l^2 |(x, y)|^2 ) | \big) \int |\om(z)| dz  \les (1 +| \log( C_l |(x, y)| ) | ) 
|| \om( 1 +|z|)^{\b} ||_{\inf} \int_{ |z| \leq S} (1 + |z|)^{-\b} d z \\
& \les || \om( 1 +|z|)^{\b} ||_{\inf} S^{2-\b} (1 + | \log(C_l |(x, y)|) | \, ).
\eal
\]
This proves Lemma \ref{lem:far}.

Next, we prove estimates in Lemmas \ref{lem:E_ker}, \ref{lem:E_vel}. 
Let $ G(y) = -\log |y|  ,  K_1 =  \f{y_1 y_2}{|y|^4}, K_2 = \f{y_1^2 - y_2^2}{2|y|^4}$.  Denote by $K_{1, s}(x, y), G_s(x, y)$ the symmetrized kernel for $K_1, G$ with odd-odd symmetry, as in \eqref{eq:kernel_du_sym}.  For $ G = -\log |y|  ,  K_1 =  \f{y_1 y_2}{|y|^4}, K_2 = \f{y_1^2 - y_2^2}{2|y|^4}$, 
and $|y| \geq 2 |x|$,  we have the following symmetrized estimates, which are proved in Part II \cite[Appendix B.1]{ChenHou2023b} using ideas similar to Taylor expansion
\begin{align}\label{eq:ker_decay1}
& |\na_x \big( G_s(x, y) - \f{8 y_1 y_2 x_1 x_2}{|y|^4} \big) | \les \f{|x|^3}{|y|^4}, \quad  | K_{1,s}(x, y) - \f{4 y_1 y_2}{|y|^4} |  \leq \f{6 |x|^2}{ \mathrm{ \Den}^2(x, y) } , \\
& | \pa_{x_1} K_{1,s}(x, y)  |  \leq \f{ 12 x_1}{ \mathrm{\Den^2(x, y)} }, \quad  | K_{1,s}(x, y) - \f{4 y_1 y_2}{|y|^4} -  2 (x_1^2 - x_2^2) \pa_{y_1}^2 K_1(y) |
\leq  \f{10 \sqrt{2} |x|^4}{ \mathrm{ \Den}^3(x, y)} ,  \notag 
\end{align}
where $\mathrm{Den}$ is defined below and satisfies $\mathrm{Den}(x, y) \asymp |y|^2$ for $|y| \geq 2|x|$
\[
\mathrm{Den}(x, y) \teq \tts{\sum}_{i=1,2} \, \min_{ |z_i| \leq x_i} |y_i- z_i|^2 = \tts{\sum}_{i=1,2} \, (\max( y_i - x_i, 0))^2, \quad \forall x, y \in \R^2_{++}.
\]

Similarly, we can obtain the following estimates for $K_{2,s}$
{\small
\beq\label{eq:ker_decay2}
\bal
&|K_{2,s}| \leq  \f{12 x_1 x_2}{ \mathrm{ \Den}^2} , 
\quad 
| \pa_{x_i} K_{2,s}(x, y) | \leq  \f{12 x_{3-i}}{ \mathrm{ \Den}^2} , \quad 
\B|K_{2,s} - 4 x_1 x_2 \cdot \f{12 y_1 y_2 (y_1^2 - y_2^2)}{ |y|^8} \B|
\leq  \f{ 40 x_1 x_2 |x|^2  }{ \mathrm{ \Den}^3} , \\
\eal
\eeq
}

\subsubsection{Proof of \eqref{eq:E_ker_est:b}, \eqref{eq:E_ker_est:c} in Lemma \ref{lem:E_ker} }\label{app:E_ker}

Below, we prove \eqref{eq:E_ker_est:b}, \eqref{eq:E_ker_est:c} in Lemma \ref{lem:E_ker}. The proof follows by decomposing the integral into the singular region $Q_2$, 
far-field region $Q_1$, and near-origin, non-singular region $Q_3$, and by symmetrizing the integral in $Q_1$.

\paragraph{\bf Proof of \eqref{eq:E_ker_est:b}}
Let  $W$ be the odd extension of $\om$ from $\R_2^+$ to $\R_2 $. 
By the Green function 
{\small 
\[
\phi = - \f{1}{2\pi} \int \log|x- y| W(y) dy, \quad \pa_{12} \phi(0) = \f{4}{\pi} \int_{\R_2^{++}} \f{y_1 y_2}{ |y|^4} W(y) dy,
\]
}
We only need to prove the estimate for the derivatives. The estimate for $\phi$ can be obtained by integration from $0$ to $|x|$.  Without loss of generality, we estimate $\pa_1 ( \phi - x_1 x_2 \pa_{12} \phi(0))$.

Similar to the proof of the first estimate in Lemma \ref{lem:E_ker}, we use the partition 
\beq\label{eq:ker_Q}
Q_1 = \{ y : |y| \geq 2|x|  \}, \quad Q_2 = \{ y: |y-x| \leq  |x| / 2 \},  \quad Q_3 = (Q_1 \cup Q_2)^c .
\eeq

Denote $M = || \om(\cdot) ( |y|^{-\al} + | y|^{\b}) ||_{\inf}$. We have 
\[
 |\om(y)| \leq M \min( |y|^{\al}, |y|^{-\b} ),  \quad \forall \, y \in \R^2_+.
 \]

In $Q_1$, we combine the estimate of $\phi$ and $-\pa_{12} \phi(0) \cdot x_1 x_2$. Symmetrizing the kernel, we get 
 \[
 \bal
&\phi - \pa_{12} \phi(0) \cdot x_1 x_2
= -\f{1}{2\pi} \int_{\R_2^{++}} K(x, y) W(y) dy , \\
& K(x, y) = \log|x-y| + \log|x+y|- \log | (x_1 -y_1, x_2 + y_2)| - \log |(x_1 + y_1, x_2 - y_2)| + \f{8 y_1 y_2 x_1 x_2}{|y|^4},
\eal
 \]
for $x, y \in \R^2$. For $\al < 2$, applying \eqref{eq:ker_decay1} $| \pa_{x_1} K(x, y)| \les  \f{|x|^3}{|y|^4} $ for $|y| \geq 2|x|$, we obtain
\beq\label{eq:E_ker2}
\bal
\B| \int_{Q_1} \pa_{x_1} K(x,y) W(y) dy \B| 
&\les M |x|^3  \int_{ |y| \geq 2|x| }  |y|^{- 4 } \min( |y|^{\al}, |y|^{-\b}) dy  \\
 &= M |x|^3 \min( |x|^{-2 + \al}, |x|^{-2 - \b} )  
  = M \min( |x|^{1 + \al} , |x|^{1 - \b}). 
 \eal
\eeq
Note that if $\al = 2$, the integral for small $y$, $\int_{ 2 |x| \leq |y| \leq 1}  |y|^{-2} dy$, leads to a $\log$ factor. 

In $Q_2, Q_3$, we estimate two integrals separately. For $\phi$ in $Q_2$, we have  $|x-y| \leq |x| / 2$ and  $|x| \asymp |y|$ in $|x-y| \leq |x| / 2$. Thus, we get 
\[
\bal
 \int_{|x-y|   \leq |x|/2  } | \pa_{x_1} \log |x- y|  W(y) | dy
 & \les M  \min( |x|^{\al} , |x|^{-\b } )  \int_{|x-y| \leq |x|/2 } |x-y|^{-1} dy  
 \les  M \min( |x|^{\al} , |x|^{-\b } )  |x|  .
 \eal
\]

In $Q_3$, we get $|x|/2 \leq |x-y| \leq 3 |x|, |y| \leq 2|x|$. It follows
\[
\bal
 \int_{ Q_3} | \pa_{x_1} \log |x- y|  W(y) | dy
 & \les M  |x|^{-1}  \int_{ |y| \leq 2 |x|} \min( |y|^{\al}, |y|^{-\b} ) dy \\
 & \les  M  |x|^{-1} \min( |x|^{2 + \al}, |x|^{2-\b})
 \les M \min( |x|^{1 + \al}, |x|^{1-\b}).
 \eal
\]

For $\phi_{xy}(0) x_1 x_2$, we have 
\[
\int_{ |y| \leq 2|x|} \B| \pa_{x_1} \f{ x_1 x_2 y_1 y_2}{ |y|^4}  W(y)\B|  dy 
\les M |x| \int_{ |y| \leq 2|x|}  |y|^{-2} \min( |y|^{\al}, |y|^{-\b})
 = M|x| \min( |x|^{\al}, |x|^{-\b}).
\]

Combining the above estimates, we prove  \eqref{eq:E_ker_est:b} in Lemma \ref{lem:E_ker}.

\paragraph{\bf Proof of \eqref{eq:E_ker_est:c}}
If $ \al \in (2, 3]$, the integrand in \eqref{eq:E_ker2} is integrable near $|x| = 0$. We yield 
\[
|\int_{Q_1} \pa_{x_1} K(x, y) W(y) dy| \les M |x|^3.
\]

Since $\al > 2, \b >0$, 
\eqref{eq:E_ker_est:c} follows from the above estimates in $Q_1, Q_2, Q_3$ and 
\[
 \min( |x|^{1 + \al}, |x|^{1-\b}) \les |x|^3 .
\]

\subsubsection{Proof of Lemma \ref{lem:E_vel}}\label{app:E_vel}

The proof is similar to Section~\ref{app:E_ker}, but now the \(X_1\)-norm provides both
weighted \(L^\infty\) control and the weighted \(C^{1/2}\) control needed in
the singular region \(Q_2\).

First, we consider the second estimate in Lemma \ref{lem:E_vel}. The bound by $ || \Om||_{X_1}$ follows from embedding. Let $W$ be the odd extension of $\Om$ from $\R_2^+$ to $\R_2$, as in \eqref{eq:ext_w_odd}. 

Recall the $X_1$-norm from \eqref{norm:X}.  We abuse the notation of $X_1$-norm for $W$ and we have
\[
 \| W \|_{X_1} =  \max( || W \vp_1 ||_{L^{\inf}(\R^2)},  || W \vp_{g,1} ||_{L^{\inf}(\R^2)} , 
[ W \psi_1 ]_{C_{g_1}^{1/2}(\R^2_+)}  ) ,
\quad 
 \| W \|_{X_1} = \| \Om \|_{X_1} . 
\]

Recall $-\D \Psi = W$ for $x \in \R^2_+$. We focus on the bound by $|x|^2 \| W \|_{X_1}$ and assume that $|x|\leq 1$. We consider the estimate for $\pa_{12} \Psi- \Psi_{12}(0)$. Firstly, we have 
\[
\pa_{12} \Psi =  \f{1}{\pi} \int K(x- y) W(y) dy , \, K(z) = \f{z_1 z_2}{|z|^4},
\quad  \pa_{12} \Psi(0) = \f{4}{\pi} \int_{\R_2^{++}} K(y) W(y) dy. 
\]
From \eqref{wg:hol}-\eqref{wg:linf_grow}, we get
\beq\label{eq:E_vel_ptest}
\max( \vp_1, \vp_{g, 1}) \gtr |x|^{-\al}, |x|^{\b}, \quad \al = 2.9, \, \b = 1/16, \quad  \psi_1 \asymp |x|^{-2}, \mathrm{ \ for \ } |x| \les 1.
\eeq

Following the standard partition \eqref{eq:ker_Q}, the symmetrization argument similar to that in the proof of Lemma \ref{lem:E_ker} or Section \ref{app:E_ker}, and the estimate \eqref{eq:ker_decay1}, we obtain 
\beq\label{eq:E_vel_far1}
 \B| \int_{|y| \geq 2|x|} K(x- y) W(y)  d y - 4 \int_{|y| \geq 2|x|, \R_2^{++}} K(y) W(y) dy\B| \les |x|^2 || W ||_{X_1}.
\eeq

For $|y |\leq 2|x|$, we have
\[
\bal
\B|  \int_{ |y|\leq 2|x|} K(y) W(y) dy \B| &\les \int_{|y| \les 2|x|}  \min(|y|^{\al}, |y|^{-\b}) |y|^{-2} dy  || W ||_{X_1} \les |x|^{\al} || W ||_{X_1} 
\les |x|^2 || W ||_{X_1}  \\
\eal
\]

In region $ Q_3 \subset \{ y: |x-y|\asymp |x|\}$, we have 
\[
\bal
|\int_{ Q_3} K(x- y) W(y )dy|
&\les |x|^{-2} \int_{ |y| \leq 3 |x|} |W(y )| dy
\les |x|^{-2} \int_{ |y| \leq 3|x|} \min( |y|^{\al} , |y|^{-\b} ) dy || W||_{X_1} \\
&\les || W||_{X_1} \min( |x|^{\al}, |x|^{-\b}).
\eal
\]

In the singular region $Q_2 = \{ y: |x-y|\leq |x|/2\}$, for any $|s|, |t| \leq |x|/2$, and $|x|\leq 1$, we have 
\[
\bal
& | W(x+ s) - W( x + t )|
= | (W \psi_1 \cdot \psi_1^{-1} )(x+ s) - (W \psi_1 \cdot  \psi_1^{-1} )(x+ t) | \\
\leq & |W \psi_1(x+s) - W \psi_1( x+ t)| \psi_1^{-1}(x+s) 
+ | \psi_1^{-1}(x+s)  - \psi_1^{-1}( x+ t) | \cdot | W \psi_1(x+t)  |  \teq I_1 + I_2.
\eal
\]

 Since $ |\na \psi_1^{-1}(x) | \les |x|$ by \eqref{wg:hol}, for $|s|, |t| \les |x|/2$, we get 
\[
| \psi_1^{-1}(x+s)  - \psi_1^{-1}( x+ t) |  \les |x| |s-t|.
\]

Since $W$ is odd extension of $\Om$ in $y$ and $\Om \psi_1 \in C_{x_1}^{1/2}(\R_2^+)$, we obtain 
$ [ W \psi_1 ]_{C_{x_1}^{1/2}(\R^2)} =  [ \Om \psi_1 ]_{C_{x_1}^{1/2}(\R_+^2)}$. 
Using the above estimates and $| \psi_1 W| \les  \f{\psi_1}{\vp_1} || W \vp_1||_{\inf} \les \f{\psi_1}{\vp_1} || W ||_{X_1}$, we yield
\[
|I_1 | \les || \Om||_{X_1} |x|^2 |s- t|^{1/2}, \quad I_2 \les |s-t| |x| \cdot \psi_1(x+t) / \vp_1(x+t)
\les |s- t| |x| \cdot |x|^{\al-2} = |s-t| |x|^{\al -1}.
\]

 Using the symmetry of the kernel that $K(s)$ is odd in $s_1, s_2$ , we get 
\[
\bal
& \B| \int_{|x-y|\leq |x|/2} K(x- y) W(y) dy \B| 
= \B| \int_{|s|\leq |x|/2, s_1 \geq 0} K(s) \big( W(x + s) - W(x + (-s_1, s_2)) \big) ds \B| \\
\les &  |x|^2 \int_{|s| \leq |x|/2} |s|^{-3/2} ds 
+ |x|^{\al-1} \int_{ |s| \leq |x|/2} |s|^{-1} ds \les |x|^{5/ 2} + |x|^{\al} \les  |x|^2.
\eal
\]
For small $|x|$, we can improve the above estimate by optimizing the window $a(x)$ for the singular region $|x-y| \leq a(x)$. We complete the estimate of the second inequality in Lemma \ref{lem:E_vel}. 

Next, we consider the first estimate in Lemma \ref{lem:E_vel}. For $W \in {X_1}$, using \eqref{eq:u_appr_near0_K}, we have 
\[
\pa_{1112} \Psi(0) =   \f{4}{\pi} \int_{\R_2^{++}}\pa_{11}K(y)  W(y) dy, \quad \pa_{11} K(y) = \f{12 y_1 y_2(y_1^2  -y_2^2) }{|y|^8 } .
\]

The proof is completely similar. We use the above partition of the domains. In $Q_1$, we use the 
symmetrizing estimates \eqref{eq:ker_decay1}, \eqref{eq:ker_decay2} to yield 
\[
 \B| \int_{|y| \geq 2|x|}  K(x- y) W(y) dy   - \int_{|y| \geq 2|x|, \R_2^{++}} ( 4 K(y)  +  2(x_1^2 - x_2^2) \pa_{11} K(y)  )W(y) dy\B| \les |x|^{5/2} || W||_{X_1}.
\]
In $Q_2, Q_3$, we estimate the integrals separately, and use the above estimates and 
\[
|x|^2\int_{|y| \les |x|} \B| \f{24 y_1 y_2(y_1^2  -y_2^2) }{|y|^8 } \cdot W(y) \B| dy
\les |x|^2 \int_{|y| \les |x| } |y|^{-4} \min( |y|^{\al}, |y|^{-\b}) dy || W ||_{X_1}
\les \min( |x|^{\al}, |x|^{-\b}  ) \| W \|_{X_1},
\]
for any $\b > 0$. Using the estimate \eqref{eq:E_vel_ptest} and the integral formulas of $\pa_{12} \Psi(0), \pa_{1112} \Psi(0) $, we get
\[ 
|\pa_{12} \Psi(0)| + |\pa_{1112} \Psi(0)|\les || W \vp_{g, 1}||_{L^{\inf}} \les || W||_{X_1}.
\]

\subsection{Ideas and examples of rigorous weighted estimates}\label{app:idea_num_est}

In this section, we outline some key methods for rigorous weighted estimates of functions and residual error using numerical analysis. More detailed methods are developed in \cite{ChenHou2023b} for rigorous numeric.

\vs{0.05in}
\paragraph{\bf Rigorous piecewise bounds}
We need to derive rigorous and tight piecewise bounds of various quantities involving the approximate steady state, singular weights, and several explicit functions. 
One of the main ideas is to use the second order error estimate 
\beq\label{est_2nd}
  \max_{ x \in [x_l, x_u] }|f(x)|  \leq \max( |f(x_l) | , |f(x_u)| ) + \tf{ 1}{8} h^2 || f_{xx}||_{L^{\inf}([x_l,x_u])}, \quad h = x_u - x_l,
  \eeq
to obtain a piecewise sharp bound of $f$ on $ [x_l, x_u]$, see e.g. \cite[Appendix {\secpiecepol}]{ChenHou2023b}. 
If we can obtain a rough bound for $f_{xx}$ on $[a, b]$, by partitioning $[a, b]$ into small intervals $[ x_{il},x_{iu}]$, evaluating $f$ on finitely many grid points $x_{il}, x_{iu}$, and using 
\eqref{est_2nd}, we obtain a tight bound of $f$ on $[a,b]$. Similarly, by estimating $ \pa_x^k f$ and applying \eqref{est_2nd} recursively, we obtain a tight bound for $\pa_x^i f$ with $i \leq k-1$. Note that for a polynomial $f$ with degree less than $d$, we have $\pa_x^k f \equiv 0 , k \geq d+1$. Our approximate steady state is represented by piecewise polynomials (See Section \ref{sec:ASS}), and we apply these estimates. 
The above estimates extend directly to multivariable functions $f(x,y)$.

For several explicit functions and more complicated functions, we bound higher derivatives inductively and use the Leibniz rule and triangle inequalities. We also derive higher-order error estimates with error terms $C_k h^k, k=3,4,5$, using numerical analysis. See more details in  Paper II \cite[Appendix C]{ChenHou2023b}. To track the round-off error in the computation, we use interval arithmetic  with  package INTLAB \cite{Ru99a} for MatLab and refer to \cite{moore2009introduction,rump2010verification} for more detailed discussions.

\vs{0.1in}

\paragraph{\bf Rigorous weighted estimates}

We next use two examples  to illustrate rigorous weighted estimates for functions represented by basis expansions. These methods are \emph{robust to numerical error} and apply to bound the weighted norm of local residual error 
$\bar \cF_{loc, i}, \cR_{loc, i}$ in \eqref{eq:lin_U} 
and the residual error $\e_1$ \eqref{eq:uerr_dec1} for Poisson equation. The key idea is to control the weighted norm
through derivative bounds. For instance, if $f=O(|x|^3)$ near $x=0$, then for
$|\xx|\le 1$,
\beq\label{eq:key_weight_est}
   |\xx|^{-3}|f(\xx)|
   \les
   \|\na^3 f\|_{L^\infty(|\xx|\le 1)} .
\eeq

\paragraph{\bf Example 1}

We consider examples related to model problem in Section \ref{sec:model_low_rank}.
Suppose that 
\[
  \overline F = \bar f \cdot \bar g ,
  \quad \bar f = c_1 x + \tts{ \sum_{i + j = 2}}  \, c_{ij} x^i y^j,
  \quad \bar g = d_1 y + \tts{ \sum_{i + j = 2} } \, d_{ij} x^i y^j ,
\]
where the coefficients $c_{ij}, d_{ij}$ are given as floating point values, and $c_1, d_1 \neq 0$. 
Clearly, we have $\pa_{xy} \overline F(0) = c_1 d_1$. We aim to bound $(\overline F - \pa_{xy} \overline F(0) \chi_\dag ) |\xx|^{-3}$ for $|\xx| < 1$, where \(\chi_\dag\) is a smooth cutoff with
\(\chi_\dag({\bf x})=xy+O(|{\bf x}|^3)\) near the origin and \(\chi_\dag\) is compactly supported.

Due to numerical error, the numerical evaluation of $\pa_{xy}\overline F(0)$ may differ from its analytic value $\pa_{xy}\overline F(0)=c_1d_1$.  However, using the \emph{analytic} definition of $\pa_{xy} \overline F(0)$, we obtain the analytic vanishing conditions  $\overline F -\pa_{xy} \overline F(0) \chi_\dag = O(|\xx|^3)$. Thus, using the Leibniz rule, we obtain
\[
  \max_{|\xx| < 1}  |\overline F -\pa_{xy} \overline F(0) \chi_\dag | \cdot |\xx|^{-3}
  \les \max_{|\xx| < 1} ( | \na^3  \overline F |  + |\pa_{xy} \overline F(0)| \cdot |\na^3 \chi_\dag| )
  \les \max_{|\xx| < 1} ( \tts{\sum}_{k + l = 3} |\na^k \bar f| \cdot |\na^l \bar g| 
    + |\pa_{xy} \overline F(0)| \cdot |\na^3 \chi_\dag| ).
\]
The upper bound is \emph{finite for any} $c_1, d_1, |\pa_{xy} \overline F(0)|$. The numerical evaluation enters the rigorous estimates by bounding the upper bounds. Using error estimates, 
e.g. the 2D analog of \eqref{est_2nd}, we obtain sharp piecewise bounds for $\na^i \bar f, 
\na^j \bar g, \na^k \chi_\dag$.
Using these bounds and the above estimate, we bound $ |\overline F -\pa_{xy} \overline F(0) \chi_\dag | \cdot |\xx|^{-3}$ rigorously by some bounded constant.

Away from $\xx=0$, we derive piecewise bounds for $\overline F - \pa_{xy} \overline F(0) \chi_\dag$ and their weighted estimates.
Although we still cannot numerically evaluate $ \pa_{xy} \overline F(0)$ exactly, the error is not  significantly amplified by the weight since the weight is less singular away from $\xx=0$. Optimizing these two estimates yields sharp piecewise bounds for 
$|\overline F -\pa_{xy} \overline F(0) \chi_\dag | \cdot |\xx|^{-3}$. 

Although numerical evaluation may not produce $\pa_{xy} \overline F(0)$ \emph{exactly}, this does \emph{not} lead to an unbounded weighted estimate. The key point is that the cancellation uses \emph{exact} analytic derivative of the represented function $\overline F$, rather than a raw floating-point approximation. By enforcing the exact vanishing order analytically,  the weighted estimate is  rigorous and bounded.

\vspace{0.1in}

\paragraph{\bf Example 2: Weighted estimates for the error}

Next, we discuss the weighted estimates for the corrected error $ (\pa_s - \cL) \wh g^{(2)}$ for $\wh g^{(2)}$ defined in \eqref{def:g2}. 
Since $\wh g^{(2)} = \wh g^{(1)} + a(t) \chi_{\dag,2}$ \eqref{def:g2}, using the vanishing \eqref{eq:g_cond} for $\wh g = \wh g^{(2)}$, \eqref{eq:key_weight_est}, 
and triangle inequality, we \emph{decouple} $a(t)$ from the numerically constructed functions
\[
\max_{|\xx| \leq 1} |\xx|^{-3} | (\pa_s - \cL) \wh g^{(2)} |
\les \max_{|\xx|\leq 1} | \na^3 (\pa_s -\cL) \wh g^{(2)} |
\les \max_{|\xx|\leq 1} ( | \na^3 (\pa_s -\cL) \wh g^{(1)} |
+ |a(t) | \cdot |\na^3 \cL \chi_{\dag,2}| 
+| a^{\pr}(t)| \cdot |\na^3 \chi_{\dag ,2}|)
\]
Since $\wh g^{(1)}, \chi_{\dag, 2}$ are represented by explicit basis functions, we bound 
$ | \na^3 (\pa_s -\cL) \wh g^{(1)} | , |a(t) | \cdot |\na^3 \cL \chi_{\dag,2}| ,| a^{\pr}(t)| \cdot |\na^3 \chi_{\dag ,2}| $ rigorously using estimates similar to \eqref{est_2nd}. We bound $a(t)$ using the \emph{analytic} formula \eqref{def:a} and the ODE \eqref{eq:ODE2}
\[
a^{\pr}(t) = \bar \lam a(t) - \pa_{xy} \cE(t, 0), 
\quad  |a(t)| \leq  e^{\bar \lam t}  |\pa_{xy} \cE_0(0)|
+ \int_0^t e^{ \bar \lam (t- s) }  |\pa_{xy} \cE(s, 0 )|  d s ,
\quad \bar \lam = \bar c(0 ) - \pa_x \overline V_1(0) - \pa_y \overline V_2(0).
\]

For the ODE system in \eqref{eq:bous_ODE_xy} and in \cite[Section 3.3]{ChenHou2023b} near $\xx=0$, which corresponds to \eqref{eq:ODE2}, we have $\bar \lam < -5$. 
Thus, by bounding $|\pa_{xy} \cE_0(0)|$ and $ |\pa_{xy} \cE(s, 0 ) |$ for $s \in \supp (\cE( \cdot, 0) )$ 
\footnote{
The approximate space-time solution $\hat G$ and its error 
$\cE(s, 0)$ in Section \ref{sec:low_rank_cor} have temporal support $s \in [0,12]$.
}
and using the exponential decay $e^{\bar \lam s}$ in $s$, we bound the factor $a(t)$ for \emph{any} $t \geq 0$.

While the above triangle inequality 
overestimates
$\nabla^3(\pa_s-\cL)\wh g^{(2)}$, since $a(t)$ is of the size
of the error, this overestimate is minor.  Away from $\xx=0$, similarly, we estimate 
\[
| |\xx|^{-3} (\pa_s - \cL) \wh g^{(2)} |
\leq | |\xx|^{-3} (\pa_s - \cL) \wh g^{(1)} |
+ |a^{\pr}(t)| \cdot |\xx|^{-3} |\chi_{\dag, 2}|  
+ |a(t) | \cdot | |\xx|^{-3}  \cL \chi_{\dag, 2} | .
\]

By deriving sharp piecewise bounds for $(\pa_s - \cL) \wh g^{(1)},  \cL \chi_{\dag, 2}, \chi_{\dag ,2}$, we obtain 
a second estimate. Optimizing the two estimates yields a sharp piecewise bound for the weighted error rigorously.  
A key point is that we \emph{do not} need to evaluate 
or approximate $a(t)$ for rigorous weighted bounds. Since $a(t)$ is very small, using the above estimates, we only 
require an \emph{upper bound} for $a(t)$.

\section{Inequalities for nonlinear stability }\label{sec:ineq}

\subsection{Inequalities for nonlinear stability}

In this Appendix, we present all the inequalities for nonlinear stability in Lemma \ref{lem:PDE_nonstab} with the final energy $E_4$ \eqref{energy4}. 
We have performed energy estimates in Section \ref{sec:EE}. 
The estimates of other nonlinear terms in the H\"older estimates are similar, and we refer them to the supplementary material I \cite[Section {\suppsecnonest}]{ChenHou2023aSupp}. Appendix~\ref{sec:ineq_linf} records the weighted
\(L^\infty\) stability inequalities.  Appendix~\ref{sec:ineq_hol} records the weighted
\(C^{1/2}\) inequalities.  Appendix~\ref{sec:ineq_cw} records the 
ODE estimates for $c_\omega$-terms and the Taylor coefficients at the origin.  Finally,
Section~\ref{sec:hol_wg_est} records the elementary estimates for the H\"older weights \(g_i\).

We verify these inequalities for nonlinear stability with computer assistance, and the codes can be found in \cite{ChenHou2023code}. The codes are implemented in MatLab with package INTLAB \cite{Ru99a} for interval arithmetic. 
The kernel-integral estimates for the weighted $L^{\infty}$ and $C^{1/2}$ estimates of nonlocal terms 
(see Section \ref{sec:u_comp_idea} for the ideas), and the construction and estimates of the approximate space-time solutions (see Sections \ref{sec:low_rank_cor}, 
\ref{sec:W2} for the ideas  and \cite[Section 3]{ChenHou2023b} for the detailed methods) are performed in parallel using the Caltech High Performance Computing 
\footnote{See more details for Caltech HPC Resources \url{https://www.hpc.caltech.edu/docs/resources.html}
}. 
Other computer-assisted estimates and main part of the verifications are done in Mac Pro (Rack, 2019) with 2.5GHz 28‑core Intel Xeon W processor and 
768GB (6x128GB) of DDR4 ECC memory.

\vs{0.05in}
\paragraph{\bf{Variables and Notations}}

Recall $W_1$ from \eqref{eq:EE_W1}, $\UU, \td \UU, \UU_A$ from \eqref{eq:U}, the energies $E_i, i\leq 4$ \eqref{energy1}, \eqref{energy2}, \eqref{energy3}, \eqref{energy4}, the weights $\rho_{10} = \rho_{01}, \rho_{ij} = \psi_1, i+j=2$ \eqref{wg:elli} for the $L^{\inf}$ estimates of $\uu_A, (\na \uu)_A$, and $\psi_u, \psi_1$ \eqref{wg:elli} for the H\"older estimate of $\psi_u \uu_A, \psi_1 (\na \uu)_A$. The notation $f_A(\om) = \td f(\om) - {\td {\hat f}}(\om) $ is introduced in \eqref{eq:vel_rem}, where $\td{\hat f}$ is the finite-rank approximation of $f(\om)$ \eqref{eq:u_appr}. 
Note that in general $\na (\UU_A) \neq (\na \UU)_A$.
The weights $\psi_i, \vp_i, \vp_{g, i}$ are defined in \eqref{wg:hol},  \eqref{wg:linf_decay}, \eqref{wg:linf_grow}.

The notation below keeps the long pointwise inequalities compact and more-readable.  We use $\cT$ for quantities associated with transport terms.  We introduce  $\cT$ to bound $ b \cdot \na f $
\beq\label{eq:tran_ln}
\cT(b, f)(x)    \teq |b_1 | \cdot |f_x| + |b_2| \cdot |f_y| = \f{1}{\rho_{10}} 
 ( |b_1 \rho_{10}| \cdot |f_x| + |b_2 \rho_{10}| \cdot |f_y|  ) . \\
\eeq
For $b = \uu_A = (u_A, v_A)$, $ \rho_{10} \cT(\uu_A, f)$ agrees with $\cT_{u}(f)$ defined in Section \ref{sec:linf_decay}. If $b = \UU$ and $f = \rho$ is a weight, 
$\cT(\UU, \rho)$ agrees with $\cT_{\mw{wg}}(\UU, \rho)$ defined in \eqref{eq:linf_non_tran2}.

We derive a piecewise weighted estimate $\rho_{ij} f_{ij}$ for a singular weight $\rho_{ij}$ associated with $f_{ij}$ and unweighted estimate for $f_{ij}$. Then, we apply two estimates to bound $\rho \cdot f_{ij}$ 
\[
 \rho f_{ij} = (\rho / \rho_{ij}) \cdot  (\rho_{ij} f_{ij}).
\]
where $f_{ij}$ denotes $u_A, v_A, u_{x,A}, v_{x,A}, u_{y,A}$ for $(i, j) = (0, 1),(1,0),(1,1),(2,0),(0, 2)$  as in \eqref{eq:u_linf_est1}, and $\rho_{10}, \rho_{20}$ are given in \eqref{wg:elli}.  The above two bounds are slightly different since $\max_{x \in Q} | f^{-1}(x)| \cdot \max_{x\in Q} | f(x) |\neq 1$. 
The second bound is useful near $x=0$ since both terms are regular. Note that we further establish weighted estimate for $\rho_4 u_A$ using $|| \om_1 \vp_1||_{\infty}$.

 Recall that we have modified the decomposition of linear and nonlinear terms in \eqref{eq:lin_U} and discussions therein to simplify the nonlocal error estimates. Below, the estimates are based on the decompositions in \eqref{eq:lin_U} and $\bar \uu^N, \bar c_{\om}^N$ etc.

In the linear $L^{\inf}(\vp)$ estimates in Section \ref{sec:linf_decay}, and $L^{\inf}(\vp_4)$ estimate in Section \ref{sec:hol_Q1}, we only use energy $E_1$ \eqref{energy1}. 
In the linear H\"older estimate in Section \ref{sec:EE_hol}, we only use 
$E_2$ \eqref{energy2}. In the $L^{\inf}(\vp_{g,i})$ estimates in Section \ref{sec:linf_grow}, we only use $E_3$ \eqref{energy3}. In the remaining energy estimates for functionals and nonlinear terms in Sections \ref{sec:W2}-\ref{sec:non}, we use the full energy $E_4$ \eqref{energy4}.

\subsubsection{Weighted $L^{\inf}$ estimate}\label{sec:ineq_linf}

This subsection collects weighted $L^{\infty}(\vp_i), L^{\infty}(\vp_4), L^{\infty}(\vp_{g, i})$ estimates from Sections~ \ref{sec:linf_decay},
\ref{sec:linf_grow}, \ref{sec:non}  with modification for equation \eqref{eq:lin_U} in Section 
\ref{sec:comb_vel_err}, and present the stability inequalities in Lemma~\ref{lem:PDE_nonstab} associated with these estimates and energy $E_4$.

\vs{0.1in}
\paragraph{\bf{Weighted $L^{\inf}(\vp_i)$ estimate}}
We establish the following linear weighted $L^{\inf}(\vp_i)$ estimates for the bad terms $B_{\mf{modi},i}$ \eqref{eq:bad_lin} in Section \ref{sec:linf_decay} 
\beq\label{linf:WG1_lin}
\bal
L_1(\vp_1) & \leq   \f{\vp_1}{\vp_2}  E_1 +  \vp_1 \cT(\UU_A , \bar \om)     , \\
L_2(\vp_2) & \leq   \f{\vp_2}{\vp_3} |\bar v^N_x| E_1 + \vp_2 \cT(\UU_A , \bar \th_x)
+ \f{\vp_2}{\rho_{20}} ( |U_{x, A} \rho_{20} \bar \th_x| 
+ |V_{x, A} \rho_{20} \bar \th_y| ) ,  \\
L_3(\vp_3) & \leq   \f{\vp_3}{\vp_2} |\bar u^N_y|  E_1
+ \vp_3 \cT(\UU_A , \bar \th_y)
+ \f{\vp_3}{ \rho_{20}}  
( |U_{y, A} \rho_{20} \bar \th_x| + |U_{x, A} \rho_{20} \bar \th_y| )  ,
\eal
\eeq
where $\psi_1 = \rho_{ij}, i+j=2$.  
These  estimates for $B_{\mf{modi}}$ 
are obtained by the same arguments as in Sections ~\ref{sec:linf_decay} 
for \eqref{eq:lin_main} with coefficients $(\bar \uu, \na \bar \uu, \bar c_{\om})$ replaced by $(\bar \uu^N, \na \bar \uu^N, \bar c_{\om}^N)$, 
and with estimates of $\uu_A, (\na \uu)_A$, replaced by  those of $\UU_A, (\na \UU)_A$. 
The same ideas apply to the remaining estimates.

We keep $ \UU_A,  (\na \UU)_A$ terms in \eqref{linf:WG1_lin} to simplify the notations. 
The terms involving $\uu_A(\om_1), (\na \uu)_A(\om_1)$ can be bounded by $ \CCs( f_{ij}, a ) E_1$
via estimates \eqref{eq:u_linf_est1}, \eqref{eq:u_linf_est12} with  some weight $a$ introduced in \eqref{eq:EE_oper} in Section \ref{sec:linf_decay}. See discussion between \eqref{eq:U} and \eqref{eq:comb_vel1} for the estimate of $\UU_A, (\na \UU)_A$. For \eqref{eq:non_dec1} and the error \eqref{eq:bous_errM}, \eqref{eq:W2_errM},
we have nonlinear energy estimates 
\bseq\label{linf:WG1_nlin}
\beq
\bal
& \cN\cF_i(\vp_i) \leq  \uds{ \eqref{eq:linf_non_tran} }{\cT_\mw{wg}(\UU, \vp_i) E_4 } +  \uds{ \eqref{def:N_W1}  }{ | N_{W_1, i} | } \cdot \vp_i  + \uds{\eqref{eq:W2_M2} }{ | N_{\hat W_2, i} | } \cdot \vp_i  +  |\bar \cF_{loc, i} \vp_i| + |\cR_{loc, i} \vp_i | , 
\\
 & | N_{W_1, 1} | \cdot \vp_1 \leq |U_x(0)| E_4, \\
& | N_{W_1, 2} |  \cdot \vp_2  \leq  \big(  |U_x(0)| + | \td U_x | + \f{\vp_2}{\vp_3} | \td V_x | \big) E_4  
 , \\
& | N_{W_1, 3} | \cdot \vp_3 \leq  \big(  3 |U_x(0)| + | \td U_y | \f{\vp_3}{\vp_2} + | \td U_x | \big) E_4   , \\
 & | N_{\hat W_2,1} | \cdot \vp_1 \leq \vp_1\cT(\td \UU, \hat \om_2)  
+  |U_x(0)| \cdot |\hat W_{2,1, M} \vp_1| , \\
& | N_{\hat W_2, 2} | \cdot \vp_2 \leq \vp_2 \cT(\td \UU, \hat \eta_2) + 
|U_x(0)| \cdot |\hat W_{2, 2,  M} \vp_2|
+ \f{\vp_2}{\rho_{20}} ( | \td U_x  \rho_{20} | \cdot |\hat \eta_2|  
+ |\td V_x \rho_{20}|\cdot |\hat \xi_2| ) , \\
& | N_{\hat W_2, 3} | \cdot \vp_3 \leq \vp_3 \cT(\td \UU, \hat \xi_2) + 
|U_x(0)| \cdot |\hat W_{2, 3,  M} \vp_3|
+ \f{\vp_3}{\rho_{20}} ( | \td U_y  \rho_{20} | \cdot |\hat \eta_2|  
+ |\td U_x \rho_{20}|\cdot |\hat \xi_2| ) , 
\eal
\eeq
\eseq 
where we have used
 \eqref{eq:linf_non_tran}, and $\wh W_{2, \cdot, M} =  (\hat \om_{2, M}, \hat \eta_{2, M}, \hat \xi_{2, M} )$ are defined in \eqref{eq:W2_M2}.
 In \eqref{linf:WG1_nlin},  $L_i(\vp_i), \cN\cF_i(\vp_i)$ are only used to indicate the weighted $L^{\inf}(\vp_i)$ estimate of linear and nonlinear terms. We have used $V_{y, A} = - U_{x, A}, \td V_y = - \td U_x$. In \eqref{linf:WG1_lin}, \eqref{linf:WG1_nlin}, we \emph{do not} multiply the terms related to $\cT_u( \cdot )$ and $(\na \UU)_A$ by $E_4$ since we will further bound it using  
\[
|\uu_A(\om_1)| \leq C_1(x) E_1, \quad |(\na \uu(\om_1))_A| \leq C_2(x) E_1, \quad E_1 \leq  E_4 \leq E_* .
\]
by some computable constant $C_i(x)$.  Under the bootstrap assumption, we can combine the estimate of $\uu_A$ and the error part in $\UU_A$ \eqref{eq:U} in \eqref{eq:comb_vel1}. 
The same reasoning applies to the $\hat W_2, U$ terms. See Sections \ref{sec:comb_vel_err} and \ref{sec:W2}.

Recall the damping coefficient  $d^{\mf{num}}_{i, L}(\cdot)$ from \eqref{eq:dp_num}.  Substituting the estimates of $\UU_A, \UU, \hat W_2$ (see also \eqref{eq:comb_vel1}), bounding $E_4$ by $E_* = 5 \cdot 10^{-6}$ in \eqref{linf:WG1_lin}, \eqref{linf:WG1_nlin}, and applying Lemma \ref{lem:PDE_nonstab}, we obtain the nonlinear stability conditions for the $L^{\inf}(\vp_i)$ estimates 
\beq\label{stab:linf_WG1}
 - d^{\mf{num}}_{i, L}(\vp_i) E_* - L_i(\vp_i) - \cN\cF_i(\vp_i) \geq \lam E_*, \quad \forall \, x \in \R_2^{++},  \ \lam > 0.
\eeq
We do the same substitution in the following stability conditions.

\vs{0.1in}
\paragraph{\bf{$L^{\inf}(\vp_4)$ estimate}}

Recall $\vp_4 = \psi_1 |x_1|^{- \f{1}{2} }$, the weight $\f{\sqrt 2}{\tau_1}$ in $E_4$ \eqref{energy1}. We have established $L^{\inf}(\vp_4)$ linear estimates in Section \ref{sec:hol_Q1} and nonlinear estimates similar to $\cN\cF_1(\vp_1)$ \eqref{linf:WG1_nlin}
\begin{align}
\f{ \sqrt{2}} {\tau_1} L_1(\vp_4) & \leq \f{\sqrt{2}}{\tau_1} ( \f{\vp_4}{\vp_2} E_1  
+ \vp_4 \cT( \UU_A, \bar \om) ),  \label{linf:WG4_lin}  \\
 \f{\sqrt{2}}{\tau_1} \cN\cF_1(\vp_4) & \leq 
\big( \uds{ \eqref{eq:linf_non_tran} }{ \cT_\mw{wg}(\UU, \vp_4) }+ |U_x(0)| \big)  E_4 + 
\f{\sqrt{2}}{\tau_1} \cdot
  \f{\vp_4}{\vp_1} \big(  |  \uds{ \eqref{eq:W2_M2}  }{ N_{\hat W_2,1} }  \vp_1 |  + |\bar \cF_{loc, 1} \vp_1| + |\cR_{loc, 1} \vp_1 | \big) .
 \notag 
\end{align}

We do not multiply $\cT_\mw{wg}(\UU, \vp_4 ) E_4$ by $ \f{ \sqrt{2}}{\tau_1}$. See estimate \eqref{eq:linf_non_tran}.
Since we have weighted $L^{\inf}(\vp_1)$ estimates of similar terms $f \vp_1$ in \eqref{linf:WG1_lin}, \eqref{linf:WG1_nlin}, e.g. $ \cT(\UU_A, \bar \om) \vp_1, \bar \cF_{loc, 1}\vp_1, \cR_{loc, 1} \vp_1, | N_{\hat W_2, 1} | \vp_1$,  we can use such estimates and estimate $\f{\vp_4}{\vp_1}$ to bound $f \vp_4$. Similar to \eqref{stab:linf_WG1}, the stability conditions read 
\beq\label{stab:linf_WG4}
  - d^{\mf{num}}_{1, L}(\vp_4) E_* - \f{\sqrt{2}}{\tau_1} ( L_1(\vp_4)+ \cN\cF_1(\vp_4) ) \geq \lam E_* , \quad \forall x \in \R_2^{++},  \ \lam > 0.
\eeq

\vs{0.1in}
\paragraph{\bf{$L^{\inf}(\vp_{g,i})$ estimate} }

In the $L^{\inf}(\vp_{g,i})$ estimate, we estimate $\mu_{g, i}|| W_{1, i} \vp_{g, i}||_{L^{\inf}}$. Recall  $\mu_g = \tau_2(\mu_4, 1, 1)$ from \eqref{energy3}, \eqref{wg:EE}. We have established linear stability estimate in Section \ref{sec:linf_grow}
\beq\label{linf:WG1g_lin}
\bal
&  \mu_{g, 1} L_1(\vp_{g, 1}) \leq \mu_{g, 1} \B( \vp_{g,1} (  \vp_2^{-1} \wedge  \mu_{g, 2}^{-1} \vp_{g,2}^{-1}  ) E_3   +   \vp_{g,1}  \cT( \UU_A ,  \bar \om) \B)  , \quad \mu_g = \tau_2(\mu_4, 1, 1), \\
& \mu_{g, 2} L_2(\vp_{g, 2}) \leq \mu_{g, 2} \B( \vp_{g,2} | \bar v^N_x | (  \vp_3^{-1} \weg  \mu_{g, 3}^{-1} \vp_{g,3}^{-1} ) E_3 
+   \vp_{g,2}\cT (\UU_A , \bar \th_x)
+  \f{\vp_{g,2}}{ \psi_1}
( | U_{x, A} \psi_1  | \cdot |\bar \th_x| + |V_{x, A} \psi_1| \cdot |\bar \th_y| ) \B) , \\
& \mu_{g, 3} L_3(\vp_{g, 3}) \leq \mu_{g, 3} \B( \vp_{g,3} | \bar u^N_y | (  \vp_2^{-1} \weg  \mu_{g, 2}^{-1} \vp_{g,2}^{-1} ) E_3
+   \vp_{g,3}\cT (\UU_A , \bar \th_y)
+  \f{\vp_{g,3}}{ \psi_1}
( | U_{y, A} \psi_1| \cdot |\bar \th_x| + |U_{x, A} \psi_1| \cdot |\bar \th_y| ) \B) . \\
\eal
\eeq

We have used $\mu_{g,i}$ to rewrite the parameters in the $L^{\inf}(\vp_{g,i})$ estimate in Section \ref{sec:linf_grow} equivalently, so that the form of estimates is more symmetric in the parameters. We also use $\psi_1 = \rho_{20}$ and keep $u_{ij, A}$ in the estimates, which can be bounded by $C_{g, ij}(x)$ defined in Section \ref{sec:linf_grow}. Denote 
\beq\label{eq:W_igj}
W_{i, g, j} = \min( \f{ \vp_{g, j}}{ \vp_i}, \f{\vp_{g, j}}{ \mu_{g, i} \vp_{g, i} })
= \vp_{g, j}  ( \vp_i^{-1} \wedge (  \mu_{g, i} \vp_{g, i})^{-1} ).
\eeq
Using $| W_{1, i} \vp_i|, \mu_{g,i}|W_{1, i} \vp_{g, i}| \leq E_3$, 
we have $ |W_{1, i} \vp_{g, j}| \leq W_{i, g, j} E_3$, which motivates the above notation. Note that $
\mu_{g, i} W_{ i, g, i} \leq 1$. We can simplify some terms in the above estimate using $W_{i, g, j}$, e.g. $\vp_{g,1} (  \vp_2^{-1} \wedge  \mu_{g, 2}^{-1} \vp_{g,2}^{-1}  )$. We have the following nonlinear estimates similar to \eqref{linf:WG1_nlin}
\begin{align}
& \mu_{g, i} \cN\cF_i( \vp_{g,i}) 
\leq \uds{ \eqref{eq:linf_non_tran} }{ \cT_\mw{wg}(\UU, \vp_{g,i} ) E_4  }
+ \mu_{g, i} \vp_{g, i} \uds{ \eqref{def:N_W1} }{ | N_{ W_1,i} | } + \f{  \mu_{g, i} \vp_{g, i}}{\vp_i} ( \uds{ \eqref{eq:W2_M2} }{  | N_{\hat W_2, i} |} \cdot \vp_i   + |\bar \cF_{loc,i } \vp_i | + |\cR_{loc, i} \vp_i | ) , \notag  \\
 & \mu_{g, 1} \vp_{g, 1} | N_{ W_1,1}|
\leq |U_x(0)| E_4, \notag  \\
& \mu_{g, 2} \vp_{g, 2} | N_{ W_1, 2} |
\leq \mu_{g, 2} \B( (|U_x(0)| + | \td U_x | ) W_{2,g,2} + | \td V_x | W_{3, g, 2} \B) E_4 , \notag \\
& \mu_{g, 3} \vp_{g, 3} | N_{ W_1, 3} |
 \leq \mu_{g, 3}\B( | \td U_y | W_{2, g, 3} + (3 |U_x(0)| + | \td U_x | ) W_{3, g, 3} \B) E_4 .
 \label{linf:WG1g_nlin}
\end{align}
Note that the terms $\cT_\mw{wg}(\UU, \vp_{g,i} ) E_4$ are not multiplied by $\mu_{g, i}$. See estimate \eqref{eq:linf_non_tran}.
 We use the weighted $L^{\inf}(\vp_1)$ estimates of similar terms $f \vp_i$ in \eqref{linf:WG1_lin}, \eqref{linf:WG1_nlin}, e.g. $ \cT(\UU_A, \bar \om) \vp_i, \bar \cF_{loc, i}\vp_i, \cR_{loc, i} \vp_i$, and  further estimate $\f{\vp_{g, i}}{\vp_i}$ to bound $ f \vp_{g, i} = ( f \vp_i) \cdot \f{\vp_{g, i} }{\vp_i}  $. The stability conditions read 
\beq\label{stab:linf_WG1g}
 - d^{\mf{num}}_{i, L}(\vp_{g,i}) E_* - \mu_{g, i} L_i(\vp_{g,i}) -  \mu_{g, i}\cN\cF_i(\vp_{g, i}) \geq \lam E_*, \quad \forall x \in \R_2^{++},  \ \lam > 0.
\eeq

We record some useful inequalities implied by \eqref{stab:linf_WG1}, \eqref{stab:linf_WG1g}.
In the stability inequalities, e.g.  \eqref{stab:linf_WG1}, \eqref{stab:linf_WG1g},
the energies are all \emph{replaced by} $E_*$. 
The term $U_x(0)$ is bounded by estimate \eqref{eq:U_bd}: $|U_x(0)| < \mu_6 E_* +| u_x(\bar \e)(0)|$. Taking $x = (x_1,0)$ with $x_1 \to \infty$, \eqref{stab:linf_WG1}, \eqref{stab:linf_WG1g} implies 
\footnote{
Since $\bar \uu^N = \na^{\perp} \bar \phi^N$ satisfies 
$\na \bar \uu^N , |\bar \uu^N| / |\xx| \to 0$ by \eqref{solu:decay}, 
we get $d_{3,L}^{\numm}( \vp_3)(x_1,0) \to  2 \bar c_{\om}^N + \bar c_l \tf17,
d_{3,L}^{\numm}( \vp_{g,3})(x_1,0) \to  2 \bar c_{\om}^N + \bar c_l \al_{g, n}$ 
as $x_1 \to \infty$ by definition \eqref{wg:linf_grow}, \eqref{wg:linf_decay}.
Moreover, as $x_1 \to \infty$,  since $ \vp_{g, 3}(x_1,0) \gg \vp_3(x_1,0)$, we get 
$\mu_{g, 3} W_{3,g,3}(x_1, 0) \to 1$ by \eqref{eq:W_igj}. 
Thus, keeping positive $-d_{i,L}^{\numm}$-term,
and negative term $-3 |U_x(0)|$ in \eqref{stab:linf_WG1}, \eqref{stab:linf_WG1g} derives \eqref{stab:infty}.
}
\beq\label{stab:infty}
-(\tf17 \bar c_l  + 2 \bar c_{\om}^N) - 3 ( \mu_6 E_* +| u_x(\bar \e)(0)| )  \geq \lam , \quad
 -( \bar c_l  \al_{g, n} + 2 \bar c_{\om}^N ) - 3 ( \mu_6 E_* +| u_x(\bar \e)(0)| )  \geq \lam .
\eeq

\subsubsection{Weighted H\"older estimates}\label{sec:ineq_hol}

 The estimate has three components.  First, \eqref{hol:lin} records the linear
H\"older estimate, separating the nonlocal terms \(L_{\rm nloc,i}\) from the
local  terms \(L_{\rm loc,i}\).  Second, \eqref{hol:nlin} records the nonlinear
H\"lder estimates.  Third, \eqref{stab:hol}  combines these estimates into the final
pointwise stability condition, which corresponds to \eqref{eq:PDE_nondiag}.

For $x, z \in \R_2^{++}$, we assume $x_1 \leq z_1$ and denote $h = x-z$. To simplify the notations, we drop the dependence of $x, z$ in $\d$. We have the following linear H\"older estimate from Section \ref{sec:EE_hol}. 
\bseq\label{hol:lin}
\beq
\bal
& \mu_{h, i} g_i(h) \d(L_i) = ( d^{\mf{num}}_{i, L}(p_{x,z}) + d^{\mf{num}}_{g,i}) \cdot \mu_{h, i} g_i(h) \d( W_{1, i} \psi_i) +  \d_{ \mf{damp} ,i}(q_{x,z}) + \mu_{h,i} g_i(h) \d( B_{\mf{modi}, i}) , \\
& |\d_{ \mf{damp}, i}(q_{x,z})| \leq \mu_{h, i} | \d( d^{\numm}_{i, L}, x, z) |g_i(h) \f{\psi_i}{\vp_i}( q_{x,z} ) E_2 ,  \\
& |\mu_{h,i} g_i(h) \d( B_{ \mf{modi}, i} )| \leq  \mu_{h, i} g_i(h) | L_{\nlocc, i} | + | L_{loc, i} |,
 \\
& |L_{\nlocc, 1}| \leq   \d_{\sq}( U_A \psi_u \cdot \bar \om_x \f{\psi_1}{\psi_u}
+  V_A \psi_u \cdot \bar \om_y \f{\psi_1}{\psi_u}
, h ) ,  \\
& |L_{\nlocc, 2}| \leq \d_{\sq}( U_A \psi_u \cdot \bar \th_{xx} \f{\psi_2}{\psi_u}
+ V_A \psi_u \cdot \bar \th_{xy} \f{ \psi_2}{\psi_u}
+ U_{x,A} \psi_1 \cdot \bar \th_{x} \f{\psi_2}{\psi_1}
+ V_{x,A} \psi_1 \cdot \bar \th_{y} \f{\psi_2}{\psi_1}, h)  , \\
& |L_{\nlocc, 3}| \leq \d_{\sq}( U_A \psi_u \cdot \bar \th_{xy} \f{\psi_2}{\psi_u}
+ V_A \psi_u \cdot \bar \th_{yy} \f{ \psi_2}{\psi_u}
+ U_{y,A} \psi_1 \cdot \bar \th_{x} \f{\psi_2}{\psi_1}
- U_{x,A} \psi_1 \cdot \bar \th_{y} \f{\psi_2}{\psi_1}, h)  ,\\
 \eal
 \eeq
and the local terms satisfy
 \beq
 \bal
 & |L_{loc, 1} | \leq 
 \min\B(  \mu_{h, 1}  g_1(x-z)( \f{\psi_1}{ \vp_2}(x) + \f{\psi_1}{\vp_2}(z) ) , 
\min_{ (p, q) = (x,z),(z,x)} ( 
 \mu_{h,1} \d_{\sq}( \f{\psi_1}{ \psi_2} , h) \f{\psi_2}{\vp_2}(p) g_1(h)
+ \f{\mu_{h,1}}{\mu_{h, 2}} \f{g_1(h)}{g_2(h)} \f{\psi_1}{ \psi_2}(q) )
\B) E_2 \\
& | L_{loc, 2} | \leq  \B( \mu_{h, 2} g_2(h) \d_{\sq}(\bar v_x^N, h) \f{\psi_2}{\vp_3}(z) 
+ \f{\mu_{h, 2}}{\mu_{h, 3}} |\bar v_x^N(x)| \B) E_2 ,  \\
& | L_{loc, 3} | \leq \B(  \mu_{h, 3} g_3(h) \d_{\sq}(\bar u_y^N, h) \f{\psi_2}{\vp_2}(z)
+ \f{\mu_{h,3}}{\mu_{h,2 }} |\bar u_y^N(x)|  \B) E_2, 
\quad 
 \mu_{h} = \tau_1^{-1}(1,  \mu_1, \mu_2 ),
\eal
\eeq
\eseq
where the coupling weights $\mu_h, \tau_1, \mu_i$ are given in \eqref{energy2}, \eqref{wg:EE}, 
notation $\d_{\sq}$ that packages basic  H\"older estimates  is defined in \eqref{eq:hol_sq}, $d^{\mf{num}}_{g_i}$ is the damping term \eqref{eq:hol_nota} from the H\"older weight $g_i$ with $b(x)$ \eqref{eq:dp} replaced by $b^N  = \bar c_l x + \bar \uu^N$, 
\footnote{
$ d^{\mf{num}}_{g_i}$ is later replaced by the full nonlinear 
damping factor $d_{g_i ,\UU}$ in \eqref{stab:hol}. 
We only use it to illustrate the linear stability estimate and is not used in the proof.
}
and $L_{\nlocc ,i}$ estimates the nonlocal terms involving $\UU_A, (\na \UU)_A$ \eqref{eq:bad_lin}, and $L_{loc, i}$ estimates the local terms 
\beq\label{eq:hol_loc2}
\bal
& 
\mu_{h, 1} g_1(h) \d( ( \psi_1 / \psi_2) \cdot \eta_1 \psi_2 ), \quad 
   \mu_{h, 2} g_2(h) \d(\bar v_x^N \cdot \xi_1 \psi_2)    , \quad 
   \mu_{h, 3} g_3(h) \d( \bar u_y^N \cdot \eta_1 \psi_2 ) ,
\eal
\eeq
where $\psi_2 = \psi_3, g_2 = g_3$. For $L_{loc, 2}, L_{loc, 3}$, we use \eqref{eq:hol_loc} with $(f, i, j)= (\bar v_x^N, 3, 2)$ and 
$( \bar u_y^N, 2, 3)$. 
Note that we assume $x_1 \leq z_1$. For the damping terms, the choice of $(p_{x,z}, q_{x,z}) = (x, z)$ or $(z,x)$ depends on the locations of $x,z$, and we have two estimates of such terms. Instead of expanding the estimates again, we refer to Section ~\ref{sec:EE_hol_bad}. We have an improved estimate for the damping coefficients $\d_{\mf{damp}, i}(q_{x, z})$, which are explicit functions, in Section {\secbddphol} in the supplementary material I \cite{ChenHou2023aSupp}. See remark \ref{rem:dp_near0}. We optimize this improved estimate and the above estimate.

We have used $\psi_2 = \psi_3, g_2 = g_3, V_{y, A} = - U_{x, A}$ and $\mu_{h, \cdot}$ \eqref{energy4} to rewrite the parameters in the estimate in Section \ref{sec:EE_hol} equivalently. In \eqref{hol:lin}, for the terms $\d_{\sq}( U_A \psi_u \cdot f_1 + V_A \psi_u \cdot f_2 + ..., \, s)$ with some functions $f_1, f_2$, we bound it using \eqref{eq:hol_sq} and  
\[
\d_i( \tts{\sum}_m \, p_m q_m, x, z) 
 \leq \tts{\sum}_m \, \d_i(p_m q_m, x, z) \leq \tts{\sum}_m \, \d_i(p_m, q_m, x, z), \quad x_{3-i} = z_{3-i}.
 \]

For example, we have 
\[
\d_i( U_A \psi_u \cdot  \bar \om_x \f{\psi_1}{\psi_u } + V_A \psi_u \cdot  \bar \om_y \f{\psi_1}{\psi_u } ,x,z)
\leq \d_i( U_A \psi_u,  \bar \om_x \f{\psi_1}{\psi_u } ,x,z)
+  \d_i( V_A \psi_u,  \bar \om_y \f{\psi_1}{\psi_u } ,x,z)
\]
and then apply \eqref{eq:hol_sq} and \eqref{eq:hol_sq2} to bound $\d_{\sq}(\cdot, x, z)$. For each term, e.g. $U_A \psi_u, \bar \om_x \f{\psi_1}{\psi_u }$, we can obtain its piecewise $C_{x_i}^{1/2}$ and $L^{\inf}$ estimate.
It  simplifies the notations and estimates. We apply the same convention for other terms and the terms below.

The estimates below separate the nonlinear terms according to whether
they come from the transport term  $\cT_{d, N}$ \eqref{eq:tran_extra}, $W_1$-nonlinearity $N_{W_1, i}$ \eqref{def:N_W1}, or $\hat W_2$-nonlinearity $N_{\hat W_2, i}$ \eqref{eq:W2_M2}. In Sections \ref{sec:non}, \ref{sec:non_hol} (see also supplementary material I \cite[Section {\suppsecnonest}]{ChenHou2023aSupp}), we establish the nonlinear estimates for \eqref{eq:lin_U}, \eqref{eq:bad_lin} (see also \eqref{eq:hol_loc})
\bseq\label{hol:nlin}
\begin{align}\label{hol:nlin:a}
& \mu_{h, i} g_i(h) |\d( \cN\cF_i( \psi_i)) | \leq \cN\cT_i +    N_{W_1, i, \mf{hol}}  +   N_{\hat W_2, i, \mf{reg}, \mf{hol}}  +  \mu_{h, i} g_i(h)  N_{ \hat W_2, i, A, \mf{hol}} \notag \\
& \qquad \qquad \qquad \qquad \qquad + \mu_{h, i} g_i(h) \d_{\sq}( \bar \cF_{loc, i} \psi_i - \cR_{loc, i} \psi_i , h ),  
\end{align}
where $\cN\cT_i$ denotes the bound for $\cT_{ d, N}(\psi_i) W_{1, i} \psi_i$ via decomposition
\eqref{eq:non_tran_hol} into
$\cT_{uA}, \cT_{uR}, \cT_{c_{\om}}$ 
\begin{align}\label{hol:nlin:b}
  & \hspace{2in}  |\mu_{h, i} g_i(h) \d( \cT_{ d, N}(\psi_i) W_{1, i} \psi_i)|  \leq  \cN\cT_i ,\\
& \cN\cT_i \leq   \B\{ |\cT_{uA}(\psi_i)(x) + \cT_{uR}(\psi_i)(x)| 
+ \mu_{h, i} g_i( h) \d_{\sq}(\cT_{uA}(\psi_i) + \cT_{uR}(\psi_i), h)  \f{\psi_i}{\vp_i}(z)  
+  \max_{p=x,z} |\cT_{c_{\om}}(p)| \notag \\
& \qquad  \quad  +  \mu_{h, i} g_i(h)\min\B( | \d(  \cT_{c_{\om}}(\psi_i) , x, z )|  \min(|x|,|z|)^{\f{1}{2}} \max_{p=x,z} \f{\psi_i(p)}{\vp_i(p) |p|^{ \f{1}{2} }}, \d_{\sq}( \cT_{c_{\om}}(\psi_i) ,h ) \max_{p=x,z} \f{\psi_i(p)}{\vp_i(p) }   \B)  \B\} E_4 , \notag
\end{align}
and  $N_{W_1, i, \mf{hol}}$ denotes the weighted $C^{1/2}$ estimate for 
$N_{W_1, i}$ in \eqref{def:N_W1}:
\beq\label{hol:nlin:b2}
\bal
& |N_{W_1, 1, \mf{hol}} | \leq |U_x(0) | E_4 , \\
& | N_{W_1, 2, \mf{hol}} |\leq  \B( |U_x(0)|  + \mu_{h, 2} g_2(h) \d_{\sq}(\td U_x, h ) \f{ \psi_2}{\vp_2}(z)
+ |\td U_x(x)|  
+ \mu_{h, 2} g_2(h) \d_{\sq}(\td V_x, h ) \f{\psi_2}{\vp_3}(z)
+ \f{\mu_{h,2} }{\mu_{h, 3}} |\td V_x(x)|
 \B)  E_4, \\
 & | N_{W_1, 3, \mf{hol}} | \leq  \B( 3 |U_x(0)|  + \mu_{h, 3} g_3(h) \d_{\sq}(\td U_y, h ) \f{ \psi_2}{\vp_2}(z)
+ \f{\mu_{h, 3}}{ \mu_{h, 2}} |\td U_y(x)|  
+ \mu_{h, 3} g_3(h) \d_{\sq}(\td U_x, h ) \f{\psi_2}{\vp_3}(z)
+ |\td U_x(x)| \B)  E_4, \\
  \eal
 \eeq
Here, we  decompose $ N_{\hat W_2, i}$ via  \eqref{eq:N_W2:hol_dec};
  $N_{\hat W_2, i, \mf{reg}, \mf{hol}}$ 
 estimates $\cB_{op, i}( ( \UU_{\FM}, (\na \UU)_{\FM} ),  \hat W_2)$ (see \eqref{eq:Blin_gen})  and $U_x(0)
\hat W_{2, i, M}$ (see \eqref{eq:W2_M2});\footnote{
We use $\mf{reg}$ for regular, since the functions in $\cB_{op, i}( ( \UU_{\FM}, (\na \UU)_{\FM} ), \hat W_2)$ and 
$U_x(0)\hat W_{2, i, M}$ are more regular.
}
 $N_{\hat W_2, i, A, \mf{hol}} $ estimates $\cB_{op, i}( ( \UU_{A}, (\na \UU)_{A} ) , \hat W_2)$ (see \eqref{eq:Blin_gen}):
\begin{align}\label{hol:nlin:c0}
 N_{\hat W_2, i, A, \mf{hol}} & = 
  | \d( \psi_i \cB_{op, i}( ( \UU_A, (\na \UU)_A ), \wh W_2) ) |   , \\
 |N_{\hat W_2, i, \mf{reg}, \mf{hol}}| & \leq  \mu_{h, i} g_i(h) \B( |U_x(0)| \d_{\sq}( \psi_i \hat W_{2, i, M} , h)
 +   \d_{\sq}( \psi_i \cB_{op, i }( ( \UU_{\FM}, (\na \UU)_{\FM} ),  \wh W_2),  h   ) \B)   ,\notag
\end{align}
The factor $\mu_{h,i}g_i(h)$ is included in
$N_{\hat W_2,i,\mf{reg},\mf{hol}}$, \emph{but not} in
$N_{\hat W_2,i,A,\mf{hol}}$. We estimate $ N_{\hat W_2, i, A, \mf{hol}}$ as
\beq\label{hol:nlin:c}
\bal
& 
 N_{\hat W_2, 1, A, \mf{hol}} 
\leq \d_{\sq}( U_A \psi_u \cdot \f{ \psi_1 \pa_x \wh W_{2, 1}  }{ \psi_u }  + 
 V_A \psi_u \cdot \f{ \psi_1 \pa_y \wh W_{2, 1}  }{ \psi_u }   ,h ) , \\
 & 
  N_{\hat W_2, 2, A, \mf{hol}} 
\leq \d_{\sq}( U_A \psi_u \cdot \f{ \psi_2 \pa_x \wh W_{2, 2}   }{ \psi_u }  + 
 V_A \psi_u \cdot \f{ \psi_2 \pa_y \wh W_{2, 2}  }{ \psi_u } 
+ U_{x,A} \psi_1 \cdot \hat \eta_2 \f{\psi_2}{\psi_1}
+ V_{x, A} \psi_1 \cdot \hat \xi_2 \f{\psi_2}{\psi_1},h )  , \\
 & 
  N_{\hat W_2, 3, A, \mf{hol}} 
\leq \d_{\sq}( U_A \psi_u \cdot \f{ \psi_2 \pa_x \wh W_{2, 3}  }{ \psi_u }  + 
 V_A \psi_u \cdot \f{ \psi_2 \pa_y \wh W_{2, 3} }{ \psi_u } 
+ U_{y, A} \psi_1 \cdot \hat \eta_2 \f{\psi_2}{ \psi_1} 
 - U_{x, A} \psi_1 \cdot \hat \xi_2 \f{ \psi_2}{\psi_1},h ) , \\
\eal
\eeq
\eseq
where we have used $\psi_2 = \psi_3, g_2 = g_3$ (defined in \eqref{wg:hol}). The explicit function  $\d( \cT_{ c_{\om} }, x, z) g_i(h) \min( |x|, |z|)^{1/2}$ 
in \eqref{hol:nlin:b} is further estimated in Supplementary material I \cite[Section {\suppsecnondphol}]{ChenHou2023aSupp} for Paper I. See Section \ref{sec:non_hol} and Remark \ref{rem:dp_near0} for motivations.

For very large $ |x|_{\infty}, |z|_{\infty} \geq y_{730}$, 
due to much larger stability factor, we  derive much simpler estimate for $\cT_{d, N}(\psi_i) W_{1, i} \psi_i$ than those in \eqref{hol:nlin:b}. We use such an estimate to verify the nonlinear stability condition \eqref{stab:hol}. See Supplementary material I \cite[Section 8.7.4]{ChenHou2023aSupp}.

For $|x-z|$ not small, using the weighted $L^{\inf}(\vp_i)$ estimates of the linear terms $L_i(\vp_i)$  in \eqref{linf:WG1_lin} and nonlinear terms \eqref{linf:WG1_nlin}, and bounding 
$\f{\psi_i}{\vp_i}$, we have a simple $L^{\inf}$ estimate 
\beq\label{eq:check_hol_lg}
\bal
 &
 |\mu_{h,i} g_i(x-z) \d(B_{\mf{modi},i} \psi_i)| \leq  \mu_{h, i} g_i(x-z) 
 \B(\f{\psi_i}{\vp_i} | L_i(\vp_i) |(x) + \f{\psi_i}{\vp_i} | L_i(\vp_i)|(z) \B), 
  \\
   & |\d( \cN\cF_i (\psi_i) )|  \leq  \sum_{y = x, z}  \f{\psi_i}{\vp_i}(y) ( \cT_\mw{wg}(\UU, \psi_i)  E_4+  | N_{W_1, i} | \cdot  \vp_i  +  | N_{\hat W_2, i} | \cdot \vp_i  +  |\bar \cF_{loc, i} \vp_i| + |\cR_{loc, i} \vp_i |  )(y) ,\\
   & |\cT_\mw{wg}(\UU, \psi_i) | \leq |\cT_{ c_{\om}}| + |\cT_{uA} |+ |\cT_{uR}|  ,\quad 
| W_{ 1,i} \psi_i | \leq   \psi_i/\vp_i  E_4,
 \eal
\eeq
where we \emph{replace} the transport part $\cT_\mw{wg}(\UU, \vp_i)$  in \eqref{linf:WG1_nlin} by $\cT_\mw{wg}( \UU, \psi_i)$ since we use weight $\psi_i$ \eqref{eq:non_dec1}. For the local terms \eqref{eq:hol_loc2}, we optimize the $C^{1/2}$ estimate of $L_{loc, i}$ in \eqref{hol:lin}, and the $L^{\inf}$ estimate in \eqref{linf:WG1_lin} with weight $ \f{\psi_i}{\vp_i}$ similar to the above.

Combining  the above estimates yields the stability conditions for the weighted H\"older estimate 
\beq\label{stab:hol}
- (d^{\mf{num}}_{i, L}(p_{x,z}) + d_{g_i, \UU} ) E_* - 
\mu_{h, i} g_i(h) \B( | \d( d^{\numm}_{i, L}, x, z) |  \f{\psi_i}{\vp_i}( q_{x,z} )
+| L_{\nlocc, i} | +    | \d( \cN\cF_i (\psi_i)  ) |  \B) - | L_{loc, i}|  \geq \lam  E_* \, ,
\eeq
for some $\lam > 0$, uniformly for any $x, z \in \R_2^{++}, x_1 \leq z_1 $, where $(p_{x,z}, q_{x,z}) = (x, z)$ or $(z,x)$ depends on the locations of $x,z$.  See Section \ref{sec:EE_hol_bad} and remark \ref{rem:dp_near0}, and Section {\secbddphol}. Here, 
$d_{g_i, \UU}$ is the nonlinear damping factor \eqref{eq:hol_nota},\eqref{eq:dp} by adding $\UU$
\[
  d_{g_i, \UU} \teq   ( b_U(x) - b_U(z) ) \cdot (\na g_i) (x-z) g_i^{-1}(x-z) , 
   \quad b_U = \bar c_l x + \bar \uu^N + \UU .
\]
 To verify the above inequalities, we follow Section \ref{sec:EE_hol_sum}.

\subsubsection{ODEs for $c_{\om}$ and $\om_{xy}(0), \th_{xxy}(0)$}\label{sec:ineq_cw}

In this subsection, we generalize the ODE estimates for $c_\omega$-terms and $\om_{xy}(0), \th_{xxy}(0)$ in Section~\ref{sec:rank1} to the modified equations \eqref{eq:lin_U}, and record the corresponding stability inequalities.
Recall the estimate \eqref{eq:w_est1} and $f_*$ \eqref{eq:green}
\[
\bal
& c_{\om}( p ) =  -\f{4}{\pi} \int_{ \R_2^{++}}  f_*(y) p(y) dy 
= -\f{4}{\pi} \la p, f_* \ra , \ 
f_*(y) = \f{y_1 y_2}{|y|^4},   \   \vp_{M, i} = \max( \vp_i, \mu_{g, i} \vp_{g, i} ), \   |W_{1,i} |  \leq \vp_{M, i}^{-1} E_4 .
\eal
\]

For \eqref{eq:lin_main_cw} and $q = 1, \chi_{\mf{ode}}$ \eqref{eq:cutoff_near0_all}, following Section \ref{sec:cw_w1}, we have the following linear estimates
\[
\bal
 |\la \G_{1, M} , q f_* \ra| 
& \leq \la   |\na \cdot ( ( \bar c_l x + \bar \uu^N)  f_* q  ) |,  \vp_{M, 1}^{-1} \ra  E_4 +\la  | U_A \bar \om_x  + V_A \bar \om_y |,  f_* q \ra , \\
   |\la \G_{2, M} , q f_* \ra| 
& \leq\la    |\na \cdot ( ( \bar c_l x + \bar \uu^N)  f_* q   ) - \bar u_x^N f_* q  | , \vp_{M, 2}^{-1}  \ra E_4 + 
\la  | U_A \bar \th_{xx}  + V_A \bar \th_{xy} |,  f_* q \ra  \\
& \quad +   \la| U_{x, A} | |\bar \th_x|  + |V_{x, A}| |\bar \th_y| ,  f_* q\ra  
+ \la |\bar v_x^N | , \vp_{M, 3}^{-1} f_* q \ra E_4 ,
\eal
\]
and nonlinear estimates 
\[
\bal
  N( c_{\om}( \om_1 q)) & \teq
\f{4}{\pi} | \la \td \cN_1 + \bar \cF_{loc, 1} -  \cR_{loc, 1} , f_* q \ra |
   \leq \g_{ 1} |U_x(0)| E_4   +  \f{4}{\pi} \B\la | \bar \cF_{loc, 1}| +| \cR_{loc, 1}| +
   \uds{ \eqref{eq:W2_M2} }{ | N_{ \hat W_2,1 } | } ,  f_* q \B\ra  \\
 & \quad + \f{4}{\pi} \B\la  |U  \pa_x( f_* q) + V \pa_y(f_* q) |, \vp_{M, 1}^{-1} \B\ra E_4,   \\
  N( c_{\om}( \eta_1 q )) &\teq \f{4}{\pi} |\la \td \cN_2 +\bar  \cF_{loc, 2} - \cR_{loc, 2} , f_* q \ra|
   \leq   \f{4}{\pi} \B\la | \bar \cF_{loc, 2}| +| \cR_{loc, 2}| + 
\uds{ \eqref{eq:W2_M2} }{ | N_{ \hat W_2, 2 }  | }
+  |\td U_x| \vp_{M, 2}^{-1} E_4 + |V_x| \vp_{M, 3}^{-1} E_4,  f_* q  \B\ra  \\
& \quad +\g_{2} |U_x(0)| E_4   +  \f{4}{\pi} \B\la  |U \pa_x( f_* q ) +  V \pa_y(f_* q ) |, \vp_{M, 2}^{-1} \B\ra E_4 ,
\eal
\]
where $ \g = ( \mu_{51}, \mu_{5 2})$ for $q = \chi_{\mf{ode}}$ \eqref{eq:cutoff_near0_all}, and $\g = (\mu_5, \mu_{62} )$ for $q =1$. 
We use $V_x = \td V_x$ since $\td V = V - V_y(0) y$ \eqref{eq:u_tilde}. We can use the weighted $L^{\inf}(\vp_i)$ estimate for $ \bar  \cF, \cR, N_{\hat W_2}$ from \eqref{linf:WG1_nlin}. 
For $q = 1$, $q f_*$ is singular near $0$. For $\uu = \UU, \bar \uu^N$, 
using a direct calculation for $f_*$,  we decompose
\[
u f_{*, x} + v f_{*, y } =\f12 (\f{ u}{x} + \f{v}{y}) (x f_{*, x} + y f_{*, y})
+ \f12 (\f{ u}{x} - \f{v}{y}) (x f_{*, x} - y f_{*, y})
=  - (\f{ u}{x} + \f{v}{y}) f_* -   (\f{ u}{x} - \f{v}{y}) \f{2 x y(x^2 - y^2)}{ |(x, y)|^6}
\]

In the first term, we exploit the cancellation near the origin:
\[
u/x + v / y = \td u /x + \td v / y, \quad \td u = u - u_x(0) x, \ \td v = v + u_x(0) y, 
\]
which vanishes $O(|x|)$ near $0$. Then we apply the piecewise bounds of $\UU, \bar \uu^N$ to estimate the integrals. 
Moreover, from \eqref{eq:ODE_IBP}, when \(q\equiv 1\), or in the region where
\(q=\chi_{\rm ode}=1\), we use
\[
\na \cdot ( x f_*(x)) = 0, \quad 
\na \cdot  ( (\bar c_l x + \bar \uu^N) q f_*  ) = \bar \uu^N \cdot \na f_*, \quad  q =1.
\]

Recall $\chi_{\mf{ode}}$ from \eqref{eq:cutoff_near0_all}, the damping terms in the ODEs $\bar c^N_{\om} \la \om_1, f_* q\ra, 2 \bar c^N_{\om} \la \eta_1, f_* q \ra$ \eqref{eq:lin_main_cw} and Section \ref{sec:rank1}. The stability conditions for $c_{\om}(\om_1 \chi_{\mf{ode}} ), c_{\om}( \eta_1 \chi_{\mf{ode}} ), c_{\om}(\eta_1)$ read
\beq\label{stab:cw1}
\bal
 & -  \bar c_{\om}^N E_* -  \mu_{51}^{-1} \big( \f{4}{\pi}  | \la \G_{1, M} , \chi_{\mf{ode}} f_* \ra| + N( c_{\om}(\om_1 \chi_{\mf{ode}} )) + \mu_{52} E_*  \big)  > \lam E_*, \\
 & - 2 \bar c_{\om}^N E_* -  \mu_{52}^{-1} \big( \f{4}{\pi}  |\la \G_{2, M} , \chi_{\mf{ode}} f_* \ra| + N(c_{\om}( \eta_1 \chi_{\mf{ode}} )  ) \big) > \lam E_*,  \\
 & -2 \bar c_{\om}^N E_* - \mu_{62}^{-1}  \big( \f{4}{\pi}  |\la \G_{2, M} , f_* \ra| + N(c_{\om}( \eta_1)) \big) > \lam E_*  ,
 \eal
\eeq
for some $\lam >0$. For the estimate of $c_{\om}( \om)$ \eqref{eq:ode_cw_all}, using the estimate in Section \ref{sec:rank1} and the above estimates, we obtain the stability conditions, 
\[
 - \bar \lam_{c_{\om}}  E_*- 
 \mu_{6}^{-1} \B( 
 \f{4}{\pi}  |\la \G_{1, M} ,  f_* \ra| + N( c_{\om}(\om_1 )) + \mu_{62} E_*
      + N_{\hat W_2, \mf{ode}} \B) > \lam E_*,
\]
for some $\lam >0$, where $N_{\hat W_2, \mf{ode}}$ bounds 
\[
\bal
N_{\hat W_2, \mf{ode}} & \leq 
 \sum_{i\geq 2}   | c_{\om}(\hat F_i(0)) a_i(W_1, \hat W_2)(t) |  + \sum_{i \geq 1} 
 \B( | a_i( W_1(t- T_i), \hat W_2(t- T_i) )  \cdot c_{\om}( \hat F_i(T_i)) |  \\
& \qquad  \qquad \qquad \qquad + \int_0^{t\wedge T_i} | a_i(W_1(t-s), \hat W_2(t-s))  | \cdot | \pa_s c_{\om} (\hat F_i(s)) - \bar \lam_{c_{\om}} c_{\om}(\hat F_i(s)) | ds \B) ,
\eal
\]
and we have used $\one_{t \geq T_i} \leq 1$ \eqref{eq:ode_cw_all}. Under the bootstrap assumptions, we bound $b_i = a_i(W_1, \hat W_2) $ by $c_i E_*$ for some constant $c_i$ \eqref{eq:W2_bi}. See Section \ref{sec:rank1}. For linear modes, $a_i$ only depends on $W_1$.

For $\om_{xy}(0), \th_{xxy}(0)$ in \eqref{eq:ODE_lin_xy}, since $|c_{\om}(\om)| \leq \mu_6 E_4, 
|\om_{xy}(0) |  \leq \mu_8 E_4, | \th_{xxy}(0) | \leq \mu_7 E_4$ \eqref{energy4},  the stability conditions read
\[
\bal
  (2 \bar c_l - \bar c_{\om} ) E_* - \mu_8^{-1} ( (\mu_7 + \mu_6 |\bar \om_{xy}(0)| ) E_*
+ \mu_6 \mu_8 E_*^2 + | \pa_{xy}\bar \cF_{1}(0)|  ) & > \lam E_* , \\
  (3\bar c_l /2 - \bar c_{\om} ) E_* - \mu_7^{-1} 
 ( \mu_6 | \bar \th_{xxy}(0)| E_* + \mu_6 \mu_7 E_*^2 + |\pa_{xy} \bar \cF_2(0)|  ) & > \lam E_* ,
\eal
\]
for some $\lam>0$, where we used $ 2 \bar c_l - 2 \bar c_{\om} + \bar u_x(0)
= \f{3}{2} \bar c_l - \bar c_{\om} $ \eqref{eq:normal}.  We verify the stronger condition 
\eqref{eq:W2_non_boot2}
\[
\bal
 &   \mu_8 \mu_6 E_* + |\pa_{xy} \bar \cF_1(0)|
 < 5 \mu_6 E_* , \quad  \mu_7 \mu_6 E_* + |\pa_{xy} \bar \cF_2(0)|
 < 10 \mu_6 E_* . 
 \eal
\]

To obtain \eqref{W_bound}, under the bootstrap assumption, we verify 
\beq\label{W2_bound}
  | \wh W_{2, i}|< 100 E_4 < 100 E_*.
\eeq

\vs{0.1in}
\paragraph{\bf{Plots of the nonlinear weighted $L^{\inf}$ estimate}}

\begin{figure}[ht]
   \centering
      \includegraphics[width = \textwidth  ]{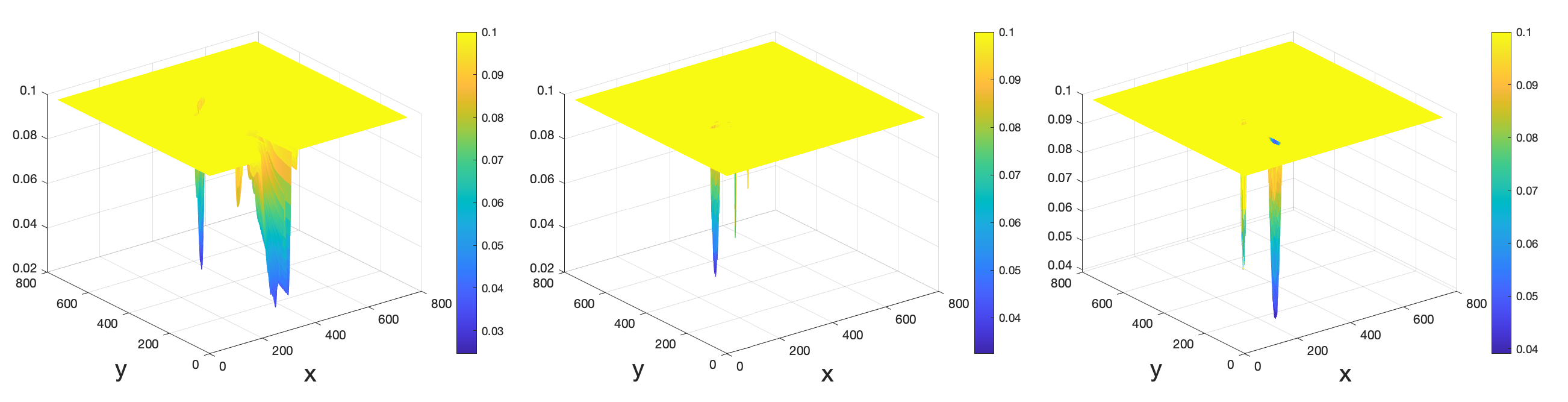}
        \captionsetup{width=0.95\textwidth} 
      \caption{ Verify nonlinear $L^{\inf}(\vp_i)$ stability conditions. 
      The strict positivity of these values verifies the condition \eqref{stab:linf_WG1} for some $\lam >0$, and for any $(x, y) \in D$.      
       }
            \label{fig:stab_WG1}
 \end{figure}

In Figure \ref{fig:stab_WG1}, we plot the rigorous piecewise lower bounds of $\min( LHS_i / E_*, 0.1)$ in a region covering $D = [0, 10^{15}]^2$, where $LHS_i$ denotes the left-hand side (LHS) of \eqref{stab:linf_WG1} in the $i$-th equation. We normalize $LHS_i$ by $1/E_*$ and take the minimum with a threshold to highlight the region with small linear damping factors. 

\begin{figure}[ht]
   \centering
      \includegraphics[width = \textwidth  ]{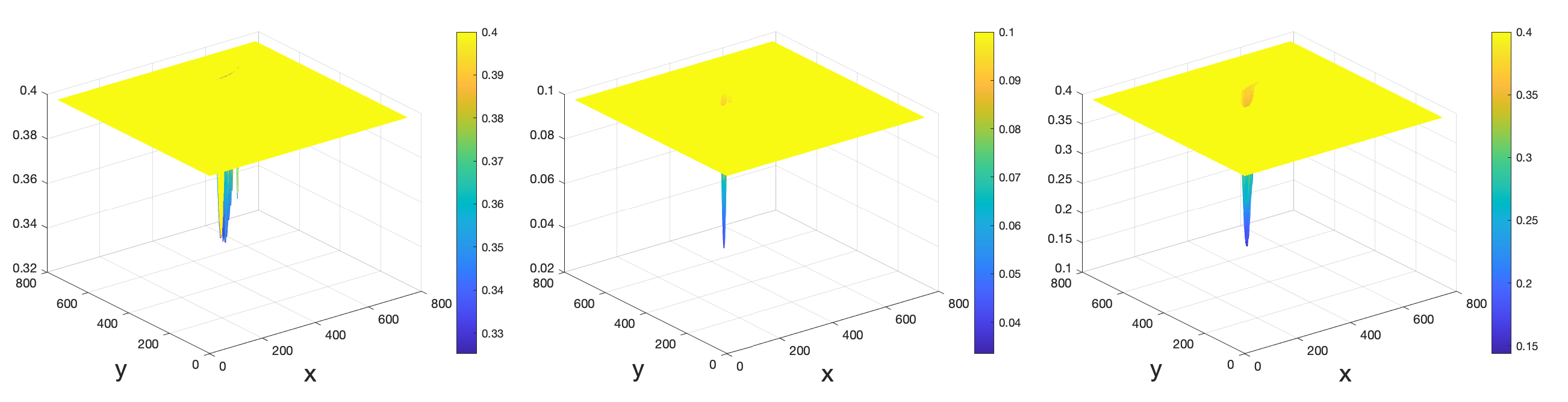}
      \captionsetup{width=0.95\textwidth} 
      \caption{ Verify nonlinear $L^{\inf}(\vp_{g,i})$ stability conditions. The strict positivity of these values verifies the condition \eqref{stab:linf_WG1g} for some $\lam >0$, and for any $(x, y) \in D$. 
       }
            \label{fig:stab_WG1g}
 \end{figure}

In Figure \ref{fig:stab_WG1g}, we plot $\min( LHS_i / E_*, c_i)$ with $c = (0.4, 0.1, 0.4)$ with $LHS$ being the left-hand side  of \eqref{stab:linf_WG1g}. All of these bounds are positive. We only use $7$ approximation terms for the velocity in $[0, 200]^2$ away from the boundary (see \eqref{eq:u_appr_2nd}, \eqref{eq:vel_cutoff_bd}). Thus the stability factor is weaker for $x$ not very large and in the bulk. 
\footnote{
Adding more velocity approximation modes would increase this stability margin in the bulk region. 
See also remark below \eqref{eq:EE_target1}. Since the current lower bounds are already strictly positive, we do not need such an improvement.
}

The weighted $L^{\inf}(\vp_4)$ estimate has a much larger stability margin $ \geq 1.5 E_*$, and thus we do not plot it. Beyond $D$, we have much larger damping factors and use the $L^{\inf}$ estimate in Section {\suppseclinffar} in the supplementary material I \cite{ChenHou2023aSupp}.

We cannot visualize the H\"older estimate condition \eqref{stab:hol}. This condition \eqref{stab:hol} depends on two spatial points and is verified by interval arithmetic in the code \cite{ChenHou2023code}.

\subsection{Estimates for the H\"older weights}\label{sec:hol_wg_est}

We record elementary bounds for the H\"older weights $g_i$ in  \eqref{wg:hol}. These weights $g = g_{\ell}$ are $(-1/2)$--homogeneous. In our energy estimates, we estimate several $0$--homogeneous functions related to $g$ for $h_1, h_2 \geq 0$:
\[
f(h) = h_k^{1/2} g(h),\  h_k \f{ (\pa_i g) }{g}(h),  \ |h| \f{ (\pa_i g) }{g} (h), \ k , i = 1,2,
\quad \mbox{or} \quad \f{ g_{i_1}(h)}{g_{i_2}(h)}, 1\leq i_1, i_2 \leq 3.
\]
Since $f(h) = f( \f{h_1}{h_2} , 1)$ for $h_2 \neq 0$ and $f(h ) = f(1, \f{ h_2}{h_1}), h_1 \neq 0$, we estimate it by partitioning $(h_1, h_2) \in [0, 1] \times \{1 \} , \{1 \} \times [0, 1] $ and 
using the monotonicity of $g, \pa_i g$. For $g = g_{\ell}$ \eqref{wg:hol}, we get 
\[
g(s) = \f{1}{ A_1(s)^{1/2} + A_2(s)^{1/2}}, \  \pa_i g = - \f{1}{2} \f{  a_{1i} A_1^{-1/2} + a_{2i} A_2^{-1/2}   }{ (A_1^{1/2} + A_2^{1/2} )^2 }, \  \f{\pa_i g}{g}
  = - \f{1}{2} \f{  a_{1i} A_1^{-1/2} + a_{2i} A_2^{-1/2}   }{ A_1^{1/2} + A_2^{1/2}  } ,
\]
where $A_m = a_{m 1} s_1 + a_{ m 2} s_2$ for some $a_{m j} >0$. Clearly, $g$ is decreasing in $|s_i|$
for $i=1,2$. For $s_1, s_2 > 0$, since $A_m$ is increasing in $s_1, s_2$, 
$\pa_i g, \f{\pa_i g}{g}$ are negative and increasing in $s_1, s_2$. It follows that $| \f{\pa_i g}{g}|$ is decreasing in $s_1, s_2$.

\vspace{0.2in}
{\bf Acknowledgments.} The research was in part supported by NSF Grants DMS-1907977 and DMS-2205590. We would like to acknowledge the generous support from Mr. K. C. Choi through the Choi Family Gift Fund and the Choi Family Postdoc Gift Fund. We would also like to thank Drs. Pengfei Liu and De Huang for a number of stimulating discussions in the early stage of this project, and Dr. Tarek Elgindi for pointing out that the Boussinesq result almost readily extends to 3D Euler. JC is grateful to Dr. Xiaoqi Chen for several suggestions on coding and the use of High Performance Computing. Part of the computation in this paper was performed using the Caltech IMSS High Performance Computing. The support from its staff is greatly appreciated.

\bibliographystyle{plain}
\bibliography{selfsimilar}

\end{document}


\begin{abstract}

In Part I \cite{ChenHou2023a_supp}, we established linear energy estimates,
nonlinear $L^{\infty}$ energy estimates, and estimates for representative
terms in the nonlinear H\"older energy estimate. In this supplementary
material, we provide additional details on the nonlinear H\"older estimate,
focusing on piecewise bounds for 
several explicit functions and explicit damping functions near the origin,
where low regularity and singular weights make the analysis delicate. We also
perform the energy estimate in the very far field using the asymptotic behavior
of the estimates.

\end{abstract}

 \maketitle


 \setcounter{section}{7}

\section{Estimates of nonlinear and error terms and piecewise bounds}

Firstly, we recall some notations from \cite{ChenHou2023a_supp,ChenHou2023b_supp}. Recall from \cite{ChenHou2023a_supp} the following variables 
\beq\label{eq:EE_W1}
W_{1, 1} = \om_1, \ W_{1, 2} = \eta_1, \  W_{1, 3} = \xi_1,
\quad \wh W_{2, 1} = \hat \om_2, \  \wh W_{2, 2} = \hat \eta_2,  \  \wh W_{2, 3} = \hat \xi_2, \quad \bar W  = (\bar \om, \bar \th_x, \bar \th_y),
\eeq
and the full energy $E_4$ for $W_1, W$ 
from \cite[Section 5.6]{ChenHou2023a_supp}, which satisfies 
\footnote{
We do not write the full $E_4$ since some norms are not used in the supplementary material.
}
\bseq\label{energy4}
\beq
\bal   
E_4(t) \geq \max\B( & \max_i || W_{1, i} \vp_i||_{\inf} , \sqrt{2} \tau_1^{-1} || \om_1 |x_1|^{-1/2} \psi_1||_{\inf},    \max_i \mu_{h, i} [ W_{1, i} \psi_i ]_{C_{g_i}^{1/2}}   ,   \\
 & \max_i  \mu_{g, i} || W_{1, i}  \vp_{g, i} ||_{\inf} , \mu_{6}^{-1} | c_{\om}(\om)|    \B) , 
 \\ 
\eal
\eeq
with 
\beq
  \mu_{h} = \tau_1^{-1}(1, \mu_1, \mu_2), \quad \mu_g = \tau_2(\mu_4, 1, 1), \quad  \psi_2 = \psi_3, \quad g_2(h) = g_3(h) ,
\eeq
\eseq
where the parameters $\tau_i, \mu_{ij}$ are given in \eqref{wg:EE}, 
and the weights $\psi, \vp$ are chosen in  \cite[Appendix {\appparawg}]{ChenHou2023a_supp}. We use $\psi$ to denote radial weights for $C^{1/2}$ estimates, and $\vp$ to denote mixed weights 
for $L^{\infty}$ estimates.
 Under the bootstrap assumption, for $\mu_5$ is given in \eqref{wg:EE}, we get 
\beq\label{eq:cw_w1_est}
  |c_{\om}(\om_1)| < \mu_5 E_4. 
\eeq

Denote the mesh for computing the approximate steady state (see Appendix {\appsolurep} in Part II \cite{ChenHou2023b_supp}) and several domains (\textbf{note that $y_i$ corresponding to $x_1(i+1)$ in MatLab indexing})
\beq\label{eq:mesh_xy}
\bal
& 0 =y_0 < y_1< y_2< .. < y_{N-1}, \quad Q_{ij}  =  [y_i, y_{i+1}] \times [y_j, y_{j+1}], \quad N = 748, \quad  N_1 = 706,  \\ 
&  N_2 = 730 < N,  \quad   R_1 = 5000,   \quad   \Om_{near} \teq [0, y_{N_1}]^2 ,
\quad y_{N-1} > 10^{15}.
\eal
\eeq

We decompose the weighted $L^{\inf}$ estimate into the region $[0, y_{N-1}]^2$ and the far-field $( [0, y_{N-1}]^2)^c$.

We decompose the H\"older estimate with $x, z$ into  five parts. In case (a.1)-(a.3), we assume that $x \in Q_{i, j}, z \in Q_{p, q}$.

{\bf (a.1)} $x$ or $z$ in the near-field and the bulk $\Om_{near}$ with $|x- z| \leq R_1$. 
Since $\Om_{near}$ is the union of $Q_{ij}, \max(i, j) \leq N_1-1$,
we check the stability inequality for 
\[
  \max( \min(i, p),  \min(j,  q) ) \leq  N_1-1, \quad  |i-p|, |j-q| \leq n_{R_1} ,
\]
where $n_{R_1}$ satisfies $y_{n_{R_1}} > R_1$. For $\max( |i-p|, |j- q| ) > n_{R_1}$, since the spacing of mesh $y_{i+1} - y_i$ is increasing, we get 
\[
|x-z| \geq \mathrm{dis}( Q_{ij}, Q_{pq}) \geq y_{n_{R_1}} > R_1.
\]

{\bf (a.2)} $x$ or $z$ in $[0, y_{N - 2}]^2 \bsh \Om_{near}$ with $|i-p|$ and $|j-q| \leq 1$. We check the stability for $|i-p|$ and $|j-q| \leq 1$ similar to the case of (a.1). 

{\bf (a.3)} $x, z \in [0, y_{N-1}]^2 \bsh [0, y_{N_1-2}]^2$ with $|i-p|$ or $|j-q| \geq 2$. In this case, 
$x$ and $z$ are at least one grid away, and thus $|x-z| \geq \min(y_{i+1}-y_i)=y_2 -y_1 > 0$. 
While the distance is not large,
since  the forcing term 
in $[0, y_{N - 2}]^2 \bsh \Om_{near}$ is very small, 
we perform the H\"older estimate using the weighted $L^{\inf}$ estimate for the forcing term. 
See Section \ref{sec:hol_case_a3}.

{\bf (b)} $x, z$ in the far-field with $|z|_{\inf} , |x|_{\inf} \geq y_{N_2-1}> 10^{12}$ and any $|x-z|\geq0$. We use the asymptotic  estimates. See Section \ref{sec:hol_far}.

{\bf (c)} 
$x \in [0, y_{N_2-1}]^2$ and $z \in \R^2_{++}$ with $|x- z| \geq R_1$. In this case, $x$ and $z$ are far apart, the H\"older estimate follows from the weighted $L^{\inf}$ estimate and triangle inequality. 
See Section \ref{sec:hol_lg}.

\vs{0.1in}

Since $y_{i+1} - y_i  > R_1$ for $i \geq 600$, 
by symmetry, cases (a.1), (a.2), (a.3), (c) cover the estimate for $x$ or $z$ in $[0, y_{N_2-1}]^2$. The remaining case $x, z \notin [0, y_{N_2-1}]^2$ is covered by case (b).

The most important case for the H\"older estimate is case (a.1) and has been established in  \cite[Section {\secEE}]{ChenHou2023a_supp}, and we have outlined the method for the verification. Estimates in cases (a.2),(a.3),(b), (c) are relatively easy since the bad term in the estimates are very small due to large $|x-z|$, or the decay of the coefficients of the nonlocal terms in the far-field. Similarly, the $L^{\inf}$ stability estimate in the far-field is simple since we have large damping factors.

\vs{0.1in}
\paragraph{\bf{Organization}}

The remaining sections are organized as follows. In Section \ref{sec:non_err}, we recall the modified 
decompositions of linear and nonlinear terms from Section {\seccombvelerr} in \cite{ChenHou2023a_supp}. In Sections \ref{sec:hol_basic}-\ref{sec:non_est}, we estimate the piecewise bounds for the damping terms near $0$ and estimate the modified nonlinear terms.
In Sections \ref{sec:linf_far}, \ref{sec:hol_far}, we discuss the energy estimate for $x$ outside the computational domain, e.g. $|x| > 10^{12}$, using the asymptotics of the coefficients in the energy estimates. In Section \ref{sec:hol_lg}, we discuss the $C^{1/2}$ estimate with large $|x-z|$ related to case (c) and estimates related to case (a.3).

\subsection{Decompositions of the nonlinear and nonlocal error terms}\label{sec:non_err}

In this section, we recall the nonlocal error terms and the modified linearization equations from \cite[Section {\seccombvelerr}]{ChenHou2023a_supp}. Then we will estimate the nonlinear terms with more details in Section \ref{sec:non_est}. 

\subsubsection{Decomposition and estimates of the velocity}\label{sec:non_dec_u}

We define
\[
\uu = (u, v ), \quad 
\td u = u - u_x(0) x, \quad \td v = v + u_x(0) y.
\]
Recall that $\uu_A(f)$ is the velocity after subtracting the approximation  $\td {\hat \uu}$ defined in Section {\secEE} \cite{ChenHou2023a_supp}
\[
\uu_A(f) \teq \td \uu(f) -  \td{ \hat{\uu} }(f) = \uu - \hat \uu , 
\quad \uu_{x, A}  \teq  \td \uu_x(f) - \td{ \hat{\uu}}_x(f), \quad 
\uu_{y, A}  \teq \td \uu_y(f) - \td{ \hat{\uu}}_y(f) \, .
\]

Firstly, we recall the following decomposition of error 
and velocity from Section {\seccombvelerr} in \cite{ChenHou2023a_supp}.  Denote by $\bar \phi^N, \hat \phi^N$ the stream function constructed numerically for the Poisson equation $-\D \bar \phi = \bar \om, -\D \hat \phi = \hat \om_2$, $\bar \uu^N, \hat \uu^N$ the associated velocity, and $\bar \e, \hat \e$ the errors
\beq\label{eq:wb_vel_err}
\bal
 \bar \uu^N = \uu( -\D \bar \phi^N) = \na^{\perp} \bar \phi^N, 
 \quad \bar \e & = \bar \om - (-\D) \bar \phi^N, 
 \quad  \ \uu(\bar \om) = \bar \uu^N + \uu(\bar \e),  \\
  \hat \e & = \hat \om_2 - (-\D) \hat \phi^N, 
\quad \uu(\hat \om_2) = \uu( (-\D) \hat \phi^N ) + \uu(\hat \e).
\eal
\eeq

The error  $ \e$ only vanishes $O(|x|^2)$ and we perform a correction near $0$. For  $\e = \bar \e$ or $\hat \e$ and some cutoff function $\chi_{ \e} = 1 +O(|x|^4)$ near $0$, we decompose 
\beq\label{eq:uerr_dec1}
\bal
  &  \e  =  \e_1 +  \e_2, \quad  \e_2 =  \e_{xy}(0) \D( \f16  x y^3   \chi_{  \e}) , \quad \uu(  \e_2) = \na^{\perp}(-\D)^{-1}  \e_2 = - \f{1}{6}   \e_{xy}(0)\na^{\perp}  ( x y^3 \chi_{  \e}) ,  \\
  & \uu(  \e)   = \uu(  \e_1) + \uu(  \e_2)  = \uu_A(  \e_1) +  \hat \uu(  \e_1 ) + \uu(  \e_2) , \quad u_x(\e)(0) = u_x(\e_1)(0), \\
 \eal
\eeq
where  
\[
 u_x(\e_2)(0) = \tf16 \e_{xy}(0) \cdot \pa_{x y} (x y^3 \chi_{\e}) |_{x, y=0} =0 .
\]
We perform a similar decomposition for $\na \uu(  \e)$. We choose $\chi_{ \e}$ in the above form such that $ \e_2 =  \e_{xy}(0) xy + O(|x|^4)$ and we can obtain $\uu( \chi_{  \e})$ explicitly. We choose $\chi_{\bar \e}$ for $\bar \e$ and $\chi_{\hat \e}$ for $\hat \e$ in
\eqref{eq:cutoff_near0_all} and they have \emph{different} parameters.

Let $\om$ be the perturbation without decomposition. We combine these errors and perturbations and perform the following decompositions 
\beq\label{eq:U}
\bal
\UU  & = (U, V) = \uu( \om + \bar \e), \qquad  \qquad \quad  \ \UU_A = \uu_A(\om_1 + \bar \e_1 + \hat \e_1) , \\
  \UU & = \UU_A + U_x(0)(x, - y) + \UU_{\FM}, 
  \qquad 
\td \UU \teq \UU - U_x(0)(x, - y) = \UU_A +\UU_{\FM},   
 \\ 
 \UU_{\FM} & = \td { \hat \uu}(\om_1 + \bar \e_1 + \hat \e_1) + \uu( \bar \e_2 + \hat \e_2) 
+ \td \uu( (- \D) \hat \phi^N ) . 
\eal
\eeq
Similarly, we decompose $ \na \UU$ and define $ (\na \UU)_A, (\na \UU)_{\FM}$.

\vs{0.1in}
\paragraph{\bf{Estimates of  velocity}}

We use  methods in \cite[Section {\secvelcomp}]{ChenHou2023b_supp} to establish the same type of weighted $L^{\inf}, C^{1/2}$ estimates for $\uu_A( \hat \e_1), \uu_A(\bar \e_1)$ as those for $u_A(\om_1)$, i.e. the weighted $L^{\inf}$ estimate for $ \uu_A, \na \uu_A$, and the H\"older estimate for $ \psi_u \uu_A, \psi_1 (\na \uu)_A$. 
For the perturbation $\uu_A(\om_1)$, we use the norm $|| \om_1 \vp_1||_{\inf}, || \om_1 \vp_{g, 1}||_{\inf}, [\om_1  \psi_1]_{ C_{x_i}^{1/2}}$ and estimate various quantities involving the velocity, e.g. 
\[
  |u_{ij, A}(x)| \leq C_{ij, 1}(x, \rho) || \om_1 \rho ||_{\inf}
  +  C_{ij, 2}(x, \rho) [ \om_1 \psi ]_{ C_x^{1/2}} + C_{ij, 3}(x, \rho) [ \om_1 \psi ]_{ C_y^{1/2}} , \quad \rho = \vp_1,  \vp_{g, 1}.
\]

For the error $\uu_A(f), f = \bar \e_1, \hat \e_1$, we use the norm $|| f \vp_{elli}||_{\inf}, [f \psi_1]_{C_{x_i}^{1/2}}$. 
We then combine the estimate of the nonlocal error term and the perturbation 
to bound $\UU_A, (\na \UU)_A$ in \eqref{eq:U} in the energy estimate in Section \ref{sec:non_dec}. 

The terms $\td \uu( (-\D) \hat \phi^N ) = \na^{\perp} \hat \phi^N +\hat \phi^N_{xy}(0)(x, -y)$ vanishes $O(|x|^3)$ near $0$, and we estimate its $C^3$ bounds following Section \ref{sec:W2_est}. 
The terms $ I_1 = \td {\hat \uu}(\om_1 + \bar \e_1 + \hat \e_1) ,  I_2 = \uu( \bar \e_2 + \hat \e_2)$ consist of finite rank operators on $\om_1, \bar \e_1, \bar \e_2, \hat \e_1, \hat \e_2  $ with smooth coefficients. Moreover, by definition, 
\[
I_1, \quad  I_2, \quad  \td \uu( (-\D) \hat \phi^N ) = O(|x|^3), \quad (\na \UU)_{\FM} = O(|x|^2)
\]
near $x=0$. The regularity of these terms are established 
\cite[Section 4.4]{ChenHou2023a_supp}.
We establish their piecewise $C^3$ estimate near the origin and $C^1$ estimate away from the origin, and combine the estimates of these terms together. See nonlinear estimates in Section {\secrankone} in \cite{ChenHou2023a_supp}. This allows us to estimate the $C^3$ bounds of $\UU_{\FM}$ and $\UU_{\FM} = O(|x|^3)$. Similarly, we estimate $(\na \UU)_{\FM}$. We factor out the term $u_x(0) (x, -y), U_x(0)(x, -y)$ in the above decomposition since the constant in our estimates for such a term is larger than others, and control it using the energy $E_4$.
Using these $C^k$ estimates of $\UU_{\FM}$, we can further obtain its $C^{1/2}$ estimate.

We can further estimate the $C^{1/2}$ norm of $ \UU, \na \UU, \td \UU, \na \td \UU $. For example, we use 
\[
 \na \td \UU = (\na \UU)_A + (\na \UU)_{\FM} , \quad 
\td U_x = (U_{x, A} \psi_1) \psi_1^{-1} +  U_{x, \FM},
\]
and the $C^{1/2}$ estimates of $U_{x, A} \psi_1, \psi_1^{-1}, U_{x, \FM}$ to obtain the $C^{1/2}$ estimate of $\td U_x$. Similarly, using the $C^1$ estimates of $\psi_1^{-1}, \psi_u^{-1}, \UU_{\FM}$ and the $C^{1/2}$ estimate of $ \UU_A \psi_u, \UU_{x, A} \psi_1 $, we can obtain the $C^{1/2}$ estimate of $\td \UU, \wt{ \na \UU}$. See Section \ref{sec:hol_basic} for basic H\"older estimates of a  product.

\subsubsection{Estimate of local and nonlocal terms}\label{sec:W2_est}

Recall from Sections {\seclinevo}, {\secvelerr} in Part II \cite{ChenHou2023b_supp} 
(see also Part I \cite[Section 4.2]{ChenHou2023a_supp} ) that we represent $ \hat W_2, \hat \phi^N $ using Duhamel's principle and piecewise fifth order polynomials and an analytic basis. 
Since the basis functions for  $\hat W_2$ are supported in the mesh domain $[0, y_{N_2-1}]^2$, 
where $N_2$ is chosen in \eqref{eq:mesh_xy}, we obtain
\beq\label{eq:W2_supp}
 \hat W_2(x, y) \equiv 0, \quad \forall  \, \max(x, y) \geq y_{N_2-1},
\eeq

 For variables that depend locally on $ \hat W_2$ \eqref{eq:EE_W1} and $\hat \phi^N$, under the bootstrap assumption $E_4(t) < E_*$, we can estimate piecewise $C^k, k \leq 4$ bounds of $ 
 \na^{\perp} \hat \phi^N, \hat W_2, \hat W_2 - (-\D) \hat \phi^N $ following Section {\secWtwo} \cite{ChenHou2023a_supp} using the energy $E_4$, e.g.
\[
 |\na^k \hat W_{2, i}(x) | \leq C_{k,  i}(x) E_4. 
\]

\subsubsection{Decomposition of the nonlinear and error terms}\label{sec:non_dec}

As in \cite[Section 5.8]{ChenHou2023a_supp}, we introduce the bilinear operator $\cB_{op,i}$ for $(\uu, M) \in \R^2 \times \R^{2 \times 2}, G = (G_1, G_2, G_3) $
\beq\label{eq:Blin_gen}
\bal
\cB_{op, 1}( (\uu, M), G )
  & = - \uu \cdot \na G_1 + M_{11}(0) G_1,  \\
\quad \cB_{op, 2}( (\uu, M), G ) & = - \uu \cdot \na G_2 +2 M_{11}(0) G_2 - M_{11} G_2 - M_{21} G_3, \\
\cB_{op, 3}( (\uu, M), G ) & = - \uu \cdot \na G_3 + 2 M_{11}(0) G_3 - M_{12} G_2 - M_{22} G_3.
\eal
\eeq
The operator is bilinear in the first input $(\uu, M)$ and the second input $G$. 
The input $G$ relates to $(\om, \eta, \xi)$ and $M$ relates to velocity gradient $\na \uu$ and its variants. 
If $M=  \na \uu, M_{11} = u_x, M_{12} = u_y, M_{21} = v_x, M_{22} = v_y$, then we drop $M$ to simplify the notation 
\beq\label{eq:Blin}
\bal
\cB_{op, 1}(\uu, G )  & = - \uu \cdot \na G_1 + u_x(0) G_1,  \\
\  \cB_{op, 2}(\uu, G) & = - \uu \cdot \na G_2 +2 u_x(0) G_2 - u_x G_2 - v_x G_3, \\
\cB_{op, 3}(\uu, G) & = - \uu \cdot \na G_3 + 2 u_x(0) G_3 - u_y G_2 - v_y G_3.
\eal
\eeq

We apply the operator $\cB_{op, i}$ to $\uu, \td \uu, \UU_A, \UU_{\FM}$.  We abuse the notation to write 
\beq\label{eq:B_abuse}
\bal
& \cB_{op, i}( \uu_A, G) = \cB_{op, i}( (\uu_A, (\na \uu)_A ), G), \quad  
 \cB_{op, i}( \UU_{\FM}, G) = \cB_{op, i}( (\UU_{\FM}, (\na \UU)_{\FM} ), G) . 
 \\
\eal
\eeq
 In $\cB_{op, i}( \uu_A, G)$ and $ \cB_{op, i}( \UU_{\FM}, G)$, we replace $\na \uu$ in \eqref{eq:Blin} by $ (\na \uu)_A, ( \na \UU)_{\FM}$. Note that $\uu_A$ does not satisfy the differential relation among $\uu_A ,(\na \uu)_A$.

Recall $W_1$ from \eqref{eq:EE_W1}. We recall the following weighted nonlinear equations 
with weights $\rho_i$ from  \cite[Section {\secnon}]{ChenHou2023a_supp} 
(the sign of $\cR_{loc, i}$ in \eqref{eq:lin_U} is \emph{minus} $-$, not plus)
\beq\label{eq:lin_U}
 \pa_t( W_{1, i} \rho_i) + (\bar c_l x + \bar \uu^N + \UU) \cdot \na ( W_{1, i} \rho_i)
  =  d^{\mf{num}}_{i, L}(\rho_i)  W_{1, i} \rho_i + B_{\mf{modi},i} \rho_i + \td \cN_i(\rho_i) + (\bar \cF_{loc, i} - \cR_{loc, i}) \rho_i ,
\eeq
where $(\bar \om, \bar \th)$ is the approximate steady state, $\bar \uu^N$ is defined 
in \eqref{eq:wb_vel_err}, $d^{\mf{num}}_{i, L}(\rho_i)$ are the damping coefficients defined in \eqref{eq:dp}, 
and $\bar c_l, \bar c_{\om}$ are determined by 
\beq\label{eq:normal}
\bar c_l = 2 \f{\bar \th_{xx}(0) }{\bar \om_x(0) }, \quad \bar c^N_{\om} = \f{1}{2} \bar c_l + \bar u^N_x(0) .
\eeq
Here $B_{\mf{modi}, i}$ denotes the linear bad terms
\beq\label{eq:bad_lin}
\bal
 B_{\mf{modi}, 1}(x) & \teq \eta_1  -  \UU_{A} \cdot \na \bar \om ,  \\
B_{\mf{modi}, 2}(x) & \teq -  \bar v^N_x \xi_1 -  \UU_A \cdot \na \bar \th_x -   \UU_{x, A} \cdot \na \bar \th , \\
 B_{\mf{modi}, 3}(x)  & \teq -  \bar u^N_y \eta_1  - \UU_A \cdot \na \bar \th_y -  \UU_{y, A} \cdot \na \bar \th,
\eal
\eeq
and $\td \cN_i$ denotes the modified nonlinear terms in the estimate of $W_{1, i}\rho_i$
\beq\label{eq:non_dec1}
\bal 
\td \cN_1(\rho_1) = &  (\UU \cdot \na \rho_1) \cdot   \om_1  + U_x(0) \om_1 \rho_1 \qquad \qquad \qquad \quad  + \cB_{op, 1} ( \td \UU, \wh W_2) \rho_1 + U_x(0) \hat W_{2, 1, M}
, 
 \\
\td \cN_2(\rho_2) = &  ( \UU \cdot \na \rho_2) \cdot \eta_1  + ( U_x(0) \eta_1 
- \td U_x \eta_1 -\td  V_x \xi_1 ) \rho_2 +  \cB_{op,2}( \td \UU, \wh W_2) \rho_2
+  U_x(0) \hat W_{2, 2, M}
, \\
 \td \cN_3(\rho_3) = &  (\UU \cdot \na \rho_3) \cdot \xi_1  + ( 3 U_x(0) \xi_1 
- \td U_y \eta_1 - \td V_y \xi_1 ) \rho_3 + \cB_{op,3}( \td \UU,\wh W_2) \rho_3 
+  U_x(0) \hat W_{2, 3, M}
\eal
\eeq
where $\cB_{op, i}$ is given in \eqref{eq:Blin}, $\hat W_{2, i, M}$ is given by
\beq\label{eq:W2_M2}
\bal
 \hat W_{2,  \cdot, M } &  =(\hat \om_2, \hat \eta_2, \hat \xi_2),  \hspace{1.3in}  \hat \om_{2, M}  \teq \hat \om_2 - x \pa_x \hat \om_2 + y \pa_y \hat \om_2 - \hat \om_{2, xy}(0) f_{\chi, 1} , \\
   \hat \eta_{2, M}  \teq & \hat \eta_2 - x \pa_x \hat \eta_2 + y \pa_y  \hat \eta_2 -  \hat \eta_{2, xy}(0) f_{\chi, 2} ,  \quad   \hat \xi_{2, M}  \teq 3 \hat \xi_2- x \pa_x \hat \xi_2 + y \pa_y \hat \xi_2 
- \hat \xi_{2, xx}(0) f_{\chi, 3},
\eal
\eeq
with $f_{\chi,i}$ defined in \eqref{eq:cutoff_near0_all}. The terms $\bar \cF_{loc, i}, \cR_{loc, i}$ are the essential local part of the residual error and residual operators. 
Since we have estimated $\bar \cF_{loc, i}, \cR_{loc, i}$ in Part I \cite{ChenHou2023a_supp} and Part II \cite{ChenHou2023b_supp}, we do not present their formulas and refer them to Section {\seccombvelerr} in \cite{ChenHou2023a_supp}.

The weighted nonlinearity $\td \cN_i(\rho)$ \eqref{eq:non_dec1} can be decomposed as 
\beq\label{def:non_wg}
\td \cN_i(\rho) = \cT_{d, N}(\rho) \cdot W_{1,i} \rho 
+ N_{W_1, i} \cdot \rho + N_{\hat W_2, i} \cdot \rho.
\eeq
where $\cT_{d, N}(\rho) $ is defined in \eqref{eq:dp}, $N_{W_1, i}$ 
denotes other nonlinear terms involving $W_1, \UU$ 
\beq\label{def:N_W1}
N_{W_1, 1}  \teq U_x(0) \om_1 , \quad 
N_{W_1, 2}  \teq  
    U_x(0) \eta_1 - \td U_x \eta_1 - \td V_x \xi_1,  \quad
N_{W_1, 3} 
 \teq 3 U_x(0) \xi_1 - \td U_y \eta_1 - \td V_y \xi_1 ,
\eeq
and $ N_{\hat W_2, i}$ denotes the last two nonlinear terms in RHS of \eqref{eq:non_dec1} of 
$\hat W_2, \UU$:
\beq\label{def:N_W2}
\bal
 N_{ \hat W_2, i} & \teq   \cB_{op, i}( \td \UU, \wh W_2) + U_x(0) \hat W_{2,  i, M}  \\
 & = \cB_{op, i}( ( \UU_A, (\na \UU)_A ) , \wh W_2)  +  \cB_{op, i}( ( \UU_{\FM}, (\na \UU)_{\FM} ), \wh W_2) 
 +  U_x(0) \hat W_{2,  i, M} .
 \eal
\eeq
Here, we have used that $\cB_{op, i}$ is bilinear 
and $\td \UU = \UU_A + \UU_{\FM}$ by \eqref{eq:U}.

We introduce notations for the damping terms $  d_{i, L}^{\numm}$  and the damping factor $d_{g_i, \UU}$ in the H\"older estimate  derived in Section {\secEEintro} \cite{ChenHou2023a_supp} 
\begin{align}\label{eq:dp}
 b_U(x) & = \bar c_l x + \bar \uu^N + \UU , \quad \quad \qquad \quad   d_{g_i, \UU} \teq   \f{ (b_U(x) - b_U(z)) \cdot (\na g_i) (x-z) }{g_i(x-z)}  ,  \\
  \cT_d^{\numm}(\rho) & = \rho^{-1}(  (\bar c_l x + \bar \uu^N) \cdot \na \rho),  
 \quad \cT_{d, N}(\rho) = \rho^{-1} (\UU \cdot  \na \rho) ,
 \quad  d_{i, L}^{\numm}(\rho)  \teq \cT^{\numm}_{d}(\rho) + d_{loc, i} , \notag \\  
 d_{loc, 1} & = \bar c^N_{\om},  \hspace{1.2in} \ \, \quad  d_{loc, 2} = 2 \bar c_{\om}^N - \bar u^N_x,  \qquad  \quad  \    d_{loc, 3} = 2 \bar c^N_{\om} + \bar u^N_x . \notag 
\end{align}

The subscripts $d, L$ are short for \textrm{damping, linear}.

\subsection{Basic H\"older estimates}\label{sec:hol_basic}

We follow \cite[Section \secholbasic]{ChenHou2023a_supp} to package several 
basic H\"older estimates. We use the following notations for the H\"older estimate
\beq\label{eq:hol_dif}
\bal
& \d_i(f, x, z) \teq \f{ |f( x ) - f(z)|}{|x-z|^{1/2}} , \quad  
 \D_i(f, x, z) \teq \f{ f(z) -f(x)}{ z_i - x_i},  \quad z_i \neq x_i, \ z_{3-i} = x_{3-i},  \\
&  \d(f)(x, z) \teq f(x)- f(z) . \\
\eal
\eeq
By definition, we have $\d_i(f, x, z) = \d_i(f, z, x), \D_i(f, x, z) = \D_i(f,z,x)$. By abusing the notations of $\d_i$, we denote by $\d_i(f, g, x, z)$ a basic $C^{1/2}$ estimate for product 
\beq\label{eq:hol_prod0}
\d_i(f, g, x, z) = \min_{(p, q) = (x,z), (z,x)} \d_i(f, x, z) |g(p)|  + \d_i(g, x, z) |f(q)  |, \quad  x_{3-i} = z_{3-i} .
\eeq
Clearly, if $g = 1$, we get $\d_i(f, g, x, z) = \d_i(f, x, z)$. We have the following basic estimates
\beq\label{eq:hol_prod}
\bal
\d(f g)(x, z) & = \d(f)(x, z) g(p) + \d(g)(x, z) f(q) , \\
\d_i(f g , x, z) & \leq \d_i(f, g, x, z), \quad x_{3-i} = z_{3-i} .
\eal
\eeq

Given the piecewise $C_x^{1/2}, C_y^{1/2}$ estimates of $f$, we use the following method for the piecewise H\"older estimate of $f$ with H\"older weight $g$ and two points $(x, z)$ 
\beq\label{eq:hol_path}
\bal
 & | (f(z) - f(x) )  g(h)| =  | ( f(z) - f(w) + f(w) - f(x) ) g(h)| \\ 
 \leq & ( \d_1( f, w, x ) |h_1|^{1/2} + \d_2(f, z, w ) |h_2|^{1/2} ) g(h) , \quad 
 h = z- x, \quad w =     (z_1, x_2) .
 \eal
\eeq

\begin{figure}[t]
\centering
\begin{subfigure}{.49\textwidth}
  \centering
  \includegraphics[width=0.45\linewidth]{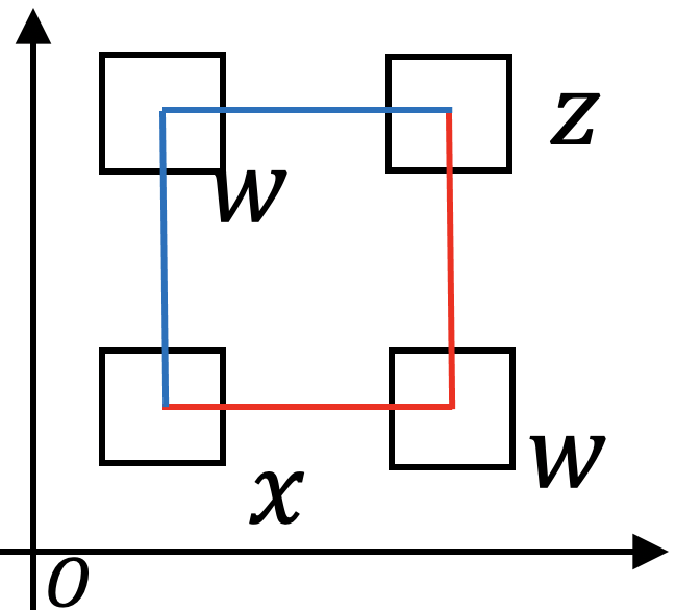}
    \end{subfigure}
    \begin{subfigure}{.49\textwidth}
  \centering
  \includegraphics[width=0.45\linewidth]{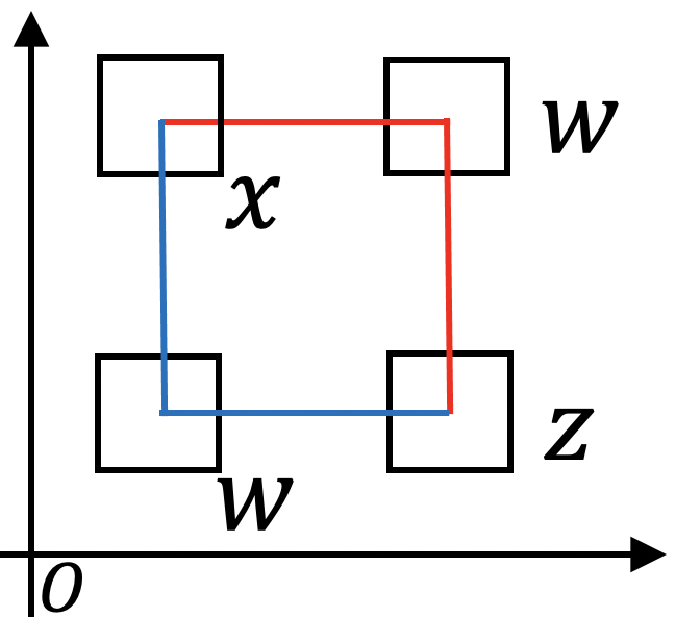}
    \end{subfigure}
\caption{Left, right figures correspond to the locations of $(x, z)$ in cases (1) $s \geq 0$, (2) $s< 0$, $s =(z_1-x_1)(z_2 - x_2)$. The red line and blue line represents two choices of the $C_{x}^{k}, C_y^{k}, k =\f{1}{2}$ or $1$ estimates used to estimate $f(x)-f(z)$. 
 }
\label{fig:hol_pa}
\vspace{-0.15in}
\end{figure}

The function $|h_i|^{1/2} g(h)$ is $0$-homogeneous and we apply the method in Section \ref{sec:hol_wg_est} to estimate it. For $\td w = (x_1, z_2) $, we derive another estimate. We optimize two estimates for $|\d(f)(x, z)| g(x-z)$. In Figure \ref{fig:hol_pa}, we illustrate the locations of $x, z$ and the $C_{x_i}^{1/2}$ ($(\d_i(f, p, q))$) estimates and the triangle inequality used to estimate $\d(f,x, z) g(h)$. We introduce $\d_{\sq}$ to denote this estimate and similar estimate for the product 
\beq\label{eq:hol_sq}
\bal
\d_{\sq}(f, x, z, s) 
& \teq \min( \d_1(f, x, (z_1, x_2) )  |s_1|^{1/2} + \d_2(f, (z_1, x_2), z ) |s_2|^{1/2},  \\
& \quad \d_2(f, x, (x_1, z_2) )  | s_2 |^{1/2} + \d_1(f, (x_1, z_2), z ) | s_1|^{1/2} ) , \\
\d_{\sq}(f, g, x, z, s) 
& \teq \min( \d_1(f, g, x, (z_1, x_2) ) | s_1 |^{1/2} + \d_2(f, g, (z_1, x_2), z ) |s_2|^{1/2},  \\
& \quad  \d_2(f, g, x, (x_1, z_2) )  | s_2|^{1/2} + \d_1(f, g, (x_1, z_2), z ) |s_1|^{1/2} ) , \\
\d_{\sq}(f, x, z, s) & \teq \d_{\sq}(f, 1, x, z, s), \\
\eal
\eeq
where $\d_i(f, g, x, z)$ is defined in \eqref{eq:hol_prod0}. We use the notation $\sq$ since it mimics Figure \ref{fig:hol_pa} and indicates that we use $C_x^{1/2}, C_y^{1/2}$ estimates to obtain the $C^{1/2}$ estimate. By definition, we get 
\beq\label{eq:hol_sq2}
|\d(f,x,z ) | \leq \d_{\sq}(f, x, z, x-z), \quad | \d(f g, x, z) | \leq \d_{\sq}(f, g, x, z, x-z).
\eeq

Note that  $\d_i(f,x,z), \d_i(f,g, x,z),  \d_{\sq}(f, x,z, s), \d_{\sq}(f, g, x, z, s)$ are symmetric in $x, z$. We introduce an extra variable $s$ to reduce bounding $\d_{\sq}(f,x,z,x-z) g(x-z)$ to estimating $\d_i(f, x, z)$ and $|x_j-z_j|^{1/2} g_i(x-z)$ separately, and $\d_{\sq}(f,x,z, s) g_i(s)$ is $0$-homogeneous in $s$.

\vs{0.1in}
\paragraph{\bf{H\"older estimate of $W_{1,i}$ terms}}

For $x_1 \leq z_1$ and $f W_{1, i} \psi_i$ with $f \in C^{1/2}$,
using the energy $||W_{1, i} \vp_i||_{\inf} , \mu_{h, i} [  W_{1, i} \psi_i]_{C_{g_i}^{1/2}} \leq E_4 $ \eqref{energy4}, \eqref{eq:hol_prod}, \eqref{eq:hol_sq2}, we perform its $C^{1/2}$ estimate as follows 
\beq\label{eq:hol_loc}
\bal
&  |\mu_{h, j} g_j(h) \d( f  W_{1, i}  \psi_i, x, z)|
 \leq  \mu_{h, j} g_j(h) ( | \d( f) W_{1, i} \psi_i(z)|
+ |  f(x) \d( W_{1, i} \psi_i)| ) \\
& \leq \f{\mu_{h, j} g_j(h)}{ \mu_{h, i} g_i(h)} | f(x) | E_4 
+ \mu_{h, j} g_j(h) \d_{\sq}( f , x, z, h) \f{\psi_i}{\vp_i}(z) E_4 , \quad h = x-z.
\eal
\eeq
If $i = j$, the first term reduces to $|f(x)| E_4$. We only pick one decomposition in \eqref{eq:hol_prod} with the coefficient of $\d( W_{1, i} \psi_i)$ evaluating at $x$, i.e. $f(x),$ to simplify the estimates. Note that $x_1 \leq z_1$.

\subsection{Additional estimates of weights}\label{sec:add_wg}

We need several asymptotic estimates of the weight and 
 the ratio among weights. 
The reader who is not interested in the rigorous estimates of explicit weights can skip to Section \ref{sec:non_est} for the estimates of nonlinear terms and error terms.

\subsubsection{Asymptotics of the weights}\label{sec:lin_asym_wg}

 Consider a radial weights 
$
\rho(x) = \rho( r) = \sum_{i \leq n} p_i r^{a_i}, 
$
with $p_i > 0$ and $a_i$ is increasing with $i$. Firstly, we have 
\[
\bal
& \pa_{i} r^{\al} =  \al \f{x_i}{r} r^{\al - 1}, \quad \pa_i \pa_j r^{\al} = \al \pa_j( x_i r^{\al-2}) =(  \al(\al-2) \f{x_i}{r} \f{x_j}{r} + \al \d_{ij}) r^{\al-2}
\teq C_{i,j}(\b, \al) r^{\al-2},  \\
\eal
\]
 where $\b = \atan (x_2 / x_1)$, and $C_{i,j}(\b, \al)$ is obtained by simplifying $ x_1 / r = \cos \b, x_2 /r= \sin \b$. For $r \geq r_*$, since $r^{a_i - a_n} \leq r_*^{a_i -a_n}$, it is easy to obtain 
\begin{align}\label{eq:lin_asym_wg1}
  & p_n r^{a_n} \leq \rho(x)  \leq r^{a_n}  (\sum_{ \ell \leq n} p_\ell r_*^{a_\ell - a_n}) ,  \
 |\pa_{x_i} \rho(x) | = |\sum_{ \ell \leq n} p_\ell a_\ell x_i r^{a_\ell - 2}|
\leq \f{x_i}{r}  r^{a_n-1} (\sum_{\ell \leq n} |p_\ell a_\ell | r_*^{a_\ell - a_n} ), \notag \\
 & |\pa_1^i \pa_2^j \rho(r)| \leq 
 r^{a_n-2} \sum\nolim_{ \ell \leq n} |p_\ell | |C_{i,j}(\b,  a_\ell )| 
 r_*^{a_\ell - a_n} \teq \rho_{ij, u}(\b) r^{a_n - 2} , \ i+ j = 2 .
\end{align}
We can bound the angular function $\rho_{ij, u}(\b)$ by partitioning $\b \in [0, \pi/2]$. 
In particular, we obtain 
\beq\label{eq:lin_asym_wg2}
  |\pa_{x_1}^i \pa_{x_2}^j \rho | \leq   \rho_{ij, u}(\b) r^{a_n - i - j} \cos \b^i \sin \b^j , \quad i+ j \leq 1 ,
\eeq
for  some constants $\rho_{ij, l}, \rho_{ij, u}$ depending on $r_*, a_i, p_i$.

For two radial weights $f, g$ and $\al$, applying the above estimates and 
\beq\label{eq:lin_asym_wg3}
  |\pa_{i}  ( \f{f}{g}  r^{\al} )| \leq |\al \f{x_i}{r} r^{\al- 1} \f{f}{g}|
  +r^{\al} ( | \f{ \pa_i f} {g} | + |\f{f \pa_i g }{g^2} |),
\eeq
we can derive the asymptotics of $ r^{\al} \f{f}{g}$ and its derivatives. Note that we can extract the common angular factor $\f{x_i}{r}$ in the above estimate. Similarly, using the triangle inequality \eqref{eq:lin_asym_wg2}-\eqref{eq:lin_asym_wg1}, we can obtain the asymptotics of
\beq\label{eq:lin_asym_wg4}
  | \pa_i \rho / \rho| \leq C(\rho, i)(\b) r^{-1} , \quad  |\pa_j ( \pa_i \rho / \rho)|
  \leq C(\rho, i, j)(\b) r^{-2}.
\eeq

\subsubsection{Piecewise bounds of the ratio of two weights}\label{sec:wg_rat}

We have estimated the piecewise bounds for the radial weights $\rho(r)$ and the mixed weight in Appendix {\appwgradial}, {\appwgmix} in Part II \cite{ChenHou2023b_supp}. Denote 
\beq\label{eq:wg_nota}
 \vp = x^{-1/2} P(r) + Q(r), \quad  P(r) = \tts{ \sum_{i\leq n_p} } p_i r^{a_i}, \   Q(r) = \tts{ \sum_{i \leq n_q} }  q_i r^{b_i}, \quad \rho(r) = \tts{ \sum_{i \leq n_s} } s_i r^{c_i} ,
\eeq
with $p_i, q_i, s_i >0$ and  $a_i, b_i, c_i$ are increasing with $i$. In all of our choices of mixed weights, $P(r)$ is more singular  than $Q$ near $r=0$ and decays faster than $Q$ for large $R$, i.e. $a_1 < b_1, a_{n_p} < b_{n_q}$.
In the energy estimates, we need several piecewise bounds of the ratio of these weights in $[x_l, x_u] \times [y_l, y_u]$. Since these weights are singular near $(0, 0)$ and along $x=0$, we need to factorize out these singularities. Due to symmetry, we focus on $x \geq 0$ and use $r = (x^2 + y^2)^{1/2}$. Consider the following decomposition for a radial weight and the mixed weight 
\beq\label{eq:wg_mix_dec1}
\bal
& P_m = \tts{ \sum_{i\leq n_p} } p_i r^{a_i - a_1},  \quad  P = r^{a_1} P_m ,  \quad  \vp = x^{-1/2} \vp_{m, h } = x^{-1/2} r^{a_1} \vp_m,  \\
& \vp_{m, h} \teq P  + x^{1/2} Q, \qquad \vp_m = r^{-a_1} \vp_{m, h} 
= P_m + r^{-a_1}  x^{1/2} Q(r).
\eal
\eeq
where \emph{m} is short for \emph{modified}. Since $r^{-a_1} Q(r)$ is not singular 
by our assumption $a_1 < b_1$, $\vp_m $ is not singular, and $\vp_{m, h}$ is only singular along $x=0$. Using the piecewise bounds for the radial weights $P, P_m, r^{-a_1} Q(r)$, we can obtain the piecewise bounds for the above modified weights. 

We introduce $D_r = r \pa_r$. Using $D_r \rho_m = \sum s_l (c_l-c_1) r^{c_l-c_1}$, we get
\beq\label{eq:wg_rad_dec1}
x_i \pa_i \rho = x_i \sum s_l c_l x_i r^{c_l-2}
 = \f{x_i^2}{r^2} \sum s_l c_l r^{c_l} , \
   \f{x_i \pa_i \rho}{ \rho} = \f{x_i^2}{r^2} (c_1 + \f{\sum s_l (c_l - c_1) r^{c_l - c_1}  }{ \sum s_l r^{c_l - c_1}}  ) 
 = \f{x_i^2}{r^2} ( c_1 + \f{ D_r \rho_m}{ \rho_m}).
\eeq
Since $c_i$ is increasing with $i$, the leading power in $D_r \rho_m$ is $r^{c_2 - c_1}$, which vanishes near $r =0$. 
For all radial weights we use, $c_2 - c_1 > 1$ and thus $x_i \pa_i \rho / \rho \in C^{1}$ locally, and we can estimate its piecewise derivatives. We will use \eqref{eq:wg_rad_dec1} in Section \ref{sec:bd_dp_hol}.

For $(x, y) \in [x_l, x_u] \times [y_l, y_u]\subset \R^2_{++}$, we have the following simple estimates 
\beq\label{eq:cos_sin}
 \f{x_l^2}{x_l^2 + y_u^2} \leq \f{x^2}{r^2} 
 \leq \f{x_u^2}{ x_u^2 + y_l^2} , 
 \quad  \f{y_l^2}{x_u^2 + y_l^2} \leq \f{y^2}{r^2} 
 \leq \f{y_u^2}{ y_u^2 + x_l^2} .
\eeq

\paragraph{\bf{Ratio between mixed weights}}
Let $\vp_i = x^{-1/2} P_i(r) + Q_i(r) $ be two mixed weights with  the leading power $r^{\al_i}$ for $ P_i(r) $. For $\al_1 \geq \al_2$,  we use \eqref{eq:wg_mix_dec1}, the decompositions
\[
 \f{\vp_1}{\vp_2 } = \f{ \vp_{1, m, h}}{ \vp_{2, m, h} } = r^{\al_1 - \al_2} \f{\vp_{1, m}}{ \vp_{2, m}},
\]
and piecewise bounds of each part to obtain three estimates and then optimize them. We use the second identity to overcome the singularity on $x=0$ and the third near $r =0$.

Similarly, for two radial weights with leading power $r^{\al_i}$, we estimate $P_1 / P_2$ using 
the following decompositions \eqref{eq:wg_mix_dec1}
\[
P_1 /P_2 = r^{\al_1 - \al_2} P_{1, m} / P_{2, m} .
\]

\vs{0.1in}
\paragraph{\bf{Ratio between the radial and the mixed weight}}
Let $\rho, \vp$ be a radial weight and a mixed weight \eqref{eq:wg_nota}. We use the following decomposition to estimate $\rho / \vp, \rho / (\vp x^{1/2}), \rho / (\vp r^{1/2})$ 
\beq\label{eq:wg_rat_C1C0}
 \f{\rho}{ \vp} = r^{ c_1 - a_1 } \f{ \rho_m x^{1/2}}{ \vp_{m}} , 
 \quad  \f{\rho}{ \vp x^{1/2} }
  = \f{\rho}{ \vp_{m , h}} = r^{c_1 - a_1} \f{ \rho_m}{\vp_m}, 
  \quad \f{\rho / r^{1/2}}{ \vp }
  = r^{c_1 - a_1} \f{\rho_m }{ \vp_m} ( \f{x}{r})^{1/2}.
\eeq
The factor $x / r$ is estimated using \eqref{eq:cos_sin}.

We remark that to obtain sharp piecewise bounds of the ratio, we can evaluate the ratio on fine grids, and then use the derivative bounds of the ratio and apply the estimate \eqref{est_2nd}. 

\subsubsection{Ratio among radial weights}\label{sec:wg_rat_rad}

In the  H\"older estimate in Section \ref{sec:non_hol}, we need to estimate $\f{ \pa_i P}{ P Q}$ for radial weights $P, Q$ \eqref{eq:wg_nota}. Using the piecewise bounds of $P, P^{-1}, Q^{-1}$, \eqref{eq:lei}, and following Appendix {\appholfunc} in Part II \cite{ChenHou2023b_supp}, we can obtain piecewise $C^1, C^{1/2}$ estimates of $\f{ \pa_i P}{ P Q}$. We need more careful estimates near $r= 0$. 
Consider the decomposition \eqref{eq:wg_nota},\eqref{eq:wg_mix_dec1} for $P, Q$. We get 
\beq\label{eq:dP_PQ}
\bal
\f{\pa_i P}{ P Q}
& = ( \f{\pa_i P_m}{ P_m} + \f{\pa_i r^{a_1}}{r^{a_1}} )  \f{r^{-b_1}}{Q_m}
= x_i  \f{  \sum_{2 \leq j \leq n_p} (a_j - a_1) p_j r^{a_j - a_1-2 -b_1}  }{ P_m Q_m} + \f{  a_1 x_i}{ r^2 } \cdot  \f{ r^{-b_1} }{Q_m} \\
& \teq I \cdot x_i + II / Q_m .
\eal
\eeq

Since $I$ is the ratio among three radial weights non-singular near $r=0$, we can estimate their $C^1$ bounds. Using \eqref{eq:lei}, we further obtain the estimate of $I x_i$. 

For $II$, since $b_1 \leq -2$, a direct calculation yields 
\[
|II | \leq |a_1| x_i^u r_u^{-b_1 -2}, \quad |\pa_{x_j} II |
\leq |a_1|  ( \d_{ij} +  |b_1 + 2| \f{x_i}{r} \f{x_j}{r} ) r^{-b_1 - 2} .
\]
Applying \eqref{eq:cos_sin} for $x_j / r$ and $r^{-b_1 - 2} \leq r_u^{-b_1-2}$, we obtain the $C^1$ estimates for $II$. Using the $C^1$ estimates of $Q_m^{-1}, II$ and \eqref{eq:lei} and following Appendix {\appholfunc} in Part II \cite{ChenHou2023b_supp}, we obtain the $C^1$ and $C^{1/2}$ estimates of $II / Q_m$. Combining two parts, we obtain the refined estimate  near $r=0$.

In Section \ref{sec:non_hol}, we  estimate $ S= \f{x_1 \pa_1 P - x_2 \pa_2 P}{P}$ for a radial weight $P$. For $r$ away from $0$, the estimates follow from the previous method. Near $r=0$, using the above computation, we get 
\[
S = (x_1^2 - x_2^2) ( \f{ \sum_{2 \leq j \leq n_p} (a_j - a_1) p_j r^{a_j - a_1 - 2}  }{P_m}   + \f{a_1}{r^2})
 = \f{x_1^2 - x_2^2}{ r^2}( \f{\td P_m}{ P_m} + a_1),
\  \td P_m = \sum_{2 \leq j \leq n_p} (a_j - a_1) p_j r^{a_j - a_1 } .
\] 

Note that $\td P_m, P_m$ are radial weights and non-singular near $r=0$, we can estimate their $C^1$ bounds and the $C^1$ bounds of $\f{\td P_m}{ P_m} + a_1 $. 
We have 
\[
\tf{x_1^2 - x_2^2}{ r^2} = 1 - 2 f(x),  \quad f(x) \teq \tf{ x_2^2}{ r^2}, 
\]
and this singular factor is not $C^{1/2}$. We estimate $\d_i(f, x, z)|w|^{1/2}$ for $w = x, z, i=1,2$ in Section \ref{sec:dp_hol_singu}. Combining these estimates, we can estimate 
\beq\label{eq:dP_P_sing}
\d_i( S, x, z) |w|^{1/2}, \quad  w = x, z, \quad i = 1, 2, \quad S = (x \pa_x P - y \pa_y P) P^{-1}.
\eeq

\subsection{Piecewise estimates of damping coefficients $d_{i, L}^{\numm}$}

Recall the damping terms $d_{i, L}^{\numm}, d_{g_i, \UU}$ from \eqref{eq:dp}. 
In this section, we derive piecewise upper bounds of $d_{i, L}^{\numm}, d_{g_i, \UU}$. 
The technicalities mainly come from the singular weights near $0$. The reader who is not interested in rigorous verification can skip to Section \ref{sec:non_est} for the estimates of nonlinear terms and error terms.

\subsubsection{Piecewise bounds of $\cT_d^{\numm}$ from the singular weights}\label{sec:bd_dp_linf}

Let $\Om = [x_l, x_u] \times [y_l, y_u] \subset \R^2_{++}$. We estimate the upper bound of the following in $\Om$
\[
\cT^{\numm}_{d}(\vp) =  ( (\bar c_l x + \bar \uu^N) \cdot \na \vp)  / \vp, \quad \vp = |x|^{-1/2} P(r) + Q(r) ,
\]
where $P, Q$ are radial weights. We assume that $\vp$ has the form \eqref{eq:wg_nota}. Using the piecewise bounds of $\pa_x^i \pa_y^j P, \pa_x^i \pa_y^j Q ,$ 
\[
|\pa_x^i x^{-1/2}| \leq |\f{ (2i-1)!!}{2^i} x_l^{-1/2 -i}|,
\]
and \eqref{eq:lei}, we obtain piecewise bounds for $\vp$. We have several estimates for $\cT_d^{\numm}(\vp)$.

(1) We evaluate $\cT_d^{\numm}(\vp)$ on the grid point $(x_{\al}, y_{\b}), \al, \b = l, u$ and estimate $ \pa_i^2 \cT_d^{\numm}(\vp)$ using the piecewise bounds for $\bar \uu^N, \vp, \vp^{-1}$ and \eqref{eq:lei}. Then we use \eqref{est_2nd} to estimate $\cT_d^{\numm}(\vp)$. 

(2) We estimate the upper bound of the numerator $  (\bar c_l x + \bar \uu^N) \cdot \na \vp$ using its grid point values, derivative bounds, and \eqref{est_2nd}. Then we estimate $\cT_d^{\numm}(\vp)$ using this upper bound, the piecewise bounds for $\vp$, and \eqref{eq:func_intval2} for the ratio $f/ g$.

Since the weight $\vp$ is singular along the line $x=0$ and near $(x, y)=0$, we need more careful estimates. 
For $(x, y) \in \R^2_+$, we introduce $B_{u,i}$ below and rewrite $\cT_d^{\numm}(\vp)$ as follows 
\beq\label{eq:bd_dp_est1}
B_{u,1} = \bar c_l + \f{\bar u^N}{x}  , \quad B_{u,2} = \bar c_l + \f{\bar v^N}{y} , \quad 
 \cT_d^{\numm}(\vp) = B_{u,1} \f{x \pa_x \vp}{\vp} + B_{u,2} \f{y \pa_y \vp }{\vp} = I_1 + I_2.
\eeq

For $\vp$ \eqref{eq:wg_nota}, using $r \pa_r = x \pa_x + y \pa_y$, we have 
\beq\label{eq:wg_id1}
\bal
 x\pa_x \vp & =  - \tf{1}{2} x^{-1/2} P + x^{1/2} \pa_x P + x \pa_x Q , \quad
 y^k \pa_y \vp = x^{-1/2}  y^k \pa_y P +  y^k \pa_y Q , \ k \geq 0 ,\\
    x \cdot \na \vp  & = - \tf{1}{2} x^{-1/2} P + x^{-1/2}  r \pa_r P + r \pa_r Q \\
  \f{x \pa_x \vp}{\vp} & = \f{ - P/2 + x \pa_x P + x^{3/2} \pa_x Q}{ P + x^{1/2} Q},
 \  \f{y \pa_y \vp}{\vp} = \f{ y P_y + x^{1/2}  y Q_y}{ P + x^{1/2} Q}.
\eal
\eeq

To overcome the singularity $|x|^{-1/2}$, we rewrite $I_1$ as follows
\[
I_1 = B_{u,1} I_3,  \quad I_3 = \f{x\pa_x \vp}{\vp} =  \f{- P/2 + x \pa_x P}{ P + x^{1/2} Q}
+   \f{ x \pa_x Q}{  x^{-1/2 } P + Q} 
=   \f{- P/2 + x \pa_x P}{ P + x^{1/2} Q}
+  \f{ x^{3/2} \pa_x Q}{   P + x^{1/2} Q} .
\] 
We estimate the upper and lower bounds of $B_{u, 1} / (P+x^{1/2} Q )$ and $B_{u, 1} / (x^{-1/2} P + Q)$, which are positive, the upper bound of $-P/2 + x\pa_x P, x\pa_x Q, x^{3/2} \pa_x Q $, and then apply \eqref{eq:func_intval2}. We use the third identity of $I_3$ to remove the singular weight $|x|^{-1/2}$, and the second identity to control $I_3$ for large $r$ and $x$ close to $0$. If we use the third identity for such an estimate, since $Q$ decays slower than $P$, the upper bound is $ \f{   (x^{3/2} \pa_x Q)_u }{P}$, which grows to $\infty$. The same argument applies to estimate the upper and lower bounds of $x \pa_x \vp / \vp$. For  $x$ close to $0$ and large $r$, we use the second identity of $I_3$ for an improved estimate. We optimize  the estimates based on these two identities. The estimates of $I_2$ \eqref{eq:bd_dp_est1}, $y \pa_y \vp / \vp$ are similar and easier.

To overcome the singularity near $(x, y) = 0$, we have an additional estimate. Using the decomposition \eqref{eq:wg_mix_dec1} for $\vp = r^{\al}( x^{-1/2} P_m + r^{-\al} Q)\teq  r^{\al} \td \vp$, we get 
\[
 \f{ \pa \vp}{\vp} = \f{ \pa \td \vp}{\td \vp} + \f{ \pa r^{\al}}{ r^{\al}} , \quad
 \f{D \vp}{\vp} = \f{ D \td \vp}{\td \vp} + \f{ D r^{\al}}{r^{\al}}, \ 
\quad  \f{ x \pa_x r^{\al}}{r^{\al}} = \al \f{x^2}{r^2} , 
\ \f{ y \pa_y r^{\al}}{r^{\al}} = \al \f{y^2}{r^2} , \quad   D = x \pa_x, y \pa_y.
\] 
We rewrite $\cT^{\numm}_{d}(\vp)$ as follows 
\[
\bal
& \cT_d^{\numm}(\vp) = \f{(\bar c_l x + \bar \uu^N) \cdot \na (r^{\al} \td \vp ) }{ r^{\al} \td \vp}
 = \f{(\bar c_l x + \bar \uu^N) \cdot \na \td \vp }{ \td \vp } + \al \f{ B_{u,1} x^2 + B_{u,2} y^2}{r^2}
 \teq J_1 + J_2 ,
 \eal
\]
where we have used $B_{u,i}$ from \eqref{eq:bd_dp_est1}. 
The first part in the above formula of $D \vp / \vp$ or $\cT_d^{\numm}$ is only singular along $x=0$ and we apply the above estimates again. For $J_2$, our weights satisfy $\al < 0$. Let $B_{u, i, l}$ be the uniform lower bounds of $B_{u,i}$ in $\Om$. Since $\al < 0$, we yield 
\[
J_2 \leq \al f , \quad  f(x, y) = \f{B_{u, 1, l} x^2 + B_{u, 2, l} y^2}{ r^2}.
\]
If $B_{u,1, l} < B_{u, 2, l}$, $f$ is increasing in $y$ and decreasing in $x$, and we get $f(x, y) \geq f(x_u, y_l)$. If $B_{u, 1, l} \geq B_{u, 2, l}$, we get $f(x, y) \geq f(x_l, y_u)$ similarly.

\vs{0.1in}
\paragraph{\bf{Estimate of $d_{i, L}^{\numm}$ for radial weights}}

We apply the above arguments to derive the piecewise upper and lower bounds for $\cT_d^{\numm}, 
d_{i, L}^{\numm}$ \eqref{eq:dp} with a radial weight $\rho$, which is simpler. We apply the above method (1) to estimate the piecewise bounds for $ \pa_j \cT_d^{\numm}(\rho), \pa_j d_{i,L}^{\numm}$. We will use these estimates in the weighted H\"older estimate. 

\vs{0.1in}
\paragraph{\bf{Estimate of $d_{i, L}^{\numm}$ for $\vp_4$}}

In our nonlinear $L^{\infty}(\vp_4)$ estimate , we need to estimate 
\[
J = | \f{U \vp_{4,x}}{\vp_4} | +| \f{V \vp_{4,y}}{\vp_4}| ,
\quad 
\vp_4 = x_1^{-1/2} \psi_1, \quad \psi_1 = |x|^{-2} + 0.6 |x|^{-1} + 0.3 |x|^{-1/6} = \tts{\sum} p_i r^{a_i} .
\]
Using \eqref{eq:wg_id1} with $(P, Q) = ( \psi_1, 0)$ we get 
\[
|J| \leq |- \f{1}{2} \f{U}{x}| + |\f{ U \psi_{1, x}}{\psi_1} | 
+ |\f{ V \psi_{1, y}}{\psi_1}| 
\leq \f{1}{2}|\f{U}{x}| + \max( | \f{U}{x}|, | \f{V}{y}|) ( | \f{x \psi_{1,x}}{\psi_1}|
+ |\f{y \psi_{1,y}}{\psi_1}| ).
\]

Since $a_i <0$, for $(x, y) \in \R^2_{++}$, we get $\psi_{1,x}, \psi_{1,y } <0$, which along with 
\eqref{eq:wg_rad_dec1} gives
\[
0 \leq J_2 =  | \f{x \psi_{1,x}}{\psi_1}|
+ |\f{y \psi_{1,y}}{\psi_1}| 
=  - \f{x \psi_{1,x}}{\psi_1} 
- \f{y \psi_{1,y}}{\psi_1}
= -\f{x_1^2 + x_2^2 }{r^2} \f{ \sum p_i a_i r^{a_i} }{ \sum p_i r^{a_i}}
= - \f{ \sum p_i a_i r^{a_i} }{ \sum p_i r^{a_i}} \leq \max_i |a_i|.
\]

Similarly, using \eqref{eq:wg_id1}, we get 
\[
 \f{ ( \bar c_l x + \bar \uu^N)\cdot \na \vp_4  }{ \vp_4}
  = - \f{1}{2} ( \bar c_l + \f{ \bar u^N}{x} )
  +  \f{ ( \bar c_l x + \bar \uu^N)\cdot \na \psi_1  }{ \psi_1 }.
\]
The second term can be estimated using the above method.

\subsubsection{Piecewise bounds of $d_{ g_i, \UU}$ from the H\"older weights}\label{sec:bd_dg}

In this section, we estimate the damping term \eqref{eq:dp} $d_{g_i, \UU}$. Using \eqref{eq:damp_id}, we only need to estimate $ \f{ (\bar \uu^N+\UU)(x) - (\bar \uu^N+\UU)(z)) \cdot (\na g)(x- z)}{ g(x-z)}$. 

Recall the notation $\D_i(f)$ from \eqref{eq:hol_dif}. We have two estimates. In the first estimate, we estimate piecewise upper and lower bounds of the Lipschitz constant in $x$ 
\[
C(f,1, l)_{ i, k, j} \leq \f{f(z) -  f(x)}{z_1 - x_1} = \D_1(f, x,z) \leq C(f, 1, u)_{ i, k, j} , \quad x_2 = z_2, \ x_1 \leq z_1 \  x \in Q_{ij}, z \in Q_{ kj}
\]
for $ f = \bar u^N + U , \bar v^N + V $ using the piecewise bounds of $ \na \bar \uu^N + \na \UU$ and the method in Appendix {\apppiecederi} in Part II \cite{ChenHou2023b_supp}. The second index in $Q_{ij}, Q_{kj}$ is the same since $x_2 = z_2$. Similarly, we estimate the piecewise bounds for the Lipschitz constant in $y$. Since we treat $\UU$ as a perturbation, we get 
\[
  \pa_i \bar u^N - |\pa_i U| \leq  \pa_i(\bar u^N + U) \leq \pa_i \bar u^N + |\pa_i U|,
\]
and similar estimate for $\bar v^N + V$. We assume $x_1 \leq z_1$ without loss of generality.  Denote 
\[
r = z - x, \quad w = (z_1, x_2), \quad s = \sgn(r_1 \cdot r_2) . \quad 
\]

We have two cases for the configuration of $x, z$: $r_1 r_2 \geq 0$ and $r_1 r_2 < 0$,
Note that $ (\pa_1 g)(h) $ is odd in $h_1$ and even in $h_2$;  $\partial_2 g$ is even in $h_1$ and odd in $h_2$. 
We write $J_1 = (u(z) - u(x)) (\pa_1 g)(z-x) / g(z-x)$ as follows 
\[
\bal
J_1 & = (  u(w) - u(x) + u(z) - u(w)  ) \f{ (\pa_1 g)(z-x) }{g(z-x)} 
= \f{u(w) - u(x) }{ z_1 - x_1 } r_1 \f{(\pa_1 g)(z-x)}{g(z-x)}
+ \f{u(z) - u(w) }{ z_2 - x_2 } r_2 \f{(\pa_1 g)(z-x)}{g(z-x)}. \\
\eal
\]
Since $r_1 \geq 0$ and $r_1 \pa_1 g( r_1, r_2) = |r_1| \pa_1 g(|r_1|, |r_2|), 
r_2 (\pa_1 g)(r_1, r_2) = |r_2| s \pa_1 g(|r_1|, |r_2|)$, we get
\[
J_1 =  \D_1( u, x, w) |r_1| \f{(\pa_1 g)( |r_1|, |r_2| )}{g( |r_1|, |r_2|)}
 + \D_2(u, w, z) |r_2| s \f{(\pa_1 g)( |r_1|, |r_2| )}{g( |r_1|, |r_2|)}.
\]
The function $|r_1| \f{(\pa_1 g)( |r_1|, |r_2| )}{g( |r_1|, |r_2|)}$ is $0-$homogeneous, and we use the method in Section \ref{sec:hol_wg_est} to estimate it by partitioning the range of $|r_1| / |r_2|, |r_2| / |r_1|$. Similarly, we have 
\[
J_2 = \f{ (v(z) - v(x)) (\pa_2 g)(z-x)}{g(z-x)}
 = \D_1(v, x, w) s |r_1| \f{ \pa_2 g(|r_1|, |r_2|)}{ g(|r_1|, |r_2| )}
 + \D_2(v, w, z)  |r_2| \f{ \pa_2 g(|r_1|, |r_2|)}{ g(|r_1|, |r_2| )}.
\]
Using the estimates of $|r_1| \f{(\pa_1 g)( |r_1|, |r_2| )}{g( |r_1|, |r_2|)}$ in Section \ref{sec:hol_wg_est}, piecewise bounds of $\D_i(f, x, z)$, we obtain piecewise estimates of $d_{g_i, \UU}$.
See the red path $w= (z_1, x_2)$ in Figure \ref{fig:hol_pa} for an illustration of using $C_{x}^{1/2}, C_y^{1/2}$ estimates of $f$ to estimate $\d(f, x, z)$.

Note that we can also choose $ w = (x_1, z_2)$ and derive another decomposition 
\[
f(z) - f(x) = f(z) - f(w) + f(w) - f(x)
= \D_1(f,z, w) r_1  + \D_2(f, w, x) r_2,
\]
for $f = \bar u^N + U, \bar v^N + V$. Using this decomposition and the above argument, we derive another estimate of $ J_1 + J_2$. We optimize both estimates to bound $d_{g_i, \UU}$ from above.

\subsubsection{Piecewise H\"older estimates of the damping terms}\label{sec:bd_dp_hol}

In this section, we estimate the piecewise bounds $\d_j(  d_{i, L }^{\numm} W_{1, i} \psi_i, x, z )$ \eqref{eq:hol_dif} for $x \in Q_{ij}, z \in Q_{p q} $. We have
\beq\label{eq:bd_dp_est2}
\bal
 & \d( d_{i, L}^{\numm} W_{1, i} \psi_i, x,z ) g_i(x- z)
 = d_{i, L}^{\numm}(p) \d( W_{1, i} \psi_i, x, z) g_i(x-z) + I,  \\
 & I =  (W_{1, i} \psi_i)(q) \d( d_{i, L}^{\numm}, x, z) g_i(x- z) , \
 |I | \leq || W_{1, i} \vp_i||_{\inf} |\d( d_{i, L}^{\numm}, x, z) g_i(x- z) | \f{\psi_i(q)}{\vp_i(q)},
 \eal
\eeq
where $(p, q) = (x, z)$ or $(z, x)$. The first term is the damping term, and we need to estimate 
the weighted difference in $|I|$.

For $x, z$ away from $0$, we use the piecewise $C^1$ bounds of $d_{i ,L}^{\numm}$
established in Section \ref{sec:bd_dp_linf}, and the method in Appendix {\appholfunc} in Part II \cite{ChenHou2023b_supp} to estimate $\d_j( d_{i, L}^{\numm}, x,z)$. For $|x-z|$ not small, we have another estimate of $\d_j( d_{i, L}^{\numm}, x, z)$ using the piecewise upper and lower bounds of $d_{i, L}^{\numm}$. 

The main difficulty is the estimate near $x=0$ since $d_{i, L}^{\numm} \notin C^{1/2}$ near $r=0$. 
We use $\psi_i / \vp_i \les r^{1/2}$ near $ r= 0$ to compensate the low regularity of $d_{i, L}^{\numm}$ near $0$. 
We focus on the estimate of $\cT_d^{\numm}$ in $d_{i, L}^{\numm}$ \eqref{eq:dp}. The remaining parts $d_{loc, i}$ in $d_{i, L}^{\numm}$, i.e. $B_{u, 4}$ below, are locally $C^1$ and their estimates follow from the previous methods. 

Recall $B_{u,1}, B_{u,2}$ from \eqref{eq:bd_dp_est1} and we further introduce $B_{u,3}, B_{u,4}$
\[
\bal
& B_{u,1} = \bar c_l + \bar u^N / x, \ B_{u,2} = \bar c_l + \bar v^N / y, \ B_{u,3} = \bar v^N / y - \bar u^N / x = B_{u,2} - B_{u,1}, \\
& d_{i, L}^{\numm}(\rho) = \cT_d^{\numm}(\rho) + B_{u, 4}, \quad
 B_{u,4} = \bar c_{\om}^N \ (i=1), \  2 \bar c_{\om}^N - \bar u_x^N \ (i=2),  \ 2 \bar c_{\om}^N + \bar u_x^N \ (i=3).
\eal
\]

Suppose that $\rho$ has the representation \eqref{eq:wg_nota}. Using the decompositions \eqref{eq:bd_dp_est1}, \eqref{eq:wg_rad_dec1}, we get 
\begin{align}\label{eq:bd_dp_est3}
\cT_d^{\numm}(\rho) & = B_{u,1} \f{x  \pa_x \rho}{\rho} + B_{u,2} \f{y \pa_y \rho}{\rho}
= \B( B_{u,1} \f{x^2}{r^2} + B_{u,2} \f{y^2}{r^2}\B) (c_1 + \f{D_r \rho_m}{\rho_m} ) 
= \B( B_{u,1} + B_{u,3} \f{y^2}{r^2}\B) (c_1 + R_{\rho}) , \notag \\
& \teq \cT_R(\rho) + \cT_S(\rho) \f{y^2}{r^2}, \notag \\
\cT_R  & \teq  B_{u,1} (c_1 + R_{\rho}), \quad  \cT_S \teq B_{u,3} (c_1 + R_{\rho}), \quad R_{\rho} \teq \f{D_r \rho_m}{\rho_m}.
\end{align}
The functions $\cT_R, \cT_S, R_{\rho}$ are locally Lipschitz, and we derive their piecewise bounds below.

From the discussion around \eqref{eq:wg_rad_dec1}, we can estimate the piecewise $C^1$ bounds of $\rho_m, D_r \rho_m$. Using \eqref{eq:lei}, we can estimate the derivatives of $R_{\rho}$ and the piecewise bounds $\d_i(R_{\rho}, x, z)$. Using $\bar u^N(0, y) =0$, the piecewise estimates in Section \ref{sec:est_f_dx}, and
\[
 \pa_x B_{u,1} = \pa_x \f{\bar u^N}{x}
  = \f{ \bar u^N_x x - \bar u^N}{x^2 }
  = \f{1}{x^2} \int_0^x \bar u^N_{xx}(z, y) z dz  , \quad \pa_y B_{u,1} = \f{\bar u^N_y}{x},
\] 
we obtain piecewise $C^1$ estimates of $ B_{u,1}$. The same argument applies to the estimates of $B_{u, 2}$. Using these $C^1$ bounds, the Leibniz rule \eqref{eq:lei}, the method in Appendix {\appholfunc} in Part II \cite{ChenHou2023b_supp}, we obtain the piecewise $C^1$ and the H\"older estimate of the regular part $\cT_R, \cT_S$. 

\vs{0.1in}
\paragraph{\bf{Weighted estimate of the singular part}}
The difficult term is $\cT_S \f{y^2}{r^2}$, which is not $C^{1/2}$ near $r=0$. We estimate $\d_i(\cT_S \f{y^2}{r^2}, x, z)$ with $\psi_i(x) / \vp_i(x) \les r^{1/2}$ \eqref{eq:bd_dp_est2}. The factor $B_{u, 3} c_1 y^2 / r^2 $  captures the anisotropic structure of the advection in $x, y$ directions near $r=0$ . See Section {\secanisoest} \cite{ChenHou2023a_supp}. In  $Q_{ij}$ \eqref{eq:mesh_xy} near $(0, 0)$, the ratio $y^2 / r^2$ can vary a lot. Thus, we introduce the angle $\b$ below to control this factor and establish the following estimate 
\beq\label{eq:dp_hol_sinb}
\bal
&\d_1( \cT_R+ \cT_S f, x , z) |w|^{1/2} \leq C_{11}(x, z) + C_{12}(x, z) \sin^2(\b),  \  x_1 \leq z_1, \  x_2 = z_2,  \\
&\d_2( \cT_R+ \cT_S f, x , z) |w|^{1/2} \leq C_{21}(x, z) + C_{22}(x, z) \sin^{3/2}(\b), \  x_1 = z_1, \ x_2 \leq z_2 ,\\ 
&  \b = \max( \atan(x_2/x_1), \atan(z_2/ z_1)) , \quad f(x, y) = y^2 / r^2 ,\\
\eal
\eeq
for $ w= x$ and $w=z$. We have estimated $\cT_R, \cT_S$ in the above. For $x \in Q_{k_1 l_1}, z \in Q_{k_2 l_2}$, we can estimate the piecewise bounds for $C_{ij}(x, z)$. For $f(z)$, we use \eqref{eq:cos_sin} to obtain its piecewise bound. We defer the estimate $\d_i(f, x, z) |w|^{1/2}$ to Section \ref{sec:dp_hol_singu}. Using 
\[
\d_i( \cT_S f ,x ,z) |w|^{1/2}
 \leq \d_i(\cT_S, x, z) f(z) |w|^{1/2}
 + \d_i( f, x, z ) \cT_S(x) |w|^{1/2} ,
\]
the piecewise estimates of $\d_i(f, x, z) |w|^{1/2}$ in \eqref{eq:dp_sin_cx}, \eqref{eq:dp_sin_cy}, we obtain the estimate of $\d_i( \cT_S f ,x ,z) |w|^{1/2}$ depending on $\b$ in \eqref{eq:dp_hol_sinb}. Using the estimate in Section \ref{sec:dp_hol_singu}, we also obtain the estimate of $\d_i(f, x, z)|w|^{1/2}, \d_1( \cT_R+ \cT_S f, x , z) |w|^{1/2}$ independent of $\b$. Using \eqref{eq:dp_hol_sinb} and \eqref{eq:bd_dp_est3}, we get 
\beq\label{eq:dp_hol_sinb2}
\cT_d^{\numm}(\rho)  = \cT_R + \cT_S (\sin \b)^2.
\eeq
Near $ r = 0$, we have $\cT_S \approx B_{u,3} c_1 , B_{u,3} \approx 5, c_1 = -2, -5/2$ \eqref{eq:bd_dp_est3} in our estimates and thus $B_{u,3} c_1, \cT_S < 0$. Although the upper bound \eqref{eq:dp_hol_sinb} becomes larger for larger $\b$, we gain a larger damping factor from $\cT_d^{\numm}$ \eqref{eq:bd_dp_est2}, which dominates the upper bounds obtained in \eqref{eq:dp_hol_sinb}.

\vs{0.1in}
\paragraph{\bf{Improved estimate of the damping terms near $0$}}
We combine \eqref{eq:bd_dp_est3}, \eqref{eq:dp_hol_sinb}, and \eqref{eq:dp_hol_sinb2} to obtain an improved estimate for $\d( d_{i, L}^{\numm} , x,z ) W_{1, i}  g_i(x- z)$ \eqref{eq:bd_dp_est2} with any $x, z \in \R^2_{++}, z \neq x$ near $0$. We assume $x_1 \leq z_1$. There are two cases: (1) $z_2 \geq x_2$, (2) $z_2 \leq x_2$. Denote 
\[
 \lam_x = \f{x_2}{|x|}, \quad \lam_z = \f{z_2}{|z|}, \quad \lam = \max(\lam_x, \lam_z).
\]

For $\b \in [\b_i, \b_{i+1}]$, we get $ \lam \in [ \sin (\b_i), \sin \b_{i+1}]$ and then derive the upper bound of \eqref{eq:dp_hol_sinb2}, \eqref{eq:dp_hol_sinb}. We choose $p = x$ or $p=z$ in \eqref{eq:bd_dp_est2} such that $\lam_p = \lam$. From \eqref{eq:dp_hol_sinb2}, this choice provides a larger damping term $\cT_d^{\numm}(\rho)$. Then we choose $q = \{ x, z \} \bsh p$. Denote 
\[
h_1 = |x_1 - z_1|, \quad h_2 = |x_2 - z_2|, \quad R_i = \f{\psi_i}{\vp_i}.
\]

Recall the notations $\d, \d_i$ from \eqref{eq:hol_dif}.

\paragraph{\bf{Case (2): $z_2 \leq x_2$}}
Since $\f{y_2}{|y|}$ is increasing in $y_2$ \eqref{eq:cos_sin}, we have $\lam_z \leq \f{ x_2}{ |(z_1, x_2)| } \leq \f{x_2}{ |x|} =\lam_x$. Thus, we get $(p, q) = (x, z), \lam = \lam_x$. Denote $w = (z_1, x_2)$. We estimate 
\[
\bal
&  |\d( d_{i, L}^{\numm} , x, z) | g_i(h) R_i(z)
 \leq  ( |\d( d_{i, L}^{\numm} , x, w ) | + |\d( d_{i, L}^{\numm} , z, w) | ) g_i(h) R_i(z)  \\
 \leq & ( \d_1( d_{i, L}^{\numm}, x, w ) |z|^{1/2} \cdot  R_i(z) |z|^{-1/2} \cdot  h_1^{1/2} 
 + \d_2( d_{i, L}^{\numm}, z, w )  |z|^{1/2} \cdot  R_i(z) |z|^{-1/2} \cdot h_2^{1/2}  )g_i(h) . 
\eal
\]

Since $|z| \leq |(z_1, x_2)| = |w|$, we use \eqref{eq:dp_hol_sinb} to bound 
\[
\d_1( d_{i, L}^{\numm}, x, w ) |z|^{1/2} 
\leq \d_1( d_{i, L}^{\numm}, x, w ) |w|^{1/2},  \quad
\d_2( d_{i, L}^{\numm}, z, w ) |z|^{1/2} . 
\]
Using $R_i(y) \sim |y_1|^{1/2} |x|^a_i, a = (0.9, 0, 0)$ near $0$, we obtain piecewise bounds for $R_i(z) |z_1|^{-1/2}$, $R_i(z) |z|^{-1/2}$. To estimate $h_j^{1/2} g_i$, we follow Section \ref{sec:hol_wg_est}. 

\begin{figure}[t]
\centering
\begin{subfigure}{.26\textwidth}
  \centering
  \includegraphics[width=0.8\linewidth]{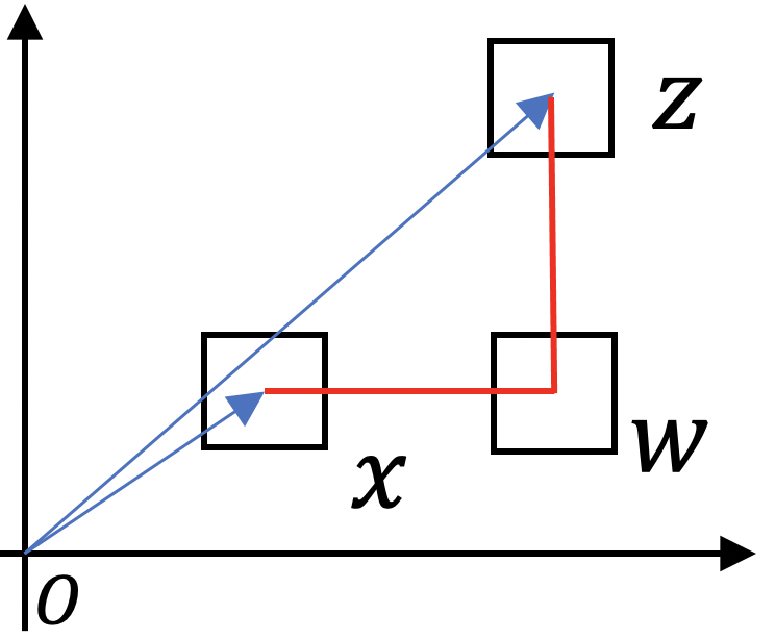}
    \end{subfigure}
    \begin{subfigure}{.36\textwidth}
  \centering
  \includegraphics[width=0.8\linewidth]{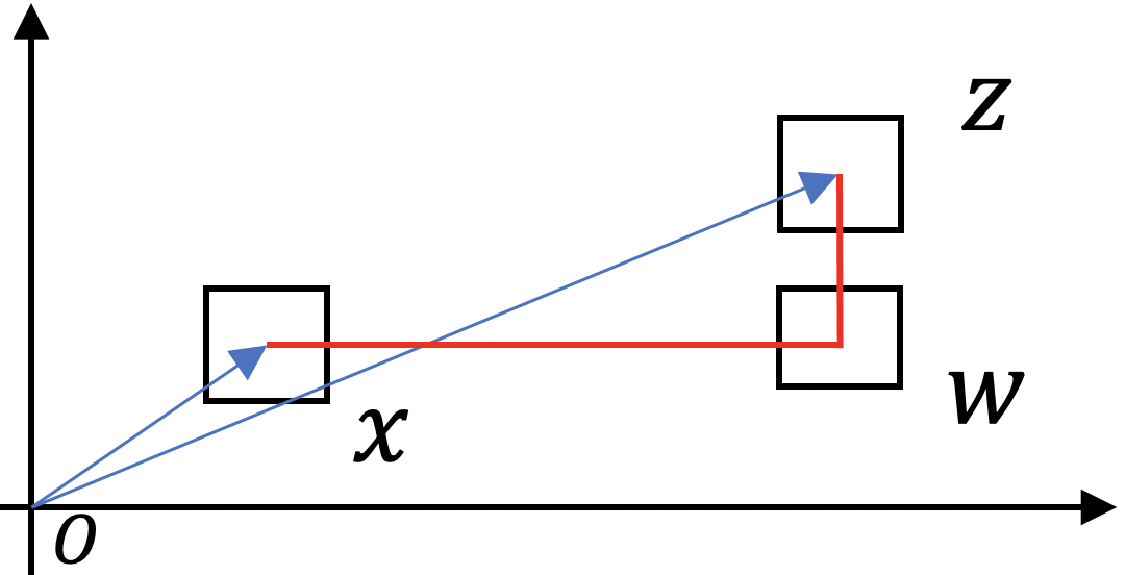}
    \end{subfigure}
        \begin{subfigure}{.36\textwidth}
  \centering
  \includegraphics[width=0.8\linewidth]{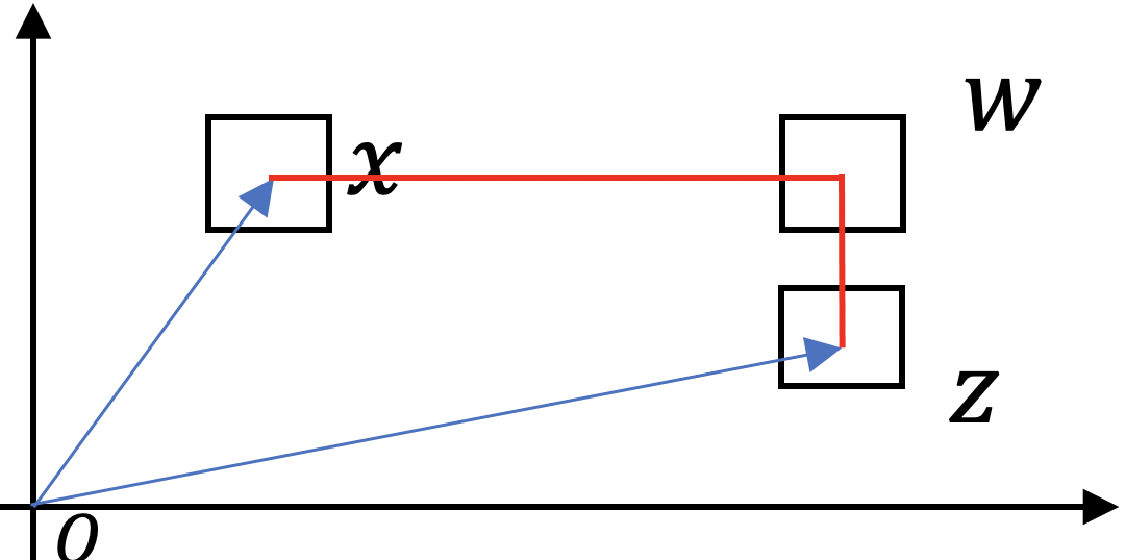}
    \end{subfigure}
\caption{Left, mid, right figure correspond to the locations of $(x, z)$ in cases (1.a), (1.b), (2). The red line represents the triangle inequalities and $C_{x_i}^{1/2}$ estimate used to estimate $|\d(d_{i, L}^{\numm}, x,z)|$.
 }
\label{fig:dp_hol}
\end{figure}

For (1) $x_2 \leq z_2$, we have two cases (1.a) $ \lam_x < \lam_z$ and (1.b) $ \lam_z < \lam_x$. 
In both cases, we choose $w = (z_1, x_2)$. In Figure \ref{fig:dp_hol}, we illustrate the locations of $(x, z)$ and the triangle inequalities used to bound $|\d( d_{i, L}^{\numm} , x, z )$ in different cases.

\paragraph{\bf{Case (1.a): $x_2 \leq z_2, \lam_x < \lam_z$}} We get $(p, q) = (z, x)$. Since $|x| \leq |w| \leq |z|$, we yield 
\[
\bal
&  |\d( d_{i, L}^{\numm} , x, z) | g_i(h) R_i(x)
 \leq  ( |\d( d_{i, L}^{\numm} , x, w ) |  + |\d( d_{i, L}^{\numm} , z, w) | ) g_i(h) R_i(x)  \\
 \leq & ( \d_1( d_{i, L}^{\numm}, x, w ) |x|^{1/2}  \cdot  h_1^{1/2} 
 + \d_2( d_{i, L}^{\numm}, w, z )  |w|^{1/2}  \cdot h_2^{1/2}  )g_i(h) 
\cdot R_i(x) |x|^{-1/2} . 
\eal
\]
We use \eqref{eq:dp_hol_sinb} to estimate $\d_1( d_{i, L}^{\numm}, x, w ) |x|^{1/2} ,\d_2( d_{i, L}^{\numm}, w, z )  |w|^{1/2}$ and estimate $R_i(x) |x|^{-\f{1}{2}}, h_j g_i(h)$ similarly.

\paragraph{\bf{Case (1.b): $x_2 \leq z_2, \lam_x > \lam_z$}} We get $(p, q) = (x, z)$ and 
\[
\bal
&  |\d( d_{i, L}^{\numm} , x, z) | g_i(h) R_i(z)
 \leq  ( |\d( d_{i, L}^{\numm} , x, w ) | + |\d( d_{i, L}^{\numm} , w, z) | ) g_i(h) R_i(z)  \\
 \leq & ( \d_1( d_{i, L}^{\numm}, x, w ) |z_1|^{1/2} \cdot  R_i |z_1|^{-1/2} \cdot  h_1^{1/2} 
 + \d_2( d_{i, L}^{\numm}, w,z )  |z|^{1/2} \cdot  R_i(z) |z|^{-1/2} \cdot h_2^{1/2}  )g_i(h) . 
\eal
\]
We estimate $R_i(z) |z|^{-\f{1}{2}}, R_i(z) |z_1|^{- \f{1}{2}}, h_j g_i(h)$ similarly. Since $|z_1| \leq |w|$,  we use \eqref{eq:dp_hol_sinb} to bound 
\[
\d_1( d_{i, L}^{\numm}, x, w ) |z_1|^{1/2} 
\leq \d_1( d_{i, L}^{\numm}, x, w ) |w|^{1/2} , \quad \d_2( d_{i, L}^{\numm}, w,z )  |z|^{1/2}.
\]

\vs{0.1in}
\paragraph{\bf{Four parameters}}
 In summary, in the improved estimate near $0$, we introduce four parameters $x, z, \lam = (\lam_x, \lam_z) , h_i = |x_i -z_i|$ and derive the bound 
\beq\label{eq:dp_hol_sinb3}
\bal
& d_{i, L}^{\numm}(p) \leq C_{31}(x, z) + C_{32}(x, z ) \lam^2, \\
& \d( d_{i, L}^{\numm} ,x,z ) g_i(x- z)   |W_{1, i} \psi_i(q)|
\leq \tts{ \sum_{l} } \, C_{l}(x, z) h_1^{a_{l1}} h^{1/2- a_{l1}}_2 g_i(h) \lam^{ a_{l3}} E_4.
\eal
\eeq
For each $x, z \in Q_{i_1 j_1}, Q_{i_2 j_2}$, we first derive the piecewise bounds for $C_l(x, z)$. By partitioning $\lam  = \sin(\b) \in [0, 1]$ (or $\b \in [0, \pi/2]$) into small intervals and using monotonicity in $\lam$, we estimate $\lam^a$.
The function in $h$ is $0$-homogeneous, which can be estimated following Section \ref{sec:hol_wg_est}.

\subsection{Estimate of nonlinear terms}\label{sec:non_est}

We estimate the nonlinear terms 
$ \td \cN_i( \rho_i)$ in \eqref{eq:lin_U}, \eqref{def:non_wg}. The estimate of each term follows from direct $L^{\inf}$ or $C^{1/2}$ estimates \eqref{eq:hol_prod0} \eqref{eq:hol_sq}. The main technicality comes from the singular weights near $r = 0$ or $x = 0$.
We have discussed the $C^k$ estimates of $\hat W_2$ in Section \ref{sec:W2_est}.

\subsubsection{Weighted H\"older estimate}\label{sec:non_hol}

In the H\"older estimate, we use radial weights $\rho_i = \psi_i$ and estimate $\mu_{h,i} [W_{1, i} \psi_i]_{ C_{g_i}^{1/2}}$ in the energy $E_4$ \eqref{energy4}.

\subsubsection{Piecewise estimates of $\cT_{d, N}(\psi_i)$}

The main technical term is the coefficient $\cT_{d, N}(\psi_i)$ (defined in \eqref{eq:dp}) from the first term $ \cT_{d, N}(\psi_i) ( W_{1, i} \psi_i)$ in \eqref{def:non_wg}. 
Since we have piecewise $C^{1/2}$ estimates of $ W_{1, i}\psi_i$ using the energy,
we only need to control $\cT_{d, N}(\psi_i)$ in order to estimate $ \cT_{d, N}(\psi_i) ( W_{1, i} \psi_i)$.
Using \eqref{eq:U}, we perform the decomposition 
\[
\cT_{d, N}(\psi_i) = \f{\UU \cdot \na \psi_i }{ \psi_i}
= U_x(0) \f{( x \pa_x \psi_i - y \pa_y \psi_i)}{ \psi_i }
 + \f{ \UU_A \cdot \na \psi_i}{\psi_i}
 + \f{ \UU_{\FM}  \cdot \na \psi_i}{\psi_i}
\teq \cT_{ c_{\om}} + \cT_{uA} + \cT_{uR}.
\]
For $\cT_{uA}$, since we have piecewise $C^{1/2}$ estimates of $\UU_A \psi_u$, we decompose it as follows 
\beq\label{eq:non_dec3}
\cT_{uA} = ( U_A  \psi_u ) \f{\pa_x \psi_i}{ \psi_i \psi_u}
+ ( V_A  \psi_u ) \f{\pa_y \psi_i}{ \psi_i \psi_u}, 
\quad \cT_{u R} =  \f{ \UU_{\FM} }{|x|^2}
 \cdot |x|^2 \f{\na \psi_i}{\psi_i}.
\eeq
Since $\psi_i, \psi_u$ are radial weights, we can estimate piecewise derivatives of $\psi_i, \psi_i^{-1}, \psi_u^{-1}$ following Appendix {\appwgradial} in Part II \cite{ChenHou2023b_supp}. See also Sections \ref{sec:wg_rat}, \ref{sec:wg_rat_rad} for the estimates of the ratio among the weights. Using these 
estimates, \eqref{eq:lei} and following Appendix {\appholfunc} in Part II \cite{ChenHou2023b_supp}, we can obtain piecewise $C^1, C^{1/2}$ estimates of 
\beq\label{eq:non_wg_est1}
\f{\pa_j \psi_i}{ \psi_i },  \quad  \f{ x \pa_x \psi_i - y \pa_y \psi_i}{ \psi_i },
\quad \f{\pa_j \psi_i } { \psi_i \psi_u } , \quad  |x|^2 \f{ \pa_j \psi_i }{ \psi_i}.
\eeq
The second ratio is estimated in Section \ref{sec:wg_rat_rad}. Using these estimates, \eqref{est_2nd}, and evaluating these functions on a fine mesh, we can refine the estimate. We use the estimates in \eqref{eq:dP_PQ} in Section \ref{sec:wg_rat_rad}
with $(P, Q) = (\psi_i, \psi_u), (\psi_i, |x|^{-2})$  to refine the estimate of $(\pa_j \psi_i ) (\psi_i \psi_u)^{-1},  |x|^2 (\pa_j \psi_i ) \psi_i^{-1} $  near $r=0$. We remark that the $C^1$ estimates of these two terms are not singular near $r=0$, while the $C^1$ estimates of the first two terms can be singular.

Using the $C^{1/2}$ estimates of $\f{\pa_j \psi_i}{ \psi_i \psi_u}$ and 
of $ u_A \psi_u, v_A \psi_u$ from the functional inequalities, \eqref{eq:hol_prod}, we obtain the $C^{1/2}$ estimate of $\cT_{d, N}$. For $\cT_{u, R}$, we can estimate 
\[
 |x|^{-k} \UU_{\FM} ,\  k=2, 3 ,  \quad 
  |x|^{-2} \pa_i  \UU_{\FM},
\]
and yield the $C^1$ estimates of  $ \f{ \UU_{\FM} }{|x|^2}$, which along with the $C^1$ estimates of $|x|^2 \f{ \pa_j \psi_i}{\psi_i}$ imply the $C^1, C^{1/2}$ estimates of $\cT_{u R}$. We have another estimate of $\cT_{u R}$ using the $C^1$ bounds of 
$ \UU_{\FM} ,  \pa_j \psi_i / \psi_i$ \eqref{eq:non_wg_est1}. We optimize these two estimates for $\cT_{u R}$.

Using the above (weighted) $C^{1/2}$ estimates for $\cT_{uA}, \cT_{uR}$, and \eqref{eq:hol_loc}, we get
\beq\label{eq:hol_non_tran1}
\bal
& \mu_{h, i} g_i(x-z) \big| \d\big( \, \big( \cT_{uA}(\psi_i) + \cT_{uR}(\psi_i) \big) ( W_{1, i} \psi_i), x, z \big)  \big| \\
\leq &  \B( |\cT_{uA}(x) + \cT_{uR}(x)|  
+ \mu_{h, i} g_i(x-z) \d_{\sq}(\cT_{uA} + \cT_{uR}, x, z, x-z) (\psi_i/ \vp_i)(z) \B) E_4 . 
\eal
\eeq

For $\cT_{c_{\om}}$ away from $r=0$, we use the $C^1$ estimate of $\cT_{c_{\om}}$ \eqref{eq:non_wg_est1}. Near $r=0$, we use the estimates in \eqref{eq:dP_P_sing} and Section \ref{sec:wg_rat_rad} to estimate $\d_i( \cT_{c_{\om}}, x, z ) \min(|x|, |z|)^{1/2}$. 
We have
\beq\label{eq:hol_non_tran_cw}
\bal
 \mu_{h, i} |\d( \cT_{c_{\om}} W_{1, i} \psi_i , x, z)  g_i & (x-z)|
 \leq \mu_{h, i} |\d(  W_{1, i} \psi_i , x, z) g_i(x-z) | \max_{p=x,z} |\cT_{c_{\om}}(p)| + \mu_{h, i} I  \\
& \qquad \quad  \ \leq  E_4 \max_{p=x,z} |\cT_{c_{\om}}(p)| + \mu_{h, i} I  ,  \\
 I &\leq |\d(  \cT_{c_{\om}} , x, z ) g_i(x-z)| \min(  |W_{1, i} \psi_i(x)|, |W_{1, i} \psi_i(z)|) \\
& \leq   | \d(  \cT_{c_{\om}} , x, z )| g_i(x-z) \min(|x|,|z|)^{1/2} \max_{p=x,z}\f{\psi_i(p)}{\vp_i(p) |p|^{1/2}}   E_4.
\eal
\eeq
The ratio $ \f{\psi_i}{\vp_i |x|^{1/2}} $ is estimated in \eqref{eq:wg_rat_C1C0}. The weighted H\"older estimate of $\cT_{c_{\om}}$ is similar to that in Section \ref{sec:bd_dp_hol}. Denote 
\[
h_i  = |x_i - z_i|, \quad M = \min(|x|, |z|) .
\]

We assume  $x_1 \leq z_1$ and consider two cases: (1) $x_2 \leq z_2$ ; (2) $x_2 \geq z_2$. Recall $\d, \d_i$ from \eqref{eq:hol_dif}.

In case (1) of $x_2 \leq z_2$, for $w = (z_1, x_2)$, we have $|x| \leq |w| \leq |z|$, and 
\[
\bal
II & \teq | \d(  \cT_{c_{\om}} , x, z )| g_i(x-z) M^{1/2}
  \leq |\d( \cT_{c_{\om}}, x, w) + \d( \cT_{c_{\om}}, w, z)| g_i(x-z) M^{1/2} \\
  & \leq  (\d_1( \cT_{c_{\om}}, x, w) |x|^{1/2}  h_1^{1/2}
  +\d_2( \cT_{c_{\om}}, w, z) |w|^{1/2} h_2^{1/2}  )  g_i(h),
  \eal
\]
and we then estimate $h_j^{1/2} g_i(h)$ following the methods in Section \ref{sec:hol_wg_est}. We can also choose $ w = ( x_1, z_2)$ and obtain another estimate. 

In case (2) $x_2 \geq z_2$, we choose $w = ( z_1, x_2)$. Then we have $|x|, |z| \leq |w|$ and 
\[
II \leq 
|\d( \cT_{c_{\om}}, x, w) + \d( \cT_{c_{\om}}, z, w)| g_i(h) M^{1/2} \\
   \leq  (\d_1( \cT_{c_{\om}}, x, w) |x|^{1/2}  h_1^{1/2}
  +\d_2( \cT_{c_{\om}}, z, w) |z|^{1/2} h_2^{1/2}  )  g_i(h).
\]

\begin{figure}[t]
\centering
\begin{subfigure}{.49\textwidth}
  \centering
  \includegraphics[width=0.45\linewidth]{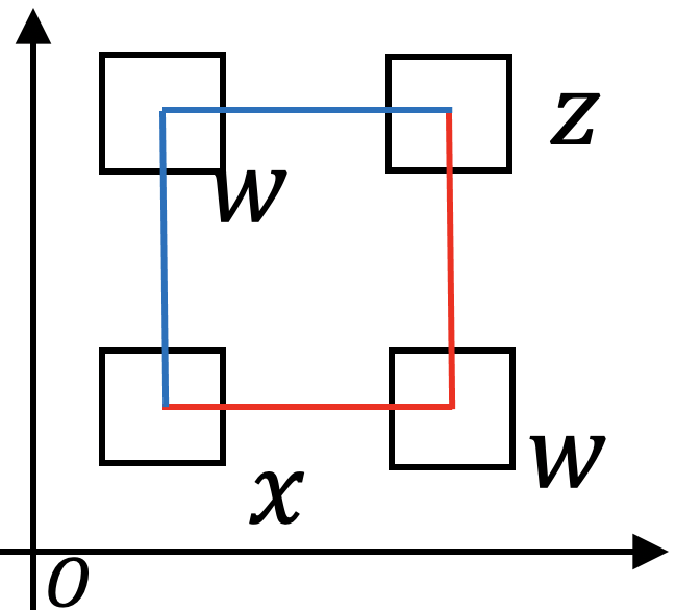}
    \end{subfigure}
    \begin{subfigure}{.49\textwidth}
  \centering
  \includegraphics[width=0.45\linewidth]{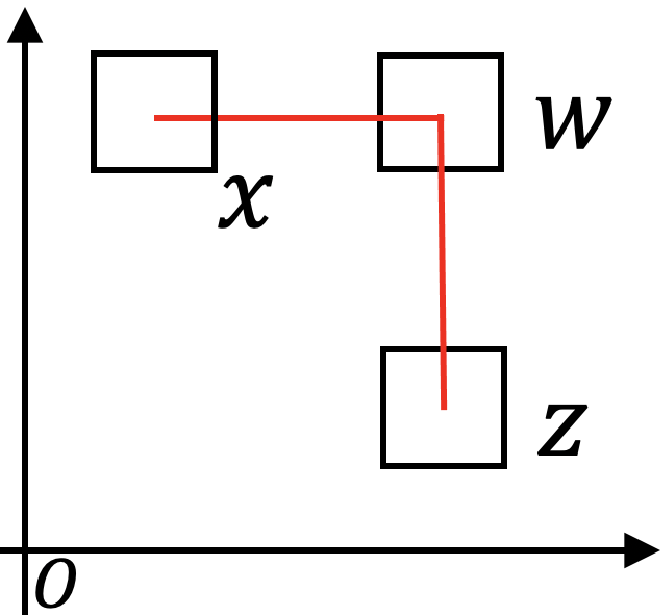}
    \end{subfigure}
\caption{Left, right figures correspond to the locations of $(x, z)$ in cases (1), (2). The red line and blue line represents the $C_{x_i}^{1/2}$ estimate 
used to estimate $|\d( \cT_{c_{\om}}, x,z)|$. 
In case (1), we have two estimates. 
 }
\label{fig:hol_cw}
\vspace{-0.1in}
\end{figure}

In the above estimates of $II$, we further use the estimate of $\d_j( \cT_{c_{\om}}, p, q) \min( |p|, |q|)^{\f{1}{2}}$. In Figure \ref{fig:hol_cw}, we illustrate $C_{x_i}^{\f{1}{2}}$ estimate and  triangle inequalities used to bound $|\d( \cT_{c_{\om}} , x, z )|$ in various cases. 

\subsubsection{Other nonlinear terms $ N_{W_1, i}, N_{\hat W_2, i}$}\label{sec:non_hol_W1W2}
The estimates of other nonlinear terms 
$ N_{W_1, i}, N_{\hat W_2, i}$ in $\td \cN_i$  \eqref{def:non_wg} are simpler.  We have discussed their estimates in the main paper \cite[Section 5.9.3]{ChenHou2023a_supp}. 
For completeness, we present them below.  We follow Sections \ref{sec:non_dec_u}, \ref{sec:W2_est} for the $L^{\inf}$ and $C^{1/2}$ estimates of $\na \UU, \UU$ \eqref{eq:U} and the 
$C^k$ estimates of $\hat W_2$. 

\vs{0.1in}
\paragraph{\bf{Nonlinearity $N_{W_1,i}$ involving $\UU, W_{1, i}$}}

We estimate nonlinear terms $N_{W_1,i}$ in \eqref{def:N_W1}, \eqref{eq:non_dec1} involving $\om_1, \eta_1, \xi_1$. 
We write down the decomposition of each term below as $A \cdot B$ and estimate $\d(f, x, z)$ for each $f = A, B$,
\beq\label{eq:non_dec_est1}
\bal
N_{W_1, 1} \psi_1   =  U_x(0) \cdot \om_1 \psi_1,  \quad
N_{W_1, 2} \psi_2  & = U_x(0) \cdot \eta_1 \psi_2  - \td U_x \cdot (\eta_1 \psi_2) - \td V_x \cdot (\xi_1 \psi_2) , \\
 N_{W_1, 3} \psi_3 & = 3 U_x(0) \cdot (\xi_1 \psi_2)
- \td U_y \cdot (\eta_1 \psi_2) - \td V_y \cdot ( \xi_1 \psi_2).
\eal
\eeq
 Since each term involves $W_{1, i} \psi_i$, we can estimate it using \eqref{eq:hol_loc}.

Note that we use the same weight $ \psi_2= \psi_3, g_2 = g_3$ for $\eta_1, \xi_1$ in the H\"older estimate.

It remains to estimate $N_{\hat W_2, i}$ in \eqref{def:non_wg},\eqref{eq:non_dec1}.
We use the decomposition \eqref{def:N_W2} and estimate three terms separately.

\vs{0.1in}
\paragraph{\bf{Nonlinear terms $\cB_{op, i}(\UU_A  , \wh W_2)$ in 
$N_{\hat W_2, i}$}}
We estimate $ \cB_{op, i}( \UU_A  ,  \wh W_2 ) \psi_i$ in 
\eqref{def:N_W2}, \eqref{eq:non_dec1}.  Similar to \eqref{eq:non_dec_est1}, we write down the decomposition of each term in $\cB_{op, i}( \UU_A ,  \wh W_2 ) \psi_i$  \eqref{eq:Blin}, \eqref{eq:B_abuse} below as $A \cdot B$ and estimate $\d_{\ell}(f, x, z)$ for each $f = A, B$ and then use \eqref{eq:hol_prod}, \eqref{eq:hol_sq2} to estimate the $AB$. 
For the advection term of $\wh W_{2, i}$ \eqref{eq:non_dec1}, \eqref{eq:B_abuse}, we use
\[
 - \UU_A \cdot \na \wh W_{2, i} \psi_i 
 = -( U_A \psi_u) \cdot  ( \f{ \pa_x \wh W_{2, i} \psi_i }{\psi_u} ) 
 - ( V_A \psi_u) \cdot  ( \f{ \pa_y \wh W_{2, i} \psi_i }{\psi_u} ) , \ i = 1,2,3.
\]
Near $0$, $\na \wh W_{2, i} \f{\psi_i}{ \psi_u}$ has a  vanishing order $O(|x|)$, and we can estimate its $C^1$ bound.

For other nonlinear terms in $\cB_{op,j}(U_A, \wh W_2), j=2,3$ involving $ (\na \UU)_A$,
\footnote{
From \eqref{eq:Blin_gen}, $\cB_{op,1}(U_A, \wh W_2)$ does not involve 
velocity-gradient  $ (\na \UU)_A$-terms. 
} 
we use 
\[
(- (\pa_i U)_{ A} \hat \eta_2 - (\pa_i V)_{ A} \hat \xi_2 )  \psi_2 
= - ( ( \pa_i U)_A \psi_1) \cdot  (\f{\hat \eta_2 \psi_2}{\psi_1} )- 
(  (\pa_i V)_A \psi_1 ) \cdot  (\f{ \hat \xi_2 \psi_2}{\psi_1} ), \quad i = 1, 2.
\]
We use $i=1, \pa_i = \pa_x$ for $\cB_{op, 2}$ and $i=2, \pa_i = \pa_y$ for $\cB_{op, 3}$. The above estimate is similar to the linear H\"older estimate in Section {\secEE} \cite{ChenHou2023a_supp} and we treat $\hat W_2$ similar to $\bar \om, \bar \th_x, \bar \th_y$.

\vs{0.1in}
\paragraph{\bf{Remaining nonlinearities of $\wh W_2$}} The remaining terms 
of $\wh W_2$ in 
$N_{\hat W_2, i}$ \eqref{eq:non_dec1},  \eqref{def:N_W2} are 
\[
  I_1 = U_x(0) \hat W_{2, i, M} , \quad 
   I_2 = \cB_{op, i}( \UU_{\FM}, \wh W_2). 
\]
With the second order correction near $x=0$, $I_1$ vanishes like $O(|x|^{3})$. By definition, we have 
\[
\wh W_{2, i} = O(|x|^2),  \quad \UU_{\FM} = O(|x|^3).
\]
It follows  that $I_2$ vanishes to the order $O(|x|^4)$ near $|x| = 0$.
For $I_1$, we estimate $C^{1/2}$ of $ \hat W_{2, i, M}  \psi_i$ by optimizing the estimate of $C^{1/2}$ estimate of $(
 \hat W_{2, i, M}  |x|^{\al_i} ) \cdot \psi_{m, i} $ and $  \hat W_{2, i, M} \psi_i$, where $\psi_i = |x|^{\al_i} \psi_{m, i} $ \eqref{eq:wg_rad_dec1}. Near $r=0$, we estimate the $C^{1/2}$ semi-norm of $
  \hat W_{2, i, M}  |x|^{p}, p = -2, -5/2 $ carefully  using the method in Section \ref{sec:hol_pow} to overcome the singularity near $r=0$.

For $I_2$, $ I_2 \psi_i = O(|x|^{ a}),a \geq \f{3}{2}$ near $r=0$,  we obtain piecewise $C^1$ estimates of $ I_2 \psi_i$. We estimate a typical term $ T = \pa_x^i \pa_y^j ( U_{\FM} \pa_x \hat \eta_2 \psi_2)$ in $\cB_{ op, 2}( \UU_{\FM}, \hat W_{2})\psi_2 $ \eqref{eq:non_dec1} as follows. Away from $0$, we estimate $T$ using \eqref{eq:lei}, the estimate of $U_{\FM}, \hat \eta_2$ and the triangle inequality. Near $r =0$, we have three terms 
\[
 D (U_{\FM}) \cdot \pa_x \hat \eta_2 \cdot \psi_2 + U_{\FM} \cdot  D \pa_x \hat \eta_2 \cdot \psi_2 + U_{\FM} \pa_x 
 \hat \eta_2 \cdot D \psi_2 = J_1 + J_2 + J_3, \quad D = \pa_x^i \pa_y^j,  \ i+j \leq 1.
\]
For $J_1$, since $U_{\FM} = O(|x|^3), \hat \eta_2 = O(|x|^2) $, we extract these vanishing orders and estimate 
\[
(D U_{\FM}   ) \cdot (\pa_x \hat \eta_2 / |x|) \cdot (\psi_2 |x|^3).
\]
The estimate of $J_2$ is similar. For $J_3$, we use \eqref{eq:wg_rad_dec1} for $\psi_2 = \psi_{m, 2} |x|^{\al_2}$. We obtain 
\[
J_3 = ( U_{\FM} \pa_x \hat \eta_2) |x|^{\al_2 -1} ) \f{ D \psi}{ |x|^{\al_2-1}}, 
\quad \f{ D \psi_2}{ |x|^{\al_2-1}} = D \psi_{m ,2} |x|
+ \al_2 \psi_{m, 2}  \f{x_1^i x_2^j}{ |x|} , \ i+ j = 1.
\]
With the above $C^{1/2}$ estimates, we can estimate $\d_{\sq}( I_k, x, z, x-z), k=1,2$ \eqref{eq:hol_sq}.

\subsection{$L^{\inf}$ stability analysis in the far-field}\label{sec:linf_far}

Recall parameters $N_i, N$ from \eqref{eq:mesh_xy}. For $|x|_{\inf} \geq y_{N_2-1}$, the piecewise polynomial part in the approximate steady states is $0$. We have
\beq\label{eq:lin_asym_ASS}
\bar \om(x) = \bar \om_1 = \bar g_1(\b) r^{\bar \al_1}, \quad  \bar \th(x) = \bar \th_1 = \bar g_2(\b) r^{ 1 + 2 \bar \al_1} , \ r = (x_1^2+ x_2^2)^{1/2}, \ \b = \atan(x_2/  x_1).
\eeq

To perform energy estimate in the far-field, we need to estimate the asymptotic behavior of $\pa_x^i \pa_y^j f(r, \b)$. Using the formulas of derivatives, 
\[
\pa_x g = (\cos \b \pa_r - \f{ \sin \b}{r} \pa_{\b} ) g, \quad  \pa_y  g= (\sin \b \pa_r + \f{ \cos \b}{r} \pa_{\b} ) g, 
\]
the triangle inequality, \eqref{eq:lin_asym_ASS}, and the piecewise bounds of $\pa_{\b}^i \bar g_j(\b)$ established using the method in Appendix {\apppiecehol}, {\secestapprfar} in Part II \cite{ChenHou2023b_supp}, we obtain 
\beq\label{eq:lin_asym_ASS2}
|\pa_x^i \pa_y^j \bar \om | \leq \bar \om_{ag, ij}(\b) r^{ \bar \al_1 - i - j }
, 
\quad 
|\pa_x^i \pa_y^j \bar \th | \leq \bar \th_{ag, ij}(\b) r^{1 + 2 \bar \al_1 - i - j }, 
\quad 
|\pa_x^i \pa_y^j \bar \phi | \leq \bar \phi_{ag, ij}(\b) r^{ 2 + \bar \al_1 - i - j }, 
\eeq
for $|x|_{\inf} \geq y_{N_2-1}$, and some piecewise constants $\bar \om_{ag, ij}(\b), \bar \th_{ag, ij}(\b), \bar \phi_{ag, ij}(\b)$.

We perform the weighted $L^{\inf}(\vp_i)$ in the far field $|x|_{\inf} \geq y_{N-1}$. 
We can decompose the weights $\vp_i$ as follows and have the following estimate uniformly for $|x| \geq y_{N-1}$ from \eqref{eq:lin_asym_wg1} 
\beq\label{eq:lin_asym_linfwg1}
\vp_i = |x|^{-1/2 } P_{1, i} + P_{2, i},  \quad | P_{j, i}(x) | \leq c_{j, i} r^{ \al_{ j, i} }, \quad  j= 1, 2, i=1,2,3,
\eeq
for some radial weights $P_{j, i},c_{j, i}$, with $P_{1, i}$ decaying faster than $P_{2, i}$. We estimate the nonlocal terms, the local terms in \eqref{eq:lin_U}, \eqref{eq:non_dec1}, and the damping terms $d_{i, L}^{\numm}$ \eqref{eq:dp} in order.

\vs{0.1in}
\paragraph{\bf{Estimate of $f/ |x|^{\al}$}}
Since the weight involves the power $|x|^{-1/2}$ singular along $x=0$, we need to estimate the asymptotic behavior of $f / |x|^{1/2}$ for $f = u_A, \bar \om_y, \na \bar \th, \bar \th_{xy}$. We estimate $ u_A / |x|^{1/2} r^{\g} $ for suitable power $\g$ at the end of Section {\seconeDlinffar} in the supplementary material of Part II \cite{ChenHou2023bSupp_supp}. For other functions $f$ depending on the profile, if $f(0, y)= 0$, or equivalently $f(r, \pi/2) = 0$, we can estimate $ f / |x|^{\al} 
= f / (r \cos \b)^{\al}$ for $\al \in [0, 1]$. First, using the polar representation \eqref{eq:lin_asym_ASS} and the method in Appendix {\secestapprfar} in Part II \cite{ChenHou2023b_supp}, we can estimate $\pa_{\b} f$. 

Consider polar coordinates $(r,\b)$ for $(x,y)\in \R^2$: $x = r \cos \b, y = r \sin \b$. Using the piecewise derivative bounds discussed in Appendix {\apppiecederi} in Part II \cite{ChenHou2023b_supp} and 
\begin{align}\label{eq:lin_asym_F_ratio}
& \f{f(r, \b)}{ \pi / 2 - \b} = - \f{1}{\pi / 2 - \b} \int_{\b}^{\pi/2}  \pa_{\tau}(f, \tau) d \tau ,   \quad \cos \b = (\f{\pi}{2} - \b) \sin \xi,  \ \xi \in [\b, \f{\pi}{2} ] ,   \\
& \f{f(r, \b)}{|x|^{\al}}
=  r^{-\al} \f{f(r, \b)}{ \f{\pi}{2} - \b} ( \f{\pi}{2}-\b)^{1-\al} (\f{\pi/2-\b}{\cos \b})^{\al}
= r^{-\al} \f{f(r, \b)}{ \f{\pi}{2} - \b} ( \f{\pi}{2}-\b)^{1-\al} ( \sin \xi )^{ - \al} , \
 \f{1}{|\sin \xi |^{ \al}} \leq \f{1}{ |\sin \b|^{\al}} , \notag
\end{align}
we can obtain the piecewise bounds of $\f{f(r, \b)}{ \pi / 2 - \b}$ and $ f(r, \b) / |x|^{\al}$.

Suppose that $f(0, y) = 0$, $ |f(r, \b)| r^b \leq C_1, |\pa_x f r^{b+1}| \leq C_2 $ for all $\max(x, y) \geq R$ and $b+1 \leq 0$.  Using $ |f_x(z, y)| \leq C_2 r^{-b-1} $ for $z \in [0, x]$, we have another simple estimate for $f / x^{1/2}$
\[
\bal
& |\f{f}{x^{1/2}} |=| \f{1}{x^{1/2}} \int_0^x f_x(z, y) dz |
\leq C_2 r^{-(b+1)} x^{1/2} = C_2 r^{-b-1/2} (\cos \b)^{1/2}
\leq C_2 r^{-b} R^{-1/2}(\cos \b)^{1/2}, 
 \ 
\b \in [ \f{\pi}{4}, \f{\pi}{2}], \\
& |\f{f}{x^{1/2}} | \leq \min( C_1 (\cos \b)^{-1/2} r^{-b-1/2}, C_1 r^{-b} R^{-1/2} )
\leq \min( C_1 (\cos \b)^{-1/2} R^{-1/2}, C_1  R^{-1/2} ) r^{-b},
\  \b \in [0, \f{\pi}{4}],
\eal
\]
where we have used $r \geq \max(x, y) \geq R$ in the second case. The above estimates also provide a uniform estimate for $f / x^{1/2} r^{b+1/2}$.

\vs{-0.1in}
\subsubsection*{\bf Estimate of the nonlocal terms}\label{sec:far_est_nloc}

First, we simplify the formula of $\hat \uu( x), \wh{ \na \uu }( x)$ for large $|x|$. From the definition of finite rank approximation in 
\cite[Section 4.3.3]{ChenHou2023a_supp} for $\hat f_1, \hat f_2$ 
and approximations  parameters in \cite[Appendix C.2]{ChenHou2023a_supp},
for $\max(x, y) \geq R_m = 2^4 \cdot 13$, 
we get 
\[
\bal
  \hat f(\om) & =  C_{f0}  I_m(\om) ,
 &  \appo{f}(\om) &= 
 C_{f0} ( I_m(\om)- u_x(\om)(0) ),
 \quad f_A = f - \hat f = f - C_{f0}  I_m(\om), 
 \eal
\]
for  $ f \in \{ u, v, u_x, u_y, v_x \}$,
where $R_m = 13 \cdot 2^{4}$, $I_m(\om) 
= - \f{4}{\pi}\int_{ \max(y_1, y_2) \geq R_m} \f{y_1 y_2}{|y|^4} \om(y) dy $ and
\beq\label{eq:u_appr_near0_coe}
C_{u0} = x , \quad C_{v0} = -y, \quad C_{u_x 0} = 1, \quad C_{u_y 0} = C_{v_x 0} = 0 .
\eeq
Bounding  $I_m(\om_1) - u_x(\om_1)(0), u_x(\om_1)(0)$ by energy $E_4$ \eqref{energy4},
and bounding $I_m(\e_1) - u_x(\e_1)(0), u_x(\e_1)(0)$ using piecewise bounds 
of error $\e = \bar \e, \hat \e$, we bound $\UU_{\FM}, \wh{\UU}$ in \eqref{eq:U}. 

We use the $L^{\inf}$ estimate of $\uu_A, (\na \uu)_A$ established in supplementary material for Part II \cite[Section {\seconeDlinffar}]{ChenHou2023bSupp_supp}. Applying  
the estimate of $\uu_A( g), (\na \uu)_A(g)$ with $g = \om_1, \bar \e_1, \hat \e_1$, 
we bound $\UU_A, ( \na \UU)_A$. 
Combining estimates for $\UU_{\FM}, \UU_A$ and using \eqref{eq:U}, we bound $\UU, \na \UU, \td \UU, \na \td \UU$.

Next, we estimate linear nonlocal terms in $B_{\mf{modi},i}$ in \eqref{eq:bad_lin}. We focus on the typical terms $ U_{x, A} \bar \th_x, U_A \bar \th_{xx}$ in the $\eta_1$ equation 
in $B_{\mf{modi},2}$ \eqref{eq:bad_lin}, \eqref{eq:lin_U}. Using \eqref{eq:lin_asym_ASS2}, \eqref{eq:lin_asym_linfwg1}, and the estimate in the paragraph above \eqref{eq:lin_asym_F_ratio}, we have 
\[
 |\bar \th_x |x_1|^{-1/2} P_{1, i} | \leq C_1(\b) r^{2 \bar \al_1 + \al_{1, i} - 1/2},
 \    |\bar \th_x P_{2, i} | \leq C_2(\b) r^{2 \bar \al_1 + \al_{2, i}},
 \  |\bar \th_{xx} P_{j, i}| \leq C_{2 + j}(\b) r^{2 \bar \al_1 -1 + \al_{j, i}},
\]
for some piecewise constant functions $C_i(\b)$. Then we apply the estimates in Section {\seconeDlinffar} in the  supplementary material of Part II \cite{ChenHou2023bSupp_supp} to estimate 
\[
U_{x, A} r^{ \al}, \quad \al = 2 \bar \al_1 + \al_{1, i} - \tf12, \ 2 \bar \al_1 + \al_{2, i}
\quad 
 U_A |x_1|^{-1/2} r^{ 2 \bar \al_1 - 1 + \al_{1,i}}, \quad U_A r^{2 \bar \al_1 -1 + \al_{2, i}}. 
\]

The estimates of other nonlocal terms, e.g. $\UU_A, (\na \UU)_A$ \eqref{eq:U} \eqref{eq:lin_U} are similar. We remark that the far field estimate of these nonlocal terms are bounded. For example, for the weights $\vp_i$ and  large $|x|$, we have 
\[
\vp_1 \sim |x|^{-1/6}, \quad  | U_{x, A} | \les \log |x| \cdot |x|^{1/6},  \quad  2 \bar \al_1  + \al_{2, i} 
= 2 \bar \al_1  + \tf{1}{7} \approx - \tf23 + \tf17,   \ i \in \{2, 3 \},
\]
and the estimate for $U_{x, A} r^{2 \bar \al_1 + 1/7 }$ decays for large $|x|$. Similar reasoning applies to other powers.

\vs{0.1in}
\paragraph{\bf{Estimate of local terms}}
The estimate of the local terms in \eqref{eq:lin_U}, e.g. $ \eta_1$ in the $\om_1$ equation, is trivial. We just need to bound $ \eta_1 \vp_1 \leq || \eta_1 \vp_2||_{\inf} | \vp_1 / \vp_2 | $ and the ratio $\vp_1 / \vp_2$. Since $\vp_2$ decays slower, $\vp_1 / \vp_2$ is bounded. 
More generally, to bound $P = |x|^{-1/2} P_1 + P_2, Q = |x|^{-1/2} Q_1  + Q_2$ with $P_i, Q_i$ radial weights and $P_i \les Q_i$, we use 
\beq\label{eq:lin_asym_linfwg2}
P / Q \leq \max( P_1/ Q_1, P_2/ Q_2),
\eeq
and the estimate in \eqref{eq:lin_asym_wg1}. The estimates of  $\bar v^N_x \eta_1, \bar u^N_y \xi_1$ in
$B_{\mf{modi}}$ \eqref{eq:lin_U}, \eqref{eq:bad_lin} follow the above estimate and \eqref{eq:lin_asym_ASS2} for $\bar v^N_x$. 
We also treat $-\bar u^N_x \eta_1,  \bar u^N_x \xi_1$ from the damping terms $d_{2, L}^{\numm} \eta_1,  d_{3, L}^{\numm} \xi_1$ in \eqref{eq:lin_U}, \eqref{eq:dp} as a perturbation using the same method.

\vs{0.1in}
\paragraph{\bf{Estimate of the damping terms $d_{i, L}^{\numm}, \cT_{d, N}$}}

We combine the estimates of linear damping term $ d_{i, L}^{\numm}(\vp_i) W_{1, i} \vp_i$ \eqref{eq:lin_U} 
and nonlinear term $ \cT_{d, N}(\vp_i) W_{1, i} \vp_i$ in \eqref{def:non_wg}.
We first estimate $\cT_{d}^{\numm}$ in \eqref{eq:dp}.
The main part is given by the $\bar c_l x  $ term since $\bar \uu^N$ grows sublinearly for large $|x|$.

As in \eqref{eq:wg_nota}, we have $\vp_i = x^{-1/2} P(r) + Q(r)$ 
for some radial weights $P, Q$.  Since $r \pa_r P, r \pa_r Q$ are radial weights, using \eqref{eq:wg_id1}, \eqref{eq:lin_asym_linfwg2}, we obtain 
\bseq\label{eq:lin_asym_linfwg3}
\beq
\bal
 & \f{x \cdot \na \vp}{\vp} \leq \f{ x^{-1/2} r \pa_r P + r \pa_r Q }{ x^{-1/2} P + Q }
 \leq \max( \f{r \pa_r P}{ P}, \f{r \pa_r Q}{Q} ),  \quad
|\f{ \pa_y \vp}{\vp}| \leq \max( |\f{ \pa_y P}{ P} | , |\f{ \pa_y Q}{ Q} | ), 
  \\
& |\f{x \pa_x \vp}{\vp}| 
\leq \max( \f{ | - 1/ 2 \cdot x^{-1/2} P + x^{1/2} \pa_x P| }{ x^{-1/2} P} 
, \f{ |x\pa_x Q| }{Q}  )
\leq \max(1/2 +   \f{ |x \pa_x P| }{  P} 
, \f{ |x\pa_x Q| }{Q}  ).
 \eal
\eeq
For radial weight $P(x) = \sum_{i\leq n} p_i r^{a_i}
= r^{a_n} \sum_{i\leq n} p_i r^{a_i - a_n} \teq r^{a_n} P_{m,n}(x)$ with $a_i$ increasing in $i$, since $a_i - a_n \leq 0$, $P_{m,n}(r)$ is decreasing in $r$, we have
\beq
\f{\pa_z H}{ H} = \f{z}{r} \cdot \f{\pa_r H}{H},
\quad 
 \f{r \pa_r P}{  P} = \f{r \pa_r r^{a_n}}{ r^{a_n}} + \f{r \pa_r P_m}{ P_m}
  = a_n +  \f{r \pa_r P_{m,n}}{ P_{m,n} } \leq a_n. 
\eeq
\eseq
for $ z \in \{x, y\}, H \in \{ P, Q\}$. The second term 
$ \f{r \pa_r P_{m,n}}{ P_{m,n} }$ has a faster decay, and $a_n$ is the main term.
Using \eqref{eq:lin_asym_ASS2}, \eqref{eq:lin_asym_linfwg1}, we derive the asymptotic estimates of $\f{\bar u^N}{x}$ and $\f{ \bar v^N}{r}$, and 
further use \eqref{eq:lin_asym_linfwg3} to  estimate $ \f{\bar \uu^N \cdot \na \vp_i}{\vp_i}$ :
\beq\label{eq:far_U_damp}
|\f{\bar \uu^N \cdot \na \vp_i}{\vp_i} | \leq |\f{\bar u^N}{x} | \cdot |\f{x \pa_x \vp_i}{\vp_i}|
+ \f{ | \bar v^N| }{r } \cdot r |\f{  \pa_y \vp_i}{\vp_i}|.
\eeq
We do not estimate $ \f{ y\pa_y \vp_i}{\vp_i}, \f{\bar v^N}{y}$ since  $\vp_i$ is not singular on $y= 0$. 
Recall $d_{loc,i}$ in $d_{i, L}^{\numm}$ from \eqref{eq:dp}
\[
 d_{loc, 1}  = \bar c^N_{\om}, \qquad  \qquad \qquad \ \, \quad  d_{loc, 2} = 2 \bar c_{\om}^N - \bar u^N_x,  \qquad  \qquad  \    d_{loc, 3} = 2 \bar c^N_{\om} + \bar u^N_x .
\]
The $\bar c_{\om}^N$ term in $d_{i, L}^{\numm}$ \eqref{eq:dp} is a damping term.  For $i=2,3$, since the local terms $\bar u^N_x$ or $-\bar u^N_x$ in $d_{loc, i}$ decay, we treat $|\bar u^N_x|$ perturbatively using \eqref{eq:lin_asym_ASS2} instead of treating 
it as a damping coefficient.
We keep the following terms as damping terms in this estimate 
\[
( (\bar c_l x \cdot \na \vp_i) / \vp_i + c_i \bar c^N_{\om} ) W_{1, i} \vp_i, \quad c = (1, 2, 2) ,
\]
whose coefficient does not decay for large $|x|$. We further bound $ \f{ \bar c_l x \cdot \na \vp_i}{\vp_i}$ from above using the first estimate in \eqref{eq:lin_asym_linfwg3}.

As in \eqref{eq:far_U_damp}, for $\cT_{d, N} W_{1, i} \vp_i$ in \eqref{eq:dp}, \eqref{eq:lin_U}, we bound 
\[
 |\cT_{d, N}(\vp_i) | \leq  \f{|U|}{|x|}  \cdot |\f{x \pa_x \vp_i}{\vp_i}|
+ \f{ | V | }{r} \cdot r |\f{\pa_y \vp_i}{\vp_i}| .
\]
We bound $\f{U}{x}, \f{|V|}{r}$ using their far-field estimates discussed in \hyperref[sec:far_est_nloc]{Nonlocal estimates},  and the estimates of weights in \eqref{eq:lin_asym_linfwg3}. 
Combining the  estimates for $d_{i, L}^{\numm}$ and $ \cT_{d, N}(\vp_i) $,  we derive lower bound for 
\[
- d_{i, L}^{\numm} -  |\cT_{d, N}(\vp_i) | ,
\]
which are used for the left hand side of stability condition. See \cite[Appendix D.1.1]{ChenHou2023a_supp}.

\vs{0.05in}

\paragraph{\bf Estimate of nonlinear terms}

For $\max(x, y) \geq y_{N_2-1}$, we have $\hat W_2 =0$ \eqref{eq:W2_supp}.
It follows 
\beq\label{eq:far_linf_non}
\td \cN_i( \rho ) = \cT_{d, N}(\rho) \cdot W_{1,i} \rho  + N_{W_1, i} \cdot \rho 
\eeq
for any weight $\rho$. 
We  estimated the local term $ \cT_{d, N}(\rho) \cdot W_{1,i} \rho$ in the above paragraph.
The $L^{\infty}$ estimates of $ N_{W_1, i}  \rho $ follows from  product estimates,
and the above local, nonlocal estimates.

Similar estimates generalize to weighted $L^{\infty}(\vp_{g, i}), L^{\infty}(\vp_4)$ estimates.

\subsubsection{Estimate the H\"older weights}\label{sec:hol_wg_est}

The H\"older weights $g = g_i$ is $-1/2$ homogeneous. 
Taking derivative on $\lam$, we yield 
\[
(x \cdot \na_x  g_i) (\lam x_1, \lam x_2) = 
\pa_{\lam} g_i(\lam x_1, \lam x_2) =\pa_{\lam} ( \lam^{-1/2} g_i( x_1, x_2) )
 = -\tf{1}{2} \lam^{-3/2}  g_i( x_1, x_2) .
\]
Choosing $\lam = 1$, we yield the following useful identity for the damping term $d_{g_i, \UU}$ \eqref{eq:dp}
\beq\label{eq:damp_id}
x \cdot \na g_i(x) =  -\, \tf{1}{2} g_i(x).
\eeq

In our energy estimates, we estimate several $0-$homogeneous functions related to $g$ for $h_i \geq 0$
\[
f(h) = h_k^{1/2} g(h),\  h_k \f{ (\pa_j g) }{g}(h),  \ |h| \f{ (\pa_j g) }{g} (h), \ k , j = 1,2,
\quad \f{ g_{i_1}(h)}{g_{i_2}(h)}, 1\leq i_1, i_2 \leq 3.
\]

Since $f(h) = f( \f{h_1}{h_2} , 1)$ for $h_2 \neq 0$ and $f(h ) = f(1, \f{ h_2}{h_1}), h_1 \neq 0$, 
 we  estimate it by partitioning $(h_1, h_2) \in [0, 1] \times \{1 \} , \{1 \} \times [0, 1] $ and using the monotonicity of $g, \pa_j g$. 
 From $g = g_i$ \eqref{wg:hol}, 
\[
g(s) = \f{1}{ A_1(s)^{1/2} + A_2(s)^{1/2}}, \  \pa_j g = - \f{1}{2} \f{  a_{1 j } A_1^{-1/2} + a_{2 j} A_2^{-1/2}   }{ (A_1^{1/2} + A_2^{1/2} )^2 }, \  \f{\pa_j g}{g}
  = - \f{1}{2} \f{  a_{1 j} A_1^{-1/2} + a_{2 j} A_2^{-1/2}   }{ A_1^{1/2} + A_2^{1/2}  } ,
\]
for $A_m = a_{m 1} s_1 + a_{m 2} s_2$ with $a_{m j } >0$. Clearly, $g(s)$ is decreasing in $|s_i|$. For $s_1, s_2 > 0$, since $A_m$ is increasing in $s_1, s_2$, 
$\pa_j g, \f{\pa_j g}{g}$ are negative and increasing in $s_1, s_2$. It follows that $| \f{\pa_j g}{g}|$ is decreasing in $s_1, s_2$.

\subsection{$C^{1/2}$ stability analysis in the far-field}\label{sec:hol_far}

We use the $C^{1/2}$ estimate of $\uu_A, (\na \uu)_A$ in the far-field established in supplementary material of Part II \cite[Section {\secuoneD}]{ChenHou2023bSupp_supp}. 
As in Section \ref{sec:non_dec_u},  we combine the nonlocal estimates for the perturbation $\om_1$ and error 
$\bar \e_1$ and $\hat \e_1$
to obtain nonlocal estimates for $\UU_A, (\na \UU)_A$.

For the $C^{1/2}$ estimate of a typical nonlocal terms in 
$B_{\mf{modi},i},  \td \cN_i$  \eqref{eq:lin_U}, 
\eqref{eq:bad_lin}, \eqref{def:non_wg} in the far-field, e.g. $U_A \bar \om_x \psi_1$, we perform the following decomposition 
\beq\label{eq:lin_asym_uwbx}
U_A \bar \om_x \psi_1 = U_A \psi_u \cdot  ( \bar \om_x \psi_1 / \psi_u).
\eeq

We need the asymptotic $L^{\inf}$ and $C^{1/2}$ estimates of $U_A \psi_u$ and $\bar \om_x \f{\psi_1}{ \psi_u }$. More generally, we need asymptotic estimate of $\UU_A \psi_u, (\na \UU)_A \psi_1$ and weighted profile $\bar S \f{f}{g} $ for two weights $f, g$.

From  the supplementary material of Part II \cite[Section {\secuoneD}]{ChenHou2023bSupp_supp}, for $|z| \geq |x|_{\inf} \geq y_{N_2-1}$, we have the weighted $C_{x_k}^{1/2}$ and $L^{\inf}$ estimates of $\uu_A(x ), \na \uu_A(x)$ uniformly:  
\beq\label{eq:lin_asym_u}
\bal
  \d_k( (\pa_1^i \pa_2^j \phi)_A \psi_u, x, z )  & \leq \CCs_{u,  ij, \mf{hol}, k}, \ i + j =1,  & 
\quad \d_k( (\pa_1^i \pa_2^j \phi)_A \psi_1, x, z ) & \leq \CCs_{u,  ij, \mf{hol}, k} ,  \ i + j = 2,  \\
  \f{ | (\pa_1^i \pa_2^j \phi)_A \psi_u(x)|  }{ |x|^{1/2} }
& \leq \CCs_{u,ij}  , \ i+ j = 1,  & 
\quad  \f{ | (\pa_1^i \pa_2^j \phi)_A \psi_1(x)|  }{ |x|^{1/2} } 
& \leq \CCs_{u,ij} ,  \ i+j=2,
 \eal
\eeq
for some computable constants $\CCs_{u,  ij, \mf{hol}, k}, \CCs_{u,ij}$. Here, we use $(\pa_1^i \pa_2^j \phi)_A$ to denote $\UU_A, (\na \UU)_A$: 
\[
-  (\pa_2 \phi)_A  = U_A, \quad  (\pa_1 \phi)_A  = V_A,
\quad -(\pa_{12} \phi)_A = U_{x, A},
\quad (\pa_{11} \phi)_A = V_{x, A},
\quad -(\pa_{22 } \phi)_A = U_{y, A} .
\quad 
\]
which are consistent with relation between velocity and stream function, e.g. $ u = -\pa_2 \phi$.

\subsubsection{Far-field H\"older product estimates}\label{sec:lin_asym_prod}

Firstly, we estimate $[AB]_{C_{x_i}^{1/2}}$ with 
\[
|\pa_x^i \pa_y^j A(r, \b)| \leq a_{ij} r^{ p - i - j}, \quad |\pa_x^i \pa_y^j B(r, \b)| \leq b_{ij} r^{ q - i - j}, \quad i + j \leq 1 ,
\]
for all $|x| \geq R$.  Using the Leibniz rule, for $i+j \leq 1$, we get 
\[
C = AB,  \quad  |\pa_x^i \pa_y^j C| \leq c_{ij} r^{ p + q - i - j}, \quad
c_{ij} = \sum\nolim_{ i_1 \leq i , j_1 \leq j }  a_{i_1 j_1} b_{i-i_1, j-j_1} ,
\]
for $|x| \geq R$, where the binomial coefficients satisfy $ \binom{ i }{i_1} =  \binom{j}{j_1 } = 1 $ since $i, j \leq 1$.

In the $C_{x_1}^{1/2}$ estimate, for $ R \leq |x| \leq |y| $ and $p + q \leq 0$, 
we have 
\[
|C(x) - C(y)| \leq 2  c_{00}|x|^{ p + q }, \quad |C(x) - C(y)| \leq  |\pa_1 C(\xi)| |x- y|
\leq c_{10} |x|^{ p+ q -1}|x-y|
\]
for some $|\xi| \geq |x|$. It follows 
\beq\label{eq:lin_asym_prod}
\f{|C(x)- C(y)|}{|x-y|^{1/2}} \leq \min ( \f{ c_{00} 2 |x|^{ p + q }  }{|x-y|^{1/2}},
 c_{10}  |x|^{ p + q - 1}  |x-y|^{1/2} ) 
\leq \sqrt{2 c_{00} c_{10}} |x|^{ p + q -1/2}.
\eeq
The $C_{x_2}^{1/2}$ estimate is similar. Using \eqref{eq:lin_asym_ASS2},\eqref{eq:lin_asym_wg2},\eqref{eq:lin_asym_wg3}, and the above estimates, we  estimate the asymptotic behavior of the $L^{\inf}$ and $C_{x_i}^{1/2}$ semi-norm of weighted profile, e.g. $\bar \om_x \f{\psi_1}{\psi_2}$ \eqref{eq:lin_asym_uwbx}.

Next, we estimate $ \d_i( u_A F, x, y)$ or $\d_i((\na u)_A F, x, y)$ with
\[
 |\d_i(F, x, y)| \leq f_i r^{\al_i}, \quad |F(x)| \leq f_0 r^{\al_0}, \quad |x|, \ |y| \geq r, 
 \quad \al_0 \leq 0, \ \al_1 + 1/2, \ \al_2 + 1/2  \leq 0.
\]

Recall the estimate of \eqref{eq:lin_asym_u}. We focus on $ U_{x, A} F$. The estimates of other terms are similar. For $|x|, |y| \geq R$, using $ |a(x) b(x) - a(y) b(y)| \leq |(a(x) - a(y)) b(x)| + |a(y) (b(x) - b(y))| $,  we obtain 
\[
\bal
| \d_i( U_{x, A} F, x, y )|
&\leq  |\d_i( U_{x, A} , x, y ) F(y)| 
+| U_{x, A}(x)   \d_i( F, x, y ) |
\leq \CCs_{u,   11, \mf{hol}, i } f_0 |x|^{\al_0 } 
+ \CCs_{u,  11} |x|^{1/2} f_i |x|^{\al_i}  \\
&\leq \CCs_{u,  11, \mf{hol} , i} f_0 R^{\al_0} + \CCs_{u,  11} f_i R^{ \al_i + 1/2}, \quad i = 1, 2.
\eal
\]
for some computable constants $ \CCs_{u,  11, \mf{hol}, i} , \CCs_{u,  11} $. Using the above estimates, we estimate \eqref{eq:lin_asym_uwbx} and other nonlocal terms in \eqref{eq:lin_U}, \eqref{eq:bad_lin} in the weighted H\"older estimate in the far-field, e.g. $ U_{x, A} \bar \th_x \psi_2, U_A \bar \th_{xx} \psi_2$.

\subsubsection{Estimates of the local terms in $B_{\mf{modi}, i}$ \eqref{eq:bad_lin} }\label{sec:lin_asym_loc}

Recall from \eqref{energy4} that the H\"older weights for $\eta_1, \xi_1$ are the same: $\psi_2 = \psi_3, g_2 = g_3$.   The H\"older estimate of the local terms in $B_{\mf{modi}, i}$ \eqref{eq:bad_lin}, e.g. $\bar v^N_x \xi_1$ in the $\eta_1$ equation, is simple. See \eqref{eq:hol_loc}. We have
 \[
|\d_i ( \bar v^N_x \xi_1  \psi_2, x, y) |
\leq  |\bar v^N_x \d_i( \xi_1 \psi_2, x, y)|
+ |\d_i(\bar v^N_x,x, y) \xi_1 \psi_2|,
\]
and use \eqref{eq:lin_asym_prod} to estimate $\d_i(\bar v^N_x,x, y)$, the energy to control $\d_i(\xi_1 \psi_2, x, y), \xi_1 \psi_2$. Note that to bound $\xi_1 \psi_2$, we use $\vp_3 \geq p_n r^{a_n}$ for some $p_n, a_n$, where $r^{a_n}$ is slowest decay power in $\vp_3$, and 
\[
|\xi_1 \psi_2 | \leq || \xi_1 \vp_3||_{\inf}  \f{\psi_2}{\vp_3} 
\leq E_4  \f{\psi_2}{ p_n r^{a_n}} .
\]
The ratio of the weight is further estimated using \eqref{eq:lin_asym_wg3}.

\subsubsection{Estimates of the damping terms and $\cT_{d,N}$}

We estimate the damping terms $d_{g_i, \UU}$ and $d_{i, L}^{\numm}$ \eqref{eq:dp} in the far-field.
We decompose $d_{g_i,\UU}$ as
\[
d_{g_i,\UU} = d_{g_i, \mf{lin}} + d_{g_i, \mf{non}},
\quad 
d_{g_i,\mf{non} } \teq \f{  (\UU(x) - \UU(z))  \cdot \na g_i(x-z) }{g_i(x-z)}, 
\]
and $d_{g_i, \mf{lin}} $ is defined below. Since $  \bar \uu^N/ |\xx|, \na \bar \uu^N$ decays for large $|\xx|$, using \eqref{eq:lin_asym_ASS2} and \eqref{eq:damp_id}, for $|x|, |z| \geq R$, we bound $ d_{g_i, \mf{lin}}$ from above as 
\[
\bal
d_{g_i, \mf{lin}} & \teq -\f{\bar c_l}{2} + \f{  (\bar \uu^N(x) - \bar \uu^N(z)) \cdot \na g_i(x - z)}{g_i(x-z)} 
 \leq - \f{ \bar c_l}{2} + R^{ \bar \al_1} 
 \B( ||\bar \phi_{ag, 11}||_{\inf} |x_1 - z_1|  \\
 & \quad +  || \bar \phi_{ag, 02}||_{\inf} |x_2 - z_2| ) \f{|\pa_1 g_i(x-z)|}{ g_i(x-z)}   
  + ( || \bar \phi_{ag, 20}||_{\inf} |x_1 - z_1| +  || \bar \phi_{ag, 11}||_{\inf} |x_2 - z_2| ) \f{|\pa_2 g_i(x-z) |}{ g_i(x-z)} \B) ,
\eal
\]
where $ \bar \phi_{ag, 11}, \bar \phi_{ag, 20}, \bar \phi_{ag, 02}$ corresponds to the far-field estimates of $u_x=-v_y,  v_x, u_y$, respectively. Similarly, using the estimate for $\na \UU$, we bound $|d_{g_i,\mf{non}}|$ perturbatively. Then, we estimate $d_{g_i, \UU}$ from above.
Since $ h_k  \f{\pa_l g_i(h)}{ g_i(h)}$ is $0-$homogeneous, we follow Section \ref{sec:hol_wg_est} to estimate it.

For $d_{i, L}^{\numm}$ \eqref{eq:dp}, let $a_{n,i}$ be the last power in the radial weight $\psi_i$. We decompose $\psi_i = r^{ a_{n, i}} \psi_{m, i}$, where $\psi_{m ,i}(x) \asymp 1$ for large $|x|$.
Since $  \f{\pa_j \psi_i}{\psi_i} =  \f{ \pa_j r^{a_{n,i} } }{ r^{ a_{n,i}}}  +  \f{ \pa_j \psi_{m, i}}{ \psi_{m, i}}$, we decompose  $d_{i, L}^{\numm}$ as follows 
\[
d_{i,L}^{\numm}(\psi_i) = \f{\bar c_l x \cdot \na |x|^{ a_{n, i}} }{ |x|^{ a_{n, i}}}
+  \f{\bar c_l x \cdot \na \psi_{m,i} }{ \psi_{m ,i}}
+ \f{ \bar \uu^N \cdot \na \psi_i }{ \psi_i}  + d_{loc, i}
= a_{n, i} \bar c_l + I_{i, 2} + I_{i, 3} +  d_{loc, i} , \\
\]

 Since $ x \cdot \na \psi_{m,i } = r \pa_r \psi_{m ,i} $ and since $\psi_{m,i}$ is a radial weight and the last power in $\psi_{m ,i}$ is $1$, 
$I_{i, 2}$ is the ratio between two radial weights with a denominator decaying faster. We use  \eqref{eq:lin_asym_wg3} and \eqref{eq:lin_asym_prod} to control $I_{i, 2}$. We use \eqref{eq:lin_asym_wg4} and \eqref{eq:lin_asym_ASS2} to estimate $ \bar \uu^N  \cdot \na \psi_i / \psi_i$ and its derivatives. 

Recall $d_{loc,i}$ from \eqref{eq:dp} 
\[
 d_{loc, 1}  = \bar c^N_{\om}, \qquad  \qquad \qquad \ \, \quad  d_{loc, 2} = 2 \bar c_{\om}^N - \bar u^N_x,  \qquad  \qquad  \    d_{loc, 3} = 2 \bar c^N_{\om} + \bar u^N_x .
\]
For $i=2,3$, since the local terms $\bar u^N_x$ or $-\bar u^N_x$ in $d_{loc, i}$ has the same decay as $I_{i, 3}$, we combine it with $I_{i, 3}$ and apply \eqref{eq:lin_asym_prod} for the H\"older estimate and bound it perturbatively. Then, we use \eqref{eq:lin_asym_ASS2}, the above estimates for $I_2, I_3$, and the estimate in Section \ref{sec:lin_asym_loc} to estimate 
\[
 (I_{1, 2} + I_{1,3} ) \om_1 \psi_1, \quad (I_{2,2} + I_{2, 3} - \bar u^N_x) \eta_1 \psi_2,
 \quad (I_{3,2} + I_{3, 3} - \bar v^N_y) \xi_1 \psi_2. 
\]

Since the coefficients in the above terms decay, we treat the above terms as perturbation. 
We keep the constant coefficients in $d_{i,L}^{\numm}$ and 
derive the main damping term in $d_{i,L}^{\numm} W_{1, i} \psi_i$ :
\beq\label{eq:damp_far_main}
d_{i, M}(\psi_i) W_{1, i} \psi_i, \quad  d_{i, M}(\psi_i)  \teq  a_{n, i}  \bar c_l   + c_i \bar c^N_{\om} , \quad  c_1 = 1, \, c_2 = 2, \, c_3 = 2. 
\eeq

\subsubsection{Estimates of nonlinear terms}

For $|x|_{\inf}, |z|_{\infty} \geq y_{N_2-1}$, 
we have $\hat W_2(x) \equiv 0$ and nonlinear terms $\td \cN_i(\psi_i)$ in \eqref{def:non_wg} 
reduces to \eqref{eq:far_linf_non}. 
The nonlinear estimates follows from  product rules. We focus on $\cT_{d, N}(\psi_i) W_{1, i} \psi_i$ 
 in \eqref{def:non_wg}. Recall $\cT_{d, N}(\psi_i)$ \eqref{eq:dp}:
\[
\cT_{d, N}(\psi_i) =  \psi_i^{-1}  (\UU \cdot \na \psi_i ).
\]
 Since we have $C^1$ far-field bounds for $\UU$, we derive the $C^1$ asymptotic bounds for $\UU$ and $ \f{\na \psi_i}{\psi_i}$ and use the product rule to bound $\cT_{d, N}(\psi_i)$.
The weighted $L^{\infty}$ and $C^{1/2}$ norm of $W_{1,i}$ are bounded by energy $E_4$. 
 Then we follow \cite[Section 5.9.3]{ChenHou2023a_supp}, the product estimates in Section \ref{sec:lin_asym_prod} and in \eqref{eq:hol_loc} to bound the nonlinear terms.

\subsection{H\"older estimates with large distance}\label{sec:hol_lg}

In this subsection, we consider  H\"older estimates with large $|x-z|$. 
We rewrite the system \eqref{eq:lin_U}  with weight $\rho_i = \psi_i$ as 
\beq\label{eq:lin_U_lg}
 \pa_t( W_{1, i} \psi_i) + (\bar c_l x + \bar \uu^N + \UU) \cdot \na ( W_{1, i} \psi_i)
   = d_{i, M} W_{1, i} \psi_i + 
   \td B_{ad, i}(\psi_i)
\eeq
where we isolate the main damping terms $d_{i, M}$ which is defined in \eqref{eq:damp_far_main}.
Here, $\td B_{ad ,i}(\cdot)$ denotes bad terms in the weighted estimates \eqref{eq:lin_U}
\[
 B_{ad, i} \teq B_{\mf{modi}, i}  \cdot \psi_i + 
 \td \cN_i(\psi_i)  + \bar \cF_{loc, i} \psi_i - \cR_{loc, i} \psi_i ,
 \quad 
\td  B_{ad, i} \teq  ( d_{i,L}^{\numm} - d_{i, M} )  W_{1, i} \psi_i  + 
  B_{ad, i}.
\]
and we \emph{include} $\psi_i$ in the definitions of $B_{ad, i}, \td  B_{ad,i}$. We have derived the piecewise $L^{\infty}$ bounds of  $B_{ad, i}$ and the damping factor $d_{i, L}^{\numm}( \psi_i)$ \eqref{eq:dp} in the mesh $Q_{i, j}$ \eqref{eq:mesh_xy} in the $C^{1/2}$ estimate of $W_{1, i} \psi_i$ with $|x-z|$ \emph{not very large}. 
For the H\"older estimate with large $|x-z|$, 
by further bounding the explicit coefficient $d_{i,L}^{\numm} - d_{i, M}$, we 
extend the piecewise $L^{\infty}$ bound of $B_{ad,i}$ to $\td  B_{ad,i}$
\[
 \max\nolim_{x \in [0,  \ y_{N-1}]^2} |\td B_{ad, i}   | \leq C_{i, 1} E_4,
\]
where $E_4$ is the energy. Since these terms decay in the far field, for $|z-x|$ large, they are very small and are treated as a small perturbation. Using these estimates and  triangle inequality, for $|x-z| \geq R_1$ with $x \in 
[0, y_{N_2-1}]^2, z \in [0, y_{N-1}]^2$, we  bound 
\beq\label{eq:hol_lg_B2}
( | \td B_{ad, i}(x)|  + | \td B_{ad,i}(z) | ) g_i(x- z)
\leq ( |\td B_{ad, i}(x)| + |\td B_{ad,i}(z) |)  R_1^{-1/2}  |x-z|^{1/2} g_i(x-z).
\eeq

Beyond the mesh $[0, y_{N-1}]^2$, from Section \ref{sec:hol_far}, we can estimate $ \f{| \td B_{ad, i}(z) | }{|z|^{1/2}}$ uniformly for $|z| \geq y_{N-1}$. Since $x \in  [0, y_{N_2-1}]^2$, we get
$|z-x| \geq |z|(1 - |x| / |z|) \geq y_{N-1}( 1 - \sqrt{2} y_{N_2-1} / y_{N-1})$ and
\beq\label{eq:hol_lg_B1}
\bal
 ( | \td B_{ad,i}(x) | +  |\td B_{ad,i}(z) | )  & g_i(x-z)
  \leq (\f{ |\td B_{ad,i}(x) | }{|z|^{1/2}} + \f{ | \td B_{ad,i}(z) | }{|z|^{1/2}} ) \f{|z|^{1/2}}{ |x-z|^{1/2}}
  |x-z|^{1/2} g_i(x-z) \\
 \leq  &  (\f{ | \td B_{ad,i}(x) | }{y_{N-1}^{1/2}}  + \f{ | \td B_{ad,i}(z) | }{|z|^{1/2}} ) \f{|y_{N-1}|^{1/2}}{( y_{N-1}- \sqrt{2} y_{N_2-1} )^{1/2}}
  |x-z|^{1/2} g_i(x-z).
\eal
\eeq
Since $  |h|^{1/2} g_i(h)$ is $0-$homogeneous, we follow Section \ref{sec:hol_wg_est} to estimate it. 

It remains to estimate the damping factor $d_{g_i, \UU }$ from  $g_i$ in \eqref{eq:dp}. We 
gain a damping factor $- \bar c_l / 2$ by \eqref{eq:damp_id} and treat $(\bar \uu^N(x) - \bar \uu^N(z)) \cdot \na g(x-z) / g(x-z)$ as perturbation. We estimate the Lipschitz norm of $\bar \uu^N$ for $|x-z| \geq R_1$. Denote $f = \bar u^N$ or $\bar v^N$.
Suppose $|x|_{\inf} \leq |z|_{\inf}$. We have 
\[
|z|_{\inf} \geq \max_i( |x|_i, |z|_i) \geq \max( |z_1 - x_1|, |z_2 - x_2|) \geq  |x-z| /\sqrt{2} \geq R_1 / \sqrt{2}.
\]

If $|x|_{\inf} \geq  \f{R_1}{2 \sqrt{2}}$, 
we take $w = (x_1, z_2)$ if $|x_1| \geq  \f{R_1}{2 \sqrt{2}}$
and $w= (z_1, x_2)$ if $|x_2| \geq  \f{R_1}{2 \sqrt{2}}$. Let $D_{R} = [0,  \f{R_1}{2 \sqrt{2}}]^2$. 
Since $|z|_{\infty}  \geq  \f{R_1}{2 \sqrt{2}}$, the line $xw$, $wz \notin D_R$. Using mean-value theorem, we get 
\begin{align}\label{eq:mean_value_lip}
|f(x) -f(z)|
 & \leq |f(x) -f(w)  | + |f(w)- f(z)|\leq |x_1 - z_1| \cdot  \| \pa_1 f \|_{L^{\infty}(D_R^c)}
+ |x_2 - z_2|  \cdot \| \pa_2 f \|_{L^{\infty}(D_R^c)} \notag \\
& \leq |x-z| ( \| \pa_1 f \|_{L^{\infty}(D_R^c)}^2 + \| \pa_2 f \|_{L^{\infty}(D_R^c)}^2)^{1/2}.
\end{align}
If $|x|_{\inf} < \f{ R_1}{  2 \sqrt{2}}$, we get 
\[
\bal
|x| \leq \sqrt{2} |x|_{\inf} \leq \f{R_1}{2} \leq \f{|x-z|}{2} , \quad |z| \leq |x-z| + |x| 
\leq  2  |x-z|  , \
\f{|f(x)-f(z)|}{|x-z|} 
\leq \f{|f(x)|}{R_1} + \f{ 2 |f(z)|}{ |z|} .
\eal
\]

By bounding $f(x)$ in $[0, R_1 / (2\sqrt{2})]^2$, and $f(z)/|z|$ outside $[0, R_1 / \sqrt{2}]^2$, we obtain the Lipschitz bound. Since $ \na \bar \uu^N$ and $\bar \uu^N(x) / |x|$ decay for large $|x|$, the above bounds are very small for large $R_1$. Denote $h_i = x_i - z_i, f = (\bar u^N + U, \bar v^N + V)$. With these estimates, we estimate 
\beq\label{eq:damp_g_pertb}
|\f{( f(x) - f(z)) \cdot \na g_i }{g_i}|
\leq  \f{|( f_1(x) - f_1(z)|}{ |h|} \cdot  | h | \B|\f{ \pa_1 g_i}{g_i}(h) \B|
+  \f{ | ( f_2(x) - f_2(z)|}{|h|} \cdot |h | \B|\f{ \pa_2 g_i}{g_i}(h) \B|.
\eeq
 Since $ |h|  \f{\pa_j g_i(h)}{g_i(h)}$ is $0-$homogeneous, we follow Section \ref{sec:hol_wg_est} to estimate it.  We treat the above terms as perturbations and we keep the following terms as damping terms 
\[
(d_{i, M} - \bar c_l / 2) \d( W_{1, i} \psi_i) g_i(x-z),  \quad   (d_{i, M} - \bar c_l / 2)< 0.
\]

\subsubsection{H\"older estimate in case (a.3)}\label{sec:hol_case_a3}
For $x ,z \in [0, y_{N-1}]^2 \bsh [0, y_{N_1 -2}]^2$ with $x \in Q_{ij}, z\in Q_{pq}$, 
$\max(|i-p| , |j-q|) \geq 2$, 
the H\"older estimate is similar since $|x|,|z|$ are very large and $\td B_{ad, i}(\cdot)$ is small. 
We use $y_0=0$, $|x-z| \geq \min_{i\geq 0} (y_{i+1} - y_i) = y_1$  and estimate 
$\td B_{ad ,i}$ following \eqref{eq:hol_lg_B2}:
\[
( | \td B_{ad, i}(x)|  + | \td B_{ad,i}(z) | ) g_i(x- z)
\leq ( |\td B_{ad, i}(x)| + |\td B_{ad,i}(z) |)  y_1^{-1/2} \cdot  |x-z|^{1/2} g_i(x-z).
\]

We do not need to use \eqref{eq:hol_lg_B1}. In the H\"older esitmate, we have 
a damping term $d_{i, M} \d( W_{1, i} \psi_i, x, z) g_i(x- z)$ as in \eqref{eq:lin_U_lg}, where $d_{i, M}$ is defined in \eqref{eq:damp_far_main}.  
For damping factor $d_{g_i, \UU }$ 
from  $g_i$   \eqref{eq:dp}, we gain an extra damping factor $- \f{\bar c_l}{2}$ by \eqref{eq:damp_id} 
and apply \eqref{eq:damp_g_pertb} 
to bound $(f(x) - f(z) ) \cdot \na g(x-z) / g(x-z)$ perturbatively, 
where $f = \bar \uu^N + \UU \in \R^2$.
Since $x, z \in [0, y_{N-1}]^2 \bsh [0, y_{N_1 -2}]^2$, we estimate the Lipschitz norms 
of $f_1, f_2$ using \eqref{eq:mean_value_lip} with domain $D_R^c$ replaced by  $([0, y_{N_1 -2}]^2)^c$.

\subsection{Piecewise bounds}\label{sec:bd_other}

In this section, we assume that the points, e.g. $x,y, z$, are in $\R^2_{++}$. 

\vs{-0.07in}
\subsubsection{Piecewise bound of integral with power}\label{sec:est_f_dx}

Given the 1D mesh $0 = x_1 <.. < x_n$ and the upper and lower bounds of in $I_i = [x_{i}, x_{i+1}]$, we estimate 
\[
I_k(f) =\f{1}{x^k} \int_0^{x} f(s) s^{k-1 } d s, \ k = 1, 2.
\]

We consider the lower bound, and the upper bound is similar. For $x \in I_i= [x_{i}, x_{i+1}]$, we get 
\[ 
I_k(f) = \f{1}{x^k}( \int_0^{x_{i}}  + \int_{x_{i}}^{x }  ) f(s)  s^{k-1} ds
\geq \f{1}{ k x^k} (  \sum_{j \leq i-1} f_{l, j} ( x_{j+1}^k - x_j^k )
+ (x^k - x_i^k) f_{l, i} )  .
\] 

Denote $S_{k, i} =  \sum_{j \leq i-1} f_{l, j} ( x_{j+1}^k - x_j^k )$. Since $x^k \in [x_i^k, x_{i+1}^k]$, we yield 
\[
I_k(f) 
 \geq \f{1}{k} f_{l, i} + \f{1}{k x^k} ( S_{k, i} - x_i^k f_{l, i})
 \geq  \f{1}{k} f_{l, i} +\f{1}{k} \min( \f{   S_{k, i} - x_i^k f_{l, i} }{x_i^k}, \f{  S_{k, i} - x_i^k f_{l, i} }{x_{i+1}^k } ).
 \]

Using a similar argument, for $f(0) =0$ and $x \in [x_k, x_{k+1} ]$, we estimate $f(x) / x$ as follows 
\[
\f{f}{x} = \f{1}{x} ( f(x_k) +\int_{x_k}^x f_x(z) dz )
\geq \f{1}{x} ( f(x_k) + (x- x_k) f_{x, l, k } )
\geq \min ( \f{f(x_k)}{x_k}, \f{f(x_k) + (x_{k+1} - x_k) f_{x, l, k}  }{ x_{k+1}}  ),
\]
where $f_{x,l,k}$ is the lower bound of $f_x$ in $[x_k, x_{k+1}]$. 

To estimate 2D functions, we first denote by $f^u_{i,j}, f^l_{i, j}$ the upper and lower bound of $f$ in a grid cell  $Q_{ij} = [x_i, x_{i+1}] \times [y_j ,y_{j+1}]$. Then for $f(0, y) =0$ with constant power and $(x, y) \in Q_{ij}  $, using the Mean-Value Theorem and $\f{A x + B}{x} \geq A+ \min( B / x_i, B/ x_{i+1})$, we have 
\beq\label{eq:f_dx_2D}
 \f{f(x, y)}{x} \geq \f{1}{x} ( f(x_i, y)  + f^l_{x, ij} (x - x_i) )
 \geq \min( \f{1}{x_i} f(x_i, y),
 \f{1}{x_{i+1}} (f(x_i, y) + f^l_{x, ij} (x_{i+1} - x_i)  ) ) .
\eeq
For $f(x_i, y)$, we further bound it from below using 
\[
f(x_i, y ) \geq  \min( f(x_i, y_j), f( x_i, y_{j+1}) )
- (y_{j+1} - y_j)^2   || f_{yy}||_{L^{\inf}(Q_{ij})} / 8.
\]
In the first interval $x \in [x_1, x_2]$ with $x_1 = 0$, we simply use 
\[
f(x, y) / x \geq f^l_{x,1 j} = \min\nolim_{(x, y)\in Q_{1j}} f_x(x, y).
\]
The estimates of the upper bound of $\f{f(x, y)}{x} $ 
and the bounds of $\f{f(x, y) }{y}$ with $f(x, 0) = 0$ are similar.

\subsubsection{H\"older estimate of the fractional power}\label{sec:dp_hol_singu}
We estimate $\d_i(f, x, z), f = x_2^2 / (x_1^2 + x_2^2)$. Denote 
\beq\label{eq:hol_sinb}
\sin \b = \max( x_2 |x|^{-1}, z_2 |z|^{-1}). 
\eeq
\paragraph{\bf{$C_x^{1/2}$ estimate of $f$ }}
In this case, since $x_2 = z_2$, we write $(x_1, x_2) = (x, y), (z_1, z_2) = (z, y), x \leq z$. Then $\sin \b = y / (x^2 + y^2)^{1/2}$. We have 
\[
\bal
&|\f{y^2}{ x^2 + y^2} - \f{y^2}{z^2 + y^2}| \f{|w|^{1/2}}{|x-z|^{1/2}}
 = \f{ y^2(x+z)|z-x|^{1/2}}{ (x^2 + y^2)(y^2 + z^2)} |w|^{1/2}
 \teq  M M_2, \\ 
 &  M = \f{y^2 |w|^{1/2}}{(x^2 + y^2) (y^2 +z^2)^{1/4}} , \quad M_2 = \f{ (x+z) |x-z|^{1/2}}{ (y^2 +z^2)^{3/4}}.
\eal
\]

For $w = (x, y)$ or $(z, y)$, we estimate $M$ using \eqref{eq:cos_sin}. Moreover, since $|w| \leq |(y, z)|$, we get
\beq\label{eq:dp_sin_cx}
M \leq y^2 / (x^2 + y^2) \leq ( \sin \b)^2,  \quad  \d_1(f, x, z) |w|^{1/2} \leq M_2 (\sin \b)^2.
\eeq

Suppose that $a \in [a_l, a_u]$ for $a = x, y, z$. To estimate the upper bound of $M_2$, we write 
$(x+z) |z-x|^{1/2} = |z^2 - x^2|^{1/2} |x+z|^{1/2}$ and have 
\[
M_2 = (\f{ z^2 - x^2}{ z^2 + y^2})^{1/2} ( \f{x}{ (y^2 + z^2)^{1/2} }  + \f{z}{ (y^2 + z^2)^{1/2}}   )^{1/2} \teq A^{1/2} ( B + C)^{1/2}.
\]
To estimate $A$, we use the fact that $A$ is increasing in $z^2$, and decreasing in $x^2, y^2$. For $C$, we use \eqref{eq:cos_sin}. For $B$, since $x \leq z$, using \eqref{eq:cos_sin}, we get 
\[
B \leq \f{x}{ (y^2 + x^2)^{1/2}} \leq  \f{x_u}{ (y_l^2 + x_u^2)^{1/2}}, \
B \leq \f{x_u}{ (y_l^2 + z_l^2)^{1/2}},
\]
and optimize two estimates. We have another estimate for $M_2$. Denote $\lam = x / z \in [0, 1]$. 
Using Young's inequality 
\beq\label{eq:dp_sin_AMGM}
(1-\lam) (1 + \lam)^2 = \f{1}{2} (2 - 2 \lam)(1 + \lam)^2
\leq \f{1}{2}  (\f{ (2 - 2 \lam) + 2 (1 + \lam) }{3})^3
= \f{1}{2} (\f{4}{3})^3 = \f{32}{27} , \quad c = ( \f{32}{27})^{1/2},
\eeq
for $\lam \in [0, 1]$, we get
\[
M_2 = \f{z^{3/2}}{ (y^2 + z^2)^{3/4} } \cdot \f{ (x+z)(z-x)^{1/2}}{ z^{3/2}} =
 \f{z^{3/2}}{ (y^2 + z^2)^{3/4} }  (1 + \lam)( 1 -\lam)^{1/2} 
\leq c \f{z^{3/2}}{ (y^2 + z^2)^{3/4} },
\]
and estimate the factor using \eqref{eq:cos_sin}.

\vs{0.1in}
\paragraph{\bf{$C_y^{1/2}$ estimate }}
In this case, we get $x_1 = z_1$. Since $f(x_2, x_1) = 1 - f(x_1, x_2), |(a_2, a_1)| =| (a_1, a_2)|$ 
\[
|f(x_1, x_2) - f(z_1, z_2)|  |x-z|^{-1/2} |w|^{1/2}= | f(x_2, x_1) - f(z_2, z_1)| \cdot |x-z|^{-1/2} |w|^{1/2}, 
\]
 the estimate of $ \d_2(f, x, z)$ follows from the above $C_x^{1/2}$ estimate. Note that we need to swap the variables $x_1, x_2$.  To bound $\d_2(f, x, z)$ using $\sin \b$ \eqref{eq:hol_sinb}, we first assume $x_2 \leq z_2$. Then we get $\sin \b = \f{z_2}{ (x_1^2 + z_2^2)^{1/2}}$. Since $z_1 = x_1$ and $|w|\leq |z|$,  we have
\beq\label{eq:dp_sin_cy}
\bal
 I &= \B|  \f{x_2^2}{x_1^2 + x_2^2} - \f{z_2^2}{x_1^2 + z_2^2 } \B| \f{|w|^{1/2}}{|x-z|^{1/2}}
 = \f{x_1^2}{x_1^2 + x_2^2} \f{(z_2-x_2)^{1/2}(z_2+x_2)}{ x_1^2 + z_2^2} |w|^{1/2} \\
&\leq  \f{x_1^2}{x_1^2 + x_2^2}  (\f{z_2}{|z|})^{3/2}
 \f{(z_2- x_2)^{1/2}(z_2 + x_2)}{ z_2^{3/2}} 
=  \f{x_1^2}{x_1^2 + x_2^2}
(\sin \b)^{3/2} (1 - \lam)^{1/2}(1+\lam )  , 
\eal
\eeq
where  $  \lam = \f{x_2}{z_2}$.  We estimate $ \f{x_1^2}{x_1^2 + x_2^2}$ using \eqref{eq:cos_sin}, and $(1- \lam)^{1/2}(\lam + 1) $ using the monotonicity of $1-\lam, 1 + \lam$, piecewise bounds of $\lam \leq 1$, and \eqref{eq:dp_sin_AMGM}.

\subsubsection{Piecewise H\"older estimates of power}\label{sec:hol_pow}

We estimate $\d_i(f, x, y)$ for $f(x) = |x|^p, p= 1/2, 1$. Since $f(x_1, x_2) = f(x_2, x_1)$, we consider $C_x^{1/2}$ estimate. The $C_y^{1/2}$ estimate can be obtained by swapping the variables $x_1, x_2$. Without loss of generality, we assume $y_1 \geq x_1, x_2 = y_2$. For each variable $a = x_1, x_2, y_1$, we assume $a \in [a^l, a^u]$.

Using $|y|^2 - |x|^2 = y_1^2 - x_1^2 = (y_1 - x_1)(y_1 + x_1)$ and  a direct computation yields 
\beq\label{eq:hol_pow_half}
\bal
 & \f{ |y|^{1/2} - |x|^{1/2}}{ |x-y|^{1/2}}
  = \f{ |y|^2 - |x|^2}{ (|y| + |x|) (|y|^{1/2} + |x|^{1/2}) |y-x|^{1/2} }
   = M M_1, 
   \ M \teq \f{y_1 + x_1}{ |x| + |y| }, \ M_1 = \f{ |y_1 - x_1|^{1/2}}{ |x|^{1/2} + |y|^{1/2}}, \\
 & \f{ |y| - |x|}{ |x-y|^{1/2}}
 = \f{ |y|^2 - |x|^2}{ |x-y|^{1/2}(|x| + |y|)}
 = \f{ |y_1-x_1|^{1/2} (x_1 + y_1)}{|x| + |y|} = |y_1-x_1|^{1/2} M
 \leq  (y_1^u - x_1^l)^{1/2} M .
  \eal
\eeq

Clearly, $M(x_1, y_1, x_2)$ is decreasing in $x_2$. A direct computation yields 
\[
\pa_{x_1} M = \f{ |x| + |y| - (x_1 + y_1) \f{x_1}{|x|} }{ (|x| + |y|)^2 } 
\geq \f{ |x| + |y| - (x_1 + y_1) } { (|x| + |y|)^2 } \geq 0.
\]
Similarly, $\pa_{y_1} M \geq 0$. We get $M(x_1, y_1, x_2) \leq M(x_1^u, y_1^u, x_2^l)$.

Since $y_1 \geq x_1$, for $M_1$, we write $M_1 = \f{ (1 - x_1 / y_1)^{1/2}}{ ( |x| / y_1)^{1/2} + (|y| / y_1)^{1/2}}$. Since $|y| / y_1$ is decreasing in $y_1$ \eqref{eq:cos_sin}, $M_1 $ is increasing in $y_1$. Clearly, $M_1(x_1, y_1, x_2)$ is decreasing in $x_2$, and decreasing in $x_1$. We get $M_1 \leq M_1(x_1^l, y_1^u, x_2^l)$. 

For $p \geq 2$ and $b \geq a$, using the convexity of $t^{p-1}$ for $p-1 \geq 1$ and Jensen's inequality, we get 
\[
b^p - a^p = p \int_a^b t^{p-1} dt 
\leq p \f{b-a}{2} (b^{p-1} + a^{p-1}).
\]

Using the above estimates, we get 
\[
 \f{|y|^p - |x|^p}{ |y-x|^{1/2}}
\leq  \f{ |y| - |x|}{ |y-x|^{1/2}}  \cdot \f{p}{2}( |x|^{p-1} + |y|^{p-1}).
\]
We have estimated $\f{ |y| - |x|}{ |y-x|^{1/2}}$ in the above. Clearly, $|x|, |y|$ are increasing in $x_1, x_2, y_1$.

\subsubsection{Piecewise H\"older estimates of weighted functions}

Given $F$ depending on the profile $\bar \om, \bar \th, \bar \uu^N$ and radial weights $f, g$ with leading power $ r^{a}, r^{b} $, respectively, we estimate $\d_i( F f/ g, x, z )$ piecewisely. Using the decomposition \eqref{eq:wg_rad_dec1}, we get 
\[
 F(x) f(r) / g(r) = F(x) f_m(r) / g_m(r) r^{a-b}.
\] 

Using the $C^1$ estimates of $f, g, F, f_m, g_m$, the Leibniz rule \eqref{eq:lei}, and \eqref{eq:hol_prod}, we obtain $C^{1/2}, L^{\inf}$ estimates of $F f /g, F f_m / g_m$. Since the weights $f_m, g_m$ are $C^1$ in our applications (see the discussions around \eqref{eq:wg_rad_dec1}), the $C^1, C^{1/2}$ estimates of $F f_m / g_m$ are not singular. For weights $f, g$ singular near $r=0$, we use the second identity and estimate the singular part $r^{a-b}$. If $a-b =0, 1/2, 1$ or $a-b \geq 2$, we use the $C^{1/2}$ estimates in Section \ref{sec:hol_pow}.

Below, we discuss three special cases $F f_m / g_m |x|^{-1/2}, F/ |x|^p, p=5/2, 2$ with $F(0, y)= 0$.

\vs{0.1in}
\paragraph{\bf{Estimate of $F  f_m / g_m/ |x|^{1/2}$}}
For $F$ odd in $x$, near $r= 0$, we have $F(x) f_m / g_m r^{a-b} \approx x r^{-1/2}$, which is $C^{1/2}$. We perform $C^{1/2}$ estimate on $F(x) r^{-1/2}$, which along with the $C^{1/2}$ estimate of $f_m / g_m$ and \eqref{eq:hol_prod} give the desired estimate. 
 Suppose that $|x| \leq |y|$ with $x_2 = y_2$ for the $C_{x_1}^{1/2}$ estimate and $x_1 = y_1$ for the $C_{x_2}^{1/2}$ estimate. A direct calculation yields 
\[
\bal
& \f{ | \f{ F(y)}{|y|^{1/2}} - \f{ F(x)}{ |x|^{1/2}} | }{ |x_i -y_i|^{1/2}}
\leq  \f{ |F(y) - F(x)|}{ |y|^{1/2} |x_i - y_i|^{1/2}}
+ |F(x)| \f{ | |x|^{-1/2} - |y|^{-1/2} | }{ |x_i - y_i |^{1/2} } 
= \f{ |F(y)  - F(x)|}{|x_i - y_i|} M_3^{1/2} +  \B| \f{F(x)}{ x_1} \B| M_4  M_5 , \\
 & M_3  = \f{ |x_i - y_i|}{ |y|}, \quad 
 M_4 \teq  \f{x_1}{ |x|^{1/2} |y|^{1/2}} , \quad M_5 \teq   \f{ | |y|^{1/2} - |x|^{1/2}| }{ |x_i-y_i|^{1/2}}, \quad T_1 =  \f{ |F(y)  - F(x)|}{|x_i - y_i|}, 
 \quad T_2 = \f{F(x)}{ x_1}.
  \eal
\]

Since $F \in C^1$, we can estimate $|\f{F(y)  - F(x)|}{|x_i - y_i|}$ piecewisely using the method in Appendix {\apppiecederi} in Part II \cite{ChenHou2023b_supp}.
We use the estimate for $M_5$ in Section \ref{sec:hol_pow} and the method in Section \ref{sec:est_f_dx} to estimate $F(x) / x_1$ piecewisely. For $M_3$ with $i=1$,
since $y_1 -x_1 \geq 0, x_2 = y_2$, $M_3$ is decreasing in $x_1, x_2$. Since $M_3 = \f{1 - x_1 / y_1 }{ (1 + (y_2 / y_1)^2)^{1/2}  }$, $M_3$ is increasing in $y_1$. Thus $M_3(x, y) \leq M_3( x_1^l, x_2^l,  y_1^u, x_2^l)$. The estimate for $M_3$ with $i=2$ is similar $M_3(x, y ) \leq M_3( x_1^l, x_2^l, x_1^l, y_2^u)$.

In the $C_{x_1}^{1/2}$ estimate, we get $x_2 = y_2$.  For $M_4$, using \eqref{eq:cos_sin} and $ |y| \geq |x|$, we have two estimates 
\[
M_4 \leq \f{x_1}{|x|} \leq \f{x_1^u}{ |( x_1^u, x_2^l)|}, 
\quad  M_4  = ( \f{x_1}{|x|})^{1/2} ( \f{x_1 }{|y|})^{1/2}
\leq ( \f{ x_1^u}{ |(x_1^u, x_2^l) |} )^{1/2}  (\f{ x_1^u}{ |(y_1^l, x_2^l)|})^{1/2}. 
\]

Near $r = 0$, we have a further improvement. Denote $x_1 / y_1 = \lam \in [0, 1]$. Since $x_1 \leq y_1 \leq |y|, |x| \leq |y|$, from the above estimate and $M M_1$ \eqref{eq:hol_pow_half}, we  bound the upper bound as follows 
\[
\bal
& M_5 = \f{|y|^{1/2} - |x|^{1/2}}{|x-y|^{1/2}} \leq M M_1 \leq M_1 \leq  \f{ (y_1 - x_1)^{1/2}}{ x_1^{1/2} + y_1^{1/2}} = \f{ (1 -\lam)^{1/2}}{1  + \lam^{1/2}},  \ 
M_3 \leq \f{y_1 -x_1}{y_1} = 1- \lam, 
 \\
& T_1 M_3^{1/2} + T_2 M_4 M_5 
\leq T_1 (1- \lam)^{1/2}
+ T_2 \f{x_1^{1/2}}{y_1^{1/2}}  M_5 
\leq (1 - \lam )^{1/2} ( T_1 +  \f{T_2 \lam^{1/2}}{1  + \lam^{1/2}} )
\leq (1 - \lam_l )^{1/2} ( T^u_1 +  \f{\lam_u^{1/2} T_2^u }{1  + \lam_l^{1/2}} ).
\eal
\]
We partition $\lam \in [0, 1]$ and use the above estimate to obtain an uniform upper bound for $\lam \in [0, 1]$.

In the $C_{x_2}^{1/2}$ estimate, we get $x_1 = y_1$ and  $M_4  = ( \f{x_1}{|x|})^{1/2} ( \f{y_1 }{|y|})^{1/2}$, which can be estimated using \eqref{eq:cos_sin}. We also have a direct estimate for $F(x) / x_1 M_4$ using the piecewise bound of $F(x)$
\[
|F(x) / x_1 M_4| \leq |F(x)| |( x_1^l, x_2^l)|^{-1/2} | (y_1^l, y_2^l)|^{-1/2}.
\]
We bound $M_5$ using previous method and then obtain the bound for $F(x)/x_1 M_4 M_5$.

\vs{0.1in}
\paragraph{\bf{Estimate of $F / |x|^p, p =2, 5/2$} }
The power $r^p,p=5/2, 2$ is the leading power in the H\"older weight $\psi_i$. In this case, our function satisfies  $F = O( r^3)$ near $r=0$ and we have piecewise $C^3$ bounds of $F$. To estimate the $C^{1/2}$ norm of $F/ r^{5/2}$, using the above estimate with $g = F / r^2$, we only need to estimate the Lipschitz norm of $g$, and the $L^{\inf}$ norm of $g / x$ and $ F / r^{5/2}$. Using the Taylor expansion and the method in Appendix {\appsingnear} in Part II \cite{ChenHou2023b_supp}, we can control $F / r^{5/2}, F/r^2 / x_1 $. To estimate the $C^1$ bound, we use
\[
 \pa_x( \f{F}{ r^2} ) = \f{ F_x }{r^2} - \f{2 F x}{r^4} ,  \quad \pa_y ( \f{F}{ r^2} ) = \f{\pa_y F}{r^2} - \f{ 2 F y}{r^4},
\]
and estimate two terms separately using the method in Appendix {\appsingnear} in Part II \cite{ChenHou2023b_supp}.

\appendix

\section{Parameters, explicit functions, and basic estimates}

\subsection{Parameters and some functions}\label{app:cutoff_near0}

Recall the following cutoff functions constructed and estimated in Appendix {\seccutoffnear} in Part II \cite{ChenHou2023b_supp}
\beq\label{eq:cutoff_near0_chi3}
\kp(x; \nu_1, \nu_2 ) = \kp_1(\f{x}{ \nu_1} ) (1 - \chi_e( \f{x}{ \nu_2} )  ) ,  \quad \kp_1(x) = \f{1}{1 + x^4}, \quad  \chi_e(x) = \B( 1 + \exp( \f{1}{x} + \f{1}{x-1} ) \B)^{-1}. 
\eeq
with $\chi_e(x) =0$ for $x\leq 0$ and $\chi_e(x) = 1$ for $x \geq 1$. For the cutoff functions in \eqref{eq:uerr_dec1} 
in Section \ref{sec:non_dec_u}, and $\chi_{j, 2}, f_{\chi, j}$, 
for $(x, y) \in \R^2$, we choose 
\beq\label{eq:cutoff_near0_all}
\bal
 \chi_{\bar \e}(x, y) &= \kp(x; \nu_{\bar \e, 1}, \nu_{\bar \e, 2}) 
\kp(y; \nu_{\bar \e, 1}, \nu_{\bar \e, 2}) ,  \quad \nu_{\bar \e, 1 } = 1/192, \quad \nu_{\bar \e,  2} = 3/2 , \\ 
 \chi_{ \hat \e}(x, y) &= \kp_*(x) \kp_*(y) , \quad  \kp_*(x) = \kp(x, 1/3, 3/2) ,
\quad  \chi_{\mf{NF}}(x, y) = \kp(x; 2,10)  \kp(y ; 2, 10 ),  \\
  f_{\chi, 1}  & = \D( \f{x y^3}{6} \chi_{ \mf{NF}}(x, y)  ) ,
\quad f_{\chi, 2} = x y \chi_{\mf{NF}}(x, y), \quad f_{\chi, 3} = \f{x^2}{2} \chi_{ \mf{NF}}(x, y). \\
\eal
\eeq
We use the following parameters for the energy \eqref{energy4}
\beq\label{wg:EE}
\bal
&\tau_1 = 5, \quad \quad \mu_1 = 0.668, \quad \mu_2  = 1.336, \quad \mu_4 = 0.065, \quad \tau_2 = 0.23 ,  \\
& \mu_5 = 76 , \quad
\mu_{51} = 61, \quad   \mu_{52} = 15, \quad  \mu_{6} = 61 ,  \quad \mu_{62} = 35.88, \\
&  \mu_7 =  9.5, \quad  \mu_8 = 4.5 ,  \quad   E_* = 5 \cdot 10^{-6}, \quad  \mu_{h} = \tau_1^{-1}(1, \mu_1, \mu_2), \quad \mu_g = \tau_2(\mu_4, 1, 1).
\eal
\eeq

Recall from Appendix {\appparawg} in Part I \cite{ChenHou2023a_supp} the following H\"older weight
\begin{align}\label{wg:hol}
   & g_{i}(h)   = g_{i0}(h) g_{i0}(1, 0)^{-1}, \quad 
 g_{i0}(h) =  ( \sqrt{ |h_1| + q_{i1} |h_2| } + q_{i3} \sqrt{ |h_2| + q_{i2} |h_1| }   )^{-1},  
 \ g_2(h) \equiv g_3(h) , \notag  \\
 & \vec q_{1, \cdot   }  = ( 0.12, 0.01, 0.25 ) , \quad  \vec q_{2, \cdot } = ( 0.14, 0.005, 0.27 ), \quad 
\vec q_{3, \cdot } = \vec q_{2,\cdot }.
\end{align}

\subsection{Basic piecewise estimates}
Denote by $f_l, f_u$ the lower and upper bound of $f$. We have 
\beq\label{eq:func_intval}
\bal
& (f-g)_{l} = f_l - g_u, \quad (f-g)_u = f_u - g_l , \quad (f+g)_{\g} = f_{\g} + g_{\g},  \\
& (f g)_l = \min( f_l g_l , f_u g_l, f_l g_u, f_u g_u),
\quad (f g)_u = \max( f_l g_l , f_u g_l, f_l g_u, f_u g_u),
\eal
\eeq
where $\g = l, u$. If $ g\geq g_l \geq 0$ or $g \geq g_l >0$, we can simplify the formula for the product 
\beq\label{eq:func_intval2}
\bal
 (f g)_{l}  & =  \min( f_l g_l, f_l g_u),  & \quad ( fg)_u  & = \max( f_u g_l, f_u g_u),  \quad & & g_l \geq 0 , \\
  (f / g)_l & = \min( f_l / g_l, f_l / g_u), & \quad   (f / g)_u & = \max( f_u / g_l, f_u / g_u),
  \quad && g_l > 0.
 \eal
\eeq

Given the piecewise bounds of the derivatives of $f, g$, using the Leibniz rule
\beq\label{eq:lei}
 \pa_x^i \pa_y^j (f g) =  \tts{ \sum\nolim_{k \leq i, l\leq j} } \binom{i}{k} \binom{j}{l} \pa_x^k \pa_y^l f \cdot 
  \pa_x^{i-k} \pa_y^{j-l} g,
\eeq
and \eqref{eq:func_intval}, \eqref{eq:func_intval2}, we can estimate the piecewise bounds of $\pa_x^i \pa_y^j(f g)$.  

We have the linear interpolation error estimates (see, e.g. Appendix {\appsoluest} in Part II 
\cite{ChenHou2023b_supp})
\beq\label{est_2nd}
\bal
 & \max_{ x \in [x_l, x_u] }|f(x)|  \leq \max( |f(x_l) | , |f(x_u)| ) + h^2 || f_{xx}||_{L^{\inf}(I)} / 8, \quad h = x_u - x_l, \\
 & \min_{\al, \b = l, u} f(x_{\al}, y_{\b}) - err_i
  \leq f(x) \leq \max_{\al, \b = l, u} f(x_{\al}, y_{\b}) + err_i  , \quad  i =1,2 ,  \\
& 
 err_1 = \tf12 ( || f_{x}||_{L^{\inf}(Q)} h_1  + || f_{y} ||_{L^{\inf}(Q)} h_2) , \quad
  err_2 = \tf18 (|| f_{xx}||_{L^{\inf}(Q)} h_1^2  + || f_{yy} ||_{L^{\inf}(Q)} h_2^2 )  .
 \eal
\eeq

\bibliographystyle{plain}